# Multi-Order Runge-Kutta Methods

## or how to numerically solve initial value problems of any order

Loris PETRONIJEVIC

## Abstract

When one wishes to numerically solve an initial value problem, it is customary to rewrite it as an equivalent first-order system to which a method, usually from the class of Runge-Kutta methods, is applied. Directly treating higher-order initial value problems without such rewriting, however, allows for significantly greater accuracy. We therefore introduce a new generalization of Runge-Kutta methods, multi-order Runge-Kutta methods, designed to solve initial value problems of arbitrary order. We establish fundamental properties of these methods, including convergence, order of consistency, and linear stability. We also analyze the structure of the system satisfied by the approximations of a method, which enables us to provide a proper definition of explicit methods and gain a finer understanding of implicit methods.

# Introduction

In the current literature on numerical methods for initial value problems, it is almost always assumed that the problem is of order one[1, 2]. This may surprise physicists, who often encounter second-order initial value problems. The usual approach is to note that any initial value problem can be rewritten as an equivalent first-order problem, whose solution can then be approximated. For example, Runge-Kutta methods apply only to first-order problems, and therefore require this rewriting, which results in a loss in accuracy.

If we consider a sufficiently smooth initial value problem of order $n \in \mathbb{N}^*$, a Runge-Kutta method of order of consistency $v$ approximates all derivatives of the solution with an error of order $\underset{h \to 0}{\mathcal{O}}\left(h^{1+v}\right)$. By contrast, a multi-order Runge-Kutta method with the same number of stages typically achieves an error of order $\underset{h \to 0}{\mathcal{O}}\left(h^{n-i+v}\right)$ for the $i^{\text{th}}$ derivative of the solution, which in the worst case $i = n - 1$ is the same as a Runge-Kutta method, but otherwise is substantially more accurate. The computational trade-off is negligible, for each stage and derivative $i$ only the computation of an additional Taylor polynomial of order $n - i - 1$ is required, which is a matter of a handful of multiplications and additions.

Another point we revisit is the usual distinction between explicit and implicit Runge-Kutta methods. Explicit methods are commonly described as those whose Butcher tableau is lower triangular. While this condition is convenient and often sufficient in practice, it does not fully express what distinguishes explicit methods from implicit ones. What truly matters is whether the computation of the stages requires solving a coupled system of equations.

To clarify this, we study the structure of the systems satisfied by the stage approximations, both for classical Runge-Kutta methods and for multi-order methods. This structural viewpoint leads to a more intrinsic understanding of explicitness and implicitness. It also gives a clearer picture of how implicit methods approximate the solution of an initial value problem.

In particular, implicitness does not necessarily mean solving a single large nonlinear system at each step. Depending on the structure of the method, the stage equations may decompose into several smaller systems, or even into independent subsystems that can be treated separately. This makes it possible to design methods whose structure is better suited to modern computational architectures, including parallel implementations.

The existing literature already contains special cases of multi-order Runge-Kutta methods such as the semi-implicit Euler method[3] and Runge-Kutta-Nyström methods[4]. These, however, are almost always restricted to order-two problems, are not well known or systematically studied, and lack a consensus on their definition. Depending on the author, additional assumptions on the initial value problem are often introduced.

For example, there is no generally accepted definition of Runge-Kutta-Nyström methods. Some authors assume the problem is independent of the derivative of the solution, while others assume that the second-order problem was obtained by differentiating a first-order one, thereby considering two differential equations simultaneously. Moreover, papers on these special cases rarely address why it may be advantageous to avoid rewriting the initial value problem into a first-order one.

To avoid any misunderstanding, we emphasize that multi-order Runge-Kutta methods are fundamentally different from Taylor series methods[1]. Taylor series methods apply only to first-order initial value problems, and higher accuracy is achieved by differentiating the differential equation function $f$ several times. In contrast, multi-order Runge-Kutta methods directly apply to higher-order initial value problems and take advantage of the fact that, in such problems, $f$ already represents a higher derivative of the unknown function $y$. This structure allows one to exploit Taylor's theorem more efficiently without ever differentiating $f$, or any function for that matter.

The study of the structural properties of Runge-Kutta methods is similarly niche. The notion of diagonally implicit Runge-Kutta methods, which relies on a particular structure, is well known[1] but remains only a special case of a more general framework developed in this paper. The operation that connects the structure of Runge-Kutta methods with the arrangement of their weights, the permutation operation, does appear in a few rare papers[5], but it has been underestimated. It has not been given a proper name, and its main significance has been overlooked, with most authors preferring adjacent notions such as reducibility[5] or B-series[1].

We also adopt a more general definition of initial value problems than is commonly used. Instead of taking the domain of an initial value problem to be the entire vector space, we only assume it is an open connected set. The reasons for this choice is explained after the formal definition. This approach requires greater rigor and leads to more technical results, but provides a deeper understanding of Runge-Kutta methods and their generalization.

Generalizing Runge-Kutta methods also requires reconstructing, from the ground up, the theory of consistency order and linear stability (including A-stability, L-stability, and related notions) for multi-order Runge-Kutta methods. Our development is strongly inspired by the classical theory of Runge-Kutta methods, but with essential additions. In particular, in the section on consistency we introduce a new discussion of rewriting initial value problems, which plays a central role. In the section on linear stability, we introduce the new concept of half-line stability, which is closely related to the notion of A($\alpha$)-stability for Runge Kutta methods.

Finally, we outline the structure of this paper:

- Section I. The first subsection gives a brief overview of Runge-Kutta methods, presenting the underlying idea that will serve as the basis for their generalization to higher-order problems. The second subsection discusses the issues with rewriting problems to first order.
- Section II. The first subsection introduces the concept of general multi-order Runge-Kutta methods and establishes conditions for the existence of their approximations. The second subsection studies the structure of the system of equations which defines the approximations of a Runge-Kutta method.
- Section III. The first subsection generalizes the notions of order of consistency and convergence. The second subsection shows how to apply multi-order Runge-Kutta methods to arbitrary initial value problems using different rewritings.
- Section IV. The first subsection compares multi-order and classical Runge-Kutta methods in terms of computational cost and order of consistency using numerical experiments. The second subsection generalizes the concepts of A-stability and L-stability, and introduces half-line stability.

# Notations and conventions

- We define $0^0 = 1$.
- We define $\mathbb{N}$ as the set of natural numbers $\{0, 1, 2, ...\}$, and $\mathbb{N}^* = \mathbb{N} \setminus \{0\}$, $\overline{\mathbb{N}} = \mathbb{N} \cup \{+\infty\}$.
- We define $\mathbb{Z}$ as the set of integers $\{... -1, 0, 1...\}$, and $\overline{\mathbb{Z}} = \mathbb{Z} \cup \{-\infty, +\infty\}$.
- We define $\mathbb{R}$ as the set of real numbers, and $\mathbb{R}^* = \mathbb{R} \setminus \{0\}$, $\mathbb{R}_- = \{x \in \mathbb{R} \mid x \leq 0\}$, $\mathbb{R}_+ = \{x \in \mathbb{R} \mid x \geq 0\}$, $\mathbb{R}_+^* = \mathbb{R}^* \cap \mathbb{R}_+$, $\mathbb{R}_-^* = \mathbb{R}^* \cap \mathbb{R}_-$, $\overline{\mathbb{R}} = \mathbb{R} \cup \{-\infty, +\infty\}$.
- We define $\boldsymbol{i}$ as the imaginary unit, $\mathbb{C}$ as the set of complex numbers. For all $z \in \mathbb{C}$, $\mathfrak{R}(z), \mathfrak{I}(z)$ are respectively the real and imaginary parts of $z$.
- Let $(a, b) \in \overline{\mathbb{R}}^2$. We define $[\![a, b]\!] = [a, b] \cap \overline{\mathbb{Z}}$, $]\!]a, b]\!] = ]a, b] \cap \overline{\mathbb{Z}}$, $[\![a, b[\![ = [a, b[ \cap \overline{\mathbb{Z}}$, $]\!]a, b[\![ = ]a, b[ \cap \overline{\mathbb{Z}}$.
- We define $\#(X)$ as the cardinal of a set $X$.
- We define $\mathcal{P}(X)$ as the power set of a set $X$.
- Let $X$, and $Y$ two non-empty sets. We define $\{X \to Y\}$ as the set of functions from $X$ to $Y$.
- Let $E, F$ two normed vector spaces and $k \in \mathbb{N}^*$. We define $D^k(E, F)$ as the set of functions from $E$ to $F$ differentiable $k$ times. We define $C^k(E, F)$ as the set of functions from $E$ to $F$ continuously differentiable $k$ times. A function in this set will be said to be of class $C^k$.
- Let $X$ a set. We define $\mathfrak{S}_X$ as the set of permutations of $X$.
- Let $X$ a set, and $\mathcal{R}$ a binary relation defined on $X$. We define, for all $x \in X$, $\mathcal{R}^{-1}(x) = \{y \in X \mid y\mathcal{R}x\}$. If $\mathcal{R}$ is an equivalence relation then $\mathcal{R}^{-1}(x)$ is the equivalence class of $x$. We define, for all $J \in \mathcal{P}(X)$, $\mathcal{R}^{-1}(J) = \{j' \in X \mid \exists j \in J, j'\mathcal{R}j\} = \bigcup_{j \in J} \mathcal{R}^{-1}(j)$
- Let $(X, d)$ a metric space, $X'$ a non-empty set, $(Y, d')$ a metric space, $U \in \mathcal{P}(X \times X')$, $a \in Y$, and $f \in \{\mathbf{U} \to Y\}$. Let $V = \{x' \in X \mid \exists x \in X', (x, x') \in U\}$, $g \in \{X \to \mathbb{R}\}$, $\tilde{x} \in X$ such that $\{\tilde{x}\} \times V \subset U$. We define a uniform big O condition as :

$$U \to Y, (x, x') \to f(x, x') = a + \underset{\substack{x \to \tilde{x} \\ x \in V}}{\mathcal{O}}(g(x))$$

$$\Leftrightarrow$$

$$\exists (\delta, C) \in \mathbb{R}_+^* \times \mathbb{R}_+^*, \ \forall (x, x') \in U, \ d(x, \tilde{x}) < \delta \Rightarrow d'(f(x, x'), a) \leq C \ |g(x)|$$

- Let $X$ a set, $p \in \mathbb{N}^*$, and $Y_1, ...Y_p$ a family of sets where no set is the empty set. We define $X^{Y_1 \times ... \times Y_p}$ the set of $p$-dimensional sequences of $X$ indexed by the elements of the sets $Y_1, ...Y_p$.
Let $v \in X^{Y_1 \times ... \times Y_p}$, and $(y_1, ...y_p) \in Y_1 \times ... \times Y_p$. We define $v_{y_1, ...y_p}$ as the entry at position $(y_1, ...y_p)$ of $v$.
Let $p' \in [\![1, p-1]\!]$. We define $v_{y_1, ...y_{p'}} \in X^{Y_{p'+1} \times ... \times Y_p}$ as the only $p - p'$-dimensional sequence which satisfies, for all $(y_{p'+1}, ...y_p) \in Y_{p'+1} \times .. \times Y_p$, the entry at position $(y_{p'+1}, ...y_p)$ of $v_{y_1, ...y_{p'}}$ is equal to $v_{y_1, ...y_p}$.
Let $x \in X$, and $j \in (\mathbb{N}^*)^{[\![1, p]\!]}$. We define $x_{j_1, ...j_p} \in X^{[\![1, j_1]\!] \times ... \times [\![1, j_p]\!]}$ as the only sequence with all of its entries equal to $x$.
We define the max norm $\|\cdot\|_{j_1,\dots j_p} \in \left\{\mathbb{R}^{[\![1,j_1]\!]\times\dots\times[\![1,j_p]\!]} \to \mathbb{R}_+\right\}$ as, for all $v \in \mathbb{R}^{[\![1,j_1]\!]\times\dots\times[\![1,j_p]\!]}$ :

$$\|v\|_{j_1,\dots j_p} = \max_{(y_1,\dots y_p)\in[\![1,j_1]\!]\times\dots\times[\![1,j_p]\!]} \left(|v_{y_1,\dots y_2}|\right)$$

We define $I_p$ as the identity matrix of dimension $p$. Let $A$ a matrix. We define $A^T$ as the transpose of $A$. Some block matrices may have in their definition rows or columns of matrices of dimension $0$, this is to avoid distinguishing cases and simply means that we should act as if these columns or rows are not here.

- We define $\otimes$ the Kronecker product as, for all $(m,n,p,q) \in (\mathbb{N}^*)^4$, for all $(A,B) \in \mathbb{C}^{[\![1,m]\!]\times[\![1,n]\!]} \times \mathbb{C}^{[\![1,p]\!]\times[\![1,q]\!]}$ :

$$A \otimes B = \begin{bmatrix} A_{1,1}B & \dots & A_{1,n}B \\ \vdots & \ddots & \vdots \\ A_{m,1}B & \dots & A_{m,n}B \end{bmatrix} \in \mathbb{C}^{[\![1,mp]\!]\times[\![1,nq]\!]}$$

To avoid using parentheses we give a higher priority to the Kronecker product than matrix multiplication. The Kronecker product will only be used as a shortcut to expand matrices :

$$A \otimes I_p = \begin{bmatrix} A_{1,1}I_p & \dots & A_{1,n}I_p \\ \vdots & \ddots & \vdots \\ A_{m,1}I_p & \dots & A_{m,n}I_p \end{bmatrix}$$

# Initial value problem and Taylor's theorem

**Definition : Initial value problem**

Let $(n, d) \in \mathbb{N}^* \times \mathbb{N}^*$, $\mathbf{U} \in \mathcal{P}\left(\mathbb{R} \times \mathbb{R}^{[\![1,n]\!] \times [\![1,d]\!]}\right)$ an open, connected non-empty set, $\left(t_0, y_{n,0}\right) \in \mathbf{U}$, and $f \in \left\{\mathbf{U} \to \mathbb{R}^{[\![1,d]\!]}\right\}$.

- An initial value problem is a tuple $\left(f, t_0, y_{n,0}\right)$.
- $\mathbf{U}$ is called the space time domain.
- $n$ and $d$ are respectively the order and the dimension of the initial value problem.
- $f$ is said to be a differential equation function of order $n$ and dimension $d$.
- $t_0$ is called the initial instant and $y_{n,0}$ is called the initial values.
- We define $\Omega = \left\{t \in \mathbb{R} \mid \exists x \in \mathbb{R}^{[\![1,n]\!] \times [\![1,d]\!]}, (t, x) \in \mathbf{U}\right\}$, the time domain, an open interval of $\mathbb{R}$.

Let $\hat{y} \in D^n\left(\Omega, \mathbb{R}^{[\![1,d]\!]}\right)$.

- $\hat{y}$ is said to be a solution candidate if and only if :

$$\forall t \in \Omega,\ \left(t, \frac{\mathrm{d}^{n-1}\hat{y}}{\mathrm{d}t^{n-1}}(t), \ldots \hat{y}(t)\right) \in \mathbf{U}$$

- $\hat{y}$ is said to satisfy the differential equation induced by $f$ if and only if $\hat{y}$ is a solution candidate and :

$$\forall t \in \Omega,\ \frac{\mathrm{d}^n\hat{y}}{\mathrm{d}t^n}(t) = f\left(t, \frac{\mathrm{d}^{n-1}\hat{y}}{\mathrm{d}t^{n-1}}(t), \ldots \hat{y}(t)\right)$$

- $\hat{y}$ is said to satisfy the initial values if and only if :

$$\forall i \in [\![0, n-1]\!],\ \frac{\mathrm{d}^i\hat{y}}{\mathrm{d}t^i}(t_0) = y_{n,0,n-i}$$

- $\hat{y}$ is said to be a solution of the initial value problem if and only if $\hat{y}$ is a solution candidate, $\hat{y}$ satisfies the differential equation induced by $f$ and $\hat{y}$ satisfies the initial values.

Let $\Omega' \in \mathcal{P}(\mathbb{R})$ an open connected subset of $\Omega$ which contains $t_0$, and $\hat{y}' \in D^n\left(\Omega', \mathbb{R}^{[\![1,d]\!]}\right)$.

- $\hat{y}'$ is said to be a maximal solution if and only if $\hat{y}'$ satisfies the initial value problem restricted to $\Omega'$, and any prolongation of $\hat{y}'$ to $\Omega''$ any open connected subset of $\Omega$ and superset of $\Omega'$, does not satisfy the initial value problem restricted to $\Omega''$.

To simplify the notations we define the jet and extended jet of a function.

**Definition : Jet & Extended jet**

Let $(n, d) \in \mathbb{N}^* \times \mathbb{N}^*$, $\Omega \in \mathcal{P}(\mathbb{R})$ an open set, and $\hat{y} \in D^n\left(\Omega \to \mathbb{R}^{[\![1,d]\!]}\right)$.

- We define $\mathcal{J}^n\hat{y} \in \left\{\Omega \to \mathbb{R}^{[\![1,n+1]\!] \times [\![1,d]\!]}\right\}$ the $n$-jet of $\hat{y}$ as :

$$\mathcal{J}^n\hat{y} = \left(\frac{\mathrm{d}^n\hat{y}}{\mathrm{d}t^n}, \ldots\ \hat{y}\right)$$

- We define $\overline{\mathcal{J}}^n\hat{y} \in \left\{\Omega \to \mathbb{R}^{\{0\} \cup ([\![1,n+1]\!] \times [\![1,d]\!])}\right\}$ the extended $n$-jet of $\hat{y}$ as :

$$\forall t \in \Omega,\ \overline{\mathcal{J}}^n\hat{y}(t) = \left(t, \frac{\mathrm{d}^n\hat{y}}{\mathrm{d}t^n}(t), \ldots\ \hat{y}(t)\right)$$

The jet and extend jet allow us to write compactly the conditions for a function $\hat{y}$ to be a solution.

**Proposition : Initial value problems & Jet**

Let $(f, t_0, y_{n,0})$ an initial value problem of order $n \in \mathbb{N}^*$ and dimension $d \in \mathbb{N}^*$, and $\hat{y} \in D^n(\Omega, \mathbb{R}^{[\![1,d]\!]})$.

- $\hat{y}$ is a solution candidate if and only if, for all $t \in \Omega,\ \overline{\mathcal{J}}^{n-1}\hat{y}(t) \in \mathbf{U}$
- $\hat{y}$ satisfies the differential equation induced by $f$ if and only if $\hat{y}$ is a solution candidate and :

$$\frac{\mathrm{d}^n \hat{y}}{\mathrm{d}t^n} = f \circ \overline{\mathcal{J}}^{n-1}\hat{y}$$

- $\hat{y}$ satisfies the initial values if and only if $\mathcal{J}^{n-1}\hat{y}(t_0) = y_{n,0}$.

If $\hat{y}$ is a solution of an initial value problem we may want to give a name to its extended jet.

**Definition : Solution curve**

Let $(f, t_0, y_{n,0})$ an initial value problem of order $n \in \mathbb{N}^*$ and dimension $d \in \mathbb{N}^*$, and $\hat{y} \in D^n(\Omega, \mathbb{R}^{[\![1,d]\!]})$ a solution of the initial value problem, if it exists.
We define the solution curve of $\hat{y}$ as the function $\overline{\mathcal{J}}^{n-1}\hat{y}$.

We now discuss this definition of an initial value problem, as everything relies on it.

Indexing the initial values by $n - i$ instead of simply $i$ may seem like a strange choice, but this change of variable is very common throughout the paper, for a good reason. The core idea is that we will eventually make $n$ tends towards $+\infty$ and it will then be necessary that the highest derivatives keep a constant numbering. For example, the derivative $n - 1$ of $\hat{y}$ is always numbered using $1$, $n - 2$ is always numbered using $2$, etc… This change of variable simplifies everything and eventually becomes intuitive.

As for $\mathbf{U}$, it needs to be an open set to have well defined derivatives, but we need to justify why $\mathbf{U}$ is assumed to be connected. First off, if $\mathbf{U}$ is connected, the domain $\Omega$ is then an open connected set of $\mathbb{R}$, hence an interval, which is necessary to have a fair definition of an initial value problem. If $\Omega$ is open but not necessarily connected, the initial value problem consists of multiple differential equations, one for each connected component, and only one is an actual initial value problem. The following paragraph proves this statement.

We define $\mathfrak{C} \in \mathcal{P}(\mathbf{U})^X$ the set of connected components of $\mathbf{U}$ indexed by a set $X$. We define $S$ the set of solutions of the initial value problem, and, for all $x \in X$, we define $S'_x$ the set of solution of the differential equation induced by $f$ restricted to $\mathfrak{C}_x$, if $(t_0, y_{n,0}) \in \mathfrak{C}_x$ the solutions must also satisfy the initial values. We now prove that the solutions of the initial value problem are the disjoint unions of $\hat{y}'_x$ with, for all $x \in X,\ \hat{y}'_x \in S'_x$.
Let $\hat{y} \in S$ a solution of the initial value problem. For all $x \in X$, the restriction of $\hat{y}$ to $\mathfrak{C}_x$, $\hat{y}'_x$, is solution of the differential equation induced by $f$ restricted to $\mathfrak{C}_x$ and satisfies the initial values if $(t_0, y_{n,0}) \in \mathfrak{C}_x$, hence $\hat{y}'_x \in S'_x$. Therefore all solutions of the initial value problem are a disjoint union of function $\hat{y}'_x$ such that, for all $x \in X,\ \hat{y}'_x \in S'_x$. We take $\hat{y}'$ such that, for all $x \in X,\ \hat{y}'_x \in S'_x$. We define $\hat{y}$ as the disjoint union of all $\hat{y}'_x$. $\hat{y}$ is a solution candidate, it satisfies the differential equation induced by $f$ on $\Omega$ and satisfies the initial values, it is thus a solution of the initial value problem.

So, if $\Omega$ is not connected we are in fact dealing with one true initial value problem and multiple differential equations. A fair definition of an initial value problem where $\Omega$ is not connected would have one initial value for each connected component, but we could then split the problem into multiple initial value problems, one for each connected component.

Since the only property that interests us is that $\Omega$ is an open interval, we could only assume that $\Omega$ is connected without assuming that $\mathbf{U}$ is connected. This would be pointless however since we can prove that all the solutions lie in the connected component of $\mathbf{U}$ which contains $\left(t_0, y_{n,0}\right)$.

To prove this, we consider a solution of the initial value problem $\hat{\mathrm{y}} \in D^n\left(\Omega, \mathbb{R}^{[\![,1,d]\!]}\right)$. The solution curve of $\hat{\mathrm{y}}$ is continuous, its image of $\Omega$ is hence connected, and since $\hat{\mathrm{y}}$ satisfies the initial values, the solution curve of $\hat{\mathrm{y}}$ contains the point $\left(t_0, y_{n,0}\right)$. The union of the images of $\Omega$ of all the solution curves is therefore connected since they all are connected and share the point $\left(t_0, y_{n,0}\right)$. This implies that the initial value problem is equivalent to one where $\mathbf{U}$ is connected.

Having $\Omega$ and $\mathbf{U}$ connected will also prevent the approximations from jumping from one connected component of $\mathbf{U}$ to another. This would be catastrophic since the approximations either jumped to another connected component of $\Omega$, and in this case we cannot approximate the value of a solution in this connected component using the values of another connected components of $\Omega$, or the approximations stayed in the same component of $\Omega$ but jumped to another connected component of $\mathbf{U}$, and in this case the approximations cannot be a value of the solution of the initial value problem.

We now state the most important theorem of this paper.

**Theorem : Taylor's theorem**

Let $(k, m, d) \in \mathbb{N}^* \times \mathbb{N}^* \times \mathbb{N}^*$, $U \in \mathcal{P}\left(\mathbb{R}^{[\![1,m]\!]}\right)$ a non-empty open set, $V$ a non-empty compact convex subset of $U$, and $(x, x') \in V^2$.

- Let $f \in D^k\left(U, \mathbb{R}^{[\![1,d]\!]}\right)$. We have :

$$f(x') = \sum_{k'=0}^{k} \frac{1}{k'!} \mathrm{d}^{k'} f(x)(x' - x) + \underset{\|x'-x\| \to 0}{o}\left(\|x' - x\|^k\right)$$

- Let $f \in C^{k+1}\left(U, \mathbb{R}^{[\![1,d]\!]}\right)$. We have :

$$f(x') = \sum_{k'=0}^{k} \frac{1}{k'!} \mathrm{d}^{k'} f(x)(x' - x) + \underset{\substack{\|x'-x\| \to 0 \\ x' \in V}}{\mathcal{O}}\left(\|x' - x\|^{k+1}\right)$$

  If $m = 1$, then :

$$f(x') = \sum_{k'=0}^{k} \frac{\left(x' - x\right)^{k'}}{k'} \frac{\mathrm{d}^{k'} f}{\mathrm{d}t^{k'}}(x) + \int_x^{x'} \frac{\left(x' - t\right)^k}{k!} \frac{\mathrm{d}^{k+1} f}{\mathrm{d}t^{k+1}}(t) \, \mathrm{d}t$$

# I Runge-Kutta methods

## I.1 Basic theory of Runge-Kutta methods

In this subsection we derive the basic theory of Runge-Kutta methods from first principles, which we will use as a guideline to generalize Runge-Kutta methods later. Since Runge-Kutta methods can only numerically solve problems of order 1, we assume $n = 1$.

The idea behind Runge-Kutta methods is to transform an initial value problem into an equivalent integration problem, then numerically approximate the value of the integral.

**I.1.1 Proposition : Integral form of an initial value problem**

Let $(f, t_0, y_{1,0,1})$ an initial value problem of order 1 with $f$ continuous, and $\hat{y} \in D^1(\Omega, \mathbb{R}^{[\![1,d]\!]})$ a solution candidate.
$\hat{y}$ is a solution of the initial value problem if and only if :

$$\forall t \in \Omega,\ \hat{y}(t) = y_{1,0,1} + \int_{t_0}^{t} f(z, \hat{y}(z))\,\mathrm{d}z$$

Proof

We prove a more general version in Section II.1.1

If we use the integral form of an initial value problem, replace $t$ with $h = t - t_0$ and substitute in the integral $x = 2\frac{z - t_0}{h} - 1$, we get :

$$\hat{y}(t_0 + h) = y_{1,0,1} + \frac{h}{2}\int_{-1}^{1} f\left(t_0 + \frac{1+x}{2}h, \hat{y}\left(t_0 + \frac{1+x}{2}h\right)\right)\mathrm{d}x$$

We can then use a Gaussian quadrature on the integral, a Gauss-Legendre quadrature[1] to be precise. A numerical quadrature, sometimes abbreviated to quadrature, is a method that numerically approximates the integral of a continuous function $g \in \{[-1, 1] \to \mathbb{R}^{[\![1,d]\!]}\}$. A Gauss-Legendre quadrature[6] $(b, c)$ with $s \in \mathbb{N}^*$ points approximates :

$$\int_{-1}^{1} g(x)\,\mathrm{d}x \approx \sum_{j=1}^{s} b_j g(c_j)$$

Where $(b, c) \in \mathbb{R}^{[\![1,s]\!]} \times \mathbb{R}^{[\![1,s]\!]}$. One could see Gaussian quadratures as a generalization of approximating an integral with rectangles. We now use a quadrature on $\hat{y}(t_0 + h)$.

$$\hat{y}(t_0 + h) \approx y_{1,0,1} + \frac{h}{2}\sum_{j=1}^{s} b_j f\left(t_0 + \frac{1+c_j}{2}h, \hat{y}\left(t_0 + \frac{1+c_j}{2}h\right)\right)$$

We define, for all $j \in [\![1, s]\!]$ :

$$w_{1,s+1,j} = \frac{b_j}{2} \in \mathbb{R},\ \tau_j = \frac{1+c_j}{2} \in \mathbb{R}$$

Hence :

$$\hat{y}(t_0 + h) \approx y_{1,0,1} + h \sum_{j=1}^{s} w_{1,s+1,j} f\big(t_0 + \tau_j h, \hat{y}(t_0 + \tau_j h)\big)$$

We can similarly approximate the values $\hat{y}(t_0 + \tau_j h)$. We note, for all $j \in [\![1, s+1]\!],\ y_{1,j,1}$ the considered approximations. We define $\tau_{s+1} = 1$ to get :

$$\forall j \in [\![1, s+1]\!],\ \hat{y}(t_0 + \tau_j h) \approx y_{1,j,1} = y_{1,0,1} + h \sum_{j'=1}^{s} w_{1,j,j'} f\big(t_0 + \tau_{j'} h, y_{1,j',1}\big)$$

With $w_1 \in \mathbb{R}^{[\![1,s+1]\!] \times [\![1,s]\!]}$. We use the same evaluations of $f$ for all quadratures because any value of $f$ used in one of them can be used by any of the other quadratures, even if it implies that some weights are 0. The vector $\tau$ and matrix $w_1$ entirely define a Runge-Kutta method. Usually these values are stored in what is called a Butcher tableau[7], but we here use a slightly different version.

**I.1.2 Definition : Runge-Kutta methods**

- A Runge-Kutta method with $s \in \mathbb{N}^*$ points is a pair $(\tau, w_1) \in \mathbb{R}^{[\![1,s]\!]} \times \mathbb{R}^{[\![1,s+1]\!] \times [\![1,s]\!]}$, with $\tau$ its nodes, and $w_1$ its weights.
- We define $\mathrm{RK}_s$ as the set of Runge-Kutta methods with $s \in \mathbb{N}^*$ points.

The vocabulary we use is non-standard. What we call a method with $s$ points is usually called a method with $s$ stages. The issue with the term stage is that it used differently by different authors. Sometimes this term is reserved for the $f\big(t_0 + \tau_j h, y_{1,j,1}\big)$ with $j \in [\![1, s]\!]$, sometimes it is used to designate the $y_{1,j,1}$ with $j \in [\![1, s]\!]$ and a distinction is made between the internal stages $y_{1,j,1}$ if $j \in [\![1, s]\!]$, and the external stage $y_{1,s+1,1}$, which implies that there are $s+1$ stages instead of $s$ stages. This situation is confusing, we thus abandon this vocabulary and use the term stage to indicate which value of $j \in [\![1, s+1]\!]$ we are considering. For example, the approximation at stage $j \in [\![1, s+1]\!]$ here designates $y_{1,j,1}$.

We can notice that a $y_{1,j,1}$ may directly or indirectly require its own value to be computed, which implies that we have to solve a system of equations to compute the approximations. To define rigorously the approximations of a Runge-Kutta method we need to introduce two closely related functions.

**I.1.3 Definition : Approximation function & Evaluation function**

Let $s \in \mathbb{N}^*$, $M = (\tau, w) \in \mathrm{RK}_s$, and $f$ a differential equation function of order 1.
We define the existence-uniqueness domain $\mathcal{U}_{1,f} \in \mathcal{P}(\mathbf{U} \times \mathbb{R})$ of $M$ as the set of $\big(t, y_{1,0,1}, h\big)$ such that the following system, called the stage system of $M$, has a unique solution :

$$\forall j \in [\![1, s+1]\!],\ y_{1,j,1} = y_{1,0,1} + h \sum_{j'=1}^{s} w_{1,j,j'} f\big(t + \tau_{j'} h, y_{1,j',1}\big)$$

$$y_{1,s+1,1} \in \mathbb{R}^{[\![1,d]\!]},\ \forall j \in [\![1, s]\!],\ \big(t + \tau_j h, y_{1,j,1}\big) \in \mathbf{U}$$

We define the approximation function $Y_{1,f} \in \big\{\mathcal{U}_{1,f} \to \mathbb{R}^{[\![1,d]\!]}\big\}^{[\![0,s+1]\!] \times \{1\} \times [\![1,d]\!]}$ and the evaluation function $F_{1,f} \in \big\{\mathcal{U}_{1,f} \to \mathbb{R}^{[\![1,d]\!]}\big\}^{[\![1,s]\!] \times [\![1,d]\!]}$ of $M$, as :

$$\forall j \in [\![0, s+1]\!],\ Y_{1,f,j,1}\big(t, y_{1,0,1}, h\big) = y_{1,j,1}$$

$$\forall j \in [\![1, s]\!],\ F_{1,f,j}\big(t, y_{1,0,1}, h\big) = f\big(t + \tau_j h, y_{1,j,1}\big)$$

There exists a class of Runge-Kutta methods that never require solving a system of equations, explicit Runge-Kutta methods.

**I.1.4 Definition : Explicit methods & Implicit methods**

Let $s \in \mathbb{N}^*$, and $M = (\tau, w_1) \in \mathrm{RK}_s$.
M is said to be explicit if and only if $w_1$ is strictly lower triangular, otherwise $M$ is said to be implicit.

In the general case however we are dealing with a system of equations.

Since $t$ and $y_{1,0,1}$ are usually constants given by an initial value problem, it is useful to define the set of $h \in \mathbb{R}$ such that $\left(t, y_{1,0,1}, h\right) \in \mathcal{U}_{1,f}$.

**I.1.5 Definition: Step sizes**

Let $s \in \mathbb{N}^*$, $M \in \mathrm{RK}_s$, $f$ a differential equation function of order 1, $\mathcal{U}_{1,f}$ the existence-uniqueness domain of $M$, and $\left(t, y_{1,0,1}\right) \in \mathbf{U}$.
We define $\mathcal{H}_{1,f,t,y_{1,0,1}}$ the set of step sizes at the point $\left(t, y_{1,0,1}\right)$ of $M$ as :

$$\mathcal{H}_{1,f,t,y_{1,0,1}} = \left\{h \in \mathbb{R} \mid \left(t, y_{1,0,1}, h\right) \in \mathcal{U}_{1,f}\right\}$$

Since implicit methods need to solve a system of equations, we need to make sure the system has a unique solution. To do so we introduce the concept of Lipschitz continuity.

**I.1.6 Definition : Lipschitz continuity**

Let $p \in \mathbb{N}^*$, $(X_k, d_k)_{k \in [\![1,p+1]\!]}$ a sequence of metric spaces, $U \in \mathcal{P}\left(X_1 \times \ldots \times X_p\right)$ an open non-empty set, $f \in \left\{U \to X_{p+1}\right\}$, and $k \in [\![1,p]\!]$.

- $f$ is said to be globally Lipschitz continuous in its $k^{\text{th}}$ variable if and only if there exists $L \in \mathbb{R}_+$ called a Lipschitz constant such that :

$$\forall x \in U,\ \forall y \in X_k, \left(x_1, \ldots, x_{k-1}, y, x_{k+1}, \ldots, x_p\right) \in \mathbf{U}$$
$$d_{p+1}\left(f\left(x_1, \ldots, x_p\right), f\left(x_1, \ldots, x_{k-1}, y, x_{k+1}, \ldots, x_p\right)\right) \leq L d_k(x_k, y)$$

  $f$ is said to be a contraction mapping in its $k^{\text{th}}$ variable if and only if $L < 1$.
- $f$ is said to be locally Lipschitz continuous in its $k^{\text{th}}$ variable if and only if, for all $x \in U$, there exists a neighborhood of $x$ contained in $U$ such that $f$ restricted to $U$ is globally Lipschitz continuous in its $k^{\text{th}}$ variable.

The following theorem is necessary to use locally Lipschitz continuous functions.

**I.1.7 Theorem : Local Lipschitz continuity & Compact sets**

Let $p \in \mathbb{N}^*$, $(X_k, d_k)_{k \in [\![1,p+1]\!]}$ a sequence of metric spaces, $U \in \mathcal{P}\left(X_1 \times \ldots X_p\right)$ an open set, and $f \in \left\{U \to X_{p+1}\right\}$.
If $f$ is locally Lipschitz continuous in its $k^{\text{th}}$ variable then, for all compact sets contained in $U$, $f$ restricted to the compact set is globally Lipschitz continuous in its $k^{\text{th}}$ variable.

A wide range of functions are Lipschitz continuous.

**I.1.8 Theorem : Lipschitz continuity & Functions of class $C^1$**

Let $(m, d) \in \mathbb{N}^* \times \mathbb{N}^*$, $U \in \mathcal{P}\left(\mathbb{R}^{[\![1,m]\!]}\right)$ an open set, and $f \in \left\{U \to R^{[\![1,d]\!]}\right\}$.

- If $f$ is of class $C^1$ then $f$ is locally Lipschitz continuous.
- If $f$ is of class $C^1$ and its partial derivatives are bounded, then $f$ is globally Lipschitz continuous.

The main interest of this notion is the study of the fixed points of a function.

### I.1.9 Definition : Fixed points

Let $U, V$ two sets such that $U \cap V \neq \emptyset$, and $f \in \{U \to V\}$.
$x \in U \cap V$ is called a fixed point of $f$ if and only if $f(x) = x$.

The most know theorem on the fixed points of a function is the following[8] :

### I.1.10 Theorem : Contraction mapping theorem

Let $(X, d)$ a non empty complete metric space, and $f \in \{X \to X\}$ a contraction mapping.
$f$ has a unique fixed point $x \in X$. For all $y \in X$, $x = \lim_{k \to +\infty} f^k(y)$.

Using this theorem, we prove that, under some conditions on $f$ and the step size, the approximations are defined.

### I.1.11 Theorem : Existence and uniqueness of approximations

Let $(s, k) \in \mathbb{N}^* \times \mathbb{N}$, $M \in \mathrm{RK}_s$, $f$ a differential equation function of order 1 and class $C^k$, locally Lipschitz continuous in its second variable, $Y_{1,f}$ the approximation function of $M$, $F_{1,f}$ the evaluation function of $M$, and $\mathcal{H}_{1,f}$ the set of step sizes of $M$.

1. For all $\left(t, y_{1,0,1}\right) \in \mathbf{U}$, there exists $h^* \in \mathbb{R}_+^*$ such that $]-h^*, h^*[ \subset \mathcal{H}_{1,f,t,y_{1,0,1}}$ and the two following functions are of class $C^k$ :

$$]-h^*, h^*[ \to \mathbb{R}^{[\![0,s+1]\!] \times \{1\} \times [\![1,d]\!]}, h \to Y_{1,f}\left(t, y_{1,0,1}, h\right)$$
$$]-h^*, h^*[ \to \mathbb{R}^{[\![1,s]\!] \times [\![1,d]\!]}, h \to F_{1,f}\left(t, y_{1,0,1}, h\right)$$

1. If $\mathbf{U} = \Omega \times \mathbb{R}^{[\![1,d]\!]}$ and $f$ is globally Lipschitz continuous in its second variable, there exists $h^* \in \mathbb{R}_+^*$ such that, for all $\left(t_0, y_{1,0,1}\right) \in \mathbf{U}$ :

$$\left\{h \in ]-h^*, h^*[ \ \mid \forall j \in [\![1,s]\!],\ t_0 + \tau_j h \in \Omega\right\} \subset \mathcal{H}_{1,f,t,y_{1,0,1}}$$

And the two following functions are of class $C^k$ :

$$\left\{h \in ]-h^*, h^*[ \ \mid \forall j \in [\![1,s]\!],\ t_0 + \tau_j h \in \Omega\right\} \to \mathbb{R}^{[\![0,s+1]\!] \times \{1\} \times [\![1,d]\!]}, h \to Y_{1,f}\left(t, y_{1,0,1}, h\right)$$
$$\left\{h \in ]-h^*, h^*[ \ \mid \forall j \in [\![1,s]\!],\ t_0 + \tau_j h \in \Omega\right\} \to \mathbb{R}^{[\![1,s]\!] \times [\![1,d]\!]}, h \to F_{1,f}\left(t, y_{1,0,1}, h\right)$$

Proof

We prove a more general version in Section II.1.6.

Now that the approximations of a Runge-Kutta methods are rigorously defined, we need to know how well a method approximates the solution of an initial value problem. To measure how close the approximation is to the solution, we first need to make sure the initial value problem we are considering has a unique solution. To do so we use the Cauchy-Lipschitz theorem, also known as the Picard-Lindelöf theorem[9].

### I.1.12 Theorem : Cauchy-Lipschitz theorem

Let $\left(f, t_0, y_{1,0,1}\right)$ an initial value problem of order 1 with $f$ continuous in its first variable and locally Lipschitz continuous in its second variable.

- There exists a unique maximal solution $\hat{y}$ to the initial value problem.
- If $\mathbf{U} = \Omega \times \mathbb{R}^{[\![1,d]\!]}$ and $f$ is globally lipschitz continuous in its second variable, then there exists a unique solution to the initial value problem.

We almost never use initial value problems with unique maximal solutions in definitions and theorems, mainly ones with unique solutions. Though it may seem more general to consider initial value problems with a unique maximal solution, we can simply restrict the time domain of the initial value problem to retrieve an initial value problem with a unique solution.

**I.1.13 Definition : Error of a Runge-Kutta method**

Let $s \in \mathbb{N}^*$, $M \in \mathrm{RK}_s$, $\left(f, t_0, y_{1,0,1}\right)$ an initial value problem with a unique solution $\hat{y}$, $Y_{1,f}$ the approximation function of $M$, and $\mathcal{U}_{1,f}$ the existence-uniqueness domain of $M$.
We define the error function $\varepsilon_{1,f,s+1,1} \in \left\{\left\{\left(t, y'_{1,0,1}, h\right) \in \mathcal{U}_{1,f} \mid t+h \in \Omega\right\} \to \mathbb{R}^{[\![1,d]\!]}\right\}$ of $M$ as, for all $\left(t, y'_{1,0,1}, h\right) \in \mathcal{U}_{1,f}$ such that $t+h \in \Omega$ :

$$\varepsilon_{1,f,s+1,1}\left(t, y'_{1,0,1}, h\right) = \hat{y}(t+h) - Y_{1,f,s+1,1}\left(t, y'_{1,0,1}, h\right)$$

We deduce that how well a method approximates $\hat{y}$ depends on the size of the error. Since no explicit general formula for the error is available, we can only estimate its behavior. Since $t$ and $y'_{1,0,1}$ are usually given constants, the only variable is $h$.

We therefore study the error as a function of h and apply Taylor's theorem in this variable to estimate how rapidly the error tends towards 0 as $h \to 0$. This approach is intuitively justified by the fact that the error vanishes at $h = 0$, and increases as the step size grows, reflecting the increasing difficulty of approximating the solution far away.

**I.1.14 Definition : Order of consistency**

Let $s \in \mathbb{N}^*$, $M \in \mathrm{RK}_s$, $\mathcal{H}$ the set of step sizes of $M$, and $\varepsilon_1$ the error function of $M$.
We define the order of consistency $v_{1,s+1,1} \in \overline{\mathbb{N}}$ of $M$ as the supremum of the set of $k \in \mathbb{N}$ such that, for all initial value problems $\left(f, t_0, y_{1,0,1}\right)$ of order 1 with a unique solution and $f$ of class $C^{\max(1,k)}$, for all $\Omega' \in \mathcal{P}(\Omega)$ any compact subset of $\Omega$, we have :

$$\left\{(t', h) \in \Omega' \times \mathbb{R} \mid h \in \mathcal{H}_{1,f,t',\hat{y}(t')}, t' + h \in \Omega\right\} \to \mathbb{R}^{[\![1,d]\!]}$$
$$(t', h) \to \varepsilon_{q,f,s+1,1}(t', \hat{y}(t'), h) = \underset{\substack{h \to 0 \\ t' \in \Omega'}}{\mathcal{O}}\left(h^{k+N}\right)$$

The higher the order of consistency, the better the method approximates the solution of the initial value problem for small step sizes.

To know the conditions on the weights for a certain order of consistency, one must use Taylor's theorem to develop $\hat{y}(t_0 + h)$ and the evaluations of $f$ in $y_{1,s+1,1}$ up to the desired order and make the coefficients cancel each other. A high order of consistency implies that the weights satisfy many conditions, the method hence needs many weights, which implies a high number of stages, which makes the method more computationally intensive.

This approach works only for small step sizes, if we wish to use a large step size, the concept of order of consistency is not useful. One way to avoid this issue is to split a step size $h'$ into multiple smaller ones $h_0, \dots h_{\tilde{q}-1}$ which satisfy $h' = \sum_{q=0}^{\tilde{q}-1} h_q$. We can then carry out a first step with $\left(t_0, h_0, y_{1,0,1}\right)$ to get the approximation $y_{1,s+1,1}$ at $t_0 + h_0$, we can then carry out another step using $\left(t_0 + h_0, h_1, y_{1,s+1,1}\right)$ to get another approximation $y_{1,2(s+1),1}$, this time at $t_0 + h_0 + h_1$, and repeat this process until $t_0 + h'$ is eventually reached.

**I.1.15 Definition : Step size sequences & Time mesh**

Let $s \in \mathbb{N}^*$, $M \in \mathrm{RK}_s$, $f$ a differential equation function of order 1, $Y_{1,f}$ the approximation function of $M$, $\mathcal{U}_{1,f}$ the existence-uniqueness domain of $M$, $\tilde{q} \in \mathbb{N}^*$, and $(t_0, y_{1,0,1}) \in \mathbf{U}$.

- $h \in \mathbb{R}^{[\![0,\tilde{q}-1]\!]}$ is said to be a step size sequence of length $\tilde{q}$ if and only if, if we note $t \in \mathbb{R}^{[\![0,\tilde{q}]\!]}$ the sequence defined as for all $q \in [\![0, \tilde{q}-1]\!]$, $t_{q+1} = t_q + h_q$, then the sequence $\left(y_{1,q(s+1),1}\right)_{q \in [\![1,\tilde{q}]\!]}$ defined as, for all $q \in [\![0, \tilde{q}-1]\!]$ :

$$y_{1,(q+1)(s+1),1} = Y_{1,f,s+1,1}\left(t_q, y_{1,q(s+1),1}, h_q\right)$$

  Exists, meaning, for all $q \in [\![0, \tilde{q}-1]\!]$, $\left(t_q, y_{1,q(s+1),1}, h_q\right) \in \mathcal{U}_{1,f}$. $t$ is called the time mesh of $(t_0, h)$.
- We define $\mathcal{H}_{\tilde{q},f,t,y_{1,0,1}} \subset \mathcal{P}\left(\mathbb{R}^{[\![0,\tilde{q}-1]\!]}\right)$ as the set of step size sequences of length $\tilde{q} \in \mathbb{N}^*$ at the point $\left(t, y_{1,0,1}\right)$.
- A sequence $h \in \mathbb{R}^{\mathbb{N}}$ is an infinite step size sequence at the point $\left(t, y_{1,0,1}\right)$ if and only if all its truncated sequences are step size sequences at the point $\left(t, y_{1,0,1}\right)$. We define the time mesh of $(t_0, h)$ as the sequence $t \in \mathbb{R}^{\mathbb{N}}$ defined as, for all $q \in \mathbb{N}$, $t_{q+1} = t_q + h_q$
- We define $\mathcal{H}_{\infty,f,t,y_{1,0,1}} \in \mathcal{P}\left(\mathbb{R}^{\mathbb{N}}\right)$ as the set of infinite step size sequences at the point $\left(t, y_{1,0,1}\right)$.

To easily use step size sequences, we expand the notations we previously defined.

**I.1.16 Definition : Prolongation of the approximation function & Evaluation function**

Let $s \in \mathbb{N}^*$, $M \in \mathrm{RK}_s$, $f$ a differential equation function of order 1, $Y_{1,f}$ the approximation function of $M$, $F_{1,f}$ the evaluation function of $M$, $\tilde{q} \in \overline{\mathbb{N}}^*$, $\left(t_0, y_{1,0,1}\right) \in \mathbf{U}$, $h \in \mathcal{H}_{\tilde{q},f,t_0,y_{1,0,1}}$ with $\mathcal{H}$ the set of step sizes of $M$, and $t$ the time mesh of $(t_0, h)$.

We define, for all $q \in \mathbb{N}$ such that $q < \tilde{q}$ :

$$\forall j \in [\![1, s+1]\!],\ Y_{1,f,q(s+1)+j,1}\left(t_0, y_{1,0,1}, h\right) = Y_{1,f,j,1}\left(t_q, Y_{1,f,q(s+1),1}\left(t_0, y_{1,0,1}, h\right), h_q\right)$$

$$\forall j \in [\![1, s]\!],\ F_{1,f,q(s+1)+j}\left(t_0, y_{1,0,1}, h\right) = f\left(t_q + \tau_j h, Y_{1,f,q(s+1)+j,1}\left(t_0, y_{1,0,1}, h\right)\right)$$

$$\varepsilon_{1,f,q(s+1),1}\left(t_0, y_{1,0,1}, h\right) = \hat{y}\left(t_q\right) - Y_{1,f,q(s+1),1}\left(t_0, y_{1,0,1}, h\right) \in \mathbb{R}^{[\![1,d]\!]}$$

If $\tilde{q} < +\infty$ and $t_{\tilde{q}} \in \Omega$, we define :

$$\varepsilon_{1,f,\tilde{q}(s+1),1}\left(t_0, y_{1,0,1}, h\right) = \hat{y}\left(t_{\tilde{q}}\right) - Y_{1,f,\tilde{q}(s+1),1}\left(t_0, y_{1,0,1}, h\right) \in \mathbb{R}^{[\![1,d]\!]}$$

**I.1.17 Proposition : Use of RK methods**

Let $s \in \mathbb{N}^*$, $M = (\tau, w) \in \mathrm{RK}_s$, $f$ differential equation function of order 1, $\tilde{q} \in \overline{\mathbb{N}}^*$, $\left(t_0, y_{1,0,1}\right) \times \mathbf{U}$, $h \in \mathcal{H}_{\tilde{q},f,t_0,y_{1,0,1}}$ with $\mathcal{H}$ the set of step sizes of $M$, and $t$ the time mesh of $(t_0, h)$.

$M$ approximates, for all $q \in \mathbb{N}$ such that $q < \tilde{q}$ :

$$\forall j \in [\![1, s+1]\!],\ y_{1,q(s+1)+j,1} = y_{1,q(s+1),1} + h_q \sum_{j'=1}^{s} w_{1,j,j'} f\left(t_q + \tau_{j'} h_q, y_{1,q(s+1)+j',1}\right)$$

$$y_{1,(q+1)(s+1),1} \in \mathbb{R}^{[\![1,d]\!]},\ \forall j \in [\![1, s]\!],\ \left(t + \tau_j h, y_{1,q(s+1)+j,1}\right) \in \mathbf{U}$$

We have to question whether or not breaking a large step size into a step size sequence actually works. We therefore introduce the most fundamental property of a Runge-Kutta method, convergence. If the considered initial value problem has a unique solution, we want the approximations to converge towards the solution as we make the time mesh finer, in a similar fashion to how Riemann's integral is defined. This implies that if we had infinite computational power, we would be able to compute the exact solution.

To make the time mesh finer we consider $(t_0, T) \in \Omega^2$ and divide $T - t_0$ into $k \in \mathbb{N}^*$ equal parts, we hence use step sizes of the form $\frac{T-t_0}{k}$.

**I.1.18 Definition : Convergence**

Let $s \in \mathbb{N}^*$, $M \in \mathrm{RK}_s$, and $\varepsilon_1$ the error function of $M$.
$M$ is said to be convergent if and only if, for all initial value problem $(f, t_0, y_{1,0,1})$ of order 1 with a unique solution such that $\mathbf{U} = \Omega \times \mathbb{R}^{[\![1,d]\!]}$ and $f$ is continuous and globally Lipschitz continuous in its second variable, for all $T \in \Omega$ :

$$\lim_{k\to\infty} \max_{q\in[\![0,k]\!]} \|\varepsilon_{1,f,q(s+1),1}\left(t_0, \left(\frac{T-t_0}{k}\right)_{q'\in[\![0,k-1]\!]}, y_{1,0,1}\right)\|_d = 0$$

This criterion is the most important one, a method which is not convergent is not useful since it does not converge towards the solution.

**I.1.19 Theorem : Characterization of convergence**

Let $s \in \mathbb{N}^*$, and $M = (\tau, w_1) \in \mathrm{RK}_s$.
$M$ is convergent if and only if :

$$\sum_{j=1}^{s} w_{1,s+1,j} = 1$$

Proof

We prove a more general version in Section III.1.36.

Convergence and a high order of consistency are important properties of a method, but with certain difficult problems (stiff problems), it is not enough. One category of initial value problem that is interesting to study are homogeneous linear differential equations with constant coefficient, thus initial value problems $(f, t_0, y_{1,0,1})$ with :

$$\lambda \in \mathbb{C},\ \forall (t,x) \in \mathbb{R} \times \mathbb{C},\ f(t,x) = \lambda x$$

The interest of this category of initial value problems stems from the linearization of the differential equation function of an initial value problem $(f, t_0, y_{1,0,1})$, with $f$ differentiable. Let $\hat{y}$ a solution candidate. Using Taylor's theorem, we approximate for small $h$ :

$$f(t_0+h, \hat{y}(t_0+h)) \approx f(t_0, \hat{y}(t_0)) + \frac{\partial f}{\partial t}(t_0, \hat{y}(t_0))h + \begin{bmatrix} \frac{\partial f_1}{\partial x_{1,1}}(t_0) & \dots & \frac{\partial f_1}{\partial x_{1,d}}(t_0) \\ \frac{\partial f_d}{\partial x_{1,1}}(t_0) & \dots & \frac{\partial f_d}{\partial x_{1,d}}(t_0) \end{bmatrix} (\hat{y}(t_0+h) - \hat{y}(t_0))$$

We thus may be interested in initial value problems $(f_l, t_0, y_{1,0,1})$ of the form :

$$\forall (t,x) \in \mathbb{R} \times \mathbb{R}^{[\![1,d]\!]},\ f_l(t,x) = C_1 + C_2 t + Ax$$

With $(C_1, C_2, A) \in \mathbb{R}^{[\![1,d]\!]} \times \mathbb{R}^{[\![1,d]\!]} \times \mathbb{R}^{[\![1,d]\!]^2}$. Classical theory tells us that the solution of this initial value problem is of the form $\hat{y}_h + \hat{y}_p$, where $\hat{y}_p$ is a solution of the differential equation induced by $f_l$, and $\hat{y}_h$ is the solution of the initial value problem $(f_h, t_0, y_{1,0,1} - \hat{y}_p(t_0))$, where, for all $(t,x) \in \mathbb{R} \times \mathbb{R}^{[\![1,d]\!]}$, $f_h(t,x) = Ax$. Assume that $A$ is a diagonalizable matrix. There hence exists $P \in \mathbb{C}^{[\![1,d]\!]^2}$ invertible, and $B \in \mathbb{C}^{[\![1,d]\!]^2}$ a diagonal matrix such that $A = P^{-1}BP$. Let $\hat{y} \in D^1(\mathbb{R}, \mathbb{R}^{[\![1,d]\!]})$. We have, for all $t \in \mathbb{R}$ :

$$\frac{\mathrm{d}\hat{y}}{\mathrm{d}t}(t) = f_h(t, \hat{y}(t)) \Leftrightarrow \frac{\mathrm{d}\hat{y}}{\mathrm{d}t}(t) = P^{-1}BP\hat{y}(t) \Leftrightarrow \frac{\mathrm{d}P\hat{y}}{\mathrm{d}t}(t) = BP\hat{y}(t)$$

The entries of $P\hat{y}$ form a system of independent linear initial value problem $\left(f_{P,k}, t_0, y'_{k,1,0,1}\right)$, with $k \in [\![1, d]\!]$. Since this system was obtained by linearizing $f$, its solution indicates how the entries of $P\hat{y}$ behave locally. We wish for this behavior to be reflected in the approximations of Runge-Kutta methods. The local behavior of $\left(f, t_0, y_{1,0,1}\right)$ dictates how the solutions of $\left(f_{P,k}, t_0, y'_{k,1,0,1}\right)$ behave asymptotically, and we can observe how a Runge-Kutta method reacts by studying how it approximates the solution of each $\left(f_{P,k}, t_0, y'_{k,1,0,1}\right)$, in particular the asymptotical behavior of its sequence of approximations.

We only consider constant step size sequences since we are only interested in the behavior of the solution for a single step. Let $h$ a step size. The last point of the mesh tends towards $+\infty$ if and only if $h > 0$, and it tends towards $-\infty$ if and only if $h < 0$.

The solution of a linear initial value problem is $t :\to y_{1,0,1}e^{\lambda(t-t_0)}$. We are interested in two behaviors, when the solution converges, and when the solutions stays bounded.

If the solution converges at one infinity then it converges towards 0. The solution tends towards 0 as $t$ tends towards $+\infty$ for all initial values if and only if $\Re(\lambda) < 0$, and the solution converges towards 0 as $t$ tends towards $-\infty$ for all initial values if and only if $\Re(\lambda) > 0$. Satisfying one of the two cases is equivalent to $\Re(h\lambda) < 0$. The solution stays bounded as $t$ tends towards $+\infty$ for all initial values if and only if $\Re(\lambda) \leq 0$, and the solution stays bounded as $t$ tends towards $-\infty$ for all initial values if and only if $\Re(\lambda) \geq 0$. Satisfying one of the two cases is equivalent to $\Re(h\lambda) \leq 0$.

We want the approximations of a Runge-Kutta method to stay bounded when the solution stays bounded and we want them to converge towards 0 when the solution converges towards 0.

**I.1.20 Definition : A-stability & Absolute A-stability**

Let $s \in \mathbb{N}^*$, $M \in \mathrm{RK}_s$, and $Y_1$ the approximation function of $M$.

- $M$ is said to be A-stable if and only if, for all infinite constant step size sequence $(h)_{q\in\mathbb{N}}$ with $h \in \mathbb{R}^*$, for all initial value problem $\left(f, t_0, y_{1,0,1}\right)$ of order 1 of the form :

$$\lambda \in \mathbb{C},\ \forall (t,x) \in \mathbb{R} \times \mathbb{C},\ f(t,x) = \lambda x$$

  If $\Re(h\lambda) \leq 0$, then $\left(Y_{1,f,q(s+1),1}\left(t_0, y_{1,0,1}, (h)_{q'\in\mathbb{N}}\right)\right)_{q\in\mathbb{N}}$ is bounded.

- $M$ is said to be absolute A-stable if and only if, for all infinite constant step size sequence $(h)_{q\in\mathbb{N}}$ with $h \in \mathbb{R}^*$, for all initial value problem $\left(f, t_0, y_{1,0,1}\right)$ of order 1 of the form :

$$\lambda \in \mathbb{C},\ \forall (t,x) \in \mathbb{R} \times \mathbb{C},\ f(t,x) = \lambda x$$

  If $\Re(h\lambda) < 0$, then $\left(Y_{1,f,q(s+1),1}\left(t_0, y_{1,0,1}, (h)_{q'\in\mathbb{N}}\right)\right)_{q\in\mathbb{N}}$ tends towards 0 as $q$ tends towards $+\infty$.

We now derive the exact expression of the approximations in order to study their convergence and boundedness. To do so, we first introduce the stability function of a method.

**I.1.21 Definition : Stability function**

Let $s \in \mathbb{N}^*$, and $M = (\tau, w_1) \in \mathrm{RK}_s$. We define :

$$E_s = \begin{bmatrix} I_s & 0_{s,1} \end{bmatrix} \in \mathbb{R}^{[\![1,s]\!]\times[\![1,s+1]\!]}$$

We define the stability function $\tilde{R}_1 \in \left\{\tilde{\mathcal{S}}_1 \to \mathbb{C}\right\}$ of $M$ with $\tilde{\mathcal{S}}_1 \in \mathcal{P}(\mathbb{C})$ its domain of definition, as :

$$\forall z \in \tilde{\mathcal{S}}_1,\ \tilde{R}_1(z) = 1 + z w_{1,s+1}^T \left(I_s - zE_s w_1\right)^{-1} 1_{s,1}$$

$$\tilde{\mathcal{S}}_1 = \{z \in \mathbb{C} \mid \det(I_s - zE_s w_1) \neq 0\}$$

### I.1.22 Proposition : Approximation of linear initial value problems

Let $s \in \mathbb{N}^*$, $M \in \mathrm{RK}_s$, $\tilde{R}_1$ the stability function of $M$, $\tilde{\mathcal{S}}_1$ its domain of definition, $f$ a linear differential equation function of order 1, with :

$$\lambda \in \mathbb{C},\ \ \forall(t,x) \in \mathbb{R} \times \mathbb{C},\ \ f(t,x) = \lambda x$$

And $\mathcal{U}_{1,f}$ the existence-uniqueness domain of $M$.
Then $\mathcal{U}_{1,f} = \left\{ \left(t, y_{1,0,1}, h\right) \in \mathbb{R} \times \mathbb{C} \times \mathbb{R} \mid h\lambda \in \tilde{\mathcal{S}}_1 \right\}$, and, for all $\left(t, y_{1,0,1}, h\right) \in \mathcal{U}_{1,f}$ :

$$Y_{1,f,s+1,1}\left(t, y_{1,0,1}, h\right) = \tilde{R}_1(h\lambda) y_{1,0,1}$$

Proof

We prove a more general version in Section IV.2.9.

We then use this expression to characterize A-stability.

### I.1.23 Theorem : Characterization of A-stability & Absolute A-stability

Let $s \in \mathbb{N}^*$, $M \in \mathrm{RK}_s$, $\tilde{R}_1$ the stability function of $M$, and $\tilde{\mathcal{S}}_1$ its domain of definition.

- $M$ is absolute A-stable if and only if, for all $z \in \tilde{\mathcal{S}}_1$, if $\mathfrak{R}(z) < 0$ then $|\tilde{R}_1(z)| < 1$.
- $M$ is A-stable if and only if, for all $z \in \mathcal{S}_1$, if $\mathfrak{R}(z) \leq 0$ then $|\tilde{R}_1(z)| \leq 1$.

Proof

We prove a more general version in Section IV.2.10.

A classical result is that explicit methods cannot be A-stable.

### I.1.24 Theorem : Second Dahlquist barrier

No explicit convergent Runge Kutta method is either A-stable or absolute A-stable.

Proof

We prove a more general version in Section IV.2.22

Although A-stability is already a strong requirement, it can be strenghtened further. The solution of a linear initial value problem is $t :\rightarrow y_{1,0,1} e^{\lambda(t-t_0)}$, hence as the real part of $\lambda$ tends towards $-\infty$, the solution converges towards 0 increasingly rapidly. If we want this behavior to show in the approximations of a Runge-Kutta method, we may use the notion of L-stability.

### I.1.25 Definition : L-stability & Absolute L-stability

Let $s \in \mathbb{N}^*$, $M \in \mathrm{RK}_s$, $\tilde{R}$ the stability function of $M$, and $\tilde{\mathcal{S}}$ its domain of definition.

- $M$ is said to be L-stable if and only if $M$ is A-stable and $\lim_{\mathfrak{R}(z) \rightarrow -\infty} \tilde{R}_1(z) = 0$.
- $M$ is said to be absolute L-stable if and only if $M$ is absolute A-stable and $\lim_{\mathfrak{R}(z) \rightarrow -\infty} \tilde{R}_1(z) = 0$.

Since A-stability or L-stability can prove to be somewhat of a heavy requirement, there exists a weaker, finer property, called A($\alpha$)-stability[1]. Instead of requiring the method to be stable for any $h\lambda$ in the left plane of $\mathbb{C}$, we only require the method to be stable for any $h\lambda$ between two half lines starting from $(0,0)$ and making an angle of $\alpha \in [0, \frac{\pi}{2}[$ with the real line :

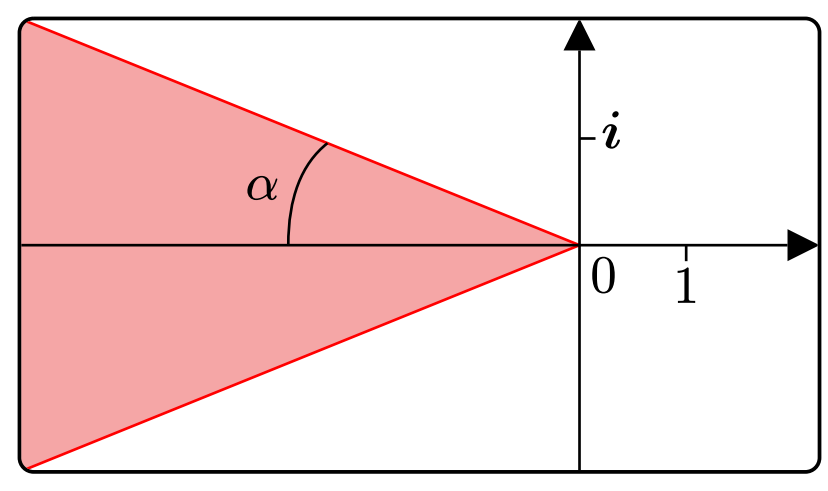

**I.1.26 Definition : A($\alpha$)-stability**

Let $s \in \mathbb{N}^*$, $M \in \mathrm{RK}_s$, $Y_1$ the approximation function of $M$, and $\alpha \in [0, \frac{\pi}{2}[$.

- $M$ is said to be A($\alpha$)-stable if and only if, for all infinite constant step size sequence $(h)_{q\in\mathbb{N}}$ with $h \in \mathbb{R}^*$, for all initial value problem $\left(f, t_0, y_{1,0,1}\right)$ of order 1 of the form :

$$\lambda \in \mathbb{C},\ \forall (t,x) \in \mathbb{R} \times \mathbb{C},\ f(t,x) = \lambda x$$

If :

$$\mathfrak{R}(h\lambda) \leq 0,\ \ \tan(\alpha)\mathfrak{R}(h\lambda) \leq \mathfrak{I}(h\lambda) \leq -\tan(\alpha)\mathfrak{R}(h\lambda)$$

Then $\left(Y_{1,f,q(s+1),1}\left(t_0, y_{1,0,1}, (h)_{q\in\mathbb{N}}\right)\right)_{q\in\mathbb{N}}$ is bounded.

- $M$ is said to be absolute A($\alpha$)-stable if and only if, for all infinite constant step size sequence $(h)_{q\in\mathbb{N}}$ with $h \in \mathbb{R}^*$, for all initial value problem $\left(f, t_0, y_{1,0,1}\right)$ of order 1 of the form :

$$\lambda \in \mathbb{C},\ \forall (t,x) \in \mathbb{R} \times \mathbb{C},\ f(t,x) = \lambda x$$

If :

$$\mathfrak{R}(h\lambda) < 0,\ \ \tan(\alpha)\mathfrak{R}(h\lambda) < \mathfrak{I}(h\lambda) < -\tan(\alpha)\mathfrak{R}(h\lambda)$$

Then $\left(Y_{1,f,q(s+1),1}\left(t_0, y_{1,0,1}, (h)_{q'\in\mathbb{N}}\right)\right)_{q\in\mathbb{N}}$ tends towards 0 as $q$ tends towards $+\infty$.

Similarly to A-stability, A($\alpha$)-stability can be characterized using the stability function.

**I.1.27 Theorem : Characterization of A($\alpha$)-stability**

Let $s \in \mathbb{N}^*$, $\alpha \in [0, \frac{\pi}{2}[$, $M \in \mathrm{RK}_s$, $\tilde{R}_1$ the stability function of $M$, and $\tilde{\mathcal{S}}_1$ its domain of definition.

- $M$ is A($\alpha$)-stable if and only if, for all $z \in \tilde{\mathcal{S}}_1$, if $\tan(\alpha)\mathfrak{R}(z) \leq \mathfrak{I}(z) \leq -\tan(\alpha)\mathfrak{R}(z),\ \mathfrak{R}(z) \leq 0$, then $|\tilde{R}_1(z)| \leq 1$.
- $M$ is absolute A($\alpha$)-stable if and only if, for all $z \in \tilde{\mathcal{S}}_1$, if $\tan(\alpha)\mathfrak{R}(z) < \mathfrak{I}(z) < -\tan(\alpha)\mathfrak{R}(z)$, $\mathfrak{R}(z) < 0$, then $|\tilde{R}_1(z)| < 1$.

Proof

We prove a more general version in Section IV.2.26

Having established the basics of the theory of Runge-Kutta methods, we now examine how Runge-Kutta methods handle initial value problems of order higher than one.

## I.2 Drawbacks of First-Order Reduction

Runge-Kutta methods can only solve initial value problems of order 1 and thus do not deal with the case $n \in \llbracket 2, +\infty \llbracket$. However, there exists a way to solve initial value problems of higher order[1].

**I.2.1 Definition : Rewriting for initial value problems of higher order**

Let $\left(f, t_0, y_{n,0}\right)$ an initial value problem of order $n \in \mathbb{N}^*$. We define :

$$\mathbf{U}' = \left\{ \left( t, \begin{bmatrix} x_1 \\ \vdots \\ x_n \end{bmatrix} \right) \in \mathbb{R} \times \mathbb{R}^{\llbracket 1, nd \rrbracket} \mid (t, x) \in \mathbf{U} \right\}$$

$$y_{1,0,1} = \begin{bmatrix} y_{n,0,1} \\ \vdots \\ y_{n,0,n} \end{bmatrix}, \ \forall \left( t, \begin{bmatrix} x_1 \\ \vdots \\ x_n \end{bmatrix} \right) \in \mathbf{U}', \ f'\left( t, \begin{bmatrix} x_1 \\ \vdots \\ x_n \end{bmatrix} \right) = \begin{bmatrix} f(t,x) \\ x_1 \\ \vdots \\ x_{n-1} \end{bmatrix} \in \mathbb{R}^{\llbracket 1, nd \rrbracket}$$

We define $\left(f', t_0, y_{1,0,1}\right)$ as the reduction to order 1 of $\left(f, t_0, y_{n,0}\right)$.

**I.2.2 Proposition : Solutions of rewritten initial value problems**

Let $\left(f, t_0, y_{n,0}\right)$ an initial value problem of order $n \in \mathbb{N}^*$, and $\left(f', t_0, y_{1,0,1}\right)$ its reduction to order 1.

- If $\hat{\mathrm{y}}$ is a solution of $\left(f, t_0, y_{n,0}\right)$ then $\mathcal{J}^{n-1}\hat{\mathrm{y}}$ is a solution of $\left(f', t_0, y_{1,0,1}\right)$.
- If $\hat{\mathrm{y}}'$ is a solution of $\left(f', t_0, y_{1,0,1}\right)$ then there exists a unique $\hat{\mathrm{y}}$ such that $\hat{\mathrm{y}}' = \mathcal{J}^{n-1}\hat{\mathrm{y}}$. $\hat{\mathrm{y}}$ is a solution of $\left(f, t_0, y_{n,0}\right)$.

Proof

Trivial.

The reduction to order 1 allows to use Runge-Kutta methods on any initial value problems.

**I.2.3 Proposition : Form of rewritten method**

Let $s \in \mathbb{N}^*$, $M = (\tau, w_1) \in \mathrm{RK}_s$, $f$ a differential equation function of order $n$, $\left(t, y_{n,0}, h\right) \in \mathbf{U} \times \mathbb{R}$, and $\left(f', t, y_{1,0,1}\right)$ the reduction to order 1 of $\left(f, t, y_{n,0}\right)$.

The stage system of $M$ for the initial value problem $\left(f', t, y_{1,0,1}\right)$ is equivalent to :

$$\forall (j, N) \in \llbracket 1, s+1 \rrbracket \times \llbracket 1, n \rrbracket$$

$$y_{n,j,N} = \sum_{N'=1}^{N} \sum_{\substack{J \in \llbracket 1, s+1 \rrbracket^{\llbracket 1, N-N'+1 \rrbracket} \\ J_1 = j \\ k \neq 1 \Rightarrow J_k \leq s}} \left( \prod_{k=1}^{N-N'} w_{1,J_k,J_{k+1}} \right) h^{N-N'} y_{n,0,N'}$$

$$+ h^N \sum_{\substack{J \in \llbracket 1, s+1 \rrbracket^{\llbracket 1, N+1 \rrbracket} \\ J_1 = j \\ k \neq 1 \Rightarrow J_k \leq s}} \left( \prod_{k=1}^{N} w_{1,J_k,J_{k+1}} \right) f\left(t + \tau_{J_{N+1}} h, y_{n,J_{N+1}}\right)$$

Where $y_{n,j,n-i}$ is the approximation of the $i^{\text{th}}$ derivative of $\hat{\mathrm{y}}$ at stage $j$.

Proof

We prove a more general version in Section III.2.29.

We can observe that this expression is similar to Taylor's theorem, for all $(i, j) \in \llbracket 0, n-1 \rrbracket \times \llbracket 1, s+1 \rrbracket$ :

$$\frac{\mathrm{d}^i \hat{\mathrm{y}}}{\mathrm{d}t^i}(t_0 + \tau_j h) = \sum_{i'=i}^{n-1} \frac{(\tau_j h)^{i'-i}}{(i'-i)!} \frac{\mathrm{d}^{i'} \hat{\mathrm{y}}}{\mathrm{d}t^{i'}}(t_0) + \underset{h \to 0}{\mathcal{O}}(h^{n-i}) = \sum_{N'=1}^{n-i} \frac{(\tau_j h)^{n-i-N'}}{(n-i-N')!} y_{n,0,N'} + \underset{h \to 0}{\mathcal{O}}(h^{n-i})$$

If the approximations were to use the initial values in the same way as Taylor's theorem, they would be more accurate.

### I.2.4 Proposition : Use of initial values of rewritten Runge-Kutta methods

Let $(s,n) \in \mathbb{N}^* \times \mathbb{N}^*$, $M = (\tau, w_1) \in \mathrm{RK}_s$, and $\tau_{s+1} = 1$.

For all differential equation function $f$ of order $n$, for all $(t, y_{n,0}, h) \in \mathbf{U} \times \mathbb{R}$, if we let $(f', t, y_{1,0,1})$ be the reduction to order 1 of $(f, t, y_{n,0})$, the stage system of $M$ for the initial value problem $(f', t, y_{1,0,1})$ is equivalent to :

$$\forall (j,N) \in [\![1, s+1]\!] \times [\![1, n]\!]$$

$$y_{n,j,N} = \sum_{N'=1}^{N} \frac{(\tau_j h)^{N-N'}}{(N-N')!} y_{n,0,N'} + h^N \sum_{\substack{J \in [\![1,s+1]\!]^{[\![1,N+1]\!]} \\ J_1 = j \\ k \neq 1 \Rightarrow J_k \leq}} \left( \prod_{k=1}^{N} w_{1,J_k,J_{k+1}} \right) f\left(t + \tau_{J_{N+1}} h, y_{n,J_{N+1}}\right)$$

If and only if :

$$w_1 \begin{bmatrix} 1 & \tau_1 & \dots & \tau_1^{n-2} \\ 1 & \tau_2 & \dots & \tau_1^{n-2} \\ \vdots & \vdots & \ddots & \vdots \\ 1 & \tau_s & \dots. & \tau_s^{n-2} \end{bmatrix} = \begin{bmatrix} \frac{\tau_1}{1} & \frac{\tau_1^2}{2} & \dots & \frac{\tau_1^{n-1}}{n-1} \\ \frac{\tau_2}{1} & \frac{\tau_2^2}{2} & \dots & \frac{\tau_2^{n-1}}{n-1} \\ \vdots & \vdots & \ddots & \vdots \\ \frac{\tau_{s+1}}{1} & \frac{\tau_{s+1}^2}{2} & \dots & \frac{\tau_{s+1}^{n-1}}{n-1} \end{bmatrix}$$

Proof

It is sufficient that :

$$\forall (j,N) \in [\![1, s+1]\!] \times [\![2, n]\!],\ \forall N' \in [\![1, N]\!], \quad \sum_{\substack{J \in [\![1,s+1]\!]^{[\![1,N-N'+1]\!]} \\ J_1 = j \\ k \neq 1 \Rightarrow J_k \leq s}} \prod_{k=1}^{N-N'} w_{1,J_k,J_{k+1}} = \frac{\tau_j^{N-N'}}{(N-N')!}$$

By taking $f$ to be identically 0 we get that this condition is also necessary. This condition is equivalent to :

$$\forall (j,N) \in [\![1, s+1]\!] \times [\![1, n-1]\!], \quad \sum_{\substack{J \in [\![1,s+1]\!]^{[\![1,N+1]\!]} \\ J_1 = j \\ k \neq 1 \Rightarrow J_k \leq s}} \prod_{k=1}^{N} w_{1,J_k,J_{k+1}} = \frac{\tau_j^N}{N!}$$

We rewrite this system as :

$$\forall j \in [\![1, s+1]\!], \quad \sum_{j'=1}^{s} w_{1,j,j'} = \tau_j$$

$$\forall (j,N) \in [\![1, s+1]\!] \times [\![2, n-1]\!], \quad \sum_{j'=1}^{s} w_{1,j,j'} \sum_{\substack{J \in [\![1,s]\!]^{[\![1,N']\!]} \\ J_1 = j'}} \prod_{k=1}^{N-1} w_{1,J_k,J_{k+1}} = \frac{\tau_j^N}{N!}$$

We replace the sum over $[\![1, s]\!]^{[\![1,N']\!]}$ with $\frac{\tau_{j'}^{N-1}}{(N-1)!}$, which yields the equivalent system of equations :

$$\forall (j, N) \in [\![1, s+1]\!] \times [\![1, n-1]\!],\ \sum_{j'=1}^{s} w_{1,j,j'} \tau_{j'}^{N-1} = \frac{\tau_j^N}{N}$$

Rewriting this system in matrix form proves the theorem.

The matrix on the left side of the condition is similar to a Vandermonde matrix, thus if it is a square matrix, which is equivalent to $n = s + 1$, this matrix is invertible if and only if the nodes are distinct. If the nodes were indeed distinct, we could multiply by the inverse of the matrix and find the expression of the weights, which implies that the weights would be entirely defined by this condition :

$$w_1 = \begin{bmatrix} \frac{\tau_1}{1} & \frac{\tau_1^2}{2} & \cdots & \frac{\tau_1^s}{s} \\ \frac{\tau_2}{1} & \frac{\tau_2^2}{2} & \cdots & \frac{\tau_2^s}{s} \\ \vdots & \vdots & \ddots & \vdots \\ \frac{\tau_{s+1}}{1} & \frac{\tau_{s+1}^2}{2} & \cdots & \frac{\tau_{s+1}^s}{s} \end{bmatrix} \begin{bmatrix} 1 & \tau_1 & \cdots & \tau_1^{s-1} \\ 1 & \tau_2 & \cdots & \tau_1^{s-1} \\ \vdots & \vdots & \ddots & \vdots \\ 1 & \tau_s & \cdots. & \tau_s^{s-1} \end{bmatrix}^{-1}$$

However, there can be nodes of equal value. There is, fortunately, a workaround. We can factorize by the nodes of equal value in the system of equations. To do so, we define $s' \in \mathbb{N}^*$ the number of distinct nodes without $\tau_{s+1}$, and $\tau'_1, \tau'_2, \dots\ \tau'_{s'}$ these distinct nodes. We also define :

$$\forall (j, j') \in [\![1, s+1]\!] \times [\![1, s']\!],\ w'_{1,j,j'} = \sum_{\substack{j'' \in [\![1,s]\!] \\ \tau'_{j'} = \tau_{j''}}} w_{1,j,j''}$$

Factorizing by $\tau'_{j'}$ yields the equivalent system :

$$\forall (N, j) \in [\![1, n-1]\!] \times [\![1, s+1]\!],\ \sum_{j'=1}^{s'} w'_{1,j,j'} \tau'^{N-1}_{j'} = \frac{\tau_j^N}{N}$$

Which in matrix form is :

$$w'_1 \begin{bmatrix} 1 & \tau'_1 & \cdots & \tau'^{n-2}_1 \\ 1 & \tau'_2 & \cdots & \tau'^{n-2}_1 \\ \vdots & \vdots & \ddots & \vdots \\ 1 & \tau'_{s'} & \cdots. & \tau'^{n-2}_{s'} \end{bmatrix} = \begin{bmatrix} \frac{\tau_1}{1} & \frac{\tau_1^2}{2} & \cdots & \frac{\tau_1^{n-1}}{n-1} \\ \frac{\tau_2}{1} & \frac{\tau_2^2}{2} & \cdots & \frac{\tau_2^{n-1}}{n-1} \\ \vdots & \vdots & \ddots & \vdots \\ \frac{\tau_{s+1}}{1} & \frac{\tau_{s+1}^2}{2} & \cdots & \frac{\tau_{s+1}^{n-1}}{n-1} \end{bmatrix}$$

We now have essentially the same system with, by definition, distinct nodes. If $n = s' + 1$, the $w'_{1,j,j'}$ are entirely defined by this equation. Since these terms are sums of weights, we actually retrieve a system of linear equations, one for each distinct node. Since none of these equations share variables, they are not redundant, which means that this condition reduces the number of free parameters by the number of distinct nodes $s'$ for each stage $j \in [\![1, s+1]\!]$.

If $n > s' + 1$, it is unlikely that this equation has a solution, thus the initial values cannot be used in accordance to Taylor's theorem. We can in fact prove that no Runge-Kutta method can satisfy this condition for $n > 2s'$.

### I.2.5 Theorem : Limit to the right use of initial values

Let $s \in \mathbb{N}^*$, $M = (\tau, w) \in \mathrm{RK}_s$, $s'$ the number of different nodes, and $\tau_{s+1} = 1$.
The highest integer $K \in \mathbb{N}^*$ such that, for all $k \in [\![1, K]\!]$ :

$$\forall j \in [\![1, s+1]\!],\ \sum_{j'=1}^{s} w_{1,j,j'} \tau_{j'}^{k-1} = \frac{\tau_j^k}{k}$$

Is inferior or equal to $2s'$.

Proof

We prove a more general version in Section III.1.21.

We here face a dilemma. If a method has many distinct nodes, it is able to make optimal use of the initial values for higher values of $n$, but the conditions leaves little to no freedom to the possible values of the weights. If a method has many equal nodes, there are more possibilities for the values of the weights, but the method uses inefficiently the initial values. None of these choices are ideal, and for any method, there exists a limit value of $n$ from which it is not possible to efficiently use the initial values. We would also need to add the conditions for the order of consistency of a method on top of this condition.

We will see in Section III.1.19 that these conditions are the solved system conditions for order of consistency $n - 1$ at rank 1 at all stages, which means that these conditions allow for a better precision, but they are too strong of a requirement.

The probability that the evaluations of $f$ are optimally used is low, but unless we assume additional conditions on the method, on its order for instance, we cannot be more precise. There is, however, an issue with whether the evaluations of $f$ and the initial values are used in the approximations of explicit methods.

### I.2.6 Theorem : Initial values and values of $f$ not used by explicit methods

Let $s \in \mathbb{N}^*$, $M \in \mathrm{RK}_s$ with $M$ explicit, $(f, t_0, y_{n,0})$ an initial value problem of order $n \in [\![2, +\infty[\![$, and $y_n$ the approximations of $M$.

1. For all $(N, j) \in [\![1, n]\!] \times [\![1, s+1]\!]$, for all $N' \in [\![1, N-j]\!]$, $y_{n,j,N}$ does not use $y_{n,0,N'}$.
2. For all $(N, j) \in [\![1, n]\!] \times [\![1, s+1]\!]$, for all $j' \in [\![1 + \max(0, j-N), s]\!]$, $y_{n,j,N}$ does not use the value $f(t_0 + \tau_{j'} h, y_{n,j'})$.

Proof

1. $y_{n,0,N'}$ is used in the calculation of $y_{n,j,N}$ if and only if :

$$\sum_{\substack{J \in [\![1,s+1]\!]^{[\![1,N-N'+1]\!]} \\ J_1 = j \\ k \neq 1 \Rightarrow J_k \leq s}} \prod_{k=1}^{N-N'} w_{1,J_k,J_{k+1}} \neq 0$$

   But :

$$\sum_{\substack{J \in [\![1,s+1]\!]^{[\![1,N-N'+1]\!]} \\ J_1 = j \\ k \neq 1 \Rightarrow J_k \leq s}} \prod_{k=1}^{N-N'} w_{1,J_k,J_{k+1}} \neq 0$$

$$\Rightarrow \exists J \in [\![1, s+1]\!]^{[\![1,N-N'+1]\!]},\ J_1 = j,\ k \neq 1 \Rightarrow J_k \leq s,\ \prod_{k=1}^{N-N'} w_{1,J_k,J_{k+1}} \neq 0$$

Since explicit methods satisfy $j \leq j' \Rightarrow w_{1,j,j'} = 0$, we have :

$$\exists J \in \llbracket 1, s+1 \rrbracket^{\llbracket 1, N-N'+1 \rrbracket},\ J_1 = j,\ k \neq 1 \Rightarrow J_k \leq s,\ \prod_{k=1}^{N-N'} w_{1,J_k,J_{k+1}} \neq 0$$

$$\Rightarrow \exists J \in \llbracket 1, s+1 \rrbracket^{\llbracket 1, N-N'+1 \rrbracket},\ J_1 = j,\ k \neq 1 \Rightarrow J_k \leq s,\ \forall k \in \llbracket 1, N-N' \rrbracket\ J_k > J_{k+1}$$

Thus, the coordinates of $J$ need to be strictly decreasing while being in between 1 and $s+1$. Since $J_1 = j$ this condition becomes : $j > J_2 > \ldots > J_{N-N'+1} > 0$. This is not possible when $N' \leq N - j$, we conclude that the values of $y_{n,0,N'}$ with $N' \in \llbracket 1, N-j \rrbracket$ are missing, even though they are necessary for the accuracy of the approximations.

2. $f\left(t_0 + \tau_{j'} h, y_{n,j'}\right)$ is used in the computation of $y_{n,j,N}$ if and only if :

$$\sum_{\substack{J \in \llbracket 1, s+1 \rrbracket^{\llbracket 1, N+1 \rrbracket} \\ J_1 = j,\ J_{N+1} = j' \\ k \neq 1 \Rightarrow J_k \leq s}} \prod_{k=1}^{N} w_{1,J_k,J_{k+1}} \neq 0$$

Using a similar reasoning to that of 1, we deduce that for explicit methods it is necessary that there exists $J \in \llbracket 1, s+1 \rrbracket^{\llbracket 1, N+1 \rrbracket}$ with strictly decreasing coordinates such that $J_1 = j,\ J_{N+1} = j'$. This is equivalent to $j > J_2 > \ldots > J_N > j'$. This is not possible when $j' \geq j + 1 - N$, thus, many evaluations of $f$ are not used, in fact if $j \leq N$, the term $y_{n,j,N}$ does not use any, and $y_{n,j,N}$ could have been computed directly using only the initial values.

These issues with explicit methods can be illustrated with the explicit Euler method, the Runge-Kutta method defined by :

$$\tau = [0] \quad w_1 = \begin{bmatrix} 0 \\ 1 \end{bmatrix}$$

Which yields when solving an initial value problem $\left(f, t_0, y_{3,0}\right)$ of order 3 :

$$\begin{aligned} y_{3,2,1} &= y_{3,0,1} + hf\left(t_0, y_{3,0}\right) \\ y_{3,2,2} &= y_{3,0,2} + hy_{3,0,1} \\ y_{3,2,3} &= y_{3,0,3} + hy_{3,0,2} \end{aligned}$$

As expected, both the initial value $y_{3,0,1}$ and the evaluation of $f$ are missing in the expression of $y_{3,2,3}$. Another issue has to do with the definition of order of consistency, the rewriting implies that :

$$\forall i \in \llbracket 0, n-1 \rrbracket,\ \frac{\mathrm{d}^i \hat{\mathrm{y}}}{\mathrm{d}t^i}(t' + h) - y_{n,s+1,n-i} = \underset{\substack{h \to 0 \\ t' \in \Omega'}}{\mathcal{O}}\left(h^{v_{1,s+1,1}+1}\right)$$

When more complex behavior can be observed. For instance, simply using Taylor's theorem gives :

$$\forall i \in \llbracket 0, n-1 \rrbracket,\ \frac{\mathrm{d}^i \hat{\mathrm{y}}}{\mathrm{d}t^i}(t' + h) - \sum_{i'=i}^{n-1} \frac{h^{i'-i}}{(i'-i)!} y_{n,0,n-i'} - \frac{h^{n-i}}{(n-i)!} f\left(t', y_{n,0}\right) = \underset{\substack{h \to 0 \\ t \in \Omega'}}{\mathcal{O}}\left(h^{n-i+1}\right)$$

We can observe that for high enough $n$, for most of the derivatives, and in particular the solution $\hat{\mathrm{y}}$, Taylor's theorem outperforms any Runge-Kutta method in terms of order of consistency.

There are plenty of issues with Runge Kutta methods, we solve them in the next sections.

# II Multi-order Runge-Kutta methods

## II.1 Derivation of multi-order Runge-Kutta methods

In this section we define a new type of method made to solve initial value problems of any order. We call this class of method multi-order Runge-Kutta methods in reference to their ability to solve initial value problems of different orders. The approach we take to define those methods is very similar to that of Runge-Kutta methods, we first transform the initial value problem into an equivalent integration problem, then approximate the integral using a quadrature.

**II.1.1 Proposition : Integral form of an initial value problem**

Let $\left(f, t_0, y_{n,0}\right)$ an initial value problem of order $n \in \mathbb{N}^*$ with $f$ continuous, and $\hat{y} \in D^n\left(\Omega, \mathbb{R}^{[\![1,d]\!]}\right)$ a solution candidate.

$\hat{y}$ is a solution of the initial value problem if and only if, for all $(i, t) \in [\![0, n-1]\!] \times \Omega$ :

$$\frac{\mathrm{d}^i \hat{y}}{\mathrm{d}t^i}(t) = \sum_{i'=i}^{n-1} \frac{\left(t-t_0\right)^{i'-i}}{(i'-i)!} \frac{\mathrm{d}^{i'} \hat{y}}{\mathrm{d}t^{i'}}(t_0) + \int_{t_0}^{t} \frac{(t-x)^{n-1-i}}{(n-1-i)!} f\left(\overline{\mathcal{J}}^{n-1} \hat{y}(x)\right) \mathrm{d}x$$

Proof

If we assume $\hat{y}$ is a solution, then its $n^{\text{th}}$ derivative is continuous, $\hat{y}$ is thus of class $C^n$. Taylor's theorem with integral remainder gives, for all $t \in \Omega$ :

$$\forall i \in [\![0, n-1]\!],\ \frac{\mathrm{d}^i \hat{y}}{\mathrm{d}t^i}(t) = \sum_{i'=i}^{n-1} \frac{\left(t-t_0\right)^{i'-i}}{(i'-i)!} \frac{\mathrm{d}^{i'} \hat{y}}{\mathrm{d}t^{i'}}(t_0) + \int_{t_0}^{t} \frac{(t-x)^{n-1-i}}{(n-1-i)!} \frac{\mathrm{d}^n \hat{y}}{\mathrm{d}t^n}(x) \, \mathrm{d}x$$

Thus, if we replace $\frac{\mathrm{d}^{i'} \hat{y}}{\mathrm{d}t^{i'}}(t_0)$ with $y_{n,0,n-i'}$ and $\frac{\mathrm{d}^n \hat{y}}{\mathrm{d}t^n}$ with $f \circ \overline{\mathcal{J}}^{n-1} \hat{y}$, we get the expression of the theorem. If $\hat{y}$ satisfies the system of equations, $\frac{\mathrm{d}^{n-1} \hat{y}}{\mathrm{d}t^{n-1}}$ is $C^1$, and by differentiating it we see that $\hat{y}$ satisfies the differential equation induced by $f$. For all $i \in [\![0, n-1]\!]$, evaluating the equation of $\frac{\mathrm{d}^i \hat{y}}{\mathrm{d}t^i}$ at $t = t_0$ yields that $\hat{y}$ satisfies the initial values. $\hat{y}$ is thus a solution of the initial value problem.

If we replace $t$ with $h = t - t_0$ and substitute in $z = 2\frac{x - t_0}{h} - 1$, we get, for all $i \in [\![0, n-1]\!]$:

$$\frac{\mathrm{d}^i \hat{y}}{\mathrm{d}t^i}(t_0 + h) = \sum_{i'=i}^{n-1} \frac{h^{i'-i}}{(i'-i)!} y_{n,0,n-i'} + \frac{h^{n-i}}{(n-1-i)! 2^{n-i}} \int_{-1}^{1} (1-z)^{n-1-i} f\left(\overline{\mathcal{J}}^{n-1} \hat{y}\left(t_0 + \frac{1+z}{2} h\right)\right) \mathrm{d}z$$

We now have an issue, Runge-Kutta methods use a Gauss-Legendre quadrature to approximate the only integral they are dealing with, but in this case there are $n$ integrals to approximate and endless possibilities. We could use a different number of points, use different nodes or even use a different type of Gaussian quadrature for each individual integral. The other type of Gaussian quadrature we could use are Gauss-Jacobi quadratures[6], which approximate, with $g \in C^0\left([-1, 1], \mathbb{R}^{[\![1,d]\!]}\right)$, $(\alpha, \beta) \in ]-1, +\infty[^2$, $s \in \mathbb{N}^*$, and $(b, c) \in \mathbb{R}^{[\![1,s]\!]} \times \mathbb{R}^{[\![1,s]\!]}$ :

$$\int_{-1}^{1} (1-x)^{\alpha} (1+x)^{\beta} g(x) \, \mathrm{d}x \approx \sum_{j=1}^{s} b_j g\left(c_j\right)$$

Here we are in the case $\beta = 0,\ \alpha = n - i - 1$. We can however observe that if an integral requires the value of $f$ at a certain node, it can be added to the calculations of the other quadratures, even if it implies that some weights are 0. All integrals thus have the same number of points and the same nodes. Using a different quadrature rule only changes the coefficients before each value of $f$, which can be put in a general weight term. We settle on Gauss-Jacobi quadratures with $s \in \mathbb{N}^*$ points :

$$\frac{\mathrm{d}^i \hat{\mathrm{y}}}{\mathrm{d}t^i}(t_0 + h) \approx \sum_{i'=i}^{n-1} \frac{h^{i'-i}}{(i'-i)!} y_{n,0,n-i'} + \frac{h^{n-i}}{(n-1-i)!2^{n-i}} \sum_{j=1}^{s} b_{i,j} f\left(\overline{\mathcal{J}}^{n-1} \hat{\mathrm{y}}\left(t_0 + \frac{1+c_j}{2} h\right)\right) \mathrm{d}x$$

With $(b, c) \in \mathbb{R}^{[\![0,n-1]\!] \times [\![1,s]\!]} \times \mathbb{R}^{[\![1,s]\!]}$. We redefine the coefficients to make the formula nicer to the eyes :

$$\forall j \in [\![1, s]\!],\ \tau_j = \frac{1 + c_j}{2} \in \mathbb{R}$$

$$\forall (i, j) \in [\![0, n-1]\!] \times [\![1, s]\!],\ w_{n-i,s+1,j} = b_{i,j} \frac{n-i}{2^{n-i}} \in \mathbb{R}$$

We get, for all $i \in [\![0, n-1]\!]$ :

$$\frac{\mathrm{d}^i \hat{\mathrm{y}}}{\mathrm{d}t^i}(t_0 + h) \approx \sum_{i'=i}^{n-1} \frac{h^{i'-i}}{(i'-i)!} y_{n,0,n-i'} + \frac{h^{n-i}}{(n-i)!} \sum_{j=1}^{s} w_{n-i,s+1,j'} f\left(\overline{\mathcal{J}}^{n-1} \hat{\mathrm{y}}\left(t_0 + \tau_j h\right)\right)$$

For all $(j, i) \in [\![1, s]\!] \times [\![0, n-1]\!]$, we approximate every derivative in $\overline{\mathcal{J}}^{n-1} \hat{\mathrm{y}}\left(t_0 + \tau_j h\right)$ using the exact same process. We also define $\tau_{s+1} = 1$ to have, for all $(j, i) \in [\![1, s+1]\!] \times [\![0, n-1]\!]$ :

$$\frac{\mathrm{d}^i \hat{\mathrm{y}}}{\mathrm{d}t^i}\left(t_0 + \tau_j h\right) \approx y_{n,j,n-i} = \sum_{i'=i}^{n-1} \frac{\left(\tau_j h\right)^{i'-i}}{(i'-i)!} y_{n,0,n-i'} + \frac{h^{n-i}}{(n-i)!} \sum_{j'=1}^{s} w_{n-i,j,j'} f\left(t_0 + \tau_{j'} h, y_{n,j'}\right)$$

All quadratures use the same values of $f$ because any value computed for any quadrature can be used by all the others. A leitmotif throughout this paper is replacing $i$ with $N = n - i$, the reasons for this change of variables have been given after the definition of an initial value problem. One fact to keep in mind is $i = n - N$. We replace $i$ with $N = n - i$ and $i'$ with $N' = n - i'$. We have, for all $(j, N) \in [\![1, s+1]\!] \times [\![1, n]\!]$ :

$$y_{n,j,N} = \sum_{N'=1}^{N} \frac{\left(\tau_j h\right)^{N-N'}}{(N-N')!} y_{n,0,N'} + \frac{h^N}{N!} \sum_{j'=1}^{s} w_{N,j,j'} f\left(t_0 + \tau_{j'} h, y_{n,j'}\right)$$

This formula derived from Taylor's theorem is not general enough, rewritten Runge-Kutta methods for instance are generally not of this form. We should hence consider a more general type of method which uses the formula, for all $(j, N) \in [\![1, s+1]\!] \times [\![1, n]\!]$ :

$$y_{n,j,N} = \sum_{N'=1}^{N} \tilde{w}_{N,N',j} h^{N-N'} y_{n,0,N'} + \frac{h^N}{N!} \sum_{j'=1}^{s} w_{N,j,j'} f\left(t_0 + \tau_{j'} h, y_{n,j'}\right)$$

Where $\tilde{w} \in \mathbb{R}^{\bigcup_{N \in [1,n]} \{N\} \times [\![1,N]\!] \times [\![1,s+1]\!]}$ is a new sequence of constants called the secondary weights to distinguish them from the main weights $w$.

The nodes, main weights and secondary weights entirely define a method. To store these values in a mathematical object, we need similar matrices to the ones used to define Runge-Kutta methods, but this time we need a finite sequence of them.

### II.1.2 Definition : General multi-order Runge-Kutta methods

Let $(n, s) \in \mathbb{N}^* \times \mathbb{N}^*$.

- We define $\mathrm{GMORK}_{n,s} = \mathbb{R}^{[\![1,s]\!]} \times \mathbb{R}^{[\![1,n]\!] \times [\![1,s+1]\!] \times [\![1,s]\!]} \times \mathbb{R}^{\cup_{N \in [1,n]} \{N\} \times [\![1,N]\!] \times [\![1,s+1]\!]}$ as the set of general multi-order Runge-Kutta methods of length $n$ with $s$ points.
- The general multi-order Runge-Kutta methods of length $n$ with $s$ points are the elements of $\mathrm{GMORK}_{n,s}$. Let $(\tau, w, \tilde{w}) \in \mathrm{GMORK}_{n,s}$. $\tau$ are the nodes, $w$ the main weights, $\tilde{w}$ the secondary weights.

We use the word length for $n$ to refer to the length of the weight sequence of the method and to avoid using the word order, which already has too many meanings.

We added the adjective general because these methods are too general to simply be called multi-order Runge-Kutta methods, their use of initial values is not sensical enough.

One type of method whose use of initial values is much more sensical is the type of methods we have previously derived from Taylor's theorem, we thus give it a name.

### II.1.3 Definition : Node-determined methods

Let $(n, s) \in \mathbb{N}^* \times \mathbb{N}^*$.

- Let $M = (\tau, w, \tilde{w}) \in \mathrm{GMORK}_{n,s}$.
  $M$ is said to be node-determined if and only if there exists $\tau_{s+1} \in \mathbb{R}$ such that :

  $$\forall (j, N) \in [\![1, s+1]\!] \times [\![1, n]\!],\ \forall N' \in [\![1, N]\!],\ \tilde{w}_{N,N',j} = \frac{\tau_j^{N-N'}}{(N-N')!}$$

  A node-determined method is entirely defined by the pair $(\tau, w)$ and can thus be identified to it, with $\tau$ its nodes, and $w$ its weight sequence.
- We define $\mathrm{NDMORK}_{n,s} = \mathbb{R}^{[\![1,s+1]\!]} \times \mathbb{R}^{[\![1,n]\!] \times [\![1,s+1]\!] \times [\![1,s]\!]}$ as the set of node-determined multi-order Runge-Kutta methods of length $n$ with $s$ points.

The adjective "node-determined" was added because the nodes dictates the coefficients of the Taylor polynomial the approximations use.

This definition is this time too restrictive, it is better to have a more general definition of multi-order Runge-Kutta methods, though this type of method will probably be the most used. To have a good definition of multi-order Runge-Kutta methods, we first need to know what is a good general multi-order Runge-Kutta method, the actual definition is hence in the section on convergence and order of consistency. Multi-order Runge-Kutta methods are somewhere between node-determined and general multi-order Runge-Kutta methods.

We dropped the requirement $\tau_{s+1} = 1$ because another condition (convergence) can enforce it, and it is useful to not be restricted to the case $\tau_{s+1} = 1$ (sub-methods).

The process of approximation of a general multi-order Runge-Kutta method seems simple enough, but it has yet to be rigorously defined. To approximate the solution of an initial value problem we use a system of equations, we therefore need to make sure that this system has a unique solution. To formalize this issue, we define two closely related functions.

**II.1.4 Definition : Approximation function & Evaluation function**

Let $(n, s) \in \mathbb{N}^* \times \mathbb{N}^*$, $M = (\tau, w, \tilde{w}) \in \mathrm{GMORK}_{n,s}$, and $f$ a differential equation function of order $n$. We define the existence-uniqueness domain $\mathcal{U}_{n,f} \in \mathcal{P}(\mathbf{U} \times \mathbb{R})$ as the set of $(t, y_{n,0}, h) \in \mathbf{U} \times \mathbb{R}$ such that the following system, called the stage system $(M, f, t, y_{n,0}, h)$, has a unique solution :

$$\forall (j, N) \in [\![1, s+1]\!] \times [\![1, n]\!],$$
$$y_{n,j,N} = \sum_{N'=1}^{N} \tilde{w}_{N,N',j} h^{N-N'} y_{n,0,N'} + \frac{h^N}{N!} \sum_{j'=1}^{s} w_{N,j,j'} f\big(t + \tau_{j'} h, y_{n,j'}\big)$$
$$y_{n,s+1} \in \mathbb{R}^{[\![1,n]\!] \times [\![1,d]\!]},\ \forall j \in [\![1, s]\!],\ \big(t + \tau_j h, y_{n,j}\big) \in \mathbf{U}$$

We define $Y_{n,f} \in \{\mathcal{U}_{n,f} \to \mathbb{R}^{[\![1,d]\!]}\}^{[\![0,s+1]\!] \times [\![1,n]\!]}$ the approximation function for length $n$ of $M$ and $F_{n,f} \in \{\mathcal{U}_{n,f} \to \mathbb{R}^{[\![1,d]\!]}\}^{[\![1,s]\!]}$ the evaluation function for length $n$ of $M$ as :

$$\forall j \in [\![0, s+1]\!],\ Y_{n,f,j}\big(t, y_{n,0}, h\big) = y_{n,j}$$
$$\forall j \in [\![1, s]\!],\ F_{n,f,j}\big(t, y_{n,0}, h\big) = f\big(t + \tau_j h, y_{n,j}\big)$$

If $(f, t, y_{n,0})$ has a unique solution $\hat{y}$, the approximation function can be seen as an approximation of the jet of the solution $\mathcal{J}^{n-1}\hat{y}$. Since $(t, y_{n,0})$ is often a given pair of constants, it is useful to define a set which contains the values of $h$ such that $(t, y_{n,0}, h) \in \mathcal{U}_{n,f}$.

**II.1.5 Definition: Step sizes**

Let $(n, s) \in \mathbb{N}^* \times \mathbb{N}^*$, $M \in \mathrm{GMORK}_{n,s}$, $f$ a differential equation function of order $n$, $(t, y_{n,0}) \in \mathbf{U}$, and $\mathcal{U}_{n,f}$ the existence-uniqueness domain of $M$.
We define $\mathcal{H}_{1,f,t,y_{n,0}}$ the set of step sizes at the point $(t, y_{n,0})$ of $M$, as :

$$\mathcal{H}_{1,f,t,y_{n,0}} = \big\{h \in \mathbb{R} \mid \big(t, y_{n,0}, h\big) \in \mathcal{U}_{n,f}\big\}$$

To end this subsection, we state a theorem on the existence of step sizes.

**II.1.6 Theorem : Existence and uniqueness of approximations**

Let $(n, s, k) \in \mathbb{N}^* \times \mathbb{N}^* \times \mathbb{N}$, $M \in \mathrm{GMORK}_{n,s}$, $f$ a differential equation function of order $n$ and class $C^k$, locally Lipschitz continuous in its second variable, $Y_{n,f}$ the approximation function for length $n$ of $M$, $F_{n,f}$ the evaluation function for length $n$ of $M$, and $\mathcal{H}$ the set of step sizes of $M$.

1. If, for all $(j, N) \in [\![1, s+1]\!] \times [\![1, n]\!]$, $\tilde{w}_{N,N,j} = 1$, then, for all $(t, y_{n,0}) \in \mathbf{U}$, there exists $h^* \in \mathbb{R}_+^*$ such that $]-h^*, h^*[ \subset \mathcal{H}_{1,f,t,y_{n,0}}$ and the two following functions are of class $C^k$ :

$$]-h^*, h^*[ \to \mathbb{R}^{[\![0,s+1]\!] \times [\![1,n]\!] \times [\![1,d]\!]}, h \to Y_{n,f}\big(t, y_{n,0}, h\big)$$
$$]-h^*, h^*[ \to \mathbb{R}^{[\![1,s]\!] \times [\![1,d]\!]}, h \to F_{n,f}\big(t, y_{n,0}, h\big)$$

2. If $\mathbf{U} = \Omega \times \mathbb{R}^{[\![1,n]\!] \times [\![1,d]\!]}$ and $f$ is globally Lipschitz continuous in its second variable, then there exists $h^* \in \mathbb{R}_+^*$ such that, for all $(t, y_{n,0}) \in \mathbf{U}$ :

$$\big\{h \in ]-h^*, h^*[\ \mid \forall j \in [\![1, s]\!],\ t + \tau_j h \in \Omega\big\} \subset \mathcal{H}_{1,f,t,y_{n,0}}$$

   And the two following functions are of class $C^k$ :

$$\big\{h \in ]-h^*, h^*[\ \mid \forall j \in [\![1, s]\!],\ t + \tau_j h \in \Omega\big\} \to \mathbb{R}^{[\![0,s+1]\!] \times [\![1,n]\!] \times [\![1,d]\!]}, h \to Y_{n,f}\big(t, y_{n,0}, h\big)$$
$$\big\{h \in ]-h^*, h^*[\ \mid \forall j \in [\![1, s]\!],\ t + \tau_j h \in \Omega\big\} \to \mathbb{R}^{[\![1,s]\!] \times [\![1,d]\!]}, h \to F_{n,f}\big(t, y_{n,0}, h\big)$$

We can note that if $f$ is of class $C^k$ with $k \in \mathbb{N}^*$, then $f$ is locally Lipschitz continuous and the second condition becomes redundant.

**II.1.7 Proposition : Value of $h^*$**

We consider the context of Section II.1.6.

- Let $C_w \in \mathbb{R}_+^*$ such that, for all $(j, N) \in [\![1, s+1]\!] \times [\![1, n]\!],\ C_w \geq \sum_{j'=1}^{s} \frac{|w_{N,j,j'}|}{N!}$.

1. - Let $r \in \mathbb{R}_+^*$ such that $\overline{\mathbf{B}}\big((t, y_{n,0}), r\big) \subset \mathbf{U}$
   - Let $L_r \in \mathbb{R}_+^*$ a Lipschitz constant of $f$ on $\overline{\mathbf{B}}\big((t, y_{n,0}), r\big)$.
   - Let $h' \in \overline{\mathbb{R}}_+^*$ such that, for all $j \in [\![1, s]\!],\ \tau_j \neq 0 \Rightarrow h' \leq \frac{r}{|\tau_j|}$.
   - Let $A = \sup_{(j,h) \in [\![1,s]\!] \times [-h', h']} f\big(t + \tau_j h, y_{n,0}\big)$
   - Let $C_{\tilde{w}}^* \in \mathbb{R}_+$ such that, for all $(j, N) \in [\![1, s+1]\!] \times [\![2, n]\!],\ C_{\tilde{w}}^* \geq \sum_{N'=1}^{N-1} |\tilde{w}_{N,N',j}|$.

   A valid choice of $h^*$ is :

   $$h^* = \min\Big(C_r,\ \sqrt[n]{C_r}, h'\Big),\ C_r = \frac{r}{(rL_r + A)C_w + C_{\tilde{w}}^*\ \|y_{n,0}\|_{n,d}}$$

2. If $\mathbf{U} = \Omega \times \mathbb{R}^{[\![1,n]\!] \times [\![1,d]\!]}$, and $f$ is globally Lipschitz continuous in its second variable with a Lipschitz constant $L \in \mathbb{R}_+^*$, then a valid choice of $h^*$ is :

   $$h^* = \min\left(\frac{1}{LC_w}, \frac{1}{\sqrt[n]{LC_w}}\right)$$

We can observe that the value of $h^*$ of the second case is the limit of $h^*$ of the first case as $r \to +\infty$.

The proof of Section II.1.6 and Section II.1.7 is somewhat involved, and the results are important in the next subsection and the section after that, we therefore break the proof down into multiple lemmas to understand precisely its structure.

To prove this theorem, we rewrite the stage system into an equivalent fixed point problem. This first requires that we define the parametrized approximation function.

**II.1.8 Definition : Parametrized approximation function**

Let $(n, s) \in \mathbb{N}^* \times \mathbb{N}^*$, $M = (\tau, w, \tilde{w}) \in \mathrm{GMORK}_{n,s}$, and $f$ a differential equation function of order $n$. We define $\mathcal{U}'_{n,f}$ the parametrized existence-uniqueness domain of $M$ and $Y'_{n,f}$ the parametrized approximation function for length $n$ of $M$ as :

$$\mathcal{U}'_{n,f} = \left\{\big(t, y_{n,0}, h, x\big) \in \mathbf{U} \times \mathbb{R} \times \mathbb{R}^{[\![1,s+1]\!] \times [\![1,n]\!] \times [\![1,d]\!]} \mid \forall j \in [\![1, s]\!],\ \big(t + \tau_j h, x_j\big) \in \mathbf{U}\right\}$$

$$\forall (j, N) \in [\![1, s+1]\!] \times [\![1, n]\!]$$

$$Y'_{n,f,j,N} : \mathcal{U}'_{n,f} \to \mathbb{R}^{[\![1,d]\!]},$$

$$\big(t, y_{n,0}, h, x\big) \to \sum_{N'=1}^{N} \tilde{w}_{N,N',j} h^{N-N'} y_{n,0,N'} + \frac{h^N}{N!} \sum_{j'=1}^{s} w_{N,j,j'} f\big(t + \tau_{j'} h, x_{j'}\big)$$

**II.1.9 Lemma : Stage system - Fixed point form**

Let $(n,s) \in \mathbb{N}^* \times \mathbb{N}^*$, $M \in \mathrm{GMORK}_{n,s}$, $f$ a differential equation function of order $n$, $\mathcal{U}'_{n,f}$ the parametrized existence-uniqueness domain of $M$, $Y'_{n,f}$ the parametrized approximation function for length $n$ of $M$, and $(t, y_{n,0}, h) \in \mathbf{U} \times \mathbb{R}$.
$y_n \in \mathbb{R}^{[\![1,s+1]\!] \times [\![1,n]\!] \times [\![1,d]\!]}$ is a solution of the stage system $(M, f, t, y_{n,0}, h)$ if and only if it is a fixed point of :

$$\left\{x \in \mathbb{R}^{[\![1,s+1]\!] \times [\![1,n]\!] \times [\![1,d]\!]} \mid (t, y_{n,0}, h, x) \in \mathcal{U}'_{n,f}\right\} \to \mathbb{R}^{[\![1,s+1]\!] \times [\![1,n]\!] \times [\![1,d]\!]},\ x \to Y'_{n,f}(t, y_{n,0}, h, x)$$

Proof

Trivial.

We wish to use the contraction mapping theorem on this fixed point problem, however, since we also want a result on the smoothness of the evaluation and approximation functions, we need a generalized version of this theorem.

**II.1.10 Theorem : Uniform contraction principle**

Let $k \in \mathbb{N}$, $(X, \|\cdot\|_X)$, $(Y, \|\cdot\|_Y)$ two Banach spaces, $V \in \mathcal{P}(X)$ a non-empty open set, $U \in \mathcal{P}(Y)$ a non empty closed set, and $f \in C^k(V \times U, U)$ a contraction mapping in its second variable.
For all $x \in V$, the function $U \to U,\ y \to f(x,y)$ has a unique fixed point noted $g(x)$. $g \in \{V \to U\}$ is of class $C^k$.

It will also prove useful to have a weaker fixed point theorem.

**II.1.11 Theorem : Brouwer's fixed point theorem**

Let $k \in \mathbb{N}^*$, $x \in \mathbb{R}^{[\![1,k]\!]}$, $r \in \mathbb{R}_+^*$, and $f \in \{\overline{\mathbf{B}}(x,r) \to \overline{\mathbf{B}}(x,r)\}$ a continuous function.
$f$ has at least a fixed point.

We now outline the structure of the proof. The main goal is to use the uniform contraction principle on the fixed point form of the stage system. The proof for existence and uniqueness also contains a proof of the existence of fixed points, which will be useful. The proof is hence divided into two parts, the existence part and the contraction part. The parametrized existence-uniqueness domain has four coordinates but we fix the first two coordinates $(t, y_{n,0})$ and make only $(h, x)$ vary.

1. Existence : We prove that there exists $h^* \in \mathbb{R}_+^*$ such that $\{(t, y_{n,0})\} \times [-h^*, h^*] \times \overline{\mathbf{B}}(y_{n,0}, r)^{[\![1,s+1]\!]}$ is a subset of $\mathcal{U}'_{n,f}$, for all $h \in [-h^*, h^*]$, the fixed points of the last coordinate $x$ are inside $\overline{\mathbf{B}}(y_{n,0}, r)$ and the image of $\{(t, y_{n,0}, h)\} \times \overline{\mathbf{B}}(y_{n,0}, r)$ by the parametrized approximation function is contained in itself. This proves the existence of fixed points thanks to Brouwer's theorem.
2. Contraction : We prove that, for all sets of the form $\{(t_0, y_{n,0})\} \times U \times V$ contained in $\mathcal{U}'_{n,d}$ such that $f$ restricted to $\bigcup_{j=1}^{s} \{(t + \tau_j h, x)\}_{(h,x) \in U \times V} \subset \mathbf{U}$ is globally Lipschitz continuous in its second variable, the parametrized approximations function is a contraction in its last coordinate for small enough $h$.

To avoid assuming $|h| \leq 1$ when we need a bound on the powers of $|h|$, we state a small lemma.

**II.1.12 Lemma : Bound on the powers of a variable**

Let $(a,b) \in \mathbb{N}^2$ such that $a \leq b$.
1. For all $h \in \mathbb{R}$, $\max_{k \in [\![a,b]\!]}(|h|^k) = \max(|h|^a, |h|^b)$.
2. For all $(h,c) \in \mathbb{R} \times \mathbb{R}_+^*$, we have $\max(|h|^a, |h|^b) \leq c \Leftrightarrow |h| \leq \min(\sqrt[a]{c}, \sqrt[b]{c})$

Proof

1. If $|h| \leq 1$, we have, for all $k \in [\![a, b-1]\!]$, $\quad |h|^{k+1} \leq |h|^k$, thus the maximum is $|h|^a$ which is equal to $\max\big(|h|^a, |h|^b\big)$. For $|h| \geq 1$ the reasoning is similar with the same conclusion.
2. It is a simple case of considering 4 cases, $c \leq 1, h \leq 1$ or $c \leq 1, h \geq 1$ or $c \geq 1, h \leq 1$ or $c \geq 1, h \geq 1$.

**II.1.13 Lemma : Fixed points of the parametrized approximation function**

Let $(n, s) \in \mathbb{N}^* \times \mathbb{N}^*$, $(\tau, w, \tilde{w}) \in \mathrm{GMORK}_{n,s}$ such that, for all $(j, N) \in [\![1, s+1]\!] \times [\![1, n]\!]$, we have $\tilde{w}_{N,N,j} = 1$, $f$ a differential equation function of order $n$, continuous and locally Lipschitz continuous in its second variable, $\big(t, y_{n,0}\big) \in \mathbf{U}$, $r \in \mathbb{R}_+^*$ such that $\overline{\mathbf{B}}\big(\big(t, y_{n,0}\big), r\big) \subset \mathbf{U}$, $\mathcal{U}'_{n,f}$ the parametrized existence-uniqueness domain of $M$, $Y'_{n,f}$ the parametrized approximation function for length $n$ of $M$, $Y_{n,f}$ the approximation function for length $n$ of $M$, and $\mathcal{H}$ the set of step sizes of $M$.
There exists $h^* \in \mathbb{R}_+^*$ such that :

$$\big\{\big(t, y_{n,0}\big)\big\} \times [-h^*, h^*] \times \overline{\mathbf{B}}\big(y_{n,0}, r\big)^{[\![1,s+1]\!]} \subset \mathcal{U}'_{n,f}$$

$$\forall N \in [\![1, n]\!],\ Y'_{n,f,j,N}\Big(t, y_{n,0}, [-h^*, h^*] \times \overline{\mathbf{B}}\big(\big(t, y_{n,0}\big), r\big)^{[\![1,s+1]\!]}\Big) \subset \overline{\mathbf{B}}\big(\big(t, y_{n,0}\big), r\big)^{[\![1,s+1]\!]}$$

And, for all $h \in [-h^*, h^*]$ the function :

$$\big\{x \in \mathbb{R}^{[\![1,s+1]\!]\times[\![1,n]\!]\times[\![1,d]\!]} \mid \big(t, y_{n,0}, h, x\big) \in \mathcal{U}'_{n,f}\big\} \to \mathbb{R}^{[\![1,s+1]\!]\times[\![1,n]\!]\times[\![1,d]\!]},\ x \to Y'_{n,f}\big(t, y_{n,0}, h, x\big)$$

Has at least a fixed point, and its fixed points are contained in $\overline{\mathbf{B}}\big(\big(t, y_{n,0}\big), r\big)^{[\![1,s+1]\!]}$. We have :

$$Y_{n,f}\Big(t, y_{n,0}, [-h^*, h^*] \cap \mathcal{H}_{1,f,t,y_{n,0}}\Big) \subset \overline{\mathbf{B}}\ \big(\big(t, y_{n,0}\big)\big)^{[\![1,s+1]\!]}$$

- Let $L_r \in \mathbb{R}_+^*$ a Lipschitz constant of $f$ on $\overline{\mathbf{B}}\big(\big(t, y_{n,0}\big), r\big)$.
- Let $h' \in \overline{\mathbb{R}}_+^*$ such that, for all $j \in [\![1, s]\!]$, $\ \tau_j \neq 0 \Rightarrow h' \leq \frac{r}{|\tau_j|}$.
- Let $A = \sup_{(j,h)\in[\![1,s]\!]\times[-h',h']} f\big(t + \tau_j h, y_{n,0}\big)$.
- Let $(C_{\tilde{w}}^*, C_w) \in \mathbb{R}_+ \times \mathbb{R}_+^*$ such that, for all $(j, N) \in [\![1, s+1]\!] \times [\![1, n]\!]$ :

$$C_{\tilde{w}}^* \geq \sum_{N'=1}^{N-1} |\tilde{w}_{N,N',j}|, \quad C_w \geq \sum_{j'=1}^{s} \frac{|w_{N,j,j'}|}{N!}$$

A valid choice of $h^*$ is :

$$h^* = \min\Big(C_r,\ \sqrt[n]{C_r}, h'\Big), \quad C_r = \frac{r}{(rL_r + A)C_w + C_{\tilde{w}}^*\ \|y_{n,0}\|_{n,d}}$$

Proof

Since $\mathbf{U}$ is an open set, there exists a closed ball around $\big(t, y_{n,0}\big)$ of radius $r$ contained in $\mathbf{U}$. We now prove that there exists $h' \in \overline{\mathbb{R}}_+^*$ such that $\big\{\big(t, y_{n,0}\big)\big\} \times [-h', h'] \times \overline{\mathbf{B}}\big(y_{n,0}, r\big)^{[\![1,s+1]\!]} \subset \mathcal{U}'_{n,f}$. Let $(h, x) \in \mathbb{R} \times \mathbb{R}^{[\![1,s+1]\!]\times[\![1,n]\!]\times[\![1,d]\!]}$. $\big(t, y_{n,0}, h, x\big) \in \mathcal{U}'_{n,f}$ if and only if, for all $j \in [\![1, s]\!]$, $\ \big(t + \tau_j h, x_j\big) \in \mathbf{U}$. It is sufficient that, for all $j \in [\![1, s]\!]$, $\ \big(t + \tau_j h, x_j\big) \in \overline{\mathbf{B}}\big(\big(t, y_{n,0}\big), r\big)$. This condition is equivalent to :

$$\|x_j - y_{n,0}\|_{n,d} \leq r, \quad |t + \tau_j h - t| \leq r$$

Which is equivalent to :

$$x_j \in \overline{\mathbf{B}}\big(y_{n,0}, r\big),\ \tau_j \neq 0 \Rightarrow |h| \leq \frac{r}{|\tau_j|}$$

Thus, if we take $h' \in \overline{\mathbb{R}}_+^*$ such that, for all $j \in [\![1, s]\!]$, $\tau_j \neq 0 \Rightarrow h' \leq \frac{r}{|\tau_j|}$, then, for all $h \in \mathbb{R}$ such that $|h| \leq h'$, for all $j \in [\![1, s]\!]$, $\left(t + \tau_j h, x_j\right) \in \overline{\mathbf{B}}\left(\left(t, y_{n,0}\right), r\right)$.

Let $h \in \mathbb{R}$ such that $|h| < h'$. We now find a bound on $\|Y'_{n,f,j,N}\left(t, y_{n,0}, h, x\right) - y_{n,0,N}\|_d$.
For all $(j, N) \in [\![1, s]\!] \times [\![1, n]\!]$ :

$$\|Y'_{n,f,j,N}\left(t, y_{n,0}, h, x\right) - y_{n,0,N}\|_d = \| \sum_{N'=1}^{N-1} \tilde{w}_{N,N',j} h^{N-N'} y_{n,0,N'} + \frac{h^N}{N!} \sum_{j'=1}^{s} w_{N,j,j'} f\left(t + \tau_{j'} h, x_{j'}\right) \|_d$$

$$\leq \sum_{N'=1}^{N-1} |\tilde{w}_{N,N',j}| \; |h|^{N-N'} \; \|y_{n,0,N'}\|_d + |h|^N \sum_{j'=1}^{s} \frac{|w_{N,j,j'}|}{N!} \; \|f\left(t + \tau_{j'} h, y_{n,0}\right)\|_d$$

$$+ \; |h|^N \sum_{j'=1}^{s} \frac{|w_{N,j,j'}|}{N!} \; \|f\left(t + \tau_{j'} h, x_{j'}\right) - f\left(t + \tau_{j'} h, y_{n,0}\right)\|_d$$

The condition $\tilde{w}_{N,N,j} = 1$ is necessary to avoid the constant term $\left(\tilde{w}_{N,N,j} - 1\right) y_{n,0,N}$. To find a bound on $\|f\left(t + \tau_{j'} h, y_{n,0}\right)\|_d$, we use the fact that $f$ is continuous. It follows from the extreme value theorem that we can take $A = \sup_{(j,h) \in [\![1,s]\!] \times [-h',h']} \|f\left(t + \tau_j h, y_{n,0}\right)\|_d$. We have taken $h$ such that, for all $j \in [\![1, s]\!]$, we have $\left(t + \tau_j h, x_j\right) \in \overline{\mathbf{B}}\left(\left(t, y_{n,0}\right), r\right)$, we can thus take $L_r \in \mathbb{R}_+^*$ a Lipschitz constant of $f$ on $\overline{\mathbf{B}}\left(\left(t, y_{n,0}\right), r\right)$. Therefore, for all $(j, N) \in [\![1, s+1]\!] \times [\![1, n]\!]$ :

$$\|Y'_{n,f,j,N}\left(t, y_{n,0}, h, x\right) - y_{n,0,N}\|_d \leq \max(|h|, |h|^n) C_{\tilde{w}}^* \; \|y_{n,0}\|_{n,d} + \max(|h|, |h|^n) \sum_{j'=1}^{s} \frac{|w_{N,j,j'}|}{N!} A$$

$$+ \max(|h|, |h|^n) \sum_{j'=1}^{s} \frac{|w_{N,j,j'}|}{N!} L_r \max_{j' \in [\![1,s]\!]} \|x_{j'} - y_{n,0}\|_{n,d}$$

$$\leq \max(|h|, |h|^n) \left( C_{\tilde{w}}^* \; \|y_{n,0}\|_{n,d} + C_w A + L_r C_w \max_{j' \in [\![1,s]\!]} \|x_{j'} - y_{n,0}\|_{n,d} \right)$$

Assume $x$ is a fixed point. We now prove $x \in \overline{\mathbf{B}}\left(y_{n,0}, r\right)$. We have $Y'_{n,f,j,N}\left(t, y_{n,0}, h, x\right) = x$, thus :

$$\max_{j \in [\![1,s]\!]} \|x_j - y_{n,0}\|_{n,d} \leq \max(|h|, |h|^n) \left( C_{\tilde{w}}^* \; \|y_{n,0}\|_{n,d} + C_w A + L_r C_w \max_{j \in [\![1,s]\!]} \|x_j - y_{n,0}\|_{n,d} \right)$$

Which is equivalent to :

$$\max_{j \in [\![1,s]\!]} \|x_j - y_{n,0}\|_{n,d} \; \left(1 - \max(|h|, |h|^n) L_r C_w\right) \leq \max(|h|, |h|^n) \left( C_{\tilde{w}}^* \; \|y_{n,0}\|_{n,d} + C_w A \right)$$

To simplify the inequality we want to find $a \in \mathbb{R}_+^*$ such that :

$$a \leq 1 - \max(|h|, |h|^n) L_r C_w$$

Which is equivalent to :

$$\max(|h|, |h|^n) \leq \frac{1 - a}{L_r C_w}$$

It is hence necessary and sufficient that $a \in ]0, 1[$. If $h$ satisfies this inequality, it is sufficient that :

$$\max_{j \in [\![1,s]\!]} \|x_j - y_{n,0}\|_{n,d} \; a \leq \max(|h|, |h|^n) \left( C_{\tilde{w}}^* \; \|y_{n,0}\|_{n,d} + C_w A \right)$$

We want $\max_{j\in[\![1,s]\!]}\|x_j - y_{n,0}\|_{n,d} \leq r$. It is sufficient that :

$$\max(|h|,|h|^n)\left(C^*_{\tilde{w}}\ \|y_{n,0}\|_{n,d} + C_w A\right) \leq ar$$

Which is equivalent to :

$$\max(|h|,|h|^n) \leq \frac{ar}{C^*_{\tilde{w}}\ \|y_{n,0}\|_{n,d} + C_w A}$$

Overall, if we want the fixed points to be in the closed ball it is sufficient that :

$$\max(|h|,|h|^n) \leq \min\left(\frac{ar}{C^*_{\tilde{w}}\ \|y_{n,0}\|_{n,d} + C_w A}, \frac{1-a}{L_r C_w}\right)$$

We now find the optimal value of $a$. Since the right term decreases as $a$ increases and the left term does the opposite, a global minimum is reached when both terms are equal, if such equality is possible.

$$\frac{ra}{C^*_{\tilde{w}}\ \|y_{n,0}\|_{n,d} + C_w A} = \frac{1-a}{L_r C_w} \Leftrightarrow a = \frac{C_w A + C^*_{\tilde{w}}\ \|y_{n,0}\|_{n,d}}{(rL_r + A)C_w + C^*_{\tilde{w}}\ \|y_{n,0}\|_{n,d}}$$

$a$ can always take this value since the numerator is strictly inferior to the denominator. We can replace $a$ in the inequality to retrieve the condition $\max(|h|,|h|^n) \leq C_r$, which is equivalent to $|h| \leq \min\left(C_r, \sqrt[n]{C_r}\right)$. If $h$ satisfies this condition, the fixed points are in $\overline{\mathbf{B}}\left(y_{n,0}, r\right)$, we thus restrict the domain of definition to this ball.

We now prove that the image of $\overline{\mathbf{B}}\left(y_{n,0}, r\right)$ is contained in itself. Let $x \in \overline{\mathbf{B}}\left(y_{n,0}, r\right)$, we get :

$$\|Y'_{n,f,j,N}\left(t, y_{n,0}, h, x\right) - y_{n,0,N}\|_d \leq \max(|h|,|h|^n)\left(C^*_{\tilde{w}}\ \|y_{n,0}\|_{n,d} + C_w A + L_r C_w r\right)$$

We want $\|Y'_{n,f,j}\left(t, y_{n,0}, h, x\right) - y_{n,0}\|_{n,d} \leq r$. It is sufficient that :

$$\max(|h|,|h|^n)\left(C^*_{\tilde{w}}\ \|y_{n,0}\|_{n,d} + (L_r r + A)C_w\right) \leq r$$

Which is equivalent to :

$$\max(|h|,|h|^n)\left(C^*_{\tilde{w}}\ \|y_{n,0}\|_{n,d} + AC_w\right) \leq r(1 - L_r C_w \max(|h|,|h|^n))$$

We want to find $a \in \mathbb{R}^*_+$ such that :

$$a \leq 1 - L_r C_w \max(|h|,|h|^n)$$

Which is equivalent to :

$$\max(|h|,|h|^n) \leq \frac{1-a}{L_r C_w}$$

It is hence necessary and sufficient that $a \in ]0,1[$. If $h$ satisfies this inequality, it sufficient that :

$$\max(|h|,|h|^n)\left(C^*_{\tilde{w}}\ \|y_{n,0}\|_{n,d} + AC_w\right) \leq ra$$

We get :

$$\max(|h|,|h|^n) \leq \frac{ra}{C^*_{\tilde{w}}\ \|y_{n,0}\|_{n,d} + AC_w}$$

The condition on $h$ is thus, with $a \in ]0,1[$ :

$$\max(|h|,|h|^n) \leq \min\left(\frac{ra}{C_{\tilde{w}}^* \ \|y_{n,0}\|_{n,d} + C_w A}, \frac{1-a}{L_r C_w}\right)$$

It is the same condition as before, it is thus sufficient that $|h| \leq \min\left(C_r, \sqrt[n]{C_r}\right)$. By adding the condition $|h| \leq h'$ we get the condition of the theorem. Brouwer's fixed point theorem directly gives us the existence of a fixed point. Since the fixed points are in $\overline{\mathbf{B}}\left(\left(t, y_{n,0}\right), r\right)^{[\![1,s+1]\!]}$, if there happens to be uniqueness of the fixed points, we deduce that the resulting approximations are also in the ball.

### II.1.14 Lemma : Contracting parametrized approximation function

Let $(n,s) \in \mathbb{N}^* \times \mathbb{N}^*$, $M = (\tau, w, \tilde{w}) \in \mathrm{GMORK}_{n,s}$, $f$ a differential equation function of order $n$, continuous and locally Lipschitz continuous in its second variable, $\mathcal{U}'_{n,f}$ the parametrized existence-uniqueness domain of $M$, $Y'_{n,f}$ the parametrized approximation function for length $n$ of $M$, $\left(t, y_{n,0}\right) \in \mathbf{U}$, and $(U,V) \in \mathcal{P}(\mathbb{R}) \times \mathcal{P}\left(\mathbb{R}^{[\![1,s+1]\!] \times [\![1,n]\!] \times [\![1,d]\!]}\right)$ two non empty sets such that $\left\{\left(t, y_{n,0}\right)\right\} \times U \times V \subset \mathcal{U}'_{n,f}$ and, if we define $\mathbf{U}' = \bigcup_{j=1}^{s} \left\{\left(t + \tau_j h, x\right)\right\}_{(h,x) \in U \times V} \subset \mathbf{U}$, then $f$ restricted to $\mathbf{U}'$ is globally Lipschitz continuous in its second variable with a constant $L \in \mathbb{R}_+^*$.
There exists $h^* \in \mathbb{R}_+^*$ such that, for all $h \in ]-h^*, h^*[ \cap U$, the function :

$$V \to \mathbb{R}^{[\![1,s+1]\!] \times [\![1,n]\!] \times [\![1,d]\!]}, \ x \to Y'_{n,f}\left(t, y_{n,0}, h, x\right)$$

Is a contraction.
Let $C_w \in \mathbb{R}_+^*$ such that, for all $(j,N) \in [\![1,s+1]\!] \times [\![1,n]\!]$, $C_w \geq \sum_{j'=1}^{s} \frac{|w_{N,j,j'}|}{N!}$.
A valid choice of $h^*$ is $\min\left(\frac{1}{LC_w}, \frac{1}{\sqrt[n]{LC_w}}\right)$.

Proof

Let $(h,x,z) \in U \times V \times V$. For all $(j,N) \in [\![1,s+1]\!] \times [\![1,n]\!]$ :

$$\begin{aligned}
&\|Y'_{n,f,j,N}\left(t, y_{n,0}, h, x\right) - Y'_{n,f,j,N}\left(t, y_{n,0}, h, z\right)\|_d \\
&= \frac{|h|^N}{N!} \ \|\sum_{j'=1}^{s} w_{N,j,j'}\left(f\left(t + \tau_{j'} h, x_{j'}\right) - f\left(t + \tau_{j'} h, z_{j'}\right)\right)\|_d \\
&\leq L \ |h|^N \sum_{j'=1}^{s} \frac{|w_{N,j,j'}|}{N!} \ \|x - z\|_{s+1,n,d} \\
&\leq \max(|h|,|h|^n) L C_w \ \|x - z\|_{s+1,n,d}
\end{aligned}$$

Thus if $\max(|h|,|h|^n) L C_w < 1$, $Y'_{n,f}$ is a contraction, which is equivalent to the condition of the lemma.

Proof of Section II.1.6

In the general case we need both lemmas of Section II.1.13 and Section II.1.14 to restrict $Y'_{n,f}$ to a closed set with its image contained in itself on which $Y'_{n,f}$ is a contraction. In the second lemma we take $U = [-h^*, h^*]$, $V = \overline{\mathbf{B}}\left(y_{n,0}, r\right)^{[\![1,s+1]\!]}$ where $h^*$ has the value of the first lemma, and for the Lipschitz constant $L$ we take $L_r$. The conditions of the second lemma is implied by the stronger condition of the first lemma :

$$\max(|h|,|h|^n) \leq \frac{1-a}{L_r C_w}$$

With $a \in ]0, 1[$, we thus only assume that the conditions of the lemma of Section II.1.13 are satisfied. It is also easy to see that if $f$ is of class $C^k$ then $Y'$ is of class $C^k$. The uniform contraction principle hence directly gives the existence of a unique fixed point and the smoothness of the approximation and evaluation functions.

For the case $\mathbf{U} = \Omega \times \mathbb{R}^{[\![1,n]\!] \times [\![1,d]\!]}$ and $f$ is globally Lipschitz continuous in its second variable, we do not need the lemma of Section II.1.13, and in the lemma of Section II.1.14 we take :

$$U = \left\{ h \in \mathbb{R} \mid \forall j \in [\![1, s]\!],\ t + \tau_j h \in \Omega \right\}$$
$$V = \mathbb{R}^{[\![1,n]\!] \times [\![1,d]\!]}$$

For all $\left(t, y_{n,0}, h\right) \in \mathbf{U} \times \mathbb{R}$, for $x \to Y'\left(t, y_{n,0}, j, x\right)$ to be a contraction in its last variable it is hence sufficient that $h \in ]-h^*, h^*[$ and, for all $j \in [\![1, s]\!],\ t + \tau_j h \in \Omega$. We notice that the choice of $h^*$ given by the lemma does not depend on $\left(t, y_{n,0}\right)$. Using the uniform contraction principle finishes the proof.

We need $f$ to be global Lipschitz continuous in the second case since otherwise we would have to restrict $f$ to a closed ball for it to be globally Lipschitz continuous, which would then force us use the first lemma.

This existence-uniqueness theorem assumes that the stage equations form a single coupled system, but in practice, however, there may be no system to solve, or the stage equations may decompose into several independent systems of smaller size. Determining the number and size of these systems is essential to understanding the computational cost of a method.

For example, methods with lower triangular weight matrix, known as diagonally implicit Runge-Kutta methods[1], require solving $s$ systems of size one instead of a single system of size $s$, which reduces significantly the computational cost.

The next section focuses on a finer study of the stage system of a method, which will lead to a definition of explicit and implicit methods and will allow us to better understand the stage system of implicit multi-order Runge-Kutta method.

## II.2 Implicit methods & Weight digraphs

The definition of explicit Runge-Kutta methods we have given in the section on Runge-Kutta methods is widely accepted, but remains a simplification. Studying why this is leads to a general framework that allows us to create methods with the structure of our choice.

This framework has both theoretical and practical applications. The notion of sub-methods yields useful results on the order of consistency. Explicit partitions of the stages are useful for studying the stability properties of implicit methods. Implicit blocks allow us to reduce the computational cost of implicit methods. Parallel blocks enable parallel implementations. It begins with a simple question, is the following Runge-Kutta method explicit ?

$$\tau = \begin{bmatrix} \frac{1}{2} \\ 1 \\ \frac{1}{2} \\ 0 \end{bmatrix},\ w_1 = \begin{bmatrix} 0 & 0 & \frac{1}{2} & 0 \\ 1 & 0 & 0 & 0 \\ 0 & 0 & 0 & \frac{1}{2} \\ 0 & 0 & 0 & 0 \\ \frac{1}{6} & \frac{1}{3} & \frac{1}{6} & \frac{1}{3} \end{bmatrix}$$

According to the definition of the section on Runge-Kutta methods, the weight matrix is not strictly lower triangular, the method is thus implicit. The only issue is that we can use this method without solving a system. We can first compute the approximation at stage 4, then stage 3, then stage 1, then stage 2 and finally stage 5. In fact, we can notice that by computing the approximations in this order we are actually using the Runge-Kutta method :

$$\tau = \begin{bmatrix} 0 \\ \frac{1}{2} \\ \frac{1}{2} \\ 1 \end{bmatrix},\ w_1 = \begin{bmatrix} 0 & 0 & 0 & 0 \\ \frac{1}{2} & 0 & 0 & 0 \\ 0 & \frac{1}{2} & 0 & 0 \\ 0 & 0 & 1 & 0 \\ \frac{1}{3} & \frac{1}{6} & \frac{1}{6} & \frac{1}{2} \end{bmatrix}$$

Which is the celebrated RK4 method, an explicit Runge-Kutta method. We here formalize this process of changing the order of the stages.

**II.2.1 Definition : Permutations of a method**

Let $(n, s) \in \mathbb{N}^* \times \mathbb{N}^*$. We define $\mathfrak{S}_s^* = \left\{\varphi \in \mathfrak{S}_{[\![1,s+1]\!]} \mid \varphi(s+1) = s+1\right\}$.
We define $*$ the permute operation as, for all $(\varphi, (\tau, w, \tilde{w})) \in \mathfrak{S}_s^* \times \mathrm{GMORK}_{n,s},\ \varphi * (\tau, w, \tilde{w}) = (\tau', w', \tilde{w}') \in \mathrm{GMORK}_{n,s}$, with :

$$\forall j \in [\![1, s]\!],\ \tau'_j = \tau_{\varphi(j)}$$
$$\forall (j, N) \in [\![1, s+1]\!] \times [\![1, n]\!],\ \forall N' \in [\![1, N]\!],\ \tilde{w}'_{N,N',j} = \tilde{w}_{N,N',\varphi(j)}$$
$$\forall (j, j', N) \in [\![1, s+1]\!] \times [\![1, s]\!] \times [\![1, n]\!],\ w'_{N,j,j'} = w_{N,\varphi(j),\varphi(j')}$$

The condition $\varphi(s+1) = s+1$ is necessary because the last stage differs from the other stages by not being used to evaluate $f$.

The permute operation satisfies the axioms of a well known object in group theory, a group action.

**II.2.2 Proposition : Permute operation is a group action**

Let $(n, s) \in \mathbb{N}^* \times \mathbb{N}^*$.
The permute operation $*$ is a left group action of the group $(\mathfrak{S}_s^*, \circ)$ on the set $\mathrm{GMORK}_{n,s}$.

Proof

Trivial.

We now state some fundamental properties of this operation.

**II.2.3 Proposition : Properties of the permute operation**

Let $(n, s) \in \mathbb{N}^* \times \mathbb{N}^*$, $M = (\tau, w, \tilde{w}) \in \mathrm{GMORK}_{n,s}$, and $\varphi \in \mathfrak{S}_s^*$.

1. Let $f$ a differential equation function of order $n$, and $\left(Y_{n,f}, \mathcal{U}_{n,f}\right), \left(Y'_{n,f}, \mathcal{U}'_{n,f}\right)$ the approximations functions and the existence-uniqueness domains of respectively $M$, $\varphi * M$.
   We have :
$$\mathcal{U}_{n,f} = \mathcal{U}'_{n,f}, \quad \forall j \in [\![1, s+1]\!],\ Y'_{n,f,j} = Y_{n,f,\varphi(j)}$$
2. Let $P_\varphi$ the permutation matrix of $\varphi$ restricted to $[\![1, s]\!]$.
   $\varphi * M = (\tau', w', \tilde{w}') \in \mathrm{GMORK}_{n,s}$, with :
$$\tau' = P_\varphi^T \tau,\ \forall N \in [\![1, n]\!],\ \tilde{w}'_N = \tilde{w}_N \begin{bmatrix} P_\varphi & 0_s \\ 0_s^T & 1 \end{bmatrix},\ w'_N = \begin{bmatrix} P_\varphi^T & 0_s \\ 0_s^T & 1 \end{bmatrix} w_N P_\varphi$$

Proof

1. Let $\left(t, y_{n,0}, h\right) \in \mathbf{U} \times \mathbb{R}$. We consider the stage system $\left(M, f, t, y_{n,0}, h\right)$. We wish to prove that it is equivalent up to a permutation to the stage system $\left(\varphi * M, f, t, y_{n,0}, h\right)$. We have :
$$\forall (j, N) \in [\![1, s+1]\!] \times [\![1, n]\!],$$
$$y_{n,j,N} = \sum_{N'=1}^{N} \tilde{w}_{N,N',j} h^{N-N'} y_{n,0,N'} + \frac{h^N}{N!} \sum_{j'=1}^{s} w_{N,j,j'} f\left(t + \tau_{j'} h, y_{n,j'}\right)$$
$$y_{n,s+1} \in \mathbb{R}^{[\![1,n]\!] \times [\![1,d]\!]},\ \forall j \in [\![1, s]\!],\ \left(t + \tau_j h, y_{n,j}\right) \in \mathbf{U}$$
   $\varphi$ is a permutation, and permuting the index $j$ does not change the system, therefore :
$$\forall (j, N) \in [\![1, s+1]\!] \times [\![1, n]\!]$$
$$y_{n,\varphi(j),N} = \sum_{N'=1}^{N} \tilde{w}_{N,N',\varphi(j)} h^{N-N'} y_{n,0,N'} + \frac{h^N}{N!} \sum_{j'=1}^{s} w_{N,\varphi(j),j'} f\left(t_0 + \tau_{j'} h, y_{n,j'}\right)$$
$$y_{n,s+1} \in \mathbb{R}^{[\![1,n]\!] \times [\![1,d]\!]},\ \forall j \in [\![1, s]\!],\ \left(t + \tau_j h, y_{n,j}\right) \in \mathbf{U}$$
   $\varphi$ is a bijection and $\varphi(\{s+1\}) = \{s+1\}$, thus $\varphi([\![1, s+1]\!] \setminus \{s+1\}) = [\![1, s+1]\!] \setminus \{s+1\}$, hence $\varphi([\![1, s]\!]) = [\![1, s]\!]$ and the restriction of $\varphi$ to $[\![1, s]\!]$ is a bijection. Applying the permutation $\varphi$ to the index of the weighted sum therefore does not change it :
$$\forall (j, N) \in [\![1, s+1]\!] \times [\![1, n]\!]$$
$$y_{n,\varphi(j),N} = \sum_{N'=1}^{N} \tilde{w}_{N,N',\varphi(j)} h^{N-N'} y_{n,0,N'} + \frac{h^N}{N!} \sum_{j'=1}^{s} w_{N,\varphi(j),\varphi(j')} f\left(t_0 + \tau_{\varphi(j')} h, y_{n,\varphi(j')}\right)$$
$$y_{n,\varphi(s+1)} \in \mathbb{R}^{[\![1,n]\!] \times [\![1,d]\!]},\ \forall j \in [\![1, s]\!],\ \left(t + \tau_{\varphi(j)} h, y_{n,\varphi(j)}\right) \in \mathbf{U}$$

We define, for all $j \in [\![1, s+1]\!],\ y'_{n,j} = y_{n,\varphi(j)}$. We get the equivalent system :

$$\forall (j, N) \in [\![1, s+1]\!] \times [\![1, n]\!]$$
$$y'_{n,j,N} = \sum_{N'=1}^{N} \tilde{w}'_{N,N',j} h^{N-N'} y_{n,0,N'} + \frac{h^N}{N!} \sum_{j'=1}^{s} w'_{N,j,j'} f\left(t_0 + \tau'_{j'} h, y'_{n,j'}\right)$$
$$y'_{n,s+1} \in \mathbb{R}^{[\![1,n]\!] \times [\![1,d]\!]},\ \forall j \in [\![1, s]\!],\ \left(t + \tau'_j h, y'_{n,j}\right) \in \mathbf{U}$$

We recovered the stage system of $\varphi * M$, which concludes the proof.

2. $P_\varphi$ is defined as :

$$P_\varphi = \left[\begin{cases} 1 \text{ if } r = \varphi(l) \\ 0 \text{ otherwise} \end{cases}\right]_{(r,l) \in [\![1,s]\!]^2} \in \mathbb{R}^{[\![1,s]\!]^2}$$

We can observe that :

$$\forall k \in \mathbb{N}^*,\ \forall M \in \mathbb{R}^{[\![1,s]\!] \times [\![1,k]\!]},\ P_\varphi^T M = \left[M_{\varphi(r),c}\right]_{(r,c) \in [\![1,s]\!] \times [\![1,k]\!]}$$
$$\forall k \in \mathbb{N}^*,\ \forall M \in \mathbb{R}^{[\![1,k]\!] \times [\![1,s]\!]},\ M P_\varphi = \left[M_{r,\varphi(c)}\right]_{(r,c) \in [\![1,k]\!] \times [\![1,s]\!]}$$

And the rest of the proof follows.

Since, as the previous proposition states, applying a permutation to a method simply relabels the stages of a method, it is hence tempting to say that $M$ is equivalent to $\varphi * M$. We thus do exactly that.

### II.2.4 Definition : Equivalence by permutation

Let $(n, s) \in \mathbb{N}^* \times \mathbb{N}^*$.
We define the relation $\mathcal{R}_p$ on $\text{GMORK}_{n,s}$ as the equivalence relation induced by the group action of $\mathfrak{S}_s^*$ on $\text{GMORK}_{n,s}$, thus, for all $M, M' \in \text{GMORK}_{n,s}^2$ :

$$M\ \mathcal{R}_p\ M' \Leftrightarrow \exists \varphi \in \mathfrak{S}_s^*,\ \varphi * M = M'$$

$M$ and $M'$ are said to be equivalent by permutation if and only if $M\ \mathcal{R}_p\ M'$. $M'$ is said to be equivalent to $M$ by permutation $\varphi \in \mathfrak{S}_s^*$ if and only if $M' = \varphi * M$.

The $p$ in $\mathcal{R}_p$ stands for permutation.

For example, the method we discussed earlier is equivalent to RK4 by permutation $(1\ 4\ 2\ 3)$. These two Runge-Kutta methods are therefore equivalent, yet one is said to be explicit while the other is said to be implicit. The real question behind the definitions of explicit and implicit methods is wether or not the approximations can be computed without solving a system. The shape of the matrix, though related, is not the essential point. What truly matters is the dependency between stages, and whether it implies that a stage requires its own value to be computed. To do so we use graph theory. Since conventions vary across the literature, we here state the definitions needed.

### II.2.5 Definition : Digraph

Let $V$ a finite set. A finite directed graph, or digraph, is a pair $(V, A)$ with $A \in \mathcal{P}(V \times V)$. An element of $V$ is a vertex (plur. vertices), and the elements of $A$ are called arcs. An arc $(v', v) \in A$ is said to be directed from $v'$ to $v$, with $v'$ its tail and $v$ its head.

**II.2.6 Definition : Diwalk**

Let $(V, A)$ a digraph, and $r \in \mathbb{N}$. $v \in V^{[\![0,r]\!]}$ is a diwalk of length $r$ if and only if, for all $k \in [\![1, r]\!]$, $(v_{k-1}, v_k) \in A$. $v$ is said to be from $v_0$ to $v_r$. $v$ is said to be closed if and only if $v_0 = v_r$.
For all $v_0 \in V$, $(v_0)$ is a valid diwalk called an empty diwalk. All empty diwalks are closed.

We here use the vertex based definition of a diwalk. It is common to use sequences of arcs of the form $((v_0, v_1), ..., (v_{r-1}, v_r))$ to define diwalks, but the concept of empty diwalks, which is very important here, does not have a natural definition in the arc based definition.

**II.2.7 Definition : Adjacency matrix**

Let $(V, A)$ a digraph. We define the adjacency matrix of $(V, A)$ as :

$$\left[\begin{cases} 1 \text{ if } (j, j') \in A \\ 0 \text{ otherwise} \end{cases}\right]_{(j,j') \in V^2}$$

**II.2.8 Definition : Contraction of a digraph**

Let $G = (V, A)$ a digraph, and $P$ a partition of $V$. Let :

$$G_{/P} = \left(P, A_{/P}\right),\ A_{/P} = \left\{(P_1, P_2) \in P^2 \mid \exists (v_1, v_2) \in (P_1 \times P_2),\ (v_1, v_2) \in A\right\}$$

We define $G_{/P}$ as the contraction of $(V, A)$ according to $P$.
Let $X$ a set, $\sigma \in \{P \to X\}$ such that $\sigma$ is a bijection. Let :

$$G_{/P,\sigma} = \left(P, A_{/P,\sigma}\right),\ A_{/P,\sigma} = \left\{(x_1, x_2) \in X^2 \mid \exists (v_1, v_2) \in \sigma^{-1}(x_1) \times \sigma^{-1}(x_2),\ (v_1, v_2) \in A\right\}$$

We define $G_{/P,\sigma}$ as the contraction of $(V, A)$ according to $P$ with labeling $\sigma$.

We are now able to define the digraph of interest.

**II.2.9 Definition : Maximum weight digraph**

Let $(n, s) \in \mathbb{N}^* \times \mathbb{N}^*$, and $M = (\tau, w, \tilde{w}) \in \mathrm{GMORK}_{n,s}$.
We define $\mathcal{G}_n$ the maximum weight digraph of $M$ as the digraph $([\![1, s+1]\!], \mathcal{A}_n)$, with :

$$\mathcal{A}_n = \left\{(j', j) \in [\![1, s]\!] \times [\![1, s+1]\!] \mid \exists N \in [\![1, n]\!],\ w_{N,j,j'} \neq 0\right\}$$

Its adjacency matrix is :

$$\left[\begin{cases} 1 \text{ if } j' \neq s+1 \text{ and } \exists N \in [\![1, n]\!],\ w_{N,j,j'} \neq 0 \\ 0 \text{ otherwise} \end{cases}\right]_{(j,j') \in [\![1,s+1]\!]^2}$$

The maximum weight digraph of a method is defined to have an arc from $j'$ to $j$ if and only if the approximation of the $j'^{\text{th}}$ stage is directly used in the $j^{\text{th}}$ stage. It follows that a stage $j'$ is used directly or indirectly in stage $j$ if and only if there exists a non-empty diwalk from $j'$ to $j$. Empty diwalks are not taken into acount since dependency is represented by arcs, and empty diwalks contain none.

We deduce that a method is implicit if and only if there exists a non-empty diwalk from a node to itself, that is if and only if the maximum weight digraph has a non-empty closed diwalk.

**II.2.10 Definition : Explicit methods & Implicit methods**

- A general multi-order Runge-Kutta method is said to be implicit if and only if its maximum weight digraph has a non-empty closed diwalk, otherwise, the method is said to be explicit.
- Let $(n, s) \in \mathbb{N}^* \times \mathbb{N}^*$.
  We define $\mathrm{EGMORK}_{n,s}$ as the set of explicit general multi-order Runge-Kutta methods of length $n$ with $s$ points, and $\mathrm{IGMORK}_{n,s}$ as the set of implicit general multi-order Runge-Kutta methods of length $n$ with $s$ points.

For instance, we draw the maximum weight digraphs of the two previous examples.

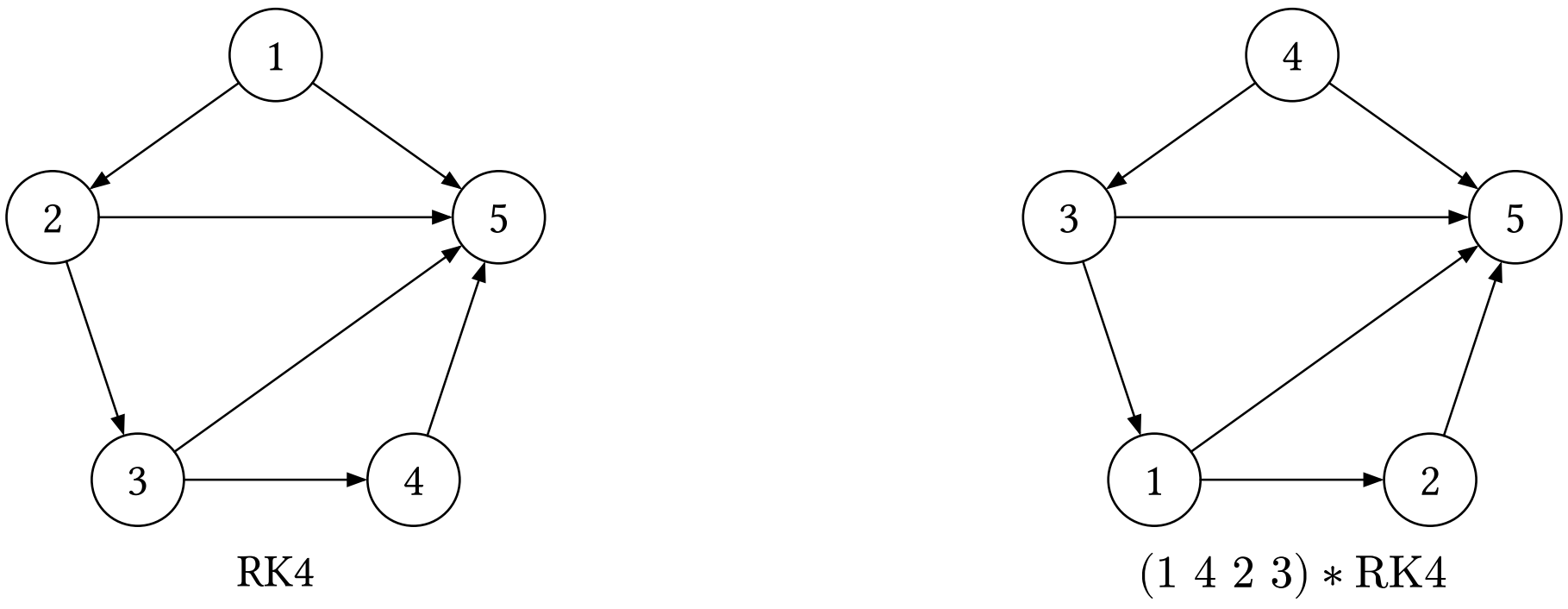

None of those digraphs have non-empty closed diwalks, RK4 and $(1\ 4\ 2\ 3) * \mathrm{RK4}$ are thus explicit methods. To characterize explicit methods, we first define a category of methods which generalizes the definition of explicit methods we gave in the section on Runge-Kutta methods.

**II.2.11 Definition : Topologically sorted methods**

A general multi-order Runge-Kutta method is said to be topologically sorted if and only if its main weight sequence is a sequence of strictly lower triangular matrices.

We now go over some properties of explicit methods.

**II.2.12 Proposition : Properties of explicit methods**

1. A general multi-order Runge-Kutta method is explicit if and only if it is equivalent by permutation to a topologically sorted method.
2. Let $n, s \in \mathbb{N}^* \times \mathbb{N}^*$, $M \in \mathrm{EGMORK}_{n,s}$, $f$ a differential equation function of order $n$, and $\mathcal{U}_{n,f}$ the existence-uniqueness domain of $M$.
   For all $\left(t, y_{n,0}, h\right) \in \mathbf{U} \times \mathbb{R}$, the stage system $\left(M, f, t, y_{n,0}, h\right)$ has at most a solution.
   If $\mathbf{U} = \Omega \times \mathbb{R}^{[\![1,n]\!] \times [\![1,d]\!]}$ then :

$$\mathcal{U}_{n,f} = \left\{\left(t, y_{n,0}, h\right) \in \mathbb{R} \times \mathbb{R}^{[\![1,n]\!] \times [\![1,d]\!]} \times \mathbb{R} \mid \forall j \in [\![1, s]\!],\ t + \tau_j h \in \Omega\right\}$$

Proof

1. Labeling the nodes of a digraph such that its adjacency matrix is strictly lower triangular is called a topological sort. A digraph can be topologically sorted if and only if it has no non-empty closed diwalk[10]. The adjacency matrix of the maximum weight digraph is strictly lower triangular if and only if the main weight sequence is a sequence of strictly lower triangular matrices.
2. Since there is no need to solve a system, the approximations are simply a sequence of computations, the only problem is if this sequence is defined or not. In the case $\mathbf{U} = \Omega \times \mathbb{R}^{[\![1,n]\!] \times [\![1,d]\!]}$, the only constraint is to stay in the interval $\Omega$.

We have not stated a result on the form of $\mathcal{U}_{n,f}$ in the general case where $\mathbf{U} \neq \Omega \times \mathbb{R}^{[\![1,n]\!] \times [\![1,d]\!]}$ because knowing that a method is explicit does not help. If a method is explicit then any fixed point is unique. This allows to skip the contraction part of the proof, which is useless since the existence part is still needed and requires stronger conditions. To remove the existence part we need $\mathbf{U} = \Omega \times \mathbb{R}^{[\![1,n]\!] \times [\![1,d]\!]}$.

The maximum weight digraph can be used for more than defining explicit methods, but to do so we first need to know which stages are needed to compute which stages.

**II.2.13 Definition : Reachability relation & Set of predecessors of a stage**

Let $(n, s) \in \mathbb{N}^* \times \mathbb{N}^*$, and $M \in \mathrm{GMORK}_{n,s}$.

1. We define $\mathcal{R}_r$ the reachability relation of $M$ on $[\![1, s+1]\!]$ as, for all $(j', j) \in [\![1, s+1]\!]^2,\ j' \mathcal{R}_r\ j$ if and only if the the maximum weight digraph of $M$ has a non-empty diwalk from $j'$ to $j$.
2. $\mathcal{R}_r^{-1}(j)$ is called the set of predecessors of the stage $j \in [\![1, s+1]\!]$, it is the set of $j' \in [\![1, s]\!]$ such that there exists a non-empty diwalk from $j'$ to $j$.

The $r$ in $\mathcal{R}_r$ stands for reachability. The stages in $\mathcal{R}_r^{-1}(j)$ are hence the ones needed to compute the stage $j \in [\![1, s+1]\!]$. A simple application of this principle gives :

**II.2.14 Definition : Useless stages & Useful stages**

Let $(n, s) \in \mathbb{N}^* \times \mathbb{N}^*$, $M \in \mathrm{GMORK}_{n,s}$, $\mathcal{R}_r$ the reachability relation of $M$, and $j \in [\![1, s]\!]$.
The stage $j$ is said to be useful if and only if $j \in \mathcal{R}_r^{-1}(s+1)$, otherwise it is said to be useless.

The stages in $[\![1, s]\!] \setminus \mathcal{R}_r^{-1}(s+1)$ are not used in the last stage and are thus useless.

We need the following lemmas for the rest of this subsection.

**II.2.15 Lemma : Some properties of the set of predecessors**

Let $(n, s) \in \mathbb{N}^* \times \mathbb{N}^*$, $M \in \mathrm{GMORK}_{n,s}$, $\varphi \in \mathfrak{S}_s^*$, $j \in [\![1, s+1]\!]$, and $\mathcal{R}_r$, $\mathcal{R}_r{}'$ the reachability relations of respectively $M$, $\varphi * M$.

1. Permutation : $\varphi(\mathcal{R}_r^{-1}(j)) = \mathcal{R}_r'^{-1}(\varphi(j))$.
2. Transitivity : $\mathcal{R}_r$ is transitive and $\mathcal{R}_r^{-1}(\mathcal{R}_r^{-1}(j)) \subset \mathcal{R}_r^{-1}(j)$
3. Unreachability : $\forall (j', j'') \in \left(\mathcal{R}_r^{-1}(j) \cup \{j\}\right) \times \left([\![1, s]\!] \setminus \mathcal{R}_r^{-1}(j)\right),\ \forall N \in [\![1, n]\!],\ w_{N,j',j''} = 0$

Proof

1. Let $j' \in [\![1, s+1]\!]$. The digraph of $M$ has a non-empty diwalk from $j'$ to $j$ if and only if the digraph of $\varphi * M$ has one from $\varphi(j')$ to $\varphi(j)$, which implies this equality.
2. Let $(j,' j'') \in [\![1, s+1]\!]^2$, if there exists a non-empty diwalk from $j'$ to $j$ and a non-empty diwalk from $j''$ to $j'$, then there exists a non-empty diwalk from $j''$ to $j$, $\mathcal{R}_r$ is thus transitive. The inclusion is directly implied by the transitivity of $\mathcal{R}_r$.
3. Let $j \in [\![1, s+1]\!]$. There can be no arc from a stage that is not in $\mathcal{R}_r^{-1}(j)$ to $j$ or any stage in $\mathcal{R}_r^{-1}(j)$, since otherwise that stage would be in $\mathcal{R}_r^{-1}(j)$.

Thanks to these lemmas we can remove useless stages.

**II.2.16 Proposition : Removing useless stages**

Let $(n, s) \in \mathbb{N}^* \times \mathbb{N}^*$, $M \in \mathrm{GMORK}_{n,s}$, $\mathcal{R}_r$ the reachability relation of $M$, $f$ a differential equation function of order $n$, $Y_{n,f}$, $\mathcal{U}_{n,f}$ the approximation function for length $n$ and the existence-uniqueness domain of $M$, $Y'_{n,f}$, $\mathcal{U}'_{n,f}$ the approximation function for length $n$ and the existence-uniqueness domain of $M$ where the useless stages have been removed, and $s' = \#\left(\mathcal{R}_r^{-1}(s+1)\right)$.
For all $\left(t, y_{n,0}, h\right) \in \mathcal{U}_{n,f} \cap \mathcal{U}'_{n,f},\ Y'_{n,f,s'+1}\left(t, y_{n,0}, h\right) = Y_{n,f,s+1}\left(t, y_{n,0}, h\right)$

Proof

Let $\left(t, y_{n,0}, h\right) \in \mathbf{U} \times \mathbb{R}$. The transitivity and unreachability lemmas give that the stage system $\left(M, f, t, y_{n,0}, h\right)$ is equivalent to :

$$\forall (j, N) \in \mathcal{R}_r^{-1}(s+1) \cup \{s+1\} \times [\![1, n]\!],$$

$$y_{n,j,N} = \sum_{N'=1}^{N} \tilde{w}_{N,N',j} h^{N-N'} y_{n,0,N'} + \frac{h^N}{N!} \sum_{j' \in \mathcal{R}_r^{-1}(s+1)} w_{N,j,j'} f\left(t + \tau_{j'} h, y_{n,j'}\right)$$

$$\forall (j, N) \in [\![1, s]\!] \setminus \mathcal{R}_r^{-1}(s+1) \times [\![1, n]\!],$$

$$y_{n,j,N} = \sum_{N'=1}^{N} \tilde{w}_{N,N',j} h^{N-N'} y_{n,0,N'} + \frac{h^N}{N!} \sum_{j'=1}^{s} w_{N,j,j'} f\left(t + \tau_{j'} h, y_{n,j'}\right)$$

$$y_{n,s+1} \in \mathbb{R}^{[\![1,n]\!] \times [\![1,d]\!]},\ \forall j \in [\![1, s]\!],\ \left(t + \tau_j h, y_{n,j}\right) \in \mathbf{U}$$

Removing the useless stages does not change the system of the stages in $\mathcal{R}_r^{-1}(s+1) \cup \{s+1\}$.

In the current literature this is referred to as 0-reducibility[5].

The useful stages form a subsystem of the stage system. More generally, we can characterize subsystems, sets of stages that only require their own values, using the following notion.

**II.2.17 Definition : Reachable-closed sets of stages**

Let $(n, s) \in \mathbb{N}^* \times \mathbb{N}^*$, $M \in \mathrm{GMORK}_{n,s}$, $\mathcal{R}_r$ the reachability relation of $M$, and $J \in \mathcal{P}([\![1, s+1]\!])$. $J$ is said to be closed under reachability, or a reachable-closed set of stages, if and only if $\mathcal{R}_r^{-1}(J) \subset J$.

We can note that $\emptyset$ is closed under reachability. If $J$ is closed under reachability, the stages of $J$ only require themselves to computed, we can hence solve independently their subsystem.

**II.2.18 Proposition : Smallest reachable-closed set of stages**

Let $(n, s) \in \mathbb{N}^* \times \mathbb{N}^*$, $M \in \mathrm{GMORK}_{n,s}$, $\mathcal{R}_r$ the reachability relation of $M$, and $J \in \mathcal{P}([\![1, s+1]\!])$. The smallest set $J' \in \mathcal{P}([\![1, s+1]\!])$ such that $J \subset J'$ and $J'$ is reachable-closed, is $J' = J \cup \mathcal{R}_r^{-1}(J)$.

Proof

Let $J' \in \mathcal{P}([\![1, s+1]\!])$ such that $J \subset J'$ and $J'$ is a reachable-closed set of stages. It is necessary that $\mathcal{R}_r^{-1}(J) \subset J'$, thus $J \cup \mathcal{R}_r^{-1}(J) \subset J'$.
We have $\mathcal{R}_r^{-1}\left(J \cup \mathcal{R}_r^{-1}(J)\right) = \mathcal{R}_r^{-1}(J) \cup \mathcal{R}_r^{-1}\left(\mathcal{R}_r^{-1}(J)\right)$. The transitivity lemma gives $\mathcal{R}_r^{-1}\left(\mathcal{R}_r^{-1}(J)\right) \subset \mathcal{R}_r^{-1}(J)$, thus $\mathcal{R}_r^{-1}\left(J \cup \mathcal{R}_r^{-1}(J)\right) \subset \mathcal{R}_r^{-1}(J) \cup \mathcal{R}_r^{-1}(J) \subset \mathcal{R}_r^{-1}(J)$. $J \cup \mathcal{R}_r^{-1}(J)$ is contained by all reachable-closed sets and it is a reachable-closed set, it is thus the smallest reachable-closed set.

The concept of reachable-closed sets is best understood through the following characterization.

**II.2.19 Theorem : Characterization of reachable-closed sets of stages**

Let $(n, s) \in \mathbb{N}^{*2}$, $M \in \mathrm{GMORK}_{n,s}$, $J \in \mathcal{P}([\![1, s]\!])$, $J \neq \emptyset$, and $\varphi \in \mathfrak{S}_s^*$ such that $\varphi(J) \subset [\![1, \#(J)]\!]$. $J$ is a reachable-closed set of stages if and only if $\varphi * M = (\tau', w', \tilde{w}') \in \mathrm{GMORK}_{n,s}$, where :

$$\forall N \in [\![1, n]\!],\ w'_N = \begin{bmatrix} B_N & 0_{\#(J), s-\#(J)} \\ * & * \end{bmatrix},\ B_N \in \mathbb{R}^{[\![1, \#(J)]\!]^2}$$

Proof

Let $\mathcal{R}_r$, ${\mathcal{R}_r}'$ the reachability relation of $M$, $\varphi * M$. Assume $J$ is a reachable-closed set.
$\mathcal{R}_r^{-1}(J) \subset J$, hence $\varphi\left(\mathcal{R}_r^{-1}(J)\right) \subset \varphi(J)$. The permutation lemma gives that $\mathcal{R}_r'^{-1}(\varphi(J)) \subset [\![1, \#(J)]\!]$, thus $\mathcal{R}_r'^{-1}([\![1, \#(J)]\!]) \subset [\![1, \#(J)]\!]$ and the unreachability lemma allows to conclude.

Assume the main weight sequence is of this form. We have $\mathcal{R}_r'^{-1}([\![1,\#(J)]\!]) \subset [\![1,\#(J)]\!]$ since otherwise there would be a non zero term in the top right. We thus have $\varphi^{-1}\big(\mathcal{R}_r'^{-1}([\![1,\#(J)]\!])\big) \subset \varphi^{-1}([\![1,\#(J)]\!])$. The permutation lemma gives $\mathcal{R}_r^{-1}\big(\varphi^{-1}([\![1,\#(J)]\!])\big) \subset J$. Finally we get $\mathcal{R}_r^{-1}(J) \subset J$.

We can observe that the previous theorem does not cover cases where $s+1 \in J$. This is because $s+1$ cannot be relabeled, thus $\varphi(s+1) \in [\![1,\#(J)]\!]$ if and only if $\#(J) = s+1$, if and only if $J = [\![1,s+1]\!]$. This is not an issue since, if the useless stages have been removed, then $\mathcal{R}_r^{-1}(s+1) = [\![1,s+1]\!]$ and this theorem is useless. The same reasoning holds for the other characterizations of this subsection.

We can observe that a reachable closed set of stages resembles a smaller multi-order Runge-Kutta method.

**II.2.20 Definition : Sub-method induced by a reachable closed set of stages**

Let $(n,s) \in \mathbb{N}^* \times \mathbb{N}^*$, $M = (\tau, w, \tilde{w}) \in \mathrm{GMORK}_{n,s}$, $J \in \mathcal{P}([\![1,s]\!])$ a non-empty reachable closed set of stages, $j_l \in [\![1,s+1]\!] \setminus J$, and $s' = \#(J)$

- Let $\sigma \in \{J \cup \{j_l\} \to [\![1,s'+1]\!]\}$.
  $\sigma$ is said to be a labeling of $(J, j_l)$ if and only if $\sigma$ is a bijection satisfying $\sigma(j_l) = s'+1$.
- Let $\sigma$ a labeling of $(J, j_l)$.
  We define $M' = (\tau', w', \tilde{w}') \in \mathrm{GMORK}_{n,s'}$ the sub-method induced by $(J, j_l, \sigma)$ as :

$$\forall j \in [\![1,s']\!],\ \tau'_j = \tau_{\sigma^{-1}(j)}$$

$$\forall (N,j,j') \in [\![1,n]\!] \times [\![1,s'+1]\!] \times [\![1,s']\!],\ w'_{N,j,j'} = w_{N,\sigma^{-1}(j),\sigma^{-1}(j')}$$

$$\forall (j,N) \in 1,s'+1]\!] \times [\![1,n]\!],\ \forall N' \in [\![1,N]\!],\ \tilde{w}_{N,N',\sigma^{-1}(j)} = \tilde{w}_{N,N',\sigma^{-1}(j)}$$

**II.2.21 Proposition : Properties of sub-methods induced by reachable-closed sets**

Let $(n,s) \in \mathbb{N}^* \times \mathbb{N}^*$, $M \in \mathrm{GMORK}_{n,s}$, $\mathcal{R}_r$ the reachability relation of $M$, $J \in \mathcal{P}([\![1,s]\!])$ a non-empty reachable closed set of stages, $j_l \in [\![1,s+1]\!] \setminus J$, and $s' = \#(J)$.

1. Let $\sigma \in \{J \cup \{j_l\} \to [\![1,s'+1]\!]\}$ a labeling of $(J, j_l)$, $M' \in \mathrm{GMORK}_{n,s}$ the sub-method induced by $(J, j_l, \sigma)$, $f$ a differential equation of order $n$, $Y_{n,f}, Y'_{n,f}$ the approximation functions for length $n$ of respectively $M, M'$, and $\mathcal{U}_{n,f}, \mathcal{U}'_{n,f}$ the existence-uniqueness domains of respectively $M, M'$.
   For all $j \in [\![1,s']\!]$, for all $\big(t, y_{n,0}, h\big) \in \mathcal{U}_{n,f} \cap \mathcal{U}'_{n,f},\ Y'_{n,f,j}\big(t, y_{n,0}, h\big) = Y_{n,f,\sigma^{-1}(j)}\big(t, y_{n,0}, h\big)$.
   If $\mathcal{R}_r^{-1}(j_l) \subset J$, then, for all $\big(t, y_{n,0}, h\big) \in \mathcal{U}_{n,f} \cap \mathcal{U}'_{n,f},\ Y'_{n,f,s'+1}\big(t, y_{n,0}, h\big) = Y_{n,f,j_l}\big(t, y_{n,0}, h\big)$
2. Let $(\sigma, \sigma') \in \{J \cup \{j_l\} \to [\![1,s'+1]\!]\}^2$ two labelings of $(J, j_l)$, and $(M', M'') \in \mathrm{GMORK}^2_{n,s}$ the sub-methods induced respectively by $(J, j_l, \sigma), (J, j_l, \sigma')$.
   There exists $\varphi \in \mathfrak{S}^*_{s'}$ such that $\varphi * M' = M''$.

Proof

1. Let $\big(t, y_{n,0}, h\big) \in \mathbf{U} \times \mathbb{R}$. The unreachability lemma implies that $\big(M, f, t, y_{n,0}, h\big)$ is equivalent to :

$$\forall (j,N) \in J \times [\![1,n]\!],$$

$$y_{n,j,N} = \sum_{N'=1}^{N} \tilde{w}_{N,N',j} h^{N-N'} y_{n,0,N'} + \frac{h^N}{N!} \sum_{j' \in J} w_{N,j,j'} f\big(t + \tau_{j'} h, y_{n,j'}\big)$$

$$\forall (j,N) \in [\![1,s+1]\!] \setminus J \times [\![1,n]\!],$$

$$y_{n,j,N} = \sum_{N'=1}^{N} \tilde{w}_{N,N',j} h^{N-N'} y_{n,0,N'} + \frac{h^N}{N!} \sum_{j'=1}^{s} w_{N,j,j'} f\big(t + \tau_{j'} h, y_{n,j'}\big)$$

$$y_{n,s+1} \in \mathbb{R}^{[\![1,n]\!] \times [\![1,d]\!]},\ \forall j \in [\![1,s]\!],\ \big(t + \tau_j h, y_{n,j}\big) \in \mathbf{U}$$

The first part of the system is equivalent to the stage system of the sub-method for the stages in $[\![1,s']\!]$, which implies, for all $j \in [\![1,s']\!],\ Y'_{n,f,j} = Y_{n,f,\sigma^{-1}(j)}$. The choice of $j_l$ does not matter for the stages in $J$, but if $\mathcal{R}_r^{-1}(j_l) \subset J$, then :

$$\forall (j,N) \in J \cup \{j_l\} \times [\![1,n]\!],$$
$$y_{n,j,N} = \sum_{N'=1}^{N} \tilde{w}_{N,N',j} h^{N-N'} y_{n,0,N'} + \frac{h^N}{N!} \sum_{j' \in J} w_{N,j,j'} f(t + \tau_{j'} h, y_{n,j'})$$
$$\forall (j,N) \in [\![1,s+1]\!] \setminus (J \cup \{j_l\}) \times [\![1,n]\!],$$
$$y_{n,j,N} = \sum_{N'=1}^{N} \tilde{w}_{N,N',j} h^{N-N'} y_{n,0,N'} + \frac{h^N}{N!} \sum_{j'=1}^{s} w_{N,j,j'} f(t + \tau_{j'} h, y_{n,j'})$$
$$y_{n,s+1} \in \mathbb{R}^{[\![1,n]\!] \times [\![1,d]\!]},\ \forall j \in [\![1,s]\!],\ (t + \tau_j h, y_{n,j}) \in \mathbf{U}$$

The first part of this system is equivalent to the stage system $(M', f, t, y_{n,0}, h)$.

2. The permutation $\varphi = \sigma^{-1} \circ \sigma'$ satisfies $\varphi * M' = M''$.

The second proposition implies that the choice of labeling $\sigma$ is meaningless.

The notion of reachable-closed sets is not precise enough since it still treats the stages in $J$ as forming a single system, when it may be further decomposed. For instance, $J$ may contain a smaller reachable-closed set of stages.

Let $j \in [\![1,s]\!]$. We can first note that if there exists a non-empty closed diwalk from $j$ to itself then $j$ requires its own value and is hence part of a system, if not, then $j$ does not require its own value to be computed and thus it is not part of a system. This allows us to make a first distinction.

### II.2.22 Definition : Implicit stages & Explicit stages

Let $(n,s) \in \mathbb{N}^* \times \mathbb{N}^*$, $M \in \mathrm{GMORK}_{n,s}$, and $\mathcal{R}_r$ the reachability relation of $M$.

- $j \in [\![1,s+1]\!]$ is said to be an implicit stage if and only if $j\ \mathcal{R}_r\ j$, otherwise it is said to be an explicit stage.
- We define $\mathcal{J} \in \mathcal{P}([\![1,s]\!])$ the set of implicit stages as $\mathcal{J} = \{j \in [\![1,s]\!] \mid j\ \mathcal{R}_r\ j\}$. $[\![1,s+1]\!] \setminus \mathcal{J}$ is the set of explicit stages.

Given an implicit stage $j \in \mathcal{J}$, we wish to know which stages are part of the system of $j$. Let $j' \in [\![1,s]\!]$. If $j\ \mathcal{R}_r\ j'$ and $j'\ \mathcal{R}_r\ j$ then these two stages need each other values to be computed and are thus part of the same system. Otherwise, at least one of the two stages does not require the other, they are hence not part of the same system. To determine whether two stages belong to the same system, we define the following relation.

### II.2.23 Definition : Dependency relation

Let $(n,s) \in \mathbb{N}^* \times \mathbb{N}^*$, $M \in \mathrm{GMORK}_{n,s}$, and $\mathcal{R}_r$ the reachability relation of $M$.
We define $\mathcal{R}_d$ the dependency relation of $M$ on $[\![1,s+1]\!]$ as, for all $(j,j') \in [\![1,s+1]\!]^2,\ j\ \mathcal{R}_d\ j'$ if and only if $j\ \mathcal{R}_r\ j'$ and $j'\ \mathcal{R}_r\ j$.

The d in $\mathcal{R}_d$ stands for dependency.

**II.2.24 Proposition : Properties of the dependency relation**

Let $(n, s) \in \mathbb{N}^* \times \mathbb{N}^*$, $M \in \mathrm{GMORK}_{n,s}$, and $\mathcal{R}_d$ the dependency relation of $M$.

1. Let $(j, j') \in [\![1, s+1]\!]^2$. $j$ and $j'$ are part of the same system if and only if $j' \mathcal{R}_d j$.
2. Let $j \in [\![1, s]\!]$ an implicit stage. The system that must be solved to compute $j$ is exactly $\mathcal{R}_d^{-1}(j)$.
3. $\mathcal{R}_d$ is transitive and symmetric.
4. $\mathcal{R}_d$ restricted to the set of implicit stages is an equivalence relation.

Proof

1. See the text above the definition of the dependency relation.
2. Since $j$ and $j'$ are part of the same system if and only if $j' \mathcal{R}_d j$, we deduce that the system of $j$ is the set of $j' \in [\![1, s]\!]$ such that $j' \mathcal{R}_d j$, thus $\mathcal{R}_d^{-1}(j)$.
3. It is clear the the dependency relation is symmetric. Let $(j, j', j'') \in [\![1, s+1]\!]^3$. If $j \mathcal{R}_d j'$ and $j' \mathcal{R}_d j''$ then $j \mathcal{R}_r j'$, $j' \mathcal{R}_r j$ and $j' \mathcal{R}_r j''$, $j'' \mathcal{R}_r j'$, which implies $j \mathcal{R}_r j''$ and $j'' \mathcal{R}_r j$, hence $j \mathcal{R}_d j''$.
4. $\mathcal{R}_d$ is reflexive on the set of implicit ranks. Since $\mathcal{R}_d$ is also transitive and symmetric on the set of implicit ranks, it is an equivalence relation.

The dependency relation allows us to find exactly the systems of a method.

**II.2.25 Theorem : Systems of a method**

Let $(n, s) \in \mathbb{N}^* \times \mathbb{N}^*$, $M \in \mathrm{IGMORK}_{n,s}$, $\mathcal{J} \in \mathcal{P}([\![1, s+1]\!])$ the implicit stages of $M$, and $\mathcal{J}_{/\mathcal{R}_d}$ the equivalence classes of the dependency relation when restricted to the implicit stages.

There are $\#\left(\mathcal{J}_{/\mathcal{R}_d}\right)$ independent systems to solve. Each $J \in \mathcal{J}_{/\mathcal{R}_d}$ is a different independent system to solve and it involves $\#(J)$ stages. $\#(J)$ is called the size of the system $J$.

Proof

It is a direct consequence of Section II.2.24.

To know how to compute a given stage we first need a small lemma.

**II.2.26 Lemma : Sets of predecessors of dependant stages**

Let $(n, s) \in \mathbb{N}^* \times \mathbb{N}^*$, $M \in \mathrm{IGMORK}_{n,s}$ and $\mathcal{R}_d, \mathcal{R}_r$ respectively the dependency and reachability relation of $M$.

For all $(j, j') \in [\![1, s+1]\!]^2$, if $j \mathcal{R}_d j'$ then $\mathcal{R}_r^{-1}(j) = \mathcal{R}_r^{-1}(j')$.

Proof

$j \mathcal{R}_d j'$ implies $j \mathcal{R}_r j'$, $j' \mathcal{R}_r j$. Let $j'' \in \mathcal{R}_r^{-1}(j')$. We have $j'' \mathcal{R}_r j'$ and $j' \mathcal{R}_r j$, thus $j'' \mathcal{R}_r j$. We deduce $\mathcal{R}_r^{-1}(j) \subset \mathcal{R}_r^{-1}(j')$, similarly we get $\mathcal{R}_r^{-1}(j') \subset \mathcal{R}_r^{-1}(j)$, thus $\mathcal{R}_r^{-1}(j) = \mathcal{R}_r^{-1}(j')$

**II.2.27 Proposition : Computing a stage**

Let $(n, s) \in \mathbb{N}^* \times \mathbb{N}^*$, $M \in \mathrm{IGMORK}_{n,s}$, $\mathcal{J} \in \mathcal{P}([\![1, s+1]\!])$ the set of implicit stages of $M$, and $j \in [\![1, s+1]\!]$.

If $j$ is an explicit stage then $j$ can be computed without solving a system using the stages in $\mathcal{R}_r^{-1}(j)$.

If $j$ is an implicit stage, it is necessary to first compute the stages in $\mathcal{R}_r^{-1}(j) \setminus \mathcal{R}_d^{-1}(j)$, then solve the system that the stages in $\mathcal{R}_d^{-1}(j)$ form.

Proof

If $j \notin \mathcal{J}$, then $j$ does not require itself to be computed and can be directly computed using the stages in $\mathcal{R}_r^{-1}(j)$. If $j \in \mathcal{J}$, then all the stages in $\mathcal{R}_d^{-1}(j)$ need to be solved for at the same time. We thus need to compute all stages in $\bigcup_{j' \in \mathcal{R}_d^{-1}(j)} \mathcal{R}_r^{-1}(j')$, the previous lemma gives that it is equal to $\mathcal{R}_r^{-1}(j)$. The stages in $\mathcal{R}_r^{-1}(j) \setminus \mathcal{R}_d^{-1}(j)$ are not part of the system since they do not need the values of $\mathcal{R}_d^{-1}(j)$ to be computed, otherwise they would be in $\mathcal{R}_d^{-1}(j)$. We thus have to compute first $\mathcal{R}_r^{-1}(j) \setminus \mathcal{R}_d^{-1}(j)$, then solve for $\mathcal{R}_d^{-1}(j)$.

The relation $\mathcal{R}_d$ is useful, but it is preferable to work with an equivalence relation defined on all stages, so as to avoid distinguishing between explicit and implicit stages. As noted earlier, $\mathcal{R}_d$ is transitive and symmetric, only reflexivity is missing.

**II.2.28 Definition : Closed equivalence**

Let $(n, s) \in \mathbb{N}^* \times \mathbb{N}^*$, $M \in \mathrm{GMORK}_{n,s}$, and $\mathcal{R}_d$ the dependency relation of $M$.
We define $\mathcal{R}_\mathrm{c}$ the closed equivalence on $[\![1, s+1]\!]$ as, for all $(j, j') \in [\![1, s+1]\!]$, $j \; \mathcal{R}_\mathrm{c} \; j'$ if and only if $j = j'$ or $j \; \mathcal{R}_d \; j'$

The c in $\mathcal{R}_\mathrm{c}$ stands for closed.

The only difference with $\mathcal{R}_d$ is that $\mathcal{R}_\mathrm{c}$ cannot distinguish empty diwalks and loops.

**II.2.29 Proposition : Properties of closed equivalence**

Let $(n, s) \in \mathbb{N}^* \times \mathbb{N}^*$, $M \in \mathrm{GMORK}_{n,s}$, and $\mathcal{R}_d$ its closed equivalence.

1. $\mathcal{R}_d$ is an equivalence relation.
2. For all $j \in [\![1, s+1]\!]$, $\mathcal{R}_\mathrm{c}^{-1}(j) = \mathcal{R}_d^{-1}(j) \cup \{j\}$

Proof

1. It is a direct consequence of the symmetry and transitivity of the dependency relation.
2. Trivial.

The closed equivalence relation is well known, though it is usually called strong connectivity. Its equivalence classes are also well known, they are called the Strongly Connected Components (SCC) of a digraph. The SCC of an implicit stage coincides with the equivalence class of $\mathcal{R}_d$ of the same stage. The SCC of a stage is hence the optimal choice for the system of this stage, but it is not the only possible choice. For example, the current literature of Runge-Kutta methods mainly assumes that the stages form one system, even though it is not always optimal. It is also easier to work with a notion that is not always optimal. We now define a generalization of a SCC.

**II.2.30 Definition : Blocks**

Let $(n, s) \in \mathbb{N}^* \times \mathbb{N}^*$, $M \in \mathrm{GMORK}_{n,s}$, $\mathcal{R}_r$ the reachability relation of $M$, and $J \in \mathcal{P}([\![1, s+1]\!])$ a non-empty set.
$J$ is said to be a block if and only if, for all $(j, j', j'') \in J \times [\![1, s]\!] \times J$, $j \; \mathcal{R}_r \; j' \; \mathcal{R}_r \; j'' \Rightarrow j' \in J$.
$\#(J)$ is called the size of the block.

The relation that defines a SCC is, for all $(j, j') \in J \times [\![1, s]\!]$, $j \; \mathcal{R}_r \; j' \; \mathcal{R}_r \; j \Rightarrow j' \in J$, we thus generalized it by making the last $j$ vary in $J$.

**II.2.31 Proposition : Properties of blocks**

Let $(n, s) \in \mathbb{N}^* \times \mathbb{N}^*$, $M \in \mathrm{GMORK}_{n,s}$, $\mathcal{R}_r$ the reachability relation of $M$, and $J \in \mathcal{P}([\![1, s+1]\!])$ a non-empty set.

1. If $J$ is a block then it is an union of SCCs.
2. If $J$ is a block then the stages of $J$ form an independent system.
3. If $J$ is a reachable-closed set then $J$ is a block.
4. $J$ is a block if and only if $\mathcal{R}_r^{-1}(J) \setminus J$ is a reachable-closed set.
5. $J$ is a block if and only if there exists $J' \in \mathcal{P}([\![1, s]\!])$ such that $J \subset J'$, $J'$ is a reachable-closed set and $J' \setminus J$ is a reachable-closed set.

Proof

1. If $J$ is a block, then, for all $j \in J$, by taking $j'' = j$ in the relation that defines a block, we deduce that $\mathcal{R}_c^{-1}(j) \subset J$. Therefore $J = \bigcup_{j \in J} j \subset \bigcup_{j \in J} \mathcal{R}_c^{-1}(j) \subset J$, thus $J = \bigcup_{j \in J} \mathcal{R}_c^{-1}(j)$.
2. This is a direct consequence of 1.
3. Let $(j, j', j'') \in J \times [\![1, s]\!] \times J$. If $j \mathcal{R}_r j' \mathcal{R}_r j''$, then $j' \in \mathcal{R}_r^{-1}(j'')$, which implies $j' \in J$.
4. Assume $J$ is block.
   We have $\mathcal{R}_r^{-1}(J) \setminus J \subset \mathcal{R}_r^{-1}(J)$, thus $\mathcal{R}_r^{-1}(\mathcal{R}_r^{-1}(J) \setminus J) \subset \mathcal{R}_r^{-1}(\mathcal{R}_r^{-1}(J))$. The transitivity lemma gives $\mathcal{R}_r^{-1}(\mathcal{R}_r^{-1}(J) \setminus J) \subset \mathcal{R}_r^{-1}(J)$.
   Let $j \in \mathcal{R}_r^{-1}(\mathcal{R}_r^{-1}(J) \setminus J)$. There exists $j' \in \mathcal{R}_r^{-1}(J) \setminus J$ such that $j \mathcal{R}_r j'$. Since $j' \in \mathcal{R}_r^{-1}(J) \setminus J$, there exists $j'' \in J$ such that $j' \mathcal{R}_r j''$. If $j \in J$ it implies that $j \mathcal{R}_r j' \mathcal{R}_r j''$ with $(j, j'') \in J^2$. Since $J$ is a block it implies that $j' \in J$, which is in contradiction with $j' \in \mathcal{R}_r^{-1}(J) \setminus J$. Thus $j \notin J$. We deduce $\mathcal{R}_r^{-1}(\mathcal{R}_r^{-1}(J) \setminus J) \subset \mathcal{R}_r^{-1}(J) \setminus J$. $\mathcal{R}_r^{-1}(J) \setminus J$ is hence a reachable-closed set.
   Assume $\mathcal{R}_r^{-1}(J) \setminus J$ is a reachable-closed set.
   Let $(j, j', j'') \in J \times [\![1, s]\!] \times J$. Assume $j \mathcal{R}_r j' \mathcal{R}_r j''$. We have $j' \in \mathcal{R}_r^{-1}(J)$. If $j' \in \mathcal{R}_r^{-1}(J) \setminus J$ we have $j \mathcal{R}_r j'$ with $j \in J$, which implies that $\mathcal{R}_r^{-1}(\mathcal{R}_r^{-1}(J) \setminus J) \cap J \neq \emptyset$, this is a contradiction with the fact that $\mathcal{R}_r^{-1}(j) \setminus J$ is closed under reachability, thus $j' \in J$. Thus $J$ is a block.
5. If $J$ is a block we can take $J' = \mathcal{R}_r^{-1}(J)$ thanks to 4.
   Let $J' \in \mathcal{P}([\![1, s]\!])$ such that $J \subset J'$ and $J'$ is a reachable-closed set such that $J' \setminus J$ is a reachable-closed set.
   Let $(j, j', j'') \in J \times [\![1, s]\!] \times J$. Assume $j \mathcal{R}_r j' \mathcal{R}_r j''$. We have $j' \in \mathcal{R}_r^{-1}(J)$, thus $j' \in \mathcal{R}_r^{-1}(J')$, which implies $j \in J'$. If $j \in J' \setminus J$, since $j \in \mathcal{R}_r^{-1}(j')$, this would imply that $\mathcal{R}_r^{-1}(J' \setminus J) \cap J \neq \emptyset$, which is impossible, thus $j' \in J$, thus $J$ is a block

We stated that blocks are independent systems, but a block may not be a system at all if all of the stages of this block are explicit stages.

**II.2.32 Definition : Implicit blocks & Explicit blocks**

Let $(n, s) \in \mathbb{N}^* \times \mathbb{N}^*$, $M \in \mathrm{GMORK}_{n,s}$, $\mathcal{I} \in \mathcal{P}([\![1, s+1]\!])$ the implicit stages of $M$, and $J \in \mathcal{P}([\![1, s+1]\!])$ a block.
$J$ is said to be implicit if and only if $J \cap \mathcal{I} \neq \emptyset$, otherwise it is said to explicit.

The fact that blocks generalize SCCs is apparent in the following proposition.

**II.2.33 Proposition : Smallest block**

Let $(n, s) \in \mathbb{N}^* \times \mathbb{N}^*$, $M \in \mathrm{GMORK}_{n,s}$, $\mathcal{R}_c$ the closed equivalence of $M$, and $j \in [\![1, s+1]\!]$.
The smallest block $J \in \mathcal{P}([\![1, s+1]\!])$ such that $j \in J$ is $\mathcal{R}_c^{-1}(j)$.

Proof

Let $J$ a block such that $j \in J$. Since $J$ is an union of SCCs, we have $\mathcal{R}_c^{-1}(j) \subset J$. Let $(j', j'', j''') \in \mathcal{R}_c^{-1}(j) \times [\![1, s]\!] \times \mathcal{R}_c^{-1}(j)$. If $j' \mathcal{R}_r j'' \mathcal{R}_r j'''$, since $j \mathcal{R}_r j'$ and $j''' \mathcal{R}_r j$ we have $j \mathcal{R}_r j'' \mathcal{R}_r j$, thus $j' \in \mathcal{R}_c^{-1}(j)$. All blocks that contain $j$ contain $\mathcal{R}_c^{-1}(j)$ and $\mathcal{R}_c^{-1}(j)$ is a block that contains $j$, $\mathcal{R}_c^{-1}(j)$ is thus the smallest block that contains $j$.

The notion of blocks is better understood through the following characterization. This characterization also explains why this notion is named block, with the right permutation, a block is a block on the diagonal of all weight matrices.

**II.2.34 Theorem : Characterization of blocks**

Let $(n, s) \in \mathbb{N}^{*2}$, $M \in \mathrm{GMORK}_{n,s}$, $J \in \mathcal{P}([\![1, s]\!])$ a non-empty set, $J' \in \mathcal{P}([\![1, s]\!])$ such that $J \subset J'$, and $\varphi \in \mathfrak{S}_s^*$ such that $\varphi(J') = [\![1, \#(J')]\!]$, $\varphi(J) = [\![1 + \#(J' \setminus J), \#(J')]\!]$.
$J$ is a block and $J'$ is a reachable-closed set such that $J' \setminus J$ is closed under reachability if and only $\varphi * M = (\tau', w', \tilde{w}') \in \mathrm{GMORK}_{n,s}$, where :

$$\forall N \in [\![1, n]\!],\ w'_N = \begin{bmatrix} C_N & 0_{\#(J' \setminus J), \#(J)} & 0_{\#(J' \setminus J), s - \#(J')} \\ * & B_N & 0_{\#(J), s - \#(J')} \\ * & * & * \end{bmatrix}$$

And $(C, B) \in \mathbb{R}^{[\![1,n]\!] \times [\![1, \#(J' \setminus J)]\!]^2} \times \mathbb{R}^{[\![1,n]\!] \times [\![1, \#(J)]\!]^2}$

Proof

Applying twice the characterization of reachable-closed sets of Section II.2.19, with $\varphi$ on $J'$ and on $J' \setminus J$, together with the fifth property of blocks of Section II.2.31 proves the proposition.

The condition on the closed reachability of $J'$ ensures that the stages needed by the block precede it. The condition that $J' \setminus J$ be closed under reachability ensures that any stage not reachable from the block does not depend on it. The existence of such $J'$ when $J$ is a block is established in the fifth proposition of Section II.2.31. The previous theorem characterizes $J$ and $J'$ simultaneously.

Since a block isolates an independent subsystem, a partition of blocks allows to decompose the stage system into multiple independent subsystems.

**II.2.35 Definition : Explicit partition of the stages & Dependency graph**

Let $(n, s) \in \mathbb{N}^* \times \mathbb{N}^*$, $M \in \mathrm{GMORK}_{n,s}$, $\mathcal{G}_n$ the maximum weight digraph of $M$, and $P$ a partition of $[\![1, s]\!]$. We define $\mathcal{G}^*_{n,P}$ as $\mathcal{G}_{n/P \cup \{\{s+1\}\}}$ without the loops.
$P$ is said to be an explicit partition of the stages if and only if $\mathcal{G}^*_{n,P}$ has no non-empty closed diwalk.
$\mathcal{G}^*_{n,P}$ is said to be the maximum dependency graph of $M$ according to $P$.

The link between blocks and explicit partitions of the stages is not clear without the following proposition.

**II.2.36 Proposition : Explicit partition of the stages & Blocks**

Let $(n, s) \in \mathbb{N}^* \times \mathbb{N}^*$, $M \in \mathrm{GMORK}_{n,s}$, and $P$ an explicit partition of the stages.
For all $J \in P$, $J$ is a block.

Proof

Let $J \in P$, and $(j, j', j'') \in J \times [\![1, s]\!] \times J$. Assume $j\ \mathcal{R}_r\ j'\ \mathcal{R}_r\ j''$. If $j' \notin J$, then there exists a non-empty closed diwalk that is not a loop between $J$ and the block which contains $j'$, which is not possible. Thus if $j\ \mathcal{R}_r\ j'\ \mathcal{R}_r\ j''$, then $j' \in J$. $J$ is hence a block.

This notion is called an explicit partition of the stages because it is a way to partition the stages in blocks such that, if we consider each block to be some kind of stage of a method, this method is explicit. The loops are removed in the definition of the dependency digraph because stages within the same block are expected to depend on each others.

**II.2.37 Proposition : Remarkable explicit partition of the stages**

Let $(n, s) \in \mathbb{N}^* \times \mathbb{N}^*$, $M \in \mathrm{GMORK}_{n,s}$, and $\mathcal{R}_c$ the closed equivalence of $M$.

1. $[\![1, s]\!]_{/\mathcal{R}_c}$, the set of SCCs without $\{s + 1\}$, is an explicit partition of the stages.
2. $\{[\![1, s]\!]\}$ is an explicit partition of the stages.

Proof

1. Let $(J, J') \in [\![1,s]\!]^2_{/\mathcal{R}_c}$ such that $J \neq J'$, and $(j, j') \in J \times J'$. If there was a closed diwalk between $J$ and $J'$, then $j \ \mathcal{R}_r \ j' \ \mathcal{R}_r \ j$, therefore $J = J'$, which is impossible. $[\![1,s]\!]_{/\mathcal{R}_c}$ is thus an explicit partition of the stages.
2. There are only two vertices in the dependency graph, $[\![1,s]\!]$ and $\{s+1\}$, and $\{s+1\}$ is not reachable by $[\![1,s]\!]$, thus the dependency graph has no non-empty closed diwalk.

Similarly to explicit methods, assuming that the dependency graph has no non-empty closed diwalks implies that there exists an order in which we can compute the blocks.

### II.2.38 Definition : Order of computation

Let $(n,s) \in \mathbb{N}^* \times \mathbb{N}^*$, $M \in \mathrm{GMORK}_{n,s}$, $\mathcal{G}_n$ the maximum weight digraph of $M$, $B \in \mathcal{P}([\![1,s]\!])^{[\![1,p]\!]}$ a partition of $[\![1,s]\!]$ with $p \in \mathbb{N}^*$, $P = \left\{B_1, ...B_p\right\}$, $\sigma \in \{P \cup \{\{s+1\}\} \to [\![1,p+1]\!]\}$ the only bijection satisfying $\sigma(\{s+1\}) = p+1$ and, for all $(J,k) \in P \times [\![1,p]\!]$, $\sigma(J) = k \Leftrightarrow J = B_k$, and let $\mathcal{G}^*_{n,B}$ be $\mathcal{G}_{n/P\cup\{\{s+1\}\},\sigma}$ without the loops.
$B$ is said to be an order of computation if and only if $P$ is an explicit partition of the stages and the adjacency matrix of $\mathcal{G}^*_{n,B}$ is strictly lower triangular. $\mathcal{G}^*_{n,B}$ is said to be the maximum dependency graph of $M$ according to $B$.

$\mathcal{G}^*_{n,B}$ is simply $\mathcal{G}^*_{n,P}$ with the labeling of the indices of $B$. This means that we can first compute $B_1$, then $B_2$, etc... This is supported by the following proposition.

### II.2.39 Proposition : Order of computations & Reachable closed set of stages

Let $(n,s) \in \mathbb{N}^* \times \mathbb{N}^*$, $M \in \mathrm{GMORK}_{n,s}$, and $B \in \mathcal{P}([\![1,s]\!])^{[\![1,p]\!]}$ an order of computation, with $p \in \mathbb{N}^*$.
For all $k \in [\![1,p]\!]$, $\bigcup_{k'=1}^{k} B_{k'}$ is a reachable-closed set of stages.

Proof

Let $k \in [\![1,p]\!]$, $j \in \bigcup_{k'=1}^{k} B_{k'}$, and $j' \in [\![1,s+1]\!]$. If $j' \ \mathcal{R}_r \ j$ then $j' \in \bigcup_{k'=1}^{k} B_{k'}$ since otherwise the adjacency matrix of the dependency graph of $B$ would not be lower triangular.

### II.2.40 Proposition : Existence of an order of computation

Let $(n,s) \in \mathbb{N}^* \times \mathbb{N}^*$, $M \in \mathrm{GMORK}_{n,s}$, and $P$ an explicit partition of the stages.
There exists an order of computation $B \in \mathcal{P}([\![1,s]\!])^{[\![1,\#(P)]\!]}$ such that $P = \left\{B_1, ...B_{\#(P)}\right\}$.

Proof

The dependency digraph has no non-empty closed diwalk, we can thus topologically sort the digraph such that its adjacency matrix is strictly lower triangular.

### II.2.41 Theorem : Characterization of order of computations

Let $(n,s) \in \mathbb{N}^* \times \mathbb{N}^*$, $M \in \mathrm{GMORK}_{n,s}$, $B \in \mathcal{P}([\![1,s]\!])^{[\![1,p]\!]}$ a partition of $[\![1,s]\!]$, $\varphi \in \mathfrak{S}^*_s$ such that, for all $k \in [\![1,p]\!]$, $\varphi(B_k) \subset [\![1 + \sum_{k'=1}^{k-1} \#(B_{k'}), \sum_{k'=1}^{k} \#(B_{k'})]\!]$
$B$ is an order of computation if and only if $\varphi * M = (\tau', w', \tilde{w}') \in \mathrm{GMORK}_{n,s}$, where :

$$\forall N \in [\![1,n]\!], \ w'_N = \begin{bmatrix} A_{1,N} & 0_{\#(B_1),\#(B_2)} & \cdots & 0_{\#(B_1),\#(B_p)} \\ * & A_{2,N} & \ddots & \vdots \\ \vdots & \ddots & \ddots & 0_{\#(B_{p-1}),\#(B_p)} \\ * & \dots & * & A_{p,N} \\ * & \dots & * & * \end{bmatrix}$$

And, for all $k \in [\![1,p]\!]$, $A_k \in \mathbb{R}^{[\![1,n]\!]\times[\![1,\#(B_k)]\!]^2}$.

Proof

The shape of the adjacency matrix is the same as that of the weight matrices.

This theorem can be proven easily, however, to describe interactions between multiple blocks or between the stages within a block, we need a general way to analyze the structure of arbitrary subset of stages. Sub-methods already allow this this to some extent, but their definition can be extended to non reachable-closed sets. This extension makes it possible to apply the previous theorems to any subset of stages.

**II.2.42 Definition : Sub-method induced by a set of stages**

Let $(n,s) \in \mathbb{N}^* \times \mathbb{N}^*$, $M = (\tau, w, \tilde{w}) \in \mathrm{GMORK}_{n,s}$, $J \in \mathcal{P}(\llbracket 1, s \rrbracket)$ a non-empty set, $s' = \#(J)$, and $j_l \in \llbracket 1, s+1 \rrbracket \setminus J$.

- Let $\sigma \in \{J \cup \{j_l\} \to \llbracket 1, s'+1 \rrbracket\}$.
  $\sigma$ is said to be a labeling of $(J, j_l)$ if and only if $\sigma$ is a bijection such that $\sigma(j_l) = s'+1$.
- Let $\sigma$ a labeling of $(J, j_l)$.
  We define $M' = (\tau', w', \tilde{w}') \in \mathrm{GMORK}_{n,s'}$ the sub-method of $M$ induced by $(J, j_l, \sigma)$ as :

$$\forall j \in \llbracket 1, s' \rrbracket,\ \tau'_j = \tau_{\sigma^{-1}(j)}$$

$$\forall (N, j, j') \in \llbracket 1, n \rrbracket \times \llbracket 1, s'+1 \rrbracket \times \llbracket 1, s' \rrbracket,\ w'_{N,j,j'} = w_{N,\sigma^{-1}(j),\sigma^{-1}(j')}$$

$$\forall (j, N) \in 1, s'+1 \rrbracket \times \llbracket 1, n \rrbracket,\ \forall N' \in \llbracket 1, N \rrbracket,\ \tilde{w}'_{N,N',\sigma^{-1}(j)} = \tilde{w}_{N,N',\sigma^{-1}(j)}$$

**II.2.43 Proposition : Properties of sub-method induced by a set of stages**

Let $(n,s) \in \mathbb{N}^* \times \mathbb{N}^*$, $M \in \mathrm{GMORK}_{n,s}$, $J \in \mathcal{P}(\llbracket 1, s \rrbracket)$ a non-empty set, $j_l \in \llbracket 1, s+1 \rrbracket \setminus J$, $s' = \#(J)$, $\sigma$ a labeling of $(J, j_l)$, and $M'$ the sub-method of $M$ induced by $(J, j_l, \sigma)$.

1. The subgraph induced by $J$ of the maximum weight digraph of $M$ with the labeling $\sigma$ is the maximum weight digraph of $M'$ where the vertex $s'+1$ has been removed.
2. Any permutation $\varphi \in \mathfrak{S}^*_{s'}$ applied to $M'$ can be applied to the stages $J$ of $M$ using the permutation $\varphi' \in \mathfrak{S}^*_s$ defined by :

$$\forall j \in \llbracket 1, s+1 \rrbracket \setminus J,\ \varphi'(j) = j$$

$$\forall j \in J,\ \varphi'(j) = \sigma^{-1}(\varphi(\sigma(j)))$$

3. Let $\sigma' \in \{J \cup \{j_l\} \to \llbracket 1, s'+1 \rrbracket\}$ a labeling of $(J, j_l)$, and $M'' \in \mathrm{GMORK}_{n,s}$ the sub-method of $M$ induced by $(J, j_l, \sigma')$.
   There exists $\varphi \in \mathfrak{S}^*_{s'}$ such that $\varphi * M' = M''$.

Proof

1. Let $(j, j') \in J$, and $\mathcal{A}_n$ the arcs of the maximum weight digraph of $M$. $(j', j) \in \mathcal{A}_n$ if and only if there exists $N \in \llbracket 1, n \rrbracket$ such that $w_{N,j,j'} \neq 0$, if and only if there exists $N \in \llbracket 1, n \rrbracket$ such that $w_{N,\sigma^{-1}(\sigma(j)),\sigma^{-1}(\sigma(j'))} \neq 0$, if and only if there exists $N \in \llbracket 1, n \rrbracket$ such that $w'_{N,\sigma(j),\sigma(j')} \neq 0$. We deduce from this the characterization of the subgraph induced by $J$.
2. Let $\varphi \in \mathfrak{S}^*_{s'}$, and $\varphi * M' = M'' = (\tau'', w'', \tilde{w}'') \in \mathrm{GMORK}_{n,s'}$, we have :

$$\forall j \in \llbracket 1, s \rrbracket,\ \tau''_j = \tau'_{\varphi(j)}$$

$$\forall (j, N) \in \llbracket 1, s+1 \rrbracket \times \llbracket 1, n \rrbracket,\ \forall N' \in \llbracket 1, N \rrbracket,\ \tilde{w}''_{N,N',j} = \tilde{w}'_{N,N',\varphi(j)}$$

$$\forall (N, j, j') \in \llbracket 1, n \rrbracket \times \llbracket 1, s+1 \rrbracket \times \llbracket 1, s \rrbracket,\ w''_{N,j,j'} = w'_{N,\varphi(j),\varphi(j')}$$

   Therefore :

$$\forall j \in J,\ \tau''_{\sigma(j)} = \tau'_{\varphi(\sigma(j))} = \tau_{\sigma^{-1}(\varphi(\sigma(j)))}$$
$$\forall (j,N) \in J \cup \{j_l\} \times [\![1,n]\!],\ \forall N' \in [\![1,N]\!],\ \tilde{w}''_{N,N',\sigma(j)} = \tilde{w}'_{N,N',\varphi(\sigma(j))} = \tilde{w}_{N,N',\sigma^{-1}(\varphi(\sigma(j)))}$$
$$\forall (N,j,j') \in [\![1,n]\!] \times J \cup \{j_l\} \times J,\ w''_{N,\sigma(j),\sigma(j')} = w'_{N,\varphi(\sigma(j)),\varphi(\sigma(j'))} = w_{N,\sigma^{-1}(\varphi(\sigma(j))),\sigma^{-1}(\varphi(\sigma(j')))}$$

We deduce the expression of $\varphi'$ from those relations.

3. The permutation $\varphi = \sigma^{-1} \circ \sigma'$ satisfies $\varphi * M' = M''$.

Using the knowledge we now have on the structure of the stage system, we can improve the contraction part of the existence-uniqueness proof. Once again, this knowledge is not useful if the existence part of the proof is needed, we hence assume $\mathbf{U} = \Omega \times \mathbb{R}^{[\![1,n]\!] \times [\![1,d]\!]}$ and $f$ is globally Lipschitz continuous.

**II.2.44 Theorem : Existence and uniqueness of approximations - Implicit blocks**

We consider the context of the theorem of Section II.1.6, and we assume the method is implicit. We consider the case 2 where $\mathbf{U} = \Omega \times \mathbb{R}^{[\![1,n]\!] \times [\![1,d]\!]}$ and $f$ is globally Lipschitz continuous with a Lipschitz constant $L$.

- Let $P$ an explicit partition of the stages.
- Let $P_\mathcal{J}$ the subset of $P$ containing only the implicit blocks.
- Let, for all $J \in P_\mathcal{J},\ C_{w,J} \in \mathbb{R}_+^*$ such that $C_{w,J} \geq \max_{(j,N) \in J \times [\![1,n]\!]} \sum_{j' \in J} \frac{|w_{N,j,j'}|}{N!}$.

A valid choice of $h^*$ is :

$$h^* = \min_{J \in P_\mathcal{J}} \min\left(\frac{1}{LC_{w,J}}, \frac{1}{\sqrt[n]{LC_{w,J}}}\right)$$

Proof

To prove this theorem we use the same reasoning as for the first proof, except it is repeated for each block of the partition. If the block $J$ considered is explicit we directly have uniqueness and existence of solutions, if $J$ is implicit we can write $\sum_{j=1}^{s}$ as $\sum_{j \in [\![1,s]\!] \setminus J}^{s} + \sum_{j \in J}^{s}$. The first part is constant since the evaluations are either not used or already known, we can then prove that the fixed point form of this system is a contraction if the step size $h$ satisfies $|h| < \min\left(\frac{1}{LC_{w,J}}, \frac{1}{\sqrt[n]{LC_{w,J}}}\right)$.

The use of this theorem is pushed to its limit with the generalization of diagonally implicit Runge-Kutta methods.

**II.2.45 Definition : Diagonally implicit methods**

Let $(n,s) \in \mathbb{N}^* \times \mathbb{N}^*$, $M \in \mathrm{GMORK}_{n,s}$, and $\mathcal{R}_c$ the closed equivalence of $M$.
$M$ is said to be diagonally implicit if and only if $M$ is implicit and $[\![1,s]\!]_{/\mathcal{R}_c} = \{\{1\}, ...\{s\}\}$.

This is equivalent to : the only non-empty closed diwalks of $M$ are loops.

**II.2.46 Proposition : Form of diagonally implicit methods**

Let $(n,s) \in \mathbb{N}^* \times \mathbb{N}^*$, and $M \in \mathrm{GMORK}_{n,s}$.
$M$ is diagonally implicit if and only if $M$ is implicit and equivalent by permutation to a general multi-order Runge-Kutta method whose main weight sequence is a sequence of lower triangular matrices.

Proof

It is a direct consequence of the characterization of order of computations.

Up to this point, implicit stages were identified by grouping together all ranks $N$ at a given stage $j \in [\![1,s+1]\!]$. However, for example, if a rank $N \in [\![1,n]\!]$ satisfies, for all $j' \in [\![1,s]\!],\ w_{N,j,j'} = 0$, the rank $N$ at stage $j$ does not depend on any stage, and no system needs to be solved to compute it. There thus exists explicit and implicit ranks.

**II.2.47 Definition : Implicit ranks & Explicit ranks**

Let $(n,s) \in \mathbb{N}^* \times \mathbb{N}^*$, $M = (\tau, w, \tilde{w}) \in \mathrm{GMORK}_{n,s}$, $(j,N) \in [\![1, s+1]\!] \times [\![1,n]\!]$, and $\mathcal{R}_{\mathrm{c}}$ the closed equivalence of $M$.
$N$ is said to be an implicit rank at stage $j$ if and only if $j$ is an implicit stage and there exists $j' \in \mathcal{R}_{\mathrm{c}}^{-1}(j)$ such that $w_{N,j,j'} \neq 0$, otherwise the rank $N$ is said to be explicit at stage $j$. We define $\mathcal{D}_{n,j}$ as the set of implicit ranks at stage $j \in [\![1, s+1]\!]$. $[\![1,n]\!] \setminus \mathcal{D}_{n,j}$ is the set of explicit ranks at stage $j$.

**II.2.48 Proposition : Use of implicit ranks & Explicit ranks**

Let $(n,s) \in \mathbb{N}^* \times \mathbb{N}^*$, $M \in \mathrm{GMORK}_{n,s}$, and $(j,N) \in [\![1, s+1]\!] \times [\![1,n]\!]$.
If $N$ is an explicit rank at stage $j$ it can be computed directly from the stages in $\mathcal{R}_r^{-1}(j) \setminus \mathcal{R}_{\mathrm{c}}^{-1}(j)$. If $N$ is an implicit rank at stage $j$ it needs to be solved simultaneously with the subsystem of $\mathcal{R}_{\mathrm{c}}^{-1}(j)$.

Proof

If $j$ is an explicit stage it does not require any value in $\mathcal{R}_{\mathrm{c}}^{-1}(j)$ and can thus be computed directly from the stages in $\mathcal{R}_r^{-1}(j) \setminus \mathcal{R}_{\mathrm{c}}^{-1}(j)$. If it requires a stage in $\mathcal{R}_{\mathrm{c}}^{-1}(j)$ then it requires itself and the rest of the subsystem to be computed.

This notion allows to improve the lower bound on $h^*$.

**II.2.49 Theorem : Existence and uniqueness of approximations - Implicit ranks**

We consider the context of the theorem of Section II.1.6, and we assume the method is implicit. We consider the case 2 where $\mathbf{U} = \Omega \times \mathbb{R}^{[\![1,n]\!] \times [\![1,d]\!]}$ and $f$ is globally Lipschitz continuous with a Lipschitz constant $L$.

- Let $P$ an explicit partition of the stages.
- Let $P_{\mathcal{J}}$ the subset of $P$ containing only the implicit blocks.
- Let $\mathcal{D}_n$ the sets of implicit ranks of $M$.
- Let, for all $J \in P_{\mathcal{J}}$, $\left(C_{w,J}, a_J, b_J\right) \in \mathbb{R}_+^* \times \mathbb{N}^* \times \mathbb{N}^*$ such that:

$$C_{w,J} \geq \max_{j \in J} \max_{N \in \mathcal{D}_{n,j}} \sum_{j' \in J} \frac{|w_{N,j,j'}|}{N!}, \ a_J \geq \max_{j \in J} \max \mathcal{D}_{n,j}, \ b_J \leq \min_{j \in J} \min \mathcal{D}_{n,j}$$

A valid choice of $h^*$ is :

$$h^* = \min_{J \in P_{\mathcal{J}}} \min\left( \frac{1}{\sqrt[b_J]{(LC_{w,J})}}, \frac{1}{\sqrt[a_J]{LC_{w,J}}} \right)$$

Proof

The proof is essentially the same as that of the previous theorem which takes into account implicit blocks, except the approximations corresponding to explicit ranks are not included in the system. We introduce $a_J$ and $b_J$ to denote the extremal exponents. Indeed, for a block $J \in P_{\mathcal{J}}$, the highest power of $h$ appearing in an implicit rank is less than or equal to $h^{a_J}$, while the lowest power is greater than or equal to $h^{b_J}$. We may use the bound $\max\left(h^{b_J}, h^{a_J}\right)$ instead of $\max(h, h^n)$, which yields a sharper estimate.

Orders of computations are useful to reduce the cost of implicit methods, but since they only allow to compute the blocks sequentially, they do not use the multithreaded architecture of modern computers. Instead of using an order of computation we can simply use a digraph traversal algorithm that makes full use of a dependency graph. To detect when two blocks do not need to be computed one after the other, we introduce the concept of parallel blocks. If one block requires another block to be computed they cannot be considered parallel, and it is sufficient that neither require the other to call them parallel blocks.

**II.2.50 Definition : Parallel blocks**

Let $(n, s) \in \mathbb{N}^* \times \mathbb{N}^*$, $M \in \text{GMORK}_{n,s}$, and $J_1, J_2 \in \mathcal{P}([\![1, s+1]\!])^2$ two disjoint blocks.
$J_1, J_2$ are said to be parallel if and only if $\mathcal{R}_r^{-1}(J_1) \cap J_2 = \emptyset$ and $\mathcal{R}_r^{-1}(J_2) \cap J_1 = \emptyset$.

Two parallel blocks do not need each others value, and thus, once the stages they depend on have been computed, we can compute them separately.

**II.2.51 Theorem : Characterization of parallel blocks**

Let $(n, s) \in \mathbb{N}^* \times \mathbb{N}^*$, $M \in \text{GMORK}_{n,s}$, $(J_1, J_2) \in \mathcal{P}([\![1, s]\!])^2$ two disjoint blocks, $J' \in \mathcal{P}([\![1, s]\!])$ such that $J_1 \cup J_2 \subset J'$, and $\varphi \in \mathfrak{S}_s^*$ such that :

$$\begin{aligned}
\varphi(J_2) &\subset [\![1 + \#(J' \setminus J_2), \#(J')]\!] \\
\varphi(J_1) &\subset [\![1 + \#(J' \setminus (J_1 \cup J_2)), \#(J' \setminus J_2)]\!] \\
\varphi(J' \setminus (J_1 \cup J_2)) &\subset [\![1, \#(J' \setminus (J_1 \cup J_2))]\!]
\end{aligned}$$

$J_1, J_2$ are independent blocks and $J'$ is a reachable-closed set such that $J' \setminus (J_1 \cup J_2)$ is a reachable-closed set if and only if $\varphi * M = (\tau', w', \tilde{w}') \in \text{GMORK}_{n,s}$, where :

$$\forall N \in [\![1, n]\!],\ w'_N = \begin{bmatrix} * & 0_{a,\#(J_1)} & 0_{a,\#(J_2)} & 0_{a,s-\#(J')} \\ * & B_N & 0_{\#(J_1),\#(J_2)} & 0_{\#(J_1),s-\#(J')} \\ * & 0_{\#(J_2),\#(J_1)} & C_N & 0_{\#(J_2),s-\#(J')} \\ * & * & * & * \end{bmatrix}$$

And $a = \#(J' \setminus (J_1 \cup J_2)),\ (B, C) \in \mathbb{R}^{[\![1,n]\!] \times [\![1,\#(J_1)]\!]^2} \times \mathbb{R}^{[\![1,n]\!] \times [\![1,\#(J_2)]\!]^2}$

Proof

Assume $J_1, J_2$ are independent blocks and $J'$ is a reachable-closed set such that $J' \setminus (J_1 \cup J_2)$ is a reachable-closed set.
$J'$ is a reachable closed set such that $J_1 \cup J_2 \subset J'$, and $\mathcal{R}_r^{-1}(J_1) \cap J_2 = \emptyset$, $\mathcal{R}_r^{-1}(J_2) \cap J_1 = \emptyset$, thus $\mathcal{R}_r^{-1}(J_1) \subset J' \setminus J_2$ and $\mathcal{R}_r^{-1}(J_2) \subset J' \setminus J_1$. $J' \setminus J_2$ is a reachable-closed set since $\mathcal{R}_r^{-1}(J' \setminus J_2) = \mathcal{R}_r^{-1}(J' \setminus (J_2 \cup J_1)) \cup \mathcal{R}_r^{-1}(J_1) \subset J' \setminus (J_2 \cup J_1) \cup J' \setminus J_2 \subset J' \setminus (J_2)$. We then use twice the characterization of blocks plus the fact that $\mathcal{R}_r^{-1}(J_2) \cap J_1 = \emptyset$ to find the form of the theorem.

Assume $\varphi * M$ has the form of the theorem.
The characterization of reachable closed sets give that $J'$, $J' \setminus J_2$ and $J' \setminus (J_1 \cup J_2)$ are reachable-closed sets. The characterization of blocks also give that $J_1$ and $J_2$ are blocks. Since $J' \setminus J_2$ contains $J_1$ and is a reachable-closed set we have $\mathcal{R}_r^{-1}(J_1) \cap J_2 = \emptyset$. There can be no diwalk from a stage of $J_2$ passing through a stage of $J' \setminus (J_1 \cup J_2)$ ending in $J_1$ since $J' \setminus (J_1 \cup J_2)$ is a reachable closed set. Thus if there exists a diwalk form a stage of $J_2$ to $J_1$ it is only composed of stages of $J_1$ and $J_2$. This diwalk cannot exist since otherwise there would be a non zero term at the left of $C_N$ in the matrix. We deduce that $\mathcal{R}_r^{-1}(J_2) \cap J_1 = \emptyset$. $J_1$ and $J_2$ are thus parallel blocks.

A stronger concept than that of parallel blocks is independent sets of stages, two parallel blocks that are each reachable-closed sets of stages. The main idea is that no stage used in the computation of one block is used in the computation of the other, they can hence be computed completely entirely independently.

**II.2.52 Definition : Independent sets of stages**

Let $(n, s) \in \mathbb{N}^* \times \mathbb{N}^*$, $M \in \mathrm{GMORK}_{n,s}$, and $(J_1, J_2) \in \mathcal{P}([\![1, s+1]\!])^2$.
$J_1, J_2$ are said to be independent sets of stages if and only if $J_1, J_2$ are reachable-closed sets of stages and parallel blocks.

There is a simpler characterization of independent sets of stages.

**II.2.53 Proposition : A simpler definition of independent sets of stages**

Let $(n, s) \in \mathbb{N}^* \times \mathbb{N}^*$, $M \in \mathrm{GMORK}_{n,s}$, and $(J_1, J_2) \in \mathcal{P}([\![1, s+1]\!])^2$.
$J_1, J_2$ are independent sets of stages if and only if $J_1, J_2$ are disjoint and reachable-closed sets of stages.

Proof

By definition $J_1, J_2$ are independent sets of stages if and only if $J_1, J_2$ are reachable-closed sets of stages, $J_1 \cap J_2 = \emptyset$, $\mathcal{R}_r^{-1}(J_1) \cap J_2 = \emptyset$, and $\mathcal{R}_r^{-1}(J_2) \cap J_1 = \emptyset$. Since $\mathcal{R}_r^{-1}(J_2) \subset J$ and $\mathcal{R}_r^{-1}(J_1) \subset J_1$, the last two conditions are redundant.

**II.2.54 Theorem : Characterization of independent sets of stages**

Let $(n, s) \in \mathbb{N}^* \times \mathbb{N}^*$, $M \in \mathrm{GMORK}_{n,s}$, $(J_1, J_2) \in \mathcal{P}([\![1, s]\!])^2$ two non-empty disjoint sets, and $\varphi \in \mathfrak{S}_s^*$ such that $\varphi(J_1) = [\![1, \#(J_1)]\!]$ and $\varphi(J_2) = [\![\#(J_1)+1, \#(J_1)+\#(J_2)]\!]$.
$J_1, J_2$ are independent sets of stages if and only if $\varphi * M = (\tau', w', \tilde{w}') \in \mathrm{GMORK}_{n,s}$, where :

$$\forall N \in [\![1, n]\!],\ w'_N = \begin{bmatrix} B_N & 0_{\#(J_1),\#(J_2)} & 0_{\#(J_1),s-\#(J_1 \cup J_2)} \\ 0_{\#(J_2),\#(J_1)} & C_N & 0_{\#(J_2),s-\#(J_1 \cup J_2)} \\ * & * & * \end{bmatrix}$$

And $(B, C) \in \mathbb{R}^{[\![1,n]\!] \times [\![1,\#(J_1)]\!]^2} \times \mathbb{R}^{[\![1,n]\!] \times [\![1,\#(J_2)]\!]^2}$

Proof

This is a special case of the characterization of parallel blocks, with $J' = J_1 \cup J_2$.

The concept of parallel blocks and independent sets of stages are most useful when the considered blocks are implicit. The strictly lower triangular form of explicit blocks implies that at least two stages use the exact same informations to make an approximation, which is not ideal. For independent explicit sets of stages this is worsened since it implies that at least two stages do not use any evaluation of $f$. It is also sensical to only parallelize computationally heavy tasks such as solving a system of equations.

Since independent sets of stages can be computed independently, it can prove useful to partition the stages into multiple independent sets of stages.

**II.2.55 Definition : Partition of independent sets of stages**

Let $(n, s) \in \mathbb{N}^* \times \mathbb{N}^*$, $M \in \mathrm{GMORK}_{n,s}$, and $P$ a partition of $[\![1, s]\!]$.
$P$ is said to be a partition of independent sets of stages if and only if, for all $(J_1, J_2) \in P^2$, if $J_1 \neq J_2$, then $J_1$ and $J_2$ are independent sets of stages.

**II.2.56 Proposition : Partitions of independent sets & Explicit partition of the stages**

Let $(n, s) \in \mathbb{N}^* \times \mathbb{N}^*$, $M \in \mathrm{GMORK}_{n,s}$, and $P$ a partition of independent sets of stages.
$P$ is an explicit partitions of the stages.

Proof

Since all the elements of $P$ are reachable closed sets, the dependency graph of $P$ only has empty diwalks.

**II.2.57 Theorem : Characterization of partitions of independent sets of stages**

Let $(n,s) \in \mathbb{N}^* \times \mathbb{N}^*$, $M \in \mathrm{GMORK}_{n,s}$, $B \in \mathcal{P}([\![1,s]\!])^{[\![1,p]\!]}$ a partition of $[\![1,s]\!]$, and $\varphi \in \mathfrak{S}_s^*$ such that, for all $k \in [\![1,p]\!],\ \varphi(B_k) \subset [\![1+\sum_{k'=1}^{k-1} \#(B_{k'}), \sum_{k'=1}^{k} \#(B_{k'})]\!]$
$B$ is an order of computation such that $\{B_1, ..., B_p\}$ is a partition of independent sets of stages if and only if $\varphi * M = (\tau', w', \tilde{w}') \in \mathrm{GMORK}_{n,s}$, where :

$$\forall N \in [\![1,n]\!],\ w'_N = \begin{bmatrix} A_{1,N} & 0_{\#(B_1),\#(B_2)} & \cdots & 0_{\#(B_1),\#(B_p)} \\ 0_{\#(B_2),\#(B_1)} & A_{2,N} & \ddots & \vdots \\ \vdots & \ddots & \ddots & 0_{\#(B_{p-1}),\#(B_p)} \\ 0_{\#(B_p),\#(B_1)} & \cdots & 0_{\#(B_p),\#(B_{p-1})} & A_{p,N} \\ * & \cdots & * & * \end{bmatrix}$$

And, for all $k \in [\![1,p]\!],\ A_k \in \mathbb{R}^{[\![1,n]\!] \times [\![1,\#(B_k)]\!]^2}$.

Proof

This is a direct consequence of the fact that the dependency graph of $P$ only has empty diwalks.

The reasoning we have followed so far applies in the general case, for special cases where we know that the differential equation function $f$ depends only on the ranks in a set $A \in \mathcal{P}([\![1,n]\!])$, we do not need to approximate intermediate values for the ranks that are not in $A$. We can thus remove the weights of the ranks that are not in $A$ for the stages $j \in [\![1,s]\!]$.

**II.2.58 Proposition : Initial value problems independent of some derivatives**

Let $(n,s) \in \mathbb{N}^* \times \mathbb{N}^*$, $M = (\tau, w, \tilde{w}) \in \mathrm{GMORK}_{n,s}$, $f$ a differential equation function such that $f$ is independent of the ranks in a set $A \in \mathcal{P}([\![1,n]\!])$, $\mathcal{U}_{n,f}$ the existence-uniqueness domain of $M$, and $Y_{n,f}$ the approximation function for length $n$ of $M$. Let $M' = (\tau, w', \tilde{w}') \in \mathrm{GMORK}_{n,s}$, where :

$$\forall N \in [\![1,n]\!] \setminus A,\ w'_N = w_N,\ \tilde{w}'_N = \tilde{w}_N$$
$$\forall N \in A,\ \forall N' \in [\![1,N]\!],\ \tilde{w}'_{N,N',s+1} = \tilde{w}_{N,N',s+1}$$
$$\forall N \in A,\ w'_{N,s+1} = w_{N,s+1}$$
$$\forall (j,N) \in [\![1,s]\!] \times A,\ \forall N' \in [\![1,N]\!],\ \tilde{w}'_{N,N',j} = 0$$
$$\forall (j,N) \in [\![1,s]\!] \times A,\ w'_{N,j} = 0_s$$

Let $\mathcal{U}'_{n,f}$ the existence-uniqueness domain of $M'$, and $Y'_{n,f}$ the approximation function for length $n$ of $M'$.
Then $\mathcal{U}_{n,f} = \mathcal{U}'_{n,f}$ and $Y_{n,f,s+1} = Y'_{n,f,s+1}$.

Proof

The ranks that are in $A$ are not used in $f$, there is hence no need to solve a system for them, and their value are not used by the other approximations. We conclude that we can replace with zeros the weights at the ranks that are in $A$ at the stages $j \in [\![1,s]\!]$ without changing the final approximations.

# III Convergence & Consistency

## III.1 Order of consistency

Now that we have rigorously defined how a general multi-order Runge-Kutta method approximates the solution of an initial value problem, we may be interested in how accurate its approximations are. To know how close the approximations are to the solution of an initial value problem, we first need to make sure the initial value problem has a unique solution.

**III.1.1 Theorem : Cauchy-Lipschitz theorem - Any order version**

Let $(f, t_0, y_{n,0})$ an initial value problem of order $n \in \mathbb{N}^*$ with $f$ continuous and locally Lipschitz continuous in its second variable.

- There exists a unique maximal solution $\hat{y} \in \{\Omega \to \mathbb{R}^d\}$ to the initial value problem.
- If $\mathbf{U} = \Omega \times \mathbb{R}^{[\![1,n]\!] \times [\![1,d]\!]}$ and $f$ is globally lipschitz continuous in its second variable, then there exists a unique solution to the initial value problem.

Proof

This theorem is implied by the Cauchy-Lipschitz theorem for equations of order 1 by rewriting the initial value problem to one of order 1.

Given an initial value problem with a unique solution, we wish to know the difference between the approximations and the actual solution. For Runge-Kutta methods we were only interested in the error of the last stage, but it will prove useful to also look into the error at the other stages. Since we cannot simply compute the error of the approximation, we fix $(t, y_{n,0}) \in \mathbf{U}$ and use Taylor's theorem to estimate how quickly the error tend towards 0 as $h$ tends towards 0.

Let $(n, s) \in \mathbb{N}^* \times \mathbb{N}^*$, $M = (\tau, w, \tilde{w}) \in \mathrm{GMORK}_{n,s}$, $f$ a differential equation function, $Y_{n,f}$ the approximation function for length $n$ of $M$, $(t, y_{n,0}, h) \in \mathbf{U} \times \mathbb{R}$, and $(j, N) \in [\![1, s+1]\!] \times [\![1, n]\!]$. Assume $(f, t, y_{n,0})$ has a unique solution $\hat{y}$. We wish to measure the difference between the approximation and the solution, but at which moment are we supposed to measure the difference ? Since multi-order methods are based on Taylor's theorem, it is safe to assume that it is at an instant of the form $t + \xi_{j,N} h$ with $h$ the step size and $\xi_{j,N} \in \mathbb{R}$ a constant. We are thus looking for a constant $\xi_{j,N}$ such that, as $h$ tends towards 0, $Y_{n,f,j,N}(t, y_{n,0}, h)$ approximates better and better the solution at $t + \xi_{j,N} h$ until at $h = 0$ it is exactly the solution. If we take $h = 0$, we have :

$$\frac{\mathrm{d}^{n-N} \hat{y}}{\mathrm{d}t^{n-N}}(t) - Y_{n,f,j,N}(t, y_{n,0}, 0) = y_{n,0,N} - \tilde{w}_{N,N,j} y_{n,0,N} = (1 - \tilde{w}_{N,N,j}) y_{n,0,N}$$

We conclude that if $\tilde{w}_{N,N,j} \neq 1$ then $Y_{n,j,N}$ does not approximate a moment. Assume $\tilde{w}_{N,N,j} = 1$ and $f$ is of class $C^1$. If $N \neq 1$, using the continuity of the approximations for small enough $h$, we have :

$$\begin{aligned}
\frac{\mathrm{d}^{n-N} \hat{y}}{\mathrm{d}t^{n-N}}(t + \xi_{j,N} h) - Y_{n,f,j,N}(t, y_{n,0}, h) &= \frac{\mathrm{d}^{n-N} \hat{y}}{\mathrm{d}t^{n-N}}(t) + \xi_{j,N} h \frac{\mathrm{d}^{n-N+1} \hat{y}}{\mathrm{d}t^{n-N+1}}(t) \\
&\quad - y_{n,0,N} - \tilde{w}_{N,N-1,j} h y_{n,0,N-1} + \underset{h \to 0}{o}(h) \\
&= (\xi_{j,N} - \tilde{w}_{N,N-1,j}) h y_{n,0,N-1} + \underset{h \to 0}{o}(h)
\end{aligned}$$

The difference is in $\underset{h\to 0}{o}(h)$ if and only if $\xi_{j,N} = \tilde{w}_{N,N-1,j}$. If $N = 1$, we have :

$$\frac{\mathrm{d}^{n-1}\hat{y}}{\mathrm{d}t^{n-1}}(t+\xi_{j,1}h) - Y_{n,f,j,1}(t, y_{n,0}, h) = \frac{\mathrm{d}^{n-1}\hat{y}}{\mathrm{d}t^{n-1}}(t) + \xi_{j,1}h\frac{\mathrm{d}^{n}\hat{y}}{\mathrm{d}t^{n}}(t)$$

$$-y_{n,0,1} - h\sum_{j'=1}^{s} w_{1,j,j'} f(t+\tau_j h, Y_{n,f,j}(t, y_{n,0}, h)) + \underset{h\to 0}{o}(h)$$

$$= \xi_{j,1}hf(t, y_{n,0}) - h \sum_{\substack{j'=1 \\ w_{1,j,j'}\neq 0}}^{s} w_{1,j,j'} f(t, \tilde{w}_{1,1,j'}y_{n,0,1}, ..., \tilde{w}_{n,n,j'}y_{n,0,n}) + \underset{h\to 0}{o}(h)$$

If, for all $j' \in [\![1, s]\!],\ w_{1,j,j'} \neq 0 \Rightarrow \forall N \in [\![1, n]\!],\ \tilde{w}_{N,N,j'} = 1$, then the method actually approximates an instant. If we assume so :

$$\frac{\mathrm{d}^{n-1}\hat{y}}{\mathrm{d}t^{n-1}}(t+\xi_{j,1}h) - Y_{n,f,j,1}(t, y_{n,0}, h) = \left(\xi_{j,1} - \sum_{j'=1}^{s} w_{1,j,j'}\right)hf(t_0, y_{n,0}) + \underset{h\to 0}{o}(h)$$

We conclude that $\xi_{j,1} = \sum_{j'=1}^{s} w_{1,j,j'}$.

### III.1.2 Definition : Approximation times

Let $(n, s) \in \mathbb{N}^* \times \mathbb{N}^*$, and $M = (\tau, w, \tilde{w}) \in \mathrm{GMORK}_{n,s}$.
We define $\xi \in \mathbb{R}^{[\![1,s+1]\!]\times[\![1,n]\!]}$ the approximation times of $M$ as :

$$\forall (j, N) \in [\![1, s+1]\!] \times [\![1, n]\!],\ \xi_{j,N} = \begin{cases} \sum_{j'=1}^{s} w_{1,j,j'} & \text{if } N = 1 \\ \tilde{w}_{N,N-1,j} & \text{otherwise} \end{cases}$$

### III.1.3 Proposition : Approximation times of node-determined methods

Let $(n, s) \in \mathbb{N}^* \times \mathbb{N}^*$, $M = (\tau, w) \in \mathrm{NDMORK}_{n,s}$, and $\xi$ the approximation times of $M$.
For all $(j, N) \in [\![1, s+1]\!] \times [\![2, n]\!],\ \xi_{j,N} = \tau_j$.

Proof

Trivial.

### III.1.4 Definition : Error function

Let $(n, s) \in \mathbb{N}^* \times \mathbb{N}^*$, $M \in \mathrm{GMORK}_{n,s}$, $Y_{n,f}$ the approximation function for length $n$ of $M$, $\mathcal{U}_{n,f}$ the existence-uniqueness domain of $M$, and $\xi$ the approximation times of $M$.
We define $\varepsilon_n$ the error function of $M$ for length $n$ as, for all $(f, t_0, y_{n,0})$ an initial value problem of order $n$ with a unique solution $\hat{y}$, for all $(j, N) \in [\![1, s+1]\!] \times [\![1, n]\!]$ :

$$\varepsilon_{n,f,j,N} : \{(t', y'_{n,0}, h) \in \mathcal{U}_{n,f} \mid t' + \xi_{j,N}h \in \Omega\} \to \mathbb{R}^{[\![1,d]\!]}$$

$$(t', y'_{n,0}, h) \to \frac{\mathrm{d}^{n-N}\hat{y}}{\mathrm{d}t^{n-N}}(t' + \xi_{j,N}h) - Y_{n,f,j,N}(t', y'_{n,0}, h)$$

We already have a theorem on the smoothness of the approximations, but we are missing a theorem on the smoothness of the solutions of an initial value problem, which is necessary to use Taylor's theorem.

### III.1.5 Proposition : Class of the solutions of an initial value problem

Let $(n, k) \in \mathbb{N}^* \times \mathbb{N}$, and $(f, t_0, y_{n,0})$ an initial value problem of order $n$.
If $f$ is of class $C^k$ then all the solutions of the initial value problem are of class $C^{n+k}$.

Proof

Let $\hat{y}$ a solution. We prove this by induction on $k$. If $\hat{y}$ is of class $C^{n+k}$ and $f$ is of class $C^{k+1}$, $\frac{d^n \hat{y}}{dt^n}$ equals $f \circ \overline{\mathcal{J}}^{n-1}\hat{y}$, which is of class $C^{k+1}$, thus $\frac{d^n \hat{y}}{dt^n}$ is of class $C^{k+1}$ and $\hat{y}$ is of class $C^{n+k+1}$.

The smoothness of the approximations and the solutions are implied by that of $f$, we can hence simply assume that $f$ is smooth enough to use Taylor's theorem. Let $(f, t_0, y_{n,0})$ an initial value problem of order $n \in \mathbb{N}^*$ with $f$ of class $C^0$ and a unique solution $\hat{y}$, $M = (\tau, w, \tilde{w}) \in \text{NDMORK}_{n,s}$, $h \in \mathcal{H}_{1,f,t_0,y_{n,0}}$, and $(j, N) \in [\![1, s+1]\!] \times [\![1, n]\!]$. $\varepsilon_n$, the error for length $n$ of $M$, satisfies thanks to Taylor's theorem :

$$\begin{aligned}
&\varepsilon_{n,f,j,N}(t_0, y_{n,0}, h) \\
&= \frac{h^N}{N!}\left(\frac{N}{2^N}\int_{-1}^{1}(1-x)^{N-1} f\left(\overline{\mathcal{J}}^{n-1}\hat{y}\left(t_0 + \frac{1+x}{2}h\right)\right) \mathrm{d}x - \sum_{j'=1}^{s} w_{N,s+1,j'} F_{n,f,j}(t_0, y_{n,0}, h)\right)
\end{aligned}$$

We deduce that the error of the approximation is the error of the term in parentheses multiplied by $h^N$.

**III.1.6 Definition : Minimum consistency matrix**

Let $(n, s) \in \mathbb{N}^* \times \mathbb{N}^*$, $M \in \text{GMORK}_{n,s}$, $\varepsilon_n$ the error function for length $n$ of $M$, and $\xi$ its approximation times.

- Let $(j, N) \in [\![1, s+1]\!] \times [\![1, n]\!]$, and $v' \in \mathbb{Z}$.
  $M$ is said to be at least of order of consistency $v'$ for length $n$ at stage $j$ and rank $N$ if and only if, for all initial value problems $(f, t_0, y_{n,0})$ of order $n$ with a unique solution and $f$ of class $C^{\max(1,v')}$, for all $\Omega' \in \mathcal{P}(\Omega)$ a compact subset of $\Omega$ :

$$\begin{aligned}
\left\{(t', h) \in \Omega' \times \mathbb{R} \mid h \in \mathcal{H}_{1,f,\overline{\mathcal{J}}^{n-1}\hat{y}(t')}, t' + \xi_{j,N} h \in \Omega\right\} &\to \mathbb{R}^{[\![1,d]\!]} \\
(t', h) \to \varepsilon_{n,f,j,N}\left(\overline{\mathcal{J}}^{n-1}\hat{y}(t'), h\right) &= \underset{\substack{h \to 0 \\ t' \in \Omega'}}{\mathcal{O}}\left(h^{v'+N}\right)
\end{aligned}$$

- We define $v_n \in \overline{\mathbb{Z}}^{[\![1,s+1]\!] \times [\![1,n]\!]}$ the minimum consistency matrix of $M$ as, for all $(j, N) \in [\![1, s+1]\!] \times [\![1, n]\!]$, $v_{n,j,N}$ is the supremum of the set of $v' \in \mathbb{Z}$ such that $M$ is at least of order of consistency $v'$ for length $n$ at stage $j$ and rank $N$. $M$ is said to be of order of consistency $v_{n,j,N}$ for length $n$ at stage $j$ and rank $N$.

The adjective for length $n$ does not imply a condition on the accuracy of a method, it simply indicates that this is the accuracy of the method applied to initial value problems of order $n$. For now, the methods only deal with a single order of initial value problem, but later $n$ will vary.

We can consider that the minimum requirement for a good general multi-order Runge-Kutta method is that all of the entries of its minimum consistency matrix are non-negative.

**III.1.7 Definition : Multi-order Runge-Kutta methods**

- A general multi-order Runge-Kutta method is said to be a multi-order Runge-Kutta method if and only if all of the entries of its minimum consistency matrix are superior or equal to 0.
- We define $\text{MORK}_{n,s}$ as the set of multi-order Runge-Kutta methods of length $n \in \mathbb{N}^*$ with $s \in \mathbb{N}^*$ points.

Thankfully, the conditions for a non-negative order of consistency are simple since we only consider the Taylor polynomials up to the $n-1^{\text{th}}$ derivative of $\hat{y}$, which do not involve $f \circ \overline{\mathcal{J}}^{n-1}\hat{y}$.

### III.1.8 Proposition : Negative order of consistency

Let $(n,s) \in \mathbb{N}^* \times \mathbb{N}^*$, $M = (\tau, w, \tilde{w}) \in \mathrm{GMORK}_{n,s}$, $v_n$ the minimum consistency matrix of $M$, and $(j, N) \in [\![1, s+1]\!] \times [\![1, n]\!]$.

1. $v_{n,j,N} \geq -N$.
2. Let $N' \in [\![1, N]\!]$. If $\tilde{w}_{N,N',j} \neq \frac{\xi_{j,N}^{N-N'}}{(N-N')!}$, then $v_{n,j,N} \leq -N'$.
3. Let $N' \in [\![1, N]\!]$. $v_{n,j,N} > -N'$ if and only if, for all $N'' \in [\![N', N]\!]$, $\tilde{w}_{N,N'',j} = \frac{\xi_{j,N}^{N-N''}}{(N-N'')!}$.

Proof

1. Since $\varepsilon_{n,f,j,N}\left(\mathcal{J}^{n-1}\hat{y}(t'), h\right) \underset{\substack{h\to 0\\ t'\in\Omega'}}{=} \mathcal{O}(1)$, $M$ is at least of order of consistency $-N$ at rank $N$.
2. It is trivial to prove if we consider the initial value problem $\left(f, 0, y_{n,0}\right)$ where $f$ is identically 0, and :

$$y_{n,0,N} = 1,\ \forall N' \in [\![1, n]\!] \setminus \{N\},\ y_{n,0,N'} = 0$$

3. The necessity of this condition is implied by 2. Taylor's theorem gives that this condition is sufficient.

This proposition allows us to characterize explicitly multi-order Runge-Kutta methods.

### III.1.9 Proposition : Characterization of multi-order Runge-Kutta methods

Let $(n,s) \in \mathbb{N}^* \times \mathbb{N}^*$, $M = (\tau, w, \tilde{w}) \in \mathrm{GMORK}_{n,s}$, and $\xi$ the approximation times of $M$.
$M$ is a multi-order Runge-Kutta method if and only if :

$$\forall (j, N) \in [\![1, s+1]\!] \times [\![1, n]\!],\ \forall N' \in [\![1, N]\!],\ \tilde{w}_{N,N',j} = \frac{\xi_{j,N}^{N-N'}}{(N-N')!}$$

Multi-order Runge-Kutta methods can thus be identified to a tuple $(\tau, w, \xi)$, with $\tau \in \mathbb{R}^{[\![1,s]\!]}$ its nodes, $w \in \mathbb{R}^{[\![1,n]\!]\times[\![1,s+1]\!]\times[\![1,s]\!]}$ its main weights, $\xi \in \mathbb{R}^{[\![1,s+1]\!]\times[\![1,n]\!]}$ its approximation times, which satisfy $\xi_{j,1} = \sum_{j'=1}^{s} w_{1,j,j'}$.

Proof

It is a direct consequence of the proposition on negative order of Section III.1.8.

### III.1.10 Proposition : Use of multi-order Runge-Kutta methods

Let $(n,s) \in \mathbb{N}^* \times \mathbb{N}^*$, $M = (\tau, w, \xi) \in \mathrm{MORK}_{n,s}$, $f$ a differential equation function of order $n$, and $\left(t, y_{n,0}, h\right) \in \mathbf{U} \times \mathbb{R}$.
The stage system $\left(M, f, t, y_{n,0}, h\right)$ is equivalent to :

$$\forall (j, N) \in [\![1, s+1]\!] \times [\![1, n]\!],$$

$$y_{n,j,N} = \sum_{N'=1}^{N} \frac{\left(\xi_{j,N} h\right)^{N-N'}}{(N-N')!} y_{n,0,N'} + \frac{h^N}{N!} \sum_{j'=1}^{s} w_{N,j,j'} f\left(t + \tau_{j'} h, y_{n,j'}\right)$$

$$y_{n,s+1} \in \mathbb{R}^{[\![1,n]\!]\times[\![1,d]\!]},\ \forall j \in [\![1, s]\!],\ \left(t + \tau_j h, y_{n,j}\right) \in \mathbf{U}$$

Proof

Trivial.

This category is the most general type of method whose use of initial value matches exactly that of a Taylor polynomial, while letting us choose the moment we wish to approximate. Node-determined methods are more restrictive.

**III.1.11 Proposition : Link between MORK and NDMORK**

Let $(n,s) \in \mathbb{N}^* \times \mathbb{N}^*$, $M = (\tau, w, \tilde{w}) \in \mathrm{GMORK}_{n,s}$, and $\xi$ the approximation times of $M$.
$M$ is a node determined method if and only if $M$ is a multi-order Runge-Kutta method and there exists $\tau_{s+1} \in \mathbb{R}$ such that :

$$\forall (j,N) \in [\![1, s+1]\!] \times [\![2, n]\!],\ \xi_{j,N} = \tau_j$$

For positive orders of consistency the conditions are more complex, we are only able to derive an equivalent system of equations.

**III.1.12 Definition : Equation for order of consistency**

Let $(n,s) \in \mathbb{N}^* \times \mathbb{N}^*$, $M = (\tau, w, \tilde{w}) \in \mathrm{GMORK}_{n,s}$, $Y_n$ the approximation function for length $n$ of $M$, $(j,N) \in [\![1, s+1]\!] \times [\![1, n]\!]$, and $(v', k) \in \mathbb{N}^* \times \mathbb{N}^*$ such that $v' \leq k$.
$M$ is said to satisfy the equation for order of consistency $v'$ of class $k$ for length $n$ at stage $j$ and rank $N$ if and only if, for all $\left(t, y_{n,0}\right) \in \mathbb{R} \times \mathbb{R}^{[\![1,n]\!] \times [\![1,d]\!]}$, for all $\mathbf{U} \in \mathcal{P}\left(\mathbb{R} \times \mathbb{R}^{[\![1,n]\!] \times [\![1,d]\!]}\right)$ a neighborhood of $\left(t, y_{n,0}\right)$, for all $f \in C^k\left(\mathbf{U}, \mathbb{R}^{[\![1,d]\!]}\right)$, with $\hat{\mathrm{y}}$ the unique maximal solution of $\left(f, t, y_{n,0}\right)$ defined on $\Omega$ :

$$\left\{h \in \mathcal{H}_{1,f,t,y_{n,0}} \mid t + \xi_{j,N} h \in \Omega\right\} \to \mathbb{R}^{[\![1,d]\!]}$$

$$h \to \sum_{j'=1}^{s} w_{N,j,j'} \sum_{k=0}^{v'-1} \frac{1}{k!} \mathrm{d}^k f\left(t, y_{n,0}\right)\left(\tau_{j'} h, Y_{n,f,j'}\left(t, y_{n,0}, h\right) - y_{n,0}\right)$$

$$- \sum_{k=0}^{v'-1} \xi_{j,N}^{N+k} h^k \frac{N!}{(k+N)!} \frac{\mathrm{d}^k f \circ \overline{\mathcal{J}}^{n-1} \hat{\mathrm{y}}}{\mathrm{d}t^k}(t) = \underset{h \to 0}{\mathcal{O}}\left(h^v\right)$$

**III.1.13 Theorem : Order of consistency & Equation for order order of consistency**

Let $(n,s) \in \mathbb{N}^* \times \mathbb{N}^*$, $M = (\tau, w, \tilde{w}) \in \mathrm{GMORK}_{n,s}$, $(j,N) \in [\![1, s+1]\!] \times [\![1, n]\!]$, and $v' \in \mathbb{N}^*$.
If, for all $(j', N') \in [\![1, s]\!] \times [\![1, n]\!]$, $\tilde{w}_{N',N',j'} = 1$, $M$ is at least of order of consistency 0 for length $n$ at stage $j$ and rank $N$, then, $M$ is at least of order of consistency $v'$ for length $n$ at stage $j$, rank $N$ if and only if $M$ satisfies the equation for order of consistency $v'$ of class $v'$ for length $n$, stage $j$, rank $N$.

Proof

Assume $M$ is at least of order of consistency $v'$ for length $n$ at stage $j$ and rank $N$. Let $\left(f, t_0, y_{n,0}\right)$ an initial value problem of order $n$ with a unique solution $\hat{\mathrm{y}}$ and $f$ of class $C^{v'}$, $\Omega' \in \mathcal{P}(\Omega)$ a compact subset of $\Omega$, and $t' \in \Omega'$. We expand the expression of $\frac{\mathrm{d}^{n-N} \hat{\mathrm{y}}}{\mathrm{d}t^{n-N}}$ :

$$\frac{\mathrm{d}^{n-N} \hat{\mathrm{y}}}{\mathrm{d}t^{n-N}}\left(t' + \xi_{j,N} h\right) = \sum_{N'=0}^{v'+N-1} \frac{\left(\xi_{j,N} h\right)^{N'}}{N'!} \frac{\mathrm{d}^{N'+n-N} \hat{\mathrm{y}}}{\mathrm{d}t^{N'+n-N}}(t') + \underset{\substack{h \to 0 \\ t' \in \Omega'}}{\mathcal{O}}\left(h^{v'+N}\right)$$

$$= \sum_{N'=0}^{N-1} \frac{\left(\xi_{j,N} h\right)^{N'}}{N'!} \frac{\mathrm{d}^{N'+n-N} \hat{\mathrm{y}}}{\mathrm{d}t^{N'+n-N}}(t') + \sum_{N'=N}^{v'+N-1} \frac{\left(\xi_{j,N} h\right)^{N'}}{N'!} \frac{\mathrm{d}^{n+N'-N} \hat{\mathrm{y}}}{\mathrm{d}t^{n+N'-N}}(t') + \underset{\substack{h \to 0 \\ t' \in \Omega'}}{\mathcal{O}}\left(h^{v'+N}\right)$$

$$= \sum_{N'=1}^{N} \frac{\left(\xi_{j,N} h\right)^{N-N'}}{(N-N')!} \frac{\mathrm{d}^{n-N'} \hat{\mathrm{y}}}{\mathrm{d}t^{n-N'}}(t') + \sum_{N'=N}^{v'+N-1} \frac{\left(\xi_{j,N} h\right)^{N'}}{N'!} \frac{\mathrm{d}^{N'-N} f \circ \overline{\mathcal{J}}^{n-1} \hat{\mathrm{y}}}{\mathrm{d}t^{N'-N}}(t') + \underset{\substack{h \to 0 \\ t' \in \Omega'}}{\mathcal{O}}\left(h^{v'+N}\right)$$

$$= \sum_{N'=1}^{N} \frac{\left(\xi_{j,N} h\right)^{N-N'}}{(N-N')!} \frac{\mathrm{d}^{n-N'} \hat{\mathrm{y}}}{\mathrm{d}t^{n-N'}}(t') + \sum_{k=0}^{v'-1} \frac{\left(\xi_{j,N} h\right)^{k+N}}{(k+N)!} \frac{\mathrm{d}^k f \circ \overline{\mathcal{J}}^{n-1} \hat{\mathrm{y}}}{\mathrm{d}t^k}(t') + \underset{\substack{h \to 0 \\ t' \in \Omega'}}{\mathcal{O}}\left(h^{v'+N}\right)$$

We also need to develop the evaluations of $f$. For all $j \in [\![1, s]\!]$ :

$$
\begin{aligned}
&F_{n,f,j}\Big(\overline{\mathcal{J}}^{n-1}\hat{\mathrm{y}}(t'), h\Big) = f\Big(t' + \tau_j h, \mathcal{J}^{n-1}\hat{\mathrm{y}}(t') + Y_{n,f,j}\Big(\overline{\mathcal{J}}^{n-1}\hat{\mathrm{y}}(t'), h\Big) - \mathcal{J}^{n-1}\hat{\mathrm{y}}(t')\Big) \\
&= \sum_{k=0}^{v'-1} \frac{1}{k!}\mathrm{d}^k f\Big(\mathcal{J}^{n-1}\hat{\mathrm{y}}(t')\Big)\Big(\tau_j h, Y_{n,f,j}\Big(\overline{\mathcal{J}}^{n-1}\hat{\mathrm{y}}(t'), h\Big) - \mathcal{J}^{n-1}\hat{\mathrm{y}}(t')\Big) \\
&\quad + \underset{\substack{\|\left(\tau_j h, Y_{n,f,j}\left(\overline{\mathcal{J}}^{n-1}\hat{\mathrm{y}}(t'), h\right) - \mathcal{J}^{n-1}\hat{\mathrm{y}}(t')\right)\| \to 0 \\ t' \in \Omega'}}{\mathcal{O}}\left(\Big\|\Big(\tau_j h, Y_{n,f,j}\Big(\overline{\mathcal{J}}^{n-1}\hat{\mathrm{y}}(t'), h\Big) - \mathcal{J}^{n-1}\hat{\mathrm{y}}(t')\Big)\Big\|^{v'}\right)
\end{aligned}
$$

For all $(j', N') \in [\![1, s]\!] \times [\![1, n]\!]$, $\tilde{w}_{N',N',j'} = 1$, thus $Y_{n,f,j'}\left(\overline{\mathcal{J}}^{n-1}\hat{\mathrm{y}}(t'), h\right) - \mathcal{J}^{n-1}\hat{\mathrm{y}}(t') = \underset{\substack{h \to 0 \\ t' \in \Omega'}}{\mathcal{O}}(h)$.

Therefore :

$$
\underset{\substack{\|\left(\tau_j h, Y_{n,f,j}\left(\overline{\mathcal{J}}^{n-1}\hat{\mathrm{y}}(t'), h\right) - \mathcal{J}^{n-1}\hat{\mathrm{y}}(t')\right)\| \to 0 \\ t' \in \Omega'}}{\mathcal{O}}\left(\Big\|\Big(\tau_j h, Y_{n,f,j}\Big(\overline{\mathcal{J}}^{n-1}\hat{\mathrm{y}}(t'), h\Big) - \mathcal{J}^{n-1}\hat{\mathrm{y}}(t')\Big)\Big\|^{v'}\right) = \underset{\substack{h \to 0 \\ t' \in \Omega'}}{\mathcal{O}}\left(h^{v'}\right)
$$

$M$ is at least of order of consistency 0 for length $n$ at stage $j$ and rank $N$, thus, for all $N' \in [\![1, N]\!]$, $\tilde{w}_{N,N',j} = \frac{\xi_{j,N}^{N-N'}}{(N-N')!}$, therefore we finally get :

$$
\begin{aligned}
&Y_{n,f,j,N}\Big(\overline{\mathcal{J}}^{n-1}\hat{\mathrm{y}}(t'), h\Big) = \sum_{N'=1}^{N} \frac{\left(\xi_{j,N} h\right)^{N-N'}}{(N-N')!} \frac{\mathrm{d}^{n-N'}\hat{\mathrm{y}}}{\mathrm{d}t^{n-N'}}(t') \\
&+ \frac{h^N}{N!} \sum_{j'=1}^{s} w_{N,j,j'} \sum_{k=0}^{v'-1} \frac{1}{k!}\mathrm{d}^k f\Big(\mathcal{J}^{n-1}\hat{\mathrm{y}}(t')\Big)\Big(\tau_{j'} h, Y_{n,f,j'}\left(\overline{\mathcal{J}}^{n-1}\hat{\mathrm{y}}(t'), h\right) - \mathcal{J}^{n-1}\hat{\mathrm{y}}(t')\Big) + \underset{\substack{h \to 0 \\ t' \in \Omega'}}{\mathcal{O}}\left(h^{v'+N}\right)
\end{aligned}
$$

The initial values in $\varepsilon_{n,f,j,N}\Big(\overline{\mathcal{J}}^{n-1}\hat{\mathrm{y}}(t'), h\Big)$ cancel out, we hence have :

$$
\begin{aligned}
& h \to \varepsilon_{n,f,j,N}\Big(\overline{\mathcal{J}}^{n-1}\hat{\mathrm{y}}(t'), h\Big) = \underset{\substack{h \to 0 \\ t' \in \Omega'}}{\mathcal{O}}\left(h^{v'+N}\right) \\
\Leftrightarrow\ & h \to \frac{h^N}{N!} \sum_{j'=1}^{s} w_{N,j,j'} \sum_{k=0}^{v'-1} \frac{1}{k!}\mathrm{d}^k f\Big(\mathcal{J}^{n-1}\hat{\mathrm{y}}(t')\Big)\Big(\tau_{j'} h, Y_{n,f,j'}\left(\overline{\mathcal{J}}^{n-1}\hat{\mathrm{y}}(t'), h\right) - \mathcal{J}^{n-1}\hat{\mathrm{y}}(t')\Big) \\
& \quad - h^N \sum_{k=0}^{v'-1} \xi_{j,N}^{N+k} \frac{h^k}{(k+N)!} \frac{\mathrm{d}^k f \circ \overline{\mathcal{J}}^{n-1}\hat{\mathrm{y}}}{\mathrm{d}t^k}(t') = h^N \underset{\substack{h \to 0 \\ t' \in \Omega'}}{\mathcal{O}}\left(h^{v'}\right) \\
\Leftrightarrow\ & h \to \sum_{j'=1}^{s} w_{N,j,j'} \sum_{k=0}^{v'-1} \frac{1}{k!}\mathrm{d}^k f\Big(\overline{\mathcal{J}}^{n-1}\hat{\mathrm{y}}(t')\Big)\Big(\tau_{j'} h, Y_{n,f,j'}\left(\overline{\mathcal{J}}^{n-1}\hat{\mathrm{y}}(t'), h\right) - \mathcal{J}^{n-1}\hat{\mathrm{y}}(t')\Big) \\
& \quad - \sum_{k=0}^{v'-1} \xi_{j,N}^{N+k} h^k \frac{N!}{(k+N)!} \frac{\mathrm{d}^k f \circ \overline{\mathcal{J}}^{n-1}\hat{\mathrm{y}}}{\mathrm{d}t^k}(t') = \underset{\substack{h \to 0 \\ t' \in \Omega'}}{\mathcal{O}}\left(h^{v'}\right)
\end{aligned}
$$

We now prove this condition is equivalent to, $M$ satisfies the equation for order of consistency $v'$ of class $v'$ for length $n$ at stage $j$ and rank $N$.

Assume the previous condition is true. We consider a case of the equation for order of consistency. Since $\hat{\mathrm{y}}$ is a unique maximal solution it is defined on a an open connected subset of $\Omega'$ which contains $t_0$. If we restrict the initial value problem to this subset, we can see the equation is satisfied by taking $t' = t_0$.

Assume the equation for order of consistency is satisfied. Since $f, t, y_{n,0}, h$ all vary, this implies that the coefficients of the Taylor polynomial cancel each others, which implies the previous condition.

Since for high orders of consistency we assume a certain smoothness of $f$, we need to make sure the approximations stay as accurate as possible when $f$ is less smooth.

**III.1.14 Proposition : Accuracy for less smooth initial value problems**

Let $(n, s) \in \mathbb{N}^* \times \mathbb{N}^*$, $M \in \mathrm{GMORK}_{n,s}$, $(j, N) \in [\![1, s+1]\!] \times [\![1, n]\!]$, $(v'_1, v'_2) \in \mathbb{Z}^{*2}$, and $(k_1, k_2) \in \overline{\mathbb{Z}}^{*2}$ such that $v'_2 \leq v'_1$, $v'_1 \leq k_1$, and $v'_2 \leq k_2$.
If $M$ is at least of order of consistency 0 for length $n$ at stage $j$ and rank $N$ and, for all $(j', N') \in [\![1, s]\!] \times [\![1, n]\!]$, $\tilde{w}_{N',N',j'} = 1$, then, if $M$ satisfies the equation for order of consistency $v'_1$ of class $k_1$ for length $n$ at stage $j$ and rank $N$, then $M$ satisfies the equation for order of consistency $v'_2$ of class $k_2$ for length $n$ at stage $j$ and rank $N$.

Proof

$M$ satisfies the equation if and only if the coefficients of the Taylor polynomials in the error cancel each other. If they do cancel each other up to degree $v'_1$, then it also does up to degree $v'_2$. The equations only need functions of high enough class to use Taylor's theorem, the actual class does not matter.

The number of terms in the equation for order of consistency explodes rapidly. One way to downplay this issue would be to only consider autonomous initial value problems, problems where the differential equation function is independent of $t$. There is a way to do so at the price of a condition on the weights at rank 1 of the method. We can rewrite any initial value problem into an equivalent autonomous initial value problem by keeping track of $t$ using $\hat{y}$ or one of its derivatives.

**III.1.15 Definition: Rewriting to autonomous initial value problems**

Let $(f, t_0, y_{n,0})$ an initial value problem of order $n \in \mathbb{N}^*$, and $N \in [\![1, n]\!]$. We define :

$$\mathbf{U}' = \left\{ \left( t, \begin{bmatrix} x_1 \\ t'_1 \end{bmatrix}, \dots, \begin{bmatrix} x_n \\ t'_n \end{bmatrix} \right) \in \mathbb{R} \times \mathbb{R}^{[\![1,n]\!] \times [\![1,d+1]\!]} \mid (t'_N, x_1,, \dots, x_n) \in \mathbf{U} \right\}$$

$$\forall \left( t, \begin{bmatrix} x_1 \\ t'_1 \end{bmatrix}, \dots, \begin{bmatrix} x_n \\ t'_n \end{bmatrix} \right) \in \mathbf{U}',\ f'\left( t, \begin{bmatrix} x_1 \\ t'_1 \end{bmatrix}, \dots, \begin{bmatrix} x_n \\ t'_n \end{bmatrix} \right) = \begin{bmatrix} f(t'_N, x_1, \dots, x_n) \\ \begin{cases} 1 \text{ if } N = 1 \\ 0 \text{ otherwise} \end{cases} \end{bmatrix}$$

$$\forall N' \in [\![1, n]\!],\ y'_{n,0,N'} = \begin{bmatrix} y_{n,0,N'} \\ \begin{cases} 0 & \text{if } N' \in [\![1, N-2]\!] \\ \frac{t_0^{N'-N+1}}{(N'-N+1)!} & \text{otherwise} \end{cases} \end{bmatrix}$$

We define $(f', t_0, y'_{n,0})$ as the rank $N$ autonomization of $(f, t_0, y_{n,0})$.

**III.1.16 Proposition : Solution of rewritten initial value problems - Autonomous rewriting**

Let $(f, t_0, y_{n,0})$ an initial value problem of order $n \in \mathbb{N}^*$, $N \in [\![1, n]\!]$, and $(f', t_0, y'_{n,0})$ the rank $N$ autonomization of $(f, t_0, y_{n,0})$.

- If $\hat{y}$ is a solution of $(f, t_0, y_{n,0})$ then $\begin{bmatrix} \hat{y} \\ \frac{t^{n-N+1}}{(n-N+1)!} \end{bmatrix}$ is a solution of $(f', t_0, y'_{n,0})$.
- If $\hat{y}'$ is a solution of $(f', t_0, y'_{n,0})$, then there exists a unique $\hat{y}$ such that $\hat{y}' = \begin{bmatrix} \hat{y} \\ \frac{t^{n-N+1}}{(n-N+1)!} \end{bmatrix}$. $\hat{y}$ is a solution of $(f, t_0, y_{n,0})$.

Proof

Trivial.

### III.1.17 Proposition : Form of rewritten method - Autonomous rewriting

Let $(n, s) \in \mathbb{N}^* \times \mathbb{N}^*$, $M = (\tau, w, \tilde{w}) \in \text{GMORK}_{n,s}$, $f$ a differential equation function of order $n$, $(t, y_{n,0}, h) \in \mathbf{U} \times \mathbb{R}$, $N \in [\![1, n]\!]$, and $(f', t, y'_{n,0})$ the rank $N$ autonomization of $(f, t, y_{n,0})$.
The stage system $(M, f', t, y'_{n,0}, h)$ is equivalent to :

$$\forall (j, N') \in [\![1, s+1]\!] \times [\![1, n]\!],$$

$$y_{n,j,N'} = \sum_{N''=1}^{N'} \tilde{w}_{N',N'',j} h^{N'-N''} y_{n,0,N''} + \frac{h^{N'}}{N'!} \sum_{j'=1}^{s} w_{N',j,j'} f\big(\tilde{w}_{N,N,j} t + \xi_{j,N} h, y_{n,j'}\big)$$

$$y_{n,s+1} \in \mathbb{R}^{[\![1,n]\!] \times [\![1,d]\!]},\ \forall j \in [\![1, s]\!],\ \big(\tilde{w}_{N,N,j} t + \xi_{j,N} h, y_{n,j}\big) \in \mathbf{U}$$

Proof

The stage system is the same as usual except $y'_{n,j,N,d+1}$ replaces $t + \tau_j h$. If $N > 1$, then $y'_{n,0,N,d+1} = t$, $y'_{n,0,N-1,d+1} = 1$ and all initial values with a lower rank are 0, hence, for all $j \in [\![1, s+1]\!]$:

$$\begin{aligned} y'_{n,j,N,d+1} &= \sum_{N'=1}^{N} \tilde{w}_{N,N',j} h^{N-N'} y'_{n,0,N',d+1} + \frac{h^N}{N!} \sum_{j'=1}^{s} w_{N,j,j'} f'_{d+1}\big(y'_{n,j'}\big) \\ &= \tilde{w}_{N,N,j} t + \tilde{w}_{N,N-1,j} h \\ &= \tilde{w}_{N,N,j} t + \xi_{j,N} h \end{aligned}$$

If $N = 1$, since $f'_{d+1} = 1$ and $y'_{n,0,1,d+1} = t$, we have, for all $j \in [\![1, s+1]\!]$ :

$$\begin{aligned} y'_{n,j,1,d+1} &= \tilde{w}_{1,1,j} y'_{n,0,1,d+1} + h \sum_{j'=1}^{s} w_{1,j,j'} f'_{d+1}\big(y'_{n,j'}\big) \\ &= \tilde{w}_{1,1,j} t + h \sum_{j'=1}^{s} w_{1,j,j'} \\ &= \tilde{w}_{1,1,j} t + \xi_{j,1} h \end{aligned}$$

We could simply accept this rewriting and apply methods only on autonomized initial value problems, but it would imply downsides similar to that of rewriting an initial value problem of order $n$ to one of order 1. We need to make sure the approximations stay the same.

### III.1.18 Proposition : Autonomous system condition

Let $(n, s) \in \mathbb{N}^* \times \mathbb{N}^*$, $M = (\tau, w, \tilde{w}) \in \text{GMORK}_{n,s}$, $f$ a differential equation function of order $n$, $(t, y_{n,0}, h) \in \mathbf{U} \times \mathbb{R}$, $N \in [\![1, n]\!]$, and $(f', t, y'_{n,0})$ the rank 1 autonomization of $(f, t, y_{n,0})$.
If, for all $j \in [\![1, s]\!]$ :

$$\tilde{w}_{1,1,j} = 1,\ \sum_{j'=1}^{s} w_{1,j,j'} = \tau_j$$

Then the stage system $(M, f', t, y'_{n,0}, h)$ is equivalent to the stage system $(M, f, t, y_{n,0}, h)$.

Proof

This is directly implied by the previous proposition.

We only consider the rank 1 autonomization because it is the only one appliable to all methods, no matter their length. This choice is necessary because in the next subsection we consider initial value problems of order inferior to the length of the method and the truncations of a method.

The autonomous condition is not necessary since, for example, if there exists $j' \in [\![1, s]\!]$ such that the evaluation of this stage is not used, thus if, for all $(j, N) \in [\![1, s+1]\!] \times [\![1, n]\!],\ w_{N,j,j'} = 0$, then we do not need this condition at stage $j'$.

This condition is well known in the current literature, and most authors assume *de facto* that the methods they consider satisfy it. This condition is sensical since it is the condition for an order of consistency 1 for length $n$ at first rank and all stages. Still, it is a loss of generality.

It is also possible to reduce the number of derivatives of $\hat{y}$ on which $f$ depends by using the rewriting of the case $n = \tilde{n}$ covered in the next subsection, though it is of limited interest and we cannot remove the dependence of $f$ on $\frac{\mathrm{d}^{n-1}\hat{y}}{\mathrm{d}t^{n-1}}$, which would be the most useful.

Since the accuracy required for a method of a certain order of consistency must hold for any initial value problem $(f, t_0, y_{n,0})$, it implies that the order of consistency must also hold if $f$ is a function of only $t$. In this case the initial value problem is simply an integration problem, the initial value problem is in some sense solved, and finding the conditions for a given order of consistency $v' \in \mathbb{N}^*$ is trivial.

**III.1.19 Theorem: Solved system conditions**

Let $(n, s) \in \mathbb{N}^* \times \mathbb{N}^*$, $M = (\tau, w, \tilde{w}) \in \mathrm{GMORK}_{n,s}$, $v' \in \mathbb{N}^*$, and $(j, N) \in [\![1, s+1]\!] \times [\![1, n]\!]$.
If $M$ is at least of order of consistency $v'$ for length $n$ at stage $j$ and rank $N$, then :

$$\forall k \in [\![0, v'-1]\!],\ \sum_{j'=1}^{s} w_{N,j,j'} \tau_{j'}^{k} = \xi_{j,N}^{N+k} \frac{k!N!}{(k+N)!}$$

Those conditions are called the solved system conditions of order of consistency $v'$ for length $n$ at stage $j$ and rank $N$.

Proof

We consider the case where $f$ depends only on $t$. $f$ is independent of the approximations, we hence do not need the condition, for all $(j, N) \in [\![1, s]\!] \times [\![1, n]\!],\ \tilde{w}_{N,N,j} = 1$, to use the equation for order of consistency $v'$, which yields :

$$\sum_{j'=1}^{s} w_{N,j,j'} \sum_{k=0}^{v'-1} \frac{1}{k!} \mathrm{d}^k f - \sum_{k=0}^{v'-1} \xi_{j,N}^{N+k} h^k \frac{N!}{(k+N)!} \frac{\mathrm{d}^k f \circ \mathcal{J}^{n-1}\hat{y}}{\mathrm{d}t^k} = \underset{h\to 0}{\mathcal{O}}\left(h^{v'}\right)$$

$$\Leftrightarrow \sum_{j'=1}^{s} w_{N,j,j'} \sum_{k=0}^{v'-1} \frac{1}{k!} \left( \frac{\partial^k f}{\partial t^k} \left(\tau_{j'} h\right)^k \right) - \sum_{k=0}^{v'-1} \xi_{j,N}^{N+k} h^k \frac{N!}{(k+N)!} \frac{\partial^k f}{\partial t^k} = \underset{h\to 0}{\mathcal{O}}\left(h^{v'}\right)$$

$$\Leftrightarrow \sum_{k=0}^{v'-1} h^k \frac{\partial^k f}{\partial t^k} \left( \left( \sum_{j'=1}^{s} w_{N,j,j'} \tau_{j'}^{k} \right) \frac{1}{k!} - \xi_{j,N}^{N+k} \frac{N!}{(k+N)!} \right) = \underset{h\to 0}{\mathcal{O}}\left(h^{v'}\right)$$

$$\Leftrightarrow \forall k \in [\![0, v'-1]\!],\ \sum_{j'=1}^{s} w_{N,j,j'} \tau_{j'}^{k} = \xi_{j,N}^{N+k} \frac{k!N!}{(k+N)!}$$

The same kind of analysis can be done in the case where $f$ is autonomous, but the conditions are not easy to find. John C. Butcher has found for Runge-Kutta methods a tree representation[1] of the order conditions in the case where $f$ is a function of only $y$, and Hairer, Wanner and Nørsett[4] generalized this tree representation to Runge-Kutta-Nyström methods (the case $n = 2$). It may be possible to do so for any $n$, especially if one uses the proposition of Section III.1.26.

The solved system conditions allow to link general multi-order Runge-Kutta methods to what motivated their definition, gaussian quadratures, which allows for an upper bound on the order of consistency of a method.

**III.1.20 Theorem : Solved system conditions & Gauss-Jacobi quadratures**

Let $(n, s) \in \mathbb{N}^* \times \mathbb{N}^*$, $M = (\tau, w, \tilde{w}) \in \mathrm{GMORK}_{n,s}$, $\xi$ its approximation times, $v' \in \mathbb{N}^*$, and $(j, N) \in [\![1, s+1]\!] \times [\![1, n]\!]$ such that $\xi_{j,N} \neq 0$. We define, for all $j' \in [\![1, s]\!]$ :

$$b_{j'} = \frac{2^N}{\xi_{j,N}^N} \frac{w_{N,j,j'}}{N}, \ c_{j'} = \tau_{j'} \frac{2}{\xi_{j',N}} - 1$$

$M$ satisfies the solved system conditions of order of consistency $v'$ for length $n$ at stage $j$ and rank $N$ if and only if $(b, c)$ is a Gauss-Jacobi quadrature of degree $v' - 1$ in the case where the power of $1 + x$ is 0 and the power of $1 - x$ is $N - 1$.

Proof

The intuition for this proof comes from the way we derived the expression of general multi-order Runge-Kutta methods. At their core general multi-order Runge-Kutta methods are Gauss-Jacobi quadratures with a few more steps. A Gauss-Jacobi quadrature with $s$ points $(b, c) \in \mathbb{R}^{[\![1,s]\!]} \times \mathbb{R}^{[\![1,s]\!]}$ approximates, with $g \in C^0\big([-1, 1], \mathbb{R}^{[\![1,d]\!]}\big)$, and $(\alpha, \beta) \in ]-1, +\infty[$ :

$$\int_{-1}^{1} (1 - x)^{\alpha} (1 + x)^{\beta} g(x) \,\mathrm{d}x \approx \sum_{j'=1}^{s} b_{j'} g\big(c_{j'}\big)$$

The main criterion for the precision of a quadrature is its degree, the highest integer $K \in \mathbb{N}$ such that, for all polynomials $P$ of degree at most $K$, the approximation is exact. The case we consider is $\beta = 0, \ \alpha = N - 1$, the condition for degree $K$ hence becomes :

$$\forall P \in \mathbb{R}_K[X], \ \sum_{j'=1}^{s} b_{j'} P\big(c_{j'}\big) = \int_{-1}^{1} (1 - x)^{N-1} P(x) \,\mathrm{d}x$$

We linearly transform this integral so that the domain of integration is $[0, 1]$.

$$\forall P \in \mathbb{R}_K[X], \ \sum_{j'=1}^{s} b_{j'} P\big(c_{j'}\big) = 2^N \int_{0}^{1} (1 - x)^{N-1} P(2x - 1) \,\mathrm{d}x$$

It is necessary and sufficient that this condition be satisfied on a basis of $\mathbb{R}_K[X]$. To simplify the integral, we use the basis $\left(\left(\frac{x+1}{2}\right)^k\right)_{k \in [\![0,K]\!]}$.

$$\forall k \in [\![0, K]\!], \ \sum_{j'=1}^{s} b_{j'} \left(\frac{c_{j'} + 1}{2}\right)^k = 2^N \int_{0}^{1} (1 - x)^{N-1} x^k \,\mathrm{d}x$$

The value of the integral is a classical result, it is exactly $\mathrm{B}(k + 1, N)$ where $\mathrm{B}$ is the beta function. Since we are dealing with integers its expression is simple, $\mathrm{B}(k + 1, N) = \frac{k!(N-1)!}{(k+N)!}$, we thus get :

$$\forall k \in [\![0, K]\!],\ \sum_{j'=1}^{s} b_{j'} \left(\frac{c_{j'}+1}{2}\right)^k = 2^N \frac{k!(N-1)!}{(k+N)!}$$

This system is very similar to the solved system conditions. In fact, we can simply multiply both sides by a non-zero constant to retrieve them.

$$\forall k \in [\![0, K]\!],\ \sum_{j'=1}^{s} N \frac{\xi_{j,N}^N}{2^N} b_{j'} \left(\xi_{j',N} \frac{c_{j'}+1}{2}\right)^k = \xi_{j,N}^{N+k} \frac{k!N!}{(k+N)!}$$

By defining :

$$K = v' - 1,\ \forall j' \in [\![1, s]\!],\ w_{N,j,j'} = N \frac{\xi_{j,N}^N}{2^N} b_{j'},\ \tau_{j'} = \xi_{j',N} \frac{c_{j'}+1}{2}$$

We exactly get the solved system conditions of order of consistency $v' - 1$ for length $n$ at rank $N$ and stage $j$. If we solve for $b, c$, we get :

$$K = v' - 1,\ \forall j' \in [\![1, s]\!],\ b_{j'} = \frac{2^N}{\xi_{j,N}^N} \frac{w_{N,j,j'}}{N},\ c_{j'} = \tau_{j'} \frac{2}{\xi_{j',N}} - 1$$

We deduce that a method satisfies the solved system conditions for order of consistency $v$ for length $n$ if and only if $(b, c)$ is a Gauss-Jacobi quadrature of degree $v - 1$ in the case $\beta = 0,\ \alpha = N - 1$.

This result is important because Gauss-Jacobi quadratures have a limit to their degree.

### III.1.21 Theorem : Upper bound on the order of consistency

Let $(n, s) \in \mathbb{N}^* \times \mathbb{N}^*$, $M \in \mathrm{GMORK}_{n,s}$, $\xi$ the approximation times of $M$, $v_n$ the minimum consistency matrix of $M$, $(j, N) \in [\![1, s+1]\!] \times [\![1, n]\!]$, and $s'$ the number of different nodes.
If $\xi_{j,N} \neq 0$ then $v_{n,j,N} \leq 2s'$.

Proof

The maximal degree of a Gauss-Jacobi quadrature with $s$ points is $2s - 1$. $v_{n,j,N} > 2s$ would require the Gauss-Jacobi quadrature of Section III.1.20 $(b, c)$ to have a degree strictly higher than $2s - 1$, which is not possible. To get the bound $2s'$ instead of $2s$, we observe that a gaussian quadrature is equivalent to the gaussian quadrature where the weights with equal nodes have been factorized. For all $(j, j') \in [\![1, s]\!]^2,\ c_j = c_{j'} \Leftrightarrow \tau_j = \tau_{j'}$. We deduce that the number of different nodes of $(b, c)$ is equal to the number of different nodes of $M$, hence the bound $2s'$.

It is possible for a method to achieve greater orders of consistency, but there are necessary conditions.

### III.1.22 Proposition : Necessary conditions for greater order of consistency

Let $(n, s) \in \mathbb{N}^* \times \mathbb{N}^*$, $M = (\tau, w, \tilde{w}) \in \mathrm{GMORK}_{n,s}$, $v_n$ the minimum consistency matrix of $M$, $(j, N) \in [\![1, s+1]\!] \times [\![1, n]\!]$, $s' \in \mathbb{N}^*$ the number of different nodes, $\tau'_1, ..., \tau'_{s'}$ those different nodes. We define $w'_{N,j} \in \mathbb{R}^{[\![1, s']\!]}$ as :

$$\forall j' \in [\![1, s']\!],\ w'_{N,j,j'} = \sum_{\substack{j'' \in [\![1, s]\!] \\ \tau'_{j'} = \tau_{j''}}} w_{N,j,j''}$$

If $v_{n,j,N} > 2s'$ then $\xi_{j,N} = 0,\ w'_{N,j} = 0_{s'}$, and, for all $N' \in [\![1, N]\!],\ \tilde{w}_{N,N',j} = 0^{N-N'}$.

Proof

We have $\xi_{j,N} = 0$, otherwise the previous theorem gives a contradiction. The conditions for negative order of consistency give, for all $N' \in [\![1, N]\!]$, $\tilde{w}_{N,N',j} = 0^{N-N'}$. The solved system conditions imply :

$$\begin{bmatrix} 1 & 1 & \dots & 1 \\ \tau'_1 & \tau'_2 & \dots & \tau'_s \\ \vdots & \vdots & \ddots & \vdots \\ \tau'^{s'-1}_1 & \tau'^{s'-1}_2 & \dots & \tau'^{s'-1}_{s'} \end{bmatrix} w'_{N,j} = 0_{s'}$$

The left term is a Vandermonde matrix where all the node are different, it is hence invertible, and by multiplying by its inverse we find $w'_{N,j} = 0_{s'}$.

The difference between $\xi_{j,N} = 0$ and $\xi_{j,N} \neq 0$ is further highlighted by the next definition and proposition.

**III.1.23 Definition : Scaling of a method**

Let $(n, s) \in \mathbb{N}^* \times \mathbb{N}^*$, $M = (\tau, w, \tilde{w}) \in \text{GMORK}_{n,s}$, $\lambda \in \mathbb{R}$, and $M' = (\tau', w', \tilde{w}') \in \text{GMORK}_{n,s}$, where :

$$\tau' = \lambda\tau,\ \forall N \in [\![1, n]\!],\ w'_N = \lambda^N w_N,\ \forall N' \in [\![1, N]\!],\ \tilde{w}'_{N,N'} = \lambda^{N-N'} \tilde{w}_{N,N'}$$

We define $M'$ as the scaling of $M$ by the factor $\lambda$.

**III.1.24 Proposition : Properties of the scaling of a method**

Let $(n, s) \in \mathbb{N}^* \times \mathbb{N}^*$, $M = (\tau, w, \tilde{w}) \in \text{GMORK}_{n,s}$, $\lambda \in \mathbb{R}$, and $M' \in \text{GMORK}_{n,s}$ the scaling of $M$ by the factor $\lambda$.

1. Let $f$ a differential equation function of order $n$, and $(t, y_{n,0}, h) \in \mathbf{U} \times \mathbb{R}$.
   The stage system $(M, f, t, y_{n,0}, \lambda h)$ is equivalent to the stage system $(M', f, t, y_{n,0}, h)$.
2. Let $\xi, \xi'$ the approximation times of respectively $M, M'$.
   $\xi' = \lambda\xi$.
3. Let $v_n, v'_n$ the minimum consistency matrices of respectively $M, M'$.
   If $\lambda \neq 0$ then $v'_n = v_n$. If $\lambda = 0$ then, for all $(j, N) \in [\![1, s+1]\!] \times [\![1, n]\!]$ :

$$v'_{n,j,N} = \begin{cases} +\infty \text{ if } \tilde{w}_{N,N,j} = 1 \\ -N \text{ otherwise} \end{cases}$$

Proof

1. $M'$ is the special case where we use $M$ with the step size multiplied by $\lambda$. The stage system $(M, f, t, y_{n,0}, \lambda h)$ is :

$$\forall (j, N) \in [\![1, s+1]\!] \times [\![1, n]\!],$$
$$y_{n,j,N} = \sum_{N'=1}^{N} \tilde{w}_{N,N',j} (h\lambda)^{N-N'} y_{n,0,N'} + \frac{(h\lambda)^N}{N!} \sum_{j'=1}^{s} w_{N,j,j'} f\big(t + \tau_{j'} h\lambda, y_{n,j'}\big)$$
$$y_{n,s+1} \in \mathbb{R}^{[\![1,n]\!] \times [\![1,d]\!]},\ \forall j \in [\![1, s]\!],\ \big(t + \tau_j h, y_{n,j}\big) \in \mathbf{U}$$

   Which is the stage system $(M', f, t, y_{n,0}, h)$.
2. Trivial.
3. Let $(j, N) \in [\![1, s+1]\!] \times [\![1, n]\!]$. If $\lambda = 0$, the approximation is exact if and only if $\tilde{w}_{N,N,j} = 1$. If true it implies an infinite order of consistency, otherwise Section III.1.8 implies that it is $-N$. If $\lambda \neq 0$, $M'$ is a special case of $M$, hence the order of consistency of $M'$ is at least that of $M$. We can observe that $M$ is the scaling of $M'$ by the factor $\frac{1}{\lambda}$, hence the order of consistency of $M$ is at least that of $M'$, which implies they are equal.

If a certain order of consistency is possible for a given approximation time then, for all $\lambda \in \mathbb{R}^*$ there exists a method with the same order of consistency and with approximation times $\lambda\xi$. We cannot escape from 0 by multiplying by a factor, the approximation time 0 is hence somewhat isolated.

If we look at the conditions for a certain of order of consistency in Section B, we may observe a pattern, only the weights at the lowest ranks appear.

**III.1.25 Proposition : Weights involved in conditions for order of consistency**

Let $(n, s) \in \mathbb{N}^* \times \mathbb{N}^*$, $M = (\tau, w, \tilde{w}) \in \mathrm{GMORK}_{n,s}$, $v' \in \mathbb{N}^*$, and $(j, N) \in [\![1, s+1]\!] \times [\![1, n]\!]$.
If, for all $(j', N') \in [\![1, s]\!] \times [\![1, n]\!]$, $\tilde{w}_{N',N',j'} = 1$, and $M$ is at least of order of consistency 0 for length $n$ at stage $j$ and rank $N$, then, the condition, $M$ satisfies the equation for order of consistency $v'$ of class $v'$ for length $n$ at stage $j$ and rank $N$, depends only on $w_{N,j}$, $\tau$, for all $N' \in [\![1, \min(n, v'-1)]\!]$, $w_{N'}$, and, for all $N' \in [\![1, n]\!]$, for all $N'' \in [\![1 + \max(0, N' - v'), N']\!]$, $\tilde{w}_{N',N''}$.

Proof

The other weights are multiplied by a power of $h$ greater than or equal to $v'$, therefore the conditions for order of consistency do not depend on these weights.

We can get more precise results by assuming a stronger condition.

**III.1.26 Proposition : Approximations involved in equation for order of consistency**

Let $(n, s) \in \mathbb{N}^* \times \mathbb{N}^*$, $M = (\tau, w, \tilde{w}) \in \mathrm{GMORK}_{n,s}$, $Y_n$ the approximation function for length $n$ of $M$, $v_n$ the minimum consistency matrix of $M$, $(j, N) \in [\![1, s+1]\!] \times [\![1, n]\!]$, $v' \in \mathbb{N}^*$, and $\mathcal{D} \in \mathcal{P}([\![1, s]\!] \times [\![1, n]\!])$ such that, for all $(j', N') \in \mathcal{D}$, $v_{n,j',N'} + N' \geq v'$.
If, for all $(j', N') \in [\![1, s]\!] \times [\![1, n]\!]$, $\tilde{w}_{N',N',j'} = 1$, and $M$ is at least of order of consistency 0 for length $n$ at stage $j$ and rank $N$, then, $M$ satisfies the equation for order of consistency $v'$ of class $v'$ for length $n$ at stage $j$ and rank $N$ if and only if, for all $(t, y_{n,0}) \in \mathbb{R} \times \mathbb{R}^{[\![1,n]\!] \times [\![1,d]\!]}$, for all neighborhoods $\mathbf{U} \in \mathcal{P}(\mathbb{R} \times \mathbb{R}^{[\![1,n]\!] \times [\![1,d]\!]})$ of $(t, y_{n,0})$, for all $f \in C^{v'}(\mathbf{U}, \mathbb{R}^{[\![1,d]\!]})$, if we note $\hat{y}$ the unique maximal solution of $(f, t, y_{n,0})$, and if we define :

$$\forall (j, N) \in [\![1, s]\!] \times [\![1, n]\!],\ Y'_{n,f,j,N}(t, y_{n,0}, h) = \begin{cases} \frac{\mathrm{d}^{n-N}\hat{y}}{\mathrm{d}t^{n-N}}(t + \xi_{j,N} h) & \text{if } (j, N) \in \mathcal{D} \\ Y_{n,f,j,N}(t, y_{n,0}, h) & \text{otherwise} \end{cases}$$

Then :

$$\begin{aligned} &\left\{h \in \mathcal{H}_{1,f,t,y_{n,0}} \mid t + \xi_{j,N} h \in \Omega\right\} \to \mathbb{R}^{[\![1,d]\!]} \\ &h \to \sum_{j'=1}^{s} w_{N,j,j'} \sum_{k=0}^{v'-1} \frac{1}{k!} \mathrm{d}^k f(t, y_{n,0})\left(\tau_{j'} h, Y'_{n,f,j'}(t, y_{n,0}, h) - y_{n,0}\right) \\ &\quad - \sum_{k=0}^{v'-1} \xi_{j,N}^{N+k} h^k \frac{N!}{(k+N)!} \frac{\mathrm{d}^k f \circ \overline{\mathcal{J}}^{n-1}\hat{y}}{\mathrm{d}t^k}(t) = \underset{h \to 0}{\mathcal{O}}(h^v) \end{aligned}$$

Proof

$Y'$ is defined to satisfy $h \to Y_{n,f}(t, y_{n,0}, h) = Y'_{f,n}(t, y_{n,0}, h) + \underset{h \to 0}{\mathcal{O}}(h^v)$, hence :

$$
\begin{aligned}
& h \to \sum_{j'=1}^{s} w_{N,j,j'} \sum_{k=0}^{v'-1} \frac{1}{k!} \mathrm{d}^k f(t, y_{n,0}) \big(\tau_{j'} h, Y_{n,f,j'}(t, y_{n,0}, h) - y_{n,0}\big) \\
& \quad - \sum_{k=0}^{v'-1} \xi_{j,N}^{N+k} h^k \frac{N!}{(k+N)!} \frac{\mathrm{d}^k f \circ \overline{\mathcal{J}}^{n-1} \hat{\mathrm{y}}}{\mathrm{d}t^k}(t) = \underset{h\to 0}{\mathcal{O}}(h^v) \\
\Leftrightarrow\ & h \to \sum_{j'=1}^{s} w_{N,j,j'} \sum_{k=0}^{v'-1} \frac{1}{k!} \mathrm{d}^k f(t, y_{n,0}) \Big(\tau_{j'} h, Y'_{f,n}(t, y_{n,0}, h) + \underset{h\to 0}{\mathcal{O}}(h^v) - y_{n,0}\Big) \\
& \quad - \sum_{k=0}^{v'-1} \xi_{j,N}^{N+k} h^k \frac{N!}{(k+N)!} \frac{\mathrm{d}^k f \circ \overline{\mathcal{J}}^{n-1} \hat{\mathrm{y}}}{\mathrm{d}t^k}(t) = \underset{h\to 0}{\mathcal{O}}(h^v) \\
\Leftrightarrow\ & h \to \sum_{j'=1}^{s} w_{N,j,j'} \sum_{k=0}^{v'-1} \frac{1}{k!} \mathrm{d}^k f(t, y_{n,0}) \big(\tau_{j'} h, Y'_{f,n}(t, y_{n,0}, h) - y_{n,0}\big) \\
& \quad - \sum_{k=0}^{v'-1} \xi_{j,N}^{N+k} h^k \frac{N!}{(k+N)!} \frac{\mathrm{d}^k f \circ \overline{\mathcal{J}}^{n-1} \hat{\mathrm{y}}}{\mathrm{d}t^k}(t) = \underset{h\to 0}{\mathcal{O}}(h^v)
\end{aligned}
$$

The previous theorems and propositions can be used on sub-methods induced by reachable-closed sets of stages, which proves particularly useful for explicit methods. We now see what informations we can gather on the order of consistency of a method by studying its permutations.

**III.1.27 Proposition : Order of consistency & Permutations**

Let $(n, s) \in \mathbb{N}^* \times \mathbb{N}^*$, $M \in \mathrm{GMORK}_{n,s}$, $\varphi \in \mathfrak{S}_s^*$, and $v_n, v'_n$ the minimum consistency matrices of respectively $M, \varphi * M$.
For all $j \in [\![1, s+1]\!]$, $v'_{n,j} = v_{n,\varphi(j)}$

Proof

A permutation simply relabels the stages of a method, which implies this equality.

**III.1.28 Proposition : Order of consistency & Sub-methods**

Let $(n, s) \in \mathbb{N}^* \times \mathbb{N}^*$, $M \in \mathrm{GMORK}_{n,s}$, $\mathcal{R}_r$ the reachability relation of $M$, $J \in \mathcal{P}([\![1, s]\!])$ a reachable closed set, $j_l \in [\![1, s+1]\!] \setminus J$, $s' = \#(J)$, $\sigma \in \{J \cup \{j_l\} \to [\![1, s'+1]\!]\}$ a labeling of $(J, j_l)$, $M' \in \mathrm{GMORK}_{n,s}$ the sub-method induced by $(J, j_l, \sigma)$, and $v_n, v_{n'}$ the minimum consistency matrixes of respectively $M, M'$.
For all $j \in [\![1, s']\!]$, $v'_{n,j} = v_{n,\sigma^{-1}(j)}$. If $\mathcal{R}_r^{-1}(j_l) \subset J$, then $v'_{n,s'+1} = v_{n,\sigma^{-1}(s'+1)}$

Proof

It is a direct consequence of Section II.2.21.

This proposition allows us to use the results on the order of consistency of a general multi-order Runge-Kutta method on stages in a reachable-closed set of stages.

The issue with the notion of order of consistency is that it is useful only when the step size is small. If we want to approximate the solution far from $t_0$, we need to break the step size $h'$ into multiple smaller steps $h \in \mathbb{R}^{[\![0, \tilde{q}-1]\!]}$ such that $h' = \sum h_k$. We can then make a first step using the step size $h_0$, then use the result of this approximation as the initial values for another step which this times uses $h_1$ as a step size, then repeat *ad nauseam*.

**III.1.29 Definition : Step size sequences & Time mesh**

Let $(n,s) \in \mathbb{N}^* \times \mathbb{N}^*$, $M \in \mathrm{GMORK}_{n,s}$, $f$ a differential equation function of order $n$, $Y_{n,f}$ the approximation function for length $n$ of $M$, $\tilde{q} \in \mathbb{N}^*$, and $\left(t_0, y_{n,0}\right) \in \mathbf{U}$.

- Let $h \in \mathbb{R}^{[\![0,\tilde{q}-1]\!]}$. $h$ is said to be a step size sequence of length $\tilde{q}$ if and only if, if we note $t \in \mathbb{R}^{[\![0,\tilde{q}]\!]}$ the sequence defined as, for all $q \in [\![0, \tilde{q}-1]\!]$, $t_{q+1} = t_q + h_q$, then the sequence $\left(y_{n,q(s+1)}\right)_{q \in [\![1,\tilde{q}]\!]}$ defined as, for all $q \in [\![0, \tilde{q}-1]\!]$, $y_{n,(q+1)(s+1)} = Y_{n,f,s+1}\left(t_q, y_{n,q(s+1)}, h_q\right)$, exists, meaning, for all $q \in [\![0, \tilde{q}-1]\!]$, $\left(t_q, y_{n,q(s+1)}, h_q\right) \in \mathcal{U}_{n,f}$. $t$ is called the time mesh of $(t_0, h)$.
- We define $\mathcal{H}_{\tilde{q},f,t_0,y_{n,0}} \subset \mathcal{P}\left(\mathbb{R}^{[\![0,\tilde{q}-1]\!]}\right)$ the set of sequences of step size sequences of length $\tilde{q} \in \mathbb{N}^*$ at the point $\left(t_0, y_{n,0}\right)$.
- A sequence $h \in \mathbb{R}^{\mathbb{N}}$ is an infinite step size sequence at the point $\left(t_0, y_{n,0}\right)$ if and only if all its truncated sequences are step size sequences at the point $\left(t_0, y_{n,0}\right)$.
- We define $\mathcal{H}_{\infty,f,t_0,y_{n,0}} \in \mathcal{P}\left(\mathbb{R}^{\mathbb{N}}\right)$ as the set of infinite sequences of step sizes at the point $\left(t_0, y_{n,0}\right)$.

To use step size sequences easily, we expand on the notations we previously defined.

**III.1.30 Definition : Prolongation of the approximation function & Evaluation function**

Let $(n,s) \in \mathbb{N}^* \times \mathbb{N}^*$, $M \in \mathrm{GMORK}_{n,s}$, $f$ a differential equation function of order $n$, $Y_{n,f}$ the approximation function for length $n$ of $M$, $F_{n,f}$ the evaluation function for length $n$ of $M$, $\left(t_0, y_{n,0}\right) \in \mathbf{U}$, $\tilde{q} \in \overline{\mathbb{N}}^*$, $h \in \mathcal{H}_{\tilde{q},f,t_0,y_{n,0}}$, and $t$ the time mesh of $(t_0, h)$.

We define, for all $q \in \mathbb{N}$ such that $q < \tilde{q}$ :

$$\forall (j,N) \in [\![1, s+1]\!] \times [\![1, n]\!]$$

$$Y_{n,f,q(s+1)+j,N}\left(t_0, y_{n,0}, h\right) = Y_{n,f,j,N}\left(t_q, Y_{n,f,q(s+1)}\left(t_0, y_{n,0}, h\right), h_q\right)$$

$$\forall j \in [\![1, s]\!],\ F_{n,f,q(s+1)+j}\left(t_0, y_{n,0}, h\right) = f\left(t_q + \tau_j h_q, Y_{n,f,q(s+1)+j}\left(t_0, y_{n,0}, h\right)\right)$$

**III.1.31 Definition : Prolongation of the error function**

Let $(n,s) \in \mathbb{N}^* \times \mathbb{N}^*$, $M \in \mathrm{GMORK}_{n,s}$, $\left(f, t_0, y_{n,0}\right)$ an initial value problem of order $n$ with a unique solution $\hat{y}$, $Y_{n,f}$ the approximation function for length $n$ of $M$, $\left(t'_0, y'_{n,0}\right) \in \mathbf{U}$, $\tilde{q} \in \overline{\mathbb{N}}^*$, $h \in \mathcal{H}_{\tilde{q},f,t'_0,y'_{n,0}}$, $t'$ the time mesh of $(t'_0, h)$, and $\varepsilon_n$ the error function of $M$ for length $n$.

We define, for all $q \in \mathbb{N}$ such that $q < \tilde{q}$, for all $(j,N) \in [\![1, s+1]\!] \times [\![1, n]\!]$, if $t'_q + \xi_{j,N} h \in \Omega$ :

$$\varepsilon_{n,f,q(s+1)+j,N}\left(t'_0, y'_{n,0}, h\right) = \frac{\mathrm{d}^{n-N}\hat{y}}{\mathrm{d}t^{n-N}}\left(t'_q + \xi_{j,N} h_q\right) - Y_{n,f,q(s+1)+j,N}\left(t'_0, y'_{n,0}, h\right)$$

We now state how a general multi-order Runge-Kutta method uses a step size sequence.

**III.1.32 Proposition : Use of general multi-order Runge-Kutta methods**

Let $(n,s) \in \mathbb{N}^* \times \mathbb{N}^*$, $M = (\tau, w, \tilde{w}) \in \mathrm{GMORK}_{n,s}$, $f$ a differential equation function of order $n$, $\left(t_0, y_{n,0}\right) \in \mathbf{U}$, $\tilde{q} \in \overline{\mathbb{N}}^*$, $h \in \mathcal{H}_{\tilde{q},f,t,y_{n,0}}$, and $t$ the time mesh of $(t_0, h)$.

$M$ approximates, for all $q \in \mathbb{N}$ such that $q < \tilde{q}$ :

$$\forall (j,N) \in [\![1, s+1]\!] \times [\![1, n]\!],$$

$$y_{n,q(s+1)+j,N} = \sum_{N'=1}^{N} \tilde{w}_{N,N',j} h_q^{N-N'} y_{n,q(s+1),N'} + \frac{h_q^N}{N!} \sum_{j'=1}^{s} w_{N,j,j'} f\left(t_q + \tau_{j'} h_q, y_{n,q(s+1)+j'}\right)$$

We need to make sure that the process of breaking a large step size into multiple step sizes actually works. We therefore introduce the most fundamental criterion of a general multi-order Runge Kutta method, convergence. In the definition of convergence of Runge Kutta methods, only constant step size sequences were considered, but here we also consider non constant step size sequences since what we really want from the concept of convergence is, if we reduce enough the size of the steps as the number of steps grows, the approximations uniformly converge towards the solution.

**III.1.33 Definition : Convergence**

Let $(n, s) \in \mathbb{N}^* \times \mathbb{N}^*$, $M \in \mathrm{GMORK}_{n,s}$, and $\varepsilon_n$ the error function of $M$ for length $n$.
$M$ is said to be convergent up to length $n$ if and only if, for all initial value problems $\left(f, t_0, y_{n,0}\right)$ of order $n$ with a unique solution and with $f$ continuous and globally Lipschitz continuous in its second variable, for all step size sequence $h \in \prod_{\tilde{q}=1}^{+\infty} \mathcal{H}_{\tilde{q},f,t_0,y_{n,0}}$, if $\tilde{q} \to \|h_{\tilde{q}}\|_{\tilde{q}} = \underset{\tilde{q}\to+\infty}{\mathcal{O}}\left(\frac{1}{\tilde{q}}\right)$ and there exists $\Omega'$ a closed subset of $\Omega$ such that, for all $\tilde{q} \in \mathbb{N}^*$, the times mesh of $\left(t_0, h_{\tilde{q}}\right)$ stays in $\Omega'$, then :

$$\lim_{\tilde{q}\to+\infty} \sup_{q\in[\![0,\tilde{q}]\!]} \|\varepsilon_{n,f,q(s+1)}\left(t_0, y_{n,0}, h_{\tilde{q}}\right)\|_{n,d} = 0$$

The big $\mathcal{O}$ condition on the step size sequence prevents the time meshes from tending towards infinity. We assume that $f$ is Lipschitz continuous because we need it to be somewhat smooth to have actual results. The time meshes must stay in a closed subset of $\Omega$ because, thanks to the big $\mathcal{O}$ condition, we get that the time meshes stay in a compact subset of $\Omega$, which is necessary to use the extreme value theorem.

$f$ is assumed to be globally Lipschitz continuous because it seems impossible to restrict $f$ to a compact set, which is necessary if we only want to assume local Lipschitz continuity. If we were to try, we would need to consider a compact set around the solution curve $\overline{\mathcal{J}}^{n-1}\hat{\mathrm{y}}$ and we would need to make sure that the approximations stay inside this compact set. The compact set cannot be a closed ball around $\left(t_0', y_{n,0}'\right)$ since, in general, the solution do not stay in such a ball, but we can still prove that such a compact set around the solution curve exists.

**III.1.34 Proposition : Thickened solution curve**

Let $n \in \mathbb{N}^*$, $\left(f, t_0, y_{n,0}\right)$ an initial value problem with $f$ continuous, $\hat{\mathrm{y}}$ a solution of the initial value problem, and $\Omega'$ a compact subset of $\Omega$.
There exists $r \in \mathbb{R}_+^*$ such that $\bigcup_{t\in\Omega'} \overline{\mathbf{B}}\left(\overline{\mathcal{J}}^{n-1}\hat{\mathrm{y}}(t), r\right) \subset \mathbf{U}$.

Proof

Let $d_{\mathbf{U}^c} : \mathbf{U} \to \mathbb{R}_+^*, x \to \inf_{y\in\mathbf{U}^c}\{d(x,y)\}$ the function which, given a point $x \in \mathbf{U}$, returns the maximum radius of an open ball around $x$. We have, for all $(x, y), \in \mathbf{U}^2$, $d_{\mathbf{U}^c}(x) \leq d(x,y) + d_{\mathbf{U}^c}(y)$. We deduce that $d_{\mathbf{U}^c}$ is Lipschitz continuous, and thus continuous. Since the solution curve is continuous, $d_{\mathbf{U}^c} \circ \overline{\mathcal{J}}^{n-1}\hat{\mathrm{y}}$ is also continuous. $d_{\mathbf{U}^c} \circ \overline{\mathcal{J}}^{n-1}\hat{\mathrm{y}}$ is continuous on the compact set $\Omega'$, it thus reaches a minimum $r \in \mathbb{R}_+$.

This minimum cannot be 0 since it would mean that there exists a sequence of points which tends towards the boundary of $\mathbf{U}$, which is not possible since we are on a compact subset of $\Omega$. This implies that, for every point $x$ of the solution curve, the closed ball of radius $\frac{r}{2}$ and center $x$ is contained in $\mathbf{U}$.

$\bigcup_{t\in\Omega'} \overline{\mathbf{B}}\left(\overline{\mathcal{J}}^{n-1}\hat{\mathrm{y}}(t), r\right)$ is a thicken line around the solution curve. Unfortunately, this result is hardly usable. It is possible to prove that the error converges point wise (for a certain step) towards 0, but proving its uniform convergence may be impossible. To keep the approximations in the thicken line we need to impose some restrictions on the size of $\|h_{\tilde{q}}\|_{\tilde{q}}$, which is only possible by increasing $\tilde{q}$, but this adds points at the end that also need to be addressed, which themselves require to shrink $\|h_{\tilde{q}}\|_{\tilde{q}}$, which adds points at the end, etc…

To prove the convergence of a method we first need the following lemma.

**III.1.35 Lemma : Lipschitz continuity of approximation function & Evaluation function**

Let $(n, s) \in \mathbb{N}^* \times \mathbb{N}^*$, $M \in \mathrm{GMORK}_{n,s}$, $f$ a differential equation function of order $n$, continuous and globally Lipschitz continuous in its second variable, $Y_{n,f}$ the approximation function for length $n$ of $M$, and $F_{n,f}$ the evaluation function for length $n$ of $M$.
There exists $h^* \in \mathbb{R}_+^*$ such that $Y_{n,f}$ and $F_{n,f}$ restricted to $\left\{\left(t, y_{n,0}, h\right) \in \mathcal{U}_{n,f} \mid h \leq h^*\right\}$ are continuous and Lipschitz continuous in their second variable.

Proof

We first prove that $Y_{n,f}$ and $F_{n,f}$ are Lipschitz continuous. Let $\left(t, y_{n,0}, z_{n,0}, h\right)$ such that $\left(t, y_{n,0}, h\right) \in \mathcal{U}_{n,f}$, and $\left(t, z_{n,0}, h\right) \in \mathcal{U}_{n,f}$. For all $j \in [\![1, s]\!]$ :

$$\begin{aligned}
&\|F_{n,f,j}\left(t, y_{n,0}, h\right) - F_{n,f,j}\left(t, z_{n,0}, h\right)\|_d \\
&= \|f\left(t + \tau_j h, Y_{n,f,j}\left(t, y_{n,0}, h\right)\right) - f\left(t + \tau_j h, Y_{n,f,j}\left(t, z_{n,0}, h\right)\right)\|_d \\
&\leq L \max_{N \in [\![1,n]\!]} \|Y_{n,f,j,N}\left(t, y_{n,0}, h\right) - Y_{n,f,j,N}\left(t, z_{n,0}, h\right))\|_d \\
&\leq L \max_{N [\![1,n]\!]} \| \sum_{N'=1}^{N} \tilde{w}_{N,N',j} h^{N-N'} \left(y_{n,0,N'} - z_{n,0,N'}\right) \|_d \\
&\quad + L \max_{N [\![1,n]\!]} \| \frac{h^N}{N!} \sum_{j'=1}^{s} w_{N,j,j'} \left(F_{n,f,j'}\left(t, y_{n,0}, h\right) - F_{n,f,j'}\left(t, z_{n,0}, h\right)\right) \|_d \\
&\leq L \max(1, |h|^{n-1}) \ \|y_{n,0} - z_{n,0}\|_{n,d} \ C_{\tilde{w}} \\
&\quad + L \max(|h|, |h|^n) \ \|F_{n,f}\left(t, y_{n,0}, h\right) - F_{n,f,1}\left(t, z_{n,0}, h\right)\|_{s,d} \ C_w
\end{aligned}$$

Thus :

$$\begin{aligned}
&(1 - L \max(|h|, |h|^n) C_w) \ \|F_{n,f}\left(t, y_{n,0}, h\right) - F_{n,f}\left(t, z_{n,0}, h\right)\|_{s,d} \\
&\leq L \max(1, |h|^{n-1}) C_{\tilde{w}} \ \|y_{n,0} - z_{n,0}\|_{n,d}
\end{aligned}$$

Let $a \in \mathbb{R}_+^*$. We want $a \leq 1 - L \max(|h|, |h|^n) C_w$, which is equivalent to $\max(|h|, |h|^n) \leq \frac{1-a}{L C_w}$. So it is necessary and sufficient that $a \in ]0, 1[$. We now restrict $h$ such that it satisfies this condition. We have :

$$\|F_{n,f}\left(t, y_{n,0}, h\right) - F_{n,f}\left(t, z_{n,0}, h\right)\|_{s,d} \leq \frac{L}{a} \max(1, |h|^{n-1}) C_{\tilde{w}} \ \|y_{n,0} - z_{n,0}\|_{n,d}$$

Thus $F_{n,f}$ is Lipschitz continuous in its second variable, thus $Y_{n,f}$ is Lipschitz continuous in its second variable.

We now prove the continuity of $Y_{n,f}$ and $F_{n,f}$. Let $\mathcal{U}'_{n,f}$ the parametrized existence-uniqueness domain of $M$ and $Y'_{n,f}$ the parametrized approximation function. Since we now assume that $f$ is globally Lipschitz continuous in its second variable, the proof of Section II.1.14 implies that the condition $h < \min\left(\frac{1}{L C_w}, \frac{1}{\sqrt[n]{L C_w}}\right)$ is sufficient for $Y'_{n,f}$ to be a contraction in its second variable with a Lipschitz constant $L' \in ]0, 1[$, for all $\left(t, y_{n,0}\right) \in \mathbf{U}$. For all $((t, x, h), (t', x', h')) \in \mathcal{U}_{n,f}^2$ such that $h < \min\left(\frac{1}{L C_w}, \frac{1}{\sqrt[n]{L C_w}}\right)$ and $h' < \min\left(\frac{1}{L C_w}, \frac{1}{\sqrt[n]{L C_w}}\right)$, since $Y_{n,f}$ is the fixed point of $Y'_{n,f}$, we have :

$$\begin{aligned}
&\|Y_{n,f}(t, x, h) - Y_{n,f}(t', x', h')\|_{s+1,n,d} \\
&= \|Y'_{n,f}\left(t, x, h, Y_{n,f}(t, x, h)\right) - Y'_{n,f}\left(t', x', h', Y_{n,f}(t', x', h')\right)\|_{s+1,n,d}
\end{aligned}$$

$$\leq \|Y'_{n,f}\big(t,x,h,Y_{n,f}(t,x,h)\big) - Y'_{n,f}\big(t,x,h,Y_{n,f}(t',x',h')\big)\|_{s+1,n,d}$$
$$+ \|Y'_{n,f}\big(t,x,h,Y_{n,f}(t',x',h')\big) - Y'_{n,f}\big(t',x',h',Y_{n,f}(t',x',h')\big)\|_{s+1,n,d}$$
$$\leq L' \ \|Y_{n,f}(t,x,h) - Y_{n,f}(t',x',h')\|_{s+1,n,d}$$
$$+ \|Y'_{n,f}\big(t,x,h,Y_{n,f}(t',x',h')\big) - Y'_{n,f}\big(t',x',h',Y_{n,f}(t',x',h')\big)\|_{s+1,n,d}$$

We deduce :

$$\|Y_{n,f}(t,x,h) - Y_{n,f}(t',x',h')\|_{s+1,n,d}$$
$$\leq \frac{1}{1-L'} \ \|Y'_{n,f}\big(t,x,h,Y_{n,f}(t',x',h')\big) - Y'_{n,f}\big(t',x',h',Y_{n,f}(t',x',h')\big)\|_{s+1,n,d}$$

$Y'_{n,f}$ is continuous in its first three variable, we deduce $\lim_{(t,x,h)\to(t',x',h')} Y_{n,f}(t,x,h) = Y_{n,f}(t',x',h')$, thus $Y_{n,f}$ is continuous, thus $F_{n,f}$ is continuous.

### III.1.36 Theorem : Conditions for convergence

Let $(n,s) \in \mathbb{N}^* \times \mathbb{N}^*$, $M = (\tau, w, \tilde{w}) \in \mathrm{GMORK}_{n,s}$, and $\xi$ the approximation times of $M$.

1. If $M$ is convergent up to length $n$, then $\xi_{s+1} = 1_n$ and, for all $N \in [\![1,n]\!]$, $\tilde{w}_{N,N,s+1} = 1$
2. Let $J = \big\{j \in [\![1,s]\!] \mid \exists N \in [\![1,n]\!],\ w_{N,s+1,j} \neq 0\big\}$. For $M$ to be convergent up to length $n$ it is sufficient that :

$$\xi_{s+1} = 1_n,\ \forall N \in [\![1,n]\!],\ \tilde{w}_{N,N,s+1} = 1,\ \forall (j,N) \in J \times [\![1,n]\!],\ \tilde{w}_{N,N,j} = 1$$

Proof

1. The condition $\forall N \in [\![1,n]\!],\ \tilde{w}_{N,N,s+1} = 1$ can be found by considering any initial value problem with initial values different from 0, and by considering only null step size sequences.
   We now prove the condition $\xi_{s+1} = 1_n$. Let $N \in [\![2,n]\!]$, and $\big(f, 0, y_{n,0}\big)$ with :
$$\forall (t,x) \in \mathbb{R} \times \mathbb{R}^{[\![1,n]\!]},\ f(t,x) = 0$$
$$\forall N' \in [\![1,n]\!],\ y_{n,N'} = \begin{cases} 1 \text{ if } N' = N-1 \\ 0 \text{ otherwise} \end{cases}$$
   The Cauchy-Lipschitz theorem implies that there exists a unique solution. The solution satisfies, for all $t \in \mathbb{R}$, $\frac{\mathrm{d}^{n-N}\hat{y}}{\mathrm{d}t^{n-N}}(t) = t$. We consider $h \in \prod_{\tilde{q}=1}^{+\infty} \mathcal{H}_{\tilde{q},f,t,y_{n,0}}$ such that, for all $\tilde{q} \in \mathbb{N}^*$, $h_{\tilde{q}} = \big(\tilde{q}^{-1}\big)_{q\in[\![0,\tilde{q}-1]\!]}$. For all $\tilde{q} \in \mathbb{N}^*$, for all $q \in [\![0,\tilde{q}]\!]$, we have $y_{n,q(s+1),N} = \frac{q}{\tilde{q}}\tilde{w}_{N,N-1,s+1}$, thus $\varepsilon_{n,f,\tilde{q}(s+1),N}\big(0, y_{n,0}, h_{\tilde{q}}\big) = 1 - \tilde{w}_{N,N-1,s+1}$, which tends towards 0 if and only if $\tilde{w}_{N,N-1,s+1} = 1$.
   The condition $\xi_{s+1,1} = 1$ is found similarly by taking $y_{n,0} = 0_n$ and $f$ identically 1.
2. Let $\big(f, t_0, y_{n,0}\big)$ an initial value problem, and $h$ a sequence of step sizes which satisfies the conditions of the definition of convergence. Let $\tilde{q} \in \mathbb{N}^*$, and $\hat{y}$ the solution of the initial value problem. Assume $\tilde{q}$ is great enough to have $\|h_{\tilde{q}}\|_{\tilde{q}} \leq 1$ and $F$ continuous and Lipschitz continuous in its second variable with a constant $L \in \mathbb{R}_+^*$. Let $(q,N) \in [\![0,\tilde{q}-1]\!] \times [\![1,n]\!]$. We have :
$$\|\varepsilon_{n,f,(q+1)(s+1),N}\big(t_0, y_{n,0}, h_{\tilde{q}}\big)\|_d = \|\frac{\mathrm{d}^{n-N}\hat{y}}{\mathrm{d}t^{n-N}}\big(t_{\tilde{q},q+1}\big) - Y_{n,f,(q+1)(s+1),N}\big(t_0, y_{n,0}, h_{\tilde{q}}\big)\|_d$$
$$\leq \|\frac{\mathrm{d}^{n-N}\hat{y}}{\mathrm{d}t^{n-N}}\big(t_{\tilde{q},q+1}\big) - Y_{n,f,s+1,N}\big(\mathcal{J}^{n-1}\hat{y}\big(t_{\tilde{q},q}\big), h_{\tilde{q},q}\big)\|_d$$
$$+ \|Y_{n,f,s+1,N}\big(\mathcal{J}^{n-1}\hat{y}\big(t_{\tilde{q},q}\big), h_{\tilde{q},q}\big) - Y_{n,f,(q+1)(s+1),N}\big(t_0, y_{n,0}, h_{\tilde{q}}\big)\|_d$$

$$
\begin{aligned}
&\leq \Big\| \frac{\mathrm{d}^{n-N}\hat{\mathrm{y}}}{\mathrm{d}t^{n-N}}\big(t_{\tilde{q},q+1}\big) - Y_{n,f,s+1,N}\Big(\overline{\mathcal{J}}^{n-1}\hat{\mathrm{y}}\big(t_{\tilde{q},q}\big), h_{\tilde{q},q}\Big)\Big\|_d \\
&\quad + \Big\| \sum_{N'=1}^{N} \tilde{w}_{N,N',s+1} h_{\tilde{q},q}^{N-N'} \left( \frac{\mathrm{d}^{n-N'}\hat{\mathrm{y}}}{\mathrm{d}t^{n-N'}}\big(t_{\tilde{q},q}\big) - Y_{n,f,q(s+1),N'}\big(t_0, y_{n,0}, h_{\tilde{q}}\big) \right) \Big\|_d \\
&\quad + \Big\| \frac{h_{\tilde{q},q}^{N}}{N!} \sum_{j=1}^{s} w_{N,s+1,j}\Big(F_{n,f,j}\big(t_{\tilde{q},q}, \mathcal{J}^{n-1}\hat{\mathrm{y}}\big(t_{\tilde{q},q}\big), h_{\tilde{q},q}\big) - F_{n,f,j}\Big(t_{\tilde{q},q}, Y_{n,f,q(s+1)}\big(t_0, y_{n,0}, h_{\tilde{q}}\big), h_{\tilde{q},q}\Big)\Big)\Big\|_d \\
&\leq \Big\| \frac{\mathrm{d}^{n-N}\hat{\mathrm{y}}}{\mathrm{d}t^{n-N}}\big(t_{\tilde{q},q+1}\big) - Y_{n,f,s+1,N}\Big(\overline{\mathcal{J}}^{n-1}\hat{\mathrm{y}}\big(t_{\tilde{q},q}\big), h_{\tilde{q},q}\Big)\Big\|_d \\
&\quad + |\tilde{w}_{N,N,s+1}| \ \Big\| \frac{\mathrm{d}^{n-N}\hat{\mathrm{y}}}{\mathrm{d}t^{n-N}}\big(t_{\tilde{q},q}\big) - Y_{n,f,q(s+1),N}\big(t_0, y_{n,0}, h_{\tilde{q}}\big)\Big\|_d \\
&\quad + \sum_{N'=1}^{N-1} |\tilde{w}_{N,N',s+1}| \ |h_{\tilde{q},q}|^{N-N'} \ \Big\| \frac{\mathrm{d}^{n-N'}\hat{\mathrm{y}}}{\mathrm{d}t^{n-N'}}\big(t_{\tilde{q},q}\big) - Y_{n,f,q(s+1),N'}\big(t_0, y_{n,0}, h_{\tilde{q}}\big)\Big\|_d \\
&\quad + |h_{\tilde{q},q}|^{N} \sum_{j=1}^{s} \frac{|w_{N,s+1,j}|}{N!} L \ \Big\| \overline{\mathcal{J}}^{n-1}\hat{\mathrm{y}}\big(t_{\tilde{q},q}\big) - Y_{n,f,q(s+1)}\big(t_0, y_{n,0}, h_{\tilde{q}}\big)\Big\|_{n,d} \\
&\leq \Big\| \frac{\mathrm{d}^{n-N}\hat{\mathrm{y}}}{\mathrm{d}t^{n-N}}\big(t_{\tilde{q},q+1}\big) - Y_{n,f,s+1,N}\Big(\overline{\mathcal{J}}^{n-1}\hat{\mathrm{y}}\big(t_{\tilde{q},q}\big), h_{\tilde{q},q}\Big)\Big\|_d + \big\|\varepsilon_{n,f,q(s+1)}\big(t_0, y_{n,0}, h_{\tilde{q}}\big)\big\|_{n,d} \\
&\quad + C_{\tilde{w}}^{*} \ |h_{\tilde{q},q}| \ \big\|\varepsilon_{n,f,q(s+1)}\big(t_0, y_{n,0}, h_{\tilde{q}}\big)\big\|_{n,d} + L \ |h_{\tilde{q},q}| \ C_w \ \big\|\varepsilon_{n,f,q(s+1)}\big(t_0, y_{n,0}, h_{\tilde{q}}\big)\big\|_{n,d} \\
&\leq \Big\| \frac{\mathrm{d}^{n-N}\hat{\mathrm{y}}}{\mathrm{d}t^{n-N}}\big(t_{\tilde{q},q} + h_{\tilde{q},q}\big) - Y_{n,f,s+1,N}\Big(\overline{\mathcal{J}}^{n-1}\hat{\mathrm{y}}\big(t_{\tilde{q},q}\big), h_{\tilde{q},q}\Big)\Big\|_d \\
&\quad + \Big(1 + (C_{\tilde{w}}^{*} + LC_w) \ \|h_{\tilde{q}}\|_{\tilde{q}}\Big) \ \big\|\varepsilon_{n,f,q(s+1)}\big(t_0, y_{n,0}, h_{\tilde{q}}\big)\big\|_{n,d}
\end{aligned}
$$

We now find a bound on the first term. Let $\Omega'$ the compact subset of $\Omega$ containing the times meshes. Let $N \in [\![2, n]\!]$. We have :

$$
\begin{aligned}
&\Big\| \frac{\mathrm{d}^{n-N}\hat{\mathrm{y}}}{\mathrm{d}t^{n-N}}\big(t'_{\tilde{q},q} + h_{\tilde{q},q}\big) - Y_{n,f,s+1,N}\Big(\overline{\mathcal{J}}^{n-1}\hat{\mathrm{y}}\big(t'_{\tilde{q},q}\big), h_{\tilde{q},q}\Big)\Big\|_d \\
&\leq \sum_{N'=1}^{N-2} \Big| \frac{\mathrm{d}^{n-N'}\hat{\mathrm{y}}}{\mathrm{d}t^{n-N'}}\big(t'_{\tilde{q},q}\big)\Big| \ |h_{\tilde{q},q}|^{N-N'} \ \Big\| \frac{1}{(N-N')!} - \tilde{w}_{N,N',s+1}\Big\|_d \\
&\quad + \frac{|h_{\tilde{q},q}|^{N}}{(N-1)!} \ \Big\| \int_0^1 x^{N-1} \frac{\mathrm{d}^{n}\hat{\mathrm{y}}}{\mathrm{d}t^{n}}\big(t'_{\tilde{q},q} + x h_{\tilde{q},q}\big)\, \mathrm{d}x \Big\|_d + |h_{\tilde{q},q}^{N}| \sum_{j=1}^{s} \frac{|w_{N,s+1,j}|}{N!} \ \Big\| F_{n,f}\Big(\overline{\mathcal{J}}^{n-1}\hat{\mathrm{y}}\big(t'_{\tilde{q},q}\big), h_{\tilde{q},q}\Big)\Big\|_{sd} \\
&\leq |h_{\tilde{q},q}|^{2} \sum_{N'=1}^{N-2} \Big| \frac{1}{(N-N')!} - \tilde{w}_{N,N',s+1}\Big| \sup_{t'' \in \Omega'} \Big( \big\| \mathcal{J}^{n-1}\hat{\mathrm{y}}(t'')\big\|_{n+1,d} \Big) \\
&\quad + |h_{\tilde{q},q}|^{2} \sup_{t'' \in \Omega'} \left( \Big\| \frac{\mathrm{d}^{n}\hat{\mathrm{y}}}{\mathrm{d}t^{n}}(t'')\Big\|_d \right) + |h_{\tilde{q},q}|^{2} \ C_w \sup_{(t'', h') \in \Omega' \times [-\|h_{\tilde{q}}\|_{\tilde{q}}, \|h_{\tilde{q}}\|_{\tilde{q}}]} \Big( F_{n,f}\Big(\overline{\mathcal{J}}^{n-1}\hat{\mathrm{y}}(t''), h'\Big)\Big) \\
&\leq \|h_{\tilde{q}}\|_{\tilde{q}}^{2} \ C'
\end{aligned}
$$

Where $C' \in \mathbb{R}_+^*$. If $N = 1$, since $\xi_{s+1,1} = 1$ and, for all $(N, j) \in [\![1, n]\!] \times J$, $\tilde{w}_{N,N,j} = 1$, we have, for all $x \in \mathbf{U}$, $\sum_{j=1}^{s} w_{1,s+1,j} F_{n,f,j}(x, 0) = \sum_{j=1}^{s} w_{1,s+1,j} f(x) = f(x)$, hence :

$$
\begin{aligned}
&\|\frac{\mathrm{d}^{n-1}\hat{y}}{\mathrm{d}t^{n-1}}\big(t'_{\tilde{q},q}+h_{\tilde{q},q}\big)-Y_{n,f,s+1,1}\Big(\overline{\mathcal{J}}^{n-1}\hat{y}\big(t'_{\tilde{q},q}\big),h_{\tilde{q},q}\Big)\|_d\\
&=\|\int_{t'_{\tilde{q},q}}^{t'_{\tilde{q},q}+h_{\tilde{q},q}}\frac{\mathrm{d}^n\hat{y}}{\mathrm{d}t^n}(x)\,\mathrm{d}x-h_{\tilde{q},q}\sum_{j=1}^{s}w_{1,s+1,j}F_{n,f,j}\Big(\overline{\mathcal{J}}^{n-1}\hat{y}\big(t'_{\tilde{q},q}\big),h_{\tilde{q},q}\Big)\|_d\\
&=\|\int_0^{h_{\tilde{q},q}}f\Big(\overline{\mathcal{J}}^{n-1}\hat{y}\big(t'_{\tilde{q},q}+x\big)\Big)-\sum_{j=1}^{s}w_{1,s+1,j}F_{n,f,j}\Big(\overline{\mathcal{J}}^{n-1}\hat{y}\big(t'_{\tilde{q},q}\big),h_{\tilde{q},q}\Big)\,\mathrm{d}x\|_d\\
&\leq\int_0^{h_{\tilde{q},q}}\|\sum_{j=1}^{s}w_{1,s+1,j}F_{n,f,j}\Big(\overline{\mathcal{J}}^{n-1}\hat{y}\big(t'_{\tilde{q},q}+x\big),0\Big)-\sum_{j=1}^{s}w_{1,s+1,j}F_{n,f,j}\Big(\overline{\mathcal{J}}^{n-1}\hat{y}\big(t'_{\tilde{q},q}\big),h_{\tilde{q},q}\Big)\ \|_d\,\mathrm{d}x\\
&\leq\int_0^{h_{\tilde{q},q}}\sum_{j=1}^{s}|w_{1,s+1,j}|\ \|F_{n,f,j}\Big(\overline{\mathcal{J}}^{n-1}\hat{y}\big(t'_{\tilde{q},q}+x\big),0\Big)-F_{n,f,j}\Big(\overline{\mathcal{J}}^{n-1}\hat{y}\big(t'_{\tilde{q},q}\big),h_{\tilde{q},q}\Big)\|_d\,\mathrm{d}x
\end{aligned}
$$

Since $F$ and $\overline{\mathcal{J}}^{n-1}\hat{y}$ are continuous, the function $\Omega'\times[-\|h_{\tilde{q}}\|_{\tilde{q}},\|h_{\tilde{q}}\|_{\tilde{q}}]\to\mathbb{R}^{[\![1,s]\!]\times[\![1,d]\!]},(t'',h')\to F_{n,f}\Big(\overline{\mathcal{J}}^{n-1}\hat{y}(t''),h'\Big)$ is also continuous. Since $\Omega'\times[-\|h_{\tilde{q}}\|_{\tilde{q}},\|h_{\tilde{q}}\|_{\tilde{q}}]$ is compact, we deduce that the function is uniformly continuous. Thus, for all $\delta\in\mathbb{R}_+^*$, there exists $\varepsilon'\in\mathbb{R}_+^*$ such that, for all $(h',t'')\in\Omega'\times[-\|h_{\tilde{q}}\|_{\tilde{q}},\|h_{\tilde{q}}\|_{\tilde{q}}]$, for all $(h'',t''')\in\Omega'\times[-\|h_{\tilde{q}}\|_{\tilde{q}},\|h_{\tilde{q}}\|_{\tilde{q}}]$ :

$$
|t''-t'''|<\varepsilon',|h'-h''|<\varepsilon'\Rightarrow\|F_{n,f}\Big(\overline{\mathcal{J}}^{n-1}\hat{y}(t''),h'\Big)-F_{n,f}\Big(\overline{\mathcal{J}}^{n-1}\hat{y}(t'''),h''\Big)\|\leq\delta
$$

We can thus define a sequence $\delta\in\mathbb{R}_+^{*\mathbb{N}^*}$ such that $\lim_{\tilde{q}\to+\infty}\delta_{\tilde{q}}=0$ and :

$$
\|\frac{\mathrm{d}^{n-1}\hat{y}}{\mathrm{d}t^{n-1}}\big(t'_{\tilde{q},q}+h_{\tilde{q},q}\big)-Y_{n,f,s+1,1}\Big(\overline{\mathcal{J}}^{n-1}\hat{y}\big(t'_{\tilde{q},q}\big),h_{\tilde{q},q}\Big)\|_d\leq\int_0^{h_{\tilde{q},q}}\sum_{j=1}^{s}|w_{1,s+1,j}|\ \delta_{\tilde{q}}\,\mathrm{d}x
$$

Hence:

$$
\|\frac{\mathrm{d}^{n-1}\hat{y}}{\mathrm{d}t^{n-1}}\big(t'_{\tilde{q},q}+h_{\tilde{q},q}\big)-Y_{n,f,s+1,N}\Big(\overline{\mathcal{J}}^{n-1}\hat{y}\big(t'_{\tilde{q},q}\big),h_{\tilde{q},q}\Big)\|_d\leq\delta_{\tilde{q}}C_w\ \|h_{\tilde{q}}\|_{\tilde{q}}
$$

Let $\delta'_{\tilde{q}}=\max\big(\|h_{\tilde{q}}\|_{\tilde{q}}\ C',\delta_{\tilde{q}}C_w\big)$. We have :

$$
\begin{aligned}
\|\varepsilon_{n,f,(q+1)(s+1)}\big(t'_0,y'_{n,0},h_{\tilde{q}}\big)\|_{n,d}&\leq\|h_{\tilde{q}}\|_{\tilde{q}}\ \delta'_{\tilde{q}}+\big(1+(C_{\tilde{w}}+LC_w)\ \|h_{\tilde{q}}\|_{\tilde{q}}\big)\ \|\varepsilon_{n,f,q(s+1)}\big(t'_0,y'_{n,0},h_{\tilde{q}}\big)\|_{n,d}\\
&\leq\|h_{\tilde{q}}\|_{\tilde{q}}\ \delta'_{\tilde{q}}+e^{(C_{\tilde{w}}+LC_w)\ \|h_{\tilde{q}}\|_{\tilde{q}}}\ \|\varepsilon_{n,f,q(s+1)}\big(t'_0,y'_{n,0},h_{\tilde{q}}\big)\|_{n,d}
\end{aligned}
$$

It is trivial to prove by induction that :

$$
\|\varepsilon_{n,f,(q+1)(s+1)}\big(t'_0,y'_{n,0},h_{\tilde{q}}\big)\|_{n,d}\leq\|h_{\tilde{q}}\|_{\tilde{q}}\ \delta'_{\tilde{q}}\left(\sum_{k=0}^{q}e^{\|h_{\tilde{q}}\|_{\tilde{q}}\ (C_{\tilde{w}}+LC_w)(q-k)}\right)
$$

Thus :

$$
\begin{aligned}
\|\varepsilon_{n,f,(q+1)(s+1)}\big(t'_0,y'_{n,0},h_{\tilde{q}}\big)\|_{n,d}&\leq\|h_{\tilde{q}}\|_{\tilde{q}}\ \delta'_{\tilde{q}}\left(\sum_{k=0}^{\tilde{q}-1}e^{\|h_{\tilde{q}}\|_{\tilde{q}}\ (C_{\tilde{w}}+LC_w)\tilde{q}}\right)\\
&\leq\|h_{\tilde{q}}\|_{\tilde{q}}\ \delta'_{\tilde{q}}\tilde{q}e^{\|h_{\tilde{q}}\|_{\tilde{q}}\ (C_{\tilde{w}}+LC_w)\tilde{q}}
\end{aligned}
$$

Since $\|h_{\tilde{q}}\|_{\tilde{q}} = \underset{\tilde{q}\to+\infty}{\mathcal{O}}\left(\frac{1}{\tilde{q}}\right)$, if we assume $\tilde{q}$ is great enough, there exists $\beta \in \mathbb{R}_+^*$ such that $\|h_{\tilde{q}}\|_{\tilde{q}} \leq \frac{\beta}{\tilde{q}}$. Hence :

$$\sup_{q\in[\![0,\tilde{q}]\!]} \|\varepsilon_{n,f,q(s+1)}\big(t_0, y_{n,0}, h_{\tilde{q}}\big)\|_{n,d} \leq \beta\delta'_{\tilde{q}} e^{\beta(C_{\tilde{w}}+LC_w)}$$

Since $\delta_{\tilde{q}}$ tends towards 0 as $\tilde{q}$ tends towards $+\infty$, the error uniformly converges towards zero.

**III.1.37 Corollary : Convergence of multi-order Runge-Kutta methods**

Let $(n, s) \in \mathbb{N}^* \times \mathbb{N}^*$, $M \in \mathrm{MORK}_{n,s}$, and $\xi$ the approximation times of $M$.
M is convergent up to length $n$ if and only if $\xi_{s+1} = 1_n$.

Proof

A multi-order Runge-Kutta method always satisfies, for all $(j, N) \in [\![1, s+1]\!] \times [\![1, n]\!]$, $\tilde{w}_{N,N,j} = 1$, $M$ thus only needs to satisfy the necessary condition $\xi_{s+1} = 1_n$.

The conditions for convergence are similar to conditions for order of consistency, which hints at a link between the order of consistency of a method and the convergence of its error towards 0.

**III.1.38 Theorem : Order of consistency & Error after the first step**

Let $(n, s, v') \in \mathbb{N}^* \times \mathbb{N}^* \times \mathbb{N}^*$, $M = (\tau, w, \tilde{w}) \in \mathrm{GMORK}_{n,s}$ such that, for all $(j, N) \in [\![1, s]\!] \times [\![1, n]\!]$, $\tilde{w}_{N,N,j} = 1$, and $\xi_{s+1} = 1_n$, where $\xi$ are the approximation times of $M$, let $v_n$ the minimum consistency matrix of $M$, and $\varepsilon_n$ the error function for length $n$ of $M$.

- Let $C_w \in \mathbb{R}_+^*$ such that, for all $(j, N) \in [\![1, s+1]\!] \times [\![1, n]\!]$, $C_w \geq \sum_{j'=1}^{s} \frac{|w_{N,j,j'}|}{N!}$.
- We define $\tilde{v}_n \in \mathbb{Z}^{[\![1,s+1]\!]\times[\![1,n]\!]}$ as, for all $(j, N) \in [\![1, s+1]\!] \times [\![1, n]\!]$, $\tilde{v}_{n,j,N} = \min\big(v_{n,j,N}, v'\big)$.
- Let $\big(f, t_0, y_{n,0}\big)$ an initial value problem of order $n$ with a unique solution and $f$ of class $C^{\max(1,v')}$ globally Lipschitz continuous with a constant $L \in \mathbb{R}_+^*$.
- Let $\tilde{q} \in \overline{\mathbb{N}}^*$, $\Omega'$ a compact subset of $\Omega$, and $(\delta, C) \in \mathbb{R}_+^{*2}$ such that $\delta < \frac{1}{LC_w}$, and for all $(j, N) \in [\![1, s+1]\!] \times [\![1, n]\!]$, $\delta$ and $C$ satisfy the uniform big $\mathcal{O}$ condition of, $M$ is at least of order of consistency $\tilde{v}_{n,j,N}$, on $\Omega'$.
- Let $h \in \mathcal{H}_{\tilde{q},f,t_0,y_{n,0}}$ such that $\|h\|_{\tilde{q}} < \delta$, $\|h\|_{\tilde{q}} \leq 1$, and the time mesh $(t_0, h)$ is contained in $\Omega'$.

1. We define $a \in \mathbb{R}_+^{[\![1,n]\!]}$ as, for all $N \in [\![1, n]\!]$ :

$$a_N = \max\left( \sqrt[N]{\frac{LC_{\tilde{w}}}{1-\delta LC_w} \sum_{j=1}^{s} |w_{N,s+1,j}|}, \sup_{N'\in[\![1,N-1]\!]} \sqrt[N-N']{|\tilde{w}_{N,N',s+1}|\ (N-N')!} \right)$$

If, for all $N \in [\![1, n]\!]$, $\tilde{w}_{N,N,s+1} = 1$, then, for all $q \in [\![0, \tilde{q}-1]\!]$ :

$$\|\varepsilon_{n,f,(q+1)(s+1)}\big(t_0, y_{n,0}, h\big)\|_{n,d} \leq C \exp\left( \|a\|_n \sum_{k'=0}^{q} |h_{k'}| \right) \sum_{k=0}^{q} |h_k|^{\min_{N\in[\![1,n]\!]} N+\tilde{v}_{n,s+1,N}}$$

2. We define $V \in \mathbb{Z}^{[\![1,n]\!]}$ as, for all $N \in [\![1, n]\!]$ :

$$V_N = \min\left( \min_{N'\in[\![1,N]\!]} \big(\tilde{v}_{n,s+1,N'}\big), \min_{N'\in[\![1,n]\!]} \big(N'-1+\tilde{v}_{n,s+1,N'}\big) \right)$$

We have, for all $q \in [\![0, \tilde{q}-1]\!]$ :

$$\|\varepsilon_{n,f,(q+1)(s+1),N}\big(t_0, y_{n,0}, h\big)\|_d \leq C \max_{k\in[\![0,q]\!]} (|h_k|)^{N+V_N} \sum_{k=0}^{q} \prod_{k'=k+1}^{q} C_{\tilde{w}} \left( 1 + \frac{LC_w}{1-\delta LC_w} |h_{k'}| \right)$$

Proof

Both statements come from the same bound. Let $q \in [\![0, \tilde{q}-1]\!]$, $Y_{n,f}$ the approximation function for length $n$ of $M$.

$$\|\varepsilon_{n,f,(q+1)(s+1),N}(t_0, y_{n,0}, h)\|_d = \|\frac{\mathrm{d}^{n-N}\hat{\mathrm{y}}}{\mathrm{d}t^{n-N}}(t_q + h_q) - Y_{n,f,(q+1)(s+1),N}(t_0, y_{n,0}, h)\|_d$$

$$\leq \|\frac{\mathrm{d}^{n-N}\hat{\mathrm{y}}}{\mathrm{d}t^{n-N}}(t_q + h_q) - Y_{n,f,s+1,N}(t_q, \mathcal{J}^{n-1}\hat{\mathrm{y}}(t_q), h_q)\|_d$$

$$+ \|Y_{n,f,s+1,N}(t_q, \mathcal{J}^{n-1}\hat{\mathrm{y}}(t_q), h_q) - Y_{n,f,(q+1)(s+1),N}(t_0, y_{n,0}, h)\|_d$$

$$\leq \|\varepsilon_{n,f,s+1,N}(\overline{\mathcal{J}}^{n-1}\hat{\mathrm{y}}(t_q), h_q)\|_d$$

$$+ \sum_{N'=1}^{N} |\tilde{w}_{N,N',s+1}| \; |h_q|^{N-N'} \; \|Y_{n,f,q(s+1),N'}(t_0, y_{n,0}, h) - \frac{\mathrm{d}^{n-N'}\hat{\mathrm{y}}}{\mathrm{d}t^{n-N'}}(t_q)\|_d$$

$$+ |h_q|^N \sum_{j=1}^{s} \frac{|w_{N,s+1,j}|}{N!} \; \|F_{n,f}(t_q, Y_{n,f,q(s+1)}(t_0, y_{n,0}, h), h_q) - F_{n,f}(t_q, \mathcal{J}^{n-1}\hat{\mathrm{y}}(t_q), h_q)\|_{s,d}$$

To bound the term in $F_{n,f}$, we can use the Lipschitz continuity of $F_{n,f}$. In the proof of Section III.1.35, we see that a sufficient condition for $F_{n,f}$ to be lipschitz is $|h_q| \leq \frac{1-a}{LC_w}$, with $a \in ]0, 1[$. We choose $a = 1 - \delta LC_w$, so that the condition becomes $|h_q| \leq \delta$, which is already verified. To have $0 < a < 1$, it is necessary and sufficient that $\delta < \frac{1}{LC_w}$. The Lipschitz constant given by the proof of Section III.1.35 is hence $\frac{LC_{\tilde{w}}}{1-\delta LC_w}$. We deduce :

$$\|\varepsilon_{n,f,(q+1)(s+1),N}(t_0, y_{n,0}, h)\|_d \leq C \; |h_q|^{N+\tilde{v}_{n,s+1,N}}$$

$$+ \sum_{N'=1}^{N} |\tilde{w}_{N,N',s+1}| \; |h_q|^{N-N'} \; \|\varepsilon_{n,f,q(s+1),N'}(t_0, y_{n,0}, h)\|_d$$

$$+ |h_q|^N \sum_{j=1}^{s} \frac{|w_{N,s+1,j}|}{N!} \frac{L}{1-\delta LC_w} C_{\tilde{w}} \; \|Y_{n,f,q(s+1)}(t_0, y_{n,0}, h) - \mathcal{J}^{n-1}\hat{\mathrm{y}}(t_q)\|_{s,d}$$

$$\leq C \; |h_q|^{N+\tilde{v}_{n,s+1,N}} + \sum_{N'=1}^{N} |\tilde{w}_{N,N',s+1}| \; |h_q|^{N-N'} \; \|\varepsilon_{n,f,q(s+1),N'}(t_0, y_{n,0}, h)\|_d$$

$$+ |h_q|^N \; \frac{LC_{\tilde{w}}}{1-\delta LC_w} \max_{N' \in [\![1,n]\!]} \|\varepsilon_{n,f,q(s+1),N'}(t_0, y_{n,0}, h)\|_d \sum_{j=1}^{s} \frac{|w_{N,s+1,j}|}{N!}$$

1. We now prove the first bound. We have :

$$\|\varepsilon_{n,f,(q+1)(s+1),N}(t_0, y_{n,0}, h)\|_d \leq C \; |h_q|^{N+\tilde{v}_{n,s+1,N}} + \Bigg(1 + \sum_{N'=1}^{N-1} |\tilde{w}_{N,N',s+1}| \; (N-N')! \frac{|h_q|^{N-N'}}{(N-N')!}$$

$$+ \frac{|h_q|^N}{N!} \frac{LC_{\tilde{w}}}{1-\delta LC_w} \sum_{j=1}^{s} |w_{N,s+1,j}| \Bigg) \; \|\varepsilon_{n,f,q(s+1)}(t_0, y_{n,0}, h)\|_{n,d}$$

$$
\begin{aligned}
&\leq C \ |h_q|^{N+\tilde{v}_{n,s+1,N}} + \left( 1 + \sum_{N'=1}^{N-1} \frac{\left( \sqrt[N-N']{|\tilde{w}_{N,N',s+1}| \ (N-N')!} \ |h_q| \right)^{N-N'}}{(N-N')!} \right. \\
&\quad \left. + \frac{\left( \sqrt[N]{\frac{LC_{\tilde{w}}}{1-\delta LC_w} \sum_{j=1}^{s} |w_{N,s+1,j}|} \ |h_q| \right)^{N}}{N!} \right) \|\varepsilon_{n,f,q(s+1)}(t_0, y_{n,0}, h)\|_{n,d} \\
&\leq C \ |h_q|^{N+\tilde{v}_{n,s+1,N}} + \left( \sum_{N'=0}^{N} \frac{\left(a_N \ |h_q|\right)^{N'}}{N'!} \right) \|\varepsilon_{n,f,q(s+1)}(t_0, y_{n,0}, h)\|_{n,d} \\
&\leq C \ |h_q|^{N+\tilde{v}_{n,s+1,N}} + \exp\left(a_N \ |h_q|\right) \ \|\varepsilon_{n,f,q(s+1)}(t_0, y_{n,0}, h)\|_{n,d}
\end{aligned}
$$

Hence :

$$
\|\varepsilon_{n,f,(q+1)(s+1)}(t_0, y_{n,0}, h)\|_{n,d} \leq C \ |h_q|^{\min\limits_{N \in [\![1,n]\!]}(N+\tilde{v}_{n,s+1,N})} + \exp\left(\|a\|_n \ |h_q|\right) \ \|\varepsilon_{n,f,q(s+1)}(t_0, y_{n,0}, h)\|_{n,d}
$$

We now prove by induction that, for all $q \in [\![0, \tilde{q}-1]\!]$ :

$$
\|\varepsilon_{n,f,(q+1)(s+1)}(t_0, y_{n,0}, h)\|_{n,d} \leq C \sum_{k=0}^{q} |h_k|^{\min\limits_{N \in [\![1,n]\!]}(N+\tilde{v}_{n,s+1,N})} \exp\left( \|a\|_n \sum_{k'=1+k}^{q} |h_{k'}| \right)
$$

The base case $q = 0$ is trivial. Assume this property is true for $q \in [\![0, \tilde{q}-2]\!]$. We have ;

$$
\begin{aligned}
&\|\varepsilon_{n,f,(q+2)(s+1)}(t_0, y_{n,0}, h)\|_{n,d} \\
&\leq C \ |h_{q+1}|^{\min\limits_{N \in [\![1,n]\!]}(N+\tilde{v}_{n,s+1,N})} + \exp\left(\|a\|_n \ |h_{q+1}|\right) \ \|\varepsilon_{n,f,(q+1)(s+1)}(t_0, y_{n,0}, h)\|_{n,d} \\
&\leq C \ |h_{q+1}|^{\min\limits_{N \in [\![1,n]\!]}(N+\tilde{v}_{n,s+1,N})} + \exp\left(\|a\|_n \ |h_{q+1}|\right) C \sum_{k=0}^{q} |h_k|^{\min\limits_{N \in [\![1,n]\!]}(N+\tilde{v}_{n,s+1,N})} \exp\left( \|a\|_n \sum_{k'=1+k}^{q} |h_{k'}| \right) \\
&\leq C \sum_{k=0}^{q+1} |h_k|^{\min\limits_{N \in [\![1,n]\!]}(N+\tilde{v}_{n,s+1,N})} \exp\left( \|a\|_n \sum_{k'=1+k}^{q+1} |h_{k'}| \right)
\end{aligned}
$$

This property is hence true for $q+1$. We can simplify this bound by making it coarser :

$$
\|\varepsilon_{n,f,(q+1)(s+1)}(t_0, y_{n,0}, h)\|_{n,d} \leq C \exp\left( \|a\|_n \sum_{k'=0}^{q} |h_{k'}| \right) \sum_{k=0}^{q} |h_k|^{\min\limits_{N \in [\![1,n]\!]}(N+\tilde{v}_{n,s+1,N})}
$$

2. We now prove by induction that, for all $q \in [\![0, \tilde{q}-1]\!]$ :

$$
\|\varepsilon_{n,f,(q+1)(s+1),N}(t_0, y_{n,0}, h)\|_d \leq C \max_{k \in [\![0,q]\!]} \left(|h_k|\right)^{N+V_N} \sum_{k=0}^{q} \prod_{k'=k+1}^{q} C_{\tilde{w}} \left( 1 + \frac{LC_w}{1-\delta LC_w} \ |h_{k'}| \right)
$$

The base case is trivial. Assume this property is true for $q \in [\![0, \tilde{q}-2]\!]$. We have :

$$
\begin{aligned}
&\|\varepsilon_{n,f,(q+2)(s+1),N}\big(t_0,y_{n,0},h\big)\|_d \le C\ |h_{q+1}|^{N+\tilde{v}_{n,s+1,N}} \\
&\quad +C_{\tilde{w}} \max_{N'\in[\![1,N]\!]}\Big(|h_{q+1}|^{N-N'}\ \|\varepsilon_{n,f,(q+1)(s+1),N'}\big(t_0,y_{n,0},h\big)\|_d\Big) \\
&\quad +\ |h_{q+1}|^{N}\ \frac{LC_{\tilde{w}}}{1-\delta LC_w} \max_{N'\in[\![1,n]\!]}\|\varepsilon_{n,f,(q+1)(s+1),N'}\big(t_0,y_{n,0},h\big)\|_d \sum_{j=1}^{s}\frac{|w_{N,s+1,j}|}{N!} \\
&\le C \max_{k\in[\![0,q+1]\!]}\big(|h_k|\big)^{N+\tilde{v}_{n,s+1,N}} + C_{\tilde{w}} \max_{N'\in[\![1,N]\!]}\Bigg[ \\
&\quad \max_{k\in[\![0,q+1]\!]}\big(|h_k|\big)^{N-N'} C \max_{k\in[\![0,q]\!]}\big(|h_k|\big)^{N'+V_{N'}} \sum_{k=0}^{q}\prod_{k'=k+1}^{q} C_{\tilde{w}}\bigg(1+\frac{LC_w}{1-\delta LC_w}\ |h_{k'}|\bigg)\Bigg] \\
&\quad +\ |h_{q+1}| \max_{k\in[\![0,q+1]\!]}\big(|h_k|\big)^{N-1} \frac{LC_w C_{\tilde{w}}}{1-\delta LC_w} \max_{N'\in[\![1,n]\!]}\Bigg[ \\
&\quad C \max_{k\in[\![0,q]\!]}\big(|h_k|\big)^{N'+V_{N'}} \sum_{k=0}^{q}\prod_{k'=k+1}^{q} C_{\tilde{w}}\bigg(1+\frac{LC_w}{1-\delta LC_w}\ |h_{k'}|\bigg)\Bigg] \\
&\le C \max_{k\in[\![0,q+1]\!]}\big(|h_k|\big)^{N}\Bigg[\max_{k\in[\![0,q+1]\!]}\big(|h_k|\big)^{V_N} \\
&\quad +C_{\tilde{w}} \max_{k\in[\![0,q+1]\!]}\big(|h_k|\big)^{\min_{N'\in[\![1,N]\!]} V_{N'}} \sum_{k=0}^{q}\prod_{k'=k+1}^{q} C_{\tilde{w}}\bigg(1+\frac{LC_w}{1-\delta LC_w}\ |h_{k'}|\bigg) \\
&\quad +\ |h_{q+1}|\ \frac{LC_w C_{\tilde{w}}}{1-\delta LC_w} \max_{k\in[\![0,q+1]\!]}\big(|h_k|\big)^{\min_{N'\in[\![1,n]\!]} N'-1+V_{N'}} \sum_{k=0}^{q}\prod_{k'=k+1}^{q} C_{\tilde{w}}\bigg(1+\frac{LC_w}{1-\delta LC_w}\ |h_{k'}|\bigg)\Bigg] \\
&\le C \max_{k\in[\![0,q+1]\!]}\big(|h_k|\big)^{N+\min\left(\min_{N'\in[\![1,N]\!]} V_{N'},\,\min_{N'\in[\![1,n]\!]} N'-1+V_{N'}\right)}\Bigg[ \\
&\quad 1+C_{\tilde{w}}\bigg(1+\frac{LC_w}{1-\delta LC_w}\ |h_{q+1}|\bigg)\sum_{k=0}^{q}\prod_{k'=k+1}^{q} C_{\tilde{w}}\bigg(1+\frac{LC_w}{1-\delta LC_w}\ |h_{k'}|\bigg)\Bigg] \\
&\le C \max_{k\in[\![0,q+1]\!]}\big(|h_k|\big)^{N+V_N} \sum_{k=0}^{q+1}\prod_{k'=k+1}^{q+1} C_{\tilde{w}}\bigg(1+\frac{LC_w}{1-\delta LC_w}\ |h_{k'}|\bigg)
\end{aligned}
$$

The property is hence true for $q+1$.

We can observe that if $M$ is a multi-order Runge-Kutta method, then for the first case, we have :

$$
a_1 = \frac{LC_{\tilde{w}}}{1-\delta LC_w}\sum_{j=1}^{s}|w_{1,s+1,j}|,\ \ \forall N\in[\![2,n]\!],\ a_N = \max\left(\sqrt[N]{\frac{LC_{\tilde{w}}}{1-\delta LC_w}\sum_{j=1}^{s}|w_{N,s+1,j}|},|\xi_{j,N}|\right)
$$

In the next section we define a category of methods called rewrite-compliant methods, and if $M$ is rewrite-compliant, we have $\min_{N\in[\![1,n]\!]} N+\tilde{v}_{n,s+1,N} = 1+\tilde{v}_{n,s+1,1}$, hence, for the second case :

$$
\forall N\in[\![1,n]\!],\ V_N = \min_{N'\in[\![1,N]\!]}\big(\tilde{v}_{n,s+1,N'}\big)
$$

The first bound of the theorem shows that the maximum error of the approximations is governed by the sum of the step sizes raised to the power of $\min_{N\in[\![1,n]\!]} N + \tilde{v}_{n,s+1,N}$. In particular, this exponent is bounded above by $1 + \tilde{v}_{n,s+1,1}$. In the next section we prove that this quantity cannot exceed the corresponding exponent for a classical Runge-Kutta method. Hence, in general, the global error cannot decay faster than that of a Runge-Kutta method with the same structure.

One advantage of multi-order methods, however, lies in the constans appearing in this bound. The additional freedom allows these constants to be reduced. This is reflected in the second bound of the theorem, which shows that the error at a given rank is governed by its own order of consistency rather than by the minimal order.

A well known advantage of this type of method is that when the initial value problem does not depend on the highest derivatives of the solution, many order conditions disappear. This leads to effective orders of consistency and convergence higher than those expected in the general case. Nevertheless, this gain in order remains bounded by $2s'$ since, even in the extreme case where the differential equation function is independent of all derivatives of the solution, the solved system conditions still need to be satisfied. As established in Section II.2.58, such initial value problem are also less computationally demanding.

We end this subsection with a generalized version of the solved system conditions.

**III.1.39 Theorem : Generalized solved system conditions**

Let $(n, s, v) \in \mathbb{N}^* \times \mathbb{N}^* \times \mathbb{N}^*$, $M = (\tau, w, \xi) \in \mathrm{MORK}_{n,s}$, $v' \in \mathbb{N}^*$, and $(j, N) \in [\![1, s+1]\!] \times [\![1, n]\!]$. If $M$ is at least of order of consistency $v'$ for length $n$ at stage $j$ and rank $N$, then :

$$\forall k \in [\![0, v'-1]\!],\ \forall \lambda \in \mathbb{N}^{[\![1,n]\!]}$$

$$|\lambda| \leq v' - 1 - k \Rightarrow \sum_{j'=1}^{s} w_{N,j,j'} \tau_{j'}^{k} \left( \prod_{N=1}^{n} \xi_{j',N}^{\lambda_N} \right) = \xi_{j,N}^{N+k+\,|\lambda|} \frac{N!(k+|\lambda|)!}{(k+|\lambda|+N)!}$$

Proof

The proof is in Section D.

We know how well a method $M$ of length $n \in \mathbb{N}^*$ approximates the solution of an initial value problem of order $n$, we are however not able to solve initial value problems of a different order using $M$. We solve this issue in the next subsection.

## III.2 Rewriting an initial value problem

Since the order of the initial value problem can be greater than 1, there are 2 cases to consider when given an initial value problem of order $n \in \mathbb{N}^*$ and a method of length $\tilde{n} \in \mathbb{N}^*$, $\tilde{n} > n$ or $\tilde{n} < n$. The first case is necessary to solve initial value problem of lower order, while the latter has downsides similar to that of the rewriting for Runge-Kutta methods. We start with the most interesting case of the two.

### The case $n = \tilde{n}$

Yes, this case is of interest because it highlights an interesting issue that needs to be addressed before tackling the two other cases. Now that there are multiple derivatives, we can take the $n - N^{\text{th}}$ derivative of $\hat{y}$, with $N \in [\![1, n]\!]$, and express it as a derivative of another function $\hat{y}'$, with $N' \in [\![1, n]\!]$ :

$$\frac{\mathrm{d}^{n-N'}\hat{y}'}{\mathrm{d}t^{n-N'}} = \frac{\mathrm{d}^{n-N}\hat{y}}{\mathrm{d}t^{n-N}}$$

With the right initial values, we get :

$$\frac{\mathrm{d}^{n}\hat{y}'}{\mathrm{d}t^{n}} = \frac{\mathrm{d}^{n+N'-N}\hat{y}}{\mathrm{d}t^{n+N'-N}}, \quad \frac{\mathrm{d}^{N-N'}\hat{y}'}{\mathrm{d}t^{N-N'}} = \hat{y}$$

We do not want the $n^{\text{th}}$ derivative of $\hat{y}'$ to be a derivative of $\hat{y}$ of order higher than $n$, hence $N' \leq N$. The case $N' = N$ is not interesting, we take $N' \neq N$. We deduce $N \neq 1$ and $N' \in [\![1, N-1]\!]$.

**III.2.1 Definition : Rewriting for initial value problems of equal order**

Let $(f, t_0, y_{n,0})$ an initial value problem of order $n \in \mathbb{N}^*$, $N \in [\![2, n]\!]$, $N' \in [\![1, N-1]\!]$. We define :

$$\mathbf{U}' = \left\{ \left(t, \begin{bmatrix} x'_1 \\ x_1 \end{bmatrix}, \dots, \begin{bmatrix} x'_n \\ x_n \end{bmatrix}\right) \in \mathbb{R} \times \mathbb{R}^{[\![1,n]\!]\times[\![1,2d]\!]} \mid (t, x_1, \dots x_{N-1}, x'_{N'}, x_{N+1}, \dots x_n) \in \mathbf{U} \right\}$$

$$\forall \left(t, \begin{bmatrix} x'_1 \\ x_1 \end{bmatrix}, \dots, \begin{bmatrix} x'_n \\ x_n \end{bmatrix}\right) \in \mathbf{U}',\ f'\left(t, \begin{bmatrix} x'_1 \\ x_1 \end{bmatrix}, \dots, \begin{bmatrix} x'_n \\ x_n \end{bmatrix}\right) = \begin{bmatrix} x_{N-N'} \\ f(t, x_1, \dots x_{N-1}, x'_{N'}, x_{N+1}, \dots x_n) \end{bmatrix}$$

$$\forall N'' \in [\![1, n]\!],\ y''_{n,0,N''} = \begin{bmatrix} y'_{n,0,N''} \\ y_{n,0,N''} \end{bmatrix}$$

$$\forall N'' \in [\![1, n]\!],\ y'_{n,0,N''} = \begin{cases} y_{n,0,N''+N-N'} & \text{if } N'' \in [\![1, N'-N+n]\!] \\ 0_d & \text{otherwise} \end{cases}$$

We define $(f', t_0, y''_{n,0})$ as the overwriting of $N$ with $N'$ of $(f, t_0, y_{n,0})$.

**III.2.2 Proposition : Solutions of rewritten initial value problems - Equal order rewriting**

Let $(f, t_0, y_{n,0})$ an initial value problem of order $n \in \mathbb{N}^*$, $N \in [\![2, n]\!]$, $N' \in [\![1, N-1]\!]$, $(f', t_0, y''_{n,0})$ the overwriting of $N$ with $N'$ of $(f, t_0, y_{n,0})$

- If $\hat{y}$ is a solution of $(f, t_0, y_{n,0})$, then, if we define $\hat{y}'$ the $(N - N')^{\text{th}}$ primitive integral of $\hat{y}$ which satisfies $\mathcal{J}^{N-N'-1}\hat{y}'(t_0) = 0_{N-N',d}$, $\begin{bmatrix} \hat{y}' \\ \hat{y} \end{bmatrix}$ is a solution of $(f', t_0, y''_{n,0})$.
- If $\hat{y}''$ is a solution of $(f', t_0, y''_{n,0})$ then there exists a unique $\hat{y}'$ such that, if we define $\hat{y} = \frac{\mathrm{d}^{N-N'}\hat{y}'}{\mathrm{d}t^{N-N'}}$, then $\hat{y}'' = \begin{bmatrix} \hat{y}' \\ \hat{y} \end{bmatrix}$. $\hat{y}$ is a solution of $(f, t_0, y_{n,0})$

Proof

Trivial.

This rewriting is named overwriting because it is essentially it does. We can then use a method on the rewritten initial value problem.

**III.2.3 Definition : Overwriting of a method**

Let $(n, s) \in \mathbb{N}^* \times \mathbb{N}^*$, $M = (\tau, w, \tilde{w}) \in \mathrm{GMORK}_{n,s}$, $N \in [\![2, n]\!]$, and $N' \in [\![1, N-1]\!]$.
We define $M' = (\tau, w', \tilde{w}') \in \mathrm{GMORK}_{n,s}$ the overwriting of $N$ with $N'$ of $M$ as :

$$\forall N'' \in [\![1, n]\!] \setminus \{N\},\ w'_{N''} = w_{N''}$$

$$\forall N'' \in [\![1, n]\!] \setminus \{N\},\ \tilde{w}'_{N''} = \tilde{w}_{N''}$$

$$\forall (j, j') \in [\![1, s+1]\!] \times [\![1, s]\!],\ w'_{N,j,j'} = N! \sum_{j''=1}^{s} \frac{w_{N',j,j''}}{N'!} \frac{w_{N-N',j'',j'}}{(N-N')!}$$

$$\forall (N'', j) \in [\![1, N]\!] \times [\![1, s+1]\!]$$

$$\tilde{w}'_{N,N'',j} = \begin{cases} \sum_{j'=1}^{s} \frac{w_{N',j,j'}}{N'!} \tilde{w}_{N-N',N'',j'} & \text{if } N'' \in [\![1, N-N']\!] \\ \tilde{w}_{N',N''+N'-N,j} & \text{otherwise} \end{cases}$$

**III.2.4 Proposition : Form of rewritten method - Equal order rewriting**

Let $(n, s) \in \mathbb{N}^* \times \mathbb{N}^*$, $M \in \mathrm{GMORK}_{n,s}$, $N \in [\![2, n]\!]$, $N' \in [\![1, N-1]\!]$, $M' \in \mathrm{GMORK}_{n,s}$ the overwriting of $N$ with $N'$ of $M$, $f$ a differential equation function of order $n$, $(t, y_{n,0}, h) \in \mathbf{U} \times \mathbb{R}$, and $(f', t, y''_{n,0})$ the overwriting of $N$ with $N'$ of $(f, t, y_{n,0})$.
The stage system $(M, f', t, y''_{n,0}, h)$ is equivalent to the stage system $(M', f, t, y_{n,0}, h)$.

Proof

The stage system $(M, f', t, y''_{n,0}, h)$ is :

$$\forall (j, N'') \in [\![1, s+1]\!] \times [\![1, n]\!]$$

$$y''_{n,j,N''} = \sum_{N'''=1}^{N''} \tilde{w}_{N'',N''',j} h^{N''-N'''} y''_{n,0,N'''} + \frac{h^{N''}}{N''!} \sum_{j'=1}^{s} w_{N'',j,j'} f'\left(t + \tau_{j'} h, y''_{n,j'}\right)$$

$$y''_{n,s+1} \in \mathbb{R}^{[\![1,n]\!] \times [\![1,2d]\!]},\ \forall j \in [\![1, s]\!],\ \left(t + \tau_j h, y''_{n,j}\right) \in \mathbf{U}'$$

We define, for all $(j, N'') \in [\![1, s+1]\!] \times [\![1, n]\!]$, $y''_{n,j,N''} = \begin{bmatrix} y'_{n,j,N''} \\ y_{n,j,N''} \end{bmatrix}$. Therefore :

$$\forall (j, N'') \in [\![1, s+1]\!] \times [\![1, n]\!]$$

$$y'_{n,j,N''} = \sum_{N'''=1}^{N''} \tilde{w}_{N'',N''',j} h^{N''-N'''} y'_{n,0,N'''} + \frac{h^{N''}}{N''!} \sum_{j'=1}^{s} w_{N'',j,j'} y_{n,j',N-N'}$$

$$y_{n,j,N''} = \sum_{N'''=1}^{N''} \tilde{w}_{N'',N''',j} h^{N''-N'''} y_{n,0,N'''}$$

$$+ \frac{h^{N''}}{N''!} \sum_{j'=1}^{s} w_{N'',j,j'} f\left(t + \tau_{j'} h, y_{n,j',1}, \dots\ y'_{n,j',N'}, \dots\ y_{n,j',n}\right)$$

$$y''_{n,s+1} \in \mathbb{R}^{[\![1,n]\!] \times [\![1,2d]\!]},\ \forall j \in [\![1, s]\!],\ \left(t + \tau_j h, y_{n,j,1}, \dots\ y'_{n,j,N'}, \dots\ y_{n,j,n}\right) \in \mathbf{U}$$

We are only interested in $y'_{n,j',N'}$ with $j \in [\![1, s+1]\!]$, $y_{n,j',N''}$ with $(j, N'') \in [\![1, s+1]\!] \times [\![1, n]\!] \setminus \{N\}$. The other variables are functions of the ones we just mentioned, we thus remove them from the system.

$$\forall j \in [\![1, s+1]\!]$$

$$y'_{n,j,N'} = \sum_{N''=1}^{N'} \tilde{w}_{N',N'',j} h^{N'-N''} y'_{n,0,N''} + \frac{h^{N'}}{N'!} \sum_{j'=1}^{s} w_{N',j,j'} \left[ \sum_{N''=1}^{N-N'} \tilde{w}_{N-N',N'',j'} h^{N-N'-N''} y_{n,0,N''} \right.$$

$$\left. + \frac{h^{N-N'}}{(N-N')!} \sum_{j''=1}^{s} w_{N-N',j',j''} f\left(t + \tau_{j''} h, y_{n,j'',1}, \dots\, y'_{n,j'',N}, \dots\, y_{n,j'',n}\right) \right]$$

$$\forall (j, N'') \in [\![1, s+1]\!] \times [\![1, n]\!] \setminus \{N\},$$

$$y_{n,j,N''} = \sum_{N'''=1}^{N''} \tilde{w}_{N'',N''',j} h^{N''-N'''} y_{n,0,N'''}$$

$$+ \frac{h^{N''}}{N''!} \sum_{j'=1}^{s} w_{N'',j,j'} f\left(t + \tau_{j'} h, y_{n,j',1}, \dots\, y'_{n,j',N}, \dots\, y_{n,j',n}\right)$$

$$y'_{n,s+1,N'} \in \mathbb{R}^{[\![1,d]\!]},\ \forall N'' \in [\![1, n]\!] \setminus \{N\},\ y_{n,s+1} \in \mathbb{R}^{[\![1,n]\!] \times [\![1,d]\!]},$$

$$\forall j \in [\![1, s]\!],\ \left(t + \tau_j h, y_{n,j,1}, \dots\, y'_{n,j,N'}, \dots\, y_{n,j,n}\right) \in \mathbf{U}$$

Let $j \in [\![1, s+1]\!]$. We now find the expression of $y'_{n,j,N'}$ :

$$y'_{n,j,N'} = \sum_{N''=1}^{N'} \tilde{w}_{N',N'',j} h^{N'-N''} y_{n,0,N-N'+N''} + \sum_{N''=1}^{N-N'} h^{N-N''} \sum_{j'=1}^{s} \frac{w_{N',j,j'}}{N'!} \tilde{w}_{N-N',N'',j'} y_{n,0,N''}$$

$$+ h^N \sum_{j'=1}^{s} \sum_{j''=1}^{s} \frac{w_{N',j,j'}}{N'!} \frac{w_{N-N',j',j''}}{(N-N')!} f\left(t + \tau_{j''} h, y_{n,j'',1}, \dots\, y'_{n,j'',N}, \dots\, y_{n,j'',n}\right)$$

$$= \sum_{N''=N-N'+1}^{N} \tilde{w}_{N',N''+N'-N,j} h^{N-N''} y_{n,0,N''} + \sum_{N''=1}^{N-N'} h^{N-N''} \sum_{j'=1}^{s} \frac{w_{N',j,j'}}{N'!} \tilde{w}_{N-N',N'',j'} y_{n,0,N''}$$

$$+ \frac{h^N}{N!} \sum_{j'=1}^{s} \left( N! \sum_{j''=1}^{s} \frac{w_{N',j,j''}}{N'!} \frac{w_{N-N',j'',j'}}{(N-N')!} \right) f\left(t + \tau_{j'} h, y_{n,j',1}, \dots\, y'_{n,j',N'}, \dots\, y_{n,j',n}\right)$$

We can see that it yields the stage system $\left(M', f, t, y_{n,0}, h\right)$.

We study the effects of this rewriting on the accuracy of a method.

**III.2.5 Proposition : Order of consistency - Equal order rewriting**

Let $(n, s) \in \mathbb{N}^* \times \mathbb{N}^*$, $M \in \mathrm{GMORK}_{n,s}$, $N \in [\![2, n]\!]$, $N' \in [\![1, N-1]\!]$, $M'$ the overwriting of $N$ with $N'$ of $M$, $\xi, \xi'$ the approximation times of respectively $M, M'$, and $v_n, v'_n$ the minimum consistency matrices of respectively $M, M'$.

If $M$ is convergent up to length $n$ then $M'$ is convergent up to length $n$. If, for all $(j, N'') \in [\![1, s]\!] \times [\![1, n]\!]$, $\tilde{w}_{N'',N'',j} = 1$, then :

$$\forall (j, N'') \in [\![1, s+1]\!] \times ([\![1, n]\!] \setminus \{N\}),\ \xi'_{j,N''} = \xi_{j,N''},\ v'_{n,j,N''} \geq v_{n,j,N''}$$

$$\forall j \in [\![1, s+1]\!],\ \xi'_{j,N} = \xi_{j,N'},\ v'_{n,j,N} \geq v_{n,j,N'} + N' - N$$

Proof

If $N' \neq 1$ the approximation times do not change. For $N' = 1$ we have :

$$\forall j \in [\![1, s+1]\!],\ \xi'_{j,N} = \tilde{w}'_{N,N-1,j} = \sum_{j'=1}^{s} w_{1,j,j'} \tilde{w}_{N-1,N-1,j'} = \sum_{j'=1}^{s} w_{1,j,j'} = \xi_{j,1}$$

$M$ and $M'$ have the same approximation times, and the rewriting is a special case of the initial value problems considered in the definition of the order of consistency, so the order of consistency of $M$ is at least that of $M'$.

With a similar reasoning we conclude that if $M$ is convergent, then $M$ also is. It follows that, if we assume the necessary conditions on the initial value problem for order of consistency $v_{n,j,N'}$, then :

$$\frac{\mathrm{d}^{n-N'}\hat{\mathrm{y}}'}{\mathrm{d}t^{n-N'}}\left(t' + \xi'_{N,j}h\right) - y'_{n,j,N'} = \underset{\substack{h \to 0 \\ t' \in \Omega'}}{\mathcal{O}}\left(h^{v_{n,j,N'}+N'}\right)$$

Therefore :

$$\frac{\mathrm{d}^{n-N}\hat{\mathrm{y}}}{\mathrm{d}t^{n-N}}\left(t_0 + \xi'_{N,j}h\right) - y'_{n,j,N'} = \underset{\substack{h \to 0 \\ t' \in \Omega'}}{\mathcal{O}}\left(h^{\left(v_{n,j,N'}+N'-N\right)+N}\right)$$

Thus $v'_{n,j,N} \geq v_{n,j,N'} + N' - N$. At the rest of the ranks and stage, the order of consistency is at least that of $M$.

The issue is that, if $v_{n,j,N'} + N' - N > v_{n,j,N}$, then $v'_{n,j,N} > v_{n,j,N}$, which means that $M'$ has a higher order of consistency for length $n$ at stage $j$ and rank $N$ than $M$. We define a category of method that cannot be outperfomed using this rewriting. Such a method must satisfy :

$$\forall (j, N) \in [\![1, s+1]\!] \times [\![2, n]\!],\ \forall N' \in [\![1, N-1]\!],\ v_{n,j,N} - v_{n,j,N'} \geq N' - N$$

We find a nicer equivalent condition.

$$\begin{aligned}
&\forall (j, N) \in [\![1, s+1]\!] \times [\![2, n]\!],\ \forall N' \in [\![1, N-1]\!],\ v_{n,j,N} - v_{n,j,N'} \geq N' - N \\
&\quad\Rightarrow \forall (j, N') \in [\![1, s+1]\!] \times [\![1, n-1]\!],\ v_{n,j,N'+1} - v_{n,j,N'} \geq -1 \\
&\Rightarrow \forall (j, N) \in [\![1, s+1]\!] \times [\![2, n]\!],\ \forall N' \in [\![1, N-1]\!],\ \sum_{N''=N'}^{N-1} \left(v_{n,j,N''+1} - v_{n,j,N''}\right) \geq -\sum_{N''=N'}^{N-1} 1 \\
&\Rightarrow \forall (j, N) \in [\![1, s+1]\!] \times [\![2, n]\!],\ \forall N' \in [\![1, N-1]\!],\ v_{n,j,N} - v_{n,j,N'} \geq N' - N
\end{aligned}$$

The previous condition is therefore equivalent to :

$$\forall (j, N) \in [\![1, s+1]\!] \times [\![1, n-1]\!],\ v_{n,j,N+1} - v_{n,j,N} \geq -1$$

This category of method is important for the two other cases $n < \tilde{n},\ n > \tilde{n}$, a method which satisfy this condition uses the other two rewritings without losing in order of consistency.

**III.2.6 Definition : Rewrite-compliant methods**

Let $(n, s) \in \mathbb{N}^* \times \mathbb{N}^*$, $M \in \mathrm{GMORK}_{n,s}$, and $v_n$ the minimum consistency matrix of $M$.
$M$ is said to be rewrite-compliant if and only if :

$$\forall (j, N) \in [\![1, s+1]\!] \times [\![1, n-1]\!],\ v_{n,j,N+1} - v_{n,j,N} \geq -1$$

This category of method cannot be outperformed in terms of order of consistency using the rewriting we have just studied, but the rewritten method may be more costly or lose some property in the process. For instance, we may ask ourselves how this rewriting affects the maximum weight digraph.

**III.2.7 Proposition : Maximum weight digraph - Equal order rewriting**

Let $(n,s) \in \mathbb{N}^* \times \mathbb{N}^*$, $M \in \mathrm{GMORK}_{n,s}$, $N \in [\![2,n]\!]$, $N' \in [\![1,N-1]\!]$, $M'$ the overwriting of $N$ with $N'$ of $M$, $\mathcal{R}_r, {\mathcal{R}_r}'$ the reachability relations of respectively $M, M'$.

1. For all $(j,j') \in [\![1,s+1]\!] \times [\![1,s]\!]$, $j' \; {\mathcal{R}_r}' \; j \Rightarrow j' \; \mathcal{R}_r \; j$.
2. Any useless stage, reachable-closed set, explicit stage, block, explicit block, explicit partition of the stages, order of computation, explicit rank, parallel blocks, independent set of stages, partition of independent sets of stages of $M$ is also one of $M'$. If $M$ is explicit $M'$ also is.

Proof

1. Let $\mathcal{G}_n = ([\![1,s+1]\!], \mathcal{A}_n)$, $\mathcal{G}'_n = ([\![1,s+1]\!], \mathcal{A}'_n)$ the weight graphs of respectively $M$, $M'$. We define :

$$A_1 = \left\{(j',j) \in [\![1,s]\!] \times [\![1,s+1]\!] \mid w_{N,j,j'} \neq 0\right\}$$
$$A_2 = \left\{(j',j) \in [\![1,s]\!] \times [\![1,s+1]\!] \mid \exists N'' \in [\![1,n]\!] \setminus \{N\},\ w_{N'',j,j'} \neq 0\right\}$$

   Thus $\mathcal{A}_n = A_1 \cup A_2$. The only changes between $M$ and $M'$ are the weights of $N$, so only $A_1$ changes and becomes :

$$A'_1 = \left\{(j',j) \in [\![1,s]\!] \times [\![1,s+1]\!] \mid w'_{N,j,j'} \neq 0\right\}$$

   We thus have $\mathcal{A}'_n = A'_1 \cup A_2$. Since both maximum weight digraph share $A_2$, only $A_1, A'_1$ are of interest. Let $(j,j') \in [\![1,s+1]\!] \times [\![1,s]\!]$, $A'_1$ has an arc from $j'$ to $j$ if and only if $w'_{N,j,j'} \neq 0$, which is equivalent to :

$$\sum_{j''=1}^{s} \frac{w_{N',j,j''}}{N'!} \frac{w_{N-N',j'',j'}}{(N-N')!} \neq 0$$

   This implies :

$$\exists j'' \in [\![1,s]\!],\ \frac{w_{N',j,j''}}{N'!} \frac{w_{N-N',j'',j'}}{(N-N')!} \neq 0$$

   We deduce that if $A'_1$ has an arc from $j'$ to $j$, then there exists $j'' \in [\![1,s]\!]$ such that $(j',j'',j)$ is a diwalk of $\mathcal{G}_n$. Hence $j' \; {\mathcal{R}_r}' \; j \Rightarrow j' \; \mathcal{R}_r \; j$.
2. Those are direct consequences of 1.

## The case $n < \tilde{n}$

To deal with this case we do not use the fact that $\hat{y}$ is vector-valued as it would come with similar downsides to what we have seen with rewritten Runge-Kutta methods and what will see in the next section. Instead, we treat $\hat{y}$ as a derivative of another function $\tilde{y}$, making the initial value problem of $\hat{y}$ artificially of greater order. In equation this can be written as :

$$\frac{\mathrm{d}^{\tilde{n}-n}\hat{y}'}{\mathrm{d}t^{\tilde{n}-n}} = \hat{y},\quad \frac{\mathrm{d}^{\tilde{n}}\hat{y}'}{\mathrm{d}t^{\tilde{n}}} = \frac{\mathrm{d}^{n}\hat{y}}{\mathrm{d}t^{n}},\ \forall N \in [\![n+1, \tilde{n}]\!],\ \frac{\mathrm{d}^{\tilde{n}-N}\hat{y}'}{\mathrm{d}t^{\tilde{n}-N}}(t_0) = 0$$

The additional initial values have been chosen to be zero, but their values do not change the outcome.

**III.2.8 Definition : Rewriting for initial value problems of lower order**

Let $(n, \tilde{n}) \in \mathbb{N}^* \times \mathbb{N}^*$ such that $n < \tilde{n}$, and $(f, t_0, y_{n,0})$ an initial value problem of order $n$. We define :

$$\begin{aligned} \mathbf{U}' &= \mathbf{U} \times \mathbb{R}^{[\![n+1,\tilde{n}]\!]} \\ \forall (t, x) \in \mathbf{U}',\ f'(t, x_1, ..., x_{\tilde{n}}) &= f(t, x_1, ..., x_n) \\ \forall N \in [\![1, n]\!],\ y_{\tilde{n},0,N} &= y_{n,0,N} \\ \forall N \in [\![n+1, \tilde{n}]\!],\ y_{\tilde{n},0,N} &= 0_d \end{aligned}$$

We define $(f', t_0, y_{\tilde{n},0})$ as the prolongation to order $\tilde{n}$ of $(f, t_0, y_{n,0})$.

**III.2.9 Proposition : Solutions of rewritten initial value problems - Lower order rewriting**

Let $(n, \tilde{n}) \in \mathbb{N}^* \times \mathbb{N}^*$ such that $n < \tilde{n}$, $(f, t_0, y_{n,0})$ an initial value problem of order $n$, and $(f', t_0, y_{\tilde{n},0})$ the prolongation to order $\tilde{n}$ of $(f, t_0, y_{n,0})$.

- If $\hat{y}$ is a solution of $(f, t_0, y_{n,0})$, then, if we define $\hat{y}'$ the $(\tilde{n} - n)^{\text{th}}$ primitive integral of $\hat{y}$ which satisfies $\mathcal{J}^{n-1}\hat{y}'(t_0) = 0_{n,d}$, then $\hat{y}'$ is a solution of $(f', t_0, y_{\tilde{n},0})$.
- If $\hat{y}'$ is a solution of $(f', t_0, y_{\tilde{n},0})$ then $\frac{\mathrm{d}^{\tilde{n}-n}\hat{y}'}{\mathrm{d}t^{\tilde{n}-n}}$ is a solution of $(f, t_0, y_{n,0})$.

Proof

Trivial.

We can then use a method on the rewritten initial value problem.

**III.2.10 Definition : Truncation of a method**

Let $(\tilde{n}, s) \in \mathbb{N}^* \times \mathbb{N}^*$, $M = (\tau, w, \tilde{w}) \in \mathrm{GMORK}_{\tilde{n},s}$, and $n \in \mathbb{N}^*$ such that $n < \tilde{n}$.
We define $M' = (\tau, w', \tilde{w}') \in \mathrm{GMORK}_{n,s}$ the truncation of $M$ at length $n$ as :

$$\tau = \tau,\ \forall N'' \in [\![1, n]\!],\ w'_N = w_N,\ \tilde{w}'_N = \tilde{w}_N$$

**III.2.11 Proposition : Form of rewritten method - Lower order rewriting**

Let $(\tilde{n}, s) \in \mathbb{N}^* \times \mathbb{N}^*$, $M \in \mathrm{GMORK}_{\tilde{n},s}$, $n \in \mathbb{N}^*$ such that $n < \tilde{n}$, $M'$ the truncation of $M$ at length $n$, $f$ a differential equation function of order $n$, $(t, y_{n,0}, h) \in \mathbf{U} \times \mathbb{R}$, and $(f', t, y_{\tilde{n},0})$ the prolongation to order $\tilde{n}$ of $(f, t, y_{n,0})$.
The stage system $(M, f', t, y_{\tilde{n},0}, h)$ is equivalent to the stage system $(M', f, t, y_{n,0}, h)$.

Proof

The stage system $(M, f', t, y_{\tilde{n},0}, h)$ is :

$$\forall (j,N) \in [\![1,s+1]\!] \times [\![1,\tilde{n}]\!],$$
$$y_{\tilde{n},j,N} = \sum_{N'=1}^{N} \tilde{w}_{N,N',j} h^{N-N'} y_{\tilde{n},0,N'} + \frac{h^N}{N!} \sum_{j'=1}^{s} w_{N,j,j'} f'\left(t + \tau_{j'} h, y_{\tilde{n},j'}\right)$$
$$y_{\tilde{n},s+1} \in \mathbb{R}^{[\![1,\tilde{n}]\!] \times [\![1,d]\!]},\ \forall j \in [\![1,s]\!],\ \left(t + \tau_j h, y_{\tilde{n},j}\right) \in \mathbf{U}'$$

We define for all $(j,N) \in [\![1,s+1]\!] \times [\![1,n]\!],\ y_{\tilde{n},j,N} = y_{n,j,N}$. The other entries of $y_{\tilde{n}}$ are functions of $y_n$, we hence remove them from the system.

$$\forall (j,N) \in [\![1,s+1]\!] \times [\![1,n]\!],$$
$$y_{n,j,N} = \sum_{N'=1}^{N} \tilde{w}_{N,N',j} h^{N-N'} y_{n,0,N'} + \frac{h^N}{N!} \sum_{j'=1}^{s} w_{N,j,j'} f\left(t + \tau_{j'} h, y_{n,j'}\right)$$
$$y_{n,s+1} \in \mathbb{R}^{[\![1,n]\!] \times [\![1,d]\!]},\ \forall j \in [\![1,s]\!],\ \left(t + \tau_j h, y_{n,j}\right) \in \mathbf{U}$$

Which is the stage system $\left(M', f, t, y_{n,0}, h\right)$.

This rewriting is simple because $n - i$ is used to index the approximations. If instead they were indexed by $i$, the rewritten approximations would be indexed by $i + \tilde{n} - n$. However, we will define methods which can be used for any $n$, which will require to let $\tilde{n} = +\infty$, and we cannot use $i + \infty - n$ as an index.

**III.2.12 Proposition : Order of consistency - Lower order rewriting**

Let $(\tilde{n}, s) \in \mathbb{N}^* \times \mathbb{N}^*$, $M \in \mathrm{GMORK}_{\tilde{n},s}$, $M'$ the truncation of $M$ at length $n$, $\xi, \xi'$ the approximation times of respectively $M, M'$, $v_n, v'_n$ the minimum consistency matrices of respectively $M, M'$.
If $M$ is convergent up to length $\tilde{n}$ then $M'$ is convergent up to length $n$. We also have :

$$\forall (j,N) \in [\![1,s+1]\!] \times [\![1,n]\!],\ \xi'_{j,N} = \xi_{j,N}$$
$$\forall (j,N) \in [\![1,s+1]\!] \times [\![1,n]\!],\ v'_{n,j,N} \geq v_{n,j,N}$$

Proof

The approximation times do not change and the prolongation of an initial value problem is a special case of the initial value problems considered in the definition of the order of consistency, thus the order of consistency at least holds, and if $M$ is convergent then $M'$ is.

If $M$ is rewrite-compliant, the rewriting for lower order equation does not make the method lose in order of consistency. Otherwise, multiple derivatives could have been tracked in the entries of ŷ to raise the order of consistency. This rewriting allows us to efficiently solve differential equations of order $n \leq \tilde{n}$. To solve initial value problems of arbitrary order, we let $\tilde{n} = +\infty$.

**III.2.13 Definition : Infinite methods**

- An infinite general multi-order method is a tuple $(\tau, w, \tilde{w})$, with $\tau \in \mathbb{R}^{[\![1,s+1]\!]}$ its nodes, $w \in \mathbb{R}^{\mathbb{N}^* \times [\![1,s+1]\!] \times [\![1,s]\!]}$ its main weights, $\tilde{w} \in \mathbb{R}^{\left(\prod_{N \in \mathbb{N}^*} \{N\} \times [\![1,N]\!]\right) \times [\![1,s+1]\!]}$ its secondary weights.
- We define $\mathrm{GMORK}_{\infty,s}$ as the set of infinite general multi-order Runge-Kutta methods with $s \in \mathbb{N}^*$ points.

We extend some of the previously introduced concepts to infinite methods and to truncations of a method. We do not discuss their applications in details since their extension is straightforward. We start with Section II.

**III.2.14 Definition : Approximation functions & Evaluation functions & Error functions**

Let $(\tilde{n}, s) \in \overline{\mathbb{N}}^* \times \mathbb{N}^*$, $M \in \mathrm{GMORK}_{\tilde{n},s}$, $n \in \mathbb{N}^*$ such that $n \leq \tilde{n}$.
We define $Y_n$ / $F_n$ / $\varepsilon_n$ the approximation / evaluation / error function for length $n$ of $M$ as the approximation / evaluation / error function for length $n$ of the truncation of $M$ at length $n$.

**III.2.15 Definition : Permutations of a method**

Let $(\tilde{n}, s) \in \overline{\mathbb{N}}^* \times \mathbb{N}^*$, and $\mathfrak{S}_s^* = \left\{\varphi \in \mathfrak{S}_{[\![1,s+1]\!]} \mid \varphi(s+1) = s+1\right\}$.
We define the permute operation $*$ as, for all $\varphi \in \mathfrak{S}_s^*, (\tau, w, \tilde{w}) \in \mathrm{GMORK}_{\tilde{n},s},\ \varphi * (\tau, w, \tilde{w}) = (\tau', w', \tilde{w}') \in \mathrm{GMORK}_{\tilde{n},s}$, where :

$$\forall j \in [\![1, s]\!],\ \tau'_j = \tau_{\varphi(j)}$$

$$\forall (j, N) \in [\![1, s+1]\!] \times \mathbb{N}^*,\ N \leq \tilde{n},\ \forall N' \in [\![1, N]\!],\ \tilde{w}'_{N,N',j} = \tilde{w}_{N,N',\varphi(j)}$$

$$\forall (j, j', N) \in [\![1, s+1]\!] \times [\![1, s]\!] \times \mathbb{N}^*,\ N \leq \tilde{n},\ w'_{N,j,j'} = w_{N,\varphi(j),\varphi(j')}$$

**III.2.16 Definition : Equivalence by permutation**

Let $(\tilde{n}, s) \in \overline{\mathbb{N}}^* \times \mathbb{N}^*$.
We define the relation $\mathcal{R}_p$ on $\mathrm{GMORK}_{\tilde{n},s}$ as the equivalence relation induced by the group action of $\mathfrak{S}_s^*$ on $\mathrm{GMORK}_{\tilde{n},s}$, thus, for all $M, M' \in \mathrm{GMORK}_{\tilde{n},s}^2$ :

$$M\ \mathcal{R}_p\ M' \Leftrightarrow \exists \varphi \in \mathfrak{S}_s^*,\ \varphi * M = M'$$

$M$ and $M'$ are said to be equivalent by permutation if and only if $M\ \mathcal{R}_p\ M'$. $M'$ is said to be equivalent to $M$ by permutation $\varphi \in \mathfrak{S}_s^*$ if and only if $M' = \varphi * M$.

**III.2.17 Definition : Weight digraphs**

Let $(\tilde{n}, s) \in \overline{\mathbb{N}}^* \times \mathbb{N}^*$, $M \in \mathrm{GMORK}_{\tilde{n},s}$.
We define the weight digraphs $\mathcal{G}$ of $M$ as, for all $n \in \mathbb{N}^*$ such that $n \leq \tilde{n}$, $\mathcal{G}_n$ is the maximum weight digraph of the truncation of $M$ at length $n$.
If $\tilde{n} = +\infty$, we define $\mathcal{G}_\infty$ the maximum weight digraphs of $M$ as :

$$\mathcal{G}_\infty = ([\![1, s+1]\!], \mathcal{A}_\infty)$$

$$\mathcal{A}_\infty = \left\{(j', j) \in [\![1, s]\!] \times [\![1, s+1]\!] \mid \exists N \in \mathbb{N}^*,\ w_{N,j,j'} \neq 0\right\}$$

There exists another definition of the maximum weight digraph of an infinite method, but to state it we need the following proposition.

**III.2.18 Proposition : Monotonicity of the weight digraphs**

Let $(\tilde{n}, s) \in \overline{\mathbb{N}}^* \times \mathbb{N}^*$, $M \in \mathrm{GMORK}_{\tilde{n},s}$, and $\mathcal{G}$ the weight digraphs of $M$.
$\mathcal{G}$ is an increasing and bounded sequence inclusion wise.

Proof

It is trivial to verify that $\mathcal{A}_n \subset \mathcal{A}_{n+1} \subset [\![1, s]\!] \times [\![1, s+1]\!]$.

The fact that $\mathcal{G}$ is an increasing sequence justifies the name maximum weight digraph.

**III.2.19 Proposition : Maximum weight digraphs as a limit - Infinite methods**

Let $s \in \mathbb{N}^*$, $M \in \mathrm{GMORK}_{\infty,s}$, and $\mathcal{G}$ the weight digraphs of $M$.
$\mathcal{G}_\infty = \lim_{n\to\infty} \mathcal{G}_n$

Proof

Trivial.

We do not generalize the remaining theorems and definitions of the second section because most are easily generalizable, especially since most are defined in terms of the maximum weight digraph. The characterizations expressed in terms of weights only require replacing $N \in [\![1, n]\!]$ with $N \in \mathbb{N}^*$. Infinite method do not produce approximations directly, rather, the relevant properties must be applied to the approximations of their truncations.

We generalize the definitions we have seen in this section.

**III.2.20 Definition : Approximation times**

Let $(\tilde{n}, s) \in \overline{\mathbb{N}}^* \times \mathbb{N}^*$, and $M = (\tau, w, \tilde{w}) \in \mathrm{GMORK}_{\tilde{n},s}$.
We define $\xi \in \mathbb{R}^{[\![1,s+1]\!] \times \mathbb{N}^*}$ the approximation times of $M$ as, for all $N \in \mathbb{N}^*$ such that $N \leq \tilde{n}$ :

$$\forall j \in [\![1, s+1]\!],\ \xi_{j,N} = \begin{cases} \sum_{j'=1}^{s} w_{1,j,j'} & \text{if } N = 1 \\ \tilde{w}_{N,N-1,j} & \text{otherwise} \end{cases}$$

**III.2.21 Definition : Consistency tensor**

Let $(\tilde{n}, s) \in \overline{\mathbb{N}}^* \times \mathbb{N}^*$, and $M \in \mathrm{GMORK}_{\tilde{n},s}$.
We define $v \in \overline{\mathbb{Z}}^{\bigcup_{n \in \mathbb{N}^*, n \leq \tilde{n}} \{n\} \times [\![1,s+1]\!] \times [\![1,n]\!]}$ the consistency tensor of $M$ as, for all $n \in \mathbb{N}^*$ such that $n \leq \tilde{n}$, $v_n$ is the minimum consistency matrix of the truncation of $M$ at length $n$. $v_n$ is said to be the consistency matrix of $M$ at length $n$.

To define the minimum consistency matrix of an infinite method we need the following proposition.

**III.2.22 Proposition : Monotonicity of the consistency tensor**

Let $(\tilde{n}, s) \in \overline{\mathbb{N}}^* \times \mathbb{N}^*$, $M \in \mathrm{GMORK}_{\tilde{n},s}$, $v$ the consistency tensor of $M$, and $(n, n') \in \mathbb{N}^* \times \mathbb{N}^*$ such that $n \leq n' \leq \tilde{n}$.
For all $(j, N) \in [\![1, s+1]\!] \times [\![1, n]\!]$, $v_{n,j,N} \geq v_{n',j,N} \geq -N$.

Proof

This is a direct consequence of Section III.2.12 and Section III.1.8.

**III.2.23 Definition : Minimum consistency matrix - Infinite methods**

Let $s \in \mathbb{N}^*$, $M \in \mathrm{GMORK}_{\infty,s}$, and $v$ the consistency tensor of $M$.
We define $v_\infty \in \overline{\mathbb{Z}}^{[\![1,s+1]\!] \times \mathbb{N}^*}$ the minimum consistency matrix of $M$ as :

$$\forall (j, N) \in [\![1, s+1]\!] \times \mathbb{N}^*,\ v_{\infty,j,N} = \lim_{n \to \infty} v_{n,j,N}$$

**III.2.24 Definition : Rewrite-compliant methods - Truncations**

Let $(\tilde{n}, s) \in \overline{\mathbb{N}}^* \times \mathbb{N}^*$, $M \in \mathrm{GMORK}_{\tilde{n},s}$, and $v_{\tilde{n}}$ the minimum consistency matrix of $M$.
$M$ is said to be rewrite-compliant if and only if :

$$\forall (j, N) \in \mathbb{N}^*,\ N \leq \tilde{n},\ v_{\tilde{n},j,N+1} - v_{\tilde{n},j,N} \geq -1$$

**III.2.25 Definition : Convergence - Truncations**

Let $(\tilde{n}, s) \in \overline{\mathbb{N}}^* \times \mathbb{N}^*$, $M \in \mathrm{GMORK}_{\tilde{n},s}$, and $n \in \overline{\mathbb{N}}^*$ such that $n \leq \tilde{n}$.
$M$ is said to be convergent up to length $n$ if and only if, for all $n' \in \mathbb{N}^*$ such that $n' \leq n$, the truncation of $M$ at length $n'$ is convergent up to length $n$.

## The case $n > \tilde{n}$

This case is not the most interesting since we now have methods capable of solving initial value problems of arbitrary order, nevertheless, we consider it for completeness.

To solve a higher order initial value problem using a rewrite-compliant method, we aim to incorporate as many derivatives of $\hat{y}$ as possible into the new function $\hat{y}'$, while keeping in mind that the $\tilde{n}^{\text{th}}$ derivative must still be availabe. We deduce :

$$\hat{y}' = \begin{bmatrix} \frac{\mathrm{d}^{n-\tilde{n}}\hat{y}}{\mathrm{d}t^{n-\tilde{n}}} \\ \vdots \\ \hat{y} \end{bmatrix} = \mathcal{J}^{n-\tilde{n}}\hat{y}$$

**III.2.26 Definition : Rewriting for initial value problems of higher order**

Let $(n, \tilde{n}) \in \mathbb{N}^* \times \mathbb{N}^*$ with $n > \tilde{n}$, and $(f, t_0, y_{n,0})$ an initial value problem of order $n$. We define :

$$\mathbf{U}' = \left\{ \left( t, \begin{bmatrix} x_{1,\tilde{n}} \\ \vdots \\ x_{1,n} \end{bmatrix}, \dots \begin{bmatrix} x_{\tilde{n},\tilde{n}} \\ \vdots \\ x_{\tilde{n},n} \end{bmatrix} \right) \in \mathbb{R} \times \mathbb{R}^{[\![1,\tilde{n}]\!] \times [\![1,(n-\tilde{n}+1)d]\!]} \mid \left(t, x_{1,\tilde{n}}, \dots\ x_{\tilde{n}-1,\tilde{n}}, x_{\tilde{n},\tilde{n}}, \dots\ x_{\tilde{n},n}\right) \in \mathbf{U} \right\}$$

$$\forall \left( t, \begin{bmatrix} x_{1,\tilde{n}} \\ \vdots \\ x_{1,n} \end{bmatrix}, \dots \begin{bmatrix} x_{\tilde{n},\tilde{n}} \\ \vdots \\ x_{\tilde{n},n} \end{bmatrix} \right) \in \mathbf{U}'$$

$$f'\left( t, \begin{bmatrix} x_{1,\tilde{n}} \\ \vdots \\ x_{1,n} \end{bmatrix}, \dots \begin{bmatrix} x_{\tilde{n},\tilde{n}} \\ \vdots \\ x_{\tilde{n},n} \end{bmatrix} \right) = \begin{bmatrix} f\left(t, x_{1,\tilde{n}}, \dots\ x_{\tilde{n}-1,\tilde{n}}, x_{\tilde{n},\tilde{n}}, \dots\ x_{\tilde{n},n}\right) \\ x_{1,\tilde{n}} \\ \vdots \\ x_{\tilde{n}-1,\tilde{n}} \\ x_{\tilde{n},\tilde{n}} \\ \vdots \\ x_{\tilde{n},n-\tilde{n}} \end{bmatrix}$$

$$\forall N \in [\![1,\tilde{n}]\!],\ y_{\tilde{n},0,N} = \begin{bmatrix} y_{n,0,N} \\ \vdots \\ y_{n,0,N+n-\tilde{n}} \end{bmatrix}$$

We define $(f', t_0, y_{\tilde{n},0})$ as the reduction to order $\tilde{n}$ of $(f, t_0, y_{n,0})$.

The visualization of $f'$ in the definition is not always exact as it requires $\tilde{n} \leq n - \tilde{n}$. In the case where this is not true, the vector would instead have $x_{n-\tilde{n},\tilde{n}}$ at the bottom.

**III.2.27 Proposition : Solution of rewritten initial value problems - Higher order rewriting**

Let $(n, \tilde{n}) \in \mathbb{N}^* \times \mathbb{N}^*$ with $n > \tilde{n}$, $(f, t_0, y_{n,0})$ an initial value problem of order $n$, and $(f', t_0, y_{\tilde{n},0})$ the reduction to order $\tilde{n}$ of $(f, t_0, y_{n,0})$

- If $\hat{y}$ is a solution of $(f, t_0, y_{n,0})$ then $\mathcal{J}^{n-\tilde{n}}\hat{y}$ is a solution of $(f', t_0, y_{\tilde{n},0})$.
- If $\hat{y}'$ is a solution of $(f', t_0, y_{\tilde{n},0})$, then there exists a unique $\hat{y}$ such that $\hat{y}' = \mathcal{J}^{n-\tilde{n}}\hat{y}$. $\hat{y}$ is a solution of $(f, t_0, y_{n,0})$.

Proof

Trivial, the basic idea is that $x_{N,N'}$ represents the $n + \tilde{n} - (N + N')^{\text{th}}$ derivative of $\hat{y}$.

### III.2.28 Definition : Extension of a method

Let $(\tilde{n}, s) \in \mathbb{N}^* \times \mathbb{N}^*$, $M = (\tau, w, \tilde{w}) \in \mathrm{GMORK}_{\tilde{n},s}$, $n \in \mathbb{N}^*$ such that $n > \tilde{n}$. We define $q \in \{\mathbb{N} \to [\![1, n]\!]\}$ and $r \in \{\mathbb{N} \to [\![1, n]\!]\}$ as, for all $N \in \mathbb{N}$, $q(N) \in \mathbb{N}^*$ is the result of the Euclidean division of $N$ by $\tilde{n}$ and $r(N)$ the remainder of the Euclidean division.
We define $M' = (\tau, w', \tilde{w}') \in \mathrm{GMORK}_{n,s}$ the extension of $M$ to length $n$ as :

$$\forall (N, j, j') \in [\![1, n]\!] \times [\![1, s+1]\!] \times [\![1, s]\!]$$

$$w'_{N,j,j'} = N! \sum_{\substack{J \in [\![1, s+1]\!]^{[\![1, q(N-1)+2]\!]} \\ J_1 = j,\ J_{q(N-1)+2} = j' \\ k \neq 1 \Rightarrow J_k \leq s}} \left( \prod_{k=1}^{q(N-1)} \frac{w_{\tilde{n}, J_k, J_{k+1}}}{\tilde{n}!} \right) \frac{w_{1+r(N-1), J_{q(N-1)+1}, J_{q(N-1)+2}}}{(1 + r(N-1))!}$$

$$\forall (j, N) \in [\![1, s+1]\!] \times [\![1, n]\!],\ \forall N' \in [\![1, N]\!]$$

$$\tilde{w}'_{N,N',j} = \sum_{\substack{J \in [\![1, s+1]\!]^{[\![1, q(N-N')+1]\!]} \\ J_1 = j,\ k \neq 1 \Rightarrow J_k \leq s}} \left( \prod_{k=1}^{q(N-N')} \frac{w_{\tilde{n}, J_k, J_{k+1}}}{\tilde{n}!} \right) \tilde{w}_{\min(N' + r(N-N'), \tilde{n}),\, \min(N', \tilde{n} - r(N-N')), J_{q(N-N')+1}}$$

### III.2.29 Proposition : Form of rewritten method - Higher order rewriting

Let $(\tilde{n}, s) \in \mathbb{N}^* \times \mathbb{N}^*$, $M \in \mathrm{GMORK}_{\tilde{n},s}$, $n \in \mathbb{N}^*$ such that $n > \tilde{n}$, $M' \in \mathrm{GMORK}_{n,s}$ the extension of $M$ to length $n$, $f$ a differential equation function of order $n$, $(t, y_{n,0}, h) \in \mathbf{U} \times \mathbb{R}$, and $(f', t, y_{\tilde{n},0})$ the reduction to order $\tilde{n}$ of $(f, t, y_{n,0})$.
The stage system $(M, f', t, y_{\tilde{n},0}, h)$ is equivalent to the stage system $(M', f, t, y_{n,0}, h)$.

Proof

Let $(\tau, w, \tilde{w}) = M$ and $(\tau', w', \tilde{w}') = M'$. The stage system $(M, f', t, y_{\tilde{n},0}, h)$ is :

$$\forall (j, N) \in [\![1, s+1]\!] \times [\![1, \tilde{n}]\!],$$

$$y_{\tilde{n},j,N} = \sum_{N'=1}^{N} \tilde{w}_{N,N',j} h^{N-N'} y_{\tilde{n},0,N'} + \frac{h^N}{N!} \sum_{j'=1}^{s} w_{N,j,j'} f'\left(t + \tau_{j'} h, y_{\tilde{n},j'}\right)$$

$$y_{\tilde{n},s+1} \in \mathbb{R}^{[\![1,\tilde{n}]\!] \times [\![1,(n-\tilde{n}+1)d]\!]},\ \forall j \in [\![1, s]\!],\ \left(t + \tau_j h, y_{\tilde{n},j}\right) \in \mathbf{U}'$$

We define :

$$\forall j \in [\![0, s+1]\!],\ y_{\tilde{n},j,N} = \begin{bmatrix} y'_{\tilde{n},j,N,\tilde{n}} \\ \vdots \\ y'_{\tilde{n},j,N,n} \end{bmatrix}$$

The approximations that interest us are the one used in the evaluation of $f$, their ranks are also given by the order we gave them in the evaluation of $f$, we hence define :

$$\begin{bmatrix} y_{n,j,\tilde{n}} \\ \vdots \\ y_{n,j,n} \end{bmatrix} = \begin{bmatrix} y'_{\tilde{n},j,\tilde{n},\tilde{n}} \\ \vdots \\ y'_{\tilde{n},j,\tilde{n},n} \end{bmatrix}, \quad \forall N \in [\![1, \tilde{n}]\!], \begin{bmatrix} y_{n,j,N} \\ * \\ \vdots \\ * \end{bmatrix} = \begin{bmatrix} y'_{\tilde{n},j,N,\tilde{n}} \\ \vdots \\ y'_{\tilde{n},j,N,n} \end{bmatrix}$$

Which can be summed up with :

$$\forall N \in [\![\tilde{n}, n]\!],\ y_{n,j,N} = y'_{\tilde{n},j,\tilde{n},N},\ \forall N \in [\![1, \tilde{n}]\!],\ y_{n,j,N} = y'_{\tilde{n},j,N,\tilde{n}}$$

We have, for all $(N, N') \in [\![1, \tilde{n}]\!] \times [\![\tilde{n}, n]\!],\ y'_{n,0,N,N'} = y_{n,0,N+N'-\tilde{n}}$. With these definitions, we get :

$$\forall (j, N) \in [\![1, s+1]\!] \times [\![1, \tilde{n}]\!],$$

$$\begin{bmatrix} y'_{\tilde{n},j,N,\tilde{n}} \\ \vdots \\ y'_{\tilde{n},j,N,n} \end{bmatrix} = \sum_{N'=1}^{N} \tilde{w}_{N,N',j} h^{N-N'} \begin{bmatrix} y'_{\tilde{n},0,N',\tilde{n}} \\ \vdots \\ y'_{\tilde{n},0,N',n} \end{bmatrix} + \frac{h^N}{N!} \sum_{j'=1}^{s} w_{N,j,j'} \begin{bmatrix} f\left(t_0 + \tau_{j'} h, y_{n,j'}\right) \\ y_{n,j',1} \\ \vdots \\ y_{n,j',n-\tilde{n}} \end{bmatrix}$$

$$y_{n,s+1} \in \mathbb{R}^{[\![1,n]\!] \times [\![1,d]\!]},\ \forall j \in [\![1, s]\!],\ \left(t + \tau_j h, y_{n,j}\right) \in \mathbf{U}$$

We deduce :

$$\forall N \in [\![\tilde{n}+1, n]\!],\ y_{n,j,N} = \sum_{N'=N-\tilde{n}+1}^{N} \tilde{w}_{\tilde{n},\tilde{n}+N'-N,j} h^{N-N'} y_{n,0,N'} + \frac{h^{\tilde{n}}}{\tilde{n}!} \sum_{j'=1}^{s} w_{\tilde{n},j,j'} y_{n,j',N-\tilde{n}}$$

$$\forall N \in [\![1, \tilde{n}]\!],\ y_{n,j,N} = \sum_{N'=1}^{N} \tilde{w}_{N,N',j} h^{N-N'} y_{n,0,N'} + \frac{h^N}{N!} \sum_{j'=1}^{s} w_{N,j,j'} f\left(t_0 + \tau_{j'} h, y_{n,j'}\right)$$

We prove that, for all $N \in [\![1, n]\!]$ :

$$y_{n,j,N} = \sum_{N'=1}^{N} \sum_{\substack{J \in [\![1,s+1]\!]^{[\![1,q(N-N')+1]\!]} \\ J_1 = j,\ k \neq 1 \Rightarrow J_k \leq s}} \left( \prod_{k'=1}^{q(N-N')} \frac{w_{\tilde{n},J_{k'},J_{k'+1}}}{\tilde{n}!} \right) \Bigg[$$

$$\tilde{w}_{\min(N'+r(N-N'),\tilde{n}),\, \min(N',\tilde{n}-r(N-N')),J_{q(N-N')+1}} h^{N-N'} y_{n,0,N'} \Bigg] + h^N \sum_{\substack{J \in [\![1,s+1]\!]^{[\![1,q(N-1)+2]\!]} \\ J_1 = j,\ k \neq 1 \Rightarrow J_k \leq s}} \Bigg[$$

$$\left( \prod_{k'=1}^{q(N-1)} \frac{w_{\tilde{n},J_{k'},J_{k'+1}}}{\tilde{n}!} \right) \frac{w_{1+r(N-1),J_{q(N-1)+1},J_{q(N-1)+2}}}{(1+r(N-1))!} f\left(t_0 + \tau_{J_{q(N-1)+2}} h, y_{n,J_{q(N-1)+2}}\right) \Bigg]$$

This can be proven by induction on $p \in \mathbb{N}$ where $p$ is the value of $q(N-1)$. The integers $N$ satisfying $q(N-1) = p$ for a given $p \in \mathbb{N}$ are those that satisfy $p\tilde{n} < N \leq (p+1)\tilde{n}$. The base case $p = 0$ has already been proven. The induction part is simply replacing $y_{N-\tilde{n}}$ with its expression and rewriting the sum of the main weights as a new coordinate of $J$. Extracting the coefficients proves the theorem.

**III.2.30 Proposition : Order of consistency - Higher order rewriting**

Let $(\tilde{n}, s) \in \mathbb{N}^* \times \mathbb{N}^*$, $M \in \mathrm{GMORK}_{\tilde{n},s}$, $v_{\tilde{n}}$ its minimum consistency matrix, $n \in \mathbb{N}^*$ such that $n > \tilde{n}$, $M'$ the extension of $M$ to length $n$, and $v_n$ its minimum consistency matrix.
If $M$ is convergent up to length $\tilde{n}$ then $M'$ is convergent up to length $n$. We have :

$$\forall (j, N) \in [\![1, s+1]\!],\ [\![\tilde{n}, n]\!],\ \xi_{j,N} = \xi_{j,\tilde{n}},\ v_{n,j,N} \geq v_{\tilde{n},j,\tilde{n}} + \tilde{n} - N$$

$$\forall (j, N) \in [\![1, s+1]\!] \times [\![1, \tilde{n}]\!],\ \xi_{j,N} = \xi'_{j,N},\ v_{n,j,N} \geq v_{\tilde{n},j,N}$$

Proof

The values of the approximation times is trivial if we stop at this steo of the stage system in the proof :

$$\forall N \in [\![\tilde{n}+1,n]\!],\ y_{n,j,N} = \sum_{N'=N-\tilde{n}+1}^{N} \tilde{w}_{\tilde{n},\tilde{n}+N'-N,j} h^{N-N'} y_{n,0,N'} + \frac{h^{\tilde{n}}}{\tilde{n}!} \sum_{j'=1}^{s} w_{\tilde{n},j,j'} y_{n,j',N-\tilde{n}}$$

$$\forall N \in [\![1,\tilde{n}]\!],\ y_{n,j,N} = \sum_{N'=1}^{N} \tilde{w}_{N,N',j} h^{N-N'} y_{n,0,N'} + \frac{h^{N}}{N!} \sum_{j'=1}^{s} w_{N,j,j'} f\big(t_0 + \tau_{j'} h, y_{n,j'}\big)$$

The approximation times do not change and the reduction to order $\tilde{n}$ of an initial value problem is a special case of the initial value problems considered in the definition of the order of consistency, thus the convergence and the the order of consistency of $M$ still holds. Since the ranks $N \in [\![\tilde{n},n]\!]$ use the order of consistency at rank $\tilde{n}$, we have $\underset{h\to 0}{\mathcal{O}}\left(h^{\tilde{n}+v_{\tilde{n},j,\tilde{n}}}\right) = \underset{h\to 0}{\mathcal{O}}\left(h^{N+v_{\tilde{n},j,\tilde{n}}+\tilde{n}-N}\right)$, we deduce $v_{n,j,N} \geq v_{\tilde{n},j,\tilde{n}} + \tilde{n} - N$.

We now examine the effects of the rewriting on the maximum weight digraph.

### III.2.31 Proposition : Maximum weight digraph - Higher order rewriting

Let $(n,\tilde{n},s) \in \mathbb{N}^* \times \mathbb{N}^* \times \mathbb{N}^*$ with $\tilde{n} < n$, $M \in \mathrm{GMORK}_{\tilde{n},s}$, $M'$ the extension of $M$ to length $n$, $\mathcal{R}_r$, $\mathcal{R}_r{}'$ the reachability relations of respectively $M, M'$.

1. For all $(j,j') \in [\![1,s+1]\!] \times [\![1,s]\!],\ j' \ \mathcal{R}_r{}' \ j \Rightarrow j' \ \mathcal{R}_r \ j$.
2. Any useless stage, reachable-closed set, explicit stage, block, explicit block, explicit partition of the stages, order of computation, explicit rank, parallel blocks, independent set of stages, partition of independent sets of stages of $M$ is also one of $M'$. If $M$ is explicit $M'$ also is.

Proof

1.
$$w'_{N,j,j'} \neq 0 \Leftrightarrow N! \sum_{\substack{J\in[\![1,s+1]\!]^{[\![1,q(N-1)+2]\!]} \\ J_1=j,\ J_{q(N-1)+2}=j' \\ k\neq 1 \Rightarrow J_k \leq s}} \left(\prod_{k'=1}^{q(N-1)} \frac{w_{\tilde{n},J_{k'},J_{k'+1}}}{\tilde{n}!}\right) \frac{w_{1+r(N-1),J_{q(N-1)+1},J_{q(N-1)+2}}}{(1+r(N-1))!} \neq 0$$

$$\Rightarrow \exists J \in [\![1,s+1]\!]^{[\![1,q(N-1)+2]\!]},\ J_1 = j,\ J_{q(N-1)+2} = j',\ k \neq 1 \Rightarrow J_k \leq s,$$

$$\left(\prod_{k'=1}^{q(N-1)} w_{\tilde{n},J_{k'},J_{k'+1}}\right) w_{1+r(N-1),J_{q(N-1)+1},J_{q(N-1)+2}} \neq 0$$

$$\Rightarrow \exists J \in [\![1,s+1]\!]^{[\![1,q(N-1)+2]\!]},\ J_1 = j,\ J_{q(N-1)+2} = j',\ k \neq 1 \Rightarrow J_k \leq s$$

$$w_{1+r(N-1),J_{q(N-1)+1},J_{q(N-1)+2}} \neq 0,\ \forall k' \in [\![1,q(N-1)+1]\!],\ w_{\tilde{n},J_{k'},J_{k'+1}} \neq 0$$

Hence if $j' \ \mathcal{R}_r{}' \ j$, then there exists $J$ a diwalk from $j$ to $j'$ in the digraph of $M$, thus $j' \ \mathcal{R}_r \ j$.

2. It is a direct consequence of 1.

The seemingly complex definition of the extension of a method arises from the fact that it conceals multiple matrix multiplications. This becomes clearer in the following lemma and proposition.

### III.2.32 Lemma: Developed matrix multiplication

Let $q \in \mathbb{N}^k$, $a \in \mathbb{N}^{*[\![1,q+1]\!]}$, $A \in \mathbb{R}^{\cup_{k=1}^{q}\{k\}\times[\![1,a_k]\!]\times[\![1,a_{k+1}]\!]}$, and $(i,j) \in [\![1,a_1]\!] \times [\![1,a_{q+1}]\!]$.

We have :

$$\left(\prod_{k=1}^{q} A_k\right)_{i,j} = \sum_{\substack{(J_1,\dots J_{q+1})\in[\![1,a_1]\!]\times\dots\times[\![1,a_{q+1}]\!] \\ J_1=i,J_{q+1}=j}} \prod_{k=1}^{q} A_{k,J_k,J_{k+1}}$$

Proof

It is trivial to prove by induction on $q$.

This sum closely resembles those appearing in the definition of the extension of a method. We therefore rewrite the expression of the extension in terms of matrix multiplication.

### III.2.33 Proposition : Matrix expression of the extension of a method

Let $(\tilde{n}, s) \in \mathbb{N}^* \times \mathbb{N}^*$, $M = (\tau, w, \tilde{w}) \in \mathrm{GMORK}_{\tilde{n},s}$, $n \in \mathbb{N}^*$ such that $n > \tilde{n}$, and $M' = (\tau', w', \tilde{w}') \in \mathrm{GMORK}_{n,s}$ the extension of $M$ to length $n$. We define $q \in \{\mathbb{N} \to [\![1, n]\!]\}$ and $r \in \{\mathbb{N} \to [\![1, n]\!]\}$ as in Section III.2.28, and $E_s = [I_s \;\; 0_s]$.

We have :

$$\forall N \in [\![1, \tilde{n}]\!],\; w'_N = \begin{cases} w_N & \text{if } N \leq \tilde{n} \\ N! \dfrac{w_{\tilde{n}}}{\tilde{n}!} \left( \dfrac{E_s w_{\tilde{n}}}{\tilde{n}!} \right)^{q(N-1)-1} \dfrac{w_{1+r(N-1)}}{(1 + r(N-1))!} & \text{otherwise} \end{cases}$$

$$\forall N \in [\![1, n]\!],\; \forall N' \in [\![1, N]\!]$$

$$\tilde{w}'_{N,N'} = \begin{cases} \tilde{w}_{\min(N,\tilde{n}),\, \min(N', \tilde{n}-N+N')} & \text{if } N' > N - \tilde{n} \\ \dfrac{w_{\tilde{n}}}{\tilde{n}!} \left( \dfrac{E_s w_{\tilde{n}}}{\tilde{n}!} \right)^{q(N-N')-1} E_s \tilde{w}_{\min(N'+r(N-N'),\tilde{n}),\, \min(N', \tilde{n}-r(N-N'))} & \text{otherwise} \end{cases}$$

Proof

It is a direct consequence of Section III.2.32.

The matrix $E_s$ is used to enforce the condition $k \neq 1 \Rightarrow J_k \leq s$, it removes the weights of the last stage :

$$\forall N \in [\![1, \tilde{n}]\!],\; E_s w_N = \begin{bmatrix} w_{N,1,1} & \cdots & w_{1,s} \\ \vdots & \ddots & \vdots \\ w_{N,s,1} & \cdots & w_{s,s} \end{bmatrix}$$

This matrix formulation is particularly useful for understanding the limitations of first-order reduction.

### III.2.34 Corollary : Extension of a Runge-Kutta method

Let $M = (\tau, w_1) \in \mathrm{RK}_s$, $n \in \mathbb{N}^*$ such that $n$, and $M' = (\tau', w', \tilde{w}') \in \mathrm{GMORK}_{n,s}$ the extension of $M$ to length $n$. We define $E_s = [I_s \;\; 0_s]$.

We have :

$$\forall N \in [\![1, \tilde{n}]\!],\; w'_N = \begin{cases} w_1 & \text{if } N = 1 \\ N! w_1 (E_s w_1)^{N-1} & \text{otherwise} \end{cases}$$

$$\forall N \in [\![1, n]\!],\; \forall N' \in [\![1, N]\!],\; \tilde{w}'_{N,N'} = \begin{cases} 1_{s+1} & \text{if } N' = N \\ w_1 (E_s w_1)^{q(N-N')-1} 1_s & \text{otherwise} \end{cases}$$

The fact that the weights of a rewritten Runge-Kutta method are given by powers of $w_1$ explains why explicit methods, whose weight matrix is strictly lower triangular and therefore nilpotent, contain many zero coefficients.

## The mixed-order case

In the beginning of this paper we stated a simplified definition of an initial value problem which assumed that the entries of $\hat{y}$ define an initial value problem of the same order, however it is not the case in general.

**III.2.35 Definition : Initial value problem - Mixed-order case**

Let $d \in \mathbb{N}^*$, $n \in \mathbb{N}^{*[\![1,d]\!]}$, $\mathbf{U} \in \mathcal{P}\left(\mathbb{R} \times \mathbb{R}^{\bigcup_{k=1}^{d} [\![1,n_k]\!] \times \{k\}}\right)$ an open connected non-empty set, $(t_0, y_{n,0}) \in \mathbf{U}$ and $f \in \{\mathbf{U} \to \mathbb{R}\}^{[\![1,d]\!]}$.

- An initial value problem is a tuple $(f, t_0, y_{n,0})$.
- $\mathbf{U}$ is called the space time domain.
- $n$ is the order vector, $d$ is the dimension and $\|n\|_d$ is the order of the initial value problem.
- $f$ is said to be a differential equation function of dimension $d$ with order vector $n$.
- $t_0$ is called the initial instant, and $y_{n,0}$ is called the initial values.
- We define $\Omega = \left\{t \in \mathbb{R} \mid \exists x \in \mathbb{R}^{\bigcup_{k=1}^{d} [\![1,n_k]\!] \times \{k\}}, (t,x) \in \mathbf{U}\right\}$ the time domain, an open interval of $\mathbb{R}$.

Let $\hat{y} \in \{\Omega \to \mathbb{R}\}^{[\![1,d]\!]}$ such that, for all $k \in [\![1,d]\!]$, $\hat{y}_k \in D^{n_k}(\Omega \to \mathbb{R})$.

- $\hat{y}$ is said to be a solution candidate if and only if :

$$\forall t \in \Omega, \ \left(t, \frac{\mathrm{d}^{n_1-1}\hat{y}_1}{\mathrm{d}t^{n_1-1}}(t), ..., \hat{y}_1(t), ..., \frac{\mathrm{d}^{n_d-1}\hat{y}_d}{\mathrm{d}t^{n_d-1}}(t), ..., \hat{y}_d(t)\right) \in \mathbf{U}$$

- $\hat{y}$ is said to satisfy the differential equation induced by $f$ if and only if $\hat{y}$ is a solution candidate and :

$$\forall (t,k) \in \Omega \times [\![1,d]\!], \ \frac{\mathrm{d}^{n_k}\hat{y}_k}{\mathrm{d}t^{n_k}}(t) = f_k\left(t, \frac{\mathrm{d}^{n_1-1}\hat{y}_1}{\mathrm{d}t^{n_1-1}}(t), ..., \hat{y}_1(t), ..., \frac{\mathrm{d}^{n_d-1}\hat{y}_d}{\mathrm{d}t^{n_d-1}}(t), ..., \hat{y}_d(t)\right)$$

- $\hat{y}$ is said to satisfy the initial values if and only if :

$$\forall k \in [\![1,d]\!], \ \forall N \in [\![1,n_k]\!], \ \frac{\mathrm{d}^{n-N}\hat{y}_k}{\mathrm{d}t^{n-N}}(t_0) = y_{n,0,N,k}$$

- $\hat{y}$ is said to be a solution of the initial value problem if and only if $\hat{y}$ is a solution candidate, $\hat{y}$ satisfies the differential equation induced by $f$ and $\hat{y}$ satisfies the initial values.

Let $\Omega' \in \mathcal{P}(\mathbb{R})$ an open connected subset of $\Omega$ which contains $t_0$, and $\hat{y}' \in \{\Omega \to \mathbb{R}\}^{[\![1,d]\!]}$ such that for all $k \in [\![1,d]\!]$, $\hat{y}'_k \in D^{n_k}(\Omega' \to \mathbb{R})$.

- $\hat{y}'$ is said to be a maximal solution if and only if $\hat{y}'$ satisfies the initial value problem restricted to $\Omega'$, and any prolongation of $\hat{y}'$ to $\Omega''$ any open connected subset of $\Omega$ and superset of $\Omega'$, does not satisfy the initial value problem restricted to $\Omega''$.

We deduce from the previous rewritings the following proposition.

**III.2.36 Proposition : Use of GMORK methods - Mixed order case**

Let $(\tilde{n}, s) \in \overline{\mathbb{N}}^* \times \mathbb{N}^*$, $(\tau, w, \tilde{w}) \in \mathrm{GMORK}_{n,s}$, $f$ a differential equation function of dimension $d \in \mathbb{N}^*$ with order vector $n \in \mathbb{N}^{*[\![1,d]\!]}$ such that $\|n\|_d \leq \tilde{n}$, $\left(t, y_{n,0}\right) \in \mathbb{R} \times \mathbf{U}$, $\tilde{q} \in \overline{\mathbb{N}}^*$, and $h \in \mathcal{H}_{\tilde{q},f,t,y_{n,0}}$.
$M$ approximates, for all $q \in \mathbb{N}$ such that $q < \tilde{q}$ :

$$\forall k \in [\![1, d]\!],\ \forall (N, j) \in [\![1, n_k]\!] \times [\![1, s+1]\!]$$

$$y_{n,q(s+1)+j,N,k} = \sum_{N'=1}^{n_k} \tilde{w}_{N,N',j} h_q^{N-N'} y_{n,q(s+1),N',k} + \frac{h_q^N}{N!} \sum_{j'=1}^{s} w_{N,j,j'} f_k\left(t_q + \tau_{j'} h_q, y_{n,j',1,1}, \dots\, y_{n,j',n_d,d}\right)$$

This expression is derived using a generalized version of the prolongation to order $\|n\|_d$, and hence satisfies similar properties. In particular, the accuracy of the approximation is governed by the consistency matrix of $M$ at length $\|n\|_d$. To deal with the case $\|n\|_d > \tilde{n}$, $M$ can be replaced with its extension to order $\|n\|_d$. For mixed-order initial values problems we can also define a mixed-order jet and extended jet.

**III.2.37 Definition : Jet & Extended jet - Mixed order case**

Let $d \in \mathbb{N}^*$, $n \in \mathbb{N}^{*[\![1,d]\!]}$, $\Omega \in \mathcal{P}(\mathbb{R})$ an open set, and $\hat{y} \in \{\Omega \to \mathbb{R}\}^{[\![1,d]\!]}$ such that, for all $k \in [\![1, d]\!]$, $\hat{y}_k \in D^{n_k}(\Omega \to \mathbb{R})$.

- We define $\mathcal{J}^n \hat{y} \in \left\{\Omega \to \mathbb{R} \times \mathbb{R}^{[\![0,n]\!] \times [\![1,d]\!]}\right\}$ the mixed-order $n$-jet of $\hat{y}$ as :

$$\mathcal{J}^n \hat{y} = \left(\frac{\mathrm{d}^{n_1} \hat{y}_1}{\mathrm{d}t^{n_1}}, \dots, \hat{y}_1, \dots, \frac{\mathrm{d}^{n_d} \hat{y}_d}{\mathrm{d}t^{n_d}}, \dots, \hat{y}_d\right) = \left(\mathcal{J}^{n_1} \hat{y}_1, \dots, \mathcal{J}^{n_d} \hat{y}_d\right)$$

- We define $\overline{\mathcal{J}}^n \hat{y} \in \left\{\Omega \to \mathbb{R} \times \mathbb{R}^{[\![0,n]\!] \times [\![1,d]\!]}\right\}$ the extended mixed-order $n$-jet of $\hat{y}$ as :

$$\forall t \in \Omega,\ \overline{\mathcal{J}}^n \hat{y}(t) = \left(t, \mathcal{J}^n \hat{y}(t)\right)$$

**III.2.38 Proposition : Initial value problem & Extended jets - Mixed order case**

Let $\left(f, t_0, y_{n,0}\right)$ an initial value problem of dimension $d \in \mathbb{N}^*$ with an order vector $n \in \mathbb{N}^{*[\![1,d]\!]}$, and $\hat{y} \in \{\Omega \to \mathbb{R}\}^{[\![1,d]\!]}$ such that, for all $k \in [\![1, d]\!]$, $\hat{y}_k \in D^{n_k}(\Omega \to \mathbb{R})$.

- $\hat{y}$ is a solution candidate if and only if, for all $t \in \Omega$, $\overline{\mathcal{J}}^{n-1_d} \hat{y}(t) \in \mathbf{U}$
- $\hat{y}$ satisfies the differential equation if and only if $\hat{y}$ is a solution candidate and $\frac{\mathrm{d}^n \hat{y}}{\mathrm{d}t^n}(t) = f \circ \overline{\mathcal{J}}^{n-1_d} \hat{y}$
- $\hat{y}$ is satisfies the initial values if and only if $\mathcal{J}^{n-1_d} \hat{y}(t_0) = y_{n,0}$.

**III.2.39 Definition : Solution curve**

Let $\left(f, t_0, y_{n,0}\right)$ an initial value problem of dimension $d \in \mathbb{N}^*$ with an order vector $n \in \mathbb{N}^{*[\![1,d]\!]}$, and $\hat{y} \in \{\Omega \to \mathbb{R}\}^{[\![1,d]\!]}$ such that, for all $k \in [\![1, d]\!]$, $\hat{y}_k \in D^{n_k}(\Omega \to \mathbb{R})$ and $\hat{y}$ is a solution of the initial value problem, if it exists.
We define the solution curve of $\hat{y}$ as the function $\overline{\mathcal{J}}^{n-1_d} \hat{y}$.

These definitions are somewhat secondary to the main development. In the special case where all entries of the order vector are equal, we do not recover the original definitions of jets and initial value problems given at the beginning of the paper. A fully satisfactory generalization would require indexing first by the component, instead of writing $y_{n_k,j,N,k}$, we would use $y_{k,n_k,j,N}$ or $y_{j,k,n_k,N}$. Since this idexing is less convenient for the setting considered in the majority of the paper, we do not adopt it.

# IV Comparisons & Stability

## IV.1 Comparisons & Numerical experiments

In this section we compare node-determined multi-order Runge-Kutta methods with Runge-Kutta methods. To make the comparison fair, we only compare methods with the same number of points $s$ and the same structure, hence the same maximum weight digraph. We first consider computational cost.

The only difference in cost is the computation of Taylor polynomials. The added cost for a single stage is $d\sum_{N=1}^{n}(N-1)$ additions and multiplications, hence $\frac{(s+1)(n-1)nd}{2}$ in total. To actually reach this number of additions and multiplication, the program needs to cache the powers of $h$ up to $n$ and the powers of $\tau_j h$ up to $n-1$, which necessitates respectively $n-1$ multiplications and $s(n-2)$ multiplications. To change the step size, a total of $(n-2)(s+1)+1$ multiplications are hence carried out.

Although this additional cost is not asymptotically negligible, initial value problems of very high order are uncommon, and even in such cases the trade-off remains reasonable. For instance, if $n=10$, the multi-order midpoint method (see below) approximates the solution $\hat{y}$ with an error of order $\underset{h\to 0}{\mathcal{O}}(h^{12})$, whereas its classical Runge-Kutta counterpart achieves only $\underset{h\to 0}{\mathcal{O}}(h^{3})$. This improvement is obtained at the cost of merely $135d$ additions and multiplications.

An explicit Runge-Kutta method would require a large number of stages to reach comparable accuracy, leading to a substantially higher computational cost, especially when evaluations of $f$ are expensive. Moreover, the numbers we gave are only upper bounds, several optimizations are possible. For instance, if a node-determined multi-order Runge-Kutta method contains repeated nodes, such as RK4, the associated Taylor polynomials and the powers of $\tau_j h$ need only be computed once for each distinct node.

There may also be a potential cost advantage for implicit multi-order Runge-Kutta methods. In these methods, the evaluations of $f$ are scaled by a factor $\frac{1}{N!}$, which is typically smaller than the corresponding coefficients appearing in classical Runge-Kutta methods. The associated Picard iterations may exhibit faster convergence toward a fixed point. This observation remains heuristic and has not yet been formally analyzed or experimentally validated.

We now present a series of numerical experiments. We begin by recalling several well-known Runge-Kutta methods and deriving explicit expression for certain node-determined multi-order Runge-Kutta methods. We then compare these methods using the implementation of Section A.

We also introduce two new methods, to the best of the author's knowledge, named here RK4b and CNb. These are the truncations at length 1 of respectively MORK4b and Multi-order CNb, two multi-order Runge-Kutta methods whose expressions are derived below.

We introduce some Runge-Kutta methods :

Explicit Euler method

$$\tau=\begin{bmatrix}0\\1\end{bmatrix},\ w_1=\begin{bmatrix}0\\1\end{bmatrix}$$

Explicit midpoint method

$$\tau=\begin{bmatrix}0\\\frac{1}{2}\\1\end{bmatrix},\ w_1=\begin{bmatrix}0&0\\\frac{1}{2}&0\\0&1\end{bmatrix}$$

Ralston's method

$$\tau = \begin{bmatrix} 0 \\ \frac{2}{3} \\ 1 \end{bmatrix},\ w_1 = \begin{bmatrix} 0 & 0 \\ \frac{2}{3} & 0 \\ \frac{1}{4} & \frac{3}{4} \end{bmatrix}$$

Heun's method

$$\tau = \begin{bmatrix} 0 \\ \frac{1}{3} \\ \frac{2}{3} \\ 1 \end{bmatrix},\ w_1 = \begin{bmatrix} 0 & 0 & 0 \\ \frac{1}{3} & 0 & 0 \\ 0 & \frac{2}{3} & 0 \\ \frac{1}{4} & 0 & \frac{3}{4} \end{bmatrix}$$

RK4

$$\tau = \begin{bmatrix} 0 \\ \frac{1}{2} \\ \frac{1}{2} \\ 1 \end{bmatrix},\ w_1 = \begin{bmatrix} 0 & 0 & 0 & 0 \\ \frac{1}{2} & 0 & 0 & 0 \\ 0 & \frac{1}{2} & 0 & 0 \\ 0 & 0 & 1 & 0 \\ \frac{1}{6} & \frac{1}{3} & \frac{1}{3} & \frac{1}{6} \end{bmatrix}$$

RK4b

$$\tau = \begin{bmatrix} 0 \\ \frac{1}{2} \\ \frac{1}{2} \\ 1 \\ 1 \end{bmatrix},\ w_1 = \begin{bmatrix} 0 & 0 & 0 & 0 \\ \frac{1}{2} & 0 & 0 & 0 \\ \frac{1}{4} & \frac{1}{4} & 0 & 0 \\ 0 & -1 & 2 & 0 \\ \frac{1}{6} & 0 & \frac{2}{3} & \frac{1}{6} \end{bmatrix}$$

Implicit Euler method

$$\tau = \begin{bmatrix} 1 \\ 1 \end{bmatrix},\ w_1 = \begin{bmatrix} 1 \\ 1 \end{bmatrix}$$

Implicit midpoint method

$$\tau = \begin{bmatrix} \frac{1}{2} \\ 1 \end{bmatrix},\ w_1 = \begin{bmatrix} \frac{1}{2} \\ 1 \end{bmatrix}$$

Crank-Nicolson method

$$\tau = \begin{bmatrix} 0 \\ 1 \\ 1 \end{bmatrix},\ w_1 = \begin{bmatrix} 0 & 0 \\ \frac{1}{2} & \frac{1}{2} \\ \frac{1}{2} & \frac{1}{2} \end{bmatrix}$$

CNb

$$\tau = \begin{bmatrix} 0 \\ \frac{2}{3} \\ 1 \end{bmatrix},\ w_1 = \begin{bmatrix} 0 & 0 \\ \frac{1}{3} & \frac{1}{3} \\ \frac{1}{4} & \frac{3}{4} \end{bmatrix}$$

Gauss-Legendre method of order 4

$$\tau = \begin{bmatrix} \frac{1}{2} - \frac{\sqrt{3}}{6} \\ \frac{1}{2} + \frac{\sqrt{3}}{6} \\ 1 \end{bmatrix},\ w_1 = \begin{bmatrix} \frac{1}{4} & \frac{1}{4} - \frac{\sqrt{3}}{6} \\ \frac{1}{4} + \frac{\sqrt{3}}{6} & \frac{1}{4} \\ \frac{1}{2} & \frac{1}{2} \end{bmatrix}$$

We now derive explicit expressions for several multi-order Runge-Kutta methods. The conditions for order of consistency at stages $j \in [\![1, s]\!]$ here play a central role since, as stated in Section III.1.25, the order conditions at the final stage do not constrain the weights of the preceding stages at sufficiently high ranks.

To avoid repeatedly repeating technical details, we adopt the following conventions. Whenever we write "the only method," this is understood modulo permutation of the stages. When we refer to a method "of order of consistency" we mean at least of order of consistency for infinite length. Likewise, "convergent" means convergent up to infinite length.

The only convergent explicit node-determined multi-order Runge-Kutta method with 1 point of order of consistency 1 at all stages and all ranks :

$$\text{Multi-order explicit Euler method :}\quad \tau = \begin{bmatrix} 0 \\ 1 \end{bmatrix}, \forall N \in \mathbb{N}^*,\ w_N = \begin{bmatrix} 0 \\ 1 \end{bmatrix}$$

For node-determined multi-order Runge Kutta method with 2 points the possibilities are endless. The only convergent explicit node-determined multi-order Runge-Kutta methods of order of consistency 2 at all ranks and last stage, of order of consistency 1 at all ranks and all stages are :

$$\tau_2 \in \mathbb{R}^*,\ \tau = \begin{bmatrix} 0 \\ \tau_2 \\ 1 \end{bmatrix},\ \forall N \in \mathbb{N}^*,\ w_N = \begin{bmatrix} 0 & 0 \\ \tau_2^N & 0 \\ 1 - \frac{1}{\tau_2(1+N)} & \frac{1}{\tau_2(1+N)} \end{bmatrix}$$

To mimick the explicit midpoint method, we take $\tau_2 = \frac{1}{2}$ :

$$\text{Multi-order explicit midpoint method :}\quad \tau = \begin{bmatrix} 0 \\ \frac{1}{2} \\ 1 \end{bmatrix},\ \forall N \in \mathbb{N}^*,\ w_N = \begin{bmatrix} 0 & 0 \\ \frac{1}{2^N} & 0 \\ \frac{N-1}{1+N} & \frac{2}{1+N} \end{bmatrix}$$

Interestingly, the case $\tau_2 = \frac{1}{2}$ is given by the solved system conditions of order of consistency 3 at second rank at last stage. If instead we choose the solved system conditions of order of consistency 3 at first rank at last stage, we get :

$$\text{Multi-order Ralston's method :}\quad \tau = \begin{bmatrix} 0 \\ \frac{2}{3} \\ 1 \end{bmatrix},\ \forall N \in \mathbb{N}^*,\ w_N = \begin{bmatrix} 0 & 0 \\ \left(\frac{2}{3}\right)^N & 0 \\ \frac{2N-1}{2(1+N)} & \frac{3}{2(1+N)} \end{bmatrix}$$

For node-determined multi-order Runge-Kutta methods with 3 points we face an issue, convergent explicit multi-order Runge-Kutta methods can achieve an order of consistency 3 at last stage at only a single rank, we thus need to choose which rank. Fortunately, the choice is easy to make, it should be the first rank. No matter the order of the initial value problem, the method benefits from this choice in the precision of its approximations, especially with initial value problems of order 1. It is also the rank which is the least advantaged by the $N$ in $\underset{h\to 0}{\mathcal{O}}\left(h^{N+v_{n,j,N}}\right)$.

The only convergent explicit node-determined multi-order Runge-Kutta methods of order of consistency 3 at first rank and last stage, of order of consistency 2 at all ranks and last stage, of order of consistency 1 at all ranks and all stages, of order of consistency 2 at stage 3 and all ranks, which satisfy the solved system condition of order of consistency 3 at last stage and all ranks are :

$$\tau_2 \in \mathbb{R}^* \setminus \left\{\frac{2}{3}\right\},\ \tau = \begin{bmatrix} 0 \\ \tau_2 \\ \frac{2}{3} \\ 1 \end{bmatrix},\ \forall N \in \mathbb{N}^*,\ w_N = \begin{bmatrix} 0 & 0 & 0 \\ \tau_2^N & 0 & 0 \\ \left(\frac{2}{3}\right)^N \frac{\tau_2(1+N)-\frac{2}{3}}{\tau_2(1+N)} & \left(\frac{2}{3}\right)^N \frac{2}{3\tau_2(1+N)} & 0 \\ 1 + \frac{6-(2+N)(3\tau_2+2)}{2\tau_2(1+N)(2+N)} & \frac{2(N-1)}{\tau_2(1+N)(2+N)(2-3\tau_2)} & \frac{9(2-\tau_2(2+N))}{2(1+N)(2+N)(2-3\tau_2)} \end{bmatrix}$$

This leads to infinitely many possible methods. The solved system conditions of order of consistency 4 cannot be satisfied at any stage and rank, therefore, to imitate Heun's method, we choose $\tau_2 = \frac{1}{3}$, which finally yields :

$$\text{Multi-order Heun's method : } \quad \tau = \begin{bmatrix} 0 \\ \frac{1}{3} \\ \frac{2}{3} \\ 1 \end{bmatrix}, \ \forall N \in \mathbb{N}^*, \ w_N = \begin{bmatrix} 0 & 0 & 0 \\ \frac{1}{3^N} & 0 & 0 \\ \left(\frac{2}{3}\right)^N \frac{N-1}{1+N} & \left(\frac{2}{3}\right)^N \frac{2}{1+N} & 0 \\ 1 - \frac{9N}{2(1+N)(2+N)} & \frac{6(N-1)}{(1+N)(2+N)} & \frac{3(4-N)}{2(1+N)(2+N)} \end{bmatrix}$$

Explicit methods with 4 points present a difficulty. Since RK4 is the most well-known Runge-Kutta method, it is natural to attempt to generalize it. However, no fully satisfactory generalization exists. All convergent explicit node-determined multi-order Runge-Kutta methods of order of consistency 4 at last stage and first rank, of order of consistency 3 at last stage and all ranks, of order of consistency 1 at all stages and all ranks, with $w_1$ and $\tau$ equal to that of RK4, must satisfy the condition $w_{2,4,1} = 2w_{2,3,2}$, which prevents the method from having its characteristic diagonal form. If we add the conditions for order of consistency 2 at stage 4 and all ranks, we get :

$$w_{1,3,2} = \frac{1}{2}, \ w_{1,4,3} = 1, \ w_{2,3,2} = \frac{1}{6}, \ w_{2,4,3} \in \mathbb{R}, \ \forall N \in [\![3, +\infty[\![, \ \left(w_{N,3,2}, w_{N,4,3}\right) \in \mathbb{R}^2$$

$$\tau = \begin{bmatrix} 0 \\ \frac{1}{2} \\ \frac{1}{2} \\ 1 \\ 1 \end{bmatrix}, \ \forall N \in \mathbb{N}^*, \ w_N = \begin{bmatrix} 0 & 0 & 0 & 0 \\ \frac{1}{2^N} & 0 & 0 & 0 \\ \frac{1}{2^N} - w_{N,3,2} & w_{N,3,2} & 0 & 0 \\ \frac{N-1}{1+N} & \frac{2}{1+N} - w_{N,4,3} & w_{N,4,3} & 0 \\ \frac{N^2}{(1+N)(2+N)} & \frac{2N}{(1+N)(2+N)} & \frac{2N}{(1+N)(2+N)} & \frac{2-N}{(1+N)(2+N)} \end{bmatrix}$$

Since no additional solved system conditions can be satisfied, we need another criterion, though somewhat arbitrary.

### IV.1.1 Definition : Polynomial-type weights

Let $(\tilde{n}, s) \in \overline{\mathbb{N}}^* \times \mathbb{N}^*$, $M = (\tau, w) \in \mathrm{NDMORK}_{\tilde{n},s}$, $j \in [\![1, s+1]\!]$, $v' \in \mathbb{N}^*$ the highest integer such that, for all $N \in \mathbb{N}^*$ with $N \leq \tilde{n}$, $M$ satisfies the solved system conditions of order of consistency $v'$ at stage $j$ and rank $N$.

$M$ is said to have polynomial-type weights at stage $j$ if and only if :

$$\forall (N, j') \in \mathbb{N}^* \times [\![1, s]\!], \ N \leq \tilde{n}, \ w_{N,j,j'} = \tau_j^N \frac{N!}{(v' - 1 + N)!} P_{j,j'}(N)$$

Where, for all $j' \in [\![1, s]\!]$, $P_{j,j'}$ is a polynomial of degree $v' - 1$ or less.

Methods with polynomial-type weights seem to arise naturally, as seen with the previous methods derived. This may be due to the solved system conditions, which in matrix form can be written as :

$$\begin{bmatrix} 1 & 1 & \dots & 1 \\ \tau_1 & \tau_2 & \dots & \tau_s \\ \vdots & \vdots & \ddots & \vdots \\ \tau_1^{v'-1} & \tau_2^{v'-1} & \dots & \tau_{s'}^{v'-1} \end{bmatrix} w_{N,j} = \tau_j^N \frac{N!}{(v' - 1 + N)!} \begin{bmatrix} \tau_j^0 0! \frac{(v'-1+N)!}{(0+N)!} \\ \tau_j^1 1! \frac{(v'-1+N)!}{(1+N)!} \\ \vdots \\ \tau_j^{v'-1} (v' - 1)! \frac{(v'-1+N)!}{(v'-1+N)!} \end{bmatrix}$$

The entries of the vector on the right are polynomials of degree $v' - 1$, which implies that $w_{N,j}$ involves this kind of expression.

If we assume the method has polynomial-type weights at stage 3, then $w_{N,3,2} = \frac{1}{2^N}\frac{a+bN}{N+1}$. The two known values, $w_{1,3,1}$ and $w_{2,3,2}$, imply $w_{N,3,2} = \frac{1}{2^{N-1}(N+1)}$. If we assume the method has polynomial-type weights at stage 4, then $w_{N,4,3} = \frac{a+bN}{N+1}$. The only known value of $w_{N,4,3}$, $w_{1,4,3}$, gives $w_{N,4,3} = \frac{2-b+bN}{N+1}$. For the sake of having $w_{N,4,1} = -w_{N,4,2}$, we take $b = 1$. We get the method :

$$\text{MORK4 :} \quad \tau = \begin{bmatrix} 0 \\ \frac{1}{2} \\ \frac{1}{2} \\ 1 \\ 1 \end{bmatrix}, \ \forall N \in \mathbb{N}^*, \ w_N = \begin{bmatrix} 0 & 0 & 0 & 0 \\ \frac{1}{2^N} & 0 & 0 & 0 \\ \frac{1}{2^N}\frac{N-1}{1+N} & \frac{1}{2^{N-1}(1+N)} & 0 & 0 \\ \frac{N-1}{1+N} & \frac{1-N}{1+N} & 1 & 0 \\ \frac{N^2}{(1+N)(2+N)} & \frac{2N}{(1+N)(2+N)} & \frac{2N}{(1+N)(2+N)} & \frac{2-N}{(1+N)(2+N)} \end{bmatrix}$$

If instead we aim for the highest order of consistency possible, there also is an infinite number of possibilities. The convergent explicit node-determined multi-order Runge-Kutta methods of order of consistency 4 at last stage and first rank, of order of consistency 3 at last stage and all ranks, of order of consistency 1 at all stage and all ranks, of order of consistency 2 at stage 3 and all ranks, of order of consistency 2 at stage 4 and all ranks, of order of consistency 3 at stage 4 and rank 2, are :

$$w_{1,4,3} = 2, \ w_{2,4,3} = \frac{2}{3}, \ \forall N \in [\![3, +\infty[\![, \ w_{N,4,3} \in \mathbb{R}$$

$$\tau = \begin{bmatrix} 0 \\ \frac{1}{2} \\ \frac{1}{2} \\ 1 \\ 1 \end{bmatrix}, \ \forall N \in \mathbb{N}^*, \ w_N = \begin{bmatrix} 0 & 0 & 0 & 0 \\ \frac{1}{2^N} & 0 & 0 & 0 \\ \frac{N}{2^N(1+N)} & \frac{1}{2^N(1+N)} & 0 & 0 \\ \frac{N-1}{1+N} & \frac{2}{1+N} - w_{N,4,3} & w_{N,4,3} & 0 \\ \frac{N^2}{(1+N)(2+N)} & 0 & \frac{4N}{(1+N)(2+N)} & \frac{2-N}{(1+N)(2+N)} \end{bmatrix}$$

No additional solved system conditions can be satisfied. We assume the method has polynomial-type weights at stage 4, thus, for all $N \in \mathbb{N}^*$, $w_{N,4,3} = \frac{aN+b}{1+N}$. The values of $w_{N,4,3}$ at $N = 1$ and $N = 2$ yield :

$$\text{MORK4b :} \quad \tau = \begin{bmatrix} 0 \\ \frac{1}{2} \\ \frac{1}{2} \\ 1 \\ 1 \end{bmatrix}, \ \forall N \in \mathbb{N}^*, \ w_N = \begin{bmatrix} 0 & 0 & 0 & 0 \\ \frac{1}{2^N} & 0 & 0 & 0 \\ \frac{N}{2^N(1+N)} & \frac{1}{2^N(1+N)} & 0 & 0 \\ \frac{N-1}{1+N} & \frac{2(N-2)}{1+N} & \frac{2(3-N)}{1+N} & 0 \\ \frac{N^2}{(1+N)(2+N)} & 0 & \frac{4N}{(1+N)(2+N)} & \frac{2-N}{(1+N)(2+N)} \end{bmatrix}$$

The only convergent node-determined multi-order Runge-Kutta method with 1 point, $\tau_1 = 1$, of order of consistency 1 at all stages and all ranks :

$$\text{Multi-order implicit Euler method :} \quad \tau = \begin{bmatrix} 1 \\ 1 \end{bmatrix}, \ w_N = \begin{bmatrix} 1 \\ 1 \end{bmatrix}$$

The only convergent node-determined multi-order Runge-Kutta methods with 1 point of order of consistency 2 at stage 1 and first rank, of order of consistency 1 at all ranks and all stages :

$$\text{Multi-order implicit midpoint method :} \quad \tau = \begin{bmatrix} \frac{1}{2} \\ 1 \end{bmatrix}, \ w_N = \begin{bmatrix} \frac{1}{2^N} \\ 1 \end{bmatrix}$$

The average of the multi-order explicit Euler method and the multi-order implicit Euler :

$$\text{Multi-order Crank-Nicolson method :}\quad \tau = \begin{bmatrix} 0 \\ 1 \\ 1 \end{bmatrix},\ w_N = \begin{bmatrix} 0 & 0 \\ \frac{1}{2} & \frac{1}{2} \\ \frac{1}{2} & \frac{1}{2} \end{bmatrix}$$

An interesting aspect of the multi-order Crank Nicolson method is that it does not stem from considerations on its order of consistency. It is only of order of consistency 2 at all ranks and last stage, which is already achievable with the multi-order implicit midpoint method. Its true interest may show in its stability properties. If we intentionally avoid the conditions for order of consistency 4, we get another method that similar to the Crank-Nicolson method in its structure. The only convergent node-determined multi-order Runge-Kutta method with 2 points, of order of consistency 3 at first rank and last stage, of order of consistency 2 at all ranks and all stages, and $\tau_1 = 0$ :

$$\text{Multi-order CNb :}\quad \tau = \begin{bmatrix} 0 \\ \frac{2}{3} \\ 1 \end{bmatrix},\ w_N = \begin{bmatrix} 0 & 0 \\ \left(\frac{2}{3}\right)^N \frac{N}{1+N} & \left(\frac{2}{3}\right)^N \frac{1}{1+N} \\ \frac{2N-1}{2(1+N)} & \frac{3}{2(1+N)} \end{bmatrix}$$

The only convergent node-determined multi-order Runge-Kutta method with 2 points, of order of consistency 4 at last stage and rank 1, of order of consistency 3 at last stage and rank 2, of order of consistency 2 at all stages and all ranks :

$$\text{Gauss-Jacobi method of order 4 :}$$

$$\tau = \begin{bmatrix} \frac{1}{2} - \frac{\sqrt{3}}{6} \\ \frac{1}{2} + \frac{\sqrt{3}}{6} \\ 1 \end{bmatrix},\ w_N = \begin{bmatrix} \frac{1}{1+N}\left(1 + \left(1+\sqrt{3}\right)\frac{N}{2}\right)\left(\frac{1}{2} - \frac{\sqrt{3}}{6}\right)^N & -\sqrt{3}\frac{N}{1+N}\left(\frac{1}{2} - \frac{\sqrt{3}}{6}\right)^{N+1} \\ \sqrt{3}\frac{N}{1+N}\left(\frac{1}{2} + \frac{\sqrt{3}}{6}\right)^{N+1} & \frac{1}{1+N}\left(1 + \left(1-\sqrt{3}\right)\frac{N}{2}\right)\left(\frac{1}{2} + \frac{\sqrt{3}}{6}\right)^N \\ \frac{1}{2} + \sqrt{3}\frac{N-1}{2(1+N)} & \frac{1}{2} - \sqrt{3}\frac{N-1}{2(1+N)} \end{bmatrix}$$

To find experimentally the order of consistency of a method, we consider an initial value problem of order $n \in \mathbb{N}^*$ we know the solution of, and plot the absolute value of the errors at each rank after a single step as a function of the step size $h$, with a logarithmic scale. The graph for each rank $N \in [\![1, n]\!]$ should be similar to a linear function. We remind the reader that the derivative index $i$ is related to the rank $N$ by $i = n - N$. The slope of those lines minus $N$ is the experimental order of consistency $v'_N$ of the method at rank $N$ and last stage. The logarithmic scale implies that the error is $\underset{h \to 0}{\mathcal{O}}\left(h^{N+v'_N}\right)$. If $f$ is a function of class $C^{\max(1, v'_N+1)}$, then, by definition, $v'_N$ is an upper bound of the actual order of consistency.

One point to keep in mind is that the principle behind the order of consistency is to cancel as many coefficients in the Taylor expansion of $f$ around $\overline{\mathcal{J}}^{n-1}\hat{y}(t_0)$. However, if some coefficients of the expansion are already zero, the experimental order of consistency can be greater than the theoretical order of consistency.

To test these methods we use a category of initial value problems that is easy to deal with, linear initial value problem, thus initial value problems with differential equation functions of the form :

$$\forall (t, x) \in \mathbb{R} \times \mathbb{R}^{[\![1,n]\!]},\ f(t, x) = \sum_{N=1}^{n} \alpha_N x_N$$

With $\alpha \in \mathbb{C}^{[\![1,n]\!]}$. This category of initial value problems is the object of study of the next section. The solution of such initial value problems depend entirely on the roots of the characteristic polynomial :

$$P_\alpha = X^n - \sum_{i=1}^{n} \alpha_N X^{n-N}$$

We here consider the particular case where the characteristic polynomial has a single root $\lambda$, hence :

$$P_\alpha = (X-\lambda)^n = \sum_{k=0}^{n} \binom{n}{k} (-\lambda)^k X^{n-k} = X^n + \sum_{k=1}^{n} \binom{n}{k} (-\lambda)^k X^{n-k}$$

We deduce :

$$\forall (t,x) \in \mathbb{R} \times \mathbb{R}^{[\![1,n]\!]},\ f(t,x) = -\sum_{N=1}^{n} \binom{n}{N} (-\lambda)^N x_N$$

The solutions are of the form $\hat{y}(t) = e^{\lambda(t-t_0)} P(t-t_0)$ with $P$ a polynomial of at most degree $n-1$. Let $a \in \mathbb{C}^{[\![0,n-1]\!]}$ the coefficients of $P$. The generalized Leibniz rule yields, for all $i \in [\![0, n-1]\!]$ :

$$\begin{aligned} \frac{\mathrm{d}^i \hat{y}}{\mathrm{d}t^i}(t) &= \sum_{k=0}^{i} \binom{i}{k} P^{(k)}(t-t_0) \lambda^{i-k} e^{\lambda(t-t_0)} \\ &= e^{\lambda(t-t_0)} \sum_{k=0}^{i} \binom{i}{k} \lambda^{i-k} P^{(k)}(t-t_0) \end{aligned}$$

We now find the coefficients of $P$ in terms of the initial values. We have $P(t) = e^{-\lambda t} \hat{y}(t+t_0)$, thus the generalized Leibniz rule yields :

$$\frac{\mathrm{d}^i P}{\mathrm{d}t^i}(t) = \sum_{k=0}^{i} \binom{i}{k} (-\lambda)^{i-k} e^{-\lambda t} \frac{\mathrm{d}^k \hat{y}}{\mathrm{d}t^k}(t+t_0)$$

Thus :

$$\frac{\mathrm{d}^i P}{\mathrm{d}t^i}(0) = \sum_{k=0}^{i} \binom{i}{k} (-\lambda)^{i-k} \frac{\mathrm{d}^k \hat{y}}{\mathrm{d}t^k}(t_0) \Leftrightarrow a_i = \sum_{k=0}^{i} \frac{(-\lambda)^{i-k}}{k!(i-k)!} y_{n,0,n-k}$$

With these formulas, we know the solution no matter the initial instant and initial values. We now plot the real part of the error after the first step of the methods in the case where $\lambda = -\frac{1}{2} + \boldsymbol{i}$, $n = 6$, $t_0 = 0$ and $P$ is identically 1. The considered initial value problem is hence :

$$\frac{\mathrm{d}^6 \hat{y}}{\mathrm{d}t^6} = 6\lambda \frac{\mathrm{d}^5 \hat{y}}{\mathrm{d}t^5} - 15\lambda^2 \frac{\mathrm{d}^4 \hat{y}}{\mathrm{d}t^4} + 20\lambda^3 \frac{\mathrm{d}^3 \hat{y}}{\mathrm{d}t^3} - 15\lambda^4 \frac{\mathrm{d}^2 \hat{y}}{\mathrm{d}t^2} + 6\lambda^5 \frac{\mathrm{d}^1 \hat{y}}{\mathrm{d}t^1} - \lambda^6 \frac{\mathrm{d}^0 \hat{y}}{\mathrm{d}t^0}$$

$$\forall i \in [\![0,5]\!],\ \frac{\mathrm{d}^i \hat{y}}{\mathrm{d}t^i}(0) = \lambda^i$$

The solution to this problem is $t \to e^{\lambda t}$. To save space, we note $y_N$ in the legends instead of $y_{n,s+1,N}$.

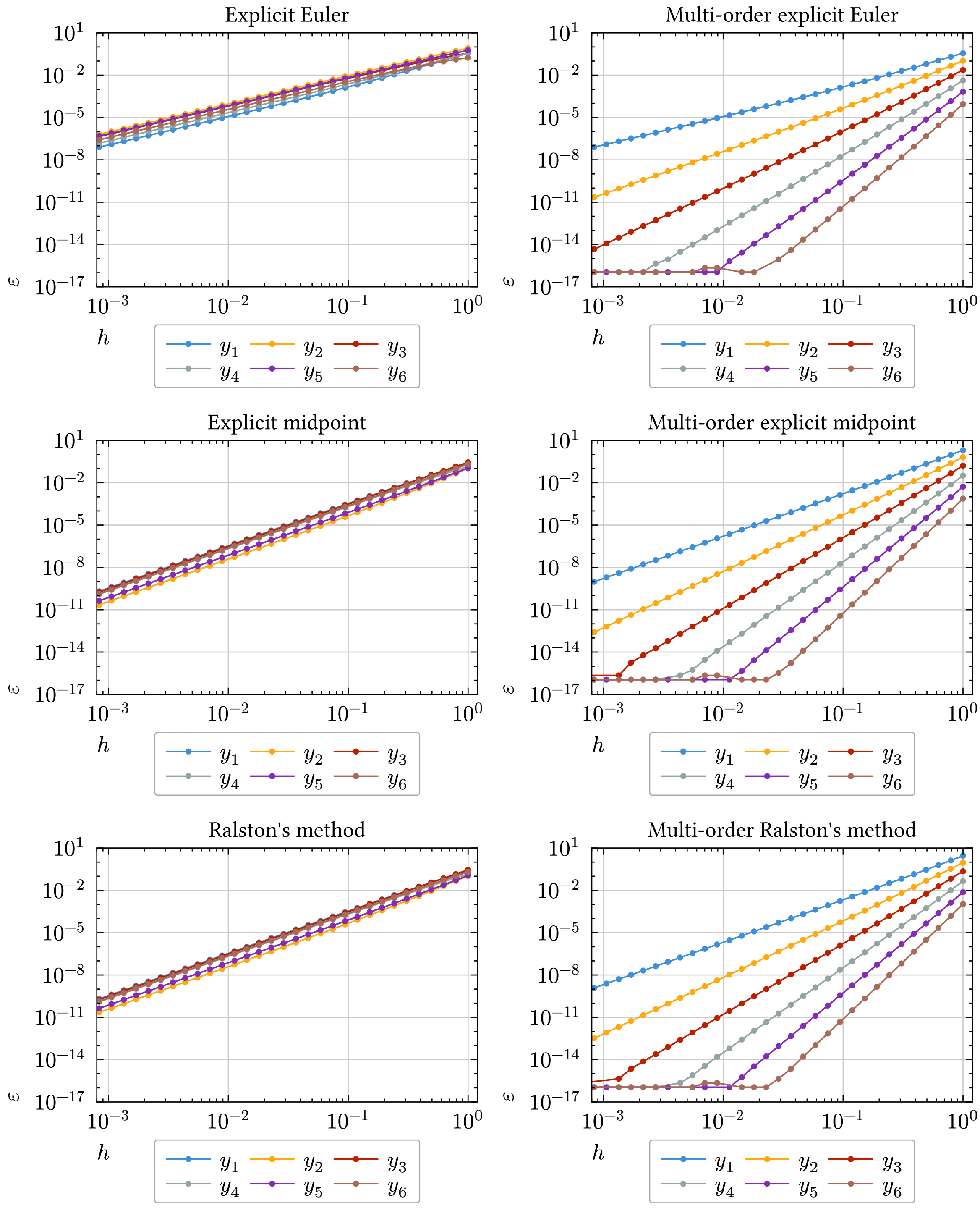

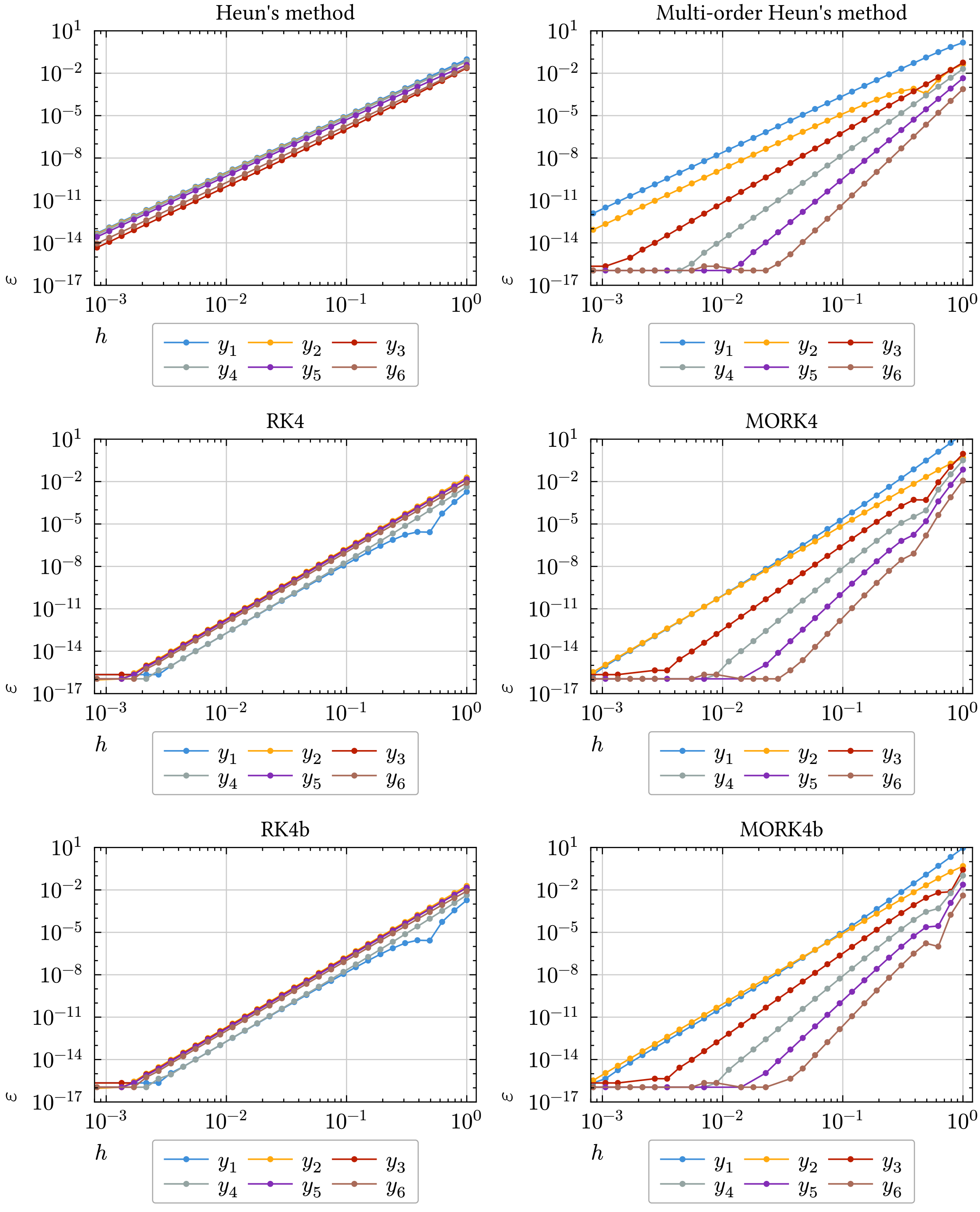

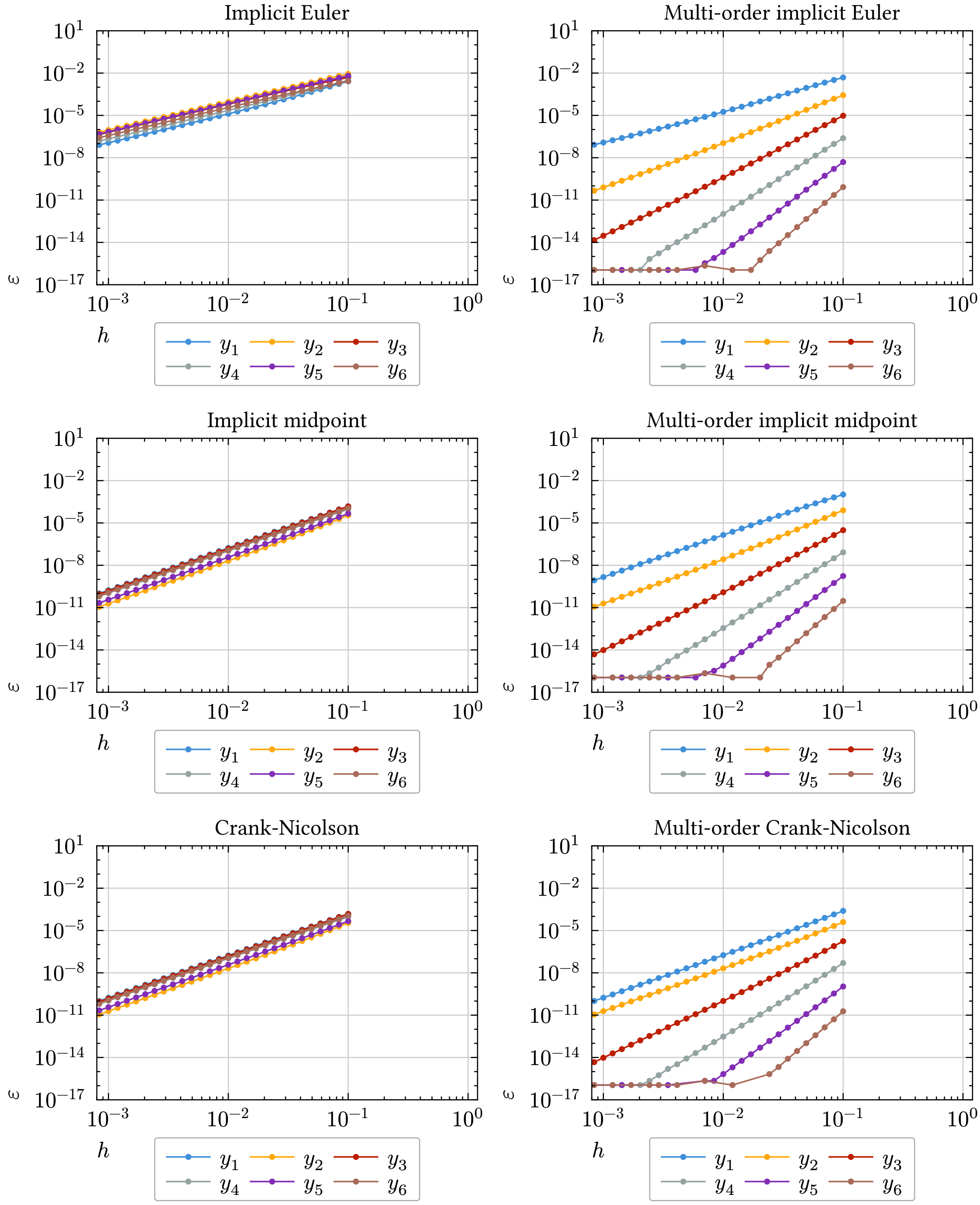

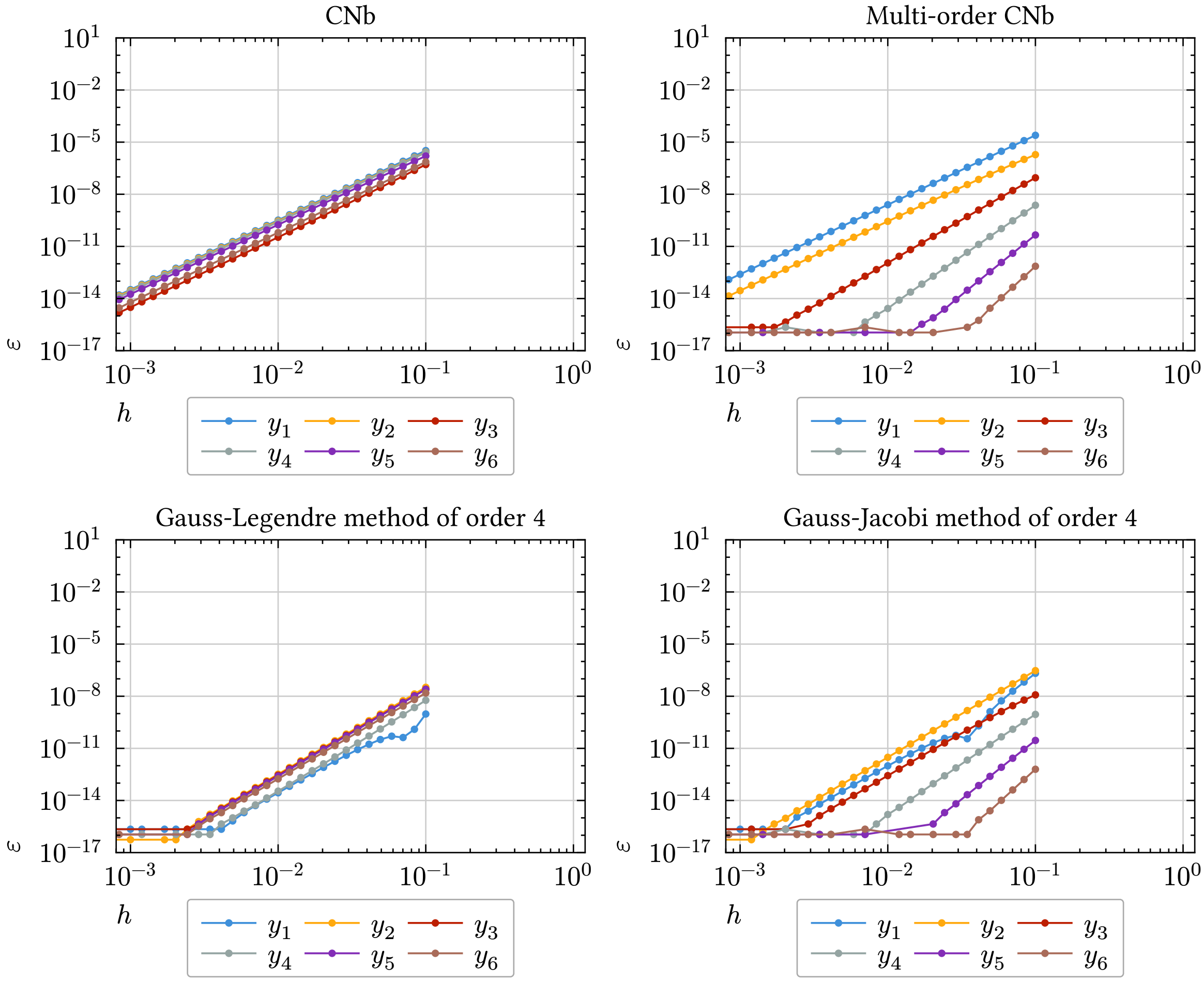

We did not plot the approximations of implicit methods for $h > 0.1$ because they all behave erratically with such step sizes, probably because Picard's iterations do not converge then. We can first observe that some curves begin with a flat line before turning into a linear function, this is due to the rounding errors made by the computer. The minimum error is around $10^{-16}$, which is consistent with the advertised 15 to 17 significant decimal digits precision of the 64 bits representation of floats we are using[11]. There also are missing points in the same flat regions, at those points the computer found the error to be 0, they hence had to be filtered out since we are using a logarithmic scale.

The experiments confirm the theoretical properties, multi-order methods are more accurate than their Runge-Kutta counterpart, except for the first rank where their errors are similar. For example, with $h = 0.1$, the error of approximation of the solution $\hat{y}$ of the multi-order explicit Euler method is a billion times smaller than the error of the explicit Euler method.

Despise this greater accuracy, Runge-Kutta methods are sometimes able to outperform multi-order Runge Kutta methods on the long term. For example, we now use the constant step size $h = \frac{1}{2}$ and plot the six first steps of each method. We first plot the solution at all points of the time mesh.

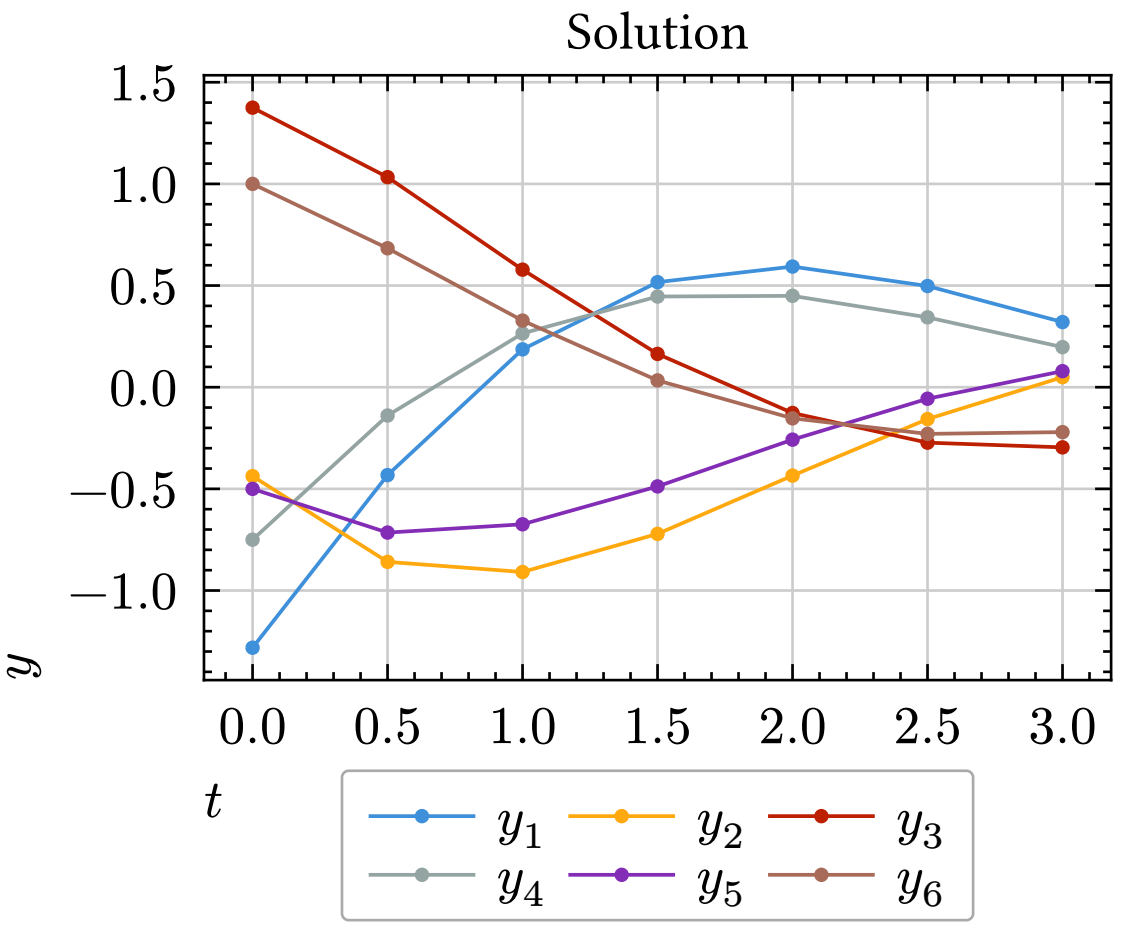

We now examine how explicit methods approximate this initial value problem.

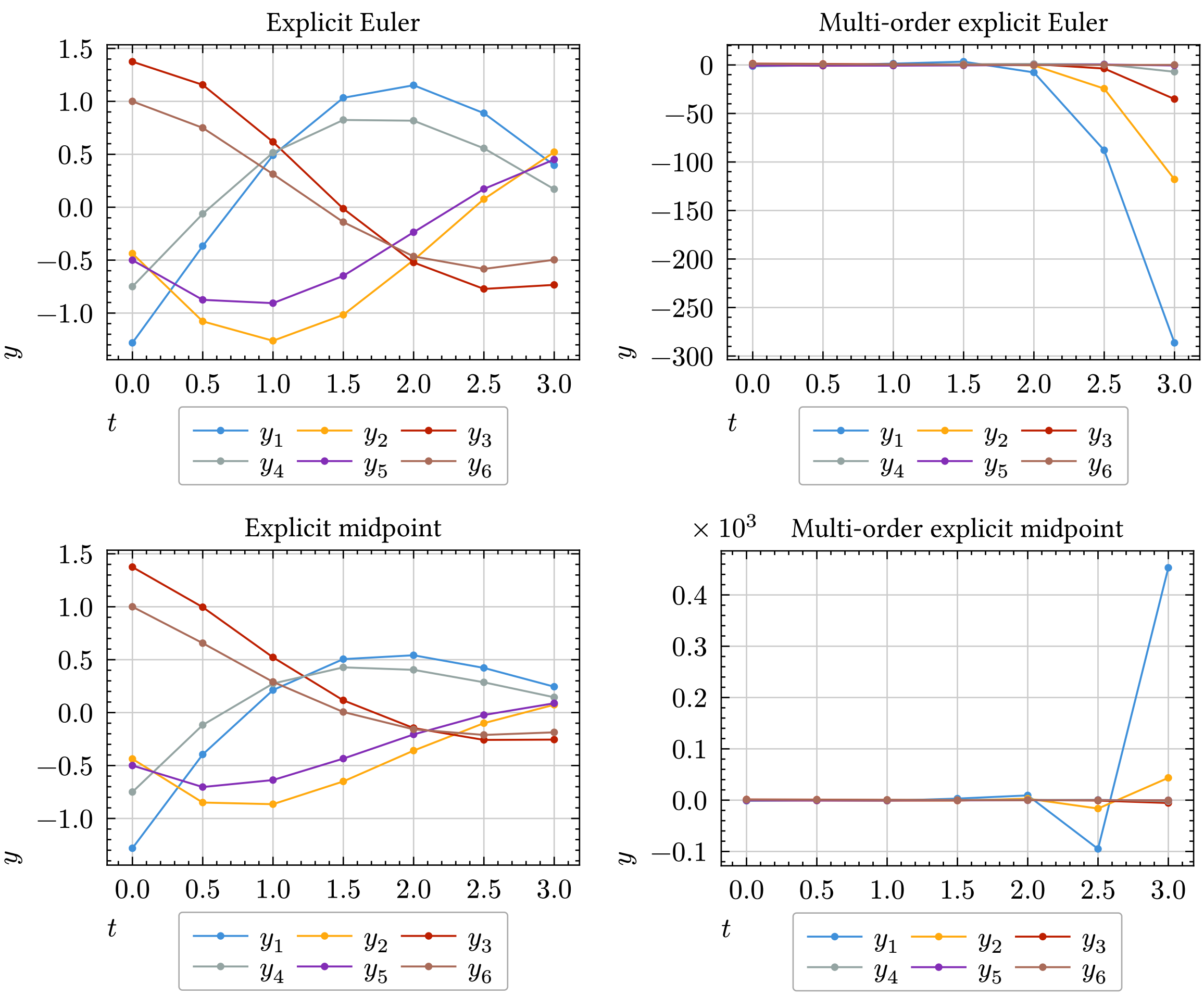

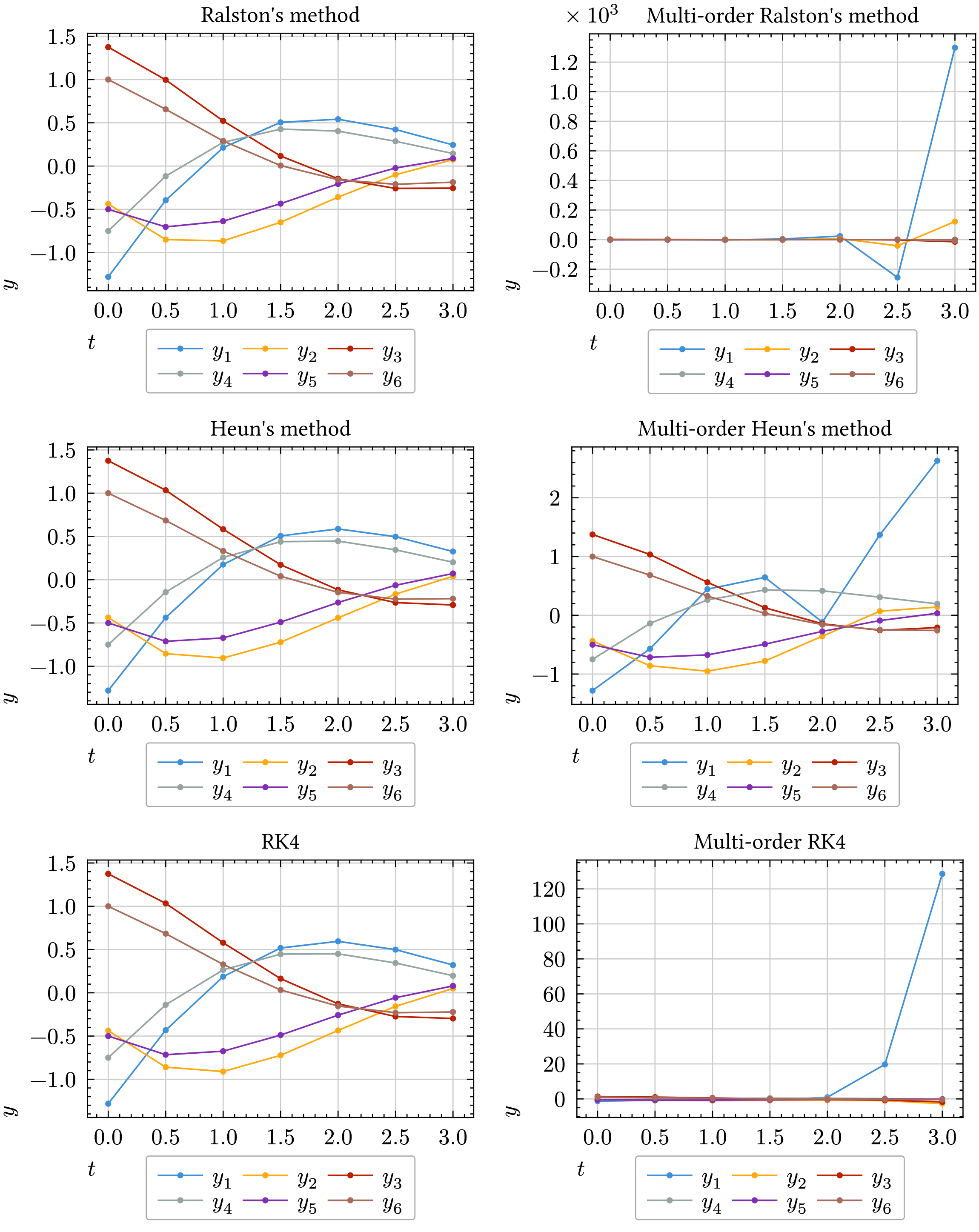

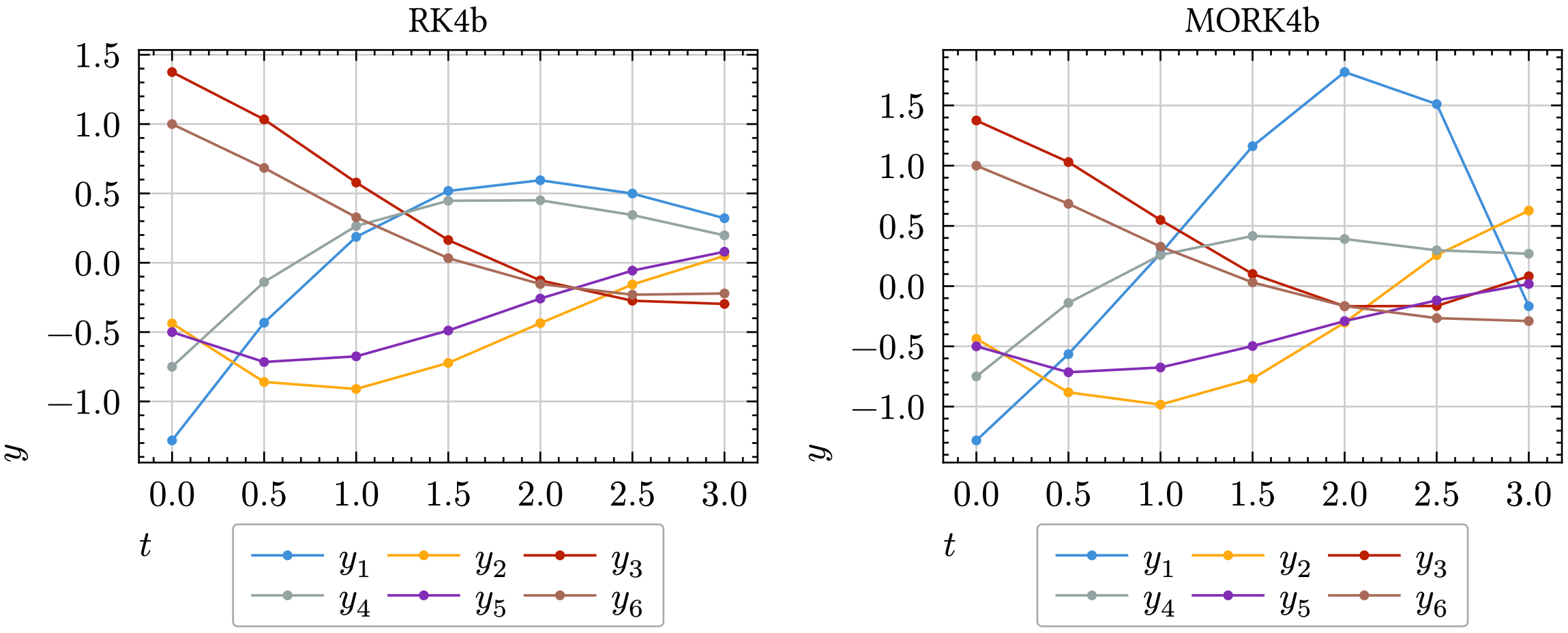

Here, the Runge-Kutta methods have better approximations. We did not plot the approximations of implicit methods because, as mentioned before, they have an erratic behavior with this step size. But if we take $h \leq 0.1$, multi-order Runge-Kutta methods approximate better the solution. For example, we now plot the error of four methods after 30 steps, with the constant step size $h$ on the x-axis :

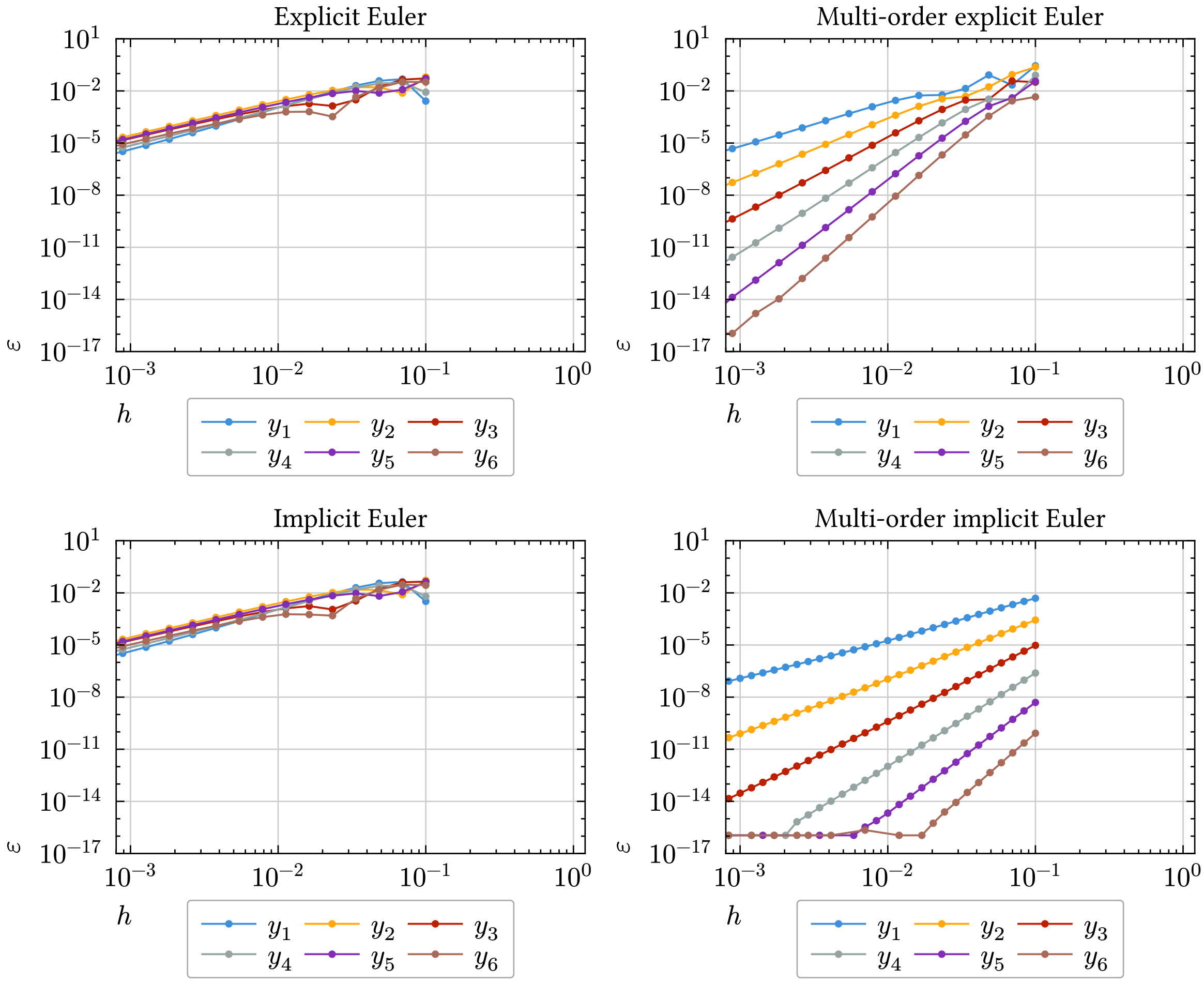

The order of consistency is hence an efficient criterion, but it requires small step sizes if one wishes to approximate the solution over long time intervals. The listed multi-order Runge-Kutta method were selected solely on the basis of their order of consistency, which describes only the local behavior of the approximation. It does not ensure long-time accuracy.

Focusing exclusively on a high order of consistency may even lead to poorer global behavior. The principle behind order of consistency is that the approximations matches a Taylor expansion up to a given degree. While higher-order Taylor expansions yield better local approximations, they do not necessarily behave well far from the expansion point. For bounded functions, higher-degree Taylor expansions diverge more rapidly as $t$ goes to $\pm\infty$ since higher powers of $t$ dominate. In this sense, increasing the order of a Taylor approximation can worsen the approximation far from the initial instant.

One way to adress this issue is to use Padé approximants[12] instead of Taylor polynomials. Padé approximants generalize Taylor expansions by using rational functions instead of polynomials as approximations. This allows the approximation to exhibit different global behavior, such as convergence as $t$ tends towards $\pm\infty$. Fortunately, Runge-Kutta methods and multi-order Runge-Kutta methods are already capable of producing Padé-type approximations.

To illustrate this, consider the linear initial value problem of order 1 :

$$\frac{\mathrm{d}\hat{\mathrm{y}}}{\mathrm{d}t} = \lambda \hat{\mathrm{y}}$$

The implicit euler method yields :

$$\hat{\mathrm{y}}(t_0 + h) \approx \frac{1}{1 - h\lambda} \hat{\mathrm{y}}(t_0)$$

Which is a rational function, not a polynomial. Linear test problems are useful for detecting whether a method uses Padé-type approximants. This observation leads naturally to the definition of A-stability, which we study in the next subsection.

Although the order of consistency is an important measure of local precision, stability is equally crucial. This is well established in the litterature on Runge-Kutta methods. For instance, Gauss-Legendre methods[1] are Runge-Kutta methods which achieve the maximal order of consistency $2s$, yet Radau or Lobatto[1] methods are also used because of their stronger stability properties, despite having lower order of consistency.

Since multi-order methods were partly motivated by extending Taylor-type constructions to higher orders, one might ask whether this emphasis on order is misguided. However, multi-order methods offer additional flexibility. They generalize classical Runge-Kutta methods and therefore provide greater freedom in the expression of the approximations.

An intresting example is the semi-implicit Euler method[3], a method widely used for hamiltonians systems. Applied to second-order initial value problems, it coincides with the node-determined multi-order Runge-Kutta method of length 2 $(\tau, w)$, where :

$$\tau = \begin{bmatrix} 0 \\ 1 \end{bmatrix}, \ w_1 = \begin{bmatrix} 0 \\ 1 \end{bmatrix}, \ w_2 = \begin{bmatrix} 0 \\ 2 \end{bmatrix}$$

Although this method is only of order of consistency 0 at last stage and rank 2, it remains extremely popular, especially in physics engines. This is due to two key properties, its simplicity and its symplecticity. A symplectic method preserves the geometric structure of Hamiltonian flows, which in practice translates into good long-time energy behavior. In this case, symplecticity matters more than high order.

Since the semi-implicit method cannot be expressed as a rewritten Runge-Kutta method, it should be regarded as a genuinely multi-order method. It is therefore plausible that there exists other multi-order methods whose orders of consistency match or exceed those of classical Runge-Kutta methods while also possessing additional desirable properties.

For instance, the MORK4b method and the multi-order Heun's method appear to perform well on the previous problem. This claim is supported by the fact that most of the weights of the method were not chosen to increase the order of consistency at last stage. As stated in Section III.1.25, the values of the weights at stages $j \in [\![1, s]\!]$ do not influence the order of consistency of the last stage at high enough ranks. This additional freedom may allow the construction of more stable multi-order methods.

The next subsection is therefore devoted to the study of the stability properties of multi-order methods.

## IV.2 A-stability, L-stability & Half-line stability

The interest of linear initial value problems stems from the linearization of differential equation functions. Let $\left(f, t_0, y_{n,0}\right)$ an initial value problem with $f$ differentiable, and $\hat{y}$ a solution candidate. For $h \in \mathbb{R}$ small enough, we approximate :

$$
\begin{aligned}
f\left(\overline{\mathcal{J}}^{n-1}\hat{y}(t+h)\right) \approx\ & f\left(\overline{\mathcal{J}}^{n-1}\hat{y}(t)\right) + \frac{\partial f}{\partial t}\left(\overline{\mathcal{J}}^{n-1}\hat{y}(t)\right)h \\
& + \sum_{i=0}^{n-1}\sum_{k=1}^{d}\frac{\partial f}{\partial x_{n-i,k}}\left(\overline{\mathcal{J}}^{n-1}\hat{y}(t)\right)\left(\frac{\mathrm{d}^i\hat{y}_k}{\mathrm{d}t^i}(t) - \frac{\mathrm{d}^i\hat{y}_k}{\mathrm{d}t^i}(t_0)\right) \\
\approx\ & f\left(\overline{\mathcal{J}}^{n-1}\hat{y}(t)\right) + \frac{\partial f}{\partial t}\left(\overline{\mathcal{J}}^{n-1}\hat{y}(t)\right)h \\
& + \sum_{i=0}^{n-1}\begin{bmatrix} \frac{\partial f_1}{\partial x_{n-i,1}}\left(\overline{\mathcal{J}}^{n-1}\hat{y}(t)\right) & \cdots & \frac{\partial f_1}{\partial x_{n-i,d}}\left(\overline{\mathcal{J}}^{n-1}\hat{y}(t)\right) \\ \vdots & \ddots & \vdots \\ \frac{\partial f_d}{\partial x_{n-i,1}}\left(\overline{\mathcal{J}}^{n-1}\hat{y}(t)\right) & \cdots & \frac{\partial f_d}{\partial x_{n-i,d}}\left(\overline{\mathcal{J}}^{n-1}\hat{y}(t)\right) \end{bmatrix}\left(\frac{\mathrm{d}^i\hat{y}}{\mathrm{d}t^i}(t) - \frac{\mathrm{d}^i\hat{y}}{\mathrm{d}t^i}(t_0)\right)
\end{aligned}
$$

So we may want to consider initial value problems $\left(f_l, t_0, y_{n,0}\right)$ with $f_l$ of the form :

$$
\forall (t,x) \in \mathbb{R} \times \mathbb{R}^{[\![1,n]\!]\times[\![1,d]\!]},\ f_l(t,x) = C_1 + C_2 t + \sum_{i=1}^{n} A_N x_N
$$

And $(C_1, C_2, A) \in \mathbb{R}^{[\![1,d]\!]} \times \mathbb{R}^{[\![1,d]\!]} \times \mathbb{R}^{[\![1,n]\!]\times[\![1,d]\!]^2}$. Solutions to this initial value problem are of the form $\hat{y}_p + \hat{y}_h$ where $\hat{y}_p$ is a particular solution of the differential equation induced by $f_l$, and $\hat{y}_h$ is a solution of the associated homogeneous initial value problem $\left(f_h, t_0, y_{n,0} - \mathcal{J}^{n-1}\hat{y}_p(t_0)\right)$, with :

$$
\forall (t,x) \in \mathbb{R} \times \mathbb{R}^{[\![1,n]\!]\times[\![1,d]\!]},\ f_h(t,x) = \sum_{i=1}^{n} A_N x_N
$$

We are here only interested in the homogeneous initial value problem. To simplify this problem we need all $A_N$ to be diagonalizable in the same basis, or equivalently, there exists $P \in \mathbb{R}^{[\![1,d]\!]^2}$ an invertible matrix such that, for all $N \in [\![1,n]\!],\ A_N = P^{-1}B_N P$ where $B \in \mathbb{C}^{[\![1,n]\!]\times[\![1,d]\!]^2}$ is a vector of diagonal matrices. This property is named simultaneous diagonalization, and it is equivalent to, all the entries of $A$ are diagonalizable and they commute. If we assume they are indeed simultaneously diagonalizable, we have, for all $t \in \mathbb{R}$ :

$$
\begin{aligned}
\frac{\mathrm{d}^n\hat{y}}{\mathrm{d}t^n}(t) = f_h\left(\overline{\mathcal{J}}^{n-1}\hat{y}(t)\right) \Leftrightarrow\ & \frac{\mathrm{d}^n\hat{y}}{\mathrm{d}t^n}(t) = \sum_{N=2}^{n} P^{-1}B_N P\frac{\mathrm{d}^{n-N}\hat{y}}{\mathrm{d}t^{n-N}}(t) \\
\Leftrightarrow\ & \frac{\mathrm{d}^n P\hat{y}}{\mathrm{d}t^n}(t) = \sum_{N=1}^{n} B_N\frac{\mathrm{d}^{n-N}P\hat{y}}{\mathrm{d}t^{n-N}}
\end{aligned}
$$

Each entry of $P\hat{y}$ is the solution of an independent homogeneous constant linear initial value problem.

### IV.2.1 Definition : Homogeneous constant linear initial value problems

Let $n \in \mathbb{N}^*$, and $f$ a differential equation function of order $n$.

- $f$ is said to be a homogeneous constant linear differential equation function of order $n$ and dimension 1 if and only if $\mathbf{U} = \mathbb{R} \times \mathbb{C}^{[\![1,n]\!]}$ and there exists $\alpha \in \mathbb{C}^{[\![1,n]\!]}$ such that :

$$\forall (t,x) \in \mathbf{U},\ f(t,x) = \sum_{N=1}^{n} \alpha_N x_N$$

Let $\left(t_0, y_{n,0}\right) \in \mathbf{U}$.

- $\left(f, t_0, y_{n,0}\right)$ is said to be a homogeneous constant linear initial value problem of order $n$ and dimension 1 if and only if $f$ is a homogeneous constant linear differential equation function of order $n$ and dimension 1.

Since we linearized $f$ to get this system, a solution of one of the linear initial value problems indicates how an entry of $P\hat{\mathrm{y}}$ behaves locally, whether it decays, grows or oscillates, and we want this behavior to show in the approximations of a method. The local behavior of $\left(f, t_0, y_{n,0}\right)$ dictates how the solution of $\left(f_h, t_0, y_{n,0} - \mathcal{J}^{n-1}\hat{\mathrm{y}}_p(t_0)\right)$ behaves asymptotically, and we can observe how a method reacts by studying how it approximates the solutions of the initial value problems of $P\hat{\mathrm{y}}$, in particular the asymptotical behavior of the sequences of approximation.

We study two behaviors, when the solution converges towards 0 (absolute A-stability) and when it stays bounded (A-stability). We are only interested in constant step size sequences since we want to make sure that the local behavior of the solution is taken into account for a single step.

Technically, general multi-order Runge-Kutta methods are defined only for real valued equations, but $\mathbb{C}$ can be seen as $\mathbb{R}^2$, and the transition from one to the other is straightforward. The Cauchy-Lipschitz theorem gives that there exists a unique solution to this initial value problem. To find the solution of a homogeneous linear problem, we use the characteristic polynomial $P_\alpha$ defined as :

$$P_\alpha = X^n - \sum_{i=0}^{n-1} \alpha_{n-i} X^i$$

Let $p \in [\![1,n]\!]$ the number of different roots of $P_\alpha$ and $\tilde{r} \in \mathbb{C}^{[\![1,p]\!]}$ these roots, with $m \in \mathbb{N}^{*[\![1,p]\!]}$ their multiplicity. The solution of this initial value problem is a linear combination of the functions $t \to t^k e^{\tilde{r}_l t}$, with $l \in [\![1,p]\!],\ k \in [\![0, m_l - 1]\!]$. Each choice of coefficients in this linear combination corresponds to a unique set of initial values.

If the solution converges then it converges towards 0. We deduce that the solution stays bounded (tends towards 0) as $t$ tends towards $+\infty$ for all initial value if and only if all roots are negative (strictly negative), and the solution stays bounded (tends towards 0) as $t$ tends toward $-\infty$ for all initial values if and only if all roots are positive (strictly positive). To avoid distinguishing between the cases $+\infty$ and $-\infty$, since we approach $\pm\infty$ by infinitely adding a constant $h \in \mathbb{R}^*$, those conditions simplifies to, for all $k \in [\![1,n]\!],\ \Re(hr_k) \leq 0$ for boundedness of the solution, and, for all $k \in [\![1,n]\!],\ \Re(hr_k) < 0$ for convergence towards 0 of the solution.

The multiplicity of the roots does not interest us, we thus define $r \in \mathbb{C}^{[\![1,n]\!]}$ the roots of $P_\alpha$ without their multiplicity. We want to be able to directly manipulate the roots, not the coefficients $\alpha$. To achieve this we use Vieta's formula [13].

**IV.2.2 Theorem : Vieta's formula**

Let $P = \sum_{i=0}^{n} c_i X^i$ a polynomial of order $n \in \mathbb{N}^*$, with $c \in \mathbb{C}^{[\![0,n]\!]}$, $c_n \neq 0$, and $r' \in \mathbb{C}^{[\![1,n]\!]}$ its roots.

$$\forall N \in [\![1, n]\!],\ c_{n-N} = c_n (-1)^N \sum_{1 \leq i_1 < \ldots < i_N \leq n} \prod_{j=1}^{N} r'_{i_j}$$

In this particular case we have :

$$c_n = 1,\ \forall N \in [\![1, n]\!],\ c_{n-N} = -\alpha_N$$

Therefore :

$$\forall N \in [\![1, n]\!],\ \alpha_N = (-1)^{N+1} \sum_{1 \leq i_1 < \ldots < i_N \leq n} \prod_{k=1}^{N} r_{i_k}$$

Though we do not yet use this formula for the sake of concision. We have all the theory needed on linear initial value problems, we now state the definition of A-stability and absolute A-stability.

**IV.2.3 Definition : A-stability & Absolute A-stability**

Let $(\tilde{n}, s) \in \overline{\mathbb{N}}^* \times \mathbb{N}^*$, $M \in \text{GMORK}_{\tilde{n},s}$, $n \in \mathbb{N}^*$ such that $n \leq \tilde{n}$, $Y_n$ the approximation function for length $n$ of $M$, and $\mathcal{H}$ the set of step sizes of $M$.

- $M$ is said to be A-stable up to length $n$ if and only if, for all initial value problem $\left(f, t_0, y_{n,0}\right)$ of the form :

$$\alpha \in \mathbb{C}^{[\![1,n]\!]},\ \left(t_0, y_{n,0}\right) \in \mathbb{R} \times \mathbb{C}^{[\![1,n]\!]}$$

$$\forall (t, x) \in \mathbb{R} \times \mathbb{C}^{[\![1,n]\!]},\ f(t, x) = \sum_{i=1}^{n} \alpha_N x_N$$

  For all non zero constant step size sequence $(h)_{q\in\mathbb{N}} \in \mathcal{H}_{\infty,f,t_0,y_{n,0}}$, if we note $r \in \mathbb{C}^{[\![1,n]\!]}$ the roots of the characteristic polynomial $X^n - \sum_{i=0}^{n-1} \alpha_{n-i} X^i$ and if, for all $k \in [\![1, n]\!]$, $\Re(hr_k) \leq 0$, then $\left(Y_{n,f,q(s+1)}\left(t_0, y_{n,0}, (h)_{q\in\mathbb{N}}\right)\right)_{q\in\mathbb{N}}$ is bounded.

- $M$ is said to be absolute A-stable for length $n$ if and only if, for all initial value problem $\left(f, t_0, y_{n,0}\right)$ of the form :

$$\alpha \in \mathbb{C}^{[\![1,n]\!]},\ \left(t_0, y_{n,0}\right) \in \mathbb{R} \times \mathbb{C}^{[\![1,n]\!]}$$

$$\forall (t, x) \in \mathbb{R} \times \mathbb{C}^{[\![1,n]\!]},\ f(t, x) = \sum_{i=1}^{n} \alpha_N x_N$$

  For all non zero constant step size sequence $(h)_{q\in\mathbb{N}} \in \mathcal{H}_{\infty,f,t_0,y_{n,0}}$, if we note $r \in \mathbb{C}^{[\![1,n]\!]}$ the roots of the characteristic polynomial $X^n - \sum_{i=0}^{n-1} \alpha_{n-i} X^i$ and if, for all $k \in [\![1, n]\!]$, $\Re(hr_k) < 0$ , then $\left(Y_{n,f,q(s+1)}\left(t_0, y_{n,0}, (h)_{q\in\mathbb{N}}\right)\right)_{q\in\mathbb{N}}$ tends towards $0_n$.

We now find out how a general multi-order Runge-Kutta method approximates the solution of linear initial value problems of order $n$.

To do so, we first rewrite the stage system in two forms, its matrix form and its evaluation form. We then solve the stage system using its evaluation form. Finally, we use the expression of the solution in the matrix form of the stage system.

**IV.2.4 Lemma : Stage system - Matrix form**

Let $(\tilde{n}, s) \in \overline{\mathbb{N}}^* \times \mathbb{N}^*$, $M = (\tau, w, \tilde{w}) \in \text{GMORK}_{\tilde{n},s}$, and $n \in \mathbb{N}^*$ such that $n \leq \tilde{n}$. We define :

$$W_{n,j} = \begin{bmatrix} w_{1,j,1} & \dots & w_{1,j,s} \\ \vdots & \ddots & \vdots \\ w_{n,j,1} & \dots & w_{n,j,s} \end{bmatrix} \in \mathbb{R}^{[\![1,n]\!] \times [\![1,s]\!]}, \ \widetilde{W}_{n,j} = \begin{bmatrix} \tilde{w}_{1,1,j} & 0 & \dots & 0 \\ \tilde{w}_{2,1,j} & \tilde{w}_{2,2,j} & \ddots & \vdots \\ \vdots & \vdots & \ddots & 0 \\ \tilde{w}_{n,1,j} & \tilde{w}_{n,2,j} & \dots & \tilde{w}_{n,n,j} \end{bmatrix} \in \mathbb{R}^{[\![1,n]\!]^2}$$

$$H_n : \mathbb{R} \to \mathbb{R}^{[\![1,n]\!]^2}, \ h \to \begin{bmatrix} h^0 & 0 & \dots & 0 \\ 0 & h^1 & \ddots & \vdots \\ \vdots & \ddots & \ddots & 0 \\ 0 & \dots & 0 & h^{n-1} \end{bmatrix}, \ D_n = \begin{bmatrix} \frac{1}{1!} & 0 & \dots & 0 \\ 0 & \frac{1}{2!} & \ddots & \vdots \\ \vdots & \ddots & \ddots & 0 \\ 0 & \dots & 0 & \frac{1}{n!} \end{bmatrix} \in \mathbb{R}^{[\![1,n]\!]^2}$$

Let $f$ a differential equation function of order $n$, and $\left(t, y_{n,0}, h\right) \in \mathbf{U} \times \mathbb{R}$.
If $h = 0$, the stage system $\left(M, f, t, y_{n,0}, h\right)$ is equivalent to :

$$\forall j \in [\![1, s+1]\!], \ \begin{bmatrix} y_{n,j,1} \\ \vdots \\ y_{n,j,n} \end{bmatrix} = \begin{bmatrix} \tilde{w}_{1,1,j} & 0 & \dots & 0 \\ 0 & \tilde{w}_{2,2,j} & \ddots & \vdots \\ \vdots & \ddots & \ddots & 0 \\ 0 & \dots & 0 & \tilde{w}_{n,n,j} \end{bmatrix} \otimes I_d \begin{bmatrix} y_{n,0,1} \\ \vdots \\ y_{n,0,n} \end{bmatrix}$$

If $h \neq 0$, the stage system $\left(M, f, t, y_{n,0}, h\right)$ is equivalent to, for all $j \in [\![1, s+1]\!]$ :

$$\begin{bmatrix} y_{n,j,1} \\ \vdots \\ y_{n,j,n} \end{bmatrix} = \left(H_n(h)\widetilde{W}_{n,j}H_n(h)^{-1}\right) \otimes I_d \begin{bmatrix} y_{n,0,1} \\ \vdots \\ y_{n,0,n} \end{bmatrix} + h\left(H_n(h)D_n W_{n,j}\right) \otimes I_d \begin{bmatrix} f\left(t + \tau_1 h, y_{n,1}\right) \\ \vdots \\ f\left(t + \tau_s h, y_{n,s}\right) \end{bmatrix}$$

$$y_{n,s+1} \in \mathbb{R}^{[\![1,n]\!] \times [\![1,d]\!]}, \ \forall j \in [\![1, s]\!], \ \left(\left(t + \tau_j h\right), y_{n,j}\right) \in \mathbf{U}$$

Proof

Trivial.

To simplify the notations we define the two following functions.

**IV.2.5 Definition : Approximation vector function & Evaluation vector function**

Let $(\tilde{n}, s) \in \overline{\mathbb{N}}^* \times \mathbb{N}^*$, $M \in \text{GMORK}_{\tilde{n},s}$, $n \in \mathbb{N}^*$ such that $n \leq \tilde{n}$, $f$ a differential equation function of order $n$, $Y_{n,f}$ the approximation function for length $n$ of $M$, $F_{n,f}$ the evaluation function for length $n$ of $M$, $\left(t_0, y_{n,0}\right) \in \mathbf{U}, \tilde{q} \in \overline{\mathbb{N}}^*, \ h \in \mathcal{H}_{\tilde{q},f,t_0,y_{n,0}}$, and $t$ the time mesh of $\left(t_0, h\right)$.
We define $\tilde{\text{Y}}_{n,f}$ the approximation vector function for length $n$ of $M$ and $\tilde{\text{F}}_{n,f}$ the evaluation vector function for length $n$ of $M$ as, for all $q \in \mathbb{N}$ with $q < \tilde{q}$ :

$$\forall j \in [\![1, s+1]\!], \ \tilde{\text{Y}}_{n,f,q(s+1)+j}\left(t_0, y_{n,0}, h\right) = \begin{bmatrix} Y_{n,f,q(s+1)+j,1}\left(t_0, y_{n,0}, h\right) \\ \vdots \\ Y_{n,f,q(s+1)+j,n}\left(t_0, y_{n,0}, h\right) \end{bmatrix}$$

$$\tilde{\text{F}}_{n,f,q}\left(t_0, y_{n,0}, h\right) = \begin{bmatrix} F_{n,f,q(s+1)+1}\left(t_0, y_{n,0}, h\right) \\ \vdots \\ F_{n,f,q(s+1)+s}\left(t_0, y_{n,0}, h\right) \end{bmatrix}$$

### IV.2.6 Lemma : Stage system - Evaluation form

Let $(\tilde{n}, s) \in \overline{\mathbb{N}}^* \times \mathbb{N}^*$, $M = (\tau, w, \tilde{w}) \in \mathrm{GMORK}_{\tilde{n},s}$, $f$ a differential equation function of order $n \in \mathbb{N}^*$ such that $n \leq \tilde{n}$, and $\left(t, y_{n,0}, h\right) \in \mathbf{U} \times \mathbb{R}$.
The stage system $\left(M, f, t, y_{n,0}, h\right)$ is equivalent to :

$$\forall j \in [\![1, s]\!],\ k_j = f\left(t + \tau_j h, y_{n,j}\right),$$

$$\forall (j, N) \in [\![1, s+1]\!] \times [\![1, n]\!],\ y_{n,j,N} = \sum_{N'=1}^{N} \tilde{w}_{N,N',j} h^{N-N'} y_{n,0,N'} + \frac{h^N}{N!} \sum_{j'=1}^{s} w_{N,j,j'} k_{j'}$$

$$y_{n,s+1} \in \mathbb{R}^{[\![1,n]\!] \times [\![1,d]\!]},\ \forall j \in [\![1, s]\!],\ \left(\left(t + \tau_j h\right), y_{n,j}\right) \in \mathbf{U}$$

Proof

Trivial.

Let $(\tilde{n}, s) \in \overline{\mathbb{N}}^* \times \mathbb{N}^*$, $M = (\tau, w, \tilde{w}) \in \mathrm{GMORK}_{\tilde{n},s}$, $n \in \mathbb{N}^*$ such that $n \leq \tilde{n}$, $f$ a linear differential equation function, and $\left(t, y_{n,0}, h\right) \in \mathbb{R} \times \mathbb{C}^{[\![1,n]\!]} \times \mathbb{R}$. The time mesh goes to infinity if and only if $h \neq 0$, we thus assume $h \neq 0$. To help solve the stage system, we define :

$$E_s = \begin{bmatrix} I_s & 0_{s,1} \end{bmatrix},\ \forall N \in [\![1, n]\!],\ T_{n,N} = \begin{bmatrix} \tilde{w}_{N,1,1} & \dots & \tilde{w}_{N,N,1} & 0_{1,n-N} \\ \vdots & \ddots & \vdots & \vdots \\ \tilde{w}_{N,1,s} & \dots & \tilde{w}_{N,N,s} & 0_{1,n-N} \end{bmatrix}$$

We now solve the stage system $\left(M, f, t, y_{n,0}, h\right)$ in its evaluation form.

$$\begin{aligned}
\begin{bmatrix} k_1 \\ \vdots \\ k_s \end{bmatrix} &= \sum_{N=1}^{n} \alpha_N \begin{bmatrix} \sum_{N'=1}^{N} h^{N-N'} \tilde{w}_{N,N',1} y_{n,0,N'} + \frac{h^N}{N!} \sum_{j'=1}^{s} w_{N,1,j'} k_{j'} \\ \vdots \\ \sum_{N'=1}^{N} h^{N-N'} \tilde{w}_{N,N',s} y_{n,0,N'} + \frac{h^N}{N!} \sum_{j'=1}^{s} w_{N,s,j'} k_{j'} \end{bmatrix} \\
&= \sum_{N=1}^{n} \alpha_N \begin{bmatrix} h^{N-1} \tilde{w}_{N,1,1} & \dots & h^0 \tilde{w}_{N,N,1} & 0_{1,n-N} \\ \vdots & \ddots & \vdots & \vdots \\ h^{N-1} \tilde{w}_{N,1,s} & \dots & h^0 \tilde{w}_{N,N,s} & 0_{1,n-N} \end{bmatrix} \begin{bmatrix} y_{n,0,1} \\ \vdots \\ y_{n,0,n} \end{bmatrix} \\
&\quad + \sum_{N=1}^{n} \alpha_N \frac{h^N}{N!} \begin{bmatrix} w_{N,1,1} & \dots & w_{N,1,s} \\ \vdots & \ddots & \vdots \\ w_{N,s,1} & \dots & w_{N,s,s} \end{bmatrix} \begin{bmatrix} k_1 \\ \vdots \\ k_s \end{bmatrix} \\
&= \sum_{N=1}^{n} \alpha_N h^{N-1} \begin{bmatrix} \tilde{w}_{N,1,1} & \dots & \tilde{w}_{N,N,1} & 0_{1,n-N} \\ \vdots & \ddots & \vdots & \vdots \\ \tilde{w}_{N,1,s} & \dots & \tilde{w}_{N,N,s} & 0_{1,n-N} \end{bmatrix} H_n(h)^{-1} \begin{bmatrix} y_{n,0,1} \\ \vdots \\ y_{n,0,n} \end{bmatrix} + \sum_{N=1}^{n} \alpha_N \frac{h^N}{N!} E_s w_N \begin{bmatrix} k_1 \\ \vdots \\ k_s \end{bmatrix} \\
&= \sum_{N=1}^{n} \alpha_N h^{N-1} T_{n,N} H_n(h)^{-1} \tilde{\mathrm{Y}}_{n,f,0}\left(t, y_{n,0}, h\right) + \sum_{N=1}^{n} \alpha_N \frac{h^N}{N!} E_s w_N \begin{bmatrix} k_1 \\ \vdots \\ k_s \end{bmatrix}
\end{aligned}$$

The stage system $\left(M, f, t, y_{n,0}, h\right)$ is thus equivalent to :

$$\left(I_s - \sum_{N=1}^{n} \alpha_N \frac{h^N}{N!} E_s w_N\right) \begin{bmatrix} k_1 \\ \vdots \\ k_s \end{bmatrix} = \sum_{N=1}^{n} \alpha_N h^{N-1} T_{n,N} H_n(h)^{-1} \tilde{\mathrm{Y}}_{n,f,0}\left(t, y_{n,0}, h\right)$$

The stage system has a unique solution if and only if the matrix to the left is invertible. If we assume so, we get the expression of the evaluation vector function.

$$\tilde{\mathrm{F}}_{n,f}\big(t, y_{n,0}, h\big) = \left(I_s - \sum_{N=1}^{n} \alpha_N \frac{h^N}{N!} E_s w_N\right)^{-1} \sum_{N=1}^{n} \alpha_N h^{N-1} T_{n,N} H_n(h)^{-1} \tilde{\Upsilon}_{n,f,0}\big(t, y_{n,0}, h\big)$$

Let $j \in [\![1, s+1]\!]$. We replace the expression of $\tilde{\mathrm{F}}_{n,f}$ in the matrix form of $\tilde{\Upsilon}_{n,f,j}$.

$$\begin{aligned}
&\tilde{\Upsilon}_{n,f,j}\big(t, y_{n,0}, h\big) = H_n(h) \widetilde{W}_{n,j} H_n(h)^{-1} \tilde{\Upsilon}_{n,f,0}\big(t, y_{n,0}, h\big) \\
&\quad + h H_n(h) D_n W_{n,j} \left(I_s - \sum_{N=1}^{n} \alpha_N \frac{h^N}{N!} E_s w_N\right)^{-1} \sum_{N=1}^{n} \alpha_N h^{N-1} T_{n,N} H_n(h)^{-1} \tilde{\Upsilon}_{n,f,0}\big(t, y_{n,0}, h\big) \\
&= H_n(h) \left(\widetilde{W}_{n,j} + D_n W_{n,j} \left(I_s - \sum_{N=1}^{n} \alpha_N h^N E_s \frac{w_N}{N!}\right)^{-1} \sum_{N=1}^{n} \alpha_N h^N T_{n,N}\right) H_n(h)^{-1} \tilde{\Upsilon}_{n,f,0}\big(t, y_{n,0}, h\big)
\end{aligned}$$

To simplify this expression, we define a function which will be useful to characterize A-stability. We first define a version of this function which uses the coefficients of the linear differential equation since it can prove useful, but we are more interested in the version which uses the roots of $P_\alpha$.

**IV.2.7 Definition : Stability functions - Coefficient form**

Let $(\tilde{n}, s) \in \overline{\mathbb{N}}^* \times \mathbb{N}^*$, $M = (\tau, w, \tilde{w}) \in \mathrm{GMORK}_{\tilde{n},s}$, and $n \in \mathbb{N}^*$ such that $n \leq \tilde{n}$. We define :

$$D_n = \begin{bmatrix} \frac{1}{1!} & 0 & \dots & 0 \\ 0 & \frac{1}{2!} & \ddots & \vdots \\ \vdots & \ddots & \ddots & 0 \\ 0 & \dots & 0 & \frac{1}{n!} \end{bmatrix} \in \mathbb{R}^{[\![1,n]\!]^2}, \ E_s = \begin{bmatrix} I_s & 0_{s,1} \end{bmatrix} \in \mathbb{R}^{[\![1,s]\!] \times [\![1,s+1]\!]}$$

$$W_{n,j} = \begin{bmatrix} w_{1,j,1} & \dots & w_{1,j,s} \\ \vdots & \ddots & \vdots \\ w_{n,j,1} & \dots & w_{n,j,s} \end{bmatrix} \in \mathbb{R}^{[\![1,n]\!] \times [\![1,s]\!]}, \ \widetilde{W}_{n,j} = \begin{bmatrix} \tilde{w}_{1,1,j} & 0 & \dots & 0 \\ \tilde{w}_{2,1,j} & \tilde{w}_{2,2,j} & \ddots & \vdots \\ \vdots & \vdots & \ddots & 0 \\ \tilde{w}_{n,1,j} & \tilde{w}_{n,2,j} & \dots & \tilde{w}_{n,n,j} \end{bmatrix} \in \mathbb{R}^{[\![1,n]\!]^2}$$

$$\forall N \in [\![1, n]\!], \ T_{n,N} = \begin{bmatrix} \tilde{w}_{N,1,1} & \dots & \tilde{w}_{N,N,1} & 0_{1,n-N} \\ \vdots & \ddots & \vdots & \vdots \\ \tilde{w}_{N,1,s} & \dots & \tilde{w}_{N,N,s} & 0_{1,n-N} \end{bmatrix} \in \mathbb{R}^{[\![1,s]\!] \times [\![1,n]\!]}$$

We define $\tilde{\mathrm{Q}}_n \in \left\{\mathbb{C}^{[\![1,n]\!]} \to \mathbb{C}^{[\![1,s]\!]^2}\right\}$ the resolvent matrix in coefficient form for length $n$, and $\tilde{\mathcal{S}}_n \in \mathcal{P}\big(\mathbb{C}^{[\![1,n]\!]}\big)$ the domain of the stability functions in coefficient form for length $n$ of $M$ as :

$$\tilde{\mathrm{Q}}_n : \mathbb{C}^{[\![1,n]\!]} \to \mathbb{C}^{[\![1,s]\!]^2}, \ z \to I_s - \sum_{N=1}^{n} \frac{z_N}{N!} E_s w_N$$

$$\tilde{\mathcal{S}}_n = \left\{ z \in \mathbb{C}^{[\![1,n]\!]} \mid \det\big(\tilde{\mathrm{Q}}_n(z)\big) \neq 0 \right\}$$

We define $\tilde{R}_n \in \left\{\tilde{\mathcal{S}}_n \to \mathbb{C}^{[\![1,n]\!]^2}\right\}$ the stability function in coefficient form for length $n$ of $M$ as :

$$\forall z \in \tilde{\mathcal{S}}_n, \ \tilde{R}_n(z) = \widetilde{W}_{n,s+1} + D_n W_{n,s+1} \tilde{\mathrm{Q}}_n(z)^{-1} \sum_{N=1}^{n} z_N T_{n,N}$$

We thus have :

$$\tilde{Y}_{n,f,s+1}(t, y_{n,0}, h) = H_n(h)\tilde{R}_n(hH_n(h)\alpha)H_n(h)^{-1}\tilde{Y}_{n,f,0}(t, y_{n,0}, h)$$

We now use Vieta's formula.

$$\tilde{Y}_{n,f,j}(t, y_{n,0}, h) = H_n(h)\Bigg(\widetilde{W}_{n,j}$$

$$+D_n W_{n,j}\left(I_s - \sum_{N=1}^{n}\left((-1)^{N+1}\sum_{1\le i_1<\ldots<i_N\le n}\prod_{k=1}^{N} r_{i_k}\right)\frac{h^N}{N!}E_s w_N\right)^{-1}\sum_{N=1}^{n}\left((-1)^{N+1}\sum_{1\le i_1<\ldots<i_N\le n}\prod_{k=1}^{N} r_{i_k}\right)h^N T_{n,N}$$

$$\Bigg)H_n(h)^{-1}\tilde{Y}_{n,f,0}(t, y_{n,0}, h)$$

$$= H_n(h)\Bigg(\widetilde{W}_{n,j}$$

$$-D_n W_{n,j}\left(I_s + \sum_{N=1}^{n}\sum_{1\le i_1<\ldots<i_N\le n}\left(\prod_{k=1}^{N}(-r_{i_k}h)\right)E_s\frac{w_N}{N!}\right)^{-1}\sum_{N=1}^{n}\sum_{1\le i_1<\ldots<i_N\le n}\left(\prod_{k=1}^{N}(-r_{i_k}h)\right)T_{n,N}$$

$$\Bigg)H_n(h)^{-1}\tilde{Y}_{n,f,0}(t, y_{n,0}, h)$$

This time, we define the root version of the stability functions, the version which actually interests us.

**IV.2.8 Definition : Stability functions**

Let $(\tilde{n}, s) \in \overline{\mathbb{N}}^* \times \mathbb{N}^*$, $M = (\tau, w, \tilde{w}) \in \mathrm{GMORK}_{\tilde{n},s}$, and $n \in \mathbb{N}^*$ such that $n \le \tilde{n}$. We consider the same matrices as in Section IV.2.7.

We define $Q_n \in \left\{\mathbb{C}^{[\![1,n]\!]} \to \mathbb{C}^{[\![1,s]\!]^2}\right\}$ the resolvent matrix for length $n$, and $\mathcal{S}_n \in \mathcal{P}\left(\mathbb{C}^{[\![1,n]\!]}\right)$ the domain of definition of the stability functions for length $n$ of $M$ as :

$$Q_n : \mathbb{C}^{[\![1,n]\!]} \to \mathbb{C}^{[\![1,s]\!]^2},\ z \to I_s + \sum_{N=1}^{n}\sum_{1\le i_1<\ldots<i_N\le n}\left(\prod_{k=1}^{N} z_{i_k}\right)E_s\frac{w_N}{N!}$$

$$\mathcal{S}_n = \left\{z \in \mathbb{C}^{[\![1,n]\!]} \mid \det(Q_n(z)) \ne 0\right\}$$

We define $R_n \in \left\{\mathcal{S}_n \to \mathbb{C}^{[\![1,n]\!]^2}\right\}$ the stability function for length $n$ of $M$ as, for all $z \in \mathcal{S}_n$ :

$$R_n(z) = \widetilde{W}_{n,s+1} - D_n W_{n,s+1} Q_n(z)^{-1}\sum_{N=1}^{n}\sum_{1\le i_1<\ldots<i_N\le n}\left(\prod_{k=1}^{N} z_{i_k}\right)T_{n,N}$$

When we refer to the stability functions we refer to the root form we just defined, not the coefficient form. We now sum up what we have seen in the previous pages.

**IV.2.9 Proposition : Approximation of linear initial value problems**

Let $(\tilde{n}, s) \in \overline{\mathbb{N}}^* \times \mathbb{N}^*$, $M = (\tau, w, \tilde{w}) \in \mathrm{GMORK}_{\tilde{n},s}$, $R$ the stability functions of $M$, $\mathcal{S}$ their domains of definition, $n \in \mathbb{N}^*$ such that $n \leq \tilde{n}$, and $(f, t_0, y_n)$ a linear initial value problem such that :

$$\alpha \in \mathbb{C}^{[\![1,n]\!]},\ \left(t_0, y_{n,0}\right) \in \mathbb{R} \times \mathbb{C}^{[\![1,n]\!]}$$

$$\forall (t,x) \in \mathbb{R} \times \mathbb{C}^{[\![1,n]\!] \times [\![1,d]\!]},\ f(t,x) = \sum_{N=1}^{n} \alpha_N x_N$$

Let $r \in \mathbb{C}^{[\![1,n]\!]}$ the roots of $-X^n + \sum_{i=0}^{n-1} \alpha_{n-i} X^i$, the characteristic polynomial of $f$, and $\tilde{\mathrm{Y}}$ the approximation vector function for length $n$ of $M$.

- We have $\mathcal{U}_{n,f} = \left\{ \left(t, h, \hat{y}_{n,0}\right) \mid -hr \in \mathcal{S}_n \right\}$.

Let $h \in \mathbb{R}$ such that $-hr \in \mathcal{S}_n$, and $\left(t, y_{n,0}\right) \in \mathbb{R} \times \mathbb{C}^{[\![1,n]\!]}$. We define :

$$H_n : \mathbb{R} \to \mathbb{R}^{[\![1,n]\!]^2},\ h \to \begin{bmatrix} h^0 & 0 & \dots & 0 \\ 0 & h^1 & \ddots & \vdots \\ \vdots & \ddots & \ddots & 0 \\ 0 & \dots & 0 & h^{n-1} \end{bmatrix}$$

- If $h \neq 0$, then $\tilde{\mathrm{Y}}_{n,f,s+1}\left(t, y_{n,0}, h\right) = H_n(h) R_n(-hr) H_n(h)^{-1} \tilde{\mathrm{Y}}_{n,f,0}\left(t, y_{n,0}, h\right)$

Proof

This proposition was proved in the text between the definition of A-stability and the present statement.

This proposition allows to characterize A-stability in terms of the stability functions.

**IV.2.10 Theorem : Characterization of A-stability & Absolute A-stability**

Let $(\tilde{n}, s) \in \overline{\mathbb{N}}^* \times \mathbb{N}^*$, $M \in \mathrm{GMORK}_{\tilde{n},s}$, $R$ the stability functions of $M$, $\mathcal{S}$ their domains of definition, $n \in \mathbb{N}^*$ such that $n \leq \tilde{n}$, and $\rho$ the spectral radius function.

1. M is absolute A-stable for length $n$ if and only if, for all $z \in \mathcal{S}_n$, if, for all $k \in [\![1, n]\!]$, $\Re(z_k) > 0$, then $\rho(R_n(z)) < 1$.
2. M is A-stable up to length $n$ if and only if, for all $z \in \mathcal{S}_n$, if, for all $k \in [\![1, n]\!]$, $\Re(z_k) \geq 0$, then $\rho(R_n(z)) \leq 1$.

Proof

Let $\tilde{\mathrm{Y}}_{n,f}$ the approximation vector function for length $n$ of $M$. Section IV.2.9 implies that a sequence $(h)_{q\in\mathbb{N}}$ with $h \in \mathbb{R}$ is a step size sequence if and only if $-hr \in \mathcal{S}_n$. If true, we have, for all $q \in \mathbb{N}$ :

$$\tilde{\mathrm{Y}}_{n,f,q(s+1)}\left(t, y_{n,0}, (h)_{q\in\mathbb{N}}\right) = H_n(h) R_n(-hr)^q H_n(h)^{-1} \tilde{\mathrm{Y}}_{n,f,0}\left(t, y_{n,0}, h\right)$$

1. Since $\tilde{\mathrm{Y}}_{n,f,q(s+1)}$ is linear in terms of the vector of initial values, it is thus equivalent to prove convergence for every vector of a basis of $\mathbb{C}^{[\![1,n]\!]}$ since the other initial values follow by linear combinations. If we choose the canonical basis, each vector of the basis gives that a different column of $H_n(h) R_n(-hr)^q H_n(h)^{-1}$ converges towards $0_n$. Therefore the approximations converges towards $0_n$ for all initial value problems if and only if $H_n(h) R_n(-hr)^q H_n(h)^{-1}$ converges towards $0_{n,n}$. Since the matrices $H_n(h)$ and $H_n(h)^{-1}$ are constant invertible matrices we multiply by their inverses in the limit to get that the approximations converges towards $0_n$ if and only if $\lim_{q\to+\infty} R_n(-hr)^q = 0_{n,n}$. A matrix satisfying this condition is called a convergent matrix, and a square matrix $A$ is convergent if and only if $\rho(A) < 1$, hence in this case, $\rho(R_n(-hr)) < 1$. By defining $z = -hr$, it is equivalent to go over all $z \in \mathcal{S}_n$ satisfying, for all $k \in [\![1, n]\!]$, $\Re(z_k) > 0$.

2. This proof is similar to 1, except each vector of the canonical basis gives that a different column of $R_n(-hr)^q$ stays bounded, and $R_n(-hr)^q$ stays bounded if and only if $\rho(R_n(-hr)) \leq 1$.

Since A-stability and absolute A-stability entirely depend on the stability functions, we now establish some properties of the stability functions, properties we later use to study the A-stability and absolute A-stability of a method.

**IV.2.11 Proposition : Properties of stability functions**

Let $(\tilde{n}, s) \in \overline{\mathbb{N}}^* \times \mathbb{N}^*$, $M \in \mathrm{GMORK}_{\tilde{n},s}$, $Q$ the resolvent matrices of $M$, $R$, the stability functions of $M$, $\mathcal{S}$ their domains of definition, and $n \in \mathbb{N}^*$ such that $n \leq \tilde{n}$.

1. $Q_n$ is a multivariate polynomial matrix of total degree of at most $n$.
2. $\mathcal{S}_n$ is an open set dense in $\mathbb{C}^{[\![1,n]\!]}$.
3. $R_n$ is a multivariate polynomial matrix of total degree inferior or equal to $ns$ divided by a multivariate polynomial of total degree inferior or equal to $ns$. We consider the matrices of Section IV.2.7. For all $z \in \mathcal{S}_n$, with Com the comatrix function :

$$R_n(z) = \frac{\widetilde{W}_{n,s+1}\det(Q_n(z)) - D_n W_{n,s+1}\mathrm{Com}(Q_n(z))^T \sum_{N=1}^{n} \sum_{1 \leq i_1 < \ldots < i_N \leq n} \left(\prod_{k=1}^{N} z_{i_k}\right) T_{n,N}}{\det(Q_n(z))}$$

Proof

1. Each entry of $Q_n$ is a multivariate polynomial of total degree of at most $n$, therefore $Q_n$ is a multivariate matrix of total degree of at most $n$.
2. We deduce from 1 that the determinant of the resolvent matrix is a multivariate polynomial of total degree of at most $ns$. $\mathcal{S}_n$ is the preimage of $\mathbb{R}^*$ under a continuous function, it is hence an open set. If we assume $\mathcal{S}_n$ is not dense, then there exists an open subset of $\mathbb{C}^{[\![1,n]\!]}$ such that $\det(Q_n)$ restricted to this subset is identically 0. Since $\det(Q_n)$ is a multivariate polynomial it would then be identically 0, but $\det(Q_n(0_n)) = 1$, which is a contradiction, thus $\mathcal{S}_n$ is dense.
3. The expression of $R_n$ comes from the fact that given an invertible matrix $A$, we have $A^{-1} = \frac{\mathrm{Com}(A)^T}{\det(A)}$. Applying this formula to $Q_n^{-1}$ in $R_n$ directly gives this expression.
   From the bound on the total degree of each entry of $Q_n$, we get that $\mathrm{Com}(Q_n)^T$ is a polynomial matrix of degree of at most $(s-1)n$. Since it is multiplied by a term with a total degree of at most $n$, the overall degree of the nominator is at most $sn$. The total degree of the denominator is at most $ns$ as well.

An important case to consider is when 0 is a root of the characteristic polynomial. Let $f$ a linear differential equation function of order $n$, and $\hat{y}$ a solution candidate. Assume 0 is a root of multiplicity $m \in \mathbb{N}^*$, which is equivalent to $\alpha_{n-m} \neq 0$ and, for all $k \in [\![n+1-m, n]\!]$, $\alpha_k = 0$. If we define $\hat{y}' = \frac{\mathrm{d}^m \hat{y}}{\mathrm{d}t^m}$, which is similar to the rewriting for initial value problem for lower order initial value problems, then :

$$\frac{\mathrm{d}^n \hat{y}}{\mathrm{d}t^n} = \sum_{N=1}^{n} \alpha_N \frac{\mathrm{d}^{n-N}\hat{y}}{\mathrm{d}t^{n-N}} \Leftrightarrow \frac{\mathrm{d}^{n-m}\hat{y}'}{\mathrm{d}t^{n-m}} = \sum_{N=1}^{n-m} \alpha_N \frac{\mathrm{d}^{n-m-N}\hat{y}'}{\mathrm{d}t^{n-m-N}}$$

The second equation is a linear differential equation of order $n-m$. This implies that the multiplicity of the root 0 dictates the actual order of the initial value problem. We now investigate the properties of the stability functions when some of the roots are 0.

**IV.2.12 Proposition : Multiplicity of the zero root & Stability functions**

Let $(\tilde{n}, s) \in \overline{\mathbb{N}}^* \times \mathbb{N}^*$, $M = (\tau, w, \tilde{w}) \in \mathrm{GMORK}_{\tilde{n},s}$, $Q$ the resolvent matrices of $M$, $R$ the stability functions of $M$, $\mathcal{S}$ the domains of the stability functions of M, $(n, n') \in \mathbb{N}^* \times \mathbb{N}^*$ such that $n \leq n' \leq \tilde{n}$. We have :

$$\forall z \in \mathbb{C}^{[\![1,n]\!]},\ Q_n(z) = Q_{n'}\left(\begin{bmatrix} z \\ 0_{n'-n} \end{bmatrix}\right)$$

$$\mathcal{S}_n = \left\{ z \in \mathbb{C}^{[\![1,n]\!]} \mid \begin{bmatrix} z \\ 0_{n'-n} \end{bmatrix} \in \mathcal{S}_{n'} \right\}$$

$$\forall z \in \mathcal{S}_n,\ R_{n'}\left(\begin{bmatrix} z \\ 0_{n'-n} \end{bmatrix}\right) = \begin{bmatrix} R_n(z) & 0_n & 0_n & \dots & 0_n \\ * & \tilde{w}_{n+1,n+1,s+1} & 0 & \dots & 0 \\ * & \tilde{w}_{n+2,n+,s+1} & \tilde{w}_{n+2,n+2,s+1} & \ddots & \vdots \\ \vdots & \vdots & \vdots & \ddots & 0 \\ * & \tilde{w}_{n',n+1,s+1} & \tilde{w}_{n',n+2,s+1} & \dots & \tilde{w}_{n',n',s+1} \end{bmatrix}$$

$$= \begin{bmatrix} R_n(z) & 0_{n,n'-n} \\ * & \begin{bmatrix} 0_{n,n'-n} & I_{n'-n} \end{bmatrix} \widetilde{W}_{n',s+1} \begin{bmatrix} 0_{n,n'-n} \\ I_{n'-n} \end{bmatrix} \end{bmatrix}$$

$$\forall z \in \mathcal{S}_n,\ \rho\left(R_{n'}\left(\begin{bmatrix} z \\ 0_{n'-n} \end{bmatrix}\right)\right) = \max\left(\rho(R_n(z)), |\tilde{w}_{n+1,n+1,s+1}|, \dots, |\tilde{w}_{n',n',s+1}|\right) \geq \rho(R_n(z))$$

Proof

We first prove the equality which involves the resolvent matrix. Let $z \in \mathbb{C}^{[\![1,n]\!]}$. We have :

$$Q_{n'}\left(\begin{bmatrix} z \\ 0_{n'-n} \end{bmatrix}\right) = I_s + \sum_{N=1}^{n'} \sum_{1 \leq i_1 < \ldots < i_N \leq n'} \left( \prod_{k=1}^{N} z_{i_k} \right) E_s \frac{w_N}{N!}$$

If the index satisfies $i_N > n$, then $z_{i_N} = 0$, and the product vanishes. We thus restrict the sum to indices satisfying $i_N \leq n$. Morevoer, if $N > n$, then the index $i$ would have $N$ distinct entries contained in $[\![1, n]\!]$, which is impossible, hence we have $N \leq n$. The second equality is a consequence of the first equality. The third equality follows from an argument similar to the one used for the first, we restrict the sum over $T_{n,N}$ to indices satisfying $N \leq n$ and $i_N \leq n$. The shape of the $T_{n,N}$ implies that the second term in the expression of $R_{n'}$ has its last $m$ columns equal to zero. The last equality follows.

The previous proposition assumes that the zero roots are labeled last, this is without loss of generality.

**IV.2.13 Proposition : Symmetry of the stability functions**

Let $(\tilde{n}, s) \in \overline{\mathbb{N}}^* \times \mathbb{N}^*$, $M \in \mathrm{GMORK}_{\tilde{n},s}$, $Q$ the resolvent matrices of $M$, $R$ the stability functions of $M$, $\mathcal{S}$, their domains of definition, and $n \in \mathbb{N}^*$ such that $n \leq \tilde{n}$.

$Q, \mathcal{S}, R$ are symmetric under coordinate permutation, hence, for all $\varphi \in \mathfrak{S}_{[\![1,n]\!]}$ :

$$\forall z \in \mathbb{C}^{[\![1,n]\!]},\ Q_n\left(\begin{bmatrix} z_{\varphi(1)} \\ \vdots \\ z_{\varphi(n)} \end{bmatrix}\right) = Q_n(z)$$

$$\forall z \in \mathcal{S}_n,\ \begin{bmatrix} z_{\varphi(1)} \\ \vdots \\ z_{\varphi(n)} \end{bmatrix} \in \mathcal{S}_n,\ R_n\left(\begin{bmatrix} z_{\varphi(1)} \\ \vdots \\ z_{\varphi(n)} \end{bmatrix}\right) = R_n(z)$$

Proof

It is a direct consequence of the fact that the way the the roots are labeled does not change the coefficients of the linear differential equation function.

Since all resolvant matrices, stability functions and domains of definition can be expressed in terms of $Q_{\tilde{n}}$, $\mathcal{S}_{\tilde{n}}$ and $R_{\tilde{n}}$, we may name these objects.

**IV.2.14 Definition : Maximum stability function - Finite methods**

Let $(\tilde{n}, s) \in \mathbb{N}^* \times \mathbb{N}^*$, $M \in \mathrm{GMORK}_{\tilde{n},s}$, $Q$ the resolvent matrices of $M$, $R$ the stability functions of $M$, and $\mathcal{S}$ their domains of definition.
We define the maximum resolvent matrix as $Q_{\tilde{n}}$, the maximum domain of the stability function as $\mathcal{S}_{\tilde{n}}$, and the maximum stability function as $R_{\tilde{n}}$.

Similar objects can also be defined for infinite methods, although this requires the use of infinite matrices. No convergence issues arise, since the series involved contain only finitely many nonzero coefficients. The corresponding infinite arguments therefore belong to $c_{00} \in \mathcal{P}(\mathbb{C}^{\mathbb{N}^*})$, where $c_{00}$ denotes the space of sequences with finite support.

**IV.2.15 Definition : Maximum stability function - Infinite methods**

Let $s \in \mathbb{N}^*$, $M = (\tau, w, \tilde{w}) \in \mathrm{GMORK}_{\infty,s}$, and $c_{00} \in \mathcal{P}(\mathbb{C}^{\mathbb{N}^*})$ the set of sequences with finite support.
We define :

$$D_\infty = \begin{bmatrix} \frac{1}{1!} & 0 & \cdots \\ 0 & \frac{1}{2!} & \ddots \\ \vdots & \ddots & \ddots \end{bmatrix} \in \mathbb{R}^{\mathbb{N}^* \times \mathbb{N}^*},\ E_s = \begin{bmatrix} I_s & 0_{s,1} \end{bmatrix} \in \mathbb{R}^{[\![1,s]\!] \times [\![1,s+1]\!]}$$

$$W_{\infty,s+1} = \begin{bmatrix} w_{1,s+1,1} & \cdots & w_{1,s+1,s} \\ w_{2,s+1,1} & \cdots & w_{2,s+1,s} \\ \vdots & \ddots & \vdots \end{bmatrix} \in \mathbb{R}^{\mathbb{N}^* \times [\![1,s]\!]},\ \widetilde{W}_{\infty,s+1} = \begin{bmatrix} \tilde{w}_{1,1,s+1} & 0 & \cdots \\ \tilde{w}_{2,1,s+1} & \tilde{w}_{2,2,s+1} & \ddots \\ \vdots & \vdots & \ddots \end{bmatrix} \in \mathbb{R}^{\mathbb{N}^* \times \mathbb{N}^*}$$

$$\forall N \in \mathbb{N}^*,\ T_{\infty,N} = \begin{bmatrix} \tilde{w}_{N,1,1} & \cdots & \tilde{w}_{N,N,1} & 0 & \cdots \\ \vdots & \ddots & \vdots & \vdots & \cdots \\ \tilde{w}_{N,1,s} & \cdots & \tilde{w}_{N,N,s} & 0 & \cdots \end{bmatrix} \in \mathbb{R}^{[\![1,s]\!] \times \mathbb{N}^*}$$

We define $Q_\infty \in \left\{c_{00} \to \mathbb{C}^{[\![1,s]\!]^2}\right\}$ the maximum resolvent matrix of $M$, $R_\infty \in \left\{\mathcal{S}_\infty \to \mathbb{C}^{\mathbb{N}^* \times \mathbb{N}^*}\right\}$ the maximum stability function of $M$, $\mathcal{S}_\infty$ its domain of definition as :

$$Q_\infty : c_{00} \to \mathbb{C}^{[\![1,s]\!]^2}, z \to I_s + \sum_{N=1}^{+\infty} \frac{1}{N!} \sum_{1 \leq i_1 < \ldots < i_N} \left( \prod_{k=1}^{N} z_{i_k} \right) E_s w_N$$

$$\mathcal{S}_\infty = \{ z \in c_{00} \mid \det(Q_\infty(z)) \neq 0 \}$$

$$R_\infty : \mathcal{S}_\infty \to \mathbb{C}^{\mathbb{N}^* \times \mathbb{N}^*}, z \to \widetilde{W}_{\infty,s+1} - D_\infty W_{\infty,s+1} Q_\infty(z)^{-1} \sum_{N=1}^{+\infty} \sum_{1 \leq i_1 < \ldots < i_N} \left( \prod_{k=1}^{N} z_{i_k} \right) T_{\infty,N}$$

We can generalize the previous theorems to the maximum stability function of infinite methods.

**IV.2.16 Proposition : Multiplicity of the zero root & Stability functions - Infinite methods**

Let $s \in \mathbb{N}^*$, $M = (\tau, w, \tilde{w}) \in \mathrm{GMORK}_{\infty,s}$, $R$ the stability functions of $M$, $\mathcal{S}$ their domains of definitions, and $n \in \mathbb{N}^*$.
We have :

$$\forall z \in \mathbb{C}^{[\![1,n]\!]},\ Q_n(z) = Q_\infty\left(\begin{bmatrix} z \\ 0_\infty \end{bmatrix}\right)$$

$$\mathcal{S}_n = \left\{ z \in \mathbb{C}^{[\![1,n]\!]} \mid \begin{bmatrix} z \\ 0_\infty \end{bmatrix} \in \mathcal{S}_\infty \right\}$$

$$\forall z \in \mathcal{S}_n,\ R_\infty\left(\begin{bmatrix} z \\ 0_\infty \end{bmatrix}\right) = \begin{bmatrix} R_n(z) & 0_n & 0_n & \cdots \\ * & \tilde{w}_{n+1,n+1,s+1} & 0 & \cdots \\ * & \tilde{w}_{n+2,n+1,s+1} & \tilde{w}_{n+2,n+2,s+1} & \ddots \\ \vdots & \vdots & \vdots & \ddots \end{bmatrix}$$

$$= \begin{bmatrix} R_n(z) & 0_{n,\infty} \\ * & \begin{bmatrix} 0_{n,\infty} & I_\infty \end{bmatrix} \widetilde{W}_{\infty,s+1} \begin{bmatrix} 0_{n,\infty} \\ I_\infty \end{bmatrix} \end{bmatrix}$$

$$\forall z \in \mathcal{S}_n,\ \rho\left(R_{n'}\left(\begin{bmatrix} z \\ 0_\infty \end{bmatrix}\right)\right) = \max\left(\rho(R_n(z)), \sup_{N \in [\![n+1,+\infty[\![} \left|\tilde{w}_{N,N,s+1}\right|\right) \geq \rho(R_n(z))$$

Proof

Similar to the proof of Section IV.2.12.

We now have all the tools needed to relate the A-stability properties of a method across different lengths.

**IV.2.17 Proposition : Different lengths of A-stability**

Let $(\tilde{n}, s) \in \overline{\mathbb{N}}^* \times \mathbb{N}^*$, $M \in \mathrm{GMORK}_{\tilde{n},s}$, and $n \in \mathbb{N}^*$ such that $n \leq \tilde{n}$.

1. If M is A-stable up to length $n$, then for all $n' \in [\![1, n]\!]$, $M$ is A-stable up to length $n'$.
2. If $M$ is absolute A-stable for length $n$ then it is A-stable up to length $n$.

Proof

1. Since, for all $z \in \mathcal{S}_{n'}$ such that, for all $k \in [\![1, n]\!]$, $\mathfrak{R}(z_k) \geq 0$ we have $\rho(R_{n'}(z)) \leq \rho\left(R_n\left(\begin{bmatrix} z \\ 0_{n-n'} \end{bmatrix}\right)\right) \leq 1$, we deduce that $M$ is A-stable up to length $n'$.
2. Let $z \in \mathcal{S}_n$ such that, for all $k \in [\![1, n]\!]$, $\mathfrak{R}(z) \geq 0$. Since $\mathcal{S}_n$ is dense there exists a sequence $z' \in \mathcal{S}_n^{\mathbb{N}}$ satisfying, for all $k \in [\![1, n]\!]$, $\mathfrak{R}(z'_k) > 0$ and $\lim_{k \to \infty} z'_k = z$. We deduce from the continuity of $\rho$ that $\rho(R_n(z)) \leq 1$

The first proposition justifies the use of up to for A-stability. It implies that it gets more and more difficult for a method to be A-stable. However, at higher ranks, we have total freedom over the weights of the stages $j \in [\![1, s]\!]$, since, as stated in Section III.1.25, they do not impact the order of consistency at last stage.

We now examine how the structure of a method affects its stability function.

**IV.2.18 Proposition : A-stability & Permutations**

Let $(\tilde{n}, s) \in \overline{\mathbb{N}}^* \times \mathbb{N}^*$, $M \in \mathrm{GMORK}_{\tilde{n},s}$, $Q$ the resolvent matrices of $M$, $R$ the stability functions of $M$, $\mathcal{S}$ their domains of definition, $\varphi \in \mathfrak{S}_s^*$, $P_\varphi$ the permutation matrix of $\varphi$ restricted to $[\![1, s]\!]$, $n \in [\![1, \tilde{n}]\!]$, $Q'$ the resolvent matrices of $\varphi * M$, $R'$ the stability functions of $\varphi * M$, and $\mathcal{S}'$ their domains of definition.
We have $Q'_n = P_\varphi Q_n P_\varphi^T,\ \mathcal{S}'_n = \mathcal{S}_n,\ R'_n = R_n$

Proof

We have $P_\varphi^{-1} = P_\varphi^T$, hence, thanks to Section II.2.3 :

$$\begin{aligned} Q'_n(z) &= I_s + \sum_{N=1}^{n} \sum_{1\le i_1<\ldots<i_N\le n} \left(\prod_{k=1}^{N} z_{i_k}\right) E_s \frac{w'_N}{N!} \\ &= P_\varphi I_s P_\varphi^T + \sum_{N=1}^{n} \sum_{1\le i_1<\ldots<i_N\le n} \left(\prod_{k=1}^{N} z_{i_k}\right) P_\varphi E_s \frac{w_N}{N!} P_\varphi^T \\ &= P_\varphi Q_n(z) P_\varphi^T \end{aligned}$$

The equality of the stability domains and stability functions is a direct consequence of the equality of the approximations at last stage proved in Section II.2.3.

We first generalize the second Dahlquist barrier. We prove that explicit methods cannot be A-stable up to length 1, which implies that they are not A-stable up to any length, which implies by contradiction that they are not absolute A-stable either. We now find the expression of the stability functions of explicit methods.

**IV.2.19 Proposition : Stability functions of explicit methods**

Let $(\tilde{n}, s) \in \overline{\mathbb{N}}^* \times \mathbb{N}^*$, $M \in \mathrm{EGMORK}_{\tilde{n},s}$, $Q$ the resolvent matrices of $M$, $\mathcal{S}$ their domains of definition, and $n \in \mathbb{N}^*$ such that $n \le \tilde{n}$.

We have :

$$\begin{aligned} \forall z \in \mathbb{C}^{[\![1,n]\!]},\ \det(Q_n(z)) &= 1 \\ \mathcal{S}_n &= \mathbb{C}^{[\![1,n]\!]} \\ \forall z \in \mathbb{C}^{[\![1,n]\!]},\ Q_n(z)^{-1} &= \sum_{k=0}^{s-1} \left( - \sum_{N=1}^{n} \sum_{1\le i_1<\ldots<i_N\le n} \left(\prod_{k=1}^{N} z_{i_k}\right) E_s \frac{w_N}{N!} \right)^k \end{aligned}$$

If $\tilde{n} = +\infty$, then :

$$\begin{aligned} \forall z \in c_{00},\ \det\left(Q_{\infty(z)}\right) &= 1 \\ \mathcal{S}_\infty &= c_{00} \\ \forall z \in c_{00},\ Q_\infty(z)^{-1} &= \sum_{k=0}^{s-1} \left( - \sum_{N=1}^{\infty} \sum_{1\le i_1<\ldots<i_N\le n} \left(\prod_{k=1}^{N} z_{i_k}\right) E_s \frac{w_N}{N!} \right)^k \end{aligned}$$

Proof

Let $A \in \mathbb{C}^{[\![1,s]\!]^2}$ a strictly lower triangular matrix. We have $\det(I_s + A) = 1$. $A$ is a nilpotent matrix of index of at most $s$ and the inverse of $I_s + A$ is a Neumann series, $(I_s + A)^{-1} = \sum_{k=0}^{s-1} (-A)^k$. $\sum_{N=1}^{n} \sum_{1\le i_1<\ldots<i_N\le n} \left(\prod_{k=1}^{N} z_{i_k}\right) E_s \frac{w_N}{N!}$ is similar to a strictly lower triangular matrix, we thus deduce the expression of its inverse and the remainder of the proposition.

**IV.2.20 Theorem : Explicit methods & A-stability**

Let $(\tilde{n}, s) \in \overline{\mathbb{N}}^* \times \mathbb{N}^*$, and $M = (\tau, w, \tilde{w}) \in \mathrm{GMORK}_{\tilde{n},s}$.

If $M$ is explicit and $\sum_{j=1}^{s} w_{1,s+1,j} \tilde{w}_{1,1,j} \neq 0$ then $M$ is neither A-stable up to any length nor absolute A-stable for any length.

Proof

$R_1$ is a scalar, its spectral radius is therefore its absolute value. Since $R_1$ is a polynomial, it is bounded by 1 if and only if it is a constant function with an absolute value inferior to 1. Section IV.2.19 gives :

$$
\begin{aligned}
R_1(z) &= \tilde{w}_{1,1,s+1} - w_{1,s+1}^T \sum_{k=0}^{s-1} (-zE_s w_1)^k z \begin{bmatrix} \tilde{w}_{1,1,1} \\ \vdots \\ \tilde{w}_{1,1,s} \end{bmatrix} \\
&= \tilde{w}_{1,1,s+1} - z w_{1,s+1}^T \begin{bmatrix} \tilde{w}_{1,1,1} \\ \vdots \\ \tilde{w}_{1,1,s} \end{bmatrix} + w_{1,s+1}^T \sum_{k=1}^{s-1} (-z)^{k+1} (E_s w_1)^k \begin{bmatrix} \tilde{w}_{1,1,1} \\ \vdots \\ \tilde{w}_{1,1,s} \end{bmatrix} \\
&= \tilde{w}_{1,1,s+1} - z \sum_{j=1}^{s} w_{1,s+1,j} \tilde{w}_{1,1,j} + w_{1,s+1}^T \sum_{k=2}^{s} (-z)^k (E_s w_1)^{k-1} \begin{bmatrix} \tilde{w}_{1,1,1} \\ \vdots \\ \tilde{w}_{1,1,s} \end{bmatrix}
\end{aligned}
$$

If $\sum_{j=1}^{s} w_{1,s+1,j} \tilde{w}_{1,1,j} \neq 0$, the polynomial is not constant, hence $M$ is not A-stable up to length 1. It follows that $M$ is not A-stable up to any length. $M$ cannot be absolute A-stable for any length, since absolute A-stability implies A-stability.

The condition $\sum_{j=1}^{s} w_{1,s+1,j} \tilde{w}_{1,1,j} \neq 0$ cannot simply be removed. For instance, the null method is both A-stable and absolutely A-stable. The null method is the explicit method whose weights and nodes are all equal to zero, and it therefore always returns the value 0. This condition is essential for explicit methods since it is a necessary condition for their convergence.

**IV.2.21 Theorem : A necessary condition for the convergence of explicit methods**

Let $(\tilde{n}, s) \in \overline{\mathbb{N}}^* \times \mathbb{N}^*$, and $M = (\tau, w, \tilde{w}) \in \mathrm{GMORK}_{\tilde{n},s}$.
If $M$ is explicit and convergent up to length 1, then $\sum_{j=1}^{s} w_{1,s+1,j} \tilde{w}_{1,1,j} = 1$.

Proof

We consider the initial value problem of order 1 $(f, y_1, t_0)$ where $f(t,x) = x$, $y_{1,0,1} = 1$, and $t_0 = 0$. The solution is $e^t$. We consider the step size sequences $\left(\frac{1}{\tilde{q}}\right)_{q \in [\![0, \tilde{q}-1]\!]}$ with $\tilde{q} \in \mathbb{N}^*$. We have :

$$
\begin{aligned}
& Y_{1,f,\tilde{q}(s+1),1}\left(t_0, y_{1,0}, \left(\frac{1}{\tilde{q}}\right)_{q \in [\![0,\tilde{q}-1]\!]}\right) = R_1\left(-\frac{1}{\tilde{q}}\right)^{\tilde{q}} y_{1,0,1} \\
&= \left(1 + \frac{1}{\tilde{q}} \sum_{j=1}^{s} w_{1,s+1,j} \tilde{w}_{1,1,j} + \mathcal{O}\left(\frac{1}{\tilde{q}^2}\right)\right)^{\tilde{q}} \\
&= \exp\left(\tilde{q} \ln\left(1 + \frac{1}{\tilde{q}} \sum_{j=1}^{s} w_{1,s+1,j} \tilde{w}_{1,1,j} + \mathcal{O}\left(\frac{1}{\tilde{q}^2}\right)\right)\right) \\
&= \exp\left(\sum_{j=1}^{s} w_{1,s+1,j} \tilde{w}_{1,1,j} + \mathcal{O}\left(\frac{1}{\tilde{q}}\right)\right) \\
&\underset{\tilde{q} \to +\infty}{\longrightarrow} \exp\left(\sum_{j=1}^{s} w_{1,s+1,j} \tilde{w}_{1,1,j}\right)
\end{aligned}
$$

If $\sum_{j=1}^{s} w_{1,s+1,j} \tilde{w}_{1,1,j} \neq 1$ then the approximations does not converge towards the solution $e$.

Together, the two previous theorems allow us to generalize the second Dahlquist barrier.

**IV.2.22 Corollary : Generalized second Dahlquist barrier**

Let $(\tilde{n}, s) \in \overline{\mathbb{N}}^* \times \mathbb{N}^*$, and $M \in \mathrm{GMORK}_{\tilde{n},s}$.
If $M$ is explicit and convergent up to length 1 then $M$ is neither A-stable up to any length nor absolute A-stable for any length.

Proof

It is a direct consequence of Section IV.2.20 and Section IV.2.21.

For explicit methods we are able to find a simple expression of the stability functions. In general there exists a weaker result which simplifies the expression of the stability functions by using an explicit partition of the stages of a method.

**IV.2.23 Theorem : A-stability & Explicit partitions of the stages**

Let $(\tilde{n}, s) \in \overline{\mathbb{N}}^* \times \mathbb{N}^*$, $M = (\tau, w, \tilde{w}) \in \mathrm{GMORK}_{\tilde{n},s}$, $Q$ the resolvent matrices of $M$, $R$ the stability functions of $M$, $\mathcal{S}$ their domains of definitions, $P$ an explicit partition of the stages of $M$, and $n \in \mathbb{N}^*$ such that $n \leq \tilde{n}$.
We have :

$$\forall z \in \mathbb{C}^{[\![1,n]\!]},\ \det(Q_n(z)) = \prod_{\substack{J \in P \\ J = \{j_1, \dots j_{\#(J)}\}}} \det\left( I_{\#(J)} + \sum_{N=1}^{n} \frac{1}{N!} \sum_{1 \leq i_1 < \dots < i_N} \left( \prod_{k=1}^{N} z_{i_k} \right) \begin{bmatrix} w_{N,j_1,j_1} & \cdots & w_{N,j_1,j_{\#(J)}} \\ \vdots & \ddots & \vdots \\ w_{N,j_{\#(J)},j_1} & \cdots & w_{N,j_{\#(J)},j_{\#(J)}} \end{bmatrix} \right)$$

$$\mathcal{S}_n = \bigcap_{\substack{J \in P \\ J = \{j_1, \dots j_{\#(J)}\}}} \left\{ z \in \mathbb{C}^{[\![1,n]\!]} \mid \det\left( I_{\#(J)} + \sum_{N=1}^{n} \frac{1}{N!} \sum_{1 \leq i_1 < \dots < i_N} \left( \prod_{k=1}^{N} z_{i_k} \right) \begin{bmatrix} w_{N,j_1,j_1} & \cdots & w_{N,j_1,j_{\#(J)}} \\ \vdots & \ddots & \vdots \\ w_{N,j_{\#(J)},j_1} & \cdots & w_{N,j_{\#(J)},j_{\#(J)}} \end{bmatrix} \right) \neq 0 \right\}$$

Proof

Let $B$ an order of computation of $P$. Using the characterization of explicit partitions, one shows that each $w_N$ is similar to a block lower triangular matrix whose diagonal blocks are precisely the blocks of $P$. It follows that the sum appearing in the resolvent matrix is also block lower triangular with the same diagonal blocks.

The determinant of the resolvent matrix is therefore the product of the determinants its diagonal blocks, which yields the desired expression. The formula for the the domain of definition follows immediately.

These expressions are independant of the labelling of the stages in each $J \in P$ since relabelling amounts to conjugation by a permutation matrix, which leaves the determinant unchanged.

This theorem implies that the implicit blocks determine entirely the domain of definitions and the denominator of the stability functions.

Although A-stability is already a strong requirement, it can be strengthened further. In the previous subsection, we observed that the entries of $P\hat{y}$ behave like linear combinations of functions of the form $t \to t^k e^{\tilde{r}_l t}$, where $\tilde{r}$ are the different roots of the characteristic polynomial. As the real parts of these roots tend to $-\infty$, the solution decays to zero increasingly rapidly. If we wish this decay property to be reflected in the approximations, we may impose the stronger requirement of L-stability.

**IV.2.24 Definition : L-stability & Absolute L-stability**

Let $(\tilde{n}, s) \in \overline{\mathbb{N}}^* \times \mathbb{N}^*$, $M \in \mathrm{GMORK}_{\tilde{n},s}$, $R$ the stability functions of $M$, $\mathcal{S}$ their domains of definition, and $n \in \mathbb{N}^*$ such that $n \leq \tilde{n}$.

- $M$ is said to be L-stable for length $n$ if and only if M is A-stable up to length $n$ and :

$$\lim_{\min_{k\in[\![1,n]\!]}(\Re(z_k))\to+\infty} R_n(z) = 0$$

- $M$ is said to be absolute L-stable for length $n$ if and only if M is absolute A-stable for length $n$ and :

$$\lim_{\min_{k\in[\![1,n]\!]}(\Re(z_k))\to+\infty} R_n(z) = 0$$

Since A-stability or L-stability can be difficult to achieve, it is natural to consider weaker notions similar to A-($\alpha$) stability. The definition of A-($\alpha$) stability however is somewhat arbitrary, we therefore introduce a weaker and more refined property. Instead of requiring stability in a sectof of the complex plane, we require stability along a single half-line. This amounts to demanding stabilty for all step sizes when when the method is applied to a fixed linear differential equation function.

**IV.2.25 Definition : Half-line stability & Absolute half-line stability**

Let $(\tilde{n}, s) \in \overline{\mathbb{N}}^* \times \mathbb{N}^*$, $M \in \mathrm{GMORK}_{\tilde{n},s}$, $n \in \mathbb{N}^*$ such that $n \leq \tilde{n}$, $Y_n$ the approximation function for length $n$ of $M$, and $\mathcal{H}$ the set of step sizes of $M$.

- Let $z \in \mathbb{C}^{[\![1,n]\!]}$ such that, for all $k \in [\![1,n]\!]$, $\Re(z_k) \geq 0$.
  $M$ is said to be half-line stable for length $n$ in the direction of $z$ if and only if, if we consider the only linear differential equation function $f$ such that its characteristic polynomial has $-r$ as its roots, then for all $(t_0, y_{n,0}) \in \mathbb{R} \times \mathbb{R}^{[\![1,n]\!]\times[\![1,d]\!]}$, for all $h \in \mathbb{R}_+^*$ such that $(h)_{k\in\mathbb{N}} \in \mathcal{H}_{\infty,f,t_0,y_{n,0}}$, $\left(Y_{n,f,q(s+1)}(t_0, y_{n,0}, (h)_{q\in\mathbb{N}})\right)_{q\in\mathbb{N}}$ is bounded.
- Let $z \in \mathbb{C}^{[\![1,n]\!]}$ such that, for all $k \in [\![1,n]\!]$, $\Re(z_k) > 0$.
  $M$ is said to be absolute half-line stable for length $n$ in the direction of $z$ if and only if, if we consider the unique linear differential equation function $f$ such that its characteristic polynomial has $-r$ as its roots, then for all $(t_0, y_{n,0}) \in \mathbb{R} \times \mathbb{R}^{[\![1,n]\!]\times[\![1,d]\!]}$, for all $h \in \mathbb{R}_+^*$ such that $(h)_{k\in\mathbb{N}} \in \mathcal{H}_{\infty,f,t_0,y_{n,0}}$, $\left(Y_{n,f,q(s+1)}(t_0, y_{n,0}, (h)_{q\in\mathbb{N}})\right)_{q\in\mathbb{N}}$ converges towards 0.

**IV.2.26 Theorem : Characterization of half-line stability & Absolute half-line stability**

Let $(\tilde{n}, s) \in \overline{\mathbb{N}}^* \times \mathbb{N}^*$, $M \in \mathrm{GMORK}_{\tilde{n},s}$, $R$ the stability functions of $M$, $\mathcal{S}$ their domains of definition, and $n \in \mathbb{N}^*$ such that $n \leq \tilde{n}$.

1. Let $z \in \mathbb{C}^{[\![1,n]\!]}$ such that, for all $z \in [\![1,n]\!]$, $\Re(z_k) \geq 0$.
   $M$ is half-line stable in the direction of $z$ for length $n$ if and only if, for all $h \in \mathbb{R}_+^*$, if $hz \in \mathcal{S}_n$, then $\rho(R_n(hz)) \leq 1$.
2. Let $z \in \mathbb{C}^{[\![1,n]\!]}$ such that, for all $k \in [\![1,n]\!]$, $\Re(z_k) > 0$.
   $M$ is absolute half-line stable in the direction of $z$ for length $n$ if and only if, for all $h \in \mathbb{R}_+^*$, if $hz \in \mathcal{S}_n$, then $\rho(R_n(hz)) < 1$.

Proof

It is a similar proof to the proof of Section IV.2.10

The function $h \to hz$ defines a half-line starting at the point $(0, 0)$ in the direction of $z$.

**IV.2.27 Proposition : Properties of half-line stability**

Let $(\tilde{n}, s) \in \overline{\mathbb{N}}^* \times \mathbb{N}^*$, $M \in \mathrm{GMORK}_{\tilde{n},s}$, $n \in \mathbb{N}^*$ such that $n \leq \tilde{n}$, and $z \in \mathbb{C}^{[\![1,n]\!]}$.

1. If $M$ is (absolute) half-line stable in the direction of $z$ for length $n$ then, for all $\lambda \in \mathbb{R}_+^*$, it is also (absolute) half-line stable in the direction of $\lambda z$ for length $n$.
2. Let $\varphi \in \mathfrak{S}_{[\![1,n]\!]}$.
   If $M$ is (absolute) half-line stable in the direction of $z$ for length $n$ then $M$ is (absolute) half-line stable in the direction of $\left(z_{\varphi(1)}, ..., z_{\varphi(n)}\right)$ for length $n$.

Proof

1. Trivial.
2. It is a direct consequence of the symmetry of the stability functions.

In some situations, it may be useful for a method to be stable not along an entire half-line, but only on a portion of it. For example, when applying the method to a linear initial value problem with the step size restricted to a given range, stability is only required on the corresponding segment of the half-line.

# Conclusion

In this paper, we depart from the usual point of view on initial value problems, where all problems are assumed to be of order one, and instead consider initial value problems of arbitrary order. We revisit the foundations of Runge-Kutta methods and examine how they rewrite any initial value problem as a first-order system to approximate its solution, and how this rewriting process hinders the accuracy of Runge-Kutta methods.

Using Taylor's theorem and Gauss-Jacobi quadratures, we construct a generalization of Runge-Kutta methods capable of directly solving initial value problems of arbitrary order, multi-order Runge-Kutta methods. To rigorously define the approximations produced by these methods, we introduce the notion of the stage system of a method, and we generalize the classical result on the well-definedness of Runge-Kutta approximations to the multi-order case.

Unlike standard Runge-Kutta theory, we then analyze the internal structure of these approximations, revealing a rich structural framework. By interpreting the weight matrices of a method as an adjacency matrix defining a directed graph, we are able to define explicit methods without the inconsistencies found in the traditional definition of explicit Runge-Kutta methods. Moreover, we study how the structure of this digraph influences the stage system, introducing the notions of implicit blocks and implicit rank, which characterize the internal dependencies of the approximations. We also define parallel blocks, which characterizes the independence between subsystems. This structural theory allows us to design methods with the structure of our choice.

We investigate the accuracy of these methods by introducing a rigorous definition of the order of consistency of a multi-order Runge-Kutta method. We characterize methods with non-negative order of consistency and derive a system of equations characterizing positive order of consistency. We show how this system can be simplified by generalizing the classical condition that allows one to restrict attention to autonomous initial value problems.

We also show that the time variable can prove useful, as it provides a natural link between multi-order Runge-Kutta methods and Gauss-Jacobi quadratures through the solved system conditions. These conditions allow us to generalize and sharpen the classical upper bound of $2s$ on the order of consistency. We also study the error after the first step by generalizing the concept of convergence and analyzing the impact of the order of consistency. We derive weaker sufficient conditions for convergence than those usually assumed.

We study the different ways in which an initial value problem can be rewritten, considering three processes: overwriting, prolongation, and reduction. The notion of overwriting allows us to define rewrite-compliant methods, methods that can be rewritten without losing order of consistency. The concepts of prolongation and extension make it possible to solve initial value problems of any order, though it is preferable to avoid extension altogether, as it generalizes the rewriting of initial problems to ones of order one and leads to a loss of accuracy. Instead, we advocate the use of infinite methods. We also investigate the effects of these rewriting on the structure of multi-order methods. Altogether, this analysis allows us to address a more general definition of initial value problems, mixed-order initial value problems.

We derive several known methods and their generalizations, as well as a two new method, MORK4b and Multi-order CNb, together with their Runge-Kutta counterpart RK4b and CNB. We compare these methods and conclude, as expected, that a multi-order Runge-Kutta method is somewhat more costly but achieves significantly smaller errors. The analysis of computational cost remains theoretical: actual speed measurements would require an optimized implementation, whereas our implementation favors flexibility for experimentation. The performance of a method also depends on context, when evaluating $f$ is computationally expensive, the additional cost of computing a few Taylor polynomials becomes negligible in comparison. Real applications of these methods are therefore needed, definitive conclusions cannot be drawn from a non-optimized program. Given the large number of parameters (the methods considered, the order of the problem, the initial values, the function $f$, and the system dimension) we restrict numerical tests to a single problem of order six.

We observe that for step sizes that are not sufficiently small, Runge-Kutta methods may outperform multi-order methods in the long term because of their stability. Although the pursuit of higher local order motivated the construction of multi-order methods, their real advantage may lie in their additional flexibility. This flexibility can potentially be used to design methods whose stability properties are improved while maintaining orders of consistency comparable to, or exceeding, that of classical Runge-Kutta methods.

We therefore extend the notion of A-stability to multi-order methods and introduce a corresponding stability function. We characterize A-stability in terms of this stability function and establish properties of the function that allow us to generalize the second Dahlquist barrier and describe how the structure of a method affects its A-stability. We further generalize the concept of L-stability and introduce the new notion of half-line stability, a weaker yet more refined version of A($\alpha$)-stability. Since multi-order Runge-Kutta methods act on $\mathbb{C}^{[\![1,n]\!]}$ rather than on $\mathbb{C}$, stability is no longer naturally described by regions of the complex plane. Instead, stability is required along individual half-lines, and more general domains can be constructed by combining half-lines.

We provide a general implementation of both multi-order and classical Runge-Kutta methods, written in pseudocode and in Rust. The implementation is currently single-threaded and does not exploit parallel blocks, but otherwise fully incorporates the structural theory developed in this paper to minimize computational cost. A multi-threaded version of this program is discussed.

Although already extensive, this paper only establishes the fundamental theory of multi-order Runge-Kutta methods, many research directions remain open. Extensive testing is necessary, and a multi-threaded implementation would be valuable. Future work includes applying the stability concepts developed here to derive more stable methods, particularly explicit ones, extending the notions of B-stability and B-series to the multi-order case, deriving conditions for higher orders of consistency, generalizing embedded methods and adaptive step-size techniques, investigating the differentiation of $f$ to increase the order of the initial value problem, and extending the framework to symplectic methods and general linear methods.

# Appendix

## A Implementation

In this section we create an implementation of Runge-Kutta methods and general multi-order Runge-Kutta methods. The algorithms are here written in pseudocode, the actual program is written in Rust and can be found on the author's GitHub [14]. The syntax of the pseudocode is inspired by the Rust programming language, we hence use enums, structs and traits. We now go over the syntax of the pseudocode.

Enums, or enumerations, is a data type that can take a finite number of predefined variants. An enum has a name and variants. A variant can have values attached to it. For example, we can define an enum named `Input` that have for variants `StringInput` or `NoInput`. We attach to the variant `StringInput` a value of type `String`, and to `NoInput`, nothing.

```
1 ENUM Input
2   StringInput(String)
3   NoInput
4 ENDENUM
```

A variable can take a variant of an enum as its value.

```
1 a := Input::StringInput("Hello World!")
2 b := Input::NoInput
```

And we can distinguish between the variants of a variable using a match statement.

```
1 MATCH a
2   CASE StringInput(string_value)
3     PRINT string_value
4   ENDCASE
5   CASE NoInput
6     PRINT "No input has been given"
7   ENDCASE
8 ENDMATCH
```

To handle variants which have not been treated in other **`CASE`** statements and which we do not need the attached value of, we can use **`CASE`** `_`.

Structs are very similar to objects in programming, in the sense that they have fields and methods. For example, we can define a 2D point :

```
1 STRUCT Point
2   x : Float
3   y : Float
4 ENDSTRUCT
```

We access a field of a struct using `struct_instance.field_name`. We also are able to add methods to a struct. The argument `self` is a reference to the instance of the struct used to access its fields.

```
1 IMPLEMENT Point
2   METHOD magnitude(self) RETURNS Float
3     RETURN sqrt(self.x**2 + self.y**2)
4   ENDMETHOD
5   METHOD translation(self, dx : Float, dy : Float)
6     self.x := self.x + dx
7     self.y := self.y + dy
8   ENDMETHOD
9 ENDIMPLEMENT
```

We create an instance of a struct using `StructName(field1 : value1,field2 : value2,...)`.

```
1 x = Point(x:0.5,y:1)
2 x.translation(-1,2)
3 x.magnitude()
```

Trait are used to define behavior shared by different structs. In practice, defining a trait means defining a set of methods that a struct with this trait must implement. For example, many different structs implement the `print` function thanks to the `Display` trait.

```
1 TRAIT Display
2   METHOD print(self) RETURNS ()
3 ENDTRAIT
```

A method or function that returns `()` returns nothing. We can implement this trait on the `Point` type :

```
1 IMPLEMENT Display FOR Point
2   METHOD print(self):
3     PRINT self.x, self.y
4   ENDMETHOD
5 ENDIMPLEMENT
```

We now state our objectives.

We wish to write an implementation of general multi-order Runge-Kutta methods and Runge-Kutta methods able to approximate any initial value problem (see mixed-order case in Section III.2.36).

To simplify the program, we only consider infinite general multi-order Runge-Kutta methods, although it would be possible to implement the rewriting for initial value problems of higher order. We implement order of computations, but not parallel blocks and multithreading to keep the program simple. We however discuss at the end how a multithreaded version would work.

To numerically solve the systems of implicit SCCs we take the simplest approach and use Picard's iteration method. More complex schemes such as Newton-based schemes can be used, they however require to know the differential of $f$, or at least an approximation, which makes their implementation difficult.

To define a general multi-order Runge-Kutta method we need :
- The nodes $\tau$.
- The secondary weights function, the function which, given a rank $N \in \mathbb{N}^*$, returns $\tilde{w}_N$.
- The main weights function, the function which, given a rank $N \in \mathbb{N}^*$, returns $w_N$.
- The maximum weight digraph.

To define a Runge-Kutta method we need :
- The nodes $\tau$.
- The weights $w_1$.

We could avoid requiring the maximum weight digraph in the definition of a multi-order method and instead construct the weight digraph at a fixed length $n \in \mathbb{N}^*$. However, if the method is then applied to an initial value problem of order greater than $n$, we would need to compute the weight digraph at a higher length, together with a corresponding order of computation. This would complicate the implementation.

We will create two different structs, `RK` for Runge-Kutta methods and `GMORK` for general multi-order Runge-Kutta methods. We want both of them to implement a trait `Solver`. This trait consists of a method named `approximate` analogous to $Y$ which takes as input an initial instant `t0`, initial values `y0`, a step size `h`, a differential equation function `f`, and outputs the approximations of the method. We now define the trait `Solver` :

```
1 TRAIT Solver
2   METHOD approximate(
3     self, t0 : Float, y0 : Float , h : Float,
4     f : (Float, Vector of Vector of Float) -> Vector of Float
5   ) RETURNS Vector of Vector of Float
6 ENDTRAIT
```

Before we implement this trait, we need to write functions related to graph theory which are necessary to find an order of computation of a method.

We first write an algorithm which detects the strongly connected components of a maximum weight digraph, then one able to contract a dependency digraph, and finally an algorithm able which topologically sorts a digraph.

To detect the SCCs of a maximum weight digraph we use Gabow's algorithm [15]. We define a function named `Gabow_SCC` which takes as input the adjacency matrix of a digraph and outputs a vector of SCC, each represented by the vector of vertices in that SCC. The algorithm has been modified to use an iterative approach instead of a recursive one, and instead of using clever tricks to know whether a vertex has already been assigned a strongly connected component or has not been visited yet or is in the current path, we simply use an enum :

```
1 ENUM NodeState
2   NotVisited
3   InPath(integer)
4   InSCC(integer)
5 ENDENUM
```

In the code an adjacency matrix is represented using a `Vector of Vector of Bool`. For an explanation of this algorithm, please refer to the original paper by Gabow [15].

```
1  FUNCTION Gabow_SCC(graph: Vector of Vector of Bool) RETURNS Vector of Vector
   of integers
2    // Initialization
3    S := empty Vector of Integer
4    B := empty Vector of Integer
5    s := length(graph)
6    IF s = 0 THEN
7      RETURN Empty Vector
8    ENDIF
9    state := Vector of size s filled with NodeState::NotVisited
10   tree := Vector containing the Vector [0,1,…,s - 1] // DFS stack
11   path := empty Vector of Integer // Currently explored path
12   scc_count := 0
13   finished := false
14
15   // Start of the algorithm
16   WHILE NOT finished DO
17     IF tree[length(tree) - 1] is not empty THEN
18       j = POP(tree[length(tree) - 1]) // next node to explore
19       MATCH state[j]
20         CASE NotVisited
21           PUSH(path, j)
22           PUSH(S,j)
23           state[j] := NodeState::InPath(length(S) - 1)
24           PUSH(B, length(S) - 1)
25           neighbors := empty Vector of Integer
26           FOR j1 FROM 0 TO s - 1 DO
27             IF graph[j][j1] THEN
28               PUSH(neighbors, j1)
29             ENDIF
30           ENDFOR
31           PUSH(tree, neighbors)
32         ENDCASE
33
34         CASE InPath(j1)
35           WHILE j1 < B[length(B) - 1] DO
36             POP(B)
37           ENDWHILE
38         ENDCASE
39
40         CASE InSCC(_)
41           // Do nothing
42         ENDCASE
43       ENDMATCH
44
45     ELSE
46       POP(tree)
47       IF path is not empty THEN
48         j := POP(path)
```

```
49       MATCH state[j]
50         CASE InPath(k1)
51           IF j1 = B[length(B) - 1] THEN
52             POP(B)
53             WHILE j1 < length(S) DO
54               state[POP(S)] = NodeState::InSCC(scc_count)
55             ENDWHILE
56             scc_count := scc_count + 1
57           ENDIF
58         ENDCASE
59         CASE _
60           // Do nothing
61         ENDCASE
62       ENDMATCH
63     ELSE
64       finished := true
65     ENDIF
66   ENDIF
67 ENDWHILE
68 // Lists the SCCs
69 scc := Vector of size scc_count filled with Empty Vector of integers
70 FOR j FROM 0 to s - 1 DO
71   MATCH state[j]
72     CASE InSCC(scc_index)
73       PUSH(scc[scc_index], j)
74     ENDCASE
75     CASE _
76       // should not happen
77     ENDCASE
78   ENDMATCH
79 ENDFOR
80 RETURN scc
81 ENDFUNCTION
```

Now that we have the list of strongly connected components, we need to contract the dependency graph. We now define a function `contract_graph` which, given an adjacency matrix and a partition of its vertices, returns the contraction of said digraph according to the partition without the loops.

```
 1 FUNCTION contract_graph(
 2   graph: Vector of Vector of Bool, partition : Vector of Vector of Integer
 3 ) RETURNS Vector of Vector of Bool
 4   s := length(graph)
 5   p := length(partition)
 6   inverse_map := Vector of size s filled with 0
 7   FOR i FROM 0 TO s - 1 DO
 8     FOR j IN partition[i] DO
 9       inverse_map[j] := i
10     ENDFOR
11   ENDFOR
```

```
12 contracted_graph := Vector of size p filled with Vector of size p filled with false
13 FOR j1 FROM 0 TO s - 1 DO
14   FOR j2 FROM 0 TO s - 1 DO
15     partition1 := inverse_map[j1]
16     partition2 := inverse_map[j2]
17     IF partition1 != partition2 AND graph[j1][j2] THEN
18       contracted_graph[partition1][partition2] := true
19     ENDIF
20   ENDFOR
21 ENDFOR
22 RETURN contracted_graph
23 ENDFUNCTION
```

We define a function `topological_sort` which topologically sorts the adjacency matrix of a digraph with no non-empty closed diwalks. This function takes as input an adjacency matrix and outputs a list of vertices in the order it must be relabeled to be strictly lower triangular.

We first define a helper function named `source_nodes` which takes the adjacency matrix of a digraph as input and outputs the nodes with no incoming arcs.

```
1 FUNCTION source_nodes(graph: Vector of Vector of Bool) RETURNS Vector of Integer
2   s := length(graph)
3   S := empty Vector of Integer
4   FOR j FROM 0 TO s - 1 DO
5     IF graph[j] contains only false THEN
6       PUSH(S,j)
7     ENDIF
8   ENDFOR
9   RETURN S
```

With the help of this function we write the `topological_sort` function [10].

```
1 FUNCTION topological_sort(dag : Vector of Vector of Bool) RETURNS Vector of Integer
2   s := length(dag)
3   SN := source_nodes(dag)
4   order := empty Vector of Integer
5   WHILE SN is not empty DO
6     node := POP(SN)
7     PUSH(order, node)
8     FOR j FROM 0 TO s - 1 DO
9       dag[j][node] := false
10    ENDFOR
11    SN := source_nodes(dag)
12    SN := SN without nodes contained in order
13  ENDWHILE
14  RETURN order
15 ENDFUNCTION
```

To distinguish between explicit and implicit SCCs we use an enum.

```
1 ENUM SCC
2   Implicit(Vector of Integer,Vector of Integer)
3   Explicit(Integer)
4 ENDENUM
```

The first argument of `Implicit` is $J$, the list of stages of the SCC, and the second one is its complementary, $[\![1, s]\!] \setminus J$. For `Explicit`, the argument is simply the only stage of the SCC.

We have everything we need to implement Runge-Kutta methods and general multi-order Runge-Kutta methods. We start with general multi-order Runge-Kutta methods. We define the `GMORK` struct.

```
 1 STRUCT GMORK
 2   nodes : Vector of Float
 3   main_weights : Vector of Vector of Vector of Float
 4   main_weights_f : Integer -> Vector of Vector of Float
 5   secondary_weights : Vector of Vector of Vector of Float
 6   secondary_weights_f : Integer -> Vector of Vector of Float
 7   computation_order : Vector of SCC
 8   s : Integer
 9   length : Integer
10   implicit_ranks : Vector of Vector of bool
11   factorial : Vector of Integer
12   h_powers : Vector of Float
13 ENDSTRUCT
```

There are many fields to this struct, we examine them one by one.

- `nodes` contains the nodes of the method.
- `main_weights` is the list of currently stored main weights.
- `main_weights_f` is the main weights function.
- `secondary_weights` is the list of currently stored secondary weights.
- `secondary_weights_f` is the secondary weight function.
- `computation_order` contains an order of computation of the method.
- `s` is the number of points of the method.
- `length` is the length of the method the struct currently stores the weights of.
- `implicit_ranks[N-1][j]` is true if and only if $N$ is an implicit rank at stage $j$.
- `factorial` stores the values of the factorial.
- `h_powers` stores the power of `h`.

We define a helper function `implicit_rank` which, given the weights of a method at a certain rank and the order of computation of the same method, returns the implicit ranks at this rank for all stages.

```
1 FUNCTION implicit_rank(weights,computation_order) RETURNS Vector of Bool
2   ranks := Vector of size length(weights) - 1 filled with false
3   FOR scc IN computation_order
4     MATCH scc
5       CASE Implicit(J)
6         FOR j IN J DO
7           FOR j1 IN J DO
8             IF main_weights[0][j][j1] != 0 THEN
```

```
 9               ranks[j] := true
10             ENDIF
11           ENDFOR
12         ENDFOR
13       ENDCASE
14       CASE Explicit
15         // do nothing
16       ENDCASE
17     ENDMATCH
18   ENDFOR
19   RETURN ranks
20 ENDFUNCTION
```

We define a new method able to create an instance of `GMORK`.

```
 1 IMPLEMENT GMORK
 2   METHOD new(
 3     main_weights_function: Integer -> Vector of Vector of Floats,
 4     secondary_weights_function: Integer -> Vector of Vector of Floats,
 5     nodes : Vector of Floats
 6     maximum_weight_graph : Vector of Vector of Bool
 7   ) RETURNS GMORK
 8     scc := Gabow_SCC(maximum_weight_graph)
 9     contracted_graph := contract_graph(maximum_weight_graph,scc)
10     order := topological_sort(contracted_graph)
11     computation_order := empty Vector of SCC
12     FOR j IN order DO
13       IF length(scc[j]) = 1 AND NOT maximum_weight_graph[scc[j][0]][scc[j][0]] THEN
14         PUSH(computation_order, SCC::Explicit(scc[j][0]))
15       ELSE
16         PUSH(computation_order, SCC::Implicit(scc[j]))
17       ENDIF
18     ENDFOR
19     RETURN GMORK(
20       nodes,
21       main_weights : [main_weights_function(1)],
22       main_weights_f : main_weights_function,
23       secondary_weights : [secondary_weights_function(1)],
24       secondary_weights_f : secondary_weights_function,
25       computation_order,
26       s : length(self.main_weights[0]) - 1,
27       length : 1,
28       implicit_ranks : [implicit_ranks(self.main_weights[0]),computation_order],
29       factorial : [1,1],
30       h_powers : [1,0],
31       h : 0,
32     )
33   ENDMETHOD
34 ENDIMPLEMENT
```

The `new` function initializes a method to a length of 1, which means it can only solve initial value problems of order 1. We thus define a method which increases the length of the multi-order method. This method named `set_minimum_length` takes an integer `n` as input, and if `method.length` is strictly lower than `n` it computes and stores the constants necessary to solve any initial value problem of order `n` or less.

```
IMPLEMENT GMORK
  METHOD set_minimum_length(self, n : Integer) RETURNS ()
    IF self.length >= n THEN
      RETURN
    ENDIF
    APPEND(self.factorial,[(self.length + 1)!, … ,n!])
    APPEND(
      self.main_weights, [self.main_weights_f(self.length+1), … ,self.main_weights_f(n)]
    )
    APPEND(
      self.secondary_weights,
      [self.secondary_weights_f(self.length+1), … ,self.secondary_weights_f(n)]
    )
    APPEND(
      self.h_powers, [self.h**(self.length+1), … ,self.h**n]
    )
    APPEND(
      self.implicit_ranks, [
        implicit_ranks(self.main_weights_f(self.length)),
        …
        implicit_ranks(self.main_weights_f(n-1))
      ]
    )
  ENDMETHOD
ENDIMPLEMENT
```

To make the code more readable we move the algorithm that implements Picard's iteration method to a helper function named `picard_iterations`. Measuring the convergence of Picard's iteration using the difference of the approximations $y_{n,j,N}$ is equivalent to measuring the convergence of $\frac{h^N}{N!} \sum_{j=1}^{s} w_{N,j,j'} f\big(t + \tau_j h, y_{n,j}\big)$. It is less costly to directly measure the convergence of the evaluations, it also has the advantage of avoiding the factor $\frac{h^N}{N!}$ which artificially lowers the difference.

```
IMPLEMENT GMORK
  METHOD picard_iterations(
    self, t : Float, h : Float, y : Vector of Vector of Vector of Float, y0 : Vector
    of Vector of Float, F : Vector of Vector of Float, f : (Float, Vector of Vector
    of Float) -> Vector of Float, J : Vector of Integer
  )
    threshold := … // Float defined by user
    minimum := … // Integer defined by user, minimum number of iterations
    maximum := … // Integer defined by user, maximum number of iterations
    constant := y
    iteration := 0
    d := threshold + 1
```

```
11     WHILE iteration < minimum OR (d > threshold AND iteration < maximum) DO
12       iteration := iteration + 1
13       d := 0
14       FOR j IN J DO
15         f_cache := f(t + self.nodes[j] * h, y[j])
16         d := MAX(d, max distance between F[j] and f_cache)
17         F[j] := f_cache
18       ENDFOR
19       FOR j IN J DO
20         FOR k FROM 0 TO length(y0) - 1 DO
21           FOR N FROM 0 TO length(y[k]) - 1 such that self.implicit_ranks[N][j] DO
22             sum := 0
23             FOR j1 IN J DO
24               sum = sum + self.main_weights[N][j][j1] * F[j1][k]
25             ENDFOR
26             y[j][k][N] = y[j][k][N] + self.h_powers[N+1] / self.factorial[N+1] * sum
27           ENDFOR
28         ENDFOR
29       ENDFOR
30     ENDWHILE
31   ENDMETHOD
32 ENDIMPLEMENT
```

We implement the Solver trait.

```
1  IMPLEMENT Solver FOR GMORK
2    METHOD approximate(
3      self,t0 : Float, y0 : Float , h : Float,
4      f : (Float, Vector of Vector of Float) -> Vector of Float
5    ) RETURNS Vector of Vector of Float
6      // Update powers of h if necessary
7      IF h != self.h_powers[1] THEN
8        self.h_powers := [h**0, … , h**self.length]
9      ENDIF
10
11     // Verifies if the stored length is higher than the order of the problem
12     FOR k FROM 0 TO length(y0) - 1 DO
13       IF length(y0[k]) > self.length THEN
14         self.set_minimum_length(length(y0[k]))
15       ENDIF
16     ENDFOR
17
18     F := Vector of size self.s filled with y0[0]
19     y := Vector of size self.s + 1 filled with y0
20
21     // Adds initial values
22     FOR j FROM 0 TO s DO
23       FOR k FROM 0 TO length(y0) - 1 DO
24         FOR N FROM 0 TO length(y0[k]) - 1 DO
```

```
25                  FOR N1 FROM 1 TO N DO
26                    y[j][k][n] := y[j][k][n] + self.secondary_weights[N][N1][j] * self.h_powers[N-
                      N1] * y0[k][N1]
27                  ENDFOR
28                ENDFOR
29              ENDFOR
30            ENDFOR
31
32            // Computes evaluations
33            FOR scc IN self.computation_order DO
34              MATCH scc
35                CASE Explicit(j)
36                  FOR k FROM 0 TO length(y0) -1 DO
37                    FOR N FROM 0 TO length(y0[k]) - 1 DO
38                      sum := 0
39                      FOR j1 FROM 0 TO self.s DO
40                        sum = sum + self.main_weights[N][j][j1] * F[j1][k]
41                      ENDFOR
42                      y[j][k][N] = y[j][k][N] + self.h_powers[N + 1] / self.factorial[N + 1] * sum
43                    ENDFOR
44                  ENDFOR
45                  IF j != self.s THEN
46                    F[j] = f(t + self.nodes[j] * h, y[j])
47                  ENDIF
48                ENDCASE
49
50                CASE Implicit(J,comp_J)
51                  FOR j IN J DO
52                    FOR k FROM 0 TO length(y0) -1 DO
53                      FOR N FROM 0 TO length(y0[k]) DO
54                        sum := 0
55                        FOR j1 IN comp_J DO
56                          sum = sum + self.main_weights[N][j][j1] * F[j1][k]
57                        ENDFOR
58                        y[j][k][N] = y[j][k][N] + self.h_powers[N + 1]/self.factorial[N + 1] * sum
59                      ENDFOR
60                    ENDFOR
61                  ENDFOR
62                  self.picard_iterations(t0, h, y, y0, F, f, J)
63                ENDCASE
64              ENDMATCH
65            ENDFOR
66            RETURN y[self.s]
67          ENDMETHOD
68        ENDIMPLEMENT
```

The implementation of general multi-order Runge-Kutta methods is done.

We implement Runge-Kutta methods. We start by defining the struct `RK`.

```
1 STRUCT RK
2   s : Integer
3   nodes : Vector of Float
4   weights : Vector of Vector of Float
5   computation_order : Vector of SCC
6 ENDSTRUCT
```

We define a new method able to create an instance of RK.

```
 1 IMPLEMENT RK
 2   METHOD new(weights: Vector of Vector of Float, nodes : Vector of Float) RETURNS RK
 3     method := new RK
 4     s := length(weights) - 1
 5     weight_graph := Vector of size method.s + 1 filled with Vector of size method.s
       + 1 filled with false
 6     FOR j FROM 0 TO method.s DO
 7       FOR j1 FROM 0 TO method.s - 1 DO
 8         IF weights[j][j1] THEN
 9           weight_graph[j][j1] := true
10         ENDIF
11       ENDFOR
12     ENDFOR
13     scc := Gabow_SCC(weight_graph)
14     contracted_graph := contract_graph(weight_graph,scc)
15     order := topological_sort(contracted_graph)
16     computation_order := empty Vector of SCC
17     FOR j IN order DO
18       IF length(scc[j]) = 1 AND NOT weight_graph[scc[j][0]][scc[j][0]] THEN
19         PUSH(computation_order, SCC::Explicit(scc[j][0]))
20       ELSE
21         PUSH(computation_order, SCC::Implicit(scc[j]))
22       ENDIF
23     ENDFOR
24     RETURN RK(
25       s,
26       nodes,
27       weights,
28       computation_order,
29     )
30   ENDMETHOD
31 ENDIMPLEMENT
```

We define a helper function which implements Picard's iteration method, picard_iterations.

```
1 IMPLEMENT RK
2   METHOD picard_iterations(t : Float, h : Float, F : Vector of Vector of Float,
    J : Vector of Integer, f : (Float, Vector of Vector of Float) -> Vector of
    Float, y : Vector of Vector of Vector of Float, y0 : Vector of Vector of Float)
3     threshold := … // Float defined by user
4     minimum := … // Integer defined by user, minimum number of iterations
```

```
5     maximum := … // Integer defined by user, maximum number of iterations
6     constant := y
7     iteration := 0
8     d := threshold + 1
9     WHILE iteration < minimum OR (d > threshold AND iteration < maximum) DO
10      iteration := iteration + 1
11      d := 0
12      FOR j IN J DO
13        f_cache := f(t + self.nodes[j] * h, y[j])
14        d := MAX(d, max distance between F[j] and cache)
15        F[j] := f_cache
16      ENDFOR
17      FOR j IN J DO
18        FOR k FROM 0 TO length(y0) - 1 DO
19          FOR N FROM length(y0[k]) - 1 TO 1 STEP −1 DO
20            sum := 0
21            FOR j1 IN J DO
22              sum = sum + self.weights[j][j1] * y[j1][k][N-1]
23            ENDFOR
24            y[j][k][N] = constant[j][k][N] + h * sum
25          ENDFOR
26          sum := 0
27          FOR j1 IN J DO
28            sum := sum + self.weights[j][j1] * F[j1][k]
29          ENDFOR
30          y[j][k][0] = constant[j][k][0] + h * sum
31        ENDFOR
32      ENDFOR
33    ENDWHILE
34  ENDMETHOD
35 ENDIMPLEMENT
```

We implement the `Solver` trait for `RK`.

```
1 IMPLEMENT Solver FOR RK
2   METHOD approximate(
3     self, t0 : Float, y0 : Float , h : Float,
4     f : (Float, Vector of Vector of Float) -> Vector of Float
5   ) RETURNS Vector of Vector of Float
6     F := Vector of size self.s filled with y0[0]
7     y := Vector of size self.s + 1 filled with y0
8     FOR scc IN self.computation_order DO
9       MATCH scc
10        CASE Explicit(j)
11          FOR k FROM 0 TO length(y0) DO
12            sum := 0
13            FOR j1 FROM 0 TO self.s - 1 DO
14              sum = sum + self.weights[j][j1] * F[j1][k]
15            ENDFOR
```

```
16              y[j][k][0] = y[j][k][0] + h * sum
17              FOR N FROM 1 T0 length(y0[k]) DO
18                sum := 0
19                FOR j1 FROM 0 TO self.s - 1 DO
20                  sum = sum + self.weights[j][j1] * y[j1][k][N - 1]
21                ENDFOR
22                y[j][k][N] = y[j][k][N] + h * sum
23              ENDFOR
24              ENDFOR
25            IF j != self.s THEN
26              F[j] = f(t + self.nodes[j] * h, y[j])
27            ENDIF
28          ENDCASE
29
30          CASE Implicit(J,comp_J)
31            FOR j IN J DO
32              FOR k FROM 0 TO length(y0) - 1 DO
33                sum := 0
34                FOR j1 IN comp_J DO
35                  sum = sum + self.weights[j][j1] * F[j1][k]
36                ENDFOR
37                y[j][k][0] = y[j][k][0] + h * sum
38                FOR N FROM 1 TO length(y0[k]) - 1 DO
39                  sum = 0
40                  FOR j1 IN comp_J DO
41                    sum = sum + self.weights[j][j1] * y[j1][k][N - 1]
42                  ENDFOR
43                  y[j][k][N] = y[j][k][N] + h * sum
44                ENDFOR
45                ENDFOR
46            ENDFOR
47            self.picard_iterations(t, h, F, J, f, y, y0)
48          ENDCASE
49        ENDMATCH
50      ENDFOR
51      RETURN y[self.s]
52    ENDMETHOD
53  ENDIMPLEMENT
```

The implementation of Runge-Kutta methods is done.

We now discuss how a multi-threaded version of this program would be implemented.

The only required modification concers how the program processes the SCCs of a method. They cannot be processes in arbitrary order since some SCCs depend on others. The procedure therefore requires a digraph traversal algorithm. The principle is straightforward. The program maintains a queue containing the SCCs assigned to threads. When a SCC has been processed, the SCCs whose dependecis have now been satisfied are added to the queue. Initially, the queue contains all SCCs with no predecessors in the dependency graph, since these can be computed directly.

We previously considered an ideal scenario where the number of available threads is equal to the number of SCCs, which is not always the case. Since we are not always able to directly process a SCC, two queues are needed, a "to process" queue and a "processing" queue. The algorithm is the same as before except that the computable SCCs are first pushed to the "to process" queue, and when a thread becomes available, the program pops a SCC from the "to process" queue and pushes it to the "processing queue".

However, some SCCs may be more important than others, and if these are processed too late, it may slow down the program. For instance, a long sequence of SCCs ($J_p$ needs $J_{p-1}$, $J_{p-1}$ needs $J_{p-2}$, ..., $J_2$ needs $J_1$) is computationally intensive, therefore $J_1$ should be processed first. If $J_1$ was processed last then, for a long time, a single thread would be working on the sequence $J_1, J_2, ..., J_p$ while the other threads are idle. With this in mind, we may want to turn the "to process" queue into a priority queue, hence, when a thread is available, the SCC with the highest priority gets pushed to the "processing" queue.

We want to assign each SCC a priority. Let $J$ a SCC. The priority of $J$ should be based on the number of sequential task that starts by computing $J$, tasks that cannot be parallelized. In terms of dependency graph, such sequences of SCCs are diwalks starting with $J$. We want to know for each SCC the maximum computational cost of a diwalk which starts with this SCC.

We need to define the computational cost of diwalk. Since a diwalk is a sequence of tasks, one way to estimate the computational cost of a diwalk is to estimate the computational cost of each SCC in the diwalk and add them together.

Finally, we evaluate the computational cost of a SCC. Adding the evaluations of $f$ already computed can be considered negligeable compared to the cost of evaluating $f$, which may be very costly depending on the context. The computational cost of a SCC should hence be based one the number of evaluations of $f$ the program has to carry to solve a SCC.

If $J$ is an explicit SCC, its cost is 1 since the program only computes the evaluation of the only stage in the SCC. If $J$ is an implicit SCC and we assume the program carries out $a \in \mathbb{N}^*$ Picard iterations to solve the system of size $\#(J)$, then the number of necessary evaluation of the function is $a\#(J)$. The cost of an implicit SCC $J$ is hence $a\#(J)$. Since $a$ cannot be known in advance, it should be replaced with an estimate.

We now seek an algorithm that finds in the dependency graph the diwalk starting at $J$ with maximal cost. Since the cost of a diwalk is the sum of the cost of each SCC, it coincides with the length of the diwalk in the dependency graph, where all the SCCs have been weighted with their cost.

A vertex is called a sink if it has no successors. Any longest diwalk necessarily ends at a sink, otherwise it could be extended further. Hence, given a SCC, we aim to determine the longest diwalk starting at $J$ and ending at a sink. This problem can be reduced to a shortest-path problem on which Dijkstra's algorithm [16] may be applied.

We first introduce an auxiliary SCC with a cost of 0 and add arcs from each sink to this new vertex. The longest diwalk ending at a sink is then equal to the longest diwalk ending at this auxiliary vertex.

Since path algorithms typically use weights on arcs rather than on vertices, we transfer the cost of each SCC to its outgoing arcs by assigning to each arc the cost of its tail.

Our goal is to compute a longest path. By multiplying all arc weights by −1, this becomes a shortest-path problem. We may then apply a modified Dijkstra algorithm that takes as input the vector of SCC costs and the adjacency matrix of the dependency graph, and returns the priority of each SCC.

```
1  FUNCTION priority(individual_cost: Vector of Integer,graph: Vector of Vector of
   Bool) RETURNS Vector of Integer
2    s := length(graph)
3    current_priority := Vector of Integer of length s filled with cost[j] at index j
     if it has no successors, 0 otherwise
4    visited := Vector of Bool of length s filled with false
5    WHILE visited contains false DO
6      scc := a non visited SCC which maximizes current_priority over non visited SCCs
7      FOR j FROM 0 TO s - 1 DO
8        IF graph[scc][j] AND individual_cost[j] + current_priority[scc] >
         current_priority[scc] THEN
9          current_priority[j] = individual_cost[j] + current_priority[scc]
10       ENDIF
11     ENDFOR
12     visited[scc] := true
13   ENDWHILE
14   RETURN current_priority
15 ENDFUNCTION
```

# B Order conditions

Let $(\tilde{n}, s) \in \overline{N}^* \times \mathbb{N}^*$, $M = (\tau, w, \tilde{w}) \in \mathrm{MORK}_{\tilde{n},s}$, $\xi$ the approximations times of $M$, $n \in \overline{\mathbb{N}}^*$ such that $n \leq \tilde{n}$, $j \in [\![1, s+1]\!]$, and $N \in \mathbb{N}^*$ such that $N \leq n$. The conditions for order of consistency 1 through 4 for length $n$ at rank $N$ and stage $j$ are :

Order of consistency 1 :

$$\sum_{j'=1}^{s} w_{N,j,j'} = \xi_{N,j}^{N} \tag{1}$$

Order of consistency 2 :

$$\sum_{j'=1}^{s} w_{N,j,j'} \tau_{j'} = \xi_{N,j}^{N+1} \frac{1}{1+N} \tag{1}$$

$$\sum_{j'=1}^{s} w_{N,j,j'} \sum_{j''=1}^{s} w_{1,j',j''} = \xi_{N,j}^{N+1} \frac{1}{1+N} \tag{2}$$

$$\forall N' \in [\![2, n]\!], \quad \sum_{j'=1}^{s} w_{N,j,j'} \xi_{N',j'} = \xi_{N,j}^{N+1} \frac{1}{1+N} \tag{3}$$

Order of consistency 3 :

$$\sum_{j'=1}^{s} w_{N,j,j'} \tau_{j'}^{2} = \xi_{N,j}^{N+2} 2 \frac{N!}{(2+N)!} \tag{1}$$

$$\sum_{j'=1}^{s} w_{N,j,j'} \tau_{j'} \sum_{j''=1}^{s} w_{1,j',j''} = \xi_{N,j}^{N+2} 2 \frac{N!}{(2+N)!} \tag{2}$$

$$\sum_{j'=1}^{s} w_{N,j,j'} \sum_{j''=1}^{s} w_{1,j',j''} \tau_{j''} = \xi_{N,j}^{N+2} \frac{N!}{(2+N)!} \tag{3}$$

$$\sum_{j'=1}^{s} w_{N,j,j'} \sum_{j''=1}^{s} w_{1,j',j''} \sum_{j'''=1}^{s} w_{1,j'',j'''} = \xi_{N,j}^{N+2} \frac{N!}{(2+N)!} \tag{4}$$

$$\sum_{j'=1}^{s} w_{N,j,j'} \left( \sum_{j''=1}^{s} w_{1,j',j''} \right)^{2} = \xi_{N,j}^{N+2} 2 \frac{N!}{(2+N)!} \tag{5}$$

$$\text{if } n > 1, \quad \sum_{j'=1}^{s} w_{N,j,j'} \sum_{j''=1}^{s} w_{2,j',j''} = \xi_{N,j}^{N+2} 2 \frac{N!}{(2+N)!} \tag{6}$$

$$\forall N' \in [\![2, n]\!], \quad \sum_{j'=1}^{s} w_{N,j,j'} \tau_{j'} \xi_{N',j'} = \xi_{N,j}^{N+2} 2 \frac{N!}{(2+N)!} \tag{7}$$

$$\forall (N', N'') \in [\![2, n]\!]^{2}, \quad \sum_{j'=1}^{s} w_{N,j,j'} \xi_{N',j'} \xi_{N'',j'} = \xi_{N,j}^{N+2} 2 \frac{N!}{(2+N)!} \tag{8}$$

$$\forall N' \in [\![2, n]\!], \quad \sum_{j'=1}^{s} w_{N,j,j'} \xi_{N',j'} \sum_{j''=1}^{s} w_{1,j',j''} = \xi_{N,j}^{N+2} 2 \frac{N!}{(2+N)!} \tag{9}$$

$$\forall N' \in [\![2,n]\!],\ \ \sum_{j'=1}^{s} w_{N,j,j'} \sum_{j''=1}^{s} w_{1,j',j''}\xi_{N',j''} = \xi_{N,j}^{N+2}\frac{N!}{(2+N)!} \tag{10}$$

Order of consistency 4 :

$$\sum_{j'=1}^{s} w_{N,j,j'}\tau_{j'}^{3} = \xi_{N,j}^{N+3}6\frac{N!}{(3+N)!} \tag{1}$$

$$\sum_{j'=1}^{s} w_{N,j,j'} \sum_{j''=1}^{s} w_{1,j',j''}\tau_{j''}^{2} = \xi_{N,j}^{N+3}2\frac{N!}{(3+N)!} \tag{2}$$

$$\sum_{j'=1}^{s} w_{N,j,j'}\tau_{j'} \sum_{j''=1}^{s} w_{1,j',j''}\tau_{j''} = \xi_{N,j}^{N+3}3\frac{N!}{(3+N)!} \tag{3}$$

$$\sum_{j'=1}^{s} w_{N,j,j'}\tau_{j'}^{2} \sum_{j''=1}^{s} w_{1,j',j''} = \xi_{N,j}^{N+3}6\frac{N!}{(3+N)!} \tag{4}$$

$$\sum_{j'=1}^{s} w_{N,j,j'} \sum_{j''=1}^{s} w_{1,j',j''} \sum_{j'''=1}^{s} w_{1,j'',j'''}\tau_{j'''} = \xi_{N,j}^{N+3}\frac{N!}{(3+N)!} \tag{5}$$

$$\sum_{j'=1}^{s} w_{N,j,j'} \sum_{j''=1}^{s} w_{1,j',j''}\tau_{j''} \sum_{j'''=1}^{s} w_{1,j'',j'''} = \xi_{N,j}^{N+3}2\frac{N!}{(3+N)!} \tag{6}$$

$$\sum_{j'=1}^{s} w_{N,j,j'}\tau_{j'} \sum_{j''=1}^{s} w_{1,j',j''} \sum_{j'''=1}^{s} w_{1,j'',j'''} = \xi_{N,j}^{N+3}3\frac{N!}{(3+N)!} \tag{7}$$

$$\sum_{j'=1}^{s} w_{N,j,j'} \left(\sum_{j''=1}^{s} w_{1,j',j''}\right) \sum_{j''=1}^{s} w_{1,j',j''}\tau_{j''} = \xi_{N,j}^{N+3}3\frac{N!}{(3+N)!} \tag{8}$$

$$\sum_{j'=1}^{s} w_{N,j,j'}\tau_{j'} \left(\sum_{j''=1}^{s} w_{1,j',j''}\right)^{2} = \xi_{N,j}^{N+3}6\frac{N!}{(3+N)!} \tag{9}$$

$$\sum_{j'=1}^{s} w_{N,j,j'} \sum_{j''=1}^{s} w_{1,j',j''} \sum_{j'''=1}^{s} w_{1,j'',j'''} \sum_{j''''=1}^{s} w_{1,j''',j''''} = \xi_{N,j}^{N+3}\frac{N!}{(3+N)!} \tag{10}$$

$$\sum_{j'=1}^{s} w_{N,j,j'} \sum_{j''=1}^{s} w_{1,j',j''} \left(\sum_{j'''=1}^{s} w_{1,j'',j'''}\right)^{2} = \xi_{N,j}^{N+3}2\frac{N!}{(3+N)!} \tag{11}$$

$$\sum_{j'=1}^{s} w_{N,j,j'} \left(\sum_{j''=1}^{s} w_{1,j',j''}\right)^{3} = \xi_{N,j}^{N+3}6\frac{N!}{(3+N)!} \tag{12}$$

$$\sum_{j'=1}^{s} w_{N,j,j'} \left(\sum_{j''=1}^{s} w_{1,j',j''}\right) \sum_{j''=1}^{s} w_{1,j',j''} \sum_{j'''=1}^{s} w_{1,j'',j'''} = \xi_{N,j}^{N+3}3\frac{N!}{(3+N)!} \tag{13}$$

$$\text{if } n > 1,\ \ \sum_{j'=1}^{s} w_{N,j,j'} \sum_{j''=1}^{s} w_{2,j',j''}\tau_{j''} = \xi_{N,j}^{N+3}2\frac{N!}{(3+N)!} \tag{14}$$

$$\text{if } n > 1,\ \ \sum_{j'=1}^{s} w_{N,j,j'} \left(\sum_{j''=1}^{s} w_{1,j',j''}\right)\left(\sum_{j''=1}^{s} w_{2,j',j''}\right) = \xi_{N,j}^{N+3}6\frac{N!}{(3+N)!} \tag{15}$$

$$\text{if } n > 1,\ \sum_{j'=1}^{s} w_{N,j,j'} \sum_{j''=1}^{s} w_{2,j',j''} \sum_{j'''=1}^{s} w_{1,j'',j'''} = \xi_{N,j}^{N+3} 2 \frac{N!}{(3+N)!} \tag{16}$$

$$\text{if } n > 1,\ \sum_{j'=1}^{s} w_{N,j,j'} \sum_{j''=1}^{s} w_{1,j',j''} \sum_{j'''=1}^{s} w_{2,j'',j'''} = \xi_{N,j}^{N+3} 2 \frac{N!}{(3+N)!} \tag{17}$$

$$\text{if } n > 1,\ \sum_{j'=1}^{s} w_{N,j,j'} \tau_{j'} \sum_{j''=1}^{s} w_{2,j',j''} = \xi_{N,j}^{N+3} 6 \frac{N!}{(3+N)!} \tag{18}$$

$$\text{if } n > 2,\ \sum_{j'=1}^{s} w_{N,j,j'} \sum_{j''=1}^{s} w_{3,j',j''} = \xi_{N,j}^{N+3} 6 \frac{N!}{(3+N)!} \tag{19}$$

$$\forall N' \in [\![2,n]\!],\ \sum_{j'=1}^{s} w_{N,j,j'} \tau_{j'}^{2} \xi_{N',j'} = \xi_{N,j}^{N+3} 6 \frac{N!}{(3+N)!} \tag{20}$$

$$\forall (N',N'') \in [\![2,n]\!]^{2},\ \sum_{j'=1}^{s} w_{N,j,j'} \tau_{j'} \xi_{N',j'} \xi_{N'',j'} = \xi_{N,j}^{N+3} 6 \frac{N!}{(3+N)!} \tag{21}$$

$$\forall N' \in [\![2,n]\!],\ \sum_{j'=1}^{s} w_{N,j,j'} \tau_{j'} \xi_{N',j'} \sum_{j''=1}^{s} w_{1,j',j''} = \xi_{N,j}^{N+3} 6 \frac{N!}{(3+N)!} \tag{22}$$

$$\forall (N',N'',N''') \in [\![2,n]\!]^{3},\ \sum_{j'=1}^{s} w_{N,j,j'} \xi_{N',j'} \xi_{N'',j'} \xi_{N''',j'} = \xi_{N,j}^{N+3} 6 \frac{N!}{(3+N)!} \tag{23}$$

$$\forall (N',N'') \in [\![2,n]\!]^{2},\ \sum_{j'=1}^{s} w_{N,j,j'} \xi_{N',j'} \xi_{N'',j'} \sum_{j''=1}^{s} w_{1,j',j''} = \xi_{N,j}^{N+3} 6 \frac{N!}{(3+N)!} \tag{24}$$

$$\forall N' \in [\![2,n]\!],\ \sum_{j'=1}^{s} w_{N,j,j'} \xi_{N',j'} \left( \sum_{j''=1}^{s} w_{1,j',j''} \right)^{2} = \xi_{N,j}^{N+3} 6 \frac{N!}{(3+N)!} \tag{25}$$

$$\forall N' \in [\![2,n]\!],\ \sum_{j'=1}^{s} w_{N,j,j'} \xi_{N',j'} \sum_{j''=1}^{s} w_{2,j',j''} = \xi_{N,j}^{N+3} 6 \frac{N!}{(3+N)!} \tag{26}$$

$$\forall N' \in [\![2,n]\!],\ \sum_{j'=1}^{s} w_{N,j,j'} \tau_{j'} \sum_{j''=1}^{s} w_{1,j',j''} \xi_{N',j''} = \xi_{N,j}^{N+3} 3 \frac{N!}{(3+N)!} \tag{27}$$

$$\forall N' \in [\![2,n]\!],\ \sum_{j'=1}^{s} w_{N,j,j'} \xi_{N',j'} \sum_{j''=1}^{s} w_{1,j',j''} \tau_{j''} = \xi_{N,j}^{N+3} 3 \frac{N!}{(3+N)!} \tag{28}$$

$$\forall (N',N'') \in [\![2,n]\!]^{2},\ \sum_{j'=1}^{s} w_{N,j,j'} \xi_{N',j'} \sum_{j''=1}^{s} w_{1,j',j''} \xi_{N'',j''} = \xi_{N,j}^{N+3} 3 \frac{N!}{(3+N)!} \tag{29}$$

$$\forall N' \in [\![2,n]\!],\ \sum_{j'=1}^{s} w_{N,j,j'} \xi_{N',j'} \sum_{j''=1}^{s} w_{1,j',j''} \sum_{j'''=1}^{s} w_{1,j'',j'''} = \xi_{N,j}^{N+3} 3 \frac{N!}{(3+N)!} \tag{30}$$

$$\forall N' \in [\![2,n]\!],\ \sum_{j'=1}^{s} w_{N,j,j'} \left( \sum_{j''=1}^{s} w_{1,j',j''} \right) \sum_{j''=1}^{s} w_{1,j',j''} \xi_{N',j''} = \xi_{N,j}^{N+3} 3 \frac{N!}{(3+N)!} \tag{31}$$

$$\forall N' \in [\![2,n]\!],\ \sum_{j'=1}^{s} w_{N,j,j'} \sum_{j''=1}^{s} w_{2,j',j''} \xi_{N',j''} = \xi_{N,j}^{N+3} 2 \frac{N!}{(3+N)!} \tag{32}$$

$$\forall N' \in [\![2,n]\!],\ \sum_{j'=1}^{s} w_{N,j,j'} \sum_{j''=1}^{s} w_{1,j',j''}\xi_{N',j''} \sum_{j'''=1}^{s} w_{1,j'',j'''} = 2\xi_{N,j}^{N+3}\frac{N!}{(3+N)!} \tag{33}$$

$$\forall N' \in [\![2,n]\!],\ \sum_{j'=1}^{s} w_{N,j,j'} \sum_{j''=1}^{s} w_{1,j',j''}\tau_{j''}\xi_{N',j''} = \xi_{N,j}^{N+3}2\frac{N!}{(3+N)!} \tag{34}$$

$$\forall (N',N'') \in [\![2,n]\!]^2,\ \sum_{j'=1}^{s} w_{N,j,j'} \sum_{j''=1}^{s} w_{1,j',j''}\xi_{N',j''}\xi_{N'',j''} = \xi_{N,j}^{N+3}2\frac{N!}{(3+N)!} \tag{35}$$

$$\forall N' \in [\![2,n]\!],\ \sum_{j'=1}^{s} w_{N,j,j'} \sum_{j''=1}^{s} w_{1,j',j''} \sum_{j'''=1}^{s} w_{1,j'',j'''}\xi_{N',j'''} = \xi_{N,j}^{N+3}\frac{N!}{(3+N)!} \tag{36}$$

For a method to be at least of order of consistency 1 it must satisfy the equations under the "order of consistency 1" title, to be at least of order of consistency 2 the method must satisfy both "order of consistency 1" and "consistency 2", etc...

To find the order of conditions up to order 4 we use the equation for order of consistency $v'$ of class 4 for length $n$ at stage $j$ and rank $N$. We first find the differentials of $f$ and the derivatives of $f \circ \overline{\mathcal{J}}^{n-1}\hat{\mathrm{y}}$ in respect to $t$ up to the third order. To do so we need the Taylor expansion of order 3 of each $y_{n,j,N}$. We then replace the expressions in the equation for order of consistency, and identify each terms. The terms with a factor $h^0$ yield the conditions for an order of consistency 1, the terms with a factor of $h^0$ or $h^1$ yield the conditions for order of consistency 2, etc...

We first find the derivatives of $f \circ \overline{\mathcal{J}}^{n-1}\hat{\mathrm{y}}$. We have :

$$\frac{\mathrm{d}f \circ \overline{\mathcal{J}}^{n-1}\hat{\mathrm{y}}}{\mathrm{d}t} = \frac{\partial f}{\partial t} + \sum_{N=1}^{n}\sum_{k=1}^{l}\frac{\partial f}{\partial x_{N,k}}\frac{\mathrm{d}^{n-N+1}\hat{\mathrm{y}}_k}{\mathrm{d}t^{n-N+1}}$$

We now differentiate twice this expression in respect to $t$. Some derivatives of $\hat{\mathrm{y}}$ are high enough to reach $f$ and its derivatives, we thus need to develop those expressions in order for the expressions to only depend on derivatives of $f$ and derivatives of $\hat{\mathrm{y}}$ of order strictly lower than $n$. For example, the developed expression of $\frac{\mathrm{d}f \circ \overline{\mathcal{J}}^{n-1}\hat{\mathrm{y}}}{\mathrm{d}t}$ is :

$$\frac{\mathrm{d}f \circ \overline{\mathcal{J}}^{n-1}\hat{\mathrm{y}}}{\mathrm{d}t} = \frac{\partial f}{\partial t} + \sum_{N=2}^{n}\sum_{k=1}^{l}\frac{\partial f}{\partial x_{N,k}}\frac{\mathrm{d}^{n-N+1}\hat{\mathrm{y}}_k}{\mathrm{d}t^{n-N+1}} + \sum_{k=1}^{l}\frac{\partial f}{\partial x_{1,k}}f_k$$

Differentiating the non-developed expression of $\frac{\mathrm{d}f \circ \overline{\mathcal{J}}^{n-1}\hat{\mathrm{y}}}{\mathrm{d}t}$ gives the non-developed expression of $\frac{\mathrm{d}^2 f \circ \overline{\mathcal{J}}^{n-1}\hat{\mathrm{y}}}{\mathrm{d}t^2}$.

$$\begin{aligned}
\frac{\mathrm{d}^2 f \circ \overline{\mathcal{J}}^{n-1}\hat{\mathrm{y}}}{\mathrm{d}t^2} &= \frac{\partial^2 f}{\partial t^2} + \sum_{N=1}^{n}\sum_{k=1}^{l}\frac{\partial^2 f}{\partial t \partial x_{N,k}}\frac{\mathrm{d}^{n-N+1}\hat{\mathrm{y}}_k}{\mathrm{d}t^{n-N+1}} + \sum_{N=1}^{n}\sum_{k=1}^{l}\frac{\partial f}{\partial x_{N,k}}\frac{\mathrm{d}^{n-N+2}\hat{\mathrm{y}}_k}{\mathrm{d}t^{n-N+2}} \\
&+ \sum_{N=1}^{n}\sum_{k=1}^{l}\left(\frac{\partial^2 f}{\partial t \partial x_{N,k}} + \sum_{N'=1}^{n}\sum_{k=1'}^{l}\frac{\partial^2 f}{\partial x_{N,k}\partial x_{N',k'}}\frac{\mathrm{d}^{n-N'+1}\hat{\mathrm{y}}_{k'}}{\mathrm{d}t^{n-N'+1}}\right)\frac{\mathrm{d}^{n-N+1}\hat{\mathrm{y}}_k}{\mathrm{d}t^{n-N+1}} \\
&= \frac{\partial^2 f}{\partial t^2} + 2\sum_{N=1}^{n}\sum_{k=1}^{l}\frac{\partial^2 f}{\partial t \partial x_{N,k}}\frac{\mathrm{d}^{n-N+1}\hat{\mathrm{y}}_k}{\mathrm{d}t^{n-N+1}} + \sum_{N=1}^{n}\sum_{k=1}^{l}\frac{\partial f}{\partial x_{N,k}}\frac{\mathrm{d}^{n-N+2}\hat{\mathrm{y}}_k}{\mathrm{d}t^{n-N+2}} \\
&+ \sum_{N=1}^{n}\sum_{k=1}^{l}\sum_{N'=1}^{n}\sum_{k'=1}^{l}\frac{\partial^2 f}{\partial x_{N,k}\partial x_{N',k'}}\frac{\mathrm{d}^{n-N'+1}\hat{\mathrm{y}}_{k'}}{\mathrm{d}t^{n-N'+1}}\frac{\mathrm{d}^{n-N+1}\hat{\mathrm{y}}_k}{\mathrm{d}t^{n-N+1}}
\end{aligned}$$

The developed expression of $\frac{\mathrm{d}^2 f \circ \overline{\mathcal{J}}^{n-1}\hat{\mathrm{y}}}{\mathrm{d}t^2}$ is :

$$
\begin{aligned}
\frac{\mathrm{d}^2 f\circ\overline{\mathscr{J}}^{n-1}\hat{\mathrm{y}}}{\mathrm{d}t^2} &= \frac{\partial^2 f}{\partial t^2} + 2\sum_{k=1}^{l}\frac{\partial^2 f}{\partial t\partial x_{1,k}}f_k + \sum_{k=1}^{l}\frac{\partial f}{\partial x_{2,k}}f_k + \sum_{k=1}^{l}\frac{\partial f}{\partial x_{1,k}}\frac{\partial f_k}{\partial t}\\
&+\sum_{k=1}^{l}\frac{\partial f}{\partial x_{1,k}}\sum_{k'=1}^{l}\frac{\partial f_k}{\partial x_{1,k'}}f_{k'} + \sum_{k=1}^{l}\sum_{k'=1}^{l}\frac{\partial^2 f}{\partial x_{1,k}\partial x_{1,k'}}f_{k'}f_k\\
&+2\sum_{N=2}^{n}\sum_{k=1}^{l}\frac{\partial^2 f}{\partial t\partial x_{N,k}}\frac{\mathrm{d}^{n-N+1}\hat{\mathrm{y}}_k}{\mathrm{d}t^{n-N+1}} + \sum_{N=3}^{n}\sum_{k=1}^{l}\frac{\partial f}{\partial x_{N,k}}\frac{\mathrm{d}^{n-N+2}\hat{\mathrm{y}}_k}{\mathrm{d}t^{n-N+2}} + \sum_{k=1}^{l}\frac{\partial f}{\partial x_{1,k}}\sum_{N=2}^{n}\sum_{k'=1}^{l}\frac{\partial f_k}{\partial x_{N,k'}}\frac{\mathrm{d}^{n-N+1}\hat{\mathrm{y}}_{k'}}{\mathrm{d}t^{n-N+1}}\\
&+\sum_{N=2}^{n}\sum_{k=1}^{l}\sum_{N'=2}^{n}\sum_{k'=1}^{l}\frac{\partial^2 f}{\partial x_{N,k}\partial x_{N',k'}}\frac{\mathrm{d}^{n-N'+1}\hat{\mathrm{y}}_{k'}}{\mathrm{d}t^{n-N'+1}}\frac{\mathrm{d}^{n-N+1}\hat{\mathrm{y}}_k}{\mathrm{d}t^{n-N+1}} + 2\sum_{N=2}^{n}\sum_{k=1}^{l}\sum_{k'=1}^{l}\frac{\partial^2 f}{\partial x_{N,k}\partial x_{1,k'}}f_{k'}\frac{\mathrm{d}^{n-N+1}\hat{\mathrm{y}}_k}{\mathrm{d}t^{n-N+1}}
\end{aligned}
$$

The non-developed expression of $\frac{\mathrm{d}^2 f\circ\overline{\mathscr{J}}^{n-1}\hat{\mathrm{y}}}{\mathrm{d}t^2}$ gives the non-developed expression of $\frac{\mathrm{d}^3 f\circ\overline{\mathscr{J}}^{n-1}\hat{\mathrm{y}}}{\mathrm{d}t^3}$.

$$
\begin{aligned}
\frac{\mathrm{d}^3 f\circ\overline{\mathscr{J}}^{n-1}\hat{\mathrm{y}}}{\mathrm{d}t^3} &= \frac{\partial^3 f}{\partial t^3} + \sum_{N=1}^{n}\sum_{k=1}^{l}\frac{\partial^3 f}{\partial t^2\partial x_{N,k}}\frac{\mathrm{d}^{n-N+1}\hat{\mathrm{y}}_k}{\mathrm{d}t^{n-N+1}}\\
&+2\sum_{N=1}^{n}\sum_{k=1}^{l}\frac{\partial^2 f}{\partial t\partial x_{N,k}}\frac{\mathrm{d}^{n-N+2}\hat{\mathrm{y}}_k}{\mathrm{d}t^{n-N+2}} + 2\sum_{N=1}^{n}\sum_{k=1}^{l}\left(\frac{\partial^3 f}{\partial t^2\partial x_{N,k}} + \sum_{N'=1}^{n}\sum_{k'=1}^{l}\frac{\partial^3 f}{\partial t\partial x_{N,k}\partial x_{N',k'}}\frac{\mathrm{d}^{n-N'+1}\hat{\mathrm{y}}_{k'}}{\mathrm{d}t^{n-N'+1}}\right)\frac{\mathrm{d}^{n-N+1}\hat{\mathrm{y}}_k}{\mathrm{d}t^{n-N+1}}\\
&+\sum_{N=1}^{n}\sum_{k=1}^{l}\frac{\partial f}{\partial x_{N,k}}\frac{\mathrm{d}^{n-N+3}\hat{\mathrm{y}}_k}{\mathrm{d}t^{n-N+3}} + \sum_{N=1}^{n}\sum_{k=1}^{l}\left(\frac{\partial^2 f}{\partial t\partial x_{N,k}} + \sum_{N'=1}^{n}\sum_{k'=1}^{l}\frac{\partial f}{\partial x_{N,k}\partial x_{N',k'}}\frac{\mathrm{d}^{n-N'+1}\hat{\mathrm{y}}_{k'}}{\mathrm{d}t^{n-N'+1}}\right)\frac{\mathrm{d}^{n-N+2}\hat{\mathrm{y}}_k}{\mathrm{d}t^{n-N+2}}\\
&+\sum_{N=1}^{n}\sum_{k=1}^{l}\sum_{N'=1}^{n}\sum_{k'=1}^{l}\left(\frac{\partial^3 f}{\partial t\partial x_{N,k}\partial x_{N',k'}} + \sum_{N''=1}^{n}\sum_{k''=1}^{l}\frac{\partial^3 f}{\partial x_{N,k}\partial x_{N',k'}\partial x_{N'',k''}}\frac{\mathrm{d}^{n-N''+1}\hat{\mathrm{y}}_{k''}}{\mathrm{d}t^{n-N''+1}}\right)\frac{\mathrm{d}^{n-N'+1}\hat{\mathrm{y}}_{k'}}{\mathrm{d}t^{n-N'+1}}\frac{\mathrm{d}^{n-N+1}\hat{\mathrm{y}}_k}{\mathrm{d}t^{n-N+1}}\\
&+2\sum_{N=1}^{n}\sum_{k=1}^{l}\sum_{N'=1}^{n}\sum_{k'=1}^{l}\frac{\partial^2 f}{\partial x_{N,k}\partial x_{N',k'}}\frac{\mathrm{d}^{n-N'+1}\hat{\mathrm{y}}_{k'}}{\mathrm{d}t^{n-N'+1}}\frac{\mathrm{d}^{n-N+2}\hat{\mathrm{y}}_k}{\mathrm{d}t^{n-N+2}}\\
&= \frac{\partial^3 f}{\partial t^3} + 3\sum_{N=1}^{n}\sum_{k=1}^{l}\frac{\partial^3 f}{\partial t^2\partial x_{N,k}}\frac{\mathrm{d}^{n-N+1}\hat{\mathrm{y}}_k}{\mathrm{d}t^{n-N+1}} + 3\sum_{N=1}^{n}\sum_{k=1}^{l}\frac{\partial^2 f}{\partial t\partial x_{N,k}}\frac{\mathrm{d}^{n-N+2}\hat{\mathrm{y}}_k}{\mathrm{d}t^{n-N+2}} + \sum_{N=1}^{n}\sum_{k=1}^{l}\frac{\partial f}{\partial x_{N,k}}\frac{\mathrm{d}^{n-N+3}\hat{\mathrm{y}}_k}{\mathrm{d}t^{n-N+3}}\\
&+3\sum_{N=1}^{n}\sum_{k=1}^{l}\sum_{N'=1}^{n}\sum_{k'=1}^{l}\frac{\partial^3 f}{\partial t\partial x_{N,k}\partial x_{N',k'}}\frac{\mathrm{d}^{n-N'+1}\hat{\mathrm{y}}_{k'}}{\mathrm{d}t^{n-N'+1}}\frac{\mathrm{d}^{n-N+1}\hat{\mathrm{y}}_k}{\mathrm{d}t^{n-N+1}} + 3\sum_{N=1}^{n}\sum_{k=1}^{l}\sum_{N'=1}^{n}\sum_{k'=1}^{l}\frac{\partial f}{\partial x_{N,k}\partial x_{N',k'}}\frac{\mathrm{d}^{n-N'+1}\hat{\mathrm{y}}_{k'}}{\mathrm{d}t^{n-N'+1}}\frac{\mathrm{d}^{n-N+2}\hat{\mathrm{y}}_k}{\mathrm{d}t^{n-N+2}}\\
&+\sum_{N=1}^{n}\sum_{k=1}^{l}\sum_{N'=1}^{n}\sum_{k'=1}^{l}\sum_{N''=1}^{n}\sum_{k''=1}^{l}\frac{\partial^3 f}{\partial x_{N,k}\partial x_{N',k'}\partial x_{N'',k''}}\frac{\mathrm{d}^{n-N''+1}\hat{\mathrm{y}}_{k''}}{\mathrm{d}t^{n-N''+1}}\frac{\mathrm{d}^{n-N'+1}\hat{\mathrm{y}}_{k'}}{\mathrm{d}t^{n-N'+1}}\frac{\mathrm{d}^{n-N+1}\hat{\mathrm{y}}_k}{\mathrm{d}t^{n-N+1}}
\end{aligned}
$$

$$
\begin{aligned}
\frac{\mathrm{d}^3 f \circ \overline{\mathcal{J}}^{n-1}\hat{\mathrm{y}}}{\mathrm{d}t^3} = {} & \frac{\partial^3 f}{\partial t^3} + 3\sum_{k=1}^{l} \frac{\partial^3 f}{\partial t^2 \partial x_{1,k}} f_k + 3\sum_{k=1}^{l} \frac{\partial^2 f}{\partial t \partial x_{2,k}} f_k + 3\sum_{k=1}^{l} \frac{\partial^2 f}{\partial t \partial x_{1,k}} \frac{\partial f_k}{\partial t} + \sum_{k=1}^{l} \frac{\partial f}{\partial x_{3,k}} f_k + \sum_{k=1}^{l} \frac{\partial f}{\partial x_{2,k}} \frac{\partial f_k}{\partial t} \\
& + \sum_{k=1}^{l} \frac{\partial f}{\partial x_{1,k}} \frac{\partial^2 f_k}{\partial t^2} + 3\sum_{k=1}^{l} \frac{\partial^2 f}{\partial t \partial x_{1,k}} \sum_{k'=1}^{l} \frac{\partial f_k}{\partial x_{1,k'}} f_{k'} + \sum_{k=1}^{l} \frac{\partial f}{\partial x_{2,k}} \sum_{k'=1}^{l} \frac{\partial f_k}{\partial x_{1,k'}} f_{k'} + 2\sum_{k=1}^{l} \frac{\partial f}{\partial x_{1,k}} \sum_{k'=1}^{l} \frac{\partial^2 f_k}{\partial t \partial x_{1,k'}} f_{k'} \\
& + \sum_{k=1}^{l} \frac{\partial f}{\partial x_{1,k}} \sum_{k'=1}^{l} \frac{\partial f_k}{\partial x_{2,k'}} f_{k'} + \sum_{k=1}^{l} \frac{\partial f}{\partial x_{1,k}} \sum_{k'=1}^{l} \frac{\partial f_k}{\partial x_{1,k'}} \frac{\partial f_{k'}}{\partial t} + 3\sum_{k=1}^{l} \sum_{k'=1}^{l} \frac{\partial^3 f}{\partial t \partial x_{1,k} \partial x_{1,k'}} f_{k'} f_k \\
& +3\sum_{k=1}^{l} \sum_{k'=1}^{l} \frac{\partial f}{\partial x_{2,k} \partial x_{1,k'}} f_{k'} f_k + 3\sum_{k=1}^{l} \sum_{k'=1}^{l} \frac{\partial f}{\partial x_{1,k} \partial x_{1,k'}} f_{k'} \frac{\partial f_k}{\partial t} + \sum_{k=1}^{l} \frac{\partial f}{\partial x_{1,k}} \sum_{k'=1}^{l} \frac{\partial f_k}{\partial x_{1,k'}} \sum_{k''=1}^{l} \frac{\partial f_{k'}}{\partial x_{1,k''}} f_{k''} \\
& + \sum_{k=1}^{l} \frac{\partial f}{\partial x_{1,k}} \sum_{k'=1}^{l} \sum_{k''=1}^{l} \frac{\partial^2 f^{[l]}}{\partial x_{1,k'} \partial x_{1,k''}} f_{k''} f_{k'} + 3\sum_{k=1}^{l} \sum_{k'=1}^{l} \frac{\partial f}{\partial x_{1,k} \partial x_{1,k'}} f_{k'} \sum_{k''=1}^{l} \frac{\partial f_k}{\partial x_{1,k''}} f_{k''} \\
& + \sum_{k=1}^{l} \sum_{k'=1}^{l} \sum_{k''=1}^{l} \frac{\partial^3 f}{\partial x_{1,k} \partial x_{1,k'} x_{1,k''}} f_{k''} f_{k'} f_k + 3\sum_{N=2}^{n} \sum_{k=1}^{l} \frac{\partial^3 f}{\partial t^2 \partial x_{N,k}} \frac{\mathrm{d}^{n-N+1}\hat{\mathrm{y}}_k}{\mathrm{d}t^{n-N+1}} + 3\sum_{N=3}^{n} \sum_{k=1}^{l} \frac{\partial^2 f}{\partial t \partial x_{N,k}} \frac{\mathrm{d}^{n-N+2}\hat{\mathrm{y}}_k}{\mathrm{d}t^{n-N+2}} \\
& + \sum_{N=4}^{n} \sum_{k=1}^{l} \frac{\partial f}{\partial x_{N,k}} \frac{\mathrm{d}^{n-N+3}\hat{\mathrm{y}}_k}{\mathrm{d}t^{n-N+3}} + 3\sum_{k=1}^{l} \frac{\partial^2 f}{\partial t \partial x_{1,k}} \sum_{N=2}^{n} \sum_{k'=1}^{l} \frac{\partial f_k}{\partial x_{N,k'}} \frac{\mathrm{d}^{n-N+1}\hat{\mathrm{y}}_{k'}}{\mathrm{d}t^{n-N+1}} + \sum_{k=1}^{l} \frac{\partial f}{\partial x_{2,k}} \sum_{N=2}^{n} \sum_{k'=1}^{l} \frac{\partial f_k}{\partial x_{N,k'}} \frac{\mathrm{d}^{n-N+1}\hat{\mathrm{y}}_{k'}}{\mathrm{d}t^{n-N+1}} \\
& +6\sum_{N=2}^{n} \sum_{k=1}^{l} \sum_{k'=1}^{l} \frac{\partial^3 f}{\partial t \partial x_{N,k} \partial x_{1,k'}} f_{k'} \frac{\mathrm{d}^{n-N+1}\hat{\mathrm{y}}_k}{\mathrm{d}t^{n-N+1}} + 2\sum_{k=1}^{l} \frac{\partial f}{\partial x_{1,k}} \sum_{N=2}^{n} \sum_{k'=1}^{l} \frac{\partial^2 f_k}{\partial t \partial x_{N,k'}} \frac{\mathrm{d}^{n-N+1}\hat{\mathrm{y}}_{k'}}{\mathrm{d}t^{n-N+1}} \\
& + \sum_{k=1}^{l} \frac{\partial f}{\partial x_{1,k}} \sum_{N=3}^{n} \sum_{k'=1}^{l} \frac{\partial f_k}{\partial x_{N,k'}} \frac{\mathrm{d}^{n-N+2}\hat{\mathrm{y}}_{k'}}{\mathrm{d}t^{n-N+2}} + 3\sum_{N=3}^{n} \sum_{k=1}^{l} \sum_{k'=1}^{l} \frac{\partial f}{\partial x_{N,k} \partial x_{1,k'}} f_{k'} \frac{\mathrm{d}^{n-N+2}\hat{\mathrm{y}}_k}{\mathrm{d}t^{n-N+2}} \\
& +3\sum_{k=1}^{l} \sum_{N=2}^{n} \sum_{k'=1}^{l} \frac{\partial f}{\partial x_{2,k} \partial x_{N,k'}} \frac{\mathrm{d}^{n-N+1}\hat{\mathrm{y}}_{k'}}{\mathrm{d}t^{n-N+1}} f_k + 3\sum_{k=1}^{l} \sum_{N=2}^{n} \sum_{k'=1}^{l} \frac{\partial f}{\partial x_{1,k} \partial x_{N,k'}} \frac{\mathrm{d}^{n-N+1}\hat{\mathrm{y}}_{k'}}{\mathrm{d}t^{n-N+1}} \frac{\partial f_k}{\partial t} \\
& +2\sum_{k=1}^{l} \frac{\partial f}{\partial x_{1,k}} \sum_{N=2}^{n} \sum_{k'=1}^{l} \sum_{k''=1}^{l} \frac{\partial^2 f_k}{\partial x_{N,k'} \partial x_{1,k''}} f_{k''} \frac{\mathrm{d}^{n-N+1}\hat{\mathrm{y}}_{k'}}{\mathrm{d}t^{n-N+1}} + \sum_{k=1}^{l} \frac{\partial f}{\partial x_{1,k}} \sum_{k'=1}^{l} \frac{\partial f_k}{\partial x_{1,k'}} \sum_{N=2}^{n} \sum_{k''=1}^{l} \frac{\partial f_{k'}}{\partial x_{N,k''}} \frac{\mathrm{d}^{i+1}\hat{\mathrm{y}}_{k''}}{\mathrm{d}t^{i+1}} \\
& +3\sum_{k=1}^{l} \sum_{N=2}^{n} \sum_{k'=1}^{l} \frac{\partial f}{\partial x_{1,k} \partial x_{N,k'}} \frac{\mathrm{d}^{n-N+1}\hat{\mathrm{y}}_{k'}}{\mathrm{d}t^{n-N+1}} \sum_{k''=1}^{l} \frac{\partial f_k}{\partial x_{1,k''}} f_{k''} + 3\sum_{k=1}^{l} \sum_{k'=1}^{l} \frac{\partial f}{\partial x_{1,k} \partial x_{1,k'}} f_{k'} \sum_{N'=2}^{n} \sum_{k''=1}^{l} \frac{\partial f_k}{\partial x_{N',k''}} \frac{\mathrm{d}^{n-N'+1}\hat{\mathrm{y}}_{k''}}{\mathrm{d}t^{n-N'+1}} \\
& +3\sum_{N=2}^{n} \sum_{k=1}^{l} \sum_{k'=1}^{l} \sum_{k''=1}^{l} \frac{\partial^3 f}{\partial x_{N,k} \partial x_{1,k'} x_{1,k''}} f_{k''} f_{k'} \frac{\mathrm{d}^{n-N+1}\hat{\mathrm{y}}_k}{\mathrm{d}t^{n-N+1}} \\
& +3\sum_{N=2}^{n} \sum_{k=1}^{l} \sum_{N'=2}^{n} \sum_{k'=1}^{l} \frac{\partial^3 f}{\partial t \partial x_{N,k} \partial x_{N',k'}} \frac{\mathrm{d}^{n-N'+1}\hat{\mathrm{y}}_{k'}}{\mathrm{d}t^{n-N'+1}} \frac{\mathrm{d}^{n-N+1}\hat{\mathrm{y}}_k}{\mathrm{d}t^{n-N+1}} + 3\sum_{N=3}^{n} \sum_{k=1}^{l} \sum_{N'=2}^{n} \sum_{k'=1}^{l} \frac{\partial f}{\partial x_{N,k} \partial x_{N',k'}} \frac{\mathrm{d}^{n-N'+1}\hat{\mathrm{y}}_{k'}}{\mathrm{d}t^{n-N'+1}} \frac{\mathrm{d}^{n-N+2}\hat{\mathrm{y}}_k}{\mathrm{d}t^{n-N+2}} \\
& + \sum_{k=1}^{l} \frac{\partial f}{\partial x_{1,k}} \sum_{N=2}^{n} \sum_{k'=1}^{l} \sum_{N'=2}^{n} \sum_{k''=1}^{l} \frac{\partial^2 f_k}{\partial x_{N,k'} \partial x_{N',k''}} \frac{\mathrm{d}^{n-N'+1}\hat{\mathrm{y}}_{k''}}{\mathrm{d}t^{n-N'+1}} \frac{\mathrm{d}^{n-N+1}\hat{\mathrm{y}}_{k'}}{\mathrm{d}t^{n-N+1}} \\
& +3\sum_{k=1}^{l} \sum_{N=2}^{n} \sum_{k'=1}^{l} \frac{\partial f}{\partial x_{1,k} \partial x_{N,k'}} \frac{\mathrm{d}^{n-N+1}\hat{\mathrm{y}}_{k'}}{\mathrm{d}t^{n-N+1}} \sum_{N'=2}^{n} \sum_{k''=1}^{l} \frac{\partial f_k}{\partial x_{N',k''}} \frac{\mathrm{d}^{n-N'+1}\hat{\mathrm{y}}_{k''}}{\mathrm{d}t^{n-N'+1}} \\
& +3\sum_{N=2}^{n} \sum_{k=1}^{l} \sum_{N'=2}^{n} \sum_{k'=1}^{l} \sum_{k''=1}^{l} \frac{\partial^3 f}{\partial x_{N,k} \partial x_{N',k'} x_{1,k''}} f_{k''} \frac{\mathrm{d}^{n-N'+1}\hat{\mathrm{y}}_{k'}}{\mathrm{d}t^{n-N'+1}} \frac{\mathrm{d}^{n-N+1}\hat{\mathrm{y}}_k}{\mathrm{d}t^{n-N+1}} \\
& + \sum_{N=2}^{n} \sum_{k=1}^{l} \sum_{N'=2}^{n} \sum_{k'=1}^{l} \sum_{N''=2}^{n} \sum_{k''=1}^{l} \frac{\partial^3 f}{\partial x_{N,k} \partial x_{N',k'} \partial x_{N'',k''}} \frac{\mathrm{d}^{n-N''+1}\hat{\mathrm{y}}_{k''}}{\mathrm{d}t^{n-N''+1}} \frac{\mathrm{d}^{n-N'+1}\hat{\mathrm{y}}_{k'}}{\mathrm{d}t^{n-N'+1}} \frac{\mathrm{d}^{n-N+1}\hat{\mathrm{y}}_k}{\mathrm{d}t^{n-N+1}}
\end{aligned}
$$

Finding the differentials of $f$ is quite simpler :

$$\mathrm{d}f = \frac{\partial f}{\partial t}\,\mathrm{d}t + \sum_{N=1}^{n}\sum_{k=1}^{l}\frac{\partial f}{\partial x_{N,k}}\,\mathrm{d}\frac{\mathrm{d}^{n-N}\hat{\mathrm{y}}_k}{\mathrm{d}t^{n-N}}$$

$$\mathrm{d}^2 f = \frac{\partial^2 f}{\partial t^2}\,\mathrm{d}t^2 + 2\sum_{N=1}^{n}\sum_{k=1}^{l}\frac{\partial^2 f}{\partial t\partial x_{N,k}}\,\mathrm{d}t\,\mathrm{d}\frac{\mathrm{d}^{n-N}\hat{\mathrm{y}}_k}{\mathrm{d}t^{n-N}} + \sum_{N=1}^{n}\sum_{k=1}^{l}\sum_{N'=1}^{n}\sum_{k'=1}^{l}\frac{\partial^2 f}{\partial x_{N,k}\partial x_{N',k'}}\,\mathrm{d}\frac{\mathrm{d}^{n-N}\hat{\mathrm{y}}_k}{\mathrm{d}t^{n-N}}\,\mathrm{d}\frac{\mathrm{d}^{n-N'}\hat{\mathrm{y}}_{k'}}{\mathrm{d}t^{n-N'}}$$

$$\mathrm{d}^3 f = \frac{\partial^3 f}{\partial t^3}\,\mathrm{d}t^3 + 3\sum_{N=1}^{n}\sum_{k=1}^{l}\frac{\partial^3 f}{\partial t^2\partial x_{N,k}}\,\mathrm{d}t^2\,\mathrm{d}\frac{\mathrm{d}^{n-N}\hat{\mathrm{y}}_k}{\mathrm{d}t^{n-N}} + 3\sum_{N=1}^{n}\sum_{k=1}^{l}\sum_{N'=1}^{n}\sum_{k'=1}^{l}\frac{\partial^3 f}{\partial t\partial x_{N,k}\partial x_{N',k'}}\,\mathrm{d}t\,\mathrm{d}\frac{\mathrm{d}^{n-N}\hat{\mathrm{y}}_k}{\mathrm{d}t^{n-N}}\,\mathrm{d}\frac{\mathrm{d}^{n-N'}\hat{\mathrm{y}}_{k'}}{\mathrm{d}t^{n-N'}}$$

$$+\sum_{N=1}^{n}\sum_{k=1}^{l}\sum_{N'=1}^{n}\sum_{k'=1}^{l}\sum_{N''=1}^{n}\sum_{k''=1}^{l}\frac{\partial^3 f}{\partial x_{N,k}\partial x_{N',k'}\partial x_{N',k''}}\,\mathrm{d}\frac{\mathrm{d}^{n-N}\hat{\mathrm{y}}_k}{\mathrm{d}t^{n-N}}\,\mathrm{d}\frac{\mathrm{d}^{n-N'}\hat{\mathrm{y}}_{k'}}{\mathrm{d}t^{n-N'}}\,\mathrm{d}\frac{\mathrm{d}^{i''}\hat{\mathrm{y}}_{k''}}{\mathrm{d}t^{i''}}$$

The less simple part is finding the Taylor expansion of degree 3 of $y_{n,j',N,k} - y_{n,0,N,k}$, with $k \in [\![1,d]\!]$ and $j' \in [\![1,s+1]\!]$:

$$\forall N \in [\![0,n-4]\!]$$

$$y_{n,j',N,k} - y_{n,0,N,k} = h\xi_{j',N}y_{n,0,N-1,k} + \frac{\left(h\xi_{j',N}\right)^2}{2}y_{n,0,N-2,k} + \frac{\left(h\xi_{j',N}\right)^3}{6}y_{n,0,N-3,k} + \underset{h\to 0}{\mathcal{O}}\left(h^4\right)$$

$$y_{n,j',3,k} - y_{n,0,3,k} = hy_{n,0,2,k}\xi_{j',3} + h^2 y_{n,0,1,k}\frac{\xi_{j',3}^2}{2} + h^3 f_k\frac{1}{6}\sum_{j''=1}^{s}w_{3,j',j''} + \underset{h\to 0}{\mathcal{O}}\left(h^4\right)$$

$$y_{n,j',2,k} - y_{n,0,2,k} = hy_{n,0,1,k}\xi_{j',2} + \frac{h^2}{2}\sum_{j''=1}^{s}w_{2,j',j''}\left(f_k + \mathrm{d}f_k\right) + \underset{h\to 0}{\mathcal{O}}\left(h^4\right)$$

$$= hy_{n,0,1,k}\xi_{j',2}$$

$$+\frac{h^2}{2}\sum_{j''=1}^{s}w_{2,j',j''}\left(f_k + \frac{\partial f_k}{\partial t}\tau_{j''}h + \sum_{N=2}^{n}\sum_{k'=1}^{l}\frac{\partial f_k}{\partial x_{N,k'}}y_{n,0,N-1,k'}h\xi_{j'',N} + \sum_{k'=1}^{l}\frac{\partial f_k}{\partial x_{1,k'}}h\sum_{j'''=1}^{s}w_{1,j'',j'''}f_{k'}\right)$$

$$+\underset{h\to 0}{\mathcal{O}}\left(h^4\right)$$

$$= hy_{n,0,1,k}\xi_{j',2} + h^2 f_k\frac{1}{2}\sum_{j''=1}^{s}w_{2,j',j''} + h^3\frac{\partial f_k}{\partial t}\frac{1}{2}\sum_{j''=1}^{s}w_{2,j',j''}\tau_{j''}$$

$$+h^3\sum_{N=2}^{n}\sum_{k'=1}^{l}\frac{\partial f_k}{\partial x_{N,k'}}y_{n,0,N-1,k'}\frac{1}{2}\sum_{j''=1}^{s}w_{2,j',j''}\xi_{j'',N} + h^3\sum_{k'=1}^{l}\frac{\partial f_k}{\partial x_{1,k'}}f_{k'}\frac{1}{2}\sum_{j''=1}^{s}w_{2,j',j''}\sum_{j'''=1}^{s}w_{1,j'',j'''} + \underset{h\to 0}{\mathcal{O}}\left(h^4\right)$$

$$
\begin{aligned}
&y_{n,j',1,k} - y_{n,0,1,k} = h \sum_{j''=1}^{s} w_{1,j',j''} \Bigg( f_k + \frac{\partial f_k}{\partial t} \tau_{j''} h \\
&+ \sum_{N=3}^{n} \sum_{k'=1}^{l} \frac{\partial f_k}{\partial x_{N,k'}} \left( y_{n,j'',N,k'} - y_{n,0,N,k'} \right) + \sum_{k'=1}^{l} \frac{\partial f_k}{\partial x_{2,k'}} \left( y_{n,j'',2,k'} - y_{n,0,2,k'} \right) + \sum_{k'=1}^{l} \frac{\partial f_k}{\partial x_{1,k'}} \left( y_{n,j'',1,k'} - y_{n,0,1,k'} \right) \\
&+ \frac{1}{2} \frac{\partial^2 f_k}{\partial t^2} \left( \tau_{j''} h \right)^2 + \sum_{N=2}^{n} \sum_{k'=1}^{l} \frac{\partial^2 f_k}{\partial t \partial x_{N,k'}} \tau_{j''} h \left( y_{n,j'',N,k'} - y_{n,0,N,k'} \right) + \sum_{k'=1}^{l} \frac{\partial^2 f_k}{\partial t \partial x_{1,k'}} \tau_{j''} h \left( y_{n,j'',1,k'} - y_{n,0,1,k'} \right) \\
&+ \frac{1}{2} \sum_{N=2}^{n} \sum_{k'=1}^{l} \sum_{N'=2}^{n} \sum_{k''=1}^{l} \frac{\partial^2 f_k}{\partial x_{N,k'} \partial x_{N',k''}} \left( y_{n,j'',N,k'} - y_{n,0,N,k'} \right) \left( y_{n,j'',N',k''} - y_{n,0,N',k''} \right) \\
&+ \sum_{N=2}^{n} \sum_{k'=1}^{l} \sum_{k''=1}^{l} \frac{\partial^2 f_k}{\partial x_{N,k'} \partial x_{1,k''}} \left( y_{n,j'',N,k'} - y_{n,0,N,k'} \right) \left( y_{n,j'',1,k''} - y_{n,0,1,k''} \right) \\
&+ \frac{1}{2} \sum_{k'=1}^{l} \sum_{k''=1}^{l} \frac{\partial^2 f_k}{\partial x_{1,k'} \partial x_{1,k''}} \left( y_{n,j'',1,k'} - y_{n,0,1,k'} \right) \left( y_{n,j'',1,k''} - y_{n,0,1,k''} \right) \\
&\Bigg) + \underset{h \to 0}{\mathcal{O}} \left( h^4 \right) \\
&= h \sum_{j''=1}^{s} w_{1,j',j''} \Bigg( f_k + \frac{\partial f_k}{\partial t} \tau_{j''} h \\
&+ \sum_{N=3}^{n} \sum_{k'=1}^{l} \frac{\partial f_k}{\partial x_{N,k'}} \left( h \xi_{j'',N} y_{n,0,N-1,k'} + \frac{\left( h \xi_{j'',N} \right)^2}{2} y_{n,0,N-2,k'} \right) \\
&+ \sum_{k'=1}^{l} \frac{\partial f_k}{\partial x_{2,k'}} \left( h \xi_{j'',2} y_{n,0,1,k'} + h^2 \frac{f_{k'}}{2} \sum_{j'''=1}^{s} w_{2,j'',j'''} \right) + \sum_{k'=1}^{l} \frac{\partial f_k}{\partial x_{1,k'}} h \sum_{j'''=1}^{s} w_{1,j'',j'''} \left( f_{k'} + \mathrm{d} f_{k'} \right) \\
&+ \frac{1}{2} \frac{\partial^2 f_k}{\partial t^2} \left( \tau_{j''} h \right)^2 + \sum_{N=2}^{n} \sum_{k'=1}^{l} \frac{\partial^2 f_k}{\partial t \partial x_{N,k'}} \tau_{j''} h^2 \xi_{j'',N} y_{n,0,N-1,k'} + \sum_{k'=1}^{l} \frac{\partial^2 f_k}{\partial t \partial x_{1,k'}} \tau_{j''} h^2 \sum_{j'''=1}^{s} w_{1,j'',j'''} f_{k'} \\
&+ \frac{1}{2} \sum_{N=2}^{n} \sum_{k'=1}^{l} \sum_{N'=2}^{n} \sum_{k''=1}^{l} \frac{\partial^2 f_k}{\partial x_{N,k'} \partial x_{N',k''}} h \xi_{j'',N} y_{n,0,N-1,k'} h \xi_{N',j''} y_{n,0,N'-1,k''} \\
&+ \sum_{N=2}^{n} \sum_{k'=1}^{l} \sum_{k''=1}^{l} \frac{\partial^2 f_k}{\partial x_{N,k'} \partial x_{1,k''}} h \xi_{j'',N} y_{n,0,N-1,k'} h \sum_{j'''=1}^{s} w_{1,j'',j'''} f_{k''} \\
&+ \frac{1}{2} \sum_{k'=1}^{l} \sum_{k''=1}^{l} \frac{\partial^2 f_k}{\partial x_{1,k'} \partial x_{1,k''}} h \left( \sum_{j'''=1}^{s} w_{1,j'',j'''} \right) f_{k'} h \left( \sum_{j'''=1}^{s} w_{1,j'',j'''} \right) f_{k''} \\
&\Bigg) + \underset{h \to 0}{\mathcal{O}} \left( h^4 \right)
\end{aligned}
$$

$$
\begin{aligned}
&= h \sum_{j''=1}^{s} w_{1,j',j''} \Bigg( f_k + \frac{\partial f_k}{\partial t}\tau_{j''} h \\
&+ \sum_{N=3}^{n}\sum_{k'=1}^{l} \frac{\partial f_k}{\partial x_{N,k'}} h\xi_{j'',N} y_{n,0,N-1,k'} + \sum_{N=3}^{n}\sum_{k'=1}^{l}\frac{\partial f_k}{\partial x_{N,k'}}\frac{\left(h\xi_{j'',N}\right)^2}{2} y_{n,0,N-2,k'} \\
&+ \sum_{k'=1}^{l}\frac{\partial f_k}{\partial x_{2,k'}} h\xi_{j'',2} y_{n,0,1,k'} + \sum_{k'=1}^{l}\frac{\partial f_k}{\partial x_{2,k'}} h^2 \frac{f_{k'}}{2}\sum_{j'''=1}^{s} w_{2,j'',j'''} \\
&+ \sum_{k'=1}^{l}\frac{\partial f_k}{\partial x_{1,k'}} h \sum_{j'''=1}^{s} w_{1,j'',j'''}\Bigg( f_{k'} + \frac{\partial f_{k'}}{\partial t}\tau_{j'''} h + \sum_{N=2}^{n}\sum_{k''=1}^{l}\frac{\partial f_{k'}}{\partial x_{N,k''}} h\xi_{N,j'''} y_{n,0,N-1,k''} \\
&\quad + \sum_{k''=1}^{l}\frac{\partial f_{k'}}{\partial x_{1,k''}} h \sum_{j''''=1}^{s} w_{1,j''',j''''} f_{k''}\Bigg) \\
&+\frac{1}{2} \frac{\partial^2 f_k}{\partial t^2}\left(\tau_{j''}h\right)^2 + \sum_{N=2}^{n}\sum_{k'=1}^{l}\frac{\partial^2 f_k}{\partial t\partial x_{N,k'}}\tau_{j''}h^2\xi_{j'',N} y_{n,0,N-1,k'} + \sum_{k'=1}^{l}\frac{\partial^2 f_k}{\partial t \partial x_{1,k'}}\tau_{j''}h^2\sum_{j'''=1}^{s} w_{1,j'',j'''} f_{k'} \\
&+\frac{1}{2} \sum_{N=2}^{n}\sum_{k'=1}^{l}\sum_{N'=2}^{n}\sum_{k''=1}^{l} \frac{\partial^2 f_k}{\partial x_{N,k'}\partial x_{N',k''}} h\xi_{j'',N} y_{n,0,N-1,k'} h\xi_{N',j''} y_{n,0,N'-1,k''} \\
&+ \sum_{N=2}^{n}\sum_{k'=1}^{l}\sum_{k''=1}^{l} \frac{\partial^2 f_k}{\partial x_{N,k'}\partial x_{1,k''}} h\xi_{j'',N} y_{n,0,N-1,k'} h \sum_{j'''=1}^{s} w_{1,j'',j'''} f_{k''} \\
&+\frac{1}{2} \sum_{k'=1}^{l}\sum_{k''=1}^{l} \frac{\partial^2 f_k}{\partial x_{1,k'}\partial x_{1,k''}} h\left(\sum_{j'''=1}^{s} w_{1,j'',j'''}\right) f_{k'} h \left(\sum_{j'''=1}^{s} w_{1,j'',j'''}\right) f_{k''} \\
&\Bigg) + \underset{h\to 0}{\mathcal{O}}\left(h^4\right)
\end{aligned}
$$

$$
\begin{aligned}
&= h f_k \sum_{j''=1}^{s} w_{1,j',j''} \\
&+h^2 \frac{\partial f_k}{\partial t}\sum_{j''=1}^{s} w_{1,j',j''}\tau_{j''} \\
&+h^2 \sum_{N=3}^{n}\sum_{k'=1}^{l} \frac{\partial f_k}{\partial x_{N,k'}} y_{n,0,N-1,k'} \sum_{j''=1}^{s} w_{1,j',j''}\xi_{j'',N} \\
&+h^2 \sum_{k'=1}^{l} \frac{\partial f_k}{\partial x_{2,k'}} y_{n,0,1,k'}\sum_{j''=1}^{s} w_{1,j',j''}\xi_{j'',2} \\
&+h^2 \sum_{k'=1}^{l} \frac{\partial f_k}{\partial x_{1,k'}} f_{k'}\sum_{j''=1}^{s} w_{1,j',j''}\sum_{j'''=1}^{s} w_{1,j'',j'''} \\
&+h^3 \sum_{N=3}^{n}\sum_{k'=1}^{l}\frac{\partial f_k}{\partial x_{N,k'}} y_{n,0,N-2,k'}\sum_{j''=1}^{s} w_{1,j',j''}\frac{\xi^2_{j'',N}}{2} \\
&+h^3 \sum_{k'=1}^{l} \frac{\partial f_k}{\partial x_{2,k'}} f_{k'}\frac{1}{2} \sum_{j''=1}^{s} w_{1,j',j''}\sum_{j'''=1}^{s} w_{2,j'',j'''}
\end{aligned}
$$

$$
\begin{aligned}
&+h^3 \sum_{k'=1}^{l} \frac{\partial f_k}{\partial x_{1,k'}} \frac{\partial f_{k'}}{\partial t} \sum_{j''=1}^{s} w_{1,j',j''} \sum_{j'''=1}^{s} w_{1,j'',j'''} \tau_{j'''} \\
&+h^3 \sum_{k'=1}^{l} \frac{\partial f_k}{\partial x_{1,k'}} \sum_{N=2}^{n} \sum_{k''=1}^{l} \frac{\partial f_{k'}}{\partial x_{N,k''}} y_{n,0,N-1,k''} \sum_{j''=1}^{s} w_{1,j',j''} \sum_{j'''=1}^{s} w_{1,j'',j'''} \xi_{N,j'''} \\
&+h^3 \sum_{k'=1}^{l} \frac{\partial f_k}{\partial x_{1,k'}} \sum_{k''=1}^{l} \frac{\partial f_{k'}}{\partial x_{1,k''}} f_{k''} \sum_{j''=1}^{s} w_{1,j',j''} \sum_{j'''=1}^{s} w_{1,j'',j'''} \sum_{j''''=1}^{s} w_{1,j''',j''''} \\
&+h^3 \frac{\partial^2 f_k}{\partial t^2} \sum_{j''=1}^{s} w_{1,j',j''} \frac{\tau_{j''}^2}{2} \\
&+h^3 \sum_{N=2}^{n} \sum_{k'=1}^{l} \frac{\partial^2 f_k}{\partial t \partial x_{N,k'}} y_{n,0,N-1,k'} \sum_{j''=1}^{s} w_{1,j',j''} \tau_{j''} \xi_{j'',N} \\
&+h^3 \sum_{k'=1}^{l} \frac{\partial^2 f_k}{\partial t \partial x_{1,k'}} \sum_{j''=1}^{s} w_{1,j',j''} \tau_{j''} \sum_{j'''=1}^{s} w_{1,j'',j'''} f_{k'} \\
&+h^3 \sum_{N=2}^{n} \sum_{k'=1}^{l} \sum_{N'=2}^{n} \sum_{k''=1}^{l} \frac{\partial^2 f_k}{\partial x_{N,k'} \partial x_{N',k''}} y_{n,0,N-1,k'} y_{n,0,N'-1,k''} \frac{1}{2} \sum_{j''=1}^{s} w_{1,j',j''} \xi_{j'',N} \xi_{N',j''} \\
&+h^3 \sum_{N=2}^{n} \sum_{k'=1}^{l} \sum_{k''=1}^{l} \frac{\partial^2 f_k}{\partial x_{N,k'} \partial x_{1,k''}} y_{n,0,N-1,k'} f_{k''} \sum_{j''=1}^{s} w_{1,j',j''} \xi_{j'',N} \sum_{j'''=1}^{s} w_{1,j'',j'''} \\
&+h^3 \sum_{k'=1}^{l} \sum_{k''=1}^{l} \frac{\partial^2 f_k}{\partial x_{1,k'} \partial x_{1,k''}} f_{k'} f_{k''} \frac{1}{2} \sum_{j''=1}^{s} w_{1,j',j''} \left( \sum_{j'''=1}^{s} w_{1,j'',j'''} \right)^2 \\
&+ \underset{h \to 0}{\mathcal{O}}\left(h^4\right)
\end{aligned}
$$

We can now plug all those values in the equation and identify the terms :

$$
\begin{aligned}
&\sum_{j'=1}^{s} w_{N,j,j'} \left( f + \mathrm{d}f + \frac{1}{2}\mathrm{d}^2 f + \frac{1}{6}\mathrm{d}^3 f \right) \\
&-\xi_{N,j}^{N} f - \xi_{N,j}^{N+1} h \frac{1}{1+N} \frac{\mathrm{d} f \circ \overline{\mathcal{J}}^{n-1}\hat{\mathrm{y}}}{\mathrm{d}t} - \xi_{N,j}^{N+2} h^2 \frac{N!}{(2+N)!} \frac{\mathrm{d}^2 f \circ \overline{\mathcal{J}}^{n-1}\hat{\mathrm{y}}}{\mathrm{d}t^2} - \xi_{N,j}^{N+3} h^3 \frac{N!}{(3+N)!} \frac{\mathrm{d}^3 f \circ \overline{\mathcal{J}}^{n-1}\hat{\mathrm{y}}}{\mathrm{d}t^3} = \underset{h \to 0}{\mathcal{O}}\left(h^4\right) \\
&\Leftrightarrow \sum_{j'=1}^{s} w_{N,j,j'} \Bigg( f + \frac{\partial f}{\partial t} \tau_{j'} h + \sum_{N=1}^{n} \sum_{k=1}^{l} \frac{\partial f}{\partial x_{N,k}} \left( y_{n,j',N,k} - y_{n,0,N,k} \right) + \frac{1}{2} \frac{\partial^2 f}{\partial t^2} \left( \tau_{j'} h \right)^2 \\
&+ \sum_{N=1}^{n} \sum_{k=1}^{l} \frac{\partial^2 f}{\partial t \partial x_{N,k}} \tau_{j'} h \left( y_{n,j',N,k} - y_{n,0,N,k} \right) \\
&+ \frac{1}{2} \sum_{N=1}^{n} \sum_{k=1}^{l} \sum_{N'=1}^{n} \sum_{k'=1}^{l} \frac{\partial^2 f}{\partial x_{N,k} \partial x_{N',k'}} \left( y_{n,j',N,k} - y_{n,0,N,k} \right) \left( y_{i',j',k'} - y_{n,0,N',k'} \right) + \frac{1}{6} \frac{\partial^3 f}{\partial t^3} \left( \tau_{j'} h \right)^3 \\
&+ \frac{1}{2} \sum_{N=1}^{n} \sum_{k=1}^{l} \frac{\partial^3 f}{\partial t^2 \partial x_{N,k}} \left( \tau_{j'} h \right)^2 \left( y_{n,j',N,k} - y_{n,0,N,k} \right) \\
&+ \frac{1}{2} \sum_{N=1}^{n} \sum_{k=1}^{l} \sum_{N'=1}^{n} \sum_{k'=1}^{l} \frac{\partial^3 f}{\partial t \partial x_{N,k} \partial x_{N',k'}} \tau_{j'} h \left( y_{n,j',N,k} - y_{n,0,N,k} \right) \left( y_{i',j',k'} - y_{n,0,N',k'} \right)
\end{aligned}
$$

$$+\frac{1}{6}\sum_{N=1}^{n}\sum_{k=1}^{l}\sum_{N'=1}^{n}\sum_{k'=1}^{l}\sum_{N''=1}^{n}\sum_{k''=1}^{l}\Bigg[$$

$$\frac{\partial^3 f}{\partial x_{N,k}\partial x_{N',k'}\partial x_{N',k''}}\big(y_{n,j',N,k}-y_{n,0,N,k}\big)\big(y_{i',j',k'}-y_{n,0,N',k'}\big)\big(y_{n,j',N'',k''}-y_{n,0,N'',k''}\big)\Bigg]\Bigg)$$

$$-\xi_{N,j}^{N}f-\xi_{N,j}^{N+1}h\frac{1}{1+N}\frac{\mathrm{d}f\circ\overline{\mathscr{J}}^{n-1}\hat{\mathrm{y}}}{\mathrm{d}t}-\xi_{N,j}^{N+2}h^2\frac{N!}{(2+N)!}\frac{\mathrm{d}^2 f\circ\overline{\mathscr{J}}^{n-1}\hat{\mathrm{y}}}{\mathrm{d}t^2}-\xi_{N,j}^{N+3}h^3\frac{N!}{(3+N)!}\frac{\mathrm{d}^3 f\circ\overline{\mathscr{J}}^{n-1}\hat{\mathrm{y}}}{\mathrm{d}t^3}\underset{h\to 0}{=}\mathcal{O}(h^4)$$

$$\Leftrightarrow\sum_{j'=1}^{s}w_{N,j,j'}(f+\frac{\partial f}{\partial t}\tau_{j'}h+\sum_{N=4}^{n}\sum_{k=1}^{l}\frac{\partial f}{\partial x_{N,k}}\big(y_{n,j',N,k}-y_{n,0,N,k}\big)$$

$$+\sum_{k=1}^{l}\frac{\partial f}{\partial x_{3,k}}\big(y_{n,j',3,k}-y_{n,0,3,k}\big)$$

$$+\sum_{k=1}^{l}\frac{\partial f}{\partial x_{2,k}}\big(y_{n,j',2,k}-y_{n,0,2,k}\big)$$

$$+\sum_{k=1}^{l}\frac{\partial f}{\partial x_{1,k}}\big(y_{n,j',1,k}-y_{n,0,1,k}\big)+\frac{1}{2}\frac{\partial^2 f}{\partial t^2}\big(\tau_{j'}h\big)^2+\sum_{N=3}^{n}\sum_{k=1}^{l}\frac{\partial^2 f}{\partial t\partial x_{N,k}}\tau_{j'}h\big(y_{n,j',N,k}-y_{n,0,N,k}\big)$$

$$+\sum_{k=1}^{l}\frac{\partial^2 f}{\partial t\partial x_{2,k}}\tau_{j'}h\big(y_{n,j',2,k}-y_{n,0,2,k}\big)+\sum_{k=1}^{l}\frac{\partial^2 f}{\partial t\partial x_{1,k}}\tau_{j'}h\big(y_{n,j',1,k}-y_{n,0,1,k}\big)$$

$$+\frac{1}{2}\sum_{N=3}^{n}\sum_{k=1}^{l}\sum_{N'=3}^{n}\sum_{k'=1}^{l}\frac{\partial^2 f}{\partial x_{N,k}\partial x_{N',k'}}\big(y_{n,j',N,k}-y_{n,0,N,k}\big)\big(y_{i',j',k'}-y_{n,0,N',k'}\big)$$

$$+\sum_{N=3}^{n}\sum_{k=1}^{l}\sum_{k'=1}^{l}\frac{\partial^2 f}{\partial x_{N,k}\partial x_{2,k'}}\big(y_{n,j',N,k}-y_{n,0,N,k}\big)\big(y_{n,j',2,k'}-y_{n,0,2,k'}\big)$$

$$+\sum_{N=3}^{n}\sum_{k=1}^{l}\sum_{k'=1}^{l}\frac{\partial^2 f}{\partial x_{N,k}\partial x_{1,k'}}\big(y_{n,j',N,k}-y_{n,0,N,k}\big)\big(y_{n,j',1,k'}-y_{n,0,1,k'}\big)$$

$$+\frac{1}{2}\sum_{k=1}^{l}\sum_{k'=1}^{l}\frac{\partial^2 f}{\partial x_{2,k}\partial x_{2,k'}}\big(y_{n,j',2,k}-y_{n,0,2,k}\big)\big(y_{n,j',2,k'}-y_{n,0,2,k'}\big)$$

$$+\sum_{k=1}^{l}\sum_{k'=1}^{l}\frac{\partial^2 f}{\partial x_{1,k}\partial x_{2,k'}}\big(y_{n,j',1,k}-y_{n,0,1,k}\big)\big(y_{n,j',2,k'}-y_{n,0,2,k'}\big)$$

$$+\frac{1}{2}\sum_{k=1}^{l}\sum_{k'=1}^{l}\frac{\partial^2 f}{\partial x_{1,k}\partial x_{1,k'}}\big(y_{n,j',1,k}-y_{n,0,1,k}\big)\big(y_{n,j',1,k'}-y_{n,0,1,k'}\big)$$

$$+\frac{1}{6}\frac{\partial^3 f}{\partial t^3}\big(\tau_{j'}h\big)^3+\frac{1}{2}\sum_{N=2}^{n}\sum_{k=1}^{l}\frac{\partial^3 f}{\partial t^2\partial x_{N,k}}\big(\tau_{j'}h\big)^2\big(y_{n,j',N,k}-y_{n,0,N,k}\big)+\frac{1}{2}\sum_{k=1}^{l}\frac{\partial^3 f}{\partial t^2\partial x_{1,k}}\big(\tau_{j'}h\big)^2\big(y_{n,j',1,k}-y_{n,0,1,k}\big)$$

$$+\frac{1}{2}\sum_{N=2}^{n}\sum_{k=1}^{l}\sum_{N'=2}^{n}\sum_{k'=1}^{l}\frac{\partial^3 f}{\partial t\partial x_{N,k}\partial x_{N',k'}}\tau_{j'}h\big(y_{n,j',N,k}-y_{n,0,N,k}\big)\big(y_{i',j',k'}-y_{n,0,N',k'}\big)$$

$$+\sum_{N=2}^{n}\sum_{k=1}^{l}\sum_{k'=1}^{l}\frac{\partial^3 f}{\partial t\partial x_{N,k}\partial x_{1,k'}}\tau_{j'}h\big(y_{n,j',N,k}-y_{n,0,N,k}\big)\big(y_{n,j',1,k'}-y_{n,0,1,k'}\big)$$

$$+\frac{1}{2}\sum_{k=1}^{l}\sum_{k'=1}^{l}\frac{\partial^3 f}{\partial t\partial x_{1,k}\partial x_{1,k'}}\tau_{j'}h\big(y_{n,j',1,k}-y_{n,0,1,k}\big)\big(y_{n,j',1,k'}-y_{n,0,1,k'}\big)+\frac{1}{6}\sum_{N=2}^{n}\sum_{k=1}^{l}\sum_{N'=2}^{n}\sum_{k'=1}^{l}\sum_{N''=2}^{n}\sum_{k''=1}^{l}\Bigg[$$

$$\frac{\partial^3 f}{\partial x_{N,k}\partial x_{N',k'}\partial x_{N'',k''}}\big(y_{n,j',N,k}-y_{n,0,N,k}\big)\big(y_{i',j',k'}-y_{n,0,N',k'}\big)\big(y_{n,j',N'',k''}-y_{n,0,N'',k''}\big)\Bigg]$$

$$+\frac{1}{2}\sum_{N=2}^{n}\sum_{k=1}^{l}\sum_{N'=2}^{n}\sum_{k'=1}^{l}\sum_{k''=1}^{l}\frac{\partial^3 f}{\partial x_{N,k}\partial x_{N',k'}\partial x_{1,k''}}\big(y_{n,j',N,k}-y_{n,0,N,k}\big)\big(y_{i',j',k'}-y_{n,0,N',k'}\big)\big(y_{n,j',1,k''}-y_{n,0,1,k''}\big)$$

$$+\frac{1}{2}\sum_{N=2}^{n}\sum_{k=1}^{l}\sum_{k'=1}^{l}\sum_{k''=1}^{l}\frac{\partial^3 f}{\partial x_{N,k}\partial x_{1,k'}\partial x_{1,k''}}\big(y_{n,j',N,k}-y_{n,0,N,k}\big)\big(y_{n,j',1,k'}-y_{n,0,1,k'}\big)\big(y_{n,j',1,k''}-y_{n,0,1,k''}\big)$$

$$+\frac{1}{6}\sum_{k=1}^{l}\sum_{k'=1}^{l}\sum_{k''=1}^{l}\frac{\partial^3 f}{\partial x_{1,k}\partial x_{1,k'}\partial x_{1,k''}}\big(y_{n,j',1,k}-y_{n,0,1,k}\big)\big(y_{n,j',1,k'}-y_{n,0,1,k'}\big)\big(y_{n,j',1,k''}-y_{n,0,1,k''}\big)\big)$$

$$-\xi_{N,j}^{N}f-\xi_{N,j}^{N+1}h\frac{1}{1+N}\frac{\mathrm{d}f\circ\overline{\mathcal{J}}^{n-1}\hat{\mathrm{y}}}{\mathrm{d}t}-\xi_{N,j}^{N+2}h^2\frac{N!}{(2+N)!}\frac{\mathrm{d}^2 f\circ\overline{\mathcal{J}}^{n-1}\hat{\mathrm{y}}}{\mathrm{d}t^2}-\xi_{N,j}^{N+3}h^3\frac{N!}{(3+N)!}\frac{\mathrm{d}^3 f\circ\overline{\mathcal{J}}^{n-1}\hat{\mathrm{y}}}{\mathrm{d}t^3}=\underset{h\to 0}{\mathcal{O}}(h^4)$$

$$\Leftrightarrow\sum_{j'=1}^{s}w_{N,j,j'}\Bigg[f$$

$$+\frac{\partial f}{\partial t}\tau_{j'}h$$

$$+\sum_{N=4}^{n}\sum_{k=1}^{l}\frac{\partial f}{\partial x_{N,k}}\left(hy_{n,0,N-1,k}\xi_{j',N}+h^2 y_{n,0,N-2,k}\frac{\xi_{j',N}^2}{2}+y_{n,0,N-3,k}h^3\frac{\xi_{j',N}^3}{6}\right)$$

$$+\sum_{k=1}^{l}\frac{\partial f}{\partial x_{3,k}}\left(hy_{n,0,2,k}\xi_{j',3}+h^2 y_{n,0,1,k}\frac{\xi_{j',3}^2}{2}+h^3 f_k\frac{1}{6}\sum_{j''=1}^{s}w_{3,j',j''}\right)$$

$$+\sum_{k=1}^{l}\frac{\partial f}{\partial x_{2,k}}\Bigg(hy_{n,0,1,k}\xi_{j',2}+h^2 f_k\frac{1}{2}\sum_{j''=1}^{s}w_{2,j',j''}+h^3\frac{\partial f_k}{\partial t}\frac{1}{2}\sum_{j''=1}^{s}w_{2,j',j''}\tau_{j''}$$

$$+h^3\sum_{N=2}^{n}\sum_{k'=1}^{l}\frac{\partial f_k}{\partial x_{N,k'}}y_{n,0,N-1,k'}\frac{1}{2}\sum_{j''=1}^{s}w_{2,j',j''}\xi_{j'',N}$$

$$+h^3\sum_{k'=1}^{l}\frac{\partial f_k}{\partial x_{1,k'}}f_{k'}\frac{1}{2}\sum_{j''=1}^{s}w_{2,j',j''}\sum_{j'''=1}^{s}w_{1,j'',j'''}\Bigg)$$

$$+\sum_{k=1}^{l}\frac{\partial f}{\partial x_{1,k}}\Bigg(hf_k\sum_{j''=1}^{s}w_{1,j',j''}$$

$$+h^2\frac{\partial f_k}{\partial t}\sum_{j''=1}^{s}w_{1,j',j''}\tau_{j''}$$

$$+h^2\sum_{N=3}^{n}\sum_{k'=1}^{l}\frac{\partial f_k}{\partial x_{N,k'}}y_{n,0,N-1,k'}\sum_{j''=1}^{s}w_{1,j',j''}\xi_{j'',N}$$

$$+h^2\sum_{k'=1}^{l}\frac{\partial f_k}{\partial x_{2,k'}}y_{n,0,1,k'}\sum_{j''=1}^{s}w_{1,j',j''}\xi_{j'',2}$$

$$+h^2 \sum_{k'=1}^{l} \frac{\partial f_k}{\partial x_{1,k'}} f_{k'} \sum_{j''=1}^{s} w_{1,j',j''} \sum_{j'''=1}^{s} w_{1,j'',j'''}$$

$$+h^3 \sum_{N=3}^{n}\sum_{k'=1}^{l} \frac{\partial f_k}{\partial x_{N,k'}} y_{n,0,N-2,k'} \sum_{j''=1}^{s} w_{1,j',j''} \frac{\xi^2_{j'',N}}{2}$$

$$+h^3 \sum_{k'=1}^{l} \frac{\partial f_k}{\partial x_{2,k'}} f_{k'} \frac{1}{2} \sum_{j''=1}^{s} w_{1,j',j''} \sum_{j'''=1}^{s} w_{2,j'',j'''}$$

$$+h^3 \sum_{k'=1}^{l} \frac{\partial f_k}{\partial x_{1,k'}} \frac{\partial f_{k'}}{\partial t} \sum_{j''=1}^{s} w_{1,j',j''} \sum_{j'''=1}^{s} w_{1,j'',j'''} \tau_{j'''}$$

$$+h^3 \sum_{k'=1}^{l} \frac{\partial f_k}{\partial x_{1,k'}} \sum_{N=2}^{n}\sum_{k''=1}^{l} \frac{\partial f_{k'}}{\partial x_{N,k''}} y_{n,0,N-1,k''} \sum_{j''=1}^{s} w_{1,j',j''} \sum_{j'''=1}^{s} w_{1,j'',j'''}\xi_{N,j'''}$$

$$+h^3 \sum_{k'=1}^{l} \frac{\partial f_k}{\partial x_{1,k'}} \sum_{k''=1}^{l} \frac{\partial f_{k'}}{\partial x_{1,k''}} f_{k''} \sum_{j''=1}^{s} w_{1,j',j''} \sum_{j'''=1}^{s} w_{1,j'',j'''} \sum_{j''''=1}^{s} w_{1,j''',j''''}$$

$$+h^3 \frac{\partial^2 f_k}{\partial t^2} \sum_{j''=1}^{s} w_{1,j',j''} \frac{\tau^2_{j''}}{2}$$

$$+h^3 \sum_{N=2}^{n}\sum_{k'=1}^{l} \frac{\partial^2 f_k}{\partial t \partial x_{N,k'}} y_{n,0,N-1,k'} \sum_{j''=1}^{s} w_{1,j',j''}\tau_{j''}\xi_{j'',N}$$

$$+h^3 \sum_{k'=1}^{l} \frac{\partial^2 f_k}{\partial t\partial x_{1,k'}} \sum_{j''=1}^{s} w_{1,j',j''}\tau_{j''} \sum_{j'''=1}^{s} w_{1,j'',j'''} f_{k'}$$

$$+h^3 \sum_{N=2}^{n}\sum_{k'=1}^{l}\sum_{N'=2}^{n}\sum_{k''=1}^{l} \frac{\partial^2 f_k}{\partial x_{N,k'}\partial x_{N',k''}} y_{n,0,N-1,k'} y_{n,0,N'-1,k''} \frac{1}{2} \sum_{j''=1}^{s} w_{1,j',j''}\xi_{j'',N}\xi_{N',j''}$$

$$+h^3 \sum_{N=2}^{n}\sum_{k'=1}^{l}\sum_{k''=1}^{l} \frac{\partial^2 f_k}{\partial x_{N,k'}\partial x_{1,k''}} y_{n,0,N-1,k'} f_{k''} \sum_{j''=1}^{s} w_{1,j',j''}\xi_{j'',N}\sum_{j'''=1}^{s} w_{1,j'',j'''}$$

$$+h^3 \sum_{k'=1}^{l}\sum_{k''=1}^{l} \frac{\partial^2 f_k}{\partial x_{1,k'}\partial x_{1,k''}} f_{k'}f_{k''}\frac{1}{2} \sum_{j''=1}^{s} w_{1,j',j''}\left(\sum_{j'''=1}^{s} w_{1,j'',j'''}\right)^2 \Bigg)$$

$$+\frac{1}{2} \frac{\partial^2 f}{\partial t^2}\left(\tau_{j'}h\right)^2$$

$$+\sum_{N=3}^{n}\sum_{k=1}^{l}\frac{\partial^2 f}{\partial t\partial x_{N,k}}\tau_{j'}h\left(hy_{n,0,N-1,k}\xi_{j',N}+h^2 y_{n,0,N-2,k}\frac{\xi^2_{j',N}}{2}\right)$$

$$+\sum_{k=1}^{l}\frac{\partial^2 f}{\partial t\partial x_{2,k}}\tau_{j'}h\left(hy_{n,0,1,k}\xi_{j',2}+h^2 f_k\frac{1}{2}\sum_{j''=1}^{s} w_{2,j',j''}\right)$$

$$+\sum_{k=1}^{l}\frac{\partial^2 f}{\partial t\partial x_{1,k}}\tau_{j'}h\Bigg(hf_k\sum_{j''=1}^{s}w_{1,j',j''}$$

$$+h^2\frac{\partial f_k}{\partial t}\sum_{j''=1}^{s}w_{1,j',j''}\tau_{j''}$$

$$+h^2 \sum_{N=3}^{n} \sum_{k'=1}^{l} \frac{\partial f_k}{\partial x_{N,k'}} y_{n,0,N-1,k'} \sum_{j''=1}^{s} w_{1,j',j''} \xi_{j'',N}$$

$$+h^2 \sum_{k'=1}^{l} \frac{\partial f_k}{\partial x_{2,k'}} y_{n,0,1,k'} \sum_{j''=1}^{s} w_{1,j',j''} \xi_{j'',2}$$

$$+h^2 \sum_{k'=1}^{l} \frac{\partial f_k}{\partial x_{1,k'}} f_{k'} \sum_{j''=1}^{s} w_{1,j',j''} \sum_{j'''=1}^{s} w_{1,j'',j'''} \Bigg)$$

$$+\frac{1}{2} \sum_{N=3}^{n} \sum_{k=1}^{l} \sum_{N'=3}^{n} \sum_{k'=1}^{l} \frac{\partial^2 f}{\partial x_{N,k} \partial x_{N',k'}} \left( h y_{n,0,N-1,k} \xi_{j',N} + h^2 y_{n,0,N-2,k} \frac{\xi_{j',N}^2}{2} \right) \left( h y_{n,0,N'-1,k'} \xi_{N',j'} + h^2 y_{n,0,N'-2,k'} \frac{\xi_{N',j'}^2}{2} \right)$$

$$+ \sum_{N=3}^{n} \sum_{k=1}^{l} \sum_{k'=1}^{l} \frac{\partial^2 f}{\partial x_{N,k} \partial x_{2,k'}} \left( h y_{n,0,N-1,k} \xi_{j',N} + h^2 y_{n,0,N-2,k} \frac{\xi_{j',N}^2}{2} \right) \left( h y_{n,0,1,k'} \xi_{j',2} + h^2 f_{k'} \frac{1}{2} \sum_{j''=1}^{s} w_{2,j',j''} \right)$$

$$+ \sum_{N=3}^{n} \sum_{k=1}^{l} \sum_{k'=1}^{l} \frac{\partial^2 f}{\partial x_{N,k} \partial x_{1,k'}} \left( h y_{n,0,N-1,k} \xi_{j',N} + h^2 y_{n,0,N-2,k} \frac{\xi_{j',N}^2}{2} \right) \Bigg( h f_{k'} \sum_{j''=1}^{s} w_{1,j',j''}$$

$$+h^2 \frac{\partial f_{k'}}{\partial t} \sum_{j''=1}^{s} w_{1,j',j''} \tau_{j''}$$

$$+h^2 \sum_{N'=3}^{n} \sum_{k''=1}^{l} \frac{\partial f_{k'}}{\partial x_{N',k''}} y_{n,0,N'-1,k''} \sum_{j''=1}^{s} w_{1,j',j''} \xi_{N',j''}$$

$$+h^2 \sum_{k''=1}^{l} \frac{\partial f_{k'}}{\partial x_{2,k''}} y_{n,0,1,k''} \sum_{j''=1}^{s} w_{1,j',j''} \xi_{j'',2}$$

$$+h^2 \sum_{k''=1}^{l} \frac{\partial f_{k'}}{\partial x_{1,k''}} f_{k''} \sum_{j''=1}^{s} w_{1,j',j''} \sum_{j'''=1}^{s} w_{1,j'',j'''} \Bigg)$$

$$+\frac{1}{2} \sum_{k=1}^{l} \sum_{k'=1}^{l} \frac{\partial^2 f}{\partial x_{2,k} \partial x_{2,k'}} \left( h y_{n,0,1,k} \xi_{j',2} + h^2 f_k \frac{1}{2} \sum_{j''=1}^{s} w_{2,j',j''} \right) \left( h y_{n,0,1,k'} \xi_{j',2} + h^2 f_{k'} \frac{1}{2} \sum_{j''=1}^{s} w_{2,j',j''} \right)$$

$$+ \sum_{k=1}^{l} \sum_{k'=1}^{l} \frac{\partial^2 f}{\partial x_{1,k} \partial x_{2,k'}} \Bigg( h f_k \sum_{j''=1}^{s} w_{1,j',j''}$$

$$+h^2 \frac{\partial f_k}{\partial t} \sum_{j''=1}^{s} w_{1,j',j''} \tau_{j''}$$

$$+h^2 \sum_{N=3}^{n} \sum_{k''=1}^{l} \frac{\partial f_k}{\partial x_{N,k''}} y_{n,0,N-1,k''} \sum_{j''=1}^{s} w_{1,j',j''} \xi_{j'',N}$$

$$+h^2 \sum_{k''=1}^{l} \frac{\partial f_k}{\partial x_{2,k''}} y_{n,0,1,k''} \sum_{j''=1}^{s} w_{1,j',j''} \xi_{j'',2}$$

$$+h^2 \sum_{k''=1}^{l} \frac{\partial f_k}{\partial x_{1,k''}} f_{k''} \sum_{j''=1}^{s} w_{1,j',j''} \sum_{j'''=1}^{s} w_{1,j'',j'''} \Bigg) \left( h y_{n,0,1,k'} \xi_{j',2} + h^2 f_{k'} \frac{1}{2} \sum_{j''=1}^{s} w_{2,j',j''} \right)$$

$$+\frac{1}{2}\sum_{k=1}^{l}\sum_{k'=1}^{l}\frac{\partial^2 f}{\partial x_{1,k}\partial x_{1,k'}}\Bigg(hf_k\sum_{j''=1}^{s} w_{1,j',j''}$$

$$+h^2\frac{\partial f_k}{\partial t}\sum_{j''=1}^{s} w_{1,j',j''}\tau_{j''}$$

$$+h^2\sum_{N=3}^{n}\sum_{k''=1}^{l}\frac{\partial f_k}{\partial x_{N,k''}}y_{n,0,N-1,k''}\sum_{j''=1}^{s} w_{1,j',j''}\xi_{j'',N}$$

$$+h^2\sum_{k''=1}^{l}\frac{\partial f_k}{\partial x_{2,k''}}y_{n,0,1,k''}\sum_{j''=1}^{s} w_{1,j',j''}\xi_{j'',2}$$

$$+h^2\sum_{k''=1}^{l}\frac{\partial f_k}{\partial x_{1,k''}}f_{k''}\sum_{j''=1}^{s} w_{1,j',j''}\sum_{j'''=1}^{s} w_{1,j'',j'''}\Bigg)\Bigg(hf_{k'}\sum_{j''=1}^{s} w_{1,j',j''}$$

$$+h^2\frac{\partial f_{k'}}{\partial t}\sum_{j''=1}^{s} w_{1,j',j''}\tau_{j''}$$

$$+h^2\sum_{N'=3}^{n}\sum_{k''=1}^{l}\frac{\partial f_{k'}}{\partial x_{N',k''}}y_{n,0,N'-1,k''}\sum_{j''=1}^{s} w_{1,j',j''}\xi_{N',j''}$$

$$+h^2\sum_{k''=1}^{l}\frac{\partial f_{k'}}{\partial x_{2,k''}}y_{n,0,1,k''}\sum_{j''=1}^{s} w_{1,j',j''}\xi_{j'',2}$$

$$+h^2\sum_{k''=1}^{l}\frac{\partial f_{k'}}{\partial x_{1,k''}}f_{k''}\sum_{j''=1}^{s} w_{1,j',j''}\sum_{j'''=1}^{s} w_{1,j'',j'''}\Bigg)$$

$$+\frac{1}{6}\frac{\partial^3 f}{\partial t^3}\big(\tau_{j'}h\big)^3$$

$$+\frac{1}{2}\sum_{N=2}^{n}\sum_{k=1}^{l}\frac{\partial^3 f}{\partial t^2\partial x_{N,k}}\big(\tau_{j'}h\big)^2\big(hy_{n,0,N-1,k}\xi_{j',N}\big)$$

$$+\frac{1}{2}\sum_{k=1}^{l}\frac{\partial^3 f}{\partial t^2\partial x_{1,k}}\big(\tau_{j'}h\big)^2\big(y_{n,j',1,k}-y_{n,0,1,k}\big)$$

$$+\frac{1}{2}\sum_{N=2}^{n}\sum_{k=1}^{l}\sum_{N'=2}^{n}\sum_{k'=1}^{l}\frac{\partial^3 f}{\partial t\partial x_{N,k}\partial x_{N',k'}}\tau_{j'}h\big(hy_{n,0,N-1,k}\xi_{j',N}\big)\big(hy_{n,0,N'-1,k'}\xi_{N',j'}\big)$$

$$+\sum_{N=2}^{n}\sum_{k=1}^{l}\sum_{k'=1}^{l}\frac{\partial^3 f}{\partial t\partial x_{N,k}\partial x_{1,k'}}\tau_{j'}h\big(hy_{n,0,N-1,k}\xi_{j',N}\big)\Bigg(hf_{k'}\sum_{j''=1}^{s} w_{1,j',j''}\Bigg)$$

$$+\frac{1}{2}\sum_{k=1}^{l}\sum_{k'=1}^{l}\frac{\partial^3 f}{\partial t\partial x_{1,k}\partial x_{1,k'}}\tau_{j'}h\Bigg(hf_k\sum_{j''=1}^{s} w_{1,j',j''}\Bigg)\Bigg(hf_{k'}\sum_{j''=1}^{s} w_{1,j',j''}\Bigg)$$

$$+\frac{1}{6}\sum_{N=2}^{n}\sum_{k=1}^{l}\sum_{N'=2}^{n}\sum_{k'=1}^{l}\sum_{N''=2}^{n}\sum_{k''=1}^{l}\frac{\partial^3 f}{\partial x_{N,k}\partial x_{N',k'}\partial x_{N'',k''}}\big(hy_{n,0,N-1,k}\xi_{j',N}\big)\big(hy_{n,0,N'-1,k'}\xi_{N',j'}\big)\big(hy_{n,0,N''-1,k''}\xi_{N'',j'}\big)$$

$$+\frac{1}{2}\sum_{N=2}^{n}\sum_{k=1}^{l}\sum_{N'=2}^{n}\sum_{k'=1}^{l}\sum_{k''=1}^{l}\frac{\partial^3 f}{\partial x_{N,k}\partial x_{N',k'}\partial x_{1,k''}}\big(hy_{n,0,N-1,k}\xi_{j',N}\big)\big(hy_{n,0,N'-1,k'}\xi_{N',j'}\big)\Bigg(hf_{k''}\sum_{j''=1}^{s} w_{1,j',j''}\Bigg)$$

$$+\frac{1}{2}\sum_{N=2}^{n}\sum_{k=1}^{l}\sum_{k'=1}^{l}\sum_{k''=1}^{l}\frac{\partial^3 f}{\partial x_{N,k}\partial x_{1,k'}\partial x_{1,k''}}\left(hy_{n,0,N-1,k}\xi_{j',N}\right)\left(hf_{k'}\sum_{j''=1}^{s}w_{1,j',j''}\right)\left(hf_{k''}\sum_{j''=1}^{s}w_{1,j',j''}\right)$$

$$+\frac{1}{6}\sum_{k=1}^{l}\sum_{k'=1}^{l}\sum_{k''=1}^{l}\frac{\partial^3 f}{\partial x_{1,k}\partial x_{1,k'}\partial x_{1,k''}}\left(hf_{k}\sum_{j''=1}^{s}w_{1,j',j''}\right)\left(hf_{k'}\sum_{j''=1}^{s}w_{1,j',j''}\right)\left(hf_{k''}\sum_{j''=1}^{s}w_{1,j',j''}\right)\Bigg]$$

$$-\xi_{N,j}^{N}f-\xi_{N,j}^{N+1}h\frac{1}{1+N}\frac{\mathrm{d}f\circ\overline{\mathcal{J}}^{n-1}\hat{\mathrm{y}}}{\mathrm{d}t}-\xi_{N,j}^{N+2}h^2\frac{N!}{(2+N)!}\frac{\mathrm{d}^2f\circ\overline{\mathcal{J}}^{n-1}\hat{\mathrm{y}}}{\mathrm{d}t^2}-\xi_{N,j}^{N+3}h^3\frac{N!}{(3+N)!}\frac{\mathrm{d}^3f\circ\overline{\mathcal{J}}^{n-1}\hat{\mathrm{y}}}{\mathrm{d}t^3}=\underset{h\to0}{\mathcal{O}}(h^4)$$

$$\Leftrightarrow\sum_{j'=1}^{s}w_{N,j,j'}\Bigg[$$

$$f$$

$$+\frac{\partial f}{\partial t}\tau_{j'}h$$

$$+\sum_{N=4}^{n}\sum_{k=1}^{l}\frac{\partial f}{\partial x_{N,k}}hy_{n,0,N-1,k}\xi_{j',N}$$

$$+\sum_{N=4}^{n}\sum_{k=1}^{l}\frac{\partial f}{\partial x_{N,k}}h^2y_{n,0,N-2,k}\frac{\xi_{j',N}^2}{2}$$

$$+\sum_{N=4}^{n}\sum_{k=1}^{l}\frac{\partial f}{\partial x_{N,k}}y_{n,0,N-3,k}h^3\frac{\xi_{j',N}^3}{6}$$

$$+\sum_{k=1}^{l}\frac{\partial f}{\partial x_{3,k}}hy_{n,0,2,k}\xi_{j',3}$$

$$+\sum_{k=1}^{l}\frac{\partial f}{\partial x_{3,k}}h^2y_{n,0,1,k}\frac{\xi_{j',3}^2}{2}$$

$$+\sum_{k=1}^{l}\frac{\partial f}{\partial x_{3,k}}h^3f_k\frac{1}{6}\sum_{j''=1}^{s}w_{3,j',j''}$$

$$+\sum_{k=1}^{l}\frac{\partial f}{\partial x_{2,k}}hy_{n,0,1,k}\xi_{j',2}$$

$$+\sum_{k=1}^{l}\frac{\partial f}{\partial x_{2,k}}h^2f_k\frac{1}{2}\sum_{j''=1}^{s}w_{2,j',j''}$$

$$+\sum_{k=1}^{l}\frac{\partial f}{\partial x_{2,k}}h^3\frac{\partial f_k}{\partial t}\frac{1}{2}\sum_{j''=1}^{s}w_{2,j',j''}\tau_{j''}$$

$$+\sum_{k=1}^{l}\frac{\partial f}{\partial x_{2,k}}h^3\sum_{N=2}^{n}\sum_{k'=1}^{l}\frac{\partial f_k}{\partial x_{N,k'}}y_{n,0,N-1,k'}\frac{1}{2}\sum_{j''=1}^{s}w_{2,j',j''}\xi_{j'',N}$$

$$+\sum_{k=1}^{l}\frac{\partial f}{\partial x_{2,k}}h^3\sum_{k'=1}^{l}\frac{\partial f_k}{\partial x_{1,k'}}f_{k'}\frac{1}{2}\sum_{j''=1}^{s}w_{2,j',j''}\sum_{j'''=1}^{s}w_{1,j'',j'''}$$

$$+\sum_{k=1}^{l}\frac{\partial f}{\partial x_{1,k}}hf_k\sum_{j''=1}^{s}w_{1,j',j''}$$

$$
\begin{aligned}
&+\sum_{k=1}^{l}\frac{\partial f}{\partial x_{1,k}}h^2\frac{\partial f_k}{\partial t}\sum_{j''=1}^{s}w_{1,j',j''}\tau_{j''}\\
&+\sum_{k=1}^{l}\frac{\partial f}{\partial x_{1,k}}h^2\sum_{N=3}^{n}\sum_{k'=1}^{l}\frac{\partial f_k}{\partial x_{N,k'}}y_{n,0,N-1,k'}\sum_{j''=1}^{s}w_{1,j',j''}\xi_{j'',N}\\
&+\sum_{k=1}^{l}\frac{\partial f}{\partial x_{1,k}}h^2\sum_{k'=1}^{l}\frac{\partial f_k}{\partial x_{2,k'}}y_{n,0,1,k'}\sum_{j''=1}^{s}w_{1,j',j''}\xi_{j'',2}\\
&+\sum_{k=1}^{l}\frac{\partial f}{\partial x_{1,k}}h^2\sum_{k'=1}^{l}\frac{\partial f_k}{\partial x_{1,k'}}f_{k'}\sum_{j''=1}^{s}w_{1,j',j''}\sum_{j'''=1}^{s}w_{1,j'',j'''}\\
&+\sum_{k=1}^{l}\frac{\partial f}{\partial x_{1,k}}h^3\sum_{N=3}^{n}\sum_{k'=1}^{l}\frac{\partial f_k}{\partial x_{N,k'}}y_{n,0,N-2,k'}\sum_{j''=1}^{s}w_{1,j',j''}\frac{\xi_{j'',N}^2}{2}\\
&+\sum_{k=1}^{l}\frac{\partial f}{\partial x_{1,k}}h^3\sum_{k'=1}^{l}\frac{\partial f_k}{\partial x_{2,k'}}f_{k'}\frac{1}{2}\sum_{j''=1}^{s}w_{1,j',j''}\sum_{j'''=1}^{s}w_{2,j'',j'''}\\
&+\sum_{k=1}^{l}\frac{\partial f}{\partial x_{1,k}}h^3\sum_{k'=1}^{l}\frac{\partial f_k}{\partial x_{1,k'}}\frac{\partial f_{k'}}{\partial t}\sum_{j''=1}^{s}w_{1,j',j''}\sum_{j'''=1}^{s}w_{1,j'',j'''}\tau_{j'''}\\
&+\sum_{k=1}^{l}\frac{\partial f}{\partial x_{1,k}}h^3\sum_{k'=1}^{l}\frac{\partial f_k}{\partial x_{1,k'}}\sum_{N=2}^{n}\sum_{k''=1}^{l}\frac{\partial f_{k'}}{\partial x_{N,k''}}y_{n,0,N-1,k''}\sum_{j''=1}^{s}w_{1,j',j''}\sum_{j'''=1}^{s}w_{1,j'',j'''}\xi_{N,j'''}\\
&+\sum_{k=1}^{l}\frac{\partial f}{\partial x_{1,k}}h^3\sum_{k'=1}^{l}\frac{\partial f_k}{\partial x_{1,k'}}\sum_{k''=1}^{l}\frac{\partial f_{k'}}{\partial x_{1,k''}}f_{k''}\sum_{j''=1}^{s}w_{1,j',j''}\sum_{j'''=1}^{s}w_{1,j'',j'''}\sum_{j''''=1}^{s}w_{1,j''',j''''}\\
&+\sum_{k=1}^{l}\frac{\partial f}{\partial x_{1,k}}h^3\frac{\partial^2 f_k}{\partial t^2}\sum_{j''=1}^{s}w_{1,j',j''}\frac{\tau_{j''}^2}{2}\\
&+\sum_{k=1}^{l}\frac{\partial f}{\partial x_{1,k}}h^3\sum_{N=2}^{n}\sum_{k'=1}^{l}\frac{\partial^2 f_k}{\partial t\partial x_{N,k'}}y_{n,0,N-1,k'}\sum_{j''=1}^{s}w_{1,j',j''}\tau_{j''}\xi_{j'',N}\\
&+\sum_{k=1}^{l}\frac{\partial f}{\partial x_{1,k}}h^3\sum_{k'=1}^{l}\frac{\partial^2 f_k}{\partial t\partial x_{1,k'}}\sum_{j''=1}^{s}w_{1,j',j''}\tau_{j''}\sum_{j'''=1}^{s}w_{1,j'',j'''}f_{k'}\\
&+\sum_{k=1}^{l}\frac{\partial f}{\partial x_{1,k}}h^3\sum_{N=2}^{n}\sum_{k'=1}^{l}\sum_{N'=2}^{n}\sum_{k''=1}^{l}\frac{\partial^2 f_k}{\partial x_{N,k'}\partial x_{N',k''}}y_{n,0,N-1,k'}y_{n,0,N'-1,k''}\frac{1}{2}\sum_{j''=1}^{s}w_{1,j',j''}\xi_{j'',N}\xi_{N',j''}\\
&+\sum_{k=1}^{l}\frac{\partial f}{\partial x_{1,k}}h^3\sum_{N=2}^{n}\sum_{k'=1}^{l}\sum_{k''=1}^{l}\frac{\partial^2 f_k}{\partial x_{N,k'}\partial x_{1,k''}}y_{n,0,N-1,k'}f_{k''}\sum_{j''=1}^{s}w_{1,j',j''}\xi_{j'',N}\sum_{j'''=1}^{s}w_{1,j'',j'''}\\
&+\sum_{k=1}^{l}\frac{\partial f}{\partial x_{1,k}}h^3\sum_{k'=1}^{l}\sum_{k''=1}^{l}\frac{\partial^2 f_k}{\partial x_{1,k'}\partial x_{1,k''}}f_{k'}f_{k''}\frac{1}{2}\sum_{j''=1}^{s}w_{1,j',j''}\left(\sum_{j'''=1}^{s}w_{1,j'',j'''}\right)^2\\
&+\frac{1}{2}\frac{\partial^2 f}{\partial t^2}\left(\tau_{j'}h\right)^2\\
&+\sum_{N=3}^{n}\sum_{k=1}^{l}\frac{\partial^2 f}{\partial t\partial x_{N,k}}\tau_{j'}h\left(hy_{n,0,N-1,k}\xi_{j',N}\right)
\end{aligned}
$$

$$+\sum_{N=3}^{n}\sum_{k=1}^{l}\frac{\partial^2 f}{\partial t\partial x_{N,k}}\tau_{j'}h\left(h^2 y_{n,0,N-2,k}\frac{\xi_{j',N}^2}{2}\right)$$

$$+\sum_{k=1}^{l}\frac{\partial^2 f}{\partial t\partial x_{2,k}}\tau_{j'}h\left(h y_{n,0,1,k}\xi_{j',2}\right)$$

$$+\sum_{k=1}^{l}\frac{\partial^2 f}{\partial t\partial x_{2,k}}\tau_{j'}h\left(h^2 f_k\frac{1}{2}\sum_{j''=1}^{s}w_{2,j',j''}\right)$$

$$+\sum_{k=1}^{l}\frac{\partial^2 f}{\partial t\partial x_{1,k}}\tau_{j'}hhf_k\sum_{j''=1}^{s}w_{1,j',j''}$$

$$+\sum_{k=1}^{l}\frac{\partial^2 f}{\partial t\partial x_{1,k}}\tau_{j'}hh^2\frac{\partial f_k}{\partial t}\sum_{j''=1}^{s}w_{1,j',j''}\tau_{j''}$$

$$+\sum_{k=1}^{l}\frac{\partial^2 f}{\partial t\partial x_{1,k}}\tau_{j'}hh^2\sum_{N=3}^{n}\sum_{k'=1}^{l}\frac{\partial f_k}{\partial x_{N,k'}}y_{n,0,N-1,k'}\sum_{j''=1}^{s}w_{1,j',j''}\xi_{j'',N}$$

$$+\sum_{k=1}^{l}\frac{\partial^2 f}{\partial t\partial x_{1,k}}\tau_{j'}hh^2\sum_{k'=1}^{l}\frac{\partial f_k}{\partial x_{2,k'}}y_{n,0,1,k'}\sum_{j''=1}^{s}w_{1,j',j''}\xi_{j'',2}$$

$$+\sum_{k=1}^{l}\frac{\partial^2 f}{\partial t\partial x_{1,k}}\tau_{j'}hh^2\sum_{k'=1}^{l}\frac{\partial f_k}{\partial x_{1,k'}}f_{k'}\sum_{j''=1}^{s}w_{1,j',j''}\sum_{j'''=1}^{s}w_{1,j'',j'''}$$

$$+\frac{1}{2}\sum_{N=3}^{n}\sum_{k=1}^{l}\sum_{N'=3}^{n}\sum_{k'=1}^{l}\frac{\partial^2 f}{\partial x_{N,k}\partial x_{N',k'}}\left(h y_{n,0,N-1,k}\xi_{j',N}\right)\left(h y_{n,0,N'-1,k'}\xi_{N',j'}\right)$$

$$+\frac{1}{2}\sum_{N=3}^{n}\sum_{k=1}^{l}\sum_{N'=3}^{n}\sum_{k'=1}^{l}\frac{\partial^2 f}{\partial x_{N,k}\partial x_{N',k'}}\left(h^2 y_{n,0,N-2,k}\frac{\xi_{j',N}^2}{2}\right)\left(h y_{n,0,N'-1,k'}\xi_{N',j'}\right)$$

$$+\frac{1}{2}\sum_{N=3}^{n}\sum_{k=1}^{l}\sum_{N'=3}^{n}\sum_{k'=1}^{l}\frac{\partial^2 f}{\partial x_{N,k}\partial x_{N',k'}}\left(h y_{n,0,N-1,k}\xi_{j',N}\right)\left(h^2 y_{n,0,N'-2,k'}\frac{\xi_{N',j'}^2}{2}\right)$$

$$+\sum_{N=3}^{n}\sum_{k=1}^{l}\sum_{k'=1}^{l}\frac{\partial^2 f}{\partial x_{N,k}\partial x_{2,k'}}\left(h y_{n,0,N-1,k}\xi_{j',N}\right)\left(h y_{n,0,1,k'}\xi_{j',2}\right)$$

$$+\sum_{N=3}^{n}\sum_{k=1}^{l}\sum_{k'=1}^{l}\frac{\partial^2 f}{\partial x_{N,k}\partial x_{2,k'}}\left(h^2 y_{n,0,N-2,k}\frac{\xi_{j',N}^2}{2}\right)\left(h y_{n,0,1,k'}\xi_{j',2}\right)$$

$$+\sum_{N=3}^{n}\sum_{k=1}^{l}\sum_{k'=1}^{l}\frac{\partial^2 f}{\partial x_{N,k}\partial x_{2,k'}}\left(h y_{n,0,N-1,k}\xi_{j',N}\right)\left(h^2 f_{k'}\frac{1}{2}\sum_{j''=1}^{s}w_{2,j',j''}\right)$$

$$+\sum_{N=3}^{n}\sum_{k=1}^{l}\sum_{k'=1}^{l}\frac{\partial^2 f}{\partial x_{N,k}\partial x_{1,k'}}\left(h y_{n,0,N-1,k}\xi_{j',N}\right)hf_{k'}\sum_{j''=1}^{s}w_{1,j',j''}$$

$$+\sum_{N=3}^{n}\sum_{k=1}^{l}\sum_{k'=1}^{l}\frac{\partial^2 f}{\partial x_{N,k}\partial x_{1,k'}}\left(h y_{n,0,N-1,k}\xi_{j',N}\right)h^2\frac{\partial f_{k'}}{\partial t}\sum_{j''=1}^{s}w_{1,j',j''}\tau_{j''}$$

$$+\sum_{N=3}^{n}\sum_{k=1}^{l}\sum_{k'=1}^{l}\frac{\partial^2 f}{\partial x_{N,k}\partial x_{1,k'}}\left(h y_{n,0,N-1,k}\xi_{j',N}\right)h^2\sum_{N'=3}^{n}\sum_{k''=1}^{l}\frac{\partial f_{k'}}{\partial x_{N',k''}}y_{n,0,N'-1,k''}\sum_{j''=1}^{s}w_{1,j',j''}\xi_{N',j''}$$

$$+\sum_{N=3}^{n}\sum_{k=1}^{l}\sum_{k'=1}^{l}\frac{\partial^2 f}{\partial x_{N,k}\partial x_{1,k'}}\left(hy_{n,0,N-1,k}\xi_{j',N}\right)h^2\sum_{k''=1}^{l}\frac{\partial f_{k'}}{\partial x_{2,k''}}y_{n,0,1,k''}\sum_{j''=1}^{s}w_{1,j',j''}\xi_{j'',2}$$

$$+\sum_{N=3}^{n}\sum_{k=1}^{l}\sum_{k'=1}^{l}\frac{\partial^2 f}{\partial x_{N,k}\partial x_{1,k'}}\left(hy_{n,0,N-1,k}\xi_{j',N}\right)h^2\sum_{k''=1}^{l}\frac{\partial f_{k'}}{\partial x_{1,k''}}f_{k''}\sum_{j''=1}^{s}w_{1,j',j''}\sum_{j'''=1}^{s}w_{1,j'',j'''}$$

$$+\sum_{N=3}^{n}\sum_{k=1}^{l}\sum_{k'=1}^{l}\frac{\partial^2 f}{\partial x_{N,k}\partial x_{1,k'}}\left(h^2y_{n,0,N-2,k}\frac{\xi^2_{j',N}}{2}\right)\left(hf_{k'}\sum_{j''=1}^{s}w_{1,j',j''}\right)$$

$$+\frac{1}{2}\sum_{k=1}^{l}\sum_{k'=1}^{l}\frac{\partial^2 f}{\partial x_{2,k}\partial x_{2,k'}}\left(hy_{n,0,1,k}\xi_{j',2}\right)\left(hy_{n,0,1,k'}\xi_{j',2}\right)$$

$$+\frac{1}{2}\sum_{k=1}^{l}\sum_{k'=1}^{l}\frac{\partial^2 f}{\partial x_{2,k}\partial x_{2,k'}}\left(h^2f_k\frac{1}{2}\sum_{j''=1}^{s}w_{2,j',j''}\right)\left(hy_{n,0,1,k'}\xi_{j',2}\right)$$

$$+\frac{1}{2}\sum_{k=1}^{l}\sum_{k'=1}^{l}\frac{\partial^2 f}{\partial x_{2,k}\partial x_{2,k'}}\left(hy_{n,0,1,k}\xi_{j',2}\right)\left(h^2f_{k'}\frac{1}{2}\sum_{j''=1}^{s}w_{2,j',j''}\right)$$

$$+\sum_{k=1}^{l}\sum_{k'=1}^{l}\frac{\partial^2 f}{\partial x_{1,k}\partial x_{2,k'}}hf_k\sum_{j''=1}^{s}w_{1,j',j''}\left(hy_{n,0,1,k'}\xi_{j',2}\right)$$

$$+\sum_{k=1}^{l}\sum_{k'=1}^{l}\frac{\partial^2 f}{\partial x_{1,k}\partial x_{2,k'}}h^2\frac{\partial f_k}{\partial t}\sum_{j''=1}^{s}w_{1,j',j''}\tau_{j''}\left(hy_{n,0,1,k'}\xi_{j',2}\right)$$

$$+\sum_{k=1}^{l}\sum_{k'=1}^{l}\frac{\partial^2 f}{\partial x_{1,k}\partial x_{2,k'}}h^2\sum_{N=3}^{n}\sum_{k''=1}^{l}\frac{\partial f_k}{\partial x_{N,k''}}y_{n,0,N-1,k''}\sum_{j''=1}^{s}w_{1,j',j''}\xi_{j'',N}\left(hy_{n,0,1,k'}\xi_{j',2}\right)$$

$$+\sum_{k=1}^{l}\sum_{k'=1}^{l}\frac{\partial^2 f}{\partial x_{1,k}\partial x_{2,k'}}h^2\sum_{k''=1}^{l}\frac{\partial f_k}{\partial x_{2,k''}}y_{n,0,1,k''}\sum_{j''=1}^{s}w_{1,j',j''}\xi_{j'',2}\left(hy_{n,0,1,k'}\xi_{j',2}\right)$$

$$+\sum_{k=1}^{l}\sum_{k'=1}^{l}\frac{\partial^2 f}{\partial x_{1,k}\partial x_{2,k'}}h^2\sum_{k''=1}^{l}\frac{\partial f_k}{\partial x_{1,k''}}f_{k''}\sum_{j''=1}^{s}w_{1,j',j''}\sum_{j'''=1}^{s}w_{1,j'',j'''}\left(hy_{n,0,1,k'}\xi_{j',2}\right)$$

$$+\sum_{k=1}^{l}\sum_{k'=1}^{l}\frac{\partial^2 f}{\partial x_{1,k}\partial x_{2,k'}}hf_k\sum_{j''=1}^{s}w_{1,j',j''}\left(h^2f_{k'}\frac{1}{2}\sum_{j''=1}^{s}w_{2,j',j''}\right)$$

$$+\sum_{k=1}^{l}\sum_{k'=1}^{l}\frac{\partial^2 f}{\partial x_{1,k}\partial x_{2,k'}}h^2\frac{\partial f_k}{\partial t}\sum_{j''=1}^{s}w_{1,j',j''}\tau_{j''}\left(h^2f_{k'}\frac{1}{2}\sum_{j''=1}^{s}w_{2,j',j''}\right)$$

$$+\sum_{k=1}^{l}\sum_{k'=1}^{l}\frac{\partial^2 f}{\partial x_{1,k}\partial x_{2,k'}}h^2\sum_{N=3}^{n}\sum_{k''=1}^{l}\frac{\partial f_k}{\partial x_{N,k''}}y_{n,0,N-1,k''}\sum_{j''=1}^{s}w_{1,j',j''}\xi_{j'',N}\left(h^2f_{k'}\frac{1}{2}\sum_{j''=1}^{s}w_{2,j',j''}\right)$$

$$+\sum_{k=1}^{l}\sum_{k'=1}^{l}\frac{\partial^2 f}{\partial x_{1,k}\partial x_{2,k'}}h^2\sum_{k''=1}^{l}\frac{\partial f_k}{\partial x_{2,k''}}y_{n,0,1,k''}\sum_{j''=1}^{s}w_{1,j',j''}\xi_{j'',2}\left(h^2f_{k'}\frac{1}{2}\sum_{j''=1}^{s}w_{2,j',j''}\right)$$

$$+\sum_{k=1}^{l}\sum_{k'=1}^{l}\frac{\partial^2 f}{\partial x_{1,k}\partial x_{2,k'}}h^2\sum_{k''=1}^{l}\frac{\partial f_k}{\partial x_{1,k''}}f_{k''}\sum_{j''=1}^{s}w_{1,j',j''}\sum_{j'''=1}^{s}w_{1,j'',j'''}\left(h^2f_{k'}\frac{1}{2}\sum_{j''=1}^{s}w_{2,j',j''}\right)$$

$$+\frac{1}{2}\sum_{k=1}^{l}\sum_{k'=1}^{l}\frac{\partial^2 f}{\partial x_{1,k}\partial x_{1,k'}}\left(hf_k\sum_{j''=1}^{s}w_{1,j',j''}\right)hf_{k'}\sum_{j''=1}^{s}w_{1,j',j''}$$

$$+\frac{1}{2}\sum_{k=1}^{l}\sum_{k'=1}^{l}\frac{\partial^2 f}{\partial x_{1,k}\partial x_{1,k'}}\left(hf_k\sum_{j''=1}^{s}w_{1,j',j''}\right)h^2\frac{\partial f_{k'}}{\partial t}\sum_{j''=1}^{s}w_{1,j',j''}\tau_{j''}$$

$$+\frac{1}{2}\sum_{k=1}^{l}\sum_{k'=1}^{l}\frac{\partial^2 f}{\partial x_{1,k}\partial x_{1,k'}}\left(hf_k\sum_{j''=1}^{s}w_{1,j',j''}\right)h^2\sum_{N=3}^{n}\sum_{k''=1}^{l}\frac{\partial f_{k'}}{\partial x_{N,k''}}y_{n,0,N-1,k''}\sum_{j''=1}^{s}w_{1,j',j''}\xi_{j'',N}$$

$$+\frac{1}{2}\sum_{k=1}^{l}\sum_{k'=1}^{l}\frac{\partial^2 f}{\partial x_{1,k}\partial x_{1,k'}}\left(hf_k\sum_{j''=1}^{s}w_{1,j',j''}\right)h^2\sum_{k''=1}^{l}\frac{\partial f_{k'}}{\partial x_{2,k''}}y_{n,0,1,k''}\sum_{j''=1}^{s}w_{1,j',j''}\xi_{j'',2}$$

$$+\frac{1}{2}\sum_{k=1}^{l}\sum_{k'=1}^{l}\frac{\partial^2 f}{\partial x_{1,k}\partial x_{1,k'}}\left(hf_k\sum_{j''=1}^{s}w_{1,j',j''}\right)h^2\sum_{k''=1}^{l}\frac{\partial f_{k'}}{\partial x_{1,k''}}f_{k''}\sum_{j''=1}^{s}w_{1,j',j''}\sum_{j'''=1}^{s}w_{1,j'',j'''}$$

$$+\frac{1}{2}\sum_{k=1}^{l}\sum_{k'=1}^{l}\frac{\partial^2 f}{\partial x_{1,k}\partial x_{1,k'}}h^2\frac{\partial f_k}{\partial t}\sum_{j''=1}^{s}w_{1,j',j''}\tau_{j''}\left(hf_{k'}\sum_{j''=1}^{s}w_{1,j',j''}\right)$$

$$+\frac{1}{2}\sum_{k=1}^{l}\sum_{k'=1}^{l}\frac{\partial^2 f}{\partial x_{1,k}\partial x_{1,k'}}h^2\sum_{N=3}^{n}\sum_{k''=1}^{l}\frac{\partial f_k}{\partial x_{N,k''}}y_{n,0,N-1,k''}\sum_{j''=1}^{s}w_{1,j',j''}\xi_{j'',N}\left(hf_{k'}\sum_{j''=1}^{s}w_{1,j',j''}\right)$$

$$+\frac{1}{2}\sum_{k=1}^{l}\sum_{k'=1}^{l}\frac{\partial^2 f}{\partial x_{1,k}\partial x_{1,k'}}h^2\sum_{k''=1}^{l}\frac{\partial f_k}{\partial x_{2,k''}}y_{n,0,1,k''}\sum_{j''=1}^{s}w_{1,j',j''}\xi_{j'',2}\left(hf_{k'}\sum_{j''=1}^{s}w_{1,j',j''}\right)$$

$$+\frac{1}{2}\sum_{k=1}^{l}\sum_{k'=1}^{l}\frac{\partial^2 f}{\partial x_{1,k}\partial x_{1,k'}}h^2\sum_{k''=1}^{l}\frac{\partial f_k}{\partial x_{1,k''}}f_{k''}\sum_{j''=1}^{s}w_{1,j',j''}\sum_{j'''=1}^{s}w_{1,j'',j'''}\left(hf_{k'}\sum_{j''=1}^{s}w_{1,j',j''}\right)$$

$$+\frac{1}{6}\frac{\partial^3 f}{\partial t^3}\left(\tau_{j'}h\right)^3$$

$$+\frac{1}{2}\sum_{N=2}^{n}\sum_{k=1}^{l}\frac{\partial^3 f}{\partial t^2\partial x_{N,k}}\left(\tau_{j'}h\right)^2\left(hy_{n,0,N-1,k}\xi_{j',N}\right)$$

$$+\frac{1}{2}\sum_{k=1}^{l}\frac{\partial^3 f}{\partial t^2\partial x_{1,k}}\left(\tau_{j'}h\right)^2\left(hf_k\sum_{j''=1}^{s}w_{1,j',j''}\right)$$

$$+\frac{1}{2}\sum_{N=2}^{n}\sum_{k=1}^{l}\sum_{N'=2}^{n}\sum_{k'=1}^{l}\frac{\partial^3 f}{\partial t\partial x_{N,k}\partial x_{N',k'}}\tau_{j'}h\left(hy_{n,0,N-1,k}\xi_{j',N}\right)\left(hy_{n,0,N'-1,k'}\xi_{N',j'}\right)$$

$$+\sum_{N=2}^{n}\sum_{k=1}^{l}\sum_{k'=1}^{l}\frac{\partial^3 f}{\partial t\partial x_{N,k}\partial x_{1,k'}}\tau_{j'}h\left(hy_{n,0,N-1,k}\xi_{j',N}\right)\left(hf_{k'}\sum_{j''=1}^{s}w_{1,j',j''}\right)$$

$$+\frac{1}{2}\sum_{k=1}^{l}\sum_{k'=1}^{l}\frac{\partial^3 f}{\partial t\partial x_{1,k}\partial x_{1,k'}}\tau_{j'}h\left(hf_k\sum_{j''=1}^{s}w_{1,j',j''}\right)\left(hf_{k'}\sum_{j''=1}^{s}w_{1,j',j''}\right)$$

$$+\frac{1}{6}\sum_{N=2}^{n}\sum_{k=1}^{l}\sum_{N'=2}^{n}\sum_{k'=1}^{l}\sum_{N''=2}^{n}\sum_{k''=1}^{l}\frac{\partial^3 f}{\partial x_{N,k}\partial x_{N',k'}\partial x_{N'',k''}}\left(hy_{n,0,N-1,k}\xi_{j',N}\right)\left(hy_{n,0,N'-1,k'}\xi_{N',j'}\right)\left(hy_{n,0,N''-1,k''}\xi_{N'',j'}\right)$$

$$+\frac{1}{2}\sum_{N=2}^{n}\sum_{k=1}^{l}\sum_{N'=2}^{n}\sum_{k'=1}^{l}\sum_{k''=1}^{l}\frac{\partial^3 f}{\partial x_{N,k}\partial x_{N',k'}\partial x_{1,k''}}\left(hy_{n,0,N-1,k}\xi_{j',N}\right)\left(hy_{n,0,N'-1,k'}\xi_{N',j'}\right)\left(hf_{k''}\sum_{j''=1}^{s}w_{1,j',j''}\right)$$

$$+\frac{1}{2}\sum_{N=2}^{n}\sum_{k=1}^{l}\sum_{k'=1}^{l}\sum_{k''=1}^{l}\frac{\partial^3 f}{\partial x_{N,k}\partial x_{1,k'}\partial x_{1,k''}}\left(hy_{n,0,N-1,k}\xi_{j',N}\right)\left(hf_{k'}\sum_{j''=1}^{s}w_{1,j',j''}\right)\left(hf_{k''}\sum_{j''=1}^{s}w_{1,j',j''}\right)$$

$$+\frac{1}{6}\sum_{k=1}^{l}\sum_{k'=1}^{l}\sum_{k''=1}^{l}\frac{\partial^3 f}{\partial x_{1,k}\partial x_{1,k'}\partial x_{1,k''}}\left(hf_{k}\sum_{j''=1}^{s}w_{1,j',j''}\right)\left(hf_{k'}\sum_{j''=1}^{s}w_{1,j',j''}\right)\left(hf_{k''}\sum_{j''=1}^{s}w_{1,j',j''}\right)$$

$$\Bigg]$$

$$-\xi_{N,j}^{N}f-\xi_{N,j}^{N+1}h\frac{1}{1+N}\frac{\mathrm{d}f\circ\overline{\mathcal{J}}^{n-1}\hat{\mathrm{y}}}{\mathrm{d}t}-\xi_{N,j}^{N+2}h^2\frac{N!}{(2+N)!}\frac{\mathrm{d}^2 f\circ\overline{\mathcal{J}}^{n-1}\hat{\mathrm{y}}}{\mathrm{d}t^2}-\xi_{N,j}^{N+3}h^3\frac{N!}{(3+N)!}\frac{\mathrm{d}^3 f\circ\overline{\mathcal{J}}^{n-1}\hat{\mathrm{y}}}{\mathrm{d}t^3}\underset{h\to 0}{=}\mathcal{O}(h^4)$$

$$\Leftrightarrow\sum_{j'=1}^{s}w_{N,j,j'}\Bigg[$$

$$f$$

$$+h\frac{\partial f}{\partial t}\tau_{j'}$$

$$+h\sum_{N=2}^{n}\sum_{k=1}^{l}\frac{\partial f}{\partial x_{N,k}}y_{n,0,N-1,k}\xi_{j',N}$$

$$+h\sum_{k=1}^{l}\frac{\partial f}{\partial x_{1,k}}f_k\sum_{j''=1}^{s}w_{1,j',j''}$$

$$+h^2\sum_{N=3}^{n}\sum_{k=1}^{l}\frac{\partial f}{\partial x_{N,k}}y_{n,0,N-2,k}\frac{\xi_{j',N}^2}{2}$$

$$+h^2\sum_{k=1}^{l}\frac{\partial f}{\partial x_{2,k}}f_k\frac{1}{2}\sum_{j''=1}^{s}w_{2,j',j''}$$

$$+h^2\sum_{k=1}^{l}\frac{\partial f}{\partial x_{1,k}}\frac{\partial f_k}{\partial t}\sum_{j''=1}^{s}w_{1,j',j''}\tau_{j''}$$

$$+h^2\sum_{k=1}^{l}\frac{\partial f}{\partial x_{1,k}}\sum_{N=2}^{n}\sum_{k'=1}^{l}\frac{\partial f_k}{\partial x_{N,k'}}y_{n,0,N-1,k'}\sum_{j''=1}^{s}w_{1,j',j''}\xi_{j'',N}$$

$$+h^2\sum_{k=1}^{l}\frac{\partial f}{\partial x_{1,k}}\sum_{k'=1}^{l}\frac{\partial f_k}{\partial x_{1,k'}}f_{k'}\sum_{j''=1}^{s}w_{1,j',j''}\sum_{j'''=1}^{s}w_{1,j'',j'''}$$

$$+h^2\frac{\partial^2 f}{\partial t^2}\frac{1}{2}\tau_{j'}^2$$

$$+h^2\sum_{N=2}^{n}\sum_{k=1}^{l}\frac{\partial^2 f}{\partial t\partial x_{N,k}}y_{n,0,N-1,k}\tau_{j'}\xi_{j',N}$$

$$+h^2\sum_{k=1}^{l}\frac{\partial^2 f}{\partial t\partial x_{1,k}}f_k\tau_{j'}\sum_{j''=1}^{s}w_{1,j',j''}$$

$$
\begin{aligned}
&+h^2 \sum_{N=2}^{n} \sum_{k=1}^{l} \sum_{N'=2}^{n} \sum_{k'=1}^{l} \frac{\partial^2 f}{\partial x_{N,k} \partial x_{N',k'}} y_{n,0,N-1,k} y_{n,0,N'-1,k'} \frac{1}{2} \xi_{j',N} \xi_{N',j'} \\
&+h^2 \sum_{N=2}^{n} \sum_{k=1}^{l} \sum_{k'=1}^{l} \frac{\partial^2 f}{\partial x_{N,k} \partial x_{1,k'}} y_{n,0,N-1,k} f_{k'} \xi_{j',N} \sum_{j''=1}^{s} w_{1,j',j''} \\
&+h^2 \sum_{k=1}^{l} \sum_{k'=1}^{l} \frac{\partial^2 f}{\partial x_{1,k} \partial x_{1,k'}} f_k f_{k'} \frac{1}{2} \left( \sum_{j''=1}^{s} w_{1,j',j''} \right) \sum_{j''=1}^{s} w_{1,j',j''} \\
&+h^3 \sum_{N=4}^{n} \sum_{k=1}^{l} \frac{\partial f}{\partial x_{N,k}} y_{n,0,N-3,k} \frac{\xi_{j',N}^3}{6} \\
&+h^3 \sum_{k=1}^{l} \frac{\partial f}{\partial x_{3,k}} f_k \frac{1}{6} \sum_{j''=1}^{s} w_{3,j',j''} \\
&+h^3 \sum_{k=1}^{l} \frac{\partial f}{\partial x_{2,k}} \frac{\partial f_k}{\partial t} \frac{1}{2} \sum_{j''=1}^{s} w_{2,j',j''} \tau_{j''} \\
&+h^3 \sum_{k=1}^{l} \frac{\partial f}{\partial x_{2,k}} \sum_{N=2}^{n} \sum_{k'=1}^{l} \frac{\partial f_k}{\partial x_{N,k'}} y_{n,0,N-1,k'} \frac{1}{2} \sum_{j''=1}^{s} w_{2,j',j''} \xi_{j'',N} \\
&+h^3 \sum_{k=1}^{l} \frac{\partial f}{\partial x_{2,k}} \sum_{k'=1}^{l} \frac{\partial f_k}{\partial x_{1,k'}} f_{k'} \frac{1}{2} \sum_{j''=1}^{s} w_{2,j',j''} \sum_{j'''=1}^{s} w_{1,j'',j'''} \\
&+h^3 \sum_{k=1}^{l} \frac{\partial f}{\partial x_{1,k}} \sum_{N=3}^{n} \sum_{k'=1}^{l} \frac{\partial f_k}{\partial x_{N,k'}} y_{n,0,N-2,k'} \sum_{j''=1}^{s} w_{1,j',j''} \frac{\xi_{j'',N}^2}{2} \\
&+h^3 \sum_{k=1}^{l} \frac{\partial f}{\partial x_{1,k}} \sum_{k'=1}^{l} \frac{\partial f_k}{\partial x_{2,k'}} f_{k'} \frac{1}{2} \sum_{j''=1}^{s} w_{1,j',j''} \sum_{j'''=1}^{s} w_{2,j'',j'''} \\
&+h^3 \sum_{k=1}^{l} \frac{\partial f}{\partial x_{1,k}} \sum_{k'=1}^{l} \frac{\partial f_k}{\partial x_{1,k'}} \frac{\partial f_{k'}}{\partial t} \sum_{j''=1}^{s} w_{1,j',j''} \sum_{j'''=1}^{s} w_{1,j'',j'''} \tau_{j'''} \\
&+h^3 \sum_{k=1}^{l} \frac{\partial f}{\partial x_{1,k}} \sum_{k'=1}^{l} \frac{\partial f_k}{\partial x_{1,k'}} \sum_{N=2}^{n} \sum_{k''=1}^{l} \frac{\partial f_{k'}}{\partial x_{N,k''}} y_{n,0,N-1,k''} \sum_{j''=1}^{s} w_{1,j',j''} \sum_{j'''=1}^{s} w_{1,j'',j'''} \xi_{N,j'''} \\
&+h^3 \sum_{k=1}^{l} \frac{\partial f}{\partial x_{1,k}} \sum_{k'=1}^{l} \frac{\partial f_k}{\partial x_{1,k'}} \sum_{k''=1}^{l} \frac{\partial f_{k'}}{\partial x_{1,k''}} f_{k''} \sum_{j''=1}^{s} w_{1,j',j''} \sum_{j'''=1}^{s} w_{1,j'',j'''} \sum_{j''''=1}^{s} w_{1,j''',j''''} \\
&+h^3 \sum_{k=1}^{l} \frac{\partial f}{\partial x_{1,k}} \frac{\partial^2 f_k}{\partial t^2} \sum_{j''=1}^{s} w_{1,j',j''} \frac{\tau_{j''}^2}{2} \\
&+h^3 \sum_{k=1}^{l} \frac{\partial f}{\partial x_{1,k}} \sum_{N=2}^{n} \sum_{k'=1}^{l} \frac{\partial^2 f_k}{\partial t \partial x_{N,k'}} y_{n,0,N-1,k'} \sum_{j''=1}^{s} w_{1,j',j''} \tau_{j''} \xi_{j'',N} \\
&+h^3 \sum_{k=1}^{l} \frac{\partial f}{\partial x_{1,k}} \sum_{k'=1}^{l} \frac{\partial^2 f_k}{\partial t \partial x_{1,k'}} f_{k'} \sum_{j''=1}^{s} w_{1,j',j''} \tau_{j''} \sum_{j'''=1}^{s} w_{1,j'',j'''} \\
&+h^3 \sum_{k=1}^{l} \frac{\partial f}{\partial x_{1,k}} \sum_{N=2}^{n} \sum_{k'=1}^{l} \sum_{N'=2}^{n} \sum_{k''=1}^{l} \frac{\partial^2 f_k}{\partial x_{N,k'} \partial x_{N',k''}} y_{n,0,N-1,k'} y_{n,0,N'-1,k''} \frac{1}{2} \sum_{j''=1}^{s} w_{1,j',j''} \xi_{j'',N} \xi_{N',j''}
\end{aligned}
$$

$$+h^3 \sum_{k=1}^{l} \frac{\partial f}{\partial x_{1,k}} \sum_{N=2}^{n} \sum_{k'=1}^{l} \sum_{k''=1}^{l} \frac{\partial^2 f_k}{\partial x_{N,k'} \partial x_{1,k''}} y_{n,0,N-1,k'} f_{k''} \sum_{j''=1}^{s} w_{1,j',j''} \xi_{j'',N} \sum_{j'''=1}^{s} w_{1,j'',j'''}$$

$$+h^3 \sum_{k=1}^{l} \frac{\partial f}{\partial x_{1,k}} \sum_{k'=1}^{l} \sum_{k''=1}^{l} \frac{\partial^2 f_k}{\partial x_{1,k'} \partial x_{1,k''}} f_{k'} f_{k''} \frac{1}{2} \sum_{j''=1}^{s} w_{1,j',j''} \left( \sum_{j'''=1}^{s} w_{1,j'',j'''} \right)^2$$

$$+h^3 \sum_{N=3}^{n} \sum_{k=1}^{l} \frac{\partial^2 f}{\partial t \partial x_{N,k}} y_{n,0,N-2,k} \tau_{j'} \frac{\xi_{j',N}^2}{2}$$

$$+h^3 \sum_{k=1}^{l} \frac{\partial^2 f}{\partial t \partial x_{2,k}} f_k \frac{1}{2} \tau_{j'} \sum_{j''=1}^{s} w_{2,j',j''}$$

$$+h^3 \sum_{k=1}^{l} \frac{\partial^2 f}{\partial t \partial x_{1,k}} \frac{\partial f_k}{\partial t} \tau_{j'} \sum_{j''=1}^{s} w_{1,j',j''} \tau_{j''}$$

$$+h^3 \sum_{k=1}^{l} \frac{\partial^2 f}{\partial t \partial x_{1,k}} \sum_{N=2}^{n} \sum_{k'=1}^{l} \frac{\partial f_k}{\partial x_{N,k'}} y_{n,0,N-1,k'} \tau_{j'} \sum_{j''=1}^{s} w_{1,j',j''} \xi_{j'',N}$$

$$+h^3 \sum_{k=1}^{l} \frac{\partial^2 f}{\partial t \partial x_{1,k}} \sum_{k'=1}^{l} \frac{\partial f_k}{\partial x_{1,k'}} f_{k'} \tau_{j'} \sum_{j''=1}^{s} w_{1,j',j''} \sum_{j'''=1}^{s} w_{1,j'',j'''}$$

$$+h^3 \sum_{N=3}^{n} \sum_{k=1}^{l} \sum_{N'=2}^{n} \sum_{k'=1}^{l} \frac{\partial^2 f}{\partial x_{N,k} \partial x_{N',k'}} y_{n,0,N-2,k} y_{n,0,N'-1,k'} \frac{1}{2} \xi_{j',N}^2 \xi_{N',j'}$$

$$+h^3 \sum_{N=2}^{n} \sum_{k=1}^{l} \sum_{k'=1}^{l} \frac{\partial^2 f}{\partial x_{N,k} \partial x_{2,k'}} y_{n,0,N-1,k} f_{k'} \frac{1}{2} \xi_{j',N} \sum_{j''=1}^{s} w_{2,j',j''}$$

$$+h^3 \sum_{N=2}^{n} \sum_{k=1}^{l} \sum_{k'=1}^{l} \frac{\partial^2 f}{\partial x_{N,k} \partial x_{1,k'}} y_{n,0,N-1,k} \frac{\partial f_{k'}}{\partial t} \xi_{j',N} \sum_{j''=1}^{s} w_{1,j',j''} \tau_{j''}$$

$$+h^3 \sum_{N=2}^{n} \sum_{k=1}^{l} \sum_{k'=1}^{l} \frac{\partial^2 f}{\partial x_{N,k} \partial x_{1,k'}} y_{n,0,N-1,k} \sum_{N'=2}^{n} \sum_{k''=1}^{l} \frac{\partial f_{k'}}{\partial x_{N',k''}} y_{n,0,N'-1,k''} \xi_{j',N} \sum_{j''=1}^{s} w_{1,j',j''} \xi_{N',j''}$$

$$+h^3 \sum_{N=2}^{n} \sum_{k=1}^{l} \sum_{k'=1}^{l} \frac{\partial^2 f}{\partial x_{N,k} \partial x_{1,k'}} y_{n,0,N-1,k} \sum_{k''=1}^{l} \frac{\partial f_{k'}}{\partial x_{1,k''}} f_{k''} \xi_{j',N} \sum_{j''=1}^{s} w_{1,j',j''} \sum_{j'''=1}^{s} w_{1,j'',j'''}$$

$$+h^3 \sum_{N=3}^{n} \sum_{k=1}^{l} \sum_{k'=1}^{l} \frac{\partial^2 f}{\partial x_{N,k} \partial x_{1,k'}} y_{n,0,N-2,k} f_{k'} \frac{\xi_{j',N}^2}{2} \sum_{j''=1}^{s} w_{1,j',j''}$$

$$+h^3 \sum_{k=1}^{l} \sum_{k'=1}^{l} \frac{\partial^2 f}{\partial x_{1,k} \partial x_{2,k'}} f_k f_{k'} \frac{1}{2} \left( \sum_{j''=1}^{s} w_{1,j',j''} \right) \left( \sum_{j''=1}^{s} w_{2,j',j''} \right)$$

$$+h^3 \sum_{k=1}^{l} \sum_{k'=1}^{l} \frac{\partial^2 f}{\partial x_{1,k} \partial x_{1,k'}} \frac{\partial f_k}{\partial t} f_{k'} \left( \sum_{j''=1}^{s} w_{1,j',j''} \right) \sum_{j''=1}^{s} w_{1,j',j''} \tau_{j''}$$

$$+h^3 \sum_{k=1}^{l} \sum_{k'=1}^{l} \frac{\partial^2 f}{\partial x_{1,k} \partial x_{1,k'}} f_k \sum_{N=2}^{n} \sum_{k''=1}^{l} \frac{\partial f_{k'}}{\partial x_{N,k''}} y_{n,0,N-1,k''} \left( \sum_{j''=1}^{s} w_{1,j',j''} \right) \sum_{j''=1}^{s} w_{1,j',j''} \xi_{j'',N}$$

$$+h^3 \sum_{k=1}^{l}\sum_{k'=1}^{l} \frac{\partial^2 f}{\partial x_{1,k}\partial x_{1,k'}} f_k \sum_{k''=1}^{l} \frac{\partial f_{k'}}{\partial x_{1,k''}} f_{k''} \left(\sum_{j''=1}^{s} w_{1,j',j''}\right) \sum_{j''=1}^{s} w_{1,j',j''} \sum_{j'''=1}^{s} w_{1,j'',j'''}$$

$$+h^3 \frac{\partial^3 f}{\partial t^3} \frac{1}{6} \tau_{j'}^3$$

$$+h^3 \sum_{N=2}^{n}\sum_{k=1}^{l} \frac{\partial^3 f}{\partial t^2 \partial x_{N,k}} y_{n,0,N-1,k} \frac{1}{2} \tau_{j'}^2 \xi_{j',N}$$

$$+h^3 \sum_{k=1}^{l} \frac{\partial^3 f}{\partial t^2 \partial x_{1,k}} f_k \frac{1}{2} \tau_{j'}^2 \sum_{j''=1}^{s} w_{1,j',j''}$$

$$+h^3 \sum_{N=2}^{n}\sum_{k=1}^{l}\sum_{N'=2}^{n}\sum_{k'=1}^{l} \frac{\partial^3 f}{\partial t \partial x_{N,k} \partial x_{N',k'}} y_{n,0,N-1,k} y_{n,0,N'-1,k'} \frac{1}{2} \tau_{j'} \xi_{j',N} \xi_{N',j'}$$

$$+h^3 \sum_{N=2}^{n}\sum_{k=1}^{l}\sum_{k'=1}^{l} \frac{\partial^3 f}{\partial t \partial x_{N,k} \partial x_{1,k'}} y_{n,0,N-1,k} f_{k'} \tau_{j'} \xi_{j',N} \sum_{j''=1}^{s} w_{1,j',j''}$$

$$+h^3 \sum_{k=1}^{l}\sum_{k'=1}^{l} \frac{\partial^3 f}{\partial t \partial x_{1,k} \partial x_{1,k'}} f_k f_{k'} \frac{1}{2} \tau_{j'} \left(\sum_{j''=1}^{s} w_{1,j',j''}\right)^2$$

$$+h^3 \sum_{N=2}^{n}\sum_{k=1}^{l}\sum_{N'=2}^{n}\sum_{k'=1}^{l}\sum_{N''=2}^{n}\sum_{k''=1}^{l} \frac{\partial^3 f}{\partial x_{N,k} \partial x_{N',k'} \partial x_{N'',k''}} y_{n,0,N-1,k} y_{n,0,N'-1,k'} y_{n,0,N''-1,k''} \frac{1}{6} \xi_{j',N} \xi_{N',j'} \xi_{N'',j'}$$

$$+h^3 \sum_{N=2}^{n}\sum_{k=1}^{l}\sum_{N'=2}^{n}\sum_{k'=1}^{l}\sum_{k''=1}^{l} \frac{\partial^3 f}{\partial x_{N,k} \partial x_{N',k'} \partial x_{1,k''}} y_{n,0,N-1,k} y_{n,0,N'-1,k'} f_{k''} \frac{1}{2} \xi_{j',N} \xi_{N',j'} \sum_{j''=1}^{s} w_{1,j',j''}$$

$$+h^3 \sum_{N=2}^{n}\sum_{k=1}^{l}\sum_{k'=1}^{l}\sum_{k''=1}^{l} \frac{\partial^3 f}{\partial x_{N,k} \partial x_{1,k'} \partial x_{1,k''}} y_{n,0,N-1,k} f_{k'} f_{k''} \frac{1}{2} \xi_{j',N} \left(\sum_{j''=1}^{s} w_{1,j',j''}\right)^2$$

$$+h^3 \sum_{k=1}^{l}\sum_{k'=1}^{l}\sum_{k''=1}^{l} \frac{\partial^3 f}{\partial x_{1,k} \partial x_{1,k'} \partial x_{1,k''}} f_k f_{k'} f_{k''} \frac{1}{6} \left(\sum_{j''=1}^{s} w_{1,j',j''}\right)^3$$

$$\Bigg]$$

$$-\xi_{N,j}^{N} f$$

$$-\xi_{N,j}^{N+1} h \frac{1}{1+N} \Bigg[$$

$$\frac{\partial f}{\partial t}$$

$$+\sum_{N=2}^{n}\sum_{k=1}^{l} \frac{\partial f}{\partial x_{N,k}} \frac{\mathrm{d}^{n-N+1} \hat{\mathrm{y}}_k}{\mathrm{d}t^{n-N+1}}$$

$$+\sum_{k=1}^{l} \frac{\partial f}{\partial x_{1,k}} f_k \Bigg]$$

$$-\xi_{N,j}^{N+2} h^2 \frac{N!}{(2+N)!} \Bigg[$$

$$\frac{\partial^2 f}{\partial t^2}$$

$$+2 \sum_{k=1}^{l} \frac{\partial^2 f}{\partial t \partial x_{1,k}} f_k$$

$$+\sum_{k=1}^{l} \frac{\partial f}{\partial x_{2,k}} f_k$$

$$+\sum_{k=1}^{l} \frac{\partial f}{\partial x_{1,k}} \frac{\partial f_k}{\partial t}$$

$$+\sum_{k=1}^{l} \frac{\partial f}{\partial x_{1,k}} \sum_{k'=1}^{l} \frac{\partial f_k}{\partial x_{1,k'}} f_{k'}$$

$$+\sum_{k=1}^{l} \sum_{k'=1}^{l} \frac{\partial^2 f}{\partial x_{1,k} \partial x_{1,k'}} f_{k'} f_k$$

$$+2 \sum_{N=2}^{n} \sum_{k=1}^{l} \frac{\partial^2 f}{\partial t \partial x_{N,k}} \frac{\mathrm{d}^{n-N+1} \hat{\mathrm{y}}_k}{\mathrm{d}t^{n-N+1}}$$

$$+\sum_{N=3}^{n} \sum_{k=1}^{l} \frac{\partial f}{\partial x_{N,k}} \frac{\mathrm{d}^{n-N+2} \hat{\mathrm{y}}_k}{\mathrm{d}t^{n-N+2}}$$

$$+\sum_{k=1}^{l} \frac{\partial f}{\partial x_{1,k}} \sum_{N=2}^{n} \sum_{k'=1}^{l} \frac{\partial f_k}{\partial x_{N,k'}} \frac{\mathrm{d}^{n-N+1} \hat{\mathrm{y}}_{k'}}{\mathrm{d}t^{n-N+1}}$$

$$+\sum_{N=2}^{n} \sum_{k=1}^{l} \sum_{N'=2}^{n} \sum_{k'=1}^{l} \frac{\partial^2 f}{\partial x_{N,k} \partial x_{N',k'}} \frac{\mathrm{d}^{n-N'+1} \hat{\mathrm{y}}_{k'}}{\mathrm{d}t^{n-N'+1}} \frac{\mathrm{d}^{n-N+1} \hat{\mathrm{y}}_k}{\mathrm{d}t^{n-N+1}}$$

$$+2 \sum_{N=2}^{n} \sum_{k=1}^{l} \sum_{k'=1}^{l} \frac{\partial^2 f}{\partial x_{N,k} \partial x_{1,k'}} f_{k'} \frac{\mathrm{d}^{n-N+1} \hat{\mathrm{y}}_k}{\mathrm{d}t^{n-N+1}} \Bigg]$$

$$-\xi_{N,j}^{N+3} h^3 \frac{N!}{(3+N)!} \Bigg[$$

$$\frac{\partial^3 f}{\partial t^3}$$

$$+3 \sum_{k=1}^{l} \frac{\partial^3 f}{\partial t^2 \partial x_{1,k}} f_k$$

$$+3 \sum_{k=1}^{l} \frac{\partial^2 f}{\partial t \partial x_{2,k}} f_k$$

$$+3 \sum_{k=1}^{l} \frac{\partial^2 f}{\partial t \partial x_{1,k}} \frac{\partial f_k}{\partial t}$$

$$+\sum_{k=1}^{l}\frac{\partial f}{\partial x_{3,k}}f_k$$

$$+\sum_{k=1}^{l}\frac{\partial f}{\partial x_{2,k}}\frac{\partial f_k}{\partial t}$$

$$+\sum_{k=1}^{l}\frac{\partial f}{\partial x_{1,k}}\frac{\partial^2 f_k}{\partial t^2}$$

$$+3\sum_{k=1}^{l}\frac{\partial^2 f}{\partial t\partial x_{1,k}}\sum_{k'=1}^{l}\frac{\partial f_k}{\partial x_{1,k'}}f_{k'}$$

$$+\sum_{k=1}^{l}\frac{\partial f}{\partial x_{2,k}}\sum_{k'=1}^{l}\frac{\partial f_k}{\partial x_{1,k'}}f_{k'}$$

$$+2\sum_{k=1}^{l}\frac{\partial f}{\partial x_{1,k}}\sum_{k'=1}^{l}\frac{\partial^2 f_k}{\partial t\partial x_{1,k'}}f_{k'}$$

$$+\sum_{k=1}^{l}\frac{\partial f}{\partial x_{1,k}}\sum_{k'=1}^{l}\frac{\partial f_k}{\partial x_{2,k'}}f_{k'}$$

$$+\sum_{k=1}^{l}\frac{\partial f}{\partial x_{1,k}}\sum_{k'=1}^{l}\frac{\partial f_k}{\partial x_{1,k'}}\frac{\partial f_{k'}}{\partial t}$$

$$+3\sum_{k=1}^{l}\sum_{k'=1}^{l}\frac{\partial^3 f}{\partial t\partial x_{1,k}\partial x_{1,k'}}f_{k'}f_k$$

$$+3\sum_{k=1}^{l}\sum_{k'=1}^{l}\frac{\partial f}{\partial x_{2,k}\partial x_{1,k'}}f_{k'}f_k$$

$$+3\sum_{k=1}^{l}\sum_{k'=1}^{l}\frac{\partial f}{\partial x_{1,k}\partial x_{1,k'}}f_{k'}\frac{\partial f_k}{\partial t}$$

$$+\sum_{k=1}^{l}\frac{\partial f}{\partial x_{1,k}}\sum_{k'=1}^{l}\frac{\partial f_k}{\partial x_{1,k'}}\sum_{k''=1}^{l}\frac{\partial f_{k'}}{\partial x_{1,k''}}f_{k''}$$

$$+\sum_{k=1}^{l}\frac{\partial f}{\partial x_{1,k}}\sum_{k'=1}^{l}\sum_{k''=1}^{l}\frac{\partial^2 f^{[l]}}{\partial x_{1,k'}\partial x_{1,k''}}f_{k''}f_{k'}$$

$$+3\sum_{k=1}^{l}\sum_{k'=1}^{l}\frac{\partial f}{\partial x_{1,k}\partial x_{1,k'}}f_{k'}\sum_{k''=1}^{l}\frac{\partial f_k}{\partial x_{1,k''}}f_{k''}$$

$$+\sum_{k=1}^{l}\sum_{k'=1}^{l}\sum_{k''=1}^{l}\frac{\partial^3 f}{\partial x_{1,k}\partial x_{1,k'}x_{1,k''}}f_{k''}f_{k'}f_k$$

$$+3\sum_{N=2}^{n}\sum_{k=1}^{l}\frac{\partial^3 f}{\partial t^2\partial x_{N,k}}\frac{\mathrm{d}^{n-N+1}\hat{\mathrm{y}}_k}{\mathrm{d}t^{n-N+1}}$$

$$+3\sum_{N=3}^{n}\sum_{k=1}^{l}\frac{\partial^2 f}{\partial t\partial x_{N,k}}\frac{\mathrm{d}^{n-N+2}\hat{\mathrm{y}}_k}{\mathrm{d}t^{n-N+2}}$$

$$+\sum_{N=4}^{n}\sum_{k=1}^{l}\frac{\partial f}{\partial x_{N,k}}\frac{\mathrm{d}^{n-N+3}\hat{y}_k}{\mathrm{d}t^{n-N+3}}$$

$$+3\sum_{k=1}^{l}\frac{\partial^2 f}{\partial t\partial x_{1,k}}\sum_{N=2}^{n}\sum_{k'=1}^{l}\frac{\partial f_k}{\partial x_{N,k'}}\frac{\mathrm{d}^{n-N+1}\hat{y}_{k'}}{\mathrm{d}t^{n-N+1}}$$

$$+\sum_{k=1}^{l}\frac{\partial f}{\partial x_{2,k}}\sum_{N=2}^{n}\sum_{k'=1}^{l}\frac{\partial f_k}{\partial x_{N,k'}}\frac{\mathrm{d}^{n-N+1}\hat{y}_{k'}}{\mathrm{d}t^{n-N+1}}$$

$$+6\sum_{N=2}^{n}\sum_{k=1}^{l}\sum_{k'=1}^{l}\frac{\partial^3 f}{\partial t\partial x_{N,k}\partial x_{1,k'}}f_{k'}\frac{\mathrm{d}^{n-N+1}\hat{y}_k}{\mathrm{d}t^{n-N+1}}$$

$$+2\sum_{k=1}^{l}\frac{\partial f}{\partial x_{1,k}}\sum_{N=2}^{n}\sum_{k'=1}^{l}\frac{\partial^2 f_k}{\partial t\partial x_{N,k'}}\frac{\mathrm{d}^{n-N+1}\hat{y}_{k'}}{\mathrm{d}t^{n-N+1}}$$

$$+\sum_{k=1}^{l}\frac{\partial f}{\partial x_{1,k}}\sum_{N=3}^{n}\sum_{k'=1}^{l}\frac{\partial f_k}{\partial x_{N,k'}}\frac{\mathrm{d}^{n-N+2}\hat{y}_{k'}}{\mathrm{d}t^{n-N+2}}$$

$$+3\sum_{N=3}^{n}\sum_{k=1}^{l}\sum_{k'=1}^{l}\frac{\partial f}{\partial x_{N,k}\partial x_{1,k'}}f_{k'}\frac{\mathrm{d}^{n-N+2}\hat{y}_k}{\mathrm{d}t^{n-N+2}}$$

$$+3\sum_{k=1}^{l}\sum_{N=2}^{n}\sum_{k'=1}^{l}\frac{\partial f}{\partial x_{2,k}\partial x_{N,k'}}\frac{\mathrm{d}^{n-N+1}\hat{y}_{k'}}{\mathrm{d}t^{n-N+1}}f_k$$

$$+3\sum_{k=1}^{l}\sum_{N=2}^{n}\sum_{k'=1}^{l}\frac{\partial f}{\partial x_{1,k}\partial x_{N,k'}}\frac{\mathrm{d}^{n-N+1}\hat{y}_{k'}}{\mathrm{d}t^{n-N+1}}\frac{\partial f_k}{\partial t}$$

$$+2\sum_{k=1}^{l}\frac{\partial f}{\partial x_{1,k}}\sum_{N=2}^{n}\sum_{k'=1}^{l}\sum_{k''=1}^{l}\frac{\partial^2 f_k}{\partial x_{N,k'}\partial x_{1,k''}}f_{k''}\frac{\mathrm{d}^{n-N+1}\hat{y}_{k'}}{\mathrm{d}t^{n-N+1}}$$

$$+\sum_{k=1}^{l}\frac{\partial f}{\partial x_{1,k}}\sum_{k'=1}^{l}\frac{\partial f_k}{\partial x_{1,k'}}\sum_{N=2}^{n}\sum_{k''=1}^{l}\frac{\partial f_{k'}}{\partial x_{N,k''}}\frac{\mathrm{d}^{i+1}\hat{y}_{k''}}{\mathrm{d}t^{i+1}}$$

$$+3\sum_{k=1}^{l}\sum_{N=2}^{n}\sum_{k'=1}^{l}\frac{\partial f}{\partial x_{1,k}\partial x_{N,k'}}\frac{\mathrm{d}^{n-N+1}\hat{y}_{k'}}{\mathrm{d}t^{n-N+1}}\sum_{k''=1}^{l}\frac{\partial f_k}{\partial x_{1,k''}}f_{k''}$$

$$+3\sum_{k=1}^{l}\sum_{k'=1}^{l}\frac{\partial f}{\partial x_{1,k}\partial x_{1,k'}}f_{k'}\sum_{N'=2}^{n}\sum_{k''=1}^{l}\frac{\partial f_k}{\partial x_{N',k''}}\frac{\mathrm{d}^{n-N'+1}\hat{y}_{k''}}{\mathrm{d}t^{n-N'+1}}$$

$$+3\sum_{N=2}^{n}\sum_{k=1}^{l}\sum_{k'=1}^{l}\sum_{k''=1}^{l}\frac{\partial^3 f}{\partial x_{N,k}\partial x_{1,k'}x_{1,k''}}f_{k''}f_{k'}\frac{\mathrm{d}^{n-N+1}\hat{y}_k}{\mathrm{d}t^{n-N+1}}$$

$$+3\sum_{N=2}^{n}\sum_{k=1}^{l}\sum_{N'=2}^{n}\sum_{k'=1}^{l}\frac{\partial^3 f}{\partial t\partial x_{N,k}\partial x_{N',k'}}\frac{\mathrm{d}^{n-N'+1}\hat{y}_{k'}}{\mathrm{d}t^{n-N'+1}}\frac{\mathrm{d}^{n-N+1}\hat{y}_k}{\mathrm{d}t^{n-N+1}}$$

$$+3\sum_{N=3}^{n}\sum_{k=1}^{l}\sum_{N'=2}^{n}\sum_{k'=1}^{l}\frac{\partial f}{\partial x_{N,k}\partial x_{N',k'}}\frac{\mathrm{d}^{n-N'+1}\hat{y}_{k'}}{\mathrm{d}t^{n-N'+1}}\frac{\mathrm{d}^{n-N+2}\hat{y}_k}{\mathrm{d}t^{n-N+2}}$$

$$+\sum_{k=1}^{l}\frac{\partial f}{\partial x_{1,k}}\sum_{N=2}^{n}\sum_{k'=1}^{l}\sum_{N'=2}^{n}\sum_{k''=1}^{l}\frac{\partial^2 f_k}{\partial x_{N,k'}\partial x_{N',k''}}\frac{\mathrm{d}^{n-N'+1}\hat{y}_{k''}}{\mathrm{d}t^{n-N'+1}}\frac{\mathrm{d}^{n-N+1}\hat{y}_{k'}}{\mathrm{d}t^{n-N+1}}$$

$$
\begin{aligned}
&+3\sum_{k=1}^{l}\sum_{N=2}^{n}\sum_{k'=1}^{l}\frac{\partial f}{\partial x_{1,k}\partial x_{N,k'}}\frac{\mathrm{d}^{n-N+1}\hat{\mathrm{y}}_{k'}}{\mathrm{d}t^{n-N+1}}\sum_{N'=2}^{n}\sum_{k''=1}^{l}\frac{\partial f_k}{\partial x_{N',k''}}\frac{\mathrm{d}^{n-N'+1}\hat{\mathrm{y}}_{k''}}{\mathrm{d}t^{n-N'+1}}\\
&+3\sum_{N=2}^{n}\sum_{k=1}^{l}\sum_{N'=2}^{n}\sum_{k'=1}^{l}\sum_{k''=1}^{l}\frac{\partial^3 f}{\partial x_{N,k}\partial x_{N',k'}x_{1,k''}}f_{k''}\frac{\mathrm{d}^{n-N'+1}\hat{\mathrm{y}}_{k'}}{\mathrm{d}t^{n-N'+1}}\frac{\mathrm{d}^{n-N+1}\hat{\mathrm{y}}_{k}}{\mathrm{d}t^{n-N+1}}\\
&+\sum_{N=2}^{n}\sum_{k=1}^{l}\sum_{N'=2}^{n}\sum_{k'=1}^{l}\sum_{N''=2}^{n}\sum_{k''=1}^{l}\frac{\partial^3 f}{\partial x_{N,k}\partial x_{N',k'}\partial x_{N'',k''}}\frac{\mathrm{d}^{n-N''+1}\hat{\mathrm{y}}_{k''}}{\mathrm{d}t^{n-N''+1}}\frac{\mathrm{d}^{n-N'+1}\hat{\mathrm{y}}_{k'}}{\mathrm{d}t^{n-N'+1}}\frac{\mathrm{d}^{n-N+1}\hat{\mathrm{y}}_{k}}{\mathrm{d}t^{n-N+1}}\\
&\Bigg]\underset{h\to 0}{=\mathcal{O}}\left(h^4\right)
\end{aligned}
$$

By equating the coefficients of the powers of $h$ and the partial derivatives, we obtain the order conditions.

# C Some solutions of the order conditions

We here list the solutions of the order conditions in some cases. The cases we consider depend on whether the method is explicit or implicit, on the number of points of the method, the stage we are solving for, and the desired order of consistency. The considered case is specified using the format KIND $s - j - v$, where KIND is either Explicit method or Implicit method, $s$ is the number of points of the method, $j$ is the stage we are solving for, and $v$ is the desired order of consistency. For explicit methods, we assume the method is topologically sorted.

For node-determined, many of the order conditions are redundant, and the only conditions we need to consider are :

- Order 1 : 1
- Order 2 : 1 and 2
- Order 3 : 1 through 6
- Order 4 : 1 through 19

To simplify the already long system of equations we only consider those equations, though we do not assume that the method is node-determined. If needed, the remainder of the equations can be added to solve for multi-order Runge-Kutta methods. Even with this simplification, the computations are long, therefore for the most complicated cases the explanations are only general guidelines.

Since there are many cases to distinguish, each case is identified using a sequence of zero and one. Each number represents a case or branching path. A 0 means that the condition is verified, a 1 means that the conditions is not verified.

### C.1 Explicit method, $s = 1$, $j = 2$, $v = 1$

The condition is $w_{N,2,1} = \xi^N_{2,N}$.

### C.2 Explicit method, $s = 2$, $j = 3$, $v = 2$

The conditions are :

Order of consistency 1 :

$$w_{N,3,1} + w_{N,3,2} = \xi_{3,N}^{N}$$

Order of consistency 2 :

$$w_{N,3,1}\tau_1 + w_{N,3,2}\tau_2 = \frac{\xi_{3,N}^{N+1}}{1+N}$$
$$w_{N,3,2}w_{1,2,1} = \frac{\xi_{3,N}^{N+1}}{1+N}$$

We rewrite this system as :

Order of consistency 1 :

$$w_{N,3,1} = \xi_{3,N}^{N} - w_{N,3,2}$$

Order of consistency 2 :

$$w_{N,3,2}(\tau_2 - \tau_1) = \xi_{3,N}^{N}\left(\frac{\xi_{3,N}}{1+N} - \tau_1\right)$$
$$w_{N,3,2}w_{1,2,1} = \frac{\xi_{3,N}^{N+1}}{1+N}$$

The solutions are :

**0-0**

$\xi_{3,N} = 0, w_{N,3,2} = 0$

$$w_{N,3,1} = 0$$

**0-1**

$\xi_{3,N} = 0, w_{N,3,2} \neq 0$

$$\tau_2 = \tau_1$$
$$w_{1,2,1} = 0$$
$$w_{N,3,1} = -w_{N,3,2}$$

**1-0**

$\xi_{3,N} \neq 0, \tau_1 = \tau_2$

$$
\begin{aligned}
w_{1,2,1} &\neq 0 \\
\tau_1 &= \frac{\xi_{3,N}}{1+N} \\
w_{N,3,1} &= \xi_{3,N}^{N} - \frac{\xi_{3,N}^{N+1}}{w_{1,2,1}(1+N)} \\
w_{N,3,2} &= \frac{\xi_{3,N}^{N+1}}{w_{1,2,1}(1+N)}
\end{aligned}
$$

**1-1**

$\xi_{3,N} \neq 0, \tau_1 \neq \tau_2$

$$
\begin{aligned}
\tau_1 &\neq \frac{\xi_{3,N}}{1+N} \\
w_{1,2,1} &= \xi_{3,N} \frac{\tau_2 - \tau_1}{\xi_{3,N} - \tau_1(1+N)} \\
w_{N,3,1} &= \frac{\xi_{3,N}^{N}}{\tau_2 - \tau_1} \left( \tau_2 - \frac{\xi_{3,N}}{1+N} \right) \\
w_{N,3,2} &= \frac{\xi_{3,N}^{N}}{\tau_2 - \tau_1} \left( \frac{\xi_{3,N}}{1+N} - \tau_1 \right)
\end{aligned}
$$

### C.3 Explicit method, $s = 3$, $j = 4$, $v = 2$

The conditions are :

Order of consistency 1 :

$$w_{N,4,1} + w_{N,4,2} + w_{N,4,3} = \xi_{4,N}^{N}$$

Order of consistency 2 :

$$w_{N,4,1}\tau_1 + w_{N,4,2}\tau_2 + w_{N,4,3}\tau_3 = \frac{\xi_{4,N}^{N+1}}{1+N}$$

$$w_{N,4,2}w_{1,2,1} + w_{N,4,3}(w_{1,3,1} + w_{1,3,2}) = \frac{\xi_{4,N}^{N+1}}{1+N}$$

We rewrite this system as :

Order of consistency 1 :

$$w_{N,4,1} = \xi_{4,N}^{N} - w_{N,4,2} - w_{N,4,3}$$

Order of consistency 2 :

$$w_{N,4,2}(\tau_2 - \tau_1) + w_{N,4,3}(\tau_3 - \tau_1) = \xi_{4,N}^{N}\left(\frac{\xi_{4,N}}{1+N} - \tau_1\right)$$

$$w_{N,4,2}w_{1,2,1} + w_{N,4,3}(w_{1,3,1} + w_{1,3,2}) = \frac{\xi_{4,N}^{N+1}}{1+N}$$

The solutions are :

**0-0-0**

$\xi_{4,N} = 0, w_{N,4,2} = 0, w_{N,4,3} = 0$

$$w_{N,4,1} = 0$$

**0-0-1**

$\xi_{4,N} = 0, w_{N,4,2} = 0, w_{N,4,3} \neq 0$

$$\tau_3 = \tau_1$$
$$w_{1,3,1} = -w_{1,3,2}$$
$$w_{N,4,1} = -w_{N,4,3}$$

**0-1-0**

$\xi_{4,N} = 0, w_{N,4,2} \neq 0, w_{N,4,3} = 0$

$$\tau_2 = \tau_1$$
$$w_{1,2,1} = 0$$
$$w_{N,4,1} = -w_{N,4,2}$$

**0-1-1-0-0**

$\xi_{4,N} = 0, w_{N,4,2} \neq 0, w_{N,4,3} \neq 0, \tau_2 = \tau_1, w_{1,2,1} = 0$

$$\tau_3 = \tau_1$$
$$w_{1,3,1} = -w_{1,3,2}$$
$$w_{N,4,1} = -w_{N,4,2} - w_{N,4,3}$$

**0-1-1-0-1**

$\xi_{4,N} = 0, w_{N,4,2} \neq 0, w_{N,4,3} \neq 0, \tau_2 = \tau_1, w_{1,2,1} \neq 0$

$$\tau_3 = \tau_1$$
$$w_{N,4,1} = w_{N,4,3} \frac{w_{1,3,1} + w_{1,3,2} - w_{1,2,1}}{w_{1,2,1}}$$
$$w_{N,4,2} = -w_{N,4,3} \frac{w_{1,3,1} + w_{1,3,2}}{w_{1,2,1}}$$

**0-1-1-1**

$\xi_{4,N} = 0, w_{N,4,2} \neq 0, w_{N,4,3} \neq 0, \tau_2 \neq \tau_1$

$$w_{1,3,1} = \frac{\tau_3 - \tau_1}{\tau_2 - \tau_1} w_{1,2,1} - w_{1,3,2}$$
$$w_{N,4,1} = w_{N,4,3} \frac{\tau_3 - \tau_2}{\tau_2 - \tau_1}$$
$$w_{N,4,2} = -w_{N,4,3} \frac{\tau_3 - \tau_1}{\tau_2 - \tau_1}$$

**1-0-0-0**

$\xi_{4,N} \neq 0, \tau_2 = \tau_1, \tau_3 = \tau_1, w_{1,2,1} = 0$

$$w_{1,3,1} + w_{1,3,2} \neq 0$$
$$\tau_1 = \frac{\xi_{4,N}}{1+N}$$
$$w_{N,4,1} = \xi_{4,N}^N - w_{N,4,2} - \frac{\xi_{4,N}^{N+1}}{(1+N)\big(w_{1,3,1} + w_{1,3,2}\big)}$$
$$w_{N,4,3} = \frac{\xi_{4,N}^{N+1}}{(1+N)\big(w_{1,3,1} + w_{1,3,2}\big)}$$

**1-0-0-1**

$\xi_{4,N} \neq 0, \tau_2 = \tau_1, \tau_3 = \tau_1, w_{1,2,1} \neq 0$

$$\tau_1 = \frac{\xi_{4,N}}{1+N}$$

$$w_{N,4,1} = \xi_{4,N}^{N} - \frac{\xi_{4,N}^{N+1}}{(1+N)w_{1,2,1}} + w_{N,4,3}\frac{w_{1,3,1} + w_{1,3,2} - w_{1,2,1}}{w_{1,2,1}}$$

$$w_{N,4,2} = \frac{\xi_{4,N}^{N+1}}{(1+N)w_{1,2,1}} - w_{N,4,3}\frac{w_{1,3,1} + w_{1,3,2}}{w_{1,2,1}}$$

**1-0-1-0**

$\xi_{4,N} \neq 0, \tau_2 = \tau_1, \tau_3 \neq \tau_1, w_{1,2,1} = 0$

$$\tau_1 \neq \frac{\xi_{4,N}}{1+N}$$

$$w_{1,3,1} = \xi_{4,N}\frac{\tau_3 - \tau_1}{\xi_{4,N} - \tau_1(1+N)} - w_{1,3,2}$$

$$w_{N,4,1} = \frac{\xi_{4,N}^{N}}{\tau_3 - \tau_1}\left(\tau_3 - \frac{\xi_{4,N}}{1+N}\right) - w_{N,4,2}$$

$$w_{N,4,3} = \frac{\xi_{4,N}^{N}}{\tau_3 - \tau_1}\left(\frac{\xi_{4,N}}{1+N} - \tau_1\right)$$

**1-0-1-1**

$\xi_{4,N} \neq 0, \tau_2 = \tau_1, \tau_3 \neq \tau_1, w_{1,2,1} \neq 0$

$$w_{N,4,1} = \xi_{4,N}^{N} - \frac{1}{w_{1,2,1}}\left(\frac{\xi_{4,N}^{N+1}}{1+N} - \frac{\xi_{4,N}^{N}}{\tau_3 - \tau_1}\left(\frac{\xi_{4,N}}{1+N} - \tau_1\right)(w_{1,3,1} + w_{1,3,2})\right) - \frac{\xi_{4,N}^{N}}{(\tau_3 - \tau_1)}\left(\frac{\xi_{4,N}}{1+N} - \tau_1\right)$$

$$w_{N,4,2} = \frac{1}{w_{1,2,1}}\left(\frac{\xi_{4,N}^{N+1}}{1+N} - \frac{\xi_{4,N}^{N}}{\tau_3 - \tau_1}\left(\frac{\xi_{4,N}}{1+N} - \tau_1\right)(w_{1,3,1} + w_{1,3,2})\right)$$

$$w_{N,4,3} = \frac{\xi_{4,N}^{N}}{(\tau_3 - \tau_1)}\left(\frac{\xi_{4,N}}{1+N} - \tau_1\right)$$

**1-1-0**

$\xi_{4,N} \neq 0, \tau_2 \neq \tau_1, w_{1,3,1} = \frac{\tau_3-\tau_1}{\tau_2-\tau_1} w_{1,2,1} - w_{1,3,2}$

$$\tau_1 \neq \frac{\xi_{4,N}}{1+N}$$

$$w_{1,2,1} = \xi_{4,N} \frac{\tau_2 - \tau_1}{\xi_{4,N} - \tau_1(1+N)}$$

$$w_{N,4,1} = \frac{\xi_{4,N}^N}{\tau_2 - \tau_1}\left(\tau_2 - \frac{\xi_{4,N}}{1+N}\right) + w_{N,4,3}\frac{\tau_3 - \tau_2}{\tau_2 - \tau_1}$$

$$w_{N,4,2} = \frac{\xi_{4,N}^N}{\tau_2 - \tau_1}\left(\frac{\xi_{4,N}}{1+N} - \tau_1\right) - w_{N,4,3}\frac{\tau_3 - \tau_1}{\tau_2 - \tau_1}$$

**1-1-1**

$\xi_{4,N} \neq 0, \tau_2 \neq \tau_1, w_{1,3,1} \neq \frac{\tau_3-\tau_1}{\tau_2-\tau_1} w_{1,2,1} - w_{1,3,2}$

$$a = \big((\tau_2 - \tau_1)\big(w_{1,3,1} + w_{1,3,2}\big) - (\tau_3 - \tau_1)w_{1,2,1}\big) \neq 0$$

$$w_{N,4,1} = \xi_{4,N}^N - \frac{\xi_{4,N}^N}{\tau_2 - \tau_1}\left(\frac{\xi_{4,N}}{1+N} - \tau_1\right) + \frac{\xi_{4,N}^N}{a}\left(\frac{\xi_{4,N}}{1+N}(\tau_2 - \tau_1) - \left(\frac{\xi_{4,N}}{1+N} - \tau_1\right)w_{1,2,1}\right)\frac{\tau_3 - \tau_2}{\tau_2 - \tau_1}$$

$$w_{N,4,2} = \frac{\xi_{4,N}^N}{\tau_2 - \tau_1}\left(\frac{\xi_{4,N}}{1+N} - \tau_1\right) - \frac{\xi_{4,N}^N}{a}\left(\frac{\xi_{4,N}}{1+N}(\tau_2 - \tau_1) - \left(\frac{\xi_{4,N}}{1+N} - \tau_1\right)w_{1,2,1}\right)\frac{\tau_3 - \tau_1}{\tau_2 - \tau_1}$$

$$w_{N,4,3} = \frac{\xi_{4,N}^N}{a}\left(\frac{\xi_{4,N}}{1+N}(\tau_2 - \tau_1) - \left(\frac{\xi_{4,N}}{1+N} - \tau_1\right)w_{1,2,1}\right)$$

### C.4 Explicit method, $s = 3, j = 4, v = 3$

The conditions are :

Order of consistency 1 :

$$w_{N,4,1} + w_{N,4,2} + w_{N,4,3} = \xi_{4,N}^{N}$$

Order of consistency 2 :

$$w_{N,4,1}\tau_1 + w_{N,4,2}\tau_2 + w_{N,4,3}\tau_3 = \frac{\xi_{4,N}^{N+1}}{1+N}$$

$$w_{N,4,2}w_{1,2,1} + w_{N,4,3}(w_{1,3,1} + w_{1,3,2}) = \frac{\xi_{4,N}^{N+1}}{1+N}$$

Order of consistency 3 :

$$w_{N,4,1}\tau_1^2 + w_{N,4,2}\tau_2^2 + w_{N,4,3}\tau_3^2 = \xi_{4,N}^{N+2} 2\frac{N!}{(2+N)!}$$

$$w_{N,4,2}\tau_2 w_{1,2,1} + w_{N,4,3}\tau_3(w_{1,3,1} + w_{1,3,2}) = \xi_{4,N}^{N+2} 2\frac{N!}{(2+N)!}$$

$$w_{N,4,2}w_{1,2,1}\tau_1 + w_{N,4,3}(w_{1,3,1}\tau_1 + w_{1,3,2}\tau_2) = \xi_{4,N}^{N+2} \frac{N!}{(2+N)!}$$

$$w_{N,4,3}w_{1,3,2}w_{1,2,1} = \xi_{4,N}^{N+2} \frac{N!}{(2+N)!}$$

$$w_{N,4,2}w_{1,2,1}^2 + w_{N,4,3}(w_{1,3,1} + w_{1,3,2})^2 = \xi_{4,N}^{N+2} 2\frac{N!}{(2+N)!}$$

$$n > 1,\; w_{N,4,2}w_{2,2,1} + w_{N,4,3}(w_{2,3,1} + w_{2,3,2}) = \xi_{4,N}^{N+2} 2\frac{N!}{(2+N)!}$$

We now list the main operations performed to rewrite this system.

1) isolate $w_{N,4,1}\tau_1$ in 2.1
2) 2.1 in 3.1
3) isolate $w_{N,4,2}w_{1,2,1}$ in 2.2
4) 2.2 in 3.2, 3.3, 3.5
5) isolate $w_{N,4,1}$ in 1.1
6) 1.1 in 2.1
7) isolate $w_{N,4,2}(\tau_2 - \tau_1)$ in 2.1
8) 2.1 in 3.1

Which yields :

Order of consistency 1 :

$$w_{N,4,1} = \xi_{4,N}^{N} - w_{N,4,2} - w_{N,4,3}$$

Order of consistency 2 :

$$w_{N,4,2}(\tau_2-\tau_1) = \xi_{4,N}^{N}\left(\frac{\xi_{4,N}}{1+N}-\tau_1\right) - w_{N,4,3}(\tau_3-\tau_1)$$
$$w_{N,4,2}w_{1,2,1} = \frac{\xi_{4,N}^{N+1}}{1+N} - w_{N,4,3}(w_{1,3,1}+w_{1,3,2})$$

Order of consistency 3 :

$$w_{N,4,3}(\tau_3-\tau_2)(\tau_3-\tau_1) = \xi_{4,N}^{N}\left(\frac{\xi_{4,N}}{1+N}\left(\xi_{4,N}\frac{2}{2+N}-\tau_1\right) - \left(\frac{\xi_{4,N}}{1+N}-\tau_1\right)\tau_2\right)$$
$$w_{N,4,3}(\tau_3-\tau_2)(w_{1,3,1}+w_{1,3,2}) = \frac{\xi_{4,N}^{N+1}}{1+N}\left(\xi_{4,N}\frac{2}{2+N}-\tau_2\right)$$
$$w_{N,4,3}w_{1,3,2}(\tau_2-\tau_1) = \frac{\xi_{4,N}^{N+1}}{1+N}\left(\frac{\xi_{4,N}}{2+N}-\tau_1\right)$$
$$w_{N,4,3}w_{1,3,2}w_{1,2,1} = \xi_{4,N}^{N+2}\frac{N!}{(2+N)!}$$
$$w_{N,4,3}(w_{1,3,1}+w_{1,3,2})(w_{1,3,1}+w_{1,3,2}-w_{1,2,1}) = \frac{\xi_{4,N}^{N+1}}{1+N}\left(\xi_{4,N}\frac{2}{2+N}-w_{1,2,1}\right)$$
$$n>1,\ w_{N,4,2}w_{2,2,1}+w_{N,4,3}(w_{2,3,1}+w_{2,3,2}) = \xi_{4,N}^{N+2}2\frac{N!}{(2+N)!}$$

The solutions are :

**0-0-0-0-0**

$\xi_{4,N}=0, \tau_3=\tau_2, \tau_2=\tau_1, w_{N,4,3}=0, w_{N,4,2}=0$

$$w_{N,4,1}=0$$

**0-0-0-0-1**

$\xi_{4,N}=0, \tau_3=\tau_2, \tau_2=\tau_1, w_{N,4,3}=0, w_{N,4,2}\neq 0$

$$w_{1,2,1}=0$$
$$n>1,\ w_{2,2,1}=0$$
$$w_{N,4,1}=-w_{N,4,2}$$

**0-0-0-1-0-0**

$\xi_{4,N}=0, \tau_3=\tau_2, \tau_2=\tau_1, w_{N,4,3}\neq 0, w_{1,3,1}+w_{1,3,2}=0, w_{1,2,1}=0$

$$n>1,\ w_{2,3,1}=-\frac{w_{N,4,2}}{w_{N,4,3}}w_{2,2,1}-w_{2,3,2}$$
$$w_{N,4,1}=-w_{N,4,2}-w_{N,4,3}$$

**0-0-0-1-0-1**

$\xi_{4,N} = 0, \tau_3 = \tau_2, \tau_2 = \tau_1, w_{N,4,3} \neq 0, w_{1,3,1} + w_{1,3,2} = 0, w_{1,2,1} \neq 0$

$$w_{1,3,2} = 0$$

$$n > 1,\ w_{2,3,1} = -\frac{w_{N,4,2}}{w_{N,4,3}} w_{2,2,1} - w_{2,3,2}$$

$$w_{N,4,1} = -w_{N,4,2} - w_{N,4,3}$$

$$w_{N,4,2} = 0$$

**0-0-0-1-1**

$\xi_{4,N} = 0, \tau_3 = \tau_2, \tau_2 = \tau_1, w_{N,4,3} \neq 0, \left(w_{1,3,1} + w_{1,3,2}\right) \neq 0$

$$w_{1,3,1} = w_{1,2,1}$$

$$n > 1,\ w_{2,3,1} = -\frac{w_{N,4,2}}{w_{N,4,3}} w_{2,2,1} - w_{2,3,2}$$

$$w_{N,4,1} = w_{1,3,2} = 0$$

$$w_{N,4,2} = -w_{N,4,3}$$

**0-0-1-0**

$\xi_{4,N} = 0, \tau_3 = \tau_2, \tau_2 \neq \tau_1, w_{N,4,3} = 0$

$$w_{N,4,1} = w_{N,4,2} = 0$$

**0-0-1-1**

$\xi_{4,N} = 0, \tau_3 = \tau_2, \tau_2 \neq \tau_1, w_{N,4,3} \neq 0$

$$w_{1,3,1} = w_{1,2,1}$$

$$n > 1,\ w_{2,3,1} = -\frac{w_{N,4,2}}{w_{N,4,3}} \left(w_{2,2,1} + \tau_2(\tau_1 - \tau_2)\right) - w_{2,3,2} - \tau_2(\tau_1 - \tau_2)$$

$$w_{N,4,1} = w_{1,3,2} = 0$$

$$w_{N,4,2} = -w_{N,4,3}$$

**0-1-0-0**

$\xi_{4,N} = 0, \tau_3 \neq \tau_2, w_{N,4,3} = 0, w_{N,4,2} = 0$

$$w_{N,4,1} = 0$$

**0-1-0-1**

$\xi_{4,N} = 0, \tau_3 \neq \tau_2, w_{N,4,3} = 0, w_{N,4,2} \neq 0$

$$\tau_2 = \tau_1$$

$$w_{1,2,1} = 0$$

$$n > 1,\ w_{2,2,1} = 0$$

$$w_{N,4,1} = -w_{N,4,2}$$

**0-1-1**

$\xi_{4,N} = 0, \tau_3 \neq \tau_2, w_{N,4,3} \neq 0$

$$\begin{aligned} \tau_3 &= \tau_1 \\ n > 1,\ w_{2,3,1} &= -w_{2,3,2} \\ w_{1,3,2} &= w_{1,3,1} = w_{N,4,2} = 0 \\ w_{N,4,1} &= -w_{N,4,3} \end{aligned}$$

**1-0**

$\xi_{4,N} \neq 0, \tau_2 = \tau_1$

We now list the main operations performed to solve this system.

1) 3.3 to find $\tau_1$, replace everywhere where added to $\xi_{4,N}$ terms
2) 2.1 in 3.1, 3.2
3) 3.1 to find $\tau_3$
4) 3.2 to find $w_{1,3,1} + w_{1,3,2}$, replace everywhere along $\tau_1$

Which yields :

$$\begin{aligned} \tau_1 &= \frac{\xi_{4,N}}{2+N} \\ \tau_3 &= \xi_{4,N} \\ w_{1,2,1} &= \frac{\xi_{4,N}}{2+N} \\ w_{1,3,1} &= -\xi_{4,N} N \\ w_{1,3,2} &= \xi_{4,N}(1+N) \\ n > 1,\ w_{2,3,1} &= \xi_{4,N}^2 2\frac{1+N}{2+N} - N(2+N)w_{2,2,1} - w_{2,3,2} \\ w_{N,4,1} &= 0 \\ w_{N,4,3} &= \frac{\xi_{4,N}^N}{(1+N)^2} \\ w_{N,4,2} &= \xi_{4,N}^N N \frac{2+N}{(1+N)^2} \end{aligned}$$

**1-1-0**

$\xi_{4,N} \neq 0, \tau_2 \neq \tau_1, \tau_3 = \tau_2$

We now list the main operations performed to solve this system.

1) multiply 2.2 by $\tau_2 - \tau_1$
2) 2.1 in 2.2
3) 3.2 in 2.2
4) multiply 3.4 by $\tau_2 - \tau_1$
5) 3.3 in 3.4
6) 3.4 in 2.2
7) $\tau_2$ in 3.2, replace everywhere
8) $\tau_1 = 0$ in 3.1, replace everywhere
9) replace everywhere $w_{1,2,1}$ from 3.4
10) 2.2 to find $w_{1,3,1} + w_{1,3,2}$, replace everywhere
$w_{1,3,2} \neq 0$

Which yields :

$$\tau_1 = 0$$

$$\tau_2 = w_{1,2,1} = \xi_{4,N} \frac{2}{2+N}$$

$$w_{1,3,1} = \xi_{4,N} \frac{2}{2+N} - w_{1,3,2}$$

$$n > 1,\ w_{2,3,1} = \frac{w_{1,3,2}}{\xi_{4,N}} \left( \xi_{4,N}^2 \frac{4}{2+N} - \left( 2 + N - \frac{\xi_{4,N}}{w_{1,3,2}} \right) w_{2,2,1} \right) - w_{2,3,2}$$

$$w_{N,4,1} = \xi_{4,N}^N \frac{N}{2(1+N)}$$

$$w_{N,4,2} = \frac{\xi_{4,N}^N}{2(1+N)} \left( 2 + N - \frac{\xi_{4,N}}{w_{1,3,2}} \right)$$

$$w_{N,4,3} = \frac{\xi_{4,N}^{N+1}}{2w_{1,3,2}(1+N)}$$

**1-1-1-0-0**

$\xi_{4,N} \neq 0, \tau_2 \neq \tau_1, \tau_3 \neq \tau_2, \tau_3 = \tau_1, \tau_1 = \xi_{4,N}$

We now list the main operations performed to solve this system.

1) multiply 2.2 by $\tau_2 - \tau_1$
2) 2.1 in 2.2
3) 3.2 in 2.2
4) multiply 3.4 by $\tau_2 - \tau_1$
5) 3.3 in 3.4
6) 3.4 in 2.2
7) multiply 2.1 by $\tau_2 - \tau_3$
8) 3.1 in 2.1
9) 3.1 in 2.2
10) 3.1 to find $\tau_2$, replace everywhere
11) multiply 3.2 by $w_{1,3,2}$
12) 3.3 in 3.2
13) 3.2 in 3.5
14) 3.3 in 3.5 twice
15) 3.4 in 3.5
16) isolate $w_{1,3,2}$ in 3.5, replace everywhere

Which yields :

$$\begin{aligned}
\tau_2 = w_{1,2,1} &= \frac{\xi_{4,N}}{2+N} \\
w_{1,3,1} &= -\xi_{4,N} N \\
w_{1,3,2} &= \xi_{4,N}(1+N) \\
n > 1,\ w_{2,3,1} &= \xi_{4,N}^2 2\frac{1+N}{2+N} - (2+N)N w_{2,2,1} - w_{2,3,2} \\
w_{N,4,1} &= 0 \\
w_{N,4,2} &= \xi_{4,N}^N \frac{(2+N)N}{(1+N)^2} \\
w_{N,4,3} &= \frac{\xi_{4,N}^N}{(1+N)^2}
\end{aligned}$$

**1-1-1-0-1**

$\xi_{4,N} \neq 0, \tau_2 \neq \tau_1, \tau_3 \neq \tau_2, \tau_3 = \tau_1, \tau_1 \neq \xi_{4,N}$

We now list the main operations performed to solve this system.

1) multiply 2.2 by $\tau_2 - \tau_1$
2) 2.1 in 2.2
3) 3.2 in 2.2
4) multiply 3.4 by $\tau_2 - \tau_1$
5) 3.3 in 3.4
6) 3.4 in 2.2
7) multiply 2.1 by $\tau_2 - \tau_3$
8) 3.1 in 2.1
9) 3.1 in 2.2
10) 2.2 to find $\tau_1 = 0$, replace everywhere
11) 3.1 to find $\tau_2$, replace everywhere
12) 3.2 to find $w_{1,3,1} + w_{1,3,2}$, replace everywhere

Which yields :

$$
\begin{aligned}
w_{1,3,2} &\neq 0 \\
\tau_1 &= 0 \\
\tau_2 = w_{1,2,1} &= \xi_{4,N} \frac{2}{2+N} \\
w_{1,3,1} &= -w_{1,3,2} \\
n > 1,\ w_{2,3,1} &= \frac{w_{1,3,2}}{\xi_{4,N}} \left( 4 \frac{\xi_{4,N}^2}{2+N} - (2+N) w_{2,2,1} \right) - w_{2,3,2} \\
w_{N,4,1} &= \frac{\xi_{4,N}^N}{2(1+N)} \left( N - \frac{\xi_{4,N}}{w_{1,3,2}} \right) \\
w_{N,4,2} &= \xi_{4,N}^N \frac{2+N}{2(1+N)} \\
w_{N,4,3} &= \frac{\xi_{4,N}^{N+1}}{2 w_{1,3,2} (1+N)}
\end{aligned}
$$

**1-1-1-1-0**

$\xi_{4,N} \neq 0, \tau_2 \neq \tau_1, \tau_3 \neq \tau_2, \tau_3 \neq \tau_1, \tau_1 = 0$

We now list the main operations performed to solve this system.

1) multiply 2.2 by $\tau_2 - \tau_1$
2) 2.1 in 2.2
3) 3.2 in 2.2
4) multiply 3.4 by $\tau_2 - \tau_1$
5) 3.3 in 3.4
6) 3.4 in 2.2
7) multiply 2.1 by $\tau_2 - \tau_3$
8) 3.1 in 2.1
9) multiply 4.5 and 4.2 by $\tau_3 - \tau_1$
10) 2.2 in 4.5
11) 3.1 in 3.2
12) 3.4 in 3.5
13) 3.2 in 3.5
14) multiply 2.2 by $\tau_3 - \tau_2$
15) 3.1 in 2.2
16) 3.4 in 2.2
17) 3.2 in 2.2
18) multiply 3.3 by $(\tau_3 - \tau_2)\tau_3$

Which yields :

$$\begin{aligned}
&\tau_2 \neq \xi_{4,N}\frac{2}{2+N} \\
&w_{1,2,1} = \tau_2 \\
&w_{1,3,1} = \tau_3 - \xi_{4,N}\frac{\tau_3(\tau_3-\tau_2)}{\tau_2\big(\xi_{4,N}2 - \tau_2(2+N)\big)} \\
&w_{1,3,2} = \xi_{4,N}\frac{\tau_3(\tau_3-\tau_2)}{\tau_2\big(\xi_{4,N}2 - \tau_2(2+N)\big)} \\
n > 1,\ &w_{2,3,1} = \frac{\tau_3}{\xi_{4,N}2 - \tau_2(2+N)}\left(\xi_{4,N}2(\tau_3-\tau_2) + \frac{\xi_{4,N}2 - \tau_3(2+N)}{\tau_2}w_{2,2,1}\right) - w_{2,3,2} \\
&w_{N,4,1} = \xi_{4,N}^{N} - \xi_{4,N}^{1+N}\frac{(\tau_3+\tau_2)(2+N) - \xi_{4,N}2}{\tau_2\tau_3(1+N)(2+N)} \\
&w_{N,4,2} = \xi_{4,N}^{1+N}\frac{\xi_{4,N}2 - \tau_3(2+N)}{\tau_2(\tau_2-\tau_3)(1+N)(2+N)} \\
&w_{N,4,3} = \xi_{4,N}^{1+N}\frac{\xi_{4,N}2 - \tau_2(2+N)}{\tau_3(\tau_3-\tau_2)(1+N)(2+N)}
\end{aligned}$$

**1-1-1-1-1**

$\xi_{4,N} \neq 0, \tau_2 \neq \tau_1, \tau_3 \neq \tau_2, \tau_3 \neq \tau_1, \tau_1 \neq 0$

We now list the main operations performed to solve this system.

1) multiply 2.2 by $\tau_2 - \tau_1$
2) 2.1 in 2.2
3) 3.2 in 2.2
4) multiply 3.4 by $\tau_2 - \tau_1$
5) 3.3 in 3.4
6) 3.4 in 2.2
7) multiply 2.1 by $\tau_2 - \tau_3$
8) 3.1 in 2.1
9) multiply 4.5 and 4.2 by $\tau_3 - \tau_1$
10) 2.2 in 4.5
11) 3.1 in 3.2
12) 3.4 in 3.5
13) 3.2 in 3.5
14) multiply 2.2 by $\tau_3 - \tau_2$
15) 3.1 in 2.2
16) 3.4 in 2.2
17) 3.2 in 2.2
18) 2.2 to find $\tau_3 = \xi_{4,N}$, replace everywhere
19) 3.5 to find $w_{1,3,1} + w_{1,3,2}$, replace everywhere
20) 3.2 to find $\tau_2$, replace everywhere
21) 3.1 in 3.3

Which yields :

$$\begin{aligned}
\tau_1 &\neq \frac{\xi_{4,N}}{2+N} \\
\tau_2 &= \frac{\xi_{4,N}}{2+N} \\
\tau_3 &= \xi_{4,N} \\
w_{1,2,1} &= \frac{\xi_{4,N}}{2+N} \\
w_{1,3,1} &= -\xi_{4,N} N \\
w_{1,3,2} &= \xi_{4,N}(1+N) \\
n > 1,\ w_{2,3,1} &= \xi_{4,N}^2 2\frac{1+N}{2+N} - N(2+N)w_{2,2,1} - w_{2,3,2} \\
w_{N,4,1} &= 0 \\
w_{N,4,2} &= \xi_{4,N}^N \frac{N(2+N)}{(1+N)^2} \\
w_{N,4,3} &= \frac{\xi_{4,N}^N}{(1+N)^2}
\end{aligned}$$

**C.5 Explicit method**, $s = 4, j = 5, v = 3$

The conditions are :

Order of consistency 1 :

$$w_{N,5,1} + w_{N,5,2} + w_{N,5,3} + w_{N,5,4} = \xi_{5,N}^{N}$$

Order of consistency 2 :

$$w_{N,5,1}\tau_1 + w_{N,5,2}\tau_2 + w_{N,5,3}\tau_3 + w_{N,5,4}\tau_4 = \frac{\xi_{5,N}^{N+1}}{1+N}$$

$$w_{N,5,2}w_{1,2,1} + w_{N,5,3}\big(w_{1,3,1} + w_{1,3,2}\big) + w_{N,5,4}\big(w_{1,4,1} + w_{1,4,2} + w_{1,4,3}\big) = \frac{\xi_{5,N}^{N+1}}{1+N}$$

Order of consistency 3 :

$$w_{N,5,1}\tau_1^2 + w_{N,5,2}\tau_2^2 + w_{N,5,3}\tau_3^2 + w_{N,5,4}\tau_4^2 = \xi_{5,N}^{N+2}2\frac{N!}{(2+N)!}$$

$$w_{N,5,2}\tau_2 w_{1,2,1} + w_{N,5,3}\tau_3\big(w_{1,3,1} + w_{1,3,2}\big) + w_{N,5,4}\tau_4\big(w_{1,4,1} + w_{1,4,2} + w_{1,4,3}\big) = \xi_{5,N}^{N+2}2\frac{N!}{(2+N)!}$$

$$w_{N,5,2}w_{1,2,1}\tau_1 + w_{N,5,3}\big(w_{1,3,1}\tau_1 + w_{1,3,2}\tau_2\big) + w_{N,5,4}\big(w_{1,4,1}\tau_1 + w_{1,4,2}\tau_2 + w_{1,4,3}\tau_3\big) = \xi_{5,N}^{N+2}\frac{N!}{(2+N)!}$$

$$w_{N,5,3}w_{1,3,2}w_{1,2,1} + w_{N,5,4}\big(w_{1,4,2}w_{1,2,1} + w_{1,4,3}\big(w_{1,3,1} + w_{1,3,2}\big)\big) = \xi_{5,N}^{N+2}\frac{N!}{(2+N)!}$$

$$w_{N,5,2}w_{1,2,1}^2 + w_{N,5,3}\big(w_{1,3,1} + w_{1,3,2}\big)^2 + w_{N,5,4}\big(w_{1,4,1} + w_{1,4,2} + w_{1,4,3}\big)^2 = \xi_{5,N}^{N+2}2\frac{N!}{(2+N)!}$$

$$n > 1,\ \ w_{N,5,2}w_{2,2,1} + w_{N,5,3}\big(w_{2,3,1} + w_{2,3,2}\big) + w_{N,5,4}\big(w_{2,4,1} + w_{2,4,2} + w_{2,4,3}\big) = \xi_{5,N}^{N+2}2\frac{N!}{(2+N)!}$$

To simplify the system of equations we assume the method satisfies the conditions for order of consistency 1 at all stages at rank 1 :

Order of consistency 1 :

$$w_{N,5,1} + w_{N,5,2} + w_{N,5,3} + w_{N,5,4} = \xi_{5,N}^{N}$$

Order of consistency 2 :

$$w_{N,5,2}\tau_2 + w_{N,5,3}\tau_3 + w_{N,5,4}\tau_4 = \frac{\xi_{5,N}^{1+N}}{1+N}$$

$$0 = 0$$

Order of consistency 3 :

$$w_{N,5,2}\tau_2^2 + w_{N,5,3}\tau_3^2 + w_{N,5,4}\tau_4^2 = 2\xi_{5,N}^{2+N}\frac{N!}{(2+N)!}$$

$$0 = 0$$

$$w_{N,5,3}w_{1,3,2}\tau_2 + w_{N,5,4}\big(w_{1,4,3}\tau_3 + w_{1,4,2}\tau_2\big) = \xi_{5,N}^{2+N}\frac{N!}{(2+N)!}$$

$$0 = 0$$
$$0 = 0$$
$$n > 1,\ w_{N,5,2}w_{2,2,1} + w_{N,5,3}\big(w_{2,3,2} + w_{2,3,1}\big) + w_{N,5,4}\big(w_{2,4,3} + w_{2,4,2} + w_{2,4,1}\big) = 2\xi_{5,N}^{2+N}\frac{N!}{(2+N)!}$$

We now list the main operations performed to rewrite this system.

1) isolate $w_{N,5,1}$ in 1.1
2) $w_{N,5,2}\tau_2$ in 2.1
3) 2.1 in 3.1
4) isolate $w_{N,5,3}\tau_3(\tau_3 - \tau_2)$ in 3.1

Which yields :

Order of consistency 1 :

$$w_{N,5,1} = \xi_{5,N}^{N} - w_{N,5,2} - w_{N,5,3} - w_{N,5,4}$$

Order of consistency 2 :

$$w_{N,5,2}\tau_2 = \frac{\xi_{5,N}^{1+N}}{1+N} - w_{N,5,3}\tau_3 - w_{N,5,4}\tau_4$$
$$0 = 0$$

Order of consistency 3 :

$$w_{N,5,3}\tau_3(\tau_3 - \tau_2) = \frac{\xi_{5,N}^{1+N}}{1+N}\left(2\frac{\xi_{5,N}}{2+N} - \tau_2\right) - w_{N,5,4}\tau_4(\tau_4 - \tau_2)$$
$$0 = 0$$
$$w_{N,5,3}w_{1,3,2}\tau_2 + w_{N,5,4}\big(w_{1,4,3}\tau_3 + w_{1,4,2}\tau_2\big) = \xi_{5,N}^{2+N}\frac{N!}{(2+N)!}$$
$$0 = 0$$
$$0 = 0$$
$$n > 1,\ w_{N,5,2}w_{2,2,1} + w_{N,5,3}\big(w_{2,3,2} + w_{2,3,1}\big) + w_{N,5,4}\big(w_{2,4,3} + w_{2,4,2} + w_{2,4,1}\big) = 2\xi_{5,N}^{2+N}\frac{N!}{(2+N)!}$$

The solutions are :

**0-0-0-0**

$\xi_{5,N} = 0, w_{N,5,4} = 0, w_{N,5,3} = 0, w_{N,5,2} = 0$

$$w_{N,5,1} = 0$$

**0-0-0-1**

$\xi_{5,N} = 0, w_{N,5,4} = 0, w_{N,5,3} = 0, w_{N,5,2} \neq 0$

$$\tau_2 = 0$$
$$n > 1,\ w_{2,2,1} = 0$$
$$w_{N,5,1} = -w_{N,5,2}$$

**0-0-1-0**

$\xi_{5,N}=0, w_{N,5,4}=0, w_{N,5,3}\neq 0, \tau_2=0$

$$\tau_3=0$$

$$n>1,\; w_{2,3,1}=-\frac{w_{N,5,2}}{w_{N,5,3}}w_{2,2,1}-w_{2,3,2}$$

$$w_{N,5,1}=-w_{N,5,2}-w_{N,5,3}$$

**0-0-1-1-0**

$\xi_{5,N}=0, w_{N,5,4}=0, w_{N,5,3}\neq 0, \tau_2\neq 0, w_{N,5,2}=0$

$$\tau_3=w_{1,3,2}=0$$

$$n>1,\; w_{2,3,1}=-w_{2,3,1}$$

$$w_{N,5,1}=-w_{N,5,2}-w_{N,5,3}$$

**0-0-1-1-1**

$\xi_{5,N}=0, w_{N,5,4}=0, w_{N,5,3}\neq 0, w_{N,5,2}\neq 0, \tau_2\neq 0$

$$\tau_2=\tau_3$$

$$w_{1,3,2}=0$$

$$n>1,\; w_{2,3,1}=-\frac{w_{N,5,2}}{w_{N,5,3}}w_{2,2,1}-w_{2,3,2}$$

$$w_{N,5,1}=-w_{N,5,2}-w_{N,5,3}$$

$$w_{N,5,2}=-w_{N,5,3}$$

**0-1-0-0-0**

$\xi_{5,N}=0, w_{N,5,4}\neq 0, \tau_2=0, \tau_4=0, \tau_3=0$

$$n>1,\; w_{2,4,1}=-\frac{w_{N,5,2}}{w_{N,5,4}}w_{2,2,1}+\frac{w_{N,5,3}}{w_{N,5,4}}\left(w_{2,3,2}+w_{2,3,1}\right)-w_{2,4,3}-w_{2,4,2}$$

$$w_{N,5,1}=-w_{N,5,2}-w_{N,5,3}-w_{N,5,4}$$

**0-1-0-0-1**

$\xi_{5,N} = 0, w_{N,5,4} \neq 0, \tau_2 = 0, \tau_4 = 0, \tau_3 \neq 0$

$$n > 1,\ w_{2,4,1} = -\frac{w_{N,5,2}}{w_{N,5,4}} w_{2,2,1} - w_{2,4,3} - w_{2,4,2}$$
$$w_{N,5,1} = -w_{N,5,2} - -w_{N,5,4}$$
$$w_{N,5,3} = w_{1,4,3} = 0$$

**0-1-0-1**

$\xi_{5,N} = 0, w_{N,5,4} \neq 0, \tau_2 = 0, \tau_4 \neq 0$

$$\tau_4 = \tau_3$$
$$w_{1,4,3} = 0$$
$$n > 1,\ w_{2,4,1} = -\frac{w_{N,5,2}}{w_{N,5,4}} w_{2,2,1} - w_{2,3,2} - w_{2,3,1} - w_{2,4,3} - w_{2,4,2}$$
$$w_{N,5,1} = -w_{N,5,2}$$
$$w_{N,5,3} = -w_{N,5,4}$$

**0-1-1-0-0-0**

$\xi_{5,N} = 0, w_{N,5,4} \neq 0, \tau_2 \neq 0, \tau_3 = 0, \tau_4 = 0$

$$w_{1,4,2} = -\frac{w_{N,5,3}}{w_{N,5,4}} w_{1,3,2}$$
$$n > 1,\ w_{2,4,1} = \frac{w_{N,5,3}}{w_{N,5,4}} \big(w_{2,3,2} + w_{2,3,1}\big) - w_{2,4,3} - w_{2,4,2}$$
$$w_{N,5,1} = -w_{N,5,3} - w_{N,5,4}$$
$$w_{N,5,2} = 0$$

**0-1-1-0-1-1**

$\xi_{5,N} = 0, w_{N,5,4} \neq 0, \tau_2 \neq 0, \tau_3 = 0, \tau_4 \neq 0$

$$\tau_2 = \tau_4$$
$$w_{1,4,2} = -\frac{w_{N,5,3}}{w_{N,5,4}} w_{1,3,2}$$
$$n > 1,\ w_{2,4,1} = w_{2,2,1} + \frac{w_{N,5,3}}{w_{N,5,4}} \big(w_{2,3,2} + w_{2,3,1}\big) - w_{2,4,3} - w_{2,4,2}$$
$$w_{N,5,1} = -w_{N,5,3}$$
$$w_{N,5,2} = -w_{N,5,4}$$

**0-1-1-1-0-0**

$\xi_{5,N} = 0, w_{N,5,4} \neq 0, \tau_2 \neq 0, \tau_3 \neq 0, \tau_3 = \tau_2, \tau_4 = 0$

$$
\begin{aligned}
w_{1,4,2} &= -\frac{w_{N,5,3}}{w_{N,5,4}} w_{1,3,2} - w_{1,4,3} \\
n > 1,\ w_{2,4,1} &= \frac{w_{N,5,3}}{w_{N,5,4}} \big(w_{2,2,1} + w_{2,3,2} + w_{2,3,1}\big) - w_{2,4,3} - w_{2,4,2} \\
w_{N,5,1} &= -w_{N,5,4} \\
w_{N,5,2} &= -w_{N,5,3}
\end{aligned}
$$

**0-1-1-1-0-1**

$\xi_{5,N} = 0, w_{N,5,4} \neq 0, \tau_2 \neq 0, \tau_3 \neq 0, \tau_3 = \tau_2, \tau_4 \neq 0$

$$
\begin{aligned}
\tau_2 &= \tau_4 \\
w_{1,4,2} &= -\frac{w_{N,5,3}}{w_{N,5,4}} w_{1,3,2} - w_{1,4,3} \\
n > 1,\ w_{2,4,1} &= \frac{w_{N,5,3} + w_{N,5,4}}{w_{N,5,4}} w_{2,2,1} + \frac{w_{N,5,3}}{w_{N,5,4}} \big(w_{2,3,2} + w_{2,3,1}\big) - w_{2,4,3} - w_{2,4,2} \\
w_{N,5,1} &= 0 \\
w_{N,5,2} &= -w_{N,5,3} - w_{N,5,4}
\end{aligned}
$$

**0-1-1-1-1**

$\xi_{5,N} = 0, w_{N,5,4} \neq 0, \tau_2 \neq 0, \tau_3 \neq 0, \tau_3 \neq \tau_2$

$$
\begin{aligned}
w_{1,4,2} &= \frac{\tau_4(\tau_2 - \tau_4)}{\tau_3(\tau_2 - \tau_3)} w_{1,3,2} - w_{1,4,3} \frac{\tau_3}{\tau_2} \\
n > 1,\ w_{2,4,1} &= \frac{\tau_4(\tau_3 - \tau_4)}{\tau_2(\tau_3 - \tau_2)} w_{2,2,1} - \frac{\tau_4(\tau_2 - \tau_4)}{\tau_3(\tau_2 - \tau_3)} \big(w_{2,3,2} + w_{2,3,1}\big) - w_{2,4,3} - w_{2,4,2} \\
w_{N,5,1} &= -w_{N,5,4} \frac{(\tau_3 - \tau_4)(\tau_2 - \tau_4)}{\tau_3 \tau_2} \\
w_{N,5,2} &= -w_{N,5,4} \frac{\tau_4(\tau_3 - \tau_4)}{\tau_2(\tau_3 - \tau_2)} \\
w_{N,5,3} &= -w_{N,5,4} \frac{\tau_4(\tau_2 - \tau_4)}{\tau_3(\tau_2 - \tau_3)}
\end{aligned}
$$

**1-0-0**

$\xi_{5,N} \neq 0, \tau_2 = 0, \tau_3 = \tau_4$

$$\tau_4 = 2\frac{\xi_{5,N}}{2+N}$$

$$n > 1,\ w_{2,4,1} = 4w_{1,4,3}\frac{\xi_{5,N}}{2+N} - 2w_{N,5,2}w_{1,4,3}w_{2,2,1}(1+N) - \left(w_{1,4,3}\frac{2+N}{\xi_{5,N}} - 1\right)(w_{2,3,2}+w_{2,3,1}) - w_{2,4,3} - w_{2,4,2}$$

$$w_{1,4,3} \neq 0$$

$$w_{N,5,1} = \xi_{5,N}^{N}\frac{N}{2(1+N)} - w_{N,5,2}$$

$$w_{N,5,3} = \frac{\xi_{5,N}^{N}}{2(1+N)}\left(2+N-\frac{\xi_{5,N}}{w_{1,4,3}}\right)$$

$$w_{N,5,4} = \frac{\xi_{5,N}^{1+N}}{2(1+N)w_{1,4,3}}$$

**1-0-1-0**

$\xi_{5,N} \neq 0, \tau_2 = 0, \tau_3 \neq \tau_4, \tau_4 = 0$

$w_{1,4,3} \neq 0$

$$\tau_3 = 2\frac{\xi_{5,N}}{2+N}$$

$$n > 1,\ w_{2,4,1} = w_{1,4,3}\left(4\frac{\xi_{5,N}}{2+N} - 2(1+N)\frac{w_{N,5,2}}{\xi_{5,N}^{1+N}}w_{2,2,1} - 4\frac{\xi_{5,N}}{\tau_3(2+N)(\tau_3-\tau_4)}(w_{2,3,2}+w_{2,3,1})\right) - w_{2,4,3} - w_{2,4,2}$$

$$w_{N,5,1} = \xi_{5,N}^{N} - w_{N,5,2} - 2\frac{\xi_{5,N}^{2+N}}{\tau_3(1+N)(2+N)(\tau_3-\tau_4)} - \frac{\xi_{5,N}^{1+N}}{2w_{1,4,3}(1+N)}$$

$$w_{N,5,3} = 2\frac{\xi_{5,N}^{2+N}}{\tau_3(1+N)(2+N)(\tau_3-\tau_4)}$$

$$w_{N,5,4} = \frac{\xi_{5,N}^{1+N}}{2w_{1,4,3}(1+N)}$$

**1-0-1-1**

$\xi_{5,N} \neq 0, \tau_2 = 0, \tau_3 \neq \tau_4, \tau_4 \neq 0$

$$w_{N,5,4} w_{1,4,3} \tau_3 \neq 0$$

$$\tau_3 \neq 2\frac{\xi_{5,N}}{2+N}$$

$$w_{1,4,3} = \xi_{5,N} \frac{\tau_4(\tau_3 - \tau_4)}{\tau_3\big(\tau_3(2+N) - 2\xi_{5,N}\big)}$$

$$n > 1,\; w_{2,4,1} = \frac{\tau_4(\tau_3 - \tau_4)}{\tau_3(2+N) - 2\xi_{5,N}} \left(2\xi_{5,N} - w_{N,5,2}(2+N)\frac{1+N}{\xi_{5,N}^{1+N}} w_{2,2,1}\right)$$
$$+ \frac{\tau_4\big(2\xi_{5,N} - \tau_4(2+N)\big)}{\tau_3\big(2\xi_{5,N} - \tau_3(2+N)\big)}\big(w_{2,3,2} + w_{2,3,1}\big) - w_{2,4,3} - w_{2,4,2}$$

$$w_{N,5,1} = \xi_{5,N}^N - w_{N,5,2} - w_{N,5,3} - w_{N,5,4}$$

$$w_{N,5,4} = \frac{\xi_{5,N}^{1+N}}{(1+N)\tau_4(\tau_3 - \tau_4)} \left(\tau_3 - 2\frac{\xi_{5,N}}{2+N}\right)$$

$$w_{N,5,3} = \frac{\xi_{5,N}^{1+N}}{\tau_3(1+N)(\tau_3 - \tau_4)} \left(2\frac{\xi_{5,N}}{2+N} - \tau_4\right)$$

**1-1-0-0-0**

$\xi_{5,N} \neq 0, \tau_2 \neq 0, \tau_3 = 0, \tau_4 = 0, w_{1,3,2} = 0$

$$w_{1,4,2} \neq 0$$

$$\tau_2 = 2\frac{\xi_{5,N}}{2+N}$$

$$n > 1,\; w_{2,2,1} = \frac{2}{(2+N)} \left(-w_{N,5,3}\frac{1+N}{\xi_{5,N}^N}\big(w_{2,3,2} + w_{2,3,1}\big) - \frac{\xi_{5,N}}{2w_{1,4,2}}\big(w_{2,4,3} + w_{2,4,2} + w_{2,4,1}\big) + 2\frac{\xi_{5,N}^2}{2+N}\right)$$

$$w_{N,5,1} = \xi_{5,N}^N \frac{N}{2(1+N)} - w_{N,5,3} - \frac{\xi_{5,N}^{1+N}}{2(1+N)w_{1,4,2}}$$

$$w_{N,5,2} = \xi_{5,N}^N \frac{2+N}{2(1+N)}$$

$$w_{N,5,4} = \frac{\xi_{5,N}^{1+N}}{2(1+N)w_{1,4,2}}$$

**1-1-0-0-1**

$\xi_{5,N} \neq 0, \tau_2 \neq 0, \tau_3 = 0, \tau_4 = 0, w_{1,3,2} \neq 0$

$$\tau_2 = 2\frac{\xi_{5,N}}{2+N}$$

$$n > 1,\ w_{2,2,1} = \frac{2}{(2+N)}\left(-\frac{1}{w_{1,3,2}}\left(\frac{\xi_{5,N}}{2} - w_{N,5,4}\frac{1+N}{\xi_{5,N}^{N}}w_{1,4,2}\right)(w_{2,3,2}+w_{2,3,1}) - w_{N,5,4}\frac{1+N}{\xi_{5,N}^{N}}(w_{2,4,3}+w_{2,4,2}+w_{2,4,1}) + 2\frac{\xi_{5,N}^{2}}{2+N}\right)$$

$$w_{N,5,1} = \xi_{5,N}^{N}\frac{N}{2(1+N)} - \frac{1}{w_{1,3,2}}\left(\frac{\xi_{5,N}^{1+N}}{2(1+N)} - w_{N,5,4}w_{1,4,2}\right) - w_{N,5,4}$$

$$w_{N,5,2} = \xi_{5,N}^{N}\frac{2+N}{2(1+N)}$$

$$w_{N,5,3} = \frac{1}{w_{1,3,2}}\left(\frac{\xi_{5,N}^{1+N}}{2(1+N)} - w_{N,5,4}w_{1,4,2}\right)$$

**1-1-0-1-0-0**

$\xi_{5,N} \neq 0, \tau_2 \neq 0, \tau_3 = 0, \tau_4 \neq 0, \tau_4 = \tau_2, w_{1,3,2} = 0$

$$w_{1,4,2} \neq 0$$

$$\tau_2 = 2\frac{\xi_{5,N}}{2+N}$$

$$n > 1,\ w_{2,4,1} = 4w_{1,4,2}\frac{\xi_{5,N}}{2+N} - 4\left((2+N)\frac{w_{1,4,2}}{\xi_{5,N}} - 1\right)w_{2,2,1} - w_{N,5,3}2w_{1,4,2}\frac{1+N}{\xi_{5,N}^{1+N}}(w_{2,3,2}+w_{2,3,1}) - w_{2,4,3} - w_{2,4,2}$$

$$w_{N,5,1} = \xi_{5,N}^{N} - 2\frac{\xi_{5,N}^{N}}{1+N}\left(2+N-\frac{\xi_{5,N}}{w_{1,4,2}}\right) - w_{N,5,3} - \frac{\xi_{5,N}^{1+N}}{2(1+N)w_{1,4,2}}$$

$$w_{N,5,2} = 2\frac{\xi_{5,N}^{N}}{1+N}\left(2+N-\frac{\xi_{5,N}}{w_{1,4,2}}\right)$$

$$w_{N,5,4} = \frac{\xi_{5,N}^{1+N}}{2(1+N)w_{1,4,2}}$$

**1-1-0-1-0-1-0**

$\xi_{5,N} \neq 0, \tau_2 \neq 0, \tau_3 = 0, \tau_4 \neq 0, \tau_4 = \tau_2, w_{1,3,2} \neq 0, w_{N,5,2} = 0$

$$\tau_2 = 2\frac{\xi_{5,N}}{2+N}$$

$$n > 1,\ w_{2,4,1} = 4\frac{\xi_{5,N}^2}{(2+N)^2} - \frac{\xi_{5,N} - (2+N)w_{1,4,2}}{w_{1,3,2}(2+N)}(w_{2,3,2} + w_{2,3,1}) - w_{2,4,3} - w_{2,4,2}$$

$$w_{N,5,1} = \xi_{5,N}^N \left( \frac{N}{2(1+N)} - \frac{\xi_{5,N} - (2+N)w_{1,4,2}}{w_{1,3,2}2(1+N)} \right)$$

$$w_{N,5,4} = \xi_{5,N}^N \frac{2+N}{2(1+N)}$$

$$w_{N,5,3} = \xi_{5,N}^N \frac{\xi_{5,N} - (2+N)w_{1,4,2}}{w_{1,3,2}2(1+N)}$$

**1-1-0-1-0-1-1**

$\xi_{5,N} \neq 0, \tau_2 \neq 0, \tau_3 = 0, \tau_4 \neq 0, \tau_4 = \tau_2, w_{1,3,2} \neq 0, w_{N,5,2} \neq 0$

$$\tau_2 = 2\frac{\xi_{5,N}}{2+N}$$

$$n > 1,\ w_{2,2,1} = \frac{1}{w_{N,5,2}} \left( -\frac{1}{w_{1,3,2}} \left( \frac{\xi_{5,N}^{1+N}}{2(1+N)} - w_{1,4,2} \left( \xi_{5,N}^N \frac{2+N}{2(1+N)} - w_{N,5,2} \right) \right) (w_{2,3,2} + w_{2,3,1}) \right.$$
$$\left. - \left( \xi_{5,N}^N \frac{2+N}{2(1+N)} - w_{N,5,2} \right) (w_{2,4,3} + w_{2,4,2} + w_{2,4,1}) + 2\xi_{5,N}^{2+N} \frac{N!}{(2+N)!} \right)$$

$$w_{N,5,1} = \xi_{5,N}^N \frac{N}{2(1+N)} - \frac{1}{w_{1,3,2}} \left( \frac{\xi_{5,N}^{1+N}}{2(1+N)} - w_{1,4,2} \left( \xi_{5,N}^N \frac{2+N}{2(1+N)} - w_{N,5,2} \right) \right)$$

$$w_{N,5,3} = \frac{1}{w_{1,3,2}} \left( \frac{\xi_{5,N}^{1+N}}{2(1+N)} - w_{1,4,2} \left( \xi_{5,N}^N \frac{2+N}{2(1+N)} - w_{N,5,2} \right) \right)$$

$$w_{N,5,4} = \xi_{5,N}^N \frac{2+N}{2(1+N)} - w_{N,5,2}$$

**1-1-0-1-1-0**

$\xi_{5,N} \neq 0, \tau_2 \neq 0, \tau_3 = 0, \tau_4 \neq 0, \tau_4 \neq \tau_2, w_{1,3,2} = 0$

$$w_{1,4,2} \neq 0$$

$$\tau_2 \neq 2\frac{\xi_{5,N}}{2+N}$$

$$w_{1,4,2} = \xi_{5,N}\frac{\tau_4(\tau_4-\tau_2)}{\tau_2\big(2\xi_{5,N}-\tau_2(2+N)\big)}$$

$$n>1,\ w_{2,4,1} = \tau_4(\tau_4-\tau_2)\big(2\xi_{5,N}-\tau_2(2+N)\big)\left(-\frac{1}{\tau_2(\tau_4-\tau_2)}\big(\tau_4(2+N)-2\xi_{5,N}\big)w_{2,2,1}\right.$$
$$\left.-(1+N)(2+N)\frac{w_{N,5,3}}{\xi_{5,N}^{1+N}}\big(w_{2,3,2}+w_{2,3,1}\big)+2\xi_{5,N}\right)-w_{2,4,3}-w_{2,4,2}$$

$$w_{N,5,1} = \xi_{5,N}^{N} - \frac{\xi_{5,N}^{1+N}}{\tau_2(1+N)(\tau_4-\tau_2)}\left(\tau_4-2\frac{\xi_{5,N}}{2+N}\right) - w_{N,5,3} - \frac{\xi_{5,N}^{1+N}}{\tau_4(1+N)(\tau_4-\tau_2)}\left(2\frac{\xi_{5,N}}{2+N}-\tau_2\right)$$

$$w_{N,5,2} = \frac{\xi_{5,N}^{1+N}}{\tau_2(1+N)(\tau_4-\tau_2)}\left(\tau_4-2\frac{\xi_{5,N}}{2+N}\right)$$

$$w_{N,5,4} = \frac{\xi_{5,N}^{1+N}}{\tau_4(1+N)(\tau_4-\tau_2)}\left(2\frac{\xi_{5,N}}{2+N}-\tau_2\right)$$

**1-1-0-1-1-1-0**

$\xi_{5,N} \neq 0, \tau_2 \neq 0, \tau_3 = 0, \tau_4 \neq 0, \tau_4 \neq \tau_2, w_{1,3,2} \neq 0, w_{N,5,2} = 0$

$$\tau_4 = 2\frac{\xi_{5,N}}{2+N}$$

$$n>1,\ w_{2,4,1} = 4\frac{\xi_{5,N}^2}{(2+N)^2} - 2\frac{1+N}{2+N}\frac{w_{N,5,2}}{\xi_{5,N}^{N}}w_{2,2,1}$$
$$-\frac{1}{w_{1,3,2}}\left(2\frac{\xi_{5,N}^2}{\tau_2(2+N)^2}-w_{1,4,2}\right)\big(w_{2,3,2}+w_{2,3,1}\big)-w_{2,4,3}-w_{2,4,2}$$

$$w_{N,5,1} = \frac{\xi_{5,N}^{N}}{2(1+N)}\left(N-\frac{\xi_{5,N}}{w_{1,3,2}}\left(2\frac{\xi_{5,N}}{\tau_2(2+N)}-\frac{2+N}{\xi_{5,N}}w_{1,4,2}\right)\right)$$

$$w_{N,5,4} = \xi_{5,N}^{N}\frac{2+N}{2(1+N)}$$

$$w_{N,5,3} = \frac{\xi_{5,N}^{1+N}}{w_{1,3,2}(1+N)(2+N)}\left(\frac{\xi_{5,N}}{\tau_2}-\frac{(2+N)^2}{2\xi_{5,N}}w_{1,4,2}\right)$$

**1-1-0-1-1-1-1**

$\xi_{5,N} \neq 0, \tau_2 \neq 0, \tau_3 = 0, \tau_4 \neq 0, \tau_4 \neq \tau_2, w_{1,3,2} \neq 0, w_{N,5,2} \neq 0$

$$n > 1,\ w_{2,2,1} = \frac{1}{\tau_4(2+N) - 2\xi_{5,N}} \left( -\frac{\tau_2}{w_{1,3,2}} \left( \xi_{5,N} \frac{\tau_4 - \tau_2}{\tau_2} - \frac{2\xi_{5,N} - \tau_2(2+N)}{\tau_4} w_{1,4,2} \right) \big(w_{2,3,2} + w_{2,3,1}\big) + \frac{\tau_2}{\tau_4}\big(\tau_2(2+N) - 2\xi_{5,N}\big)\big(w_{2,4,3} + w_{2,4,2} + w_{2,4,1}\big) + 2\tau_2(\tau_4 - \tau_2)\xi_{5,N} \right)$$

$$w_{N,5,1} = \xi_{5,N}^N - \frac{\xi_{5,N}^{1+N}}{(1+N)(2+N)} \left( \frac{1}{\tau_2(\tau_4 - \tau_2)} \big(\tau_4(2+N) - 2\xi_{5,N}\big) + \frac{1}{\tau_4(\tau_4 - \tau_2)} \big(2\xi_{5,N} - \tau_2(2+N)\big) + \frac{1}{w_{1,3,2}} \left( \frac{\xi_{5,N}}{\tau_2} - \frac{2\xi_{5,N} - \tau_2(2+N)}{\tau_4(\tau_4 - \tau_2)} w_{1,4,2} \right) \right)$$

$$w_{N,5,2} = \frac{\xi_{5,N}^{1+N}}{\tau_2(1+N)(\tau_4 - \tau_2)} \left( \tau_4 - 2\frac{\xi_{5,N}}{2+N} \right)$$

$$w_{N,5,3} = \frac{\xi_{5,N}^{1+N}}{w_{1,3,2}(1+N)(2+N)} \left( \frac{\xi_{5,N}}{\tau_2} - \frac{2\xi_{5,N} - \tau_2(2+N)}{\tau_4(\tau_4 - \tau_2)} w_{1,4,2} \right)$$

$$w_{N,5,4} = \frac{\xi_{5,N}^{1+N}}{\tau_4(1+N)(\tau_4 - \tau_2)} \left( 2\frac{\xi_{5,N}}{2+N} - \tau_2 \right)$$

**1-1-1-0-0-0**

$\xi_{5,N} \neq 0, \tau_2 \neq 0, \tau_3 \neq 0, \tau_3 = \tau_2, \tau_2 = \tau_4, w_{1,3,2} = 0$

$$w_{1,4,3} \neq -w_{1,4,2}$$

$$\tau_2 = 2\frac{\xi_{5,N}}{2+N}$$

$$n > 1,\ w_{2,4,1} = \frac{w_{1,4,3} + w_{1,4,2}}{\xi_{5,N}} \left( \xi_{5,N}^2 \frac{4}{2+N} - \left( 2 + N - \frac{\xi_{5,N}}{w_{1,4,3} + w_{1,4,2}} \right) w_{2,2,1} - w_{N,5,3} \frac{2(1+N)}{\xi_{5,N}^N} \big(w_{2,3,2} + w_{2,3,1} - w_{2,2,1}\big) \right) - w_{2,4,3} - w_{2,4,2}$$

$$w_{N,5,1} = \xi_{5,N}^N \frac{N}{2(1+N)}$$

$$w_{N,5,2} = \frac{\xi_{5,N}^N}{2(1+N)} \left( 2 + N - \frac{\xi_{5,N}}{w_{1,4,3} + w_{1,4,2}} \right) - w_{N,5,3}$$

$$w_{N,5,4} = \frac{\xi_{5,N}^{1+N}}{2(1+N)\big(w_{1,4,3} + w_{1,4,2}\big)}$$

**1-1-1-0-0-1-0**

$\xi_{5,N} \neq 0, \tau_2 \neq 0, \tau_3 \neq 0, \tau_3 = \tau_2, \tau_2 = \tau_4, w_{1,3,2} \neq 0, w_{N,5,4} = 0$

$$\tau_2 = 2\frac{\xi_{5,N}}{2+N}$$

$$n > 1,\ w_{2,3,1} = 2w_{1,3,2}\left(2\frac{\xi_{5,N}}{2+N} - \frac{1+N}{\xi_{5,N}^{1+N}}w_{N,5,2}w_{2,2,1}\right) - w_{2,3,2}$$

$$w_{N,5,1} = \xi_{5,N}^{N}\frac{N}{2(1+N)}$$

$$w_{N,5,2} = \frac{\xi_{5,N}^{N}}{2(1+N)}\left(2+N-\frac{\xi_{5,N}}{w_{1,3,2}}\right)$$

$$w_{N,5,3} = \frac{\xi_{5,N}^{1+N}}{2w_{1,3,2}(1+N)}$$

**1-1-1-0-0-1-1**

$\xi_{5,N} \neq 0, \tau_2 \neq 0, \tau_3 \neq 0, \tau_3 = \tau_2, \tau_2 = \tau_4, w_{1,3,2} \neq 0, w_{N,5,4} \neq 0$

$$\tau_2 = 2\frac{\xi_{5,N}}{2+N}$$

$$n > 1,\ w_{2,4,1} = \frac{1}{w_{N,5,4}}\left(2\xi_{5,N}^{2+N}\frac{N!}{(2+N)!} - w_{N,5,2}w_{2,2,1} - w_{N,5,3}\left(w_{2,3,2}+w_{2,3,1}\right)\right) - w_{2,4,3} - w_{2,4,2}$$

$$w_{N,5,1} = \xi_{5,N}^{N}\frac{N}{2(1+N)}$$

$$w_{N,5,2} = \frac{\xi_{5,N}^{N}}{2(1+N)}\left(2+N-\frac{\xi_{5,N}}{w_{1,3,2}}\right) + \frac{w_{N,5,4}}{w_{1,3,2}}\left(w_{1,4,3}+w_{1,4,2}-w_{1,3,2}\right)$$

$$w_{N,5,3} = \frac{1}{w_{1,3,2}}\left(\frac{\xi_{5,N}^{1+N}}{2(1+N)} - w_{N,5,4}\left(w_{1,4,3}+w_{1,4,2}\right)\right)$$

**1-1-1-0-1-0-0**

$\xi_{5,N} \neq 0, \tau_2 \neq 0, \tau_3 \neq 0, \tau_3 = \tau_2, \tau_4 \neq \tau_2, \tau_4 = 0, w_{1,3,2} = 0$

$$\tau_2 = 2\frac{\xi_{5,N}}{2+N}$$

$$n > 1,\ w_{2,4,1} = 2\big(w_{1,4,3} + w_{1,4,2}\big)\left(2\frac{\xi_{5,N}}{2+N} - \frac{2+N}{2\xi_{5,N}}w_{2,2,1} - w_{N,5,3}\frac{1+N}{\xi_{5,N}^{1+N}}\big(w_{2,3,2} + w_{2,3,1} - w_{2,2,1}\big)\right) - w_{2,4,3} - w_{2,4,2}$$

$$w_{1,4,3} \neq -w_{1,4,2}$$

$$w_{N,5,1} = \xi_{5,N}^{N}\frac{N}{2(1+N)}\left(1 - \frac{\xi_{5,N}}{2\big(w_{1,4,3} + w_{1,4,2}\big)}\right)$$

$$w_{N,5,2} = \xi_{5,N}^{N}\frac{2+N}{2(1+N)} - w_{N,5,3}$$

$$w_{N,5,4} = \frac{\xi_{5,N}^{1+N}}{2\big(w_{1,4,3} + w_{1,4,2}\big)(1+N)}$$

**1-1-1-0-1-0-1-0**

$\xi_{5,N} \neq 0, \tau_2 \neq 0, \tau_3 \neq 0, \tau_3 = \tau_2, \tau_4 \neq \tau_2, \tau_4 = 0, w_{1,3,2} \neq 0, w_{N,5,4} = 0$

$$\tau_2 = 2\frac{\xi_{5,N}}{2+N}$$

$$n > 1,\ w_{2,3,1} = 4w_{1,3,2}\frac{\xi_{5,N}}{2+N} - w_{2,2,1}\left(\frac{2+N}{\xi_{5,N}}w_{1,3,2} - 1\right) - w_{2,3,2}$$

$$w_{N,5,1} = \xi_{5,N}^{N}\frac{N}{2(1+N)}$$

$$w_{N,5,2} = \frac{\xi_{5,N}^{N}}{2(1+N)}\left(2 + N - \frac{\xi_{5,N}}{w_{1,3,2}}\right)$$

$$w_{N,5,3} = \frac{\xi_{5,N}^{1+N}}{2w_{1,3,2}(1+N)}$$

**1-1-1-0-1-0-1-1**

$\xi_{5,N} \neq 0, \tau_2 \neq 0, \tau_3 \neq 0, \tau_3 = \tau_2, \tau_4 \neq \tau_2, \tau_4 = 0, w_{1,3,2} \neq 0, w_{N,5,4} \neq 0$

$$\tau_2 = 2\frac{\xi_{5,N}}{2+N}$$

$$n > 1,\ w_{2,4,1} = \frac{1}{w_{N,5,4}}\left(2\xi_{5,N}^{2+N}\frac{N!}{(2+N)!} - w_{2,2,1}\xi_{5,N}^{N}\frac{2+N}{2(1+N)}\right.$$
$$\left.-\frac{1}{w_{1,3,2}}\left(\frac{\xi_{5,N}^{1+N}}{2(1+N)} - w_{N,5,4}\big(w_{1,4,3} + w_{1,4,2}\big)\right)\big(w_{2,3,2} + w_{2,3,1} - w_{2,2,1}\big)\right) - w_{2,4,3} - w_{2,4,2}$$

$$w_{N,5,1} = \xi_{5,N}^{N}\frac{N}{2(1+N)} - w_{N,5,4}$$

$$w_{N,5,2} = \xi_{5,N}^{N}\frac{2+N}{2(1+N)} - \frac{1}{w_{1,3,2}}\left(\frac{\xi_{5,N}^{1+N}}{2(1+N)} - w_{N,5,4}\big(w_{1,4,3} + w_{1,4,2}\big)\right)$$

$$w_{N,5,3} = \frac{1}{w_{1,3,2}}\left(\frac{\xi_{5,N}^{1+N}}{2(1+N)} - w_{N,5,4}\big(w_{1,4,3} + w_{1,4,2}\big)\right)$$

**1-1-1-0-1-1-0**

$\xi_{5,N} \neq 0, \tau_2 \neq 0, \tau_3 \neq 0, \tau_3 = \tau_2, \tau_4 \neq \tau_2, \tau_4 \neq 0, w_{1,3,2} = 0$

$$\tau_2 \neq 2\frac{\xi_{5,N}}{2+N}$$

$$w_{1,4,2} = \frac{\xi_{5,N}}{\tau_2}\tau_4\frac{\tau_4 - \tau_2}{2\xi_{5,N} - \tau_2(2+N)} - w_{1,4,3}$$

$$n > 1,\ w_{2,4,1} = \frac{\tau_4\big(2\xi_{5,N} - \tau_4(2+N)\big)}{\tau_2\big(2\xi_{5,N} - \tau_2(2+N)\big)}w_{2,2,1} - w_{2,4,3} - w_{2,4,2}$$
$$+\tau_4\frac{\tau_4 - \tau_2}{2\xi_{5,N} - \tau_2(2+N)}\left(2\xi_{5,N} - w_{N,5,3}(1+N)\frac{2+N}{\xi_{5,N}^{1+N}}\big(w_{2,3,2} + w_{2,3,1} - w_{2,2,1}\big)\right)$$

$$w_{N,5,1} = \xi_{5,N}^{N} - w_{N,5,2} - w_{N,5,3} - w_{N,5,4}$$

$$w_{N,5,2} = \frac{\xi_{5,N}^{1+N}}{\tau_2(\tau_4 - \tau_2)(1+N)}\left(\tau_4 - 2\frac{\xi_{5,N}}{2+N}\right) - w_{N,5,3}$$

$$w_{N,5,4} = \frac{\xi_{5,N}^{1+N}}{\tau_4(\tau_4 - \tau_2)(1+N)}\left(2\frac{\xi_{5,N}}{2+N} - \tau_2\right)$$

**1-1-1-0-1-1-1-0**

$\xi_{5,N} \neq 0, \tau_2 \neq 0, \tau_3 \neq 0, \tau_3 = \tau_2, \tau_4 \neq \tau_2, \tau_4 \neq 0, w_{1,3,2} \neq 0, w_{N,5,4} = 0$

$$\tau_2 = 2\frac{\xi_{5,N}}{2+N}$$

$$n > 1,\ w_{2,3,1} = 4\xi_{5,N}\frac{w_{1,3,2}}{2+N} - \left(w_{1,3,2}\frac{2+N}{\xi_{5,N}} - 1\right)w_{2,2,1} - w_{2,3,2}$$

$$w_{N,5,1} = \xi_{5,N}^{N} - \frac{\xi_{5,N}^{1+N}}{(\tau_4-\tau_2)(1+N)(2+N)}\left(\frac{\tau_4(2+N)-2\xi_{5,N}}{\tau_2} + \frac{2\xi_{5,N}-\tau_2(2+N)}{\tau_4}\right)$$

$$w_{N,5,2} = \frac{\xi_{5,N}^{N}}{2(1+N)}\left(2+N-\frac{\xi_{5,N}}{w_{1,3,2}}\right)$$

$$w_{N,5,3} = \frac{\xi_{5,N}^{1+N}}{w_{1,3,2}2(1+N)}$$

**1-1-1-0-1-1-1-1**

$\xi_{5,N} \neq 0, \tau_2 \neq 0, \tau_3 \neq 0, \tau_3 = \tau_2, \tau_4 \neq \tau_2, \tau_4 \neq 0, w_{1,3,2} \neq 0, w_{N,5,4} \neq 0$

$$\tau_2 \neq 2\frac{\xi_{5,N}}{2+N}$$

$$n > 1,\ w_{2,4,1} = \frac{\tau_4(\tau_4-\tau_2)}{2\xi_{5,N}-\tau_2(2+N)}\left(2\xi_{5,N} - w_{2,2,1}\frac{\tau_4(2+N)-2\xi_{5,N}}{\tau_2(\tau_4-\tau_2)}\right.$$
$$\left. - \left(\frac{\xi_{5,N}}{\tau_2} - \frac{2\xi_{5,N}-\tau_2(2+N)}{\tau_4(\tau_4-\tau_2)}(w_{1,4,3}+w_{1,4,2})\right)\frac{w_{2,3,2}+w_{2,3,1}-w_{2,2,1}}{w_{1,3,2}}\right) - w_{2,4,3} - w_{2,4,2}$$

$$w_{N,5,1} = \xi_{5,N}^{N} - \frac{\xi_{5,N}^{1+N}}{(\tau_4-\tau_2)(1+N)(2+N)}\left(\frac{\tau_4(2+N)-2\xi_{5,N}}{\tau_2} + \frac{2\xi_{5,N}-\tau_2(2+N)}{\tau_4}\right)$$

$$w_{N,5,2} = \frac{\xi_{5,N}^{1+N}}{(1+N)(2+N)}\left(\frac{\tau_4(2+N)-2\xi_{5,N}}{\tau_2(\tau_4-\tau_2)} - \frac{1}{w_{1,3,2}}\left(\frac{\xi_{5,N}}{\tau_2} - \frac{2\xi_{5,N}-\tau_2(2+N)}{\tau_4(\tau_4-\tau_2)}(w_{1,4,3}+w_{1,4,2})\right)\right)$$

$$w_{N,5,3} = \frac{\xi_{5,N}^{1+N}}{w_{1,3,2}(1+N)(2+N)}\left(\frac{\xi_{5,N}}{\tau_2} - \frac{2\xi_{5,N}-\tau_2(2+N)}{\tau_4(\tau_4-\tau_2)}(w_{1,4,3}+w_{1,4,2})\right)$$

$$w_{N,5,4} = \frac{\xi_{5,N}^{1+N}\left(2\xi_{5,N}-\tau_2(2+N)\right)}{\tau_4(\tau_4-\tau_2)(1+N)(2+N)}$$

**1-1-1-1-0-0**

$\xi_{5,N} \neq 0, \tau_2 \neq 0, \tau_3 \neq 0, \tau_3 \neq \tau_2, w_{1,4,2} = w_{1,3,2}\frac{\tau_4(\tau_4-\tau_2)}{\tau_3(\tau_3-\tau_2)} - w_{1,4,3}\frac{\tau_3}{\tau_2}, w_{N,5,4} = 0$

$$\tau_2 \neq 2\frac{\xi_{5,N}}{2+N}$$

$$w_{1,3,2} = \xi_{5,N}\frac{\tau_3(\tau_3-\tau_2)}{\tau_2\big(2\xi_{5,N}-\tau_2(2+N)\big)}$$

$$n>1,\ w_{2,3,1} = \frac{\tau_3}{2\xi_{5,N}-\tau_2(2+N)}\left(2\xi_{5,N}(\tau_3-\tau_2) + \frac{2\xi_{5,N}-\tau_3(2+N)}{\tau_2}w_{2,2,1}\right) - w_{2,3,2}$$

$$w_{N,5,1} = \xi_{5,N}^N - \frac{\xi_{5,N}^{1+N}}{(\tau_3-\tau_2)(1+N)(2+N)}\left(\frac{\tau_3(2+N)-2\xi_{5,N}}{\tau_2} + \frac{2\xi_{5,N}-\tau_2(2+N)}{\tau_3}\right)$$

$$w_{N,5,2} = \frac{\xi_{5,N}^{1+N}}{\tau_2(1+N)(\tau_3-\tau_2)}\left(\tau_3 - 2\frac{\xi_{5,N}}{2+N}\right)$$

$$w_{N,5,3} = \xi_{5,N}^{1+N}\frac{2\xi_{5,N}-\tau_2(2+N)}{\tau_3(\tau_3-\tau_2)(1+N)(2+N)}$$

**1-1-1-1-0-1**

$\xi_{5,N} \neq 0, \tau_2 \neq 0, \tau_3 \neq 0, \tau_3 \neq \tau_2, w_{1,4,2} = w_{1,3,2}\frac{\tau_4(\tau_4-\tau_2)}{\tau_3(\tau_3-\tau_2)} - w_{1,4,3}\frac{\tau_3}{\tau_2}, w_{N,5,4} \neq 0$

$$w_{1,3,2} = \xi_{5,N}\frac{\tau_3(\tau_3-\tau_2)}{\tau_2\big(2\xi_{5,N}-\tau_2(2+N)\big)}$$

$$\begin{aligned} n>1,\ w_{2,4,1} = \frac{1}{w_{N,5,4}}\Bigg( & \\ & 2\xi_{5,N}^{2+N}\frac{N!}{(2+N)!} - \frac{\xi_{5,N}^{1+N}}{\tau_2(1+N)(\tau_3-\tau_2)}\left(\tau_3 - 2\frac{\xi_{5,N}}{2+N}\right) + w_{N,5,4}\frac{\tau_4(\tau_4-\tau_3)}{\tau_2(\tau_3-\tau_2)}w_{2,2,1} \\ & -\frac{1}{\tau_3(\tau_3-\tau_2)}\left(\frac{\xi_{5,N}^{1+N}}{1+N}\left(2\frac{\xi_{5,N}}{2+N}-\tau_2\right) - w_{N,5,4}\tau_4(\tau_4-\tau_2)\right)(w_{2,3,2}+w_{2,3,1}) \\ & \Bigg) - w_{2,4,3} - w_{2,4,2} \end{aligned}$$

$$w_{N,5,1} = \xi_{5,N}^N - w_{N,5,2} - w_{N,5,3} - w_{N,5,4}$$

$$w_{N,5,2} = \frac{\xi_{5,N}^{1+N}}{\tau_2(1+N)(\tau_3-\tau_2)}\left(\tau_3 - 2\frac{\xi_{5,N}}{2+N}\right) + w_{N,5,4}\frac{\tau_4(\tau_4-\tau_3)}{\tau_2(\tau_3-\tau_2)}$$

$$w_{N,5,3} = \frac{1}{\tau_3(\tau_3-\tau_2)}\left(\frac{\xi_{5,N}^{1+N}}{1+N}\left(2\frac{\xi_{5,N}}{2+N}-\tau_2\right) - w_{N,5,4}\tau_4(\tau_4-\tau_2)\right)$$

**1-1-1-1-1-0**

$\xi_{5,N} \neq 0, \tau_2 \neq 0, \tau_3 \neq 0, \tau_3 \neq \tau_2, w_{1,4,2} \neq w_{1,3,2}\frac{\tau_4(\tau_4-\tau_2)}{\tau_3(\tau_3-\tau_2)} - w_{1,4,3}\frac{\tau_3}{\tau_2}, w_{N,5,4} = 0$

$$\tau_2 \neq 0$$

$$\tau_2 \neq 2\frac{\xi_{5,N}}{2+N}$$

$$w_{1,3,2} = \frac{\xi_{5,N}\tau_3(\tau_3-\tau_2)}{\tau_2\big(2\xi_{5,N}-\tau_2(2+N)\big)}$$

$$n>1,\ w_{2,3,1} = \frac{\tau_3(\tau_3-\tau_2)}{2\xi_{5,N}-\tau_2(2+N)}\left(2\xi_{5,N} - \frac{\tau_3(2+N)-2\xi_{5,N}}{\tau_2(\tau_3-\tau_2)}w_{2,2,1}\right) - w_{2,3,2}$$

$$w_{N,5,1} = \xi_{5,N}^{N} - \frac{\xi_{5,N}^{1+N}}{\tau_2\tau_3(1+N)(2+N)}\big((\tau_3+\tau_2)(2+N)-2\xi_{5,N}\big)$$

$$w_{N,5,2} = \frac{\xi_{5,N}^{1+N}}{(1+N)(2+N)\tau_2(\tau_3-\tau_2)}\big(\tau_3(2+N)-2\xi_{5,N}\big)$$

$$w_{N,5,3} = \frac{\xi_{5,N}^{1+N}}{\tau_3(\tau_3-\tau_2)(1+N)(2+N)}\big(2\xi_{5,N}-\tau_2(2+N)\big)$$

**1-1-1-1-1-1**

$\xi_{5,N} \neq 0, \tau_2 \neq 0, \tau_3 \neq 0, \tau_3 \neq \tau_2, w_{1,4,2} \neq w_{1,3,2}\frac{\tau_4(\tau_4-\tau_2)}{\tau_3(\tau_3-\tau_2)} - w_{1,4,3}\frac{\tau_3}{\tau_2}, w_{N,5,4} \neq 0$

$$a = w_{1,4,3}\tau_3^2(\tau_3-\tau_2) + \tau_2\big(w_{1,4,2}\tau_3(\tau_3-\tau_2) - w_{1,3,2}\tau_4(\tau_4-\tau_2)\big)$$

$$b = \xi_{5,N}\tau_3(\tau_3-\tau_2) - w_{1,3,2}\tau_2\big(2\xi_{5,N}-\tau_2(2+N)\big) \neq 0$$

$$n>1,\ w_{2,4,1} = \frac{a}{b}\left(2\xi_{5,N} - \frac{\tau_3(2+N)-2\xi_{5,N}}{\tau_2(\tau_3-\tau_2)}w_{2,2,1} - \frac{2\xi_{5,N}-\tau_2(2+N)}{\tau_3(\tau_3-\tau_2)}\big(w_{2,3,2}+w_{2,3,1}\big)\right)$$

$$+\frac{\tau_4(\tau_4-\tau_2)}{\tau_3(\tau_3-\tau_2)}\big(w_{2,3,2}+w_{2,3,1}-w_{2,2,1}\big) - w_{2,4,3} - w_{2,4,2}$$

$$w_{N,5,1} = \xi_{5,N}^{N} - \frac{\xi_{5,N}^{1+N}}{\tau_2\tau_3(1+N)(2+N)}\left((\tau_3+\tau_2)(2+N) - 2\xi_{5,N} + (\tau_4-\tau_2)(\tau_4-\tau_3)\frac{b}{a}\right)$$

$$w_{N,5,2} = \frac{\xi_{5,N}^{1+N}}{(1+N)(2+N)\tau_2(\tau_3-\tau_2)}\left(\tau_3(2+N) - 2\xi_{5,N} + \frac{b}{a}\tau_4(\tau_4-\tau_3)\right)$$

$$w_{N,5,3} = \frac{\xi_{5,N}^{1+N}}{\tau_3(\tau_3-\tau_2)(1+N)(2+N)}\left(2\xi_{5,N} - \tau_2(2+N) - \frac{b}{a}\tau_4(\tau_4-\tau_2)\right)$$

$$w_{N,5,4} = \xi_{5,N}^{1+N}\frac{b}{(1+N)(2+N)a}$$

**C.6 Explicit**, $s=4$, $j=5$, $v=4$

The conditions are :

Order of consistency 1 :

$$w_{N,5,1}+w_{N,5,2}+w_{N,5,3}+w_{N,5,4}=\xi_{5,N}^{N}$$

Order of consistency 2 :

$$w_{N,5,1}\tau_1+w_{N,5,2}\tau_2+w_{N,5,3}\tau_3+w_{N,5,4}\tau_4=\frac{\xi_{5,N}^{N+1}}{1+N}$$

$$w_{N,5,2}w_{1,2,1}+w_{N,5,3}\left(w_{1,3,1}+w_{1,3,2}\right)+w_{N,5,4}\left(w_{1,4,1}+w_{1,4,2}+w_{1,4,3}\right)=\frac{\xi_{5,N}^{N+1}}{1+N}$$

Order of consistency 3 :

$$w_{N,5,1}\tau_1^2+w_{N,5,2}\tau_2^2+w_{N,5,3}\tau_3^2+w_{N,5,4}\tau_4^2=\xi_{5,N}^{N+2}2\frac{N!}{(2+N)!}$$

$$w_{N,5,2}\tau_2w_{1,2,1}+w_{N,5,3}\tau_3\left(w_{1,3,1}+w_{1,3,2}\right)+w_{N,5,4}\tau_4\left(w_{1,4,1}+w_{1,4,2}+w_{1,4,3}\right)=\xi_{5,N}^{N+2}2\frac{N!}{(2+N)!}$$

$$w_{N,5,2}w_{1,2,1}\tau_1+w_{N,5,3}\left(w_{1,3,1}\tau_1+w_{1,3,2}\tau_2\right)+w_{N,5,4}\left(w_{1,4,1}\tau_1+w_{1,4,2}\tau_2+w_{1,4,3}\tau_3\right)=\xi_{5,N}^{N+2}\frac{N!}{(2+N)!}$$

$$w_{N,5,3}w_{1,3,2}w_{1,2,1}+w_{N,5,4}\left(w_{1,4,2}w_{1,2,1}+w_{1,4,3}\left(w_{1,3,1}+w_{1,3,2}\right)\right)=\xi_{5,N}^{N+2}\frac{N!}{(2+N)!}$$

$$w_{N,5,2}w_{1,2,1}^2+w_{N,5,3}\left(w_{1,3,1}+w_{1,3,2}\right)^2+w_{N,5,4}\left(w_{1,4,1}+w_{1,4,2}+w_{1,4,3}\right)^2=\xi_{5,N}^{N+2}2\frac{N!}{(2+N)!}$$

$$n>1,\;\; w_{N,5,2}w_{2,2,1}+w_{N,5,3}\left(w_{2,3,1}+w_{2,3,2}\right)+w_{N,5,4}\left(w_{2,4,1}+w_{2,4,2}+w_{2,4,3}\right)=\xi_{5,N}^{N+2}2\frac{N!}{(2+N)!}$$

Order of consistency 4 :

$$w_{N,5,1}\tau_1^3+w_{N,5,2}\tau_2^3+w_{N,5,3}\tau_3^3+w_{N,5,4}\tau_4^3=\xi_{5,N}^{N+3}6\frac{N!}{(3+N)!}$$

$$w_{N,5,2}w_{1,2,1}\tau_1^2+w_{N,5,3}\left(w_{1,3,1}\tau_1^2+w_{1,3,2}\tau_2^2\right)+w_{N,5,4}\left(w_{1,4,1}\tau_1^2+w_{1,4,2}\tau_2^2+w_{1,4,3}\tau_3^2\right)=\xi_{5,N}^{N+3}2\frac{N!}{(3+N)!}$$

$$w_{N,5,2}\tau_2w_{1,2,1}\tau_1+w_{N,5,3}\tau_3\left(w_{1,3,1}\tau_1+w_{1,3,2}\tau_2\right)+w_{N,5,4}\tau_4\left(w_{1,4,1}\tau_1+w_{1,4,2}\tau_2+w_{1,4,3}\tau_3\right)=\xi_{5,N}^{N+3}3\frac{N!}{(3+N)!}$$

$$w_{N,5,2}\tau_2^2w_{1,2,1}+w_{N,5,3}\tau_3^2\left(w_{1,3,1}+w_{1,3,2}\right)+w_{N,5,4}\tau_4^2\left(w_{1,4,1}+w_{1,4,2}+w_{1,4,3}\right)=\xi_{5,N}^{N+3}6\frac{N!}{(3+N)!}$$

$$w_{N,5,3}w_{1,3,2}w_{1,2,1}\tau_1+w_{N,5,4}\left(w_{1,4,2}w_{1,2,1}\tau_1+w_{1,4,3}\left(w_{1,3,1}\tau_1+w_{1,3,2}\tau_2\right)\right)=\xi_{5,N}^{N+3}\frac{N!}{(3+N)!}$$

$$w_{N,5,3}w_{1,3,2}\tau_2w_{1,2,1}+w_{N,5,4}\left(w_{1,4,2}\tau_2w_{1,2,1}+w_{1,4,3}\tau_3\left(w_{1,3,1}+w_{1,3,2}\right)\right)=\xi_{5,N}^{N+3}2\frac{N!}{(3+N)!}$$

$$w_{N,5,3}\tau_3w_{1,3,2}w_{1,2,1}+w_{N,5,4}\tau_4\left(w_{1,4,2}w_{1,2,1}+w_{1,4,3}\left(w_{1,3,1}+w_{1,3,2}\right)\right)=\xi_{5,N}^{N+3}3\frac{N!}{(3+N)!}$$

$$w_{N,5,2}w_{1,2,1}w_{1,2,1}\tau_1 + w_{N,5,3}\left(w_{1,3,1} + w_{1,3,2}\right)\left(w_{1,3,1}\tau_1 + w_{1,3,2}\tau_2\right)$$

$$+w_{N,5,4}\left(w_{1,4,1} + w_{1,4,2} + w_{1,4,3}\right)\left(w_{1,4,1}\tau_1 + w_{1,4,2}\tau_2 + w_{1,4,3}\tau_3\right) = \xi_{5,N}^{N+3}3\frac{N!}{(3+N)!}$$

$$w_{N,5,2}\tau_2 w_{1,2,1}^2 + w_{N,5,3}\tau_3\left(w_{1,3,1} + w_{1,3,2}\right)^2 + w_{N,5,4}\tau_4\left(w_{1,4,1} + w_{1,4,2} + w_{1,4,3}\right)^2 = \xi_{5,N}^{N+3}6\frac{N!}{(3+N)!}$$

$$w_{N,5,4}w_{1,4,3}w_{1,3,2}w_{1,2,1} = \xi_{5,N}^{N+3}\frac{N!}{(3+N)!}$$

$$w_{N,5,3}w_{1,3,2}w_{1,2,1}^2 + w_{N,5,4}\left(w_{1,4,2}w_{1,2,1}^2 + w_{1,4,3}\left(w_{1,3,1} + w_{1,3,2}\right)^2\right) = \xi_{5,N}^{N+3}2\frac{N!}{(3+N)!}$$

$$w_{N,5,3}\left(w_{1,3,1} + w_{1,3,2}\right)w_{1,3,2}w_{1,2,1}$$

$$+w_{N,5,4}\left(w_{1,4,1} + w_{1,4,2} + w_{1,4,3}\right)\left(w_{1,4,2}w_{1,2,1} + w_{1,4,3}\left(w_{1,3,1} + w_{1,3,2}\right)\right) = \xi_{5,N}^{N+3}3\frac{N!}{(3+N)!}$$

$$w_{N,5,2}w_{1,2,1}^3 + w_{N,5,3}\left(w_{1,3,1} + w_{1,3,2}\right)^3 + w_{N,5,4}\left(w_{1,4,1} + w_{1,4,2} + w_{1,4,3}\right)^3 = \xi_{5,N}^{N+3}6\frac{N!}{(3+N)!}$$

$$n > 1,\ w_{N,5,2}w_{2,2,1}\tau_1 + w_{N,5,3}\left(w_{2,3,1}\tau_1 + w_{2,3,2}\tau_2\right) + w_{N,5,4}\left(w_{2,4,1}\tau_1 + w_{2,4,2}\tau_2 + w_{2,4,3}\tau_3\right) = \xi_{5,N}^{N+3}2\frac{N!}{(3+N)!}$$

$$n > 1,\ w_{N,5,2}w_{1,2,1}w_{2,2,1} + w_{N,5,3}\left(w_{1,3,1} + w_{1,3,2}\right)\left(w_{2,3,1} + w_{2,3,2}\right)$$

$$+w_{N,5,4}\left(w_{1,4,1} + w_{1,4,2} + w_{1,4,3}\right)\left(w_{2,4,1} + w_{2,4,2} + w_{2,4,3}\right) = \xi_{5,N}^{N+3}6\frac{N!}{(3+N)!}$$

$$n > 1,\ w_{N,5,3}w_{2,3,2}w_{1,2,1} + w_{N,5,4}\left(w_{2,4,2}w_{1,2,1} + w_{2,4,3}\left(w_{1,3,1} + w_{1,3,2}\right)\right) = \xi_{5,N}^{N+3}2\frac{N!}{(3+N)!}$$

$$n > 1,\ w_{N,5,3}w_{1,3,2}w_{2,2,1} + w_{N,5,4}\left(w_{1,4,2}w_{2,2,1} + w_{1,4,3}\left(w_{2,3,1} + w_{2,3,2}\right)\right) = \xi_{5,N}^{N+3}2\frac{N!}{(3+N)!}$$

$$n > 1,\ w_{N,5,2}\tau_2 w_{2,2,1} + w_{N,5,3}\tau_3\left(w_{2,3,1} + w_{2,3,2}\right) + w_{N,5,4}\tau_4\left(w_{2,4,1} + w_{2,4,2} + w_{2,4,3}\right) = \xi_{5,N}^{N+3}6\frac{N!}{(3+N)!}$$

$$n > 2,\ w_{N,5,2}w_{3,2,1} + w_{N,5,3}\left(w_{3,3,1} + w_{3,3,2}\right) + w_{N,5,4}\left(w_{3,4,1} + w_{3,4,2} + w_{3,4,3}\right) = \xi_{5,N}^{N+3}6\frac{N!}{(3+N)!}$$

We now list the main operations performed to rewrite this system.

1) isolate $w_{N,5,1}\tau_1^2$ in 3.1
2) 3.1 in 4.1
3) isolate $w_{N,5,2}\tau_2 w_{1,2,1}$ in 3.2
4) 3.2 in 4.4
5) isolate $w_{N,5,2}w_{1,2,1}\tau_1$ in 3.3
6) 3.3 in 4.2, 4.3
7) isolate $w_{N,5,2}w_{1,2,1}^2$ in 3.5
8) 3.5 in 4.8, 4.9, 4.13
9) isolate $w_{N,5,2}w_{2,2,1}$ in 3.6
10) 3.6 in 4.14, 4.15, 4.18
11) isolate $w_{N,5,3}w_{1,3,2}w_{1,2,1}$ in 3.4
12) 3.4 in 4.5, 4.6, 4.7, 4.11, 4.12
13) isolate $w_{N,5,1}\tau_1$ in 2.1
14) 2.1 in 3.1
15) isolate $w_{N,5,2}w_{1,2,1}$ in 2.2
16) 2.2 in 3.2, 3.3, 3.5
17) isolate $w_{N,5,2}\tau_2(\tau_2-\tau_1)$ in 3.1
18) 3.1 in 4.1
19) isolate $w_{N,5,3}(\tau_3-\tau_2)\big(w_{1,3,1}+w_{1,3,2}\big)$ in 3.2
20) 3.2 in 4.3, 4.4, 4.9
21) isolate $w_{N,5,3}w_{1,3,2}(\tau_2-\tau_1)$ in 3.3
22) 3.3 in 4.2, 4.3, 4.8
23) isolate $w_{N,5,3}\big(w_{1,3,1}+w_{1,3,2}\big)\big(w_{1,3,1}+w_{1,3,2}-w_{1,2,1}\big)$ in 3.5
24) 3.5 in 4.13
25) isolate $w_{N,5,1}$ in 1.1
26) 1.1 in 2.1
27) isolate $w_{N,5,2}(\tau_2-\tau_1)$ in 2.1
28) 2.1 in 3.1
29) isolate $w_{N,5,3}(\tau_3-\tau_2)(\tau_3-\tau_1)$ in 3.1
30) 3.1 in 4.1

Order of consistency 1 :

$$w_{N,5,1}=\xi_{5,N}^{N}-w_{N,5,2}-w_{N,5,3}-w_{N,5,4}$$

Order of consistency 2 :

$$w_{N,5,2}(\tau_2-\tau_1)=\xi_{5,N}^{N}\left(\frac{\xi_{5,N}}{1+N}-\tau_1\right)-w_{N,5,3}(\tau_3-\tau_1)-w_{N,5,4}(\tau_4-\tau_1)$$

$$w_{N,5,2}w_{1,2,1}=\frac{\xi_{5,N}^{N+1}}{1+N}-w_{N,5,3}\big(w_{1,3,1}+w_{1,3,2}\big)-w_{N,5,4}\big(w_{1,4,1}+w_{1,4,2}+w_{1,4,3}\big)$$

Order of consistency 3 :

$$w_{N,5,3}(\tau_3-\tau_2)(\tau_3-\tau_1)=\xi_{5,N}^{N}\left(\frac{\xi_{5,N}}{1+N}\left(\xi_{5,N}\frac{2}{2+N}-\tau_1\right)-\left(\frac{\xi_{5,N}}{1+N}-\tau_1\right)\tau_2\right)-w_{N,5,4}(\tau_4-\tau_2)(\tau_4-\tau_1)$$

$$w_{N,5,3}(\tau_3-\tau_2)(w_{1,3,1}+w_{1,3,2}) = -w_{N,5,4}(\tau_4-\tau_2)(w_{1,4,1}+w_{1,4,2}+w_{1,4,3}) + \frac{\xi_{5,N}^{N+1}}{1+N}\left(\xi_{5,N}\frac{2}{2+N}-\tau_2\right)$$

$$w_{N,5,3}w_{1,3,2}(\tau_2-\tau_1) = -w_{N,5,4}\big(w_{1,4,2}(\tau_2-\tau_1)+w_{1,4,3}(\tau_3-\tau_1)\big) + \frac{\xi_{5,N}^{N+1}}{1+N}\left(\frac{\xi_{5,N}}{2+N}-\tau_1\right)$$

$$w_{N,5,3}w_{1,3,2}w_{1,2,1} = -w_{N,5,4}\big(w_{1,4,2}w_{1,2,1}+w_{1,4,3}(w_{1,3,1}+w_{1,3,2})\big) + \xi_{5,N}^{N+2}\frac{N!}{(2+N)!}$$

$$w_{N,5,3}(w_{1,3,1}+w_{1,3,2})(w_{1,3,1}+w_{1,3,2}-w_{1,2,1}) = -w_{N,5,4}(w_{1,4,1}+w_{1,4,2}+w_{1,4,3})(w_{1,4,1}+w_{1,4,2}+w_{1,4,3}-w_{1,2,1}) + \frac{\xi_{5,N}^{N+1}}{1+N}\left(\xi_{5,N}\frac{2}{2+N}-w_{1,2,1}\right)$$

$$n>1,\ w_{N,5,2}w_{2,2,1} = -w_{N,5,3}(w_{2,3,1}+w_{2,3,2}) - w_{N,5,4}(w_{2,4,1}+w_{2,4,2}+w_{2,4,3}) + \xi_{5,N}^{N+2}2\frac{N!}{(2+N)!}$$

Order of consistency 4 :

$$w_{N,5,4}(\tau_4-\tau_3)(\tau_4-\tau_2)(\tau_4-\tau_1)$$
$$= \xi_{5,N}^{N}\left(\frac{\xi_{5,N}}{1+N}\left(\xi_{5,N}\frac{2}{2+N}\left(\xi_{5,N}\frac{3}{3+N}-\tau_1\right)-\left(\xi_{5,N}\frac{2}{2+N}-\tau_1\right)\tau_2\right)-\left(\frac{\xi_{5,N}}{1+N}\left(\xi_{5,N}\frac{2}{2+N}-\tau_1\right)-\left(\frac{\xi_{5,N}}{1+N}-\tau_1\right)\tau_2\right)\tau_3\right)$$

$$w_{N,5,4}w_{1,4,3}(\tau_3-\tau_2)(\tau_3-\tau_1) = \frac{\xi_{5,N}^{N+1}}{1+N}\left(\frac{\xi_{5,N}}{2+N}\left(\xi_{5,N}\frac{2}{3+N}-\tau_1\right)-\left(\frac{\xi_{5,N}}{2+N}-\tau_1\right)\tau_2\right)$$

$$w_{N,5,4}(\tau_4-\tau_3)\big(w_{1,4,2}(\tau_2-\tau_1)+w_{1,4,3}(\tau_3-\tau_1)\big) = \frac{\xi_{5,N}^{N+1}}{1+N}\left(\frac{\xi_{5,N}}{2+N}\left(\xi_{5,N}\frac{3}{3+N}-2\tau_1\right)-\left(\frac{\xi_{5,N}}{2+N}-\tau_1\right)\tau_3\right)$$

$$w_{N,5,4}(\tau_4-\tau_3)(\tau_4-\tau_2)(w_{1,4,1}+w_{1,4,2}+w_{1,4,3}) = \frac{\xi_{5,N}^{N+1}}{1+N}\left(\xi_{5,N}\frac{2}{2+N}\left(\xi_{5,N}\frac{3}{3+N}-\tau_2\right)-\left(\xi_{5,N}\frac{2}{2+N}-\tau_2\right)\tau_3\right)$$

$$w_{N,5,4}w_{1,4,3}w_{1,3,2}(\tau_2-\tau_1) = \xi_{5,N}^{N+2}\frac{N!}{(2+N)!}\left(\frac{\xi_{5,N}}{3+N}-\tau_1\right)$$

$$w_{N,5,4}w_{1,4,3}(\tau_3-\tau_2)(w_{1,3,1}+w_{1,3,2}) = \xi_{5,N}^{N+2}\frac{N!}{(2+N)!}\left(\xi_{5,N}\frac{2}{3+N}-\tau_2\right)$$

$$w_{N,5,4}(\tau_4-\tau_3)\big(w_{1,4,2}w_{1,2,1}+w_{1,4,3}(w_{1,3,1}+w_{1,3,2})\big) = \xi_{5,N}^{N+2}\frac{N!}{(2+N)!}\left(\xi_{5,N}\frac{3}{3+N}-\tau_3\right)$$

$$w_{N,5,4}\big(w_{1,4,1}+w_{1,4,2}+w_{1,4,3}-(w_{1,3,1}+w_{1,3,2})\big)\big(w_{1,4,2}(\tau_2-\tau_1)+w_{1,4,3}(\tau_3-\tau_1)\big)$$
$$=$$
$$\frac{\xi_{5,N}^{N+1}}{1+N}\left(\frac{\xi_{5,N}}{2+N}\left(\xi_{5,N}\frac{3}{3+N}-2\tau_1\right)-\left(\frac{\xi_{5,N}}{2+N}-\tau_1\right)(w_{1,3,1}+w_{1,3,2})\right)$$

$$w_{N,5,4}(\tau_4-\tau_2)(w_{1,4,1}+w_{1,4,2}+w_{1,4,3})\big(w_{1,4,1}+w_{1,4,2}+w_{1,4,3}-(w_{1,3,1}+w_{1,3,2})\big)$$
$$=$$
$$\frac{\xi_{5,N}^{N+1}}{1+N}\left(\xi_{5,N}\frac{2}{2+N}\left(\xi_{5,N}\frac{3}{3+N}-\tau_2\right)-\left(\xi_{5,N}\frac{2}{2+N}-\tau_2\right)(w_{1,3,1}+w_{1,3,2})\right)$$

$$w_{N,5,4}w_{1,4,3}w_{1,3,2}w_{1,2,1} = \xi_{5,N}^{N+3}\frac{N!}{(3+N)!}$$

$$w_{N,5,4}w_{1,4,3}\left(w_{1,3,1}+w_{1,3,2}\right)\left(w_{1,3,1}+w_{1,3,2}-w_{1,2,1}\right)=\xi_{5,N}^{N+2}\frac{N!}{(2+N)!}\left(\xi_{5,N}\frac{2}{3+N}-w_{1,2,1}\right)$$

$$w_{N,5,4}\left(w_{1,4,1}+w_{1,4,2}+w_{1,4,3}-\left(w_{1,3,1}+w_{1,3,2}\right)\right)\left(w_{1,4,2}w_{1,2,1}+w_{1,4,3}\left(w_{1,3,1}+w_{1,3,2}\right)\right)$$
$$=\xi_{5,N}^{N+2}\frac{N!}{(2+N)!}\left(\xi_{5,N}\frac{3}{3+N}-\left(w_{1,3,1}+w_{1,3,2}\right)\right)$$

$$w_{N,5,4}\left(w_{1,4,1}+w_{1,4,2}+w_{1,4,3}\right)\left(w_{1,4,1}+w_{1,4,2}+w_{1,4,3}-\left(w_{1,3,1}+w_{1,3,2}\right)\right)\left(w_{1,4,1}+w_{1,4,2}+w_{1,4,3}-w_{1,2,1}\right)$$
$$=\frac{\xi_{5,N}^{N+1}}{1+N}\left(\xi_{5,N}\frac{2}{2+N}\left(\xi_{5,N}\frac{3}{3+N}-w_{1,2,1}\right)-\left(\xi_{5,N}\frac{2}{2+N}-w_{1,2,1}\right)\left(w_{1,3,1}+w_{1,3,2}\right)\right)$$

$$n>1,\ w_{N,5,3}w_{2,3,2}(\tau_2-\tau_1)+w_{N,5,4}\left(w_{2,4,2}(\tau_2-\tau_1)+w_{2,4,3}(\tau_3-\tau_1)\right)=\xi_{5,N}^{N+3}2\frac{N!}{(3+N)!}-\xi_{5,N}^{N+2}2\frac{N!}{(2+N)!}\tau_1$$

$$n>1,\ w_{N,5,3}\left(w_{1,3,1}+w_{1,3,2}-w_{1,2,1}\right)\left(w_{2,3,1}+w_{2,3,2}\right)$$
$$+w_{N,5,4}\left(w_{1,4,1}+w_{1,4,2}+w_{1,4,3}-w_{1,2,1}\right)\left(w_{2,4,1}+w_{2,4,2}+w_{2,4,3}\right)=\xi_{5,N}^{N+3}6\frac{N!}{(3+N)!}-\xi_{5,N}^{N+2}2\frac{N!}{(2+N)!}w_{1,2,1}$$

$$n>1,\ w_{N,5,3}w_{2,3,2}w_{1,2,1}+w_{N,5,4}\left(w_{2,4,2}w_{1,2,1}+w_{2,4,3}\left(w_{1,3,1}+w_{1,3,2}\right)\right)=\xi_{5,N}^{N+3}2\frac{N!}{(3+N)!}$$

$$n>1,\ w_{N,5,3}w_{1,3,2}w_{2,2,1}+w_{N,5,4}\left(w_{1,4,2}w_{2,2,1}+w_{1,4,3}\left(w_{2,3,1}+w_{2,3,2}\right)\right)=\xi_{5,N}^{N+3}2\frac{N!}{(3+N)!}$$

$$n>1,\ w_{N,5,3}(\tau_3-\tau_2)\left(w_{2,3,1}+w_{2,3,2}\right)+w_{N,5,4}(\tau_4-\tau_2)\left(w_{2,4,1}+w_{2,4,2}+w_{2,4,3}\right)=\xi_{5,N}^{N+3}6\frac{N!}{(3+N)!}-\xi_{5,N}^{N+2}2\frac{N!}{(2+N)!}\tau_2$$

$$n>2,\ w_{N,5,2}w_{3,2,1}+w_{N,5,3}\left(w_{3,3,1}+w_{3,3,2}\right)+w_{N,5,4}\left(w_{3,4,1}+w_{3,4,2}+w_{3,4,3}\right)=\xi_{5,N}^{N+3}6\frac{N!}{(3+N)!}$$

The solutions are :

**0-0-0-0-0-0-0-0**

$\xi_{5,N}=0,\tau_4=\tau_1,\tau_1=\tau_3,\tau_1=\tau_2,\tau_2=0,w_{N,5,4}=0,w_{N,5,3}=0,w_{N,5,2}=0$

$$w_{N,5,1}=0$$

**0-0-0-0-0-0-0-1**

$\xi_{5,N}=0,\tau_4=\tau_1,\tau_1=\tau_3,\tau_1=\tau_2,\tau_2=0,w_{N,5,4}=0,w_{N,5,3}=0,w_{N,5,2}\neq 0$

$$w_{1,2,1}=0$$
$$n>1,\ w_{2,2,1}=0$$
$$n>2,\ w_{3,2,1}=0$$
$$w_{N,5,1}=-w_{N,5,2}$$

**0-0-0-0-0-0-1-0-0**

$\xi_{5,N} = 0, \tau_4 = \tau_1, \tau_1 = \tau_3, \tau_1 = \tau_2, \tau_2 = 0, w_{N,5,4} = 0, w_{N,5,3} \neq 0, w_{1,2,1} = 0, w_{1,3,2} = 0$

$$w_{1,3,1} = 0$$
$$n > 1,\ w_{2,3,1} = -\frac{w_{N,5,2}}{w_{N,5,3}} w_{2,2,1} - w_{2,3,2}$$
$$n > 2,\ w_{3,3,1} = -\frac{w_{N,5,2}}{w_{N,5,3}} w_{3,2,1} - w_{3,3,2}$$
$$w_{N,5,1} = -w_{N,5,2} - w_{N,5,3}$$

**0-0-0-0-0-0-1-0-1**

$\xi_{5,N} = 0, \tau_4 = \tau_1, \tau_1 = \tau_3, \tau_1 = \tau_2, \tau_2 = 0, w_{N,5,4} = 0, w_{N,5,3} \neq 0, w_{1,2,1} = 0, w_{1,3,2} \neq 0$

$$w_{1,3,1} = -w_{1,3,2}$$
$$n > 1,\ w_{2,3,1} = -w_{2,3,2}$$
$$n > 1,\ w_{2,2,1} = 0$$
$$n > 2,\ w_{3,3,1} = -\frac{w_{N,5,2}}{w_{N,5,3}} w_{3,2,1} - w_{3,3,2}$$
$$w_{N,5,1} = -w_{N,5,2} - w_{N,5,3}$$

**0-0-0-0-0-0-1-1-0**

$\xi_{5,N} = 0, \tau_4 = \tau_1, \tau_1 = \tau_3, \tau_1 = \tau_2, \tau_2 = 0, w_{N,5,4} = 0, w_{N,5,3} \neq 0, w_{1,2,1} \neq 0, w_{N,5,2} \neq 0$

$$w_{1,3,1} = w_{1,2,1}$$
$$n > 1,\ w_{2,3,1} = w_{2,2,1}$$
$$n > 1,\ w_{2,3,2} = 0$$
$$n > 2,\ w_{3,3,1} = w_{3,2,1} - w_{3,3,2}$$
$$w_{N,5,1} = w_{1,3,2} = 0$$
$$w_{N,5,2} = -w_{N,5,3}$$

**0-0-0-0-0-0-1-1-1**

$\xi_{5,N} = 0, \tau_4 = \tau_1, \tau_1 = \tau_3, \tau_1 = \tau_2, \tau_2 = 0, w_{N,5,4} = 0, w_{N,5,3} \neq 0, w_{1,2,1} \neq 0, w_{N,5,2} = 0$

$$w_{1,3,1} = w_{1,3,2} = 0$$
$$n > 1,\ w_{2,3,1} = w_{2,3,2} = 0$$
$$n > 2,\ w_{3,3,1} = -w_{3,3,2}$$
$$w_{N,5,1} = -w_{N,5,3}$$

**0-0-0-0-0-1-0-0-0-0-0**

$\xi_{5,N} = 0, \tau_4 = \tau_1, \tau_1 = \tau_3, \tau_1 = \tau_2, \tau_2 = 0, w_{N,5,4} \neq 0, w_{1,4,3} = 0, w_{1,2,1} = 0, w_{1,4,1} = -w_{1,4,2}, w_{1,3,1} = -w_{1,3,2}, w_{1,4,2} = -\frac{w_{N,5,3}}{w_{N,5,4}} w_{1,3,2}$

$$n > 1,\ w_{2,4,1} = -\frac{1}{w_{N,5,4}}\big(w_{N,5,3}\big(w_{2,3,1} + w_{2,3,2}\big) + w_{N,5,2}w_{2,2,1}\big) - w_{2,4,2} - w_{2,4,3}$$

$$n > 2,\ w_{3,4,1} = -\frac{1}{w_{N,5,4}}\big(w_{N,5,2}w_{3,2,1} + w_{N,5,3}\big(w_{3,3,1} + w_{3,3,2}\big)\big) - w_{3,4,2} - w_{3,4,3}$$

$$w_{N,5,1} = -w_{N,5,2} - w_{N,5,3} - w_{N,5,4}$$

**0-0-0-0-0-1-0-0-0-0-1**

$\xi_{5,N} = 0, \tau_4 = \tau_1, \tau_1 = \tau_3, \tau_1 = \tau_2, \tau_2 = 0, w_{N,5,4} \neq 0, w_{1,4,3} = 0, w_{1,2,1} = 0, w_{1,4,1} = -w_{1,4,2}, w_{1,3,1} = -w_{1,3,2}, w_{1,4,2} \neq -\frac{w_{N,5,3}}{w_{N,5,4}} w_{1,3,2}$

$$n > 1,\ w_{2,4,1} = -\frac{w_{N,5,3}}{w_{N,5,4}}\big(w_{2,3,1} + w_{2,3,2}\big) - w_{2,4,2} - w_{2,4,3}$$

$$n > 1,\ w_{2,2,1} = 0$$

$$n > 2,\ w_{3,4,1} = -\frac{1}{w_{N,5,4}}\big(w_{N,5,2}w_{3,2,1} + w_{N,5,3}\big(w_{3,3,1} + w_{3,3,2}\big)\big) - w_{3,4,2} - w_{3,4,3}$$

$$w_{N,5,1} = -w_{N,5,2} - w_{N,5,3} - w_{N,5,4}$$

**0-0-0-0-0-1-0-0-0-1-0**

$\xi_{5,N} = 0, \tau_4 = \tau_1, \tau_1 = \tau_3, \tau_1 = \tau_2, \tau_2 = 0, w_{N,5,4} \neq 0, w_{1,4,3} = 0, w_{1,2,1} = 0, w_{1,4,1} = -w_{1,4,2}, w_{1,3,1} \neq -w_{1,3,2}, w_{1,4,2} = 0$

$$n > 1,\ w_{2,4,1} = -\frac{w_{N,5,2}}{w_{N,5,4}} w_{2,2,1} - w_{2,4,2}$$

$$n > 1,\ w_{2,4,3} = 0$$

$$n > 2,\ w_{3,4,1} = -\frac{w_{N,5,2}}{w_{N,5,4}} w_{3,2,1} - w_{3,4,2} - w_{3,4,3}$$

$$w_{N,5,1} = -w_{N,5,2} - w_{N,5,4}$$

$$w_{N,5,3} = 0$$

**0-0-0-0-0-1-0-0-0-1-1**

$\xi_{5,N} = 0, \tau_4 = \tau_1, \tau_1 = \tau_3, \tau_1 = \tau_2, \tau_2 = 0, w_{N,5,4} \neq 0, w_{1,4,3} = 0, w_{1,2,1} = 0, w_{1,4,1} = -w_{1,4,2}, w_{1,3,1} \neq -w_{1,3,2}, w_{1,4,2} \neq 0$

$$n > 1,\ w_{2,4,1} = -w_{2,4,2}$$

$$n > 1,\ w_{2,4,3} = 0$$

$$n > 1,\ w_{2,2,1} = 0$$

$$n > 2,\ w_{3,4,1} = -\frac{w_{N,5,2}}{w_{N,5,4}} w_{3,2,1} - w_{3,4,2} - w_{3,4,3}$$

$$w_{N,5,1} = -w_{N,5,2} - w_{N,5,4}$$

$$w_{N,5,3} = 0$$

**0-0-0-0-0-1-0-0-1-0-0**

$\xi_{5,N} = 0, \tau_4 = \tau_1, \tau_1 = \tau_3, \tau_1 = \tau_2, \tau_2 = 0, w_{N,5,4} \neq 0, w_{1,4,3} = 0, w_{1,2,1} = 0, w_{1,4,1} \neq -w_{1,4,2}, w_{N,5,2} = 0, w_{1,4,2} = w_{1,3,2}$

$$\begin{aligned}
&w_{1,3,1} = w_{1,4,1}\\
n > 1,\ &w_{2,4,3} = 0\\
n > 1,\ &w_{2,4,1} = w_{2,3,1} + w_{2,3,2} - w_{2,4,2}\\
n > 2,\ &w_{3,3,1} = \frac{w_{N,5,2}}{w_{N,5,4}} w_{3,2,1} - w_{3,3,2} + w_{3,4,1} + w_{3,4,2} + w_{3,4,3}\\
&w_{N,5,1} = -w_{N,5,2}\\
&w_{N,5,3} = -w_{N,5,4}
\end{aligned}$$

**0-0-0-0-0-1-0-0-1-0-1**

$\xi_{5,N} = 0, \tau_4 = \tau_1, \tau_1 = \tau_3, \tau_1 = \tau_2, \tau_2 = 0, w_{N,5,4} \neq 0, w_{1,4,3} = 0, w_{1,2,1} = 0, w_{1,4,1} \neq -w_{1,4,2}, w_{N,5,2} = 0, w_{1,4,2} \neq w_{1,3,2}$

$$\begin{aligned}
&w_{1,3,1} = w_{1,4,1} + w_{1,4,2} - w_{1,3,2}\\
n > 1,\ &w_{2,4,3} = w_{2,2,1} = 0\\
n > 1,\ &w_{2,4,1} = w_{2,3,1} + w_{2,3,2} - w_{2,4,2}\\
n > 2,\ &w_{3,3,1} = \frac{w_{N,5,2}}{w_{N,5,4}} w_{3,2,1} - w_{3,3,2} + w_{3,4,1} + w_{3,4,2} + w_{3,4,3}\\
&w_{N,5,1} = -w_{N,5,2}\\
&w_{N,5,3} = -w_{N,5,4}
\end{aligned}$$

**0-0-0-0-0-1-0-0-1-1**

$\xi_{5,N} = 0, \tau_4 = \tau_1, \tau_1 = \tau_3, \tau_1 = \tau_2, \tau_2 = 0, w_{N,5,4} \neq 0, w_{1,4,3} = 0, w_{1,2,1} = 0, w_{1,4,1} \neq -w_{1,4,2}, w_{N,5,2} \neq 0$

$$\begin{aligned}
&w_{1,3,1} = w_{1,4,1} + w_{1,4,2} - w_{1,3,2}\\
n > 1,\ &w_{2,2,1} = w_{2,4,3} = 0\\
n > 1,\ &w_{2,4,1} = w_{2,3,1} + w_{2,3,2} - w_{2,4,2}\\
n > 2,\ &w_{3,3,1} = \frac{w_{N,5,2}}{w_{N,5,4}} w_{3,2,1} - w_{3,3,2} + w_{3,4,1} + w_{3,4,2} + w_{3,4,3}\\
&w_{N,5,1} = -w_{N,5,2}\\
&w_{N,5,3} = -w_{N,5,4}
\end{aligned}$$

**0-0-0-0-0-1-0-1-0-0-0**

$\xi_{5,N} = 0, \tau_4 = \tau_1, \tau_1 = \tau_3, \tau_1 = \tau_2, \tau_2 = 0, w_{N,5,4} \neq 0, w_{1,4,3} = 0, w_{1,2,1} \neq 0, w_{1,4,2} = 0, w_{N,5,3} = 0, w_{1,4,1} = 0$

$$
\begin{aligned}
n > 1,\ w_{2,4,1} &= \frac{w_{2,4,3}}{w_{1,2,1}}\big(w_{1,3,1} + w_{1,3,2}\big) - w_{2,4,3} \\
n > 1,\ w_{2,4,2} &= -\frac{w_{2,4,3}}{w_{1,2,1}}\big(w_{1,3,1} + w_{1,3,2}\big) \\
n > 2,\ w_{3,4,1} &= -w_{3,4,2} - w_{3,4,3} \\
w_{N,5,1} &= -w_{N,5,4} \\
w_{N,5,2} &= 0
\end{aligned}
$$

**0-0-0-0-0-1-0-1-0-0-1**

$\xi_{5,N} = 0, \tau_4 = \tau_1, \tau_1 = \tau_3, \tau_1 = \tau_2, \tau_2 = 0, w_{N,5,4} \neq 0, w_{1,4,3} = 0, w_{1,2,1} \neq 0, w_{1,4,2} = 0, w_{N,5,3} = 0, w_{1,4,1} \neq 0$

$$
\begin{aligned}
w_{1,4,1} &= w_{1,2,1} \\
n > 1,\ w_{2,4,1} &= w_{2,2,1} + \frac{w_{2,4,3}}{w_{1,2,1}}\big(w_{1,3,1} + w_{1,3,2}\big) - w_{2,4,3} \\
n > 1,\ w_{2,4,2} &= -\frac{w_{2,4,3}}{w_{1,2,1}}\big(w_{1,3,1} + w_{1,3,2}\big) \\
n > 2,\ w_{3,4,1} &= w_{3,2,1} - w_{3,4,2} - w_{3,4,3} \\
w_{N,5,1} &= 0 \\
w_{N,5,2} &= -w_{N,5,4}
\end{aligned}
$$

**0-0-0-0-0-1-0-1-0-1-0-0**

$\xi_{5,N} = 0, \tau_4 = \tau_1, \tau_1 = \tau_3, \tau_1 = \tau_2, \tau_2 = 0, w_{N,5,4} \neq 0, w_{1,4,3} = 0, w_{1,2,1} \neq 0, w_{1,4,2} = 0, w_{N,5,3} \neq 0, w_{1,4,1} = 0, w_{1,3,1} = 0$

$$
\begin{aligned}
n > 1,\ w_{2,4,1} &= \frac{w_{N,5,3}}{w_{N,5,4}}\big(w_{2,3,2} - \big(w_{2,3,1} + w_{2,3,2}\big)\big) - w_{2,4,3} \\
n > 1,\ w_{2,4,2} &= -\frac{w_{N,5,3}}{w_{N,5,4}} w_{2,3,2} \\
n > 2,\ w_{3,4,1} &= -\frac{w_{N,5,3}}{w_{N,5,4}}\big(w_{3,3,1} + w_{3,3,2}\big) - w_{3,4,2} - w_{3,4,3} \\
w_{N,5,1} &= -w_{N,5,3} - w_{N,5,4} \\
w_{N,5,2} &= w_{1,3,2} = 0
\end{aligned}
$$

**0-0-0-0-0-1-0-1-0-1-0-1**

$\xi_{5,N} = 0, \tau_4 = \tau_1, \tau_1 = \tau_3, \tau_1 = \tau_2, \tau_2 = 0, w_{N,5,4} \neq 0, w_{1,4,3} = 0, w_{1,2,1} \neq 0, w_{1,4,2} = 0, w_{N,5,3} \neq 0, w_{1,4,1} = 0, w_{1,3,1} \neq 0$

$$
\begin{aligned}
& w_{1,3,2} = 0 \\
& w_{1,3,1} = w_{1,2,1} \\
n > 1,\ & w_{2,4,1} = \frac{w_{N,5,3}}{w_{N,5,4}} w_{2,3,2} \\
n > 1,\ & w_{2,4,2} = -\frac{w_{N,5,3}}{w_{N,5,4}} w_{2,3,2} - w_{2,4,3} \\
n > 1,\ & w_{2,2,1} = w_{2,3,1} + w_{2,3,2} \\
n > 2,\ & w_{3,4,1} = -\frac{w_{N,5,3}}{w_{N,5,4}} \left( w_{3,3,1} + w_{3,3,2} - w_{3,2,1} \right) - w_{3,4,2} - w_{3,4,3} \\
& w_{N,5,1} = -w_{N,5,4} \\
& w_{N,5,2} = -w_{N,5,3}
\end{aligned}
$$

**0-0-0-0-0-1-0-1-0-1-1-0**

$\xi_{5,N} = 0, \tau_4 = \tau_1, \tau_1 = \tau_3, \tau_1 = \tau_2, \tau_2 = 0, w_{N,5,4} \neq 0, w_{1,4,3} = 0, w_{1,2,1} \neq 0, w_{1,4,2} = 0, w_{N,5,3} \neq 0, w_{1,4,1} \neq 0, w_{1,3,1} = 0$

$$
\begin{aligned}
& w_{1,3,2} = 0 \\
& w_{1,4,1} = w_{1,2,1} \\
n > 1,\ & w_{2,2,1} = w_{2,4,1} + w_{2,4,3} - \frac{w_{N,5,3}}{w_{N,5,4}} w_{2,3,2} \\
n > 1,\ & w_{2,4,2} = -\frac{w_{N,5,3}}{w_{N,5,4}} w_{2,3,2} \\
n > 1,\ & w_{2,3,1} = -w_{2,3,2} \\
n > 2,\ & w_{3,2,1} = \frac{w_{N,5,3}}{w_{N,5,4}} \left( w_{3,3,1} + w_{3,3,2} \right) + w_{3,4,1} + w_{3,4,2} + w_{3,4,3} \\
& w_{N,5,1} = -w_{N,5,3} \\
& w_{N,5,2} = -w_{N,5,4}
\end{aligned}
$$

**0-0-0-0-0-1-0-1-0-1-1-1-0**

$\xi_{5,N} = 0, \tau_4 = \tau_1, \tau_1 = \tau_3, \tau_1 = \tau_2, \tau_2 = 0, w_{N,5,4} \neq 0, w_{1,4,3} = 0, w_{1,2,1} \neq 0, w_{1,4,2} = 0, w_{N,5,3} \neq 0, w_{1,4,1} \neq 0, w_{1,3,1} \neq 0, w_{1,3,1} = w_{1,2,1}$

$$
\begin{aligned}
& w_{1,3,2} = 0 \\
& w_{1,4,1} = w_{1,2,1} \\
n > 1,\ & w_{2,4,1} = -\frac{w_{N,5,3}}{w_{N,5,4}} w_{2,3,1} + \frac{w_{N,5,3} + w_{N,5,4}}{w_{N,5,4}} w_{2,2,1} \\
n > 1,\ & w_{2,4,2} = -\frac{w_{N,5,3}}{w_{N,5,4}} w_{2,3,2} - w_{2,4,3} \\
n > 2,\ & w_{3,4,1} = \frac{w_{N,5,3} + w_{N,5,4}}{w_{N,5,4}} w_{3,2,1} - \frac{w_{N,5,3}}{w_{N,5,4}} \left( w_{3,3,1} + w_{3,3,2} \right) - w_{3,4,2} - w_{3,4,3} \\
& w_{N,5,1} = 0 \\
& w_{N,5,2} = -w_{N,5,3} - w_{N,5,4}
\end{aligned}
$$

**0-0-0-0-0-1-0-1-0-1-1-1-1**

$\xi_{5,N} = 0, \tau_4 = \tau_1, \tau_1 = \tau_3, \tau_1 = \tau_2, \tau_2 = 0, w_{N,5,4} \neq 0, w_{1,4,3} = 0, w_{1,2,1} \neq 0, w_{1,4,2} = 0, w_{N,5,3} \neq 0, w_{1,4,1} \neq 0, w_{1,3,1} \neq 0, w_{1,3,1} \neq w_{1,2,1}$

$$
\begin{aligned}
& w_{1,4,1} = w_{1,3,1} \\
n > 1,\ & w_{2,3,2} = w_{2,4,2} + w_{2,4,3} \frac{w_{1,3,1}}{w_{1,2,1}} \\
n > 1,\ & w_{2,3,1} = w_{2,4,1} + w_{2,4,3} \left( 1 - \frac{w_{1,3,1}}{w_{1,2,1}} \right) \\
n > 2,\ & w_{3,3,1} = w_{3,4,1} + w_{3,4,2} + w_{3,4,3} - w_{3,3,2} \\
& w_{N,5,1} = w_{N,5,2} = w_{1,3,2} = 0 \\
& w_{N,5,3} = -w_{N,5,4}
\end{aligned}
$$

**0-0-0-0-0-1-0-1-1-0**

$\xi_{5,N} = 0, \tau_4 = \tau_1, \tau_1 = \tau_3, \tau_1 = \tau_2, \tau_2 = 0, w_{N,5,4} \neq 0, w_{1,4,3} = 0, w_{1,2,1} \neq 0, w_{1,4,2} \neq 0, w_{1,2,1} = w_{1,3,1} + w_{1,3,2}$

$$w_{N,5,3} w_{1,3,2} \neq 0$$
$$w_{1,4,1} = w_{1,3,1} + w_{1,3,2} - w_{1,4,2}$$
$$n > 1,\ w_{2,3,1} = \frac{w_{1,3,2}}{w_{1,4,2}} w_{2,4,1} + \frac{w_{1,4,2} - w_{1,3,2}}{w_{1,4,2}} w_{2,2,1}$$
$$n > 1,\ w_{2,3,2} = \frac{w_{1,3,2}}{w_{1,4,2}} \left( w_{2,4,2} + w_{2,4,3} \right)$$
$$n > 2,\ w_{3,3,1} = \frac{w_{1,4,2} - w_{1,3,2}}{w_{1,4,2}} w_{3,2,1} + \frac{w_{1,3,2}}{w_{1,4,2}} \left( w_{3,4,1} + w_{3,4,2} + w_{3,4,3} \right) - w_{3,3,2}$$
$$w_{N,5,1} = 0$$
$$w_{N,5,2} = -w_{N,5,3} \frac{w_{1,4,2} - w_{1,3,2}}{w_{1,4,2}}$$
$$w_{N,5,4} = -w_{N,5,3} \frac{w_{1,3,2}}{w_{1,4,2}}$$

**0-0-0-0-0-1-0-1-1-1-0**

$\xi_{5,N} = 0, \tau_4 = \tau_1, \tau_1 = \tau_3, \tau_1 = \tau_2, \tau_2 = 0, w_{N,5,4} \neq 0, w_{1,4,3} = 0, w_{1,2,1} \neq 0, w_{1,4,2} \neq 0, w_{1,2,1} \neq w_{1,3,1} + w_{1,3,2}, w_{1,4,2} = w_{1,3,2}$

$$w_{N,5,3} w_{1,3,2} \neq 0$$
$$w_{1,4,1} = w_{1,3,1}$$
$$n > 1,\ w_{2,3,1} = w_{2,4,1} + w_{2,4,3} \left( 1 - \frac{w_{1,3,1} + w_{1,3,2}}{w_{1,2,1}} \right)$$
$$n > 1,\ w_{2,3,2} = w_{2,4,2} + w_{2,4,3} \frac{w_{1,3,1} + w_{1,3,2}}{w_{1,2,1}}$$
$$n > 2,\ w_{3,3,1} = w_{3,4,1} + w_{3,4,2} + w_{3,4,3} - w_{3,3,2}$$
$$w_{N,5,1} = w_{N,5,2} = 0$$
$$w_{N,5,4} = -w_{N,5,3}$$

**0-0-0-0-0-1-0-1-1-1-1**

$\xi_{5,N} = 0, \tau_4 = \tau_1, \tau_1 = \tau_3, \tau_1 = \tau_2, \tau_2 = 0, w_{N,5,4} \neq 0, w_{1,4,3} = 0, w_{1,2,1} \neq 0, w_{1,4,2} \neq 0, w_{1,2,1} \neq w_{1,3,1} + w_{1,3,2}, w_{1,4,2} \neq w_{1,3,2}$

$$
\begin{aligned}
w_{N,5,3}w_{1,3,2} &\neq 0 \\
w_{1,3,1} &= -w_{1,3,2} \\
w_{1,4,1} &= -w_{1,4,2} \\
n > 1,\ w_{2,3,1} &= \frac{w_{1,3,2}}{w_{1,4,2}}\left(w_{2,4,1} + w_{2,4,3}\right) \\
n > 1,\ w_{2,3,2} &= \frac{w_{1,3,2}}{w_{1,4,2}}w_{2,4,2} \\
n > 2,\ w_{3,3,1} &= \frac{w_{1,3,2}}{w_{1,4,2}}\left(w_{3,4,1} + w_{3,4,2} + w_{3,4,3}\right) - w_{3,3,2} \\
w_{N,5,1} &= w_{N,5,3}\left(\frac{w_{1,3,2}}{w_{1,4,2}} - 1\right) \\
w_{N,5,2} &= 0 \\
w_{N,5,4} &= -w_{N,5,3}\frac{w_{1,3,2}}{w_{1,4,2}}
\end{aligned}
$$

**0-0-0-0-0-1-1-0**

$\xi_{5,N} = 0, \tau_4 = \tau_1, \tau_1 = \tau_3, \tau_1 = \tau_2, \tau_2 = 0, w_{N,5,4} \neq 0, w_{1,4,3} \neq 0, w_{1,2,1} = 0$

$$
\begin{aligned}
w_{1,4,1} &= -w_{1,4,2} - w_{1,4,3} \\
w_{1,3,1} &= -w_{1,3,2} \\
n > 1,\ w_{2,4,1} &= \left(\left(w_{N,5,3}w_{1,3,2} + w_{N,5,4}w_{1,4,2}\right)\frac{w_{N,5,3}}{w_{N,5,4}w_{1,4,3}} - w_{N,5,2}\right)\frac{w_{2,2,1}}{w_{N,5,4}} - w_{2,4,2} - w_{2,4,3} \\
n > 1,\ w_{2,3,1} &= -\left(w_{N,5,3}w_{1,3,2} - w_{N,5,4}w_{1,4,2}\right)\frac{w_{2,2,1}}{w_{1,4,3}w_{N,5,4}} - w_{2,3,2} \\
n > 2,\ w_{3,4,1} &= -\frac{w_{N,5,2}}{w_{N,5,4}}w_{3,2,1} - \frac{w_{N,5,3}}{w_{N,5,4}}\left(w_{3,3,1} + w_{3,3,2}\right) - w_{3,4,2} - w_{3,4,3} \\
w_{N,5,1} &= -w_{N,5,2} - w_{N,5,3} - w_{N,5,4}
\end{aligned}
$$

**0-0-0-0-0-1-1-1-0-0**

$\xi_{5,N} = 0, \tau_4 = \tau_1, \tau_1 = \tau_3, \tau_1 = \tau_2, \tau_2 = 0, w_{N,5,4} \neq 0, w_{1,4,3} \neq 0, w_{1,2,1} \neq 0, w_{1,4,2} = 0, w_{1,2,1} = w_{1,4,1} + w_{1,4,3}$

$$\begin{aligned} w_{1,3,1} &= w_{1,3,2} = 0 \\ n > 1,\ w_{2,2,1} &= \frac{w_{N,5,3}}{w_{N,5,4}} w_{2,3,1} + w_{2,4,1} + w_{2,4,3} \\ n > 1,\ w_{2,4,2} &= -\frac{w_{N,5,3}}{w_{N,5,4}} w_{2,3,2} \\ n > 1,\ w_{2,3,1} &= -w_{2,3,2} \\ n > 2,\ w_{3,4,1} &= w_{3,2,1} - \frac{w_{N,5,3}}{w_{N,5,4}} \big(w_{3,3,1} + w_{3,3,2}\big) - w_{3,4,2} - w_{3,4,3} \\ w_{N,5,1} &= -w_{N,5,3} \\ w_{N,5,2} &= -w_{N,5,4} \end{aligned}$$

**0-0-0-0-0-1-1-1-0-1**

$\xi_{5,N} = 0, \tau_4 = \tau_1, \tau_1 = \tau_3, \tau_1 = \tau_2, \tau_2 = 0, w_{N,5,4} \neq 0, w_{1,4,3} \neq 0, w_{1,2,1} \neq 0, w_{1,4,2} = 0, w_{1,2,1} \neq w_{1,4,1} + w_{1,4,3}$

$$\begin{aligned} w_{1,4,1} &= -w_{1,4,3} \\ n > 1,\ w_{2,4,1} &= \frac{w_{N,5,3}}{w_{N,5,4}} w_{2,3,2} - w_{2,4,3} \\ n > 1,\ w_{2,4,2} &= -\frac{w_{N,5,3}}{w_{N,5,4}} w_{2,3,2} \\ n > 1,\ w_{2,3,1} &= -w_{2,3,2} \\ n > 2,\ w_{3,4,1} &= -\frac{w_{N,5,3}}{w_{N,5,4}} \big(w_{3,3,1} + w_{3,3,2}\big) - w_{3,4,2} - w_{3,4,3} \\ w_{N,5,1} &= -w_{N,5,3} - w_{N,5,4} \\ w_{N,5,2} &= w_{1,3,1} = w_{1,3,2} = 0 \end{aligned}$$

**0-0-0-0-0-1-1-1-1-0**

$\xi_{5,N} = 0, \tau_4 = \tau_1, \tau_1 = \tau_3, \tau_1 = \tau_2, \tau_2 = 0, w_{N,5,4} \neq 0, w_{1,4,3} \neq 0, w_{1,2,1} \neq 0, w_{1,4,2} \neq 0, w_{1,4,1} = 0$

$$
\begin{aligned}
& w_{1,3,1} = w_{1,2,1} \\
& w_{1,3,2} = 0 \\
& w_{1,4,2} = -w_{1,4,3} \\
n > 1,\ & w_{2,4,1} = \frac{w_{N,5,3}}{w_{N,5,4}} w_{2,3,2} \\
n > 1,\ & w_{2,4,2} = -w_{2,4,3} - \frac{w_{N,5,3}}{w_{N,5,4}} w_{2,3,2} \\
n > 1,\ & w_{2,2,1} = w_{2,3,1} + w_{2,3,2} \\
n > 2,\ & w_{3,4,1} = \frac{w_{N,5,3}}{w_{N,5,4}} \big(w_{3,2,1} - \big(w_{3,3,1} + w_{3,3,2}\big)\big) - w_{3,4,2} - w_{3,4,3} \\
& w_{N,5,1} = -w_{N,5,4} \\
& w_{N,5,2} = -w_{N,5,3}
\end{aligned}
$$

**0-0-0-0-0-1-1-1-1-1**

$\xi_{5,N} = 0, \tau_4 = \tau_1, \tau_1 = \tau_3, \tau_1 = \tau_2, \tau_2 = 0, w_{N,5,4} \neq 0, w_{1,4,3} \neq 0, w_{1,2,1} \neq 0, w_{1,4,2} \neq 0, w_{1,4,1} \neq 0$

$$
\begin{aligned}
& w_{1,4,1} = w_{1,3,1} = w_{1,2,1} \\
& w_{1,4,2} = -w_{1,4,3} \\
n > 1,\ & w_{2,4,1} = w_{2,2,1} + \frac{w_{N,5,3}}{w_{N,5,4}} w_{2,3,2} \\
n > 1,\ & w_{2,4,2} = -w_{2,4,3} - \frac{w_{N,5,3}}{w_{N,5,4}} w_{2,3,2} \\
n > 1,\ & w_{2,2,1} = w_{2,3,1} + w_{2,3,2} \\
n > 2,\ & w_{3,4,1} = -\frac{w_{N,5,3}}{w_{N,5,4}} \big(w_{3,3,1} + w_{3,3,2} - w_{3,2,1}\big) - w_{3,4,2} - w_{3,4,3} + w_{3,2,1} \\
& w_{N,5,1} = w_{1,3,2} = 0 \\
& w_{N,5,2} = -w_{N,5,3} - w_{N,5,4}
\end{aligned}
$$

**0-0-0-0-1-0-0-0-0-0**

$\xi_{5,N} = 0, \tau_4 = \tau_1, \tau_1 = \tau_3, \tau_1 = \tau_2 \neq 0, w_{1,2,1} = 0, w_{1,4,3} = 0, w_{N,5,4} = 0, w_{N,5,3} = 0, w_{N,5,2} = 0$

$$
w_{N,5,1} = 0
$$

**0-0-0-0-1-0-0-0-0-1**

$\xi_{5,N} = 0, \tau_4 = \tau_1, \tau_1 = \tau_3, \tau_1 = \tau_2 \neq 0, w_{1,2,1} = 0, w_{1,4,3} = 0, w_{N,5,4} = 0, w_{N,5,3} = 0, w_{N,5,2} \neq 0$

$$
\begin{aligned}
n > 1,\ & w_{2,2,1} = 0 \\
n > 2,\ & w_{3,2,1} = 0 \\
& w_{N,5,1} = -w_{N,5,2}
\end{aligned}
$$

**0-0-0-0-1-0-0-0-1-0**

$\xi_{5,N} = 0, \tau_4 = \tau_1, \tau_1 = \tau_3, \tau_1 = \tau_2 \neq 0, w_{1,2,1} = 0, w_{1,4,3} = 0, w_{N,5,4} = 0, w_{N,5,3} \neq 0, w_{1,3,2} = 0$

$$w_{1,3,1} = 0$$

$$n > 1,\ w_{2,3,1} = -\frac{w_{N,5,2}}{w_{N,5,3}} w_{2,2,1} - w_{2,3,2}$$

$$n > 2,\ w_{3,3,1} = -\frac{w_{N,5,2}}{w_{N,5,3}} w_{3,2,1} - w_{3,3,2}$$

$$w_{N,5,1} = -w_{N,5,2} - w_{N,5,3}$$

**0-0-0-0-1-0-0-0-1-1**

$\xi_{5,N} = 0, \tau_4 = \tau_1, \tau_1 = \tau_3, \tau_1 = \tau_2 \neq 0, w_{1,2,1} = 0, w_{1,4,3} = 0, w_{N,5,4} = 0, w_{N,5,3} \neq 0, w_{1,3,2} \neq 0$

$$w_{1,3,1} = -w_{1,3,2}$$

$$n > 1,\ w_{2,2,1} = 0$$

$$n > 1,\ w_{2,3,1} = -w_{2,3,2}$$

$$n > 2,\ w_{3,3,1} = -\frac{w_{N,5,2}}{w_{N,5,3}} w_{3,2,1} - w_{3,3,2}$$

$$w_{N,5,1} = -w_{N,5,2} - w_{N,5,3}$$

**0-0-0-0-1-0-0-1-0-0-0**

$\xi_{5,N} = 0, \tau_4 = \tau_1, \tau_1 = \tau_3, \tau_1 = \tau_2 \neq 0, w_{1,2,1} = 0, w_{1,4,3} = 0, w_{N,5,4} \neq 0, w_{1,4,1} = -w_{1,4,2}, w_{1,3,1} = -w_{1,3,2}, w_{1,4,2} = -\frac{w_{N,5,3}}{w_{N,5,4}} w_{1,3,2}$

$$n > 1,\ w_{2,4,1} = -\frac{w_{N,5,2}}{w_{N,5,4}} w_{2,2,1} - \frac{w_{N,5,3}}{w_{N,5,4}} \big(w_{2,3,1} + w_{2,3,2}\big) - w_{2,4,2} - w_{2,4,3}$$

$$n > 2,\ w_{3,4,1} = -\frac{w_{N,5,2}}{w_{N,5,4}} w_{3,2,1} - \frac{w_{N,5,3}}{w_{N,5,4}} \big(w_{3,3,1} + w_{3,3,2}\big) - w_{3,4,2} - w_{3,4,3}$$

$$w_{N,5,1} = -w_{N,5,2} - w_{N,5,3} - w_{N,5,4}$$

**0-0-0-0-1-0-0-1-0-0-1**

$\xi_{5,N} = 0, \tau_4 = \tau_1, \tau_1 = \tau_3, \tau_1 = \tau_2 \neq 0, w_{1,2,1} = 0, w_{1,4,3} = 0, w_{N,5,4} \neq 0, w_{1,4,1} = -w_{1,4,2}, w_{1,3,1} = -w_{1,3,2}, w_{1,4,2} \neq -\frac{w_{N,5,3}}{w_{N,5,4}} w_{1,3,2}$

$$n > 1,\ w_{2,4,1} = -\frac{w_{N,5,3}}{w_{N,5,4}} \big(w_{2,3,1} + w_{2,3,2}\big) - w_{2,4,2} - w_{2,4,3}$$

$$n > 1,\ w_{2,2,1} = 0$$

$$n > 2,\ w_{3,4,1} = -\frac{w_{N,5,2}}{w_{N,5,4}} w_{3,2,1} - \frac{w_{N,5,3}}{w_{N,5,4}} \big(w_{3,3,1} + w_{3,3,2}\big) - w_{3,4,2} - w_{3,4,3}$$

$$w_{N,5,1} = -w_{N,5,2} - w_{N,5,3} - w_{N,5,4}$$

**0-0-0-0-1-0-0-1-0-1-0**

$\xi_{5,N} = 0, \tau_4 = \tau_1, \tau_1 = \tau_3, \tau_1 = \tau_2 \neq 0, w_{1,2,1} = 0, w_{1,4,3} = 0, w_{N,5,4} \neq 0, w_{1,4,1} = -w_{1,4,2}, w_{1,3,1} \neq -w_{1,3,2}, w_{1,4,2} = 0$

$$
\begin{aligned}
n > 1,\ w_{2,4,1} &= -\frac{w_{N,5,2}}{w_{N,5,4}} w_{2,2,1} - w_{2,4,2} \\
n > 1,\ w_{2,4,3} &= 0 \\
n > 2,\ w_{3,4,1} &= -\frac{w_{N,5,2}}{w_{N,5,4}} w_{3,2,1} - w_{3,4,2} - w_{3,4,3} \\
w_{N,5,1} &= -w_{N,5,2} - w_{N,5,4} \\
w_{N,5,3} &= 0
\end{aligned}
$$

**0-0-0-0-1-0-0-1-0-1-1**

$\xi_{5,N} = 0, \tau_4 = \tau_1, \tau_1 = \tau_3, \tau_1 = \tau_2 \neq 0, w_{1,2,1} = 0, w_{1,4,3} = 0, w_{N,5,4} \neq 0, w_{1,4,1} = -w_{1,4,2}, w_{1,3,1} \neq -w_{1,3,2}, w_{1,4,2} \neq 0$

$$
\begin{aligned}
n > 1,\ w_{2,4,1} &= -w_{2,4,2} \\
n > 1,\ w_{2,4,3} &= w_{2,2,1} = 0 \\
n > 2,\ w_{3,4,1} &= -\frac{w_{N,5,2}}{w_{N,5,4}} w_{3,2,1} - w_{3,4,2} - w_{3,4,3} \\
w_{N,5,1} &= -w_{N,5,2} - w_{N,5,4} \\
w_{N,5,3} &= 0
\end{aligned}
$$

**0-0-0-0-1-0-0-1-1-0-0**

$\xi_{5,N} = 0, \tau_4 = \tau_1, \tau_1 = \tau_3, \tau_1 = \tau_2 \neq 0, w_{1,2,1} = 0, w_{1,4,3} = 0, w_{N,5,4} \neq 0, w_{1,4,1} \neq -w_{1,4,2}, w_{N,5,2} = 0, w_{1,3,2} = w_{1,4,2}$

$$
\begin{aligned}
w_{1,3,1} &= w_{1,4,1} \\
n > 1,\ w_{2,3,1} &= w_{2,4,1} + w_{2,4,2} - w_{2,3,2} \\
n > 1,\ w_{2,4,3} &= 0 \\
n > 2,\ w_{3,3,1} &= w_{3,4,1} + w_{3,4,2} + w_{3,4,3} + \frac{w_{N,5,2}}{w_{N,5,4}} w_{3,2,1} - w_{3,3,2} \\
w_{N,5,1} &= 0 \\
w_{N,5,3} &= -w_{N,5,4}
\end{aligned}
$$

**0-0-0-0-1-0-0-1-1-0-1**

$\xi_{5,N} = 0, \tau_4 = \tau_1, \tau_1 = \tau_3, \tau_1 = \tau_2 \neq 0, w_{1,2,1} = 0, w_{1,4,3} = 0, w_{N,5,4} \neq 0, w_{1,4,1} \neq -w_{1,4,2}, w_{N,5,2} = 0, w_{1,3,2} \neq w_{1,4,2}$

$$
\begin{aligned}
w_{N,5,1} &= 0 \\
w_{N,5,3} &= -w_{N,5,4} \\
w_{1,3,1} &= w_{1,4,1} + w_{1,4,2} - w_{1,3,2} \\
n > 1,\ w_{2,3,1} &= w_{2,4,1} + w_{2,4,2} - w_{2,3,2} \\
n > 1,\ w_{2,4,3} &= w_{2,2,1} = 0 \\
n > 2,\ w_{3,3,1} &= w_{3,4,1} + w_{3,4,2} + w_{3,4,3} + \frac{w_{N,5,2}}{w_{N,5,4}} w_{3,2,1} - w_{3,3,2}
\end{aligned}
$$

**0-0-0-0-1-0-0-1-1-1**

$\xi_{5,N} = 0, \tau_4 = \tau_1, \tau_1 = \tau_3, \tau_1 = \tau_2, \tau_2 \neq 0, w_{1,2,1} = 0, w_{1,4,3} = 0, w_{N,5,4} \neq 0, w_{1,4,1} \neq -w_{1,4,2}, w_{N,5,2} \neq 0$

$$
\begin{aligned}
w_{N,5,1} &= -w_{N,5,2} \\
w_{N,5,3} &= -w_{N,5,4} \\
w_{1,3,1} &= w_{1,4,1} + w_{1,4,2} - w_{1,3,2} \\
n > 1,\ w_{2,2,1} &= w_{2,4,3} = 0 \\
n > 1,\ w_{2,3,1} &= w_{2,4,1} + w_{2,4,2} - w_{2,3,2} \\
n > 2,\ w_{3,3,1} &= w_{3,4,1} + w_{3,4,2} + w_{3,4,3} + \frac{w_{N,5,2}}{w_{N,5,4}} w_{3,2,1} - w_{3,3,2}
\end{aligned}
$$

**0-0-0-0-1-0-1-0-0-0**

$\xi_{5,N} = 0, \tau_4 = \tau_1, \tau_1 = \tau_3, \tau_1 = \tau_2, \tau_2 \neq 0, w_{1,2,1} = 0, w_{1,4,3} \neq 0, w_{N,5,4} = 0, w_{N,5,3} = 0, w_{N,5,2} = 0$

$$w_{N,5,1} = 0$$

**0-0-0-0-1-0-1-0-0-1**

$\xi_{5,N} = 0, \tau_4 = \tau_1, \tau_1 = \tau_3, \tau_1 = \tau_2, \tau_2 \neq 0, w_{1,2,1} = 0, w_{1,4,3} \neq 0, w_{N,5,4} = 0, w_{N,5,3} = 0, w_{N,5,2} \neq 0$

$$
\begin{aligned}
w_{N,5,1} &= -w_{N,5,2} \\
n > 1,\ w_{2,2,1} &= w_{3,2,1} = 0
\end{aligned}
$$

**0-0-0-0-1-0-1-0-1-0**

$\xi_{5,N} = 0, \tau_4 = \tau_1, \tau_1 = \tau_3, \tau_1 = \tau_2, \tau_2 \neq 0, w_{1,2,1} = 0, w_{1,4,3} \neq 0, w_{N,5,4} = 0, w_{N,5,3} \neq 0, w_{1,3,2} = 0$

$$
\begin{aligned}
w_{N,5,1} &= -w_{N,5,2} - w_{N,5,3} \\
w_{1,3,1} &= 0 \\
n > 1,\ w_{2,3,1} &= -\frac{w_{N,5,2}}{w_{N,5,3}} w_{2,2,1} - w_{2,3,2} \\
n > 2,\ w_{3,3,1} &= -\frac{w_{N,5,2}}{w_{N,5,3}} w_{3,2,1} - w_{3,3,2}
\end{aligned}
$$

**0-0-0-0-1-0-1-0-1-1**

$\xi_{5,N} = 0, \tau_4 = \tau_1, \tau_1 = \tau_3, \tau_1 = \tau_2, \tau_2 \neq 0, w_{1,2,1} = 0, w_{1,4,3} \neq 0, w_{N,5,4} = 0, w_{N,5,3} \neq 0, w_{1,3,2} \neq 0$

$$w_{N,5,1} = -w_{N,5,2} - w_{N,5,3}$$
$$w_{1,3,1} = -w_{1,3,2}$$
$$n > 1,\ w_{2,3,1} = -w_{2,3,2}$$
$$n > 1,\ w_{2,2,1} = 0$$
$$n > 2,\ w_{3,3,1} = -\frac{w_{N,5,2}}{w_{N,5,3}} w_{3,2,1} - w_{3,3,2}$$

**0-0-0-0-1-0-1-1**

$\xi_{5,N} = 0, \tau_4 = \tau_1, \tau_1 = \tau_3, \tau_1 = \tau_2, \tau_2 \neq 0, w_{1,2,1} = 0, w_{1,4,3} \neq 0, w_{N,5,4} \neq 0$

$$w_{N,5,1} = -w_{N,5,2} - w_{N,5,3} - w_{N,5,4}$$
$$w_{1,3,1} = -w_{1,3,2}$$
$$w_{1,4,1} = -w_{1,4,2} - w_{1,4,3}$$
$$n > 1,\ w_{2,4,1} = -\left( w_{N,5,2} + \frac{w_{N,5,3}}{w_{1,4,3}} \left( \frac{w_{N,5,3}}{w_{N,5,4}} w_{1,3,2} + w_{1,4,2} \right) \right) \frac{w_{2,2,1}}{w_{N,5,4}} - w_{2,4,2} - w_{2,4,3}$$
$$n > 1,\ w_{2,3,1} = -\left( \frac{w_{N,5,3}}{w_{N,5,4}} w_{1,3,2} + w_{1,4,2} \right) \frac{w_{2,2,1}}{w_{1,4,3}} - w_{2,3,2}$$
$$n > 2,\ w_{3,4,1} = -\frac{w_{N,5,2}}{w_{N,5,4}} w_{3,2,1} - \frac{w_{N,5,3}}{w_{N,5,4}} (w_{3,3,1} + w_{3,3,2}) - w_{3,4,2} - w_{3,4,3}$$

**0-0-0-0-1-1-0-0-0**

$\xi_{5,N} = 0, \tau_4 = \tau_1, \tau_1 = \tau_3, \tau_1 = \tau_2, \tau_2 \neq 0, w_{1,2,1} \neq 0, w_{N,5,4} = 0, w_{N,5,3} = 0, w_{N,5,2} = 0$

$$w_{N,5,1} = 0$$

**0-0-0-0-1-1-0-0-1**

$\xi_{5,N} = 0, \tau_4 = \tau_1, \tau_1 = \tau_3, \tau_1 = \tau_2, \tau_2 \neq 0, w_{1,2,1} \neq 0, w_{N,5,4} = 0, w_{N,5,3} = 0, w_{N,5,2} \neq 0$

$$w_{N,5,1} = -w_{N,5,2}$$
$$w_{1,2,1} = 0$$
$$n > 1,\ w_{2,2,1} = 0$$
$$n > 2,\ w_{3,2,1} = 0$$

**0-0-0-0-1-1-0-1-0**

$\xi_{5,N} = 0, \tau_4 = \tau_1, \tau_1 = \tau_3, \tau_1 = \tau_2, \tau_2 \neq 0, w_{1,2,1} \neq 0, w_{N,5,4} = 0, w_{N,5,3} \neq 0, w_{1,3,1} = 0$

$$w_{N,5,1} = -w_{N,5,3}$$
$$w_{N,5,2} = w_{1,3,2} = 0$$
$$n > 1,\ w_{2,3,1} = w_{2,3,2} = 0$$
$$n > 2,\ w_{3,3,1} = -w_{3,3,2}$$

**0-0-0-0-1-1-0-1-1**

$\xi_{5,N}=0, \tau_4=\tau_1, \tau_1=\tau_3, \tau_1=\tau_2, \tau_2\neq 0, w_{1,2,1}\neq 0, w_{N,5,4}=0, w_{N,5,3}\neq 0, w_{1,3,1}\neq 0$

$$\begin{aligned} w_{N,5,1} &= w_{1,3,2}=0 \\ w_{N,5,2} &= -w_{N,5,3} \\ w_{1,2,1} &= w_{1,3,1} \\ n>1,\ w_{2,3,1} &= w_{2,2,1} \\ n>1,\ w_{2,3,2} &= 0 \\ n>2,\ w_{3,2,1} &= w_{3,3,1}+w_{3,3,2} \end{aligned}$$

**0-0-0-0-1-1-1-0-0-0**

$\xi_{5,N}=0, \tau_4=\tau_1, \tau_1=\tau_3, \tau_1=\tau_2, \tau_2\neq 0, w_{1,2,1}\neq 0, w_{N,5,4}\neq 0, w_{1,4,3}=0, w_{N,5,3}=0, w_{1,4,1}=0$

$$\begin{aligned} w_{N,5,1} &= -w_{N,5,4} \\ w_{N,5,2} &= w_{1,4,2}=0 \\ n>1,\ w_{2,4,1} &= \frac{w_{2,4,3}}{w_{1,2,1}}\left(w_{1,3,1}+w_{1,3,2}-w_{1,2,1}\right) \\ n>1,\ w_{2,4,2} &= -w_{2,4,3}\frac{w_{1,3,1}+w_{1,3,2}}{w_{1,2,1}} \\ n>2,\ w_{3,4,1} &= -\frac{w_{N,5,2}}{w_{N,5,4}}w_{3,2,1}-w_{3,4,2}-w_{3,4,3} \end{aligned}$$

**0-0-0-0-1-1-1-0-0-1**

$\xi_{5,N}=0, \tau_4=\tau_1, \tau_1=\tau_3, \tau_1=\tau_2, \tau_2\neq 0, w_{1,2,1}\neq 0, w_{N,5,4}\neq 0, w_{1,4,3}=0, w_{N,5,3}=0, w_{1,4,1}\neq 0$

$$\begin{aligned} w_{N,5,1} &= w_{1,4,2}=0 \\ w_{N,5,2} &= -w_{N,5,4} \\ w_{1,4,1} &= w_{1,2,1} \\ n>1,\ w_{2,4,1} &= w_{2,2,1}-\frac{w_{2,4,3}}{w_{1,2,1}}\left(w_{1,2,1}-\left(w_{1,3,1}+w_{1,3,2}\right)\right) \\ n>1,\ w_{2,4,2} &= -w_{2,4,3}\frac{w_{1,3,1}+w_{1,3,2}}{w_{1,2,1}} \\ n>2,\ w_{3,4,1} &= -\frac{w_{N,5,2}}{w_{N,5,4}}w_{3,2,1}-w_{3,4,2}-w_{3,4,3} \end{aligned}$$

**0-0-0-0-1-1-1-0-1-0-0-0**

$\xi_{5,N} = 0, \tau_4 = \tau_1, \tau_1 = \tau_3, \tau_1 = \tau_2, \tau_2 \neq 0, w_{1,2,1} \neq 0, w_{N,5,4} \neq 0, w_{1,4,3} = 0, w_{N,5,3} \neq 0, w_{1,3,2} = 0, w_{1,3,1} = 0, w_{1,4,1} = 0$

$$\begin{aligned}
w_{N,5,1} &= -w_{N,5,3} - w_{N,5,4} \\
w_{N,5,2} &= 0 \\
w_{1,4,2} &= 0 \\
n > 1,\ w_{2,4,1} &= -\frac{w_{N,5,3}}{w_{N,5,4}} w_{2,3,1} - w_{2,4,3} \\
n > 1,\ w_{2,4,2} &= -\frac{w_{N,5,3}}{w_{N,5,4}} w_{2,3,2} \\
n > 2,\ w_{3,4,1} &= -\frac{w_{N,5,3}}{w_{N,5,4}} \big(w_{3,3,1} + w_{3,3,2}\big) - w_{3,4,2} - w_{3,4,3}
\end{aligned}$$

**0-0-0-0-1-1-1-0-1-0-0-1**

$\xi_{5,N} = 0, \tau_4 = \tau_1, \tau_1 = \tau_3, \tau_1 = \tau_2, \tau_2 \neq 0, w_{1,2,1} \neq 0, w_{N,5,4} \neq 0, w_{1,4,3} = 0, w_{N,5,3} \neq 0, w_{1,3,2} = 0, w_{1,3,1} = 0, w_{1,4,1} \neq 0$

$$\begin{aligned}
w_{N,5,1} &= -w_{N,5,3} \\
w_{N,5,2} &= -w_{N,5,4} \\
w_{1,4,2} &= 0 \\
w_{1,4,1} &= w_{1,2,1} \\
n > 1,\ w_{2,3,1} &= -w_{2,3,2} \\
n > 1,\ w_{2,4,1} &= w_{2,2,1} + \frac{w_{N,5,3}}{w_{N,5,4}} w_{2,3,2} - w_{2,4,3} \\
n > 1,\ w_{2,4,2} &= -\frac{w_{N,5,3}}{w_{N,5,4}} w_{2,3,2} \\
n > 2,\ w_{3,4,1} &= w_{3,2,1} - \frac{w_{N,5,3}}{w_{N,5,4}} \big(w_{3,3,1} + w_{3,3,2}\big) - w_{3,4,2} - w_{3,4,3}
\end{aligned}$$

**0-0-0-0-1-1-1-0-1-0-1-0-0**

$\xi_{5,N} = 0, \tau_4 = \tau_1, \tau_1 = \tau_3, \tau_1 = \tau_2, \tau_2 \neq 0, w_{1,2,1} \neq 0, w_{N,5,4} \neq 0, w_{1,4,3} = 0, w_{N,5,3} \neq 0, w_{1,3,2} = 0, w_{1,3,1} \neq 0, w_{1,3,1} = w_{1,2,1}, w_{1,4,1} = 0$

$$\begin{aligned}
w_{N,5,1} &= -w_{N,5,4} \\
w_{N,5,2} &= -w_{N,5,3} \\
w_{1,4,2} &= 0 \\
n > 1,\ w_{2,4,1} &= \frac{w_{N,5,3}}{w_{N,5,4}} w_{2,3,2} \\
n > 1,\ w_{2,4,2} &= -\frac{w_{N,5,3}}{w_{N,5,4}} w_{2,3,2} - w_{2,4,3} \\
n > 1,\ w_{2,3,1} &= w_{2,2,1} - w_{2,3,2} \\
n > 2,\ w_{3,4,1} &= -\frac{w_{N,5,2}}{w_{N,5,4}} w_{3,2,1} - \frac{w_{N,5,3}}{w_{N,5,4}} \big(w_{3,3,1} + w_{3,3,2}\big) - w_{3,4,2} - w_{3,4,3}
\end{aligned}$$

**0-0-0-0-1-1-1-0-1-0-1-0-1**

$\xi_{5,N} = 0, \tau_4 = \tau_1, \tau_1 = \tau_3, \tau_1 = \tau_2, \tau_2 \neq 0, w_{1,2,1} \neq 0, w_{N,5,4} \neq 0, w_{1,4,3} = 0, w_{N,5,3} \neq 0, w_{1,3,2} = 0, w_{1,3,1} \neq 0, w_{1,3,1} = w_{1,2,1}, w_{1,4,1} \neq 0$

$$\begin{aligned}
w_{N,5,1} &= 0 \\
w_{N,5,2} &= -w_{N,5,3} - w_{N,5,4} \\
w_{1,4,2} &= 0 \\
w_{1,4,1} &= w_{1,2,1} \\
n > 1,\ w_{2,4,1} &= \frac{w_{N,5,3}}{w_{N,5,4}} \big(w_{2,2,1} - w_{2,3,1}\big) + w_{2,2,1} \\
n > 1,\ w_{2,4,2} &= -\frac{w_{N,5,3}}{w_{N,5,4}} w_{2,3,2} - w_{2,4,3} \\
n > 2,\ w_{3,4,1} &= \frac{w_{N,5,3} + w_{N,5,4}}{w_{N,5,4}} w_{3,2,1} - \frac{w_{N,5,3}}{w_{N,5,4}} \big(w_{3,3,1} + w_{3,3,2}\big) - w_{3,4,2} - w_{3,4,3}
\end{aligned}$$

**0-0-0-0-1-1-1-0-1-0-1-1**

$\xi_{5,N} = 0, \tau_4 = \tau_1, \tau_1 = \tau_3, \tau_1 = \tau_2, \tau_2 \neq 0, w_{1,2,1} \neq 0, w_{N,5,4} \neq 0, w_{1,4,3} = 0, w_{N,5,3} \neq 0, w_{1,3,2} = 0, w_{1,3,1} \neq 0, w_{1,3,1} \neq w_{1,2,1}$

$$\begin{aligned}
w_{N,5,1} &= w_{N,5,2} = w_{1,4,2} = 0 \\
w_{N,5,3} &= -w_{N,5,4} \\
w_{1,4,1} &= w_{1,3,1} \\
n > 1,\ w_{2,3,1} &= w_{2,4,1} + \frac{w_{2,4,3}}{w_{1,2,1}} \big(w_{1,2,1} - w_{1,3,1}\big) \\
n > 1,\ w_{2,3,2} &= w_{2,4,2} + w_{2,4,3} \frac{w_{1,3,1}}{w_{1,2,1}} \\
n > 2,\ w_{3,4,1} &= w_{3,3,1} + w_{3,3,2} - w_{3,4,2} - w_{3,4,3}
\end{aligned}$$

**0-0-0-0-1-1-1-0-1-1-0**

$\xi_{5,N} = 0, \tau_4 = \tau_1, \tau_1 = \tau_3, \tau_1 = \tau_2, \tau_2 \neq 0, w_{1,2,1} \neq 0, w_{N,5,4} \neq 0, w_{1,4,3} = 0, w_{N,5,3} \neq 0, w_{1,3,2} \neq 0, w_{N,5,3} = -w_{N,5,4}$

$$w_{N,5,1} = -w_{N,5,3} - w_{N,5,4}$$
$$w_{N,5,2} = 0$$
$$w_{1,4,2} = w_{1,3,2}$$
$$w_{1,4,1} = w_{1,3,1}$$
$$n > 1,\ w_{2,3,1} = w_{2,4,1} + w_{2,4,3}\frac{w_{1,2,1} - \left(w_{1,3,1} + w_{1,3,2}\right)}{w_{1,2,1}}$$
$$n > 1,\ w_{2,3,2} = w_{2,4,2} + w_{2,4,3}\frac{w_{1,3,1} + w_{1,3,2}}{w_{1,2,1}}$$
$$n > 2,\ w_{3,3,1} = w_{3,4,1} + w_{3,4,2} + w_{3,4,3} - w_{3,3,2}$$

**0-0-0-0-1-1-1-0-1-1-1-0**

$\xi_{5,N} = 0, \tau_4 = \tau_1, \tau_1 = \tau_3, \tau_1 = \tau_2, \tau_2 \neq 0, w_{1,2,1} \neq 0, w_{N,5,4} \neq 0, w_{1,4,3} = 0, w_{N,5,3} \neq 0, w_{1,3,2} \neq 0, w_{N,5,3} \neq -w_{N,5,4}, w_{1,4,1} = -w_{1,4,2}$

$$w_{N,5,1} = -w_{N,5,3} - w_{N,5,4}$$
$$w_{N,5,2} = 0$$
$$w_{1,4,2} = -\frac{w_{N,5,3}}{w_{N,5,4}}w_{1,3,2}$$
$$w_{1,3,1} = -w_{1,3,2}$$
$$n > 1,\ w_{2,4,1} = -\frac{w_{N,5,3}}{w_{N,5,4}}w_{2,3,1} - w_{2,4,3}$$
$$n > 1,\ w_{2,4,2} = -\frac{w_{N,5,3}}{w_{N,5,4}}w_{2,3,2}$$
$$n > 2,\ w_{3,4,1} = -\frac{w_{N,5,3}}{w_{N,5,4}}\left(w_{3,3,1} + w_{3,3,2}\right) - w_{3,4,2} - w_{3,4,3}$$

**0-0-0-0-1-1-1-0-1-1-1-1**

$\xi_{5,N} = 0, \tau_4 = \tau_1, \tau_1 = \tau_3, \tau_1 = \tau_2, \tau_2 \neq 0, w_{1,2,1} \neq 0, w_{N,5,4} \neq 0, w_{1,4,3} = 0, w_{N,5,3} \neq 0, w_{1,3,2} \neq 0, w_{N,5,3} \neq -w_{N,5,4}, w_{1,4,1} \neq -w_{1,4,2}$

$$
\begin{aligned}
w_{N,5,1} &= 0 \\
w_{N,5,2} &= -w_{N,5,3} - w_{N,5,4} \\
w_{1,4,2} &= -\frac{w_{N,5,3}}{w_{N,5,4}} w_{1,3,2} \\
w_{1,4,1} &= w_{1,2,1} + \frac{w_{N,5,3}}{w_{N,5,4}} w_{1,3,2} \\
w_{1,3,1} &= w_{1,2,1} - w_{1,3,2} \\
n > 1,\ w_{2,4,1} &= \frac{w_{N,5,3}}{w_{N,5,4}} \big( w_{2,2,1} - w_{2,3,1} \big) + w_{2,2,1} \\
n > 1,\ w_{2,4,2} &= -\frac{w_{N,5,3}}{w_{N,5,4}} w_{2,3,2} - w_{2,4,3} \\
n > 2,\ w_{3,4,1} &= \frac{w_{N,5,3} + w_{N,5,4}}{w_{N,5,4}} w_{3,2,1} - \frac{w_{N,5,3}}{w_{N,5,4}} \big( w_{3,3,1} + w_{3,3,2} \big) - w_{3,4,2} - w_{3,4,3}
\end{aligned}
$$

**0-0-0-0-1-1-1-1-0-0**

$\xi_{5,N} = 0, \tau_4 = \tau_1, \tau_1 = \tau_3, \tau_1 = \tau_2, \tau_2 \neq 0, w_{1,2,1} \neq 0, w_{N,5,4} \neq 0, w_{1,4,3} \neq 0, w_{1,4,2} = 0, w_{1,4,1} = -w_{1,4,3}$

$$
\begin{aligned}
w_{N,5,1} &= -w_{N,5,3} - w_{N,5,4} \\
w_{N,5,2} &= w_{1,3,1} = w_{1,3,2} = 0 \\
n > 1,\ w_{2,4,2} &= -\frac{w_{N,5,3}}{w_{N,5,4}} w_{2,3,2} \\
n > 1,\ w_{2,3,1} &= -w_{2,3,2} \\
n > 1,\ w_{2,4,1} &= \frac{w_{N,5,3}}{w_{N,5,4}} w_{2,3,2} - w_{2,4,3} \\
n > 2,\ w_{3,4,1} &= -\frac{w_{N,5,3}}{w_{N,5,4}} \big( w_{3,3,1} + w_{3,3,2} \big) - w_{3,4,2} - w_{3,4,3}
\end{aligned}
$$

**0-0-0-0-1-1-1-1-0-1**

$\xi_{5,N} = 0, \tau_4 = \tau_1, \tau_1 = \tau_3, \tau_1 = \tau_2, \tau_2 \neq 0, w_{1,2,1} \neq 0, w_{N,5,4} \neq 0, w_{1,4,3} \neq 0, w_{1,4,2} = 0, w_{1,4,1} \neq -w_{1,4,3}$

$$
\begin{aligned}
w_{N,5,1} &= -w_{N,5,3} \\
w_{N,5,2} &= -w_{N,5,4} \\
w_{1,3,1} &= w_{1,3,2} = 0 \\
w_{1,4,1} &= w_{1,2,1} - w_{1,4,3} \\
n > 1,\ w_{2,4,2} &= -\frac{w_{N,5,3}}{w_{N,5,4}} w_{2,3,2} \\
n > 1,\ w_{2,3,1} &= -w_{2,3,2} \\
n > 1,\ w_{2,4,1} &= w_{2,2,1} + \frac{w_{N,5,3}}{w_{N,5,4}} w_{2,3,2} - w_{2,4,3} \\
n > 2,\ w_{3,2,1} &= \frac{w_{N,5,3}}{w_{N,5,4}} \left(w_{3,3,1} + w_{3,3,2}\right) + w_{3,4,1} + w_{3,4,2} + w_{3,4,3}
\end{aligned}
$$

**0-0-0-0-1-1-1-1-1-0**

$\xi_{5,N} = 0, \tau_4 = \tau_1, \tau_1 = \tau_3, \tau_1 = \tau_2, \tau_2 \neq 0, w_{1,2,1} \neq 0, w_{N,5,4} \neq 0, w_{1,4,3} \neq 0, w_{1,4,2} \neq 0, w_{1,4,1} = 0$

$$
\begin{aligned}
w_{N,5,1} &= -w_{N,5,4} \\
w_{N,5,2} &= -w_{N,5,3} \\
w_{1,4,3} &= -w_{1,4,2} \\
w_{1,3,2} &= 0 \\
w_{1,3,1} &= w_{1,2,1} \\
n > 1,\ w_{2,4,2} &= -w_{N,5,3} \frac{w_{2,3,2}}{w_{N,5,4}} - w_{2,4,3} \\
n > 1,\ w_{2,3,1} &= w_{2,2,1} - w_{2,3,2} \\
n > 1,\ w_{2,4,1} &= w_{N,5,3} \frac{w_{2,3,2}}{w_{N,5,4}} \\
n > 2,\ w_{3,4,1} &= \frac{w_{N,5,3}}{w_{N,5,4}} \left(w_{3,2,1} - \left(w_{3,3,1} + w_{3,3,2}\right)\right) - w_{3,4,2} = w_{3,4,3}
\end{aligned}
$$

**0-0-0-0-1-1-1-1-1-1**

$\xi_{5,N} = 0, \tau_4 = \tau_1, \tau_1 = \tau_3, \tau_1 = \tau_2, \tau_2 \neq 0, w_{1,2,1} \neq 0, w_{N,5,4} \neq 0, w_{1,4,3} \neq 0, w_{1,4,2} \neq 0, w_{1,4,1} \neq 0$

$$w_{N,5,1} = w_{1,3,2} = 0$$
$$w_{N,5,2} = -w_{N,5,3} - w_{N,5,4}$$
$$w_{1,4,2} = -w_{1,4,3}$$
$$w_{1,4,1} = w_{1,3,1} = w_{1,2,1}$$
$$n > 1,\ w_{2,4,1} = \frac{w_{N,5,3}}{w_{N,5,4}} w_{2,3,2} + w_{2,2,1}$$
$$n > 1,\ w_{2,4,2} = -\frac{w_{N,5,3}}{w_{N,5,4}} w_{2,3,2} - w_{2,4,3}$$
$$n > 1,\ w_{2,3,1} = w_{2,2,1} - w_{2,3,2}$$
$$n > 2,\ w_{3,4,1} = -\frac{w_{N,5,3}}{w_{N,5,4}} \big(w_{3,3,1} + w_{3,3,2} - w_{3,2,1}\big) + w_{3,2,1} - w_{3,4,2} - w_{3,4,3}$$

**0-0-0-1-0-0**

$\xi_{5,N} = 0, \tau_4 = \tau_1, \tau_1 = \tau_3, \tau_1 \neq \tau_2, w_{N,5,4} = 0, w_{N,5,3} = 0$

$$w_{N,5,1} = w_{N,5,2} = 0$$

**0-0-0-1-0-1**

$\xi_{5,N} = 0, \tau_4 = \tau_1, \tau_1 = \tau_3, \tau_1 \neq \tau_2, w_{N,5,4} = 0, w_{N,5,3} \neq 0$

$$w_{N,5,1} = -w_{N,5,3}$$
$$w_{N,5,2} = w_{1,3,1} = w_{1,3,2} = 0$$
$$n > 1,\ w_{2,3,1} = w_{2,3,2} = 0$$
$$n > 2,\ w_{3,3,1} = -w_{3,3,2}$$

**0-0-0-1-1-0-0-0**

$\xi_{5,N} = 0, \tau_4 = \tau_1, \tau_1 = \tau_3, \tau_1 \neq \tau_2, w_{N,5,4} \neq 0, w_{1,4,3} = 0, w_{N,5,3} = 0, w_{1,3,1} = -w_{1,3,2}$

$$w_{N,5,1} = -w_{N,5,4}$$
$$w_{N,5,2} = w_{1,4,1} = w_{1,4,2} = 0$$
$$n > 1,\ w_{2,4,1} = -w_{2,4,3}$$
$$n > 1,\ w_{2,4,2} = 0$$
$$n > 2,\ w_{3,4,1} = -w_{3,4,2} - w_{3,4,3}$$

**0-0-0-1-1-0-0-1**

$\xi_{5,N} = 0, \tau_4 = \tau_1, \tau_1 = \tau_3, \tau_1 \neq \tau_2, w_{N,5,4} \neq 0, w_{1,4,3} = 0, w_{N,5,3} = 0, w_{1,3,1} \neq -w_{1,3,2}$

$$w_{N,5,1} = -w_{N,5,4}$$
$$w_{N,5,2} = w_{1,4,1} = w_{1,4,2} = 0$$
$$n > 1,\ w_{2,4,1} = w_{2,4,2} = w_{2,4,3} = 0$$
$$n > 2,\ w_{3,4,1} = -w_{3,4,2} - w_{3,4,3}$$

**0-0-0-1-1-0-1-0**

$\xi_{5,N} = 0, \tau_4 = \tau_1, \tau_1 = \tau_3, \tau_1 \neq \tau_2, w_{N,5,4} \neq 0, w_{1,4,3} = 0, w_{N,5,3} \neq 0, w_{1,3,1} = -w_{1,3,2}$

$$w_{N,5,1} = -w_{N,5,3} - w_{N,5,4}$$
$$w_{N,5,2} = 0$$
$$w_{1,4,1} = \frac{w_{N,5,3}}{w_{N,5,4}} w_{1,3,2}$$
$$w_{1,4,2} = -\frac{w_{N,5,3}}{w_{N,5,4}} w_{1,3,2}$$
$$n > 1,\ w_{2,4,1} = -\frac{w_{N,5,3}}{w_{N,5,4}} w_{2,3,1} - w_{2,4,3}$$
$$n > 1,\ w_{2,4,2} = -\frac{w_{N,5,3}}{w_{N,5,4}} w_{2,3,2}$$
$$n > 2,\ w_{3,4,1} = -\frac{w_{N,5,3}}{w_{N,5,4}} \left(w_{3,3,1} + w_{3,3,2}\right) - w_{3,4,2} - w_{3,4,3}$$

**0-0-0-1-1-0-1-1**

$\xi_{5,N} = 0, \tau_4 = \tau_1, \tau_1 = \tau_3, \tau_1 \neq \tau_2, w_{N,5,4} \neq 0, w_{1,4,3} = 0, w_{N,5,3} \neq 0, w_{1,3,1} \neq -w_{1,3,2}$

$$w_{N,5,1} = w_{N,5,2} = 0$$
$$w_{N,5,3} = -w_{N,5,4}$$
$$w_{1,3,2} = w_{1,4,2}$$
$$w_{1,4,1} = w_{1,3,1}$$
$$n > 1,\ w_{2,4,1} = w_{2,3,1}$$
$$n > 1,\ w_{2,3,2} = w_{2,4,2}$$
$$n > 1,\ w_{2,4,3} = 0$$
$$n > 2,\ w_{3,3,1} = w_{3,4,1} + w_{3,4,2} + w_{3,4,3} - w_{3,3,2}$$

**0-0-0-1-1-1**

$\xi_{5,N} = 0, \tau_4 = \tau_1, \tau_1 = \tau_3, \tau_1 \neq \tau_2, w_{N,5,4} \neq 0, w_{1,4,3} \neq 0$

$$w_{N,5,1} = -w_{N,5,3} - w_{N,5,4}$$
$$w_{N,5,2} = w_{1,4,2} = 0$$
$$w_{1,4,1} = -w_{1,4,3}$$
$$w_{1,3,1} = w_{1,3,2} = 0$$
$$n > 1,\ w_{2,4,1} = \frac{w_{N,5,3}}{w_{N,5,4}} w_{2,3,2} - w_{2,4,3}$$
$$n > 1,\ w_{2,4,2} = -\frac{w_{N,5,3}}{w_{N,5,4}} w_{2,3,2}$$
$$n > 1,\ w_{2,3,1} = -w_{2,3,2}$$
$$n > 2,\ w_{3,4,1} = -\frac{w_{N,5,3}}{w_{N,5,4}} \left(w_{3,3,1} + w_{3,3,2}\right) - w_{3,4,2} - w_{3,4,3}$$

**0-0-1-0-0-0**

$\xi_{5,N} = 0, \tau_4 = \tau_1, \tau_1 \neq \tau_3, \tau_3 = \tau_2, w_{N,5,4} = 0, w_{N,5,2} = 0$

$$w_{N,5,1} = w_{N,5,3} = 0$$

**0-0-1-0-0-1**

$\xi_{5,N} = 0, \tau_4 = \tau_1, \tau_1 \neq \tau_3, \tau_3 = \tau_2, w_{N,5,4} = 0, w_{N,5,2} \neq 0$

$$w_{N,5,1} = w_{1,3,2} = 0$$
$$w_{N,5,3} = -w_{N,5,2}$$
$$w_{1,3,1} = w_{1,2,1}$$
$$n > 1,\ w_{2,2,1} = w_{2,3,1}$$
$$n > 1,\ w_{2,3,2} = 0$$
$$n > 2,\ w_{3,2,1} = w_{3,3,1} + w_{3,3,2}$$

**0-0-1-0-1-0-0-0**

$\xi_{5,N} = 0, \tau_4 = \tau_1, \tau_1 \neq \tau_3, \tau_3 = \tau_2, w_{N,5,4} \neq 0, w_{N,5,2} = 0, w_{1,4,3} = 0, w_{1,2,1} = w_{1,3,1} + w_{1,3,2}$

$$w_{N,5,1} = -w_{N,5,4}$$
$$w_{N,5,3} = w_{1,4,1} = w_{1,4,2} = 0$$
$$n > 1,\ w_{2,4,1} = 0$$
$$n > 1,\ w_{2,4,2} = -w_{2,4,3}$$
$$n > 2,\ w_{3,4,1} = -w_{3,4,2} - w_{3,4,3}$$

**0-0-1-0-1-0-0-1**

$\xi_{5,N} = 0, \tau_4 = \tau_1, \tau_1 \neq \tau_3, \tau_3 = \tau_2, w_{N,5,4} \neq 0, w_{N,5,2} = 0, w_{1,4,3} = 0, w_{1,2,1} \neq w_{1,3,1} + w_{1,3,2}$

$$\begin{aligned}
w_{N,5,1} &= -w_{N,5,4} \\
w_{N,5,3} &= w_{1,4,1} = w_{1,4,2} = 0 \\
n > 1,\ w_{2,4,1} &= w_{2,4,2} = w_{2,4,3} = 0 \\
n > 2,\ w_{3,4,1} &= -w_{3,4,2} - w_{3,4,3}
\end{aligned}$$

**0-0-1-0-1-0-1**

$\xi_{5,N} = 0, \tau_4 = \tau_1, \tau_1 \neq \tau_3, \tau_3 = \tau_2, w_{N,5,4} \neq 0, w_{N,5,2} = 0, w_{1,4,3} \neq 0$

$$\begin{aligned}
w_{N,5,1} &= -w_{N,5,4} \\
w_{N,5,3} &= w_{1,4,1} = w_{1,3,2} = 0 \\
w_{1,4,2} &= -w_{1,4,3} \\
w_{1,3,1} &= w_{1,2,1} \\
n > 1,\ w_{2,4,1} &= 0 \\
n > 1,\ w_{2,4,2} &= -w_{2,4,3} \\
n > 1,\ w_{2,2,1} &= w_{2,3,1} + w_{2,3,2} \\
n > 2,\ w_{3,4,1} &= -w_{3,4,2} - w_{3,4,3}
\end{aligned}$$

**0-0-1-0-1-1**

$\xi_{5,N} = 0, \tau_4 = \tau_1, \tau_1 \neq \tau_3, \tau_3 = \tau_2, w_{N,5,4} \neq 0, w_{N,5,2} \neq 0$

$$\begin{aligned}
w_{N,5,1} &= -w_{N,5,4} \\
w_{N,5,3} &= -w_{N,5,2} \\
w_{1,3,2} &= w_{1,4,1} = 0 \\
w_{1,3,1} &= w_{1,2,1} \\
w_{1,4,2} &= -w_{1,4,3} \\
n > 1,\ w_{2,4,1} &= -\frac{w_{N,5,2}}{w_{N,5,4}} w_{2,3,2} \\
n > 1,\ w_{2,4,2} &= \frac{w_{N,5,2}}{w_{N,5,4}} w_{2,3,2} - w_{2,4,3} \\
n > 1,\ w_{2,2,1} &= w_{2,3,1} + w_{2,3,2} \\
n > 2,\ w_{3,4,1} &= \frac{w_{N,5,2}}{w_{N,5,4}} \left( w_{3,3,1} + w_{3,3,2} - w_{3,2,1} \right) - w_{3,4,2} - w_{3,4,3}
\end{aligned}$$

**0-0-1-1-0-0-0**

$\xi_{5,N} = 0, \tau_4 = \tau_1, \tau_1 \neq \tau_3, \tau_2 \neq \tau_3, \tau_2 = \tau_1, w_{N,5,4} = 0, w_{N,5,2} = 0$

$$w_{N,5,1} = w_{N,5,3} = 0$$

**0-0-1-1-0-0-1**

$\xi_{5,N} = 0, \tau_4 = \tau_1, \tau_1 \neq \tau_3, \tau_2 \neq \tau_3, \tau_2 = \tau_1, w_{N,5,4} = 0, w_{N,5,2} \neq 0$

$$w_{N,5,1} = -w_{N,5,2}$$
$$w_{N,5,3} = w_{1,2,1} = 0$$
$$n > 1,\ w_{2,2,1} = 0$$
$$n > 2,\ w_{3,2,1} = 0$$

**0-0-1-1-0-1-0-0**

$\xi_{5,N} = 0, \tau_4 = \tau_1, \tau_1 \neq \tau_3, \tau_2 \neq \tau_3, \tau_2 = \tau_1, w_{N,5,4} \neq 0, w_{1,2,1} = 0, w_{1,4,2} = 0$

$$w_{N,5,1} = -w_{N,5,2} - w_{N,5,4}$$
$$w_{N,5,3} = w_{1,4,3} = w_{1,4,1} = 0$$
$$n > 1,\ w_{2,4,1} = -\frac{w_{N,5,2}}{w_{N,5,4}} w_{2,2,1} - w_{2,4,2}$$
$$n > 1,\ w_{2,4,3} = 0$$
$$n > 2,\ w_{3,4,1} = -\frac{w_{N,5,2}}{w_{N,5,4}} w_{3,2,1} - w_{3,4,2} - w_{3,4,3}$$

**0-0-1-1-0-1-0-1**

$\xi_{5,N} = 0, \tau_4 = \tau_1, \tau_1 \neq \tau_3, \tau_2 \neq \tau_3, \tau_2 = \tau_1, w_{N,5,4} \neq 0, w_{1,2,1} = 0, w_{1,4,2} \neq 0$

$$w_{N,5,1} = -w_{N,5,2} - w_{N,5,4}$$
$$w_{N,5,3} = w_{1,4,3} = 0$$
$$w_{1,4,1} = -w_{1,4,2}$$
$$n > 1,\ w_{2,4,1} = -w_{2,4,2}$$
$$n > 1,\ w_{2,4,3} = w_{2,2,1} = 0$$
$$n > 2,\ w_{3,4,1} = -\frac{w_{N,5,2}}{w_{N,5,4}} w_{3,2,1} - w_{3,4,2} - w_{3,4,3}$$

**0-0-1-1-0-1-1-0**

$\xi_{5,N} = 0, \tau_4 = \tau_1, \tau_1 \neq \tau_3, \tau_2 \neq \tau_3, \tau_2 = \tau_1, w_{N,5,4} \neq 0, w_{1,2,1} \neq 0, w_{1,4,1} = 0$

$$w_{N,5,1} = -w_{N,5,4}$$
$$w_{N,5,3} = w_{N,5,2} = w_{1,4,3} = w_{1,4,2} = 0$$
$$n > 1,\ w_{2,4,1} = -\frac{w_{N,5,2}}{w_{N,5,4}} w_{2,2,1}$$
$$n > 1,\ w_{2,4,3} = w_{2,4,2} = 0$$
$$n > 2,\ w_{3,4,1} = -w_{3,4,2} - w_{3,4,3}$$

**0-0-1-1-0-1-1-1**

$\xi_{5,N} = 0, \tau_4 = \tau_1, \tau_1 \neq \tau_3, \tau_2 \neq \tau_3, \tau_2 = \tau_1, w_{N,5,4} \neq 0, w_{1,2,1} \neq 0, w_{1,4,1} \neq 0$

$$\begin{aligned} w_{N,5,1} &= w_{N,5,3} = 0 \\ w_{N,5,2} &= -w_{N,5,4} \\ w_{1,4,3} &= w_{1,4,2} = 0 \\ w_{1,4,1} &= w_{1,2,1} \\ n > 1,\ w_{2,4,1} &= w_{2,2,1} \\ n > 1,\ w_{2,4,3} &= w_{2,4,2} = 0 \\ n > 2,\ w_{3,2,1} &= w_{3,4,1} + w_{3,4,2} + w_{3,4,3} \end{aligned}$$

**0-0-1-1-1-0**

$\xi_{5,N} = 0, \tau_4 = \tau_1, \tau_1 \neq \tau_3, \tau_2 \neq \tau_3, \tau_2 \neq \tau_1, w_{N,5,4} = 0$

$$w_{N,5,1} = w_{N,5,3} = w_{N,5,2} = 0$$

**0-0-1-1-1-1-0**

$\xi_{5,N} = 0, \tau_4 = \tau_1, \tau_1 \neq \tau_3, \tau_2 \neq \tau_3, \tau_2 \neq \tau_1, w_{N,5,4} \neq 0, w_{1,2,1} = \frac{\tau_2 - \tau_1}{\tau_3 - \tau_1}\left(w_{1,3,1} + w_{1,3,2}\right)$

$$\begin{aligned} w_{N,5,1} &= -w_{N,5,4} \\ w_{N,5,3} &= w_{N,5,2} = 0 \\ w_{1,4,1} &= w_{1,4,3} = w_{1,4,2} = 0 \\ n > 1,\ w_{2,4,1} &= -w_{2,4,2}\frac{\tau_3 - \tau_2}{\tau_3 - \tau_1} \\ n > 1,\ w_{2,4,3} &= -w_{2,4,2}\frac{\tau_2 - \tau_1}{\tau_3 - \tau_1} \\ n > 2,\ w_{3,4,1} &= -w_{3,4,2} - w_{3,4,3} \end{aligned}$$

**0-0-1-1-1-1-1**

$\xi_{5,N} = 0, \tau_4 = \tau_1, \tau_1 \neq \tau_3, \tau_2 \neq \tau_3, \tau_2 \neq \tau_1, w_{N,5,4} \neq 0, w_{1,2,1} \neq \frac{\tau_2 - \tau_1}{\tau_3 - \tau_1}\left(w_{1,3,1} + w_{1,3,2}\right)$

$$\begin{aligned} w_{N,5,1} &= -w_{N,5,4} \\ w_{N,5,3} &= w_{N,5,2} = 0 \\ w_{1,4,1} &= w_{1,4,3} = w_{1,4,2} = 0 \\ n > 1,\ w_{2,4,1} &= w_{2,4,2} = w_{2,4,3} = 0 \\ n > 2,\ w_{3,4,1} &= -w_{3,4,2} - w_{3,4,3} \end{aligned}$$

**0-1-0-0-0-0**

$\xi_{5,N} = 0, \tau_4 \neq \tau_1, \tau_4 = \tau_3, \tau_3 = \tau_2, w_{N,5,4} = 0, w_{N,5,3} = 0$

$$w_{N,5,1} = w_{N,5,2} = 0$$

**0-1-0-0-0-1**

$\xi_{5,N} = 0, \tau_4 \neq \tau_1, \tau_4 = \tau_3, \tau_3 = \tau_2, w_{N,5,4} = 0, w_{N,5,3} \neq 0$

$$w_{N,5,1} = w_{1,3,2} = 0$$
$$w_{N,5,2} = -w_{N,5,3}$$
$$w_{1,3,1} = w_{1,2,1}$$
$$n > 1,\ w_{2,3,1} = w_{2,2,1}$$
$$n > 1,\ w_{2,3,2} = 0$$
$$n > 2,\ w_{3,2,1} = w_{3,3,1} + w_{3,3,2}$$

**0-1-0-0-1-0-0-0-0**

$\xi_{5,N} = 0, \tau_4 \neq \tau_1, \tau_4 = \tau_3, \tau_3 = \tau_2, w_{N,5,4} \neq 0, w_{1,4,3} = 0, w_{1,3,1} = w_{1,2,1}, w_{N,5,2} = 0, w_{1,3,2} = 0$

$$w_{N,5,1} = 0$$
$$w_{N,5,3} = -w_{N,5,4}$$
$$w_{1,4,1} = w_{1,2,1}$$
$$w_{1,4,2} = w_{1,3,2}$$
$$n > 1,\ w_{2,3,1} = w_{2,4,1}$$
$$n > 1,\ w_{2,3,2} = w_{2,4,2} + w_{2,4,3}$$
$$n > 2,\ w_{3,3,1} = w_{3,4,1} + w_{3,4,2} + w_{3,4,3} - w_{3,3,2}$$

**0-1-0-0-1-0-0-0-1**

$\xi_{5,N} = 0, \tau_4 \neq \tau_1, \tau_4 = \tau_3, \tau_3 = \tau_2, w_{N,5,4} \neq 0, w_{1,4,3} = 0, w_{1,3,1} = w_{1,2,1}, w_{N,5,2} = 0, w_{1,3,2} \neq 0$

$$w_{N,5,1} = 0$$
$$w_{N,5,3} = -w_{N,5,4}$$
$$w_{1,4,1} = w_{1,2,1}$$
$$w_{1,4,2} = w_{1,3,2}$$
$$n > 1,\ w_{2,3,1} = w_{2,4,1}$$
$$n > 1,\ w_{2,3,2} = w_{2,4,2}$$
$$n > 1,\ w_{2,4,3} = 0$$
$$n > 2,\ w_{3,3,1} = w_{3,4,1} + w_{3,4,2} + w_{3,4,3} - w_{3,3,2}$$

**0-1-0-0-1-0-0-1-0**

$\xi_{5,N} = 0, \tau_4 \neq \tau_1, \tau_4 = \tau_3, \tau_3 = \tau_2, w_{N,5,4} \neq 0, w_{1,4,3} = 0, w_{1,3,1} = w_{1,2,1}, w_{N,5,2} \neq 0, w_{1,3,2} = 0$

$$
\begin{aligned}
w_{N,5,1} &= w_{1,4,2} = 0 \\
w_{N,5,2} &= -\left(w_{N,5,3} + w_{N,5,4}\right) \\
w_{1,4,1} &= w_{1,2,1} \\
n > 1,\ w_{2,4,1} &= w_{2,2,1} - \frac{w_{N,5,3}}{w_{N,5,4}}\left(w_{2,3,1} - w_{2,2,1}\right) \\
n > 1,\ w_{2,4,2} &= -\frac{w_{N,5,3}}{w_{N,5,4}} w_{2,3,2} - w_{2,4,3} \\
n > 2,\ w_{3,4,1} &= w_{3,2,1} - \frac{w_{N,5,3}}{w_{N,5,4}}\left(w_{3,3,1} + w_{3,3,2} - w_{3,2,1}\right) - w_{3,4,2} - w_{3,4,3}
\end{aligned}
$$

**0-1-0-0-1-0-0-1-1**

$\xi_{5,N} = 0, \tau_4 \neq \tau_1, \tau_4 = \tau_3, \tau_3 = \tau_2, w_{N,5,4} \neq 0, w_{1,4,3} = 0, w_{1,3,1} = w_{1,2,1}, w_{N,5,2} \neq 0, w_{1,3,2} \neq 0$

$$
\begin{aligned}
w_{N,5,1} &= w_{N,5,3} = w_{1,4,2} = 0 \\
w_{N,5,2} &= -w_{N,5,4} \\
w_{1,4,1} &= w_{1,2,1} \\
n > 1,\ w_{2,4,1} &= w_{2,2,1} \\
n > 1,\ w_{2,4,2} &= w_{2,4,3} = 0 \\
n > 2,\ w_{3,2,1} &= w_{3,4,1} + w_{3,4,2} + w_{3,4,3}
\end{aligned}
$$

**0-1-0-0-1-0-1-0**

$\xi_{5,N} = 0, \tau_4 \neq \tau_1, \tau_4 = \tau_3, \tau_3 = \tau_2, w_{N,5,4} \neq 0, w_{1,4,3} = 0, w_{1,3,1} \neq w_{1,2,1}, w_{1,2,1} = w_{1,3,1} + w_{1,3,2}$

$$
\begin{aligned}
w_{N,5,1} &= 0 \\
w_{N,5,2} &= -w_{N,5,3} - w_{N,5,4} \\
w_{1,4,1} &= \frac{w_{N,5,3}}{w_{N,5,4}} w_{1,3,2} + w_{1,2,1} \\
w_{1,4,2} &= -\frac{w_{N,5,3}}{w_{N,5,4}} w_{1,3,2} \\
n > 1,\ w_{2,4,1} &= -\frac{w_{N,5,3}}{w_{N,5,4}}\left(w_{2,3,1} - w_{2,2,1}\right) + w_{2,2,1} \\
n > 1,\ w_{2,4,2} &= -\frac{w_{N,5,3}}{w_{N,5,4}} w_{2,3,2} - w_{2,4,3} \\
n > 2,\ w_{3,4,1} &= w_{3,2,1} - \frac{w_{N,5,3}}{w_{N,5,4}}\left(w_{3,3,1} + w_{3,3,2} - w_{3,2,1}\right) - w_{3,4,2} - w_{3,4,3}
\end{aligned}
$$

**0-1-0-0-1-0-1-1-0**

$\xi_{5,N} = 0, \tau_4 \neq \tau_1, \tau_4 = \tau_3, \tau_3 = \tau_2, w_{N,5,4} \neq 0, w_{1,4,3} = 0, w_{1,3,1} \neq w_{1,2,1}, w_{1,2,1} \neq w_{1,3,1} + w_{1,3,2}, w_{N,5,3} = 0$

$$
\begin{aligned}
w_{N,5,1} &= 0 \\
w_{N,5,2} &= -w_{N,5,4} \\
w_{1,4,1} &= w_{1,2,1} \\
w_{1,4,2} &= 0 \\
n > 1,\ w_{2,4,1} &= w_{2,2,1} \\
n > 1,\ w_{2,4,2} &= w_{2,4,3} = 0 \\
n > 2,\ w_{3,2,1} &= w_{3,4,1} + w_{3,4,2} + w_{3,4,3}
\end{aligned}
$$

**0-1-0-0-1-0-1-1-1**

$\xi_{5,N} = 0, \tau_4 \neq \tau_1, \tau_4 = \tau_3, \tau_3 = \tau_2, w_{N,5,4} \neq 0, w_{1,4,3} = 0, w_{1,3,1} \neq w_{1,2,1}, w_{1,2,1} \neq w_{1,3,1} + w_{1,3,2}, w_{N,5,3} \neq 0$

$$
\begin{aligned}
w_{N,5,1} &= w_{N,5,2} = 0 \\
w_{1,4,1} &= w_{1,3,1} \\
w_{1,4,2} &= w_{1,3,2} \\
w_{N,5,3} &= -w_{N,5,4} \\
n > 1,\ w_{2,3,1} &= w_{2,4,1} \\
n > 1,\ w_{2,3,2} &= w_{2,4,2} \\
n > 1,\ w_{2,4,3} &= 0 \\
n > 2,\ w_{3,3,1} &= w_{3,4,1} + w_{3,4,2} + w_{3,4,3} - w_{3,3,2}
\end{aligned}
$$

**0-1-0-0-1-1**

$\xi_{5,N} = 0, \tau_4 \neq \tau_1, \tau_4 = \tau_3, \tau_3 = \tau_2, w_{N,5,4} \neq 0, w_{1,4,3} \neq 0$

$$
\begin{aligned}
w_{N,5,1} &= w_{1,3,2} = 0 \\
w_{N,5,2} &= -\left(w_{N,5,3} + w_{N,5,4}\right) \\
w_{1,4,1} &= w_{1,3,1} = w_{1,2,1} \\
w_{1,4,2} &= -w_{1,4,3} \\
n > 1,\ w_{2,4,1} &= w_{2,3,1} + w_{2,3,2}\frac{w_{N,5,4} + w_{N,5,3}}{w_{N,5,4}} \\
n > 1,\ w_{2,4,2} &= -\frac{w_{N,5,3}}{w_{N,5,4}}w_{2,3,2} - w_{2,4,3} \\
n > 1,\ w_{2,2,1} &= w_{2,3,1} + w_{2,3,2} \\
n > 2,\ w_{3,4,1} &= -\frac{w_{N,5,3}}{w_{N,5,4}}\left(w_{3,3,1} + w_{3,3,2} - w_{3,2,1}\right) + w_{3,2,1} - w_{3,4,2} - w_{3,4,3}
\end{aligned}
$$

**0-1-0-1-0-0**

$\xi_{5,N} = 0, \tau_4 \neq \tau_1, \tau_4 = \tau_3, \tau_3 \neq \tau_2, w_{N,5,4} = 0, w_{N,5,2} = 0$

$$w_{N,5,1} = w_{N,5,3} = 0$$

**0-1-0-1-0-1**

$\xi_{5,N} = 0, \tau_4 \neq \tau_1, \tau_4 = \tau_3, \tau_3 \neq \tau_2, w_{N,5,4} = 0, w_{N,5,2} \neq 0$

$$w_{N,5,1} = -w_{N,5,2}$$
$$\tau_2 = \tau_1$$
$$w_{1,2,1} = w_{N,5,3} = 0$$
$$n > 1,\ w_{2,2,1} = w_{3,2,1} = 0$$

**0-1-0-1-1-0-0-0**

$\xi_{5,N} = 0, \tau_4 \neq \tau_1, \tau_4 = \tau_3, \tau_3 \neq \tau_2, w_{N,5,4} \neq 0, w_{1,3,2} = w_{1,4,2}, w_{N,5,2} = 0, w_{1,2,1} = \frac{\tau_2-\tau_1}{\tau_3-\tau_1}\big(w_{1,3,1} + w_{1,3,2}\big)$

$$w_{N,5,1} = w_{1,4,3} = 0$$
$$w_{N,5,3} = -w_{N,5,4}$$
$$w_{1,3,1} = w_{1,4,1}$$
$$n > 1,\ w_{2,3,1} = w_{2,4,1} + \big(w_{2,4,2} - w_{2,3,2}\big)\frac{\tau_3 - \tau_2}{\tau_3 - \tau_1}$$
$$n > 1,\ w_{2,4,3} = -\big(w_{2,4,2} - w_{2,3,2}\big)\frac{\tau_2 - \tau_1}{\tau_3 - \tau_1}$$
$$n > 2,\ w_{3,3,1} = -w_{3,3,2} + w_{3,4,1} + w_{3,4,2} + w_{3,4,3}$$

**0-1-0-1-1-0-0-1**

$\xi_{5,N} = 0, \tau_4 \neq \tau_1, \tau_4 = \tau_3, \tau_3 \neq \tau_2, w_{N,5,4} \neq 0, w_{1,3,2} = w_{1,4,2}, w_{N,5,2} = 0, w_{1,2,1} \neq \frac{\tau_2-\tau_1}{\tau_3-\tau_1}\big(w_{1,3,1} + w_{1,3,2}\big)$

$$w_{N,5,1} = w_{1,4,3} = 0$$
$$w_{N,5,3} = -w_{N,5,4}$$
$$w_{1,3,1} = w_{1,4,1}$$
$$n > 1,\ w_{2,3,1} = w_{2,4,1}$$
$$n > 1,\ w_{2,4,3} = 0$$
$$n > 1,\ w_{2,4,2} = w_{2,3,2}$$
$$n > 2,\ w_{3,3,1} = -w_{3,3,2} + w_{3,4,1} + w_{3,4,2} + w_{3,4,3}$$

**0-1-0-1-1-0-1-0-0**

$\xi_{5,N} = 0, \tau_4 \neq \tau_1, \tau_4 = \tau_3, \tau_3 \neq \tau_2, w_{N,5,4} \neq 0, w_{1,3,2} = w_{1,4,2}, w_{N,5,2} \neq 0, w_{1,3,1} = -w_{1,3,2}, \tau_3 = 0$

$$
\begin{aligned}
w_{N,5,1} &= -w_{N,5,2} \\
\tau_2 &= \tau_1 \\
w_{1,2,1} &= w_{1,4,3} = 0 \\
w_{N,5,3} &= -w_{N,5,4} \\
w_{1,3,1} &= w_{1,4,1} \\
n > 1,\ w_{2,3,1} &= w_{2,4,1} + w_{2,4,2} + \frac{w_{N,5,2}}{w_{N,5,4}} w_{2,2,1} - w_{2,3,2} \\
n > 1,\ w_{2,4,3} &= 0 \\
n > 2,\ w_{3,3,1} &= \frac{w_{N,5,2}}{w_{N,5,4}} w_{3,2,1} - w_{3,3,2} + \left(w_{3,4,1} + w_{3,4,2} + w_{3,4,3}\right)
\end{aligned}
$$

**0-1-0-1-1-0-1-0-1**

$\xi_{5,N} = 0, \tau_4 \neq \tau_1, \tau_4 = \tau_3, \tau_3 \neq \tau_2, w_{N,5,4} \neq 0, w_{1,3,2} = w_{1,4,2}, w_{N,5,2} \neq 0, w_{1,3,1} = -w_{1,3,2}, \tau_3 \neq 0$

$$
\begin{aligned}
w_{N,5,1} &= -w_{N,5,2} \\
\tau_2 &= \tau_1 \\
w_{1,2,1} &= w_{1,4,3} = 0 \\
w_{N,5,3} &= -w_{N,5,4} \\
w_{1,3,1} &= w_{1,4,1} \\
n > 1,\ w_{2,3,1} &= w_{2,4,1} + w_{2,4,2} - w_{2,3,2} \\
n > 1,\ w_{2,4,3} &= w_{2,2,1} = 0 \\
n > 2,\ w_{3,3,1} &= \frac{w_{N,5,2}}{w_{N,5,4}} w_{3,2,1} - w_{3,3,2} + \left(w_{3,4,1} + w_{3,4,2} + w_{3,4,3}\right)
\end{aligned}
$$

**0-1-0-1-1-0-1-1**

$\xi_{5,N} = 0, \tau_4 \neq \tau_1, \tau_4 = \tau_3, \tau_3 \neq \tau_2, w_{N,5,4} \neq 0, w_{1,3,2} = w_{1,4,2}, w_{N,5,2} \neq 0, w_{1,3,1} \neq -w_{1,3,2}$

$$
\begin{aligned}
w_{N,5,1} &= -w_{N,5,2} \\
\tau_2 &= \tau_1 \\
w_{1,2,1} &= w_{1,4,3} = 0 \\
w_{N,5,3} &= -w_{N,5,4} \\
w_{1,3,1} &= w_{1,4,1} \\
n > 1,\ w_{2,3,1} &= w_{2,4,1} + w_{2,4,2} - w_{2,3,2} \\
n > 1,\ w_{2,4,3} &= w_{2,2,1} = 0 \\
n > 2,\ w_{3,3,1} &= \frac{w_{N,5,2}}{w_{N,5,4}} w_{3,2,1} - w_{3,3,2} + \left(w_{3,4,1} + w_{3,4,2} + w_{3,4,3}\right)
\end{aligned}
$$

**0-1-0-1-1-1**

$\xi_{5,N} = 0, \tau_4 \neq \tau_1, \tau_4 = \tau_3, \tau_3 \neq \tau_2, w_{N,5,4} \neq 0, w_{1,3,2} \neq w_{1,4,2}$

$$
\begin{aligned}
w_{N,5,1} &= -w_{N,5,2} \\
w_{N,5,3} &= -w_{N,5,4} \\
w_{1,4,1} &= -w_{1,4,2} + w_{1,3,2} + w_{1,3,1} \\
w_{1,4,3} &= w_{1,2,1} = 0 \\
\tau_2 &= \tau_1 \\
n > 1,\ w_{2,3,1} &= w_{2,4,1} + w_{2,4,2} - w_{2,3,2} \\
n > 1,\ w_{2,4,3} &= w_{2,2,1} = 0 \\
n > 2,\ w_{3,3,1} &= \frac{w_{N,5,2}}{w_{N,5,4}} w_{3,2,1} - w_{3,3,2} + w_{3,4,1} + w_{3,4,2} + w_{3,4,3}
\end{aligned}
$$

**0-1-1-0-0-0**

$\xi_{5,N} = 0, \tau_4 \neq \tau_1, \tau_4 \neq \tau_3, \tau_2 = \tau_1, w_{N,5,3} = 0, w_{N,5,2} = 0$

$$
w_{N,5,1} = w_{N,5,4} = 0
$$

**0-1-1-0-0-1**

$\xi_{5,N} = 0, \tau_4 \neq \tau_1, \tau_4 \neq \tau_3, \tau_2 = \tau_1, w_{N,5,3} = 0, w_{N,5,2} \neq 0$

$$
\begin{aligned}
w_{N,5,1} &= -w_{N,5,2} \\
w_{1,2,1} &= w_{N,5,4} = 0 \\
n > 1,\ w_{2,2,1} &= 0 \\
n > 2,\ w_{3,2,1} &= 0
\end{aligned}
$$

**0-1-1-0-1-0-0-0**

$\xi_{5,N} = 0, \tau_4 \neq \tau_1, \tau_4 \neq \tau_3, \tau_2 = \tau_1, w_{N,5,3} \neq 0, w_{1,3,1} = -w_{1,3,2}, w_{1,2,1} = 0, w_{1,3,2} = 0$

$$
\begin{aligned}
w_{N,5,1} &= -w_{N,5,2} - w_{N,5,3} \\
\tau_3 &= \tau_1 \\
w_{N,5,4} &= 0 \\
n > 1,\ w_{2,3,1} &= -\frac{w_{N,5,2}}{w_{N,5,3}} w_{2,2,1} - w_{2,3,2} \\
n > 2,\ w_{3,3,1} &= -\frac{w_{N,5,2}}{w_{N,5,3}} w_{3,2,1} - w_{3,3,2}
\end{aligned}
$$

**0-1-1-0-1-0-0-1**

$\xi_{5,N} = 0, \tau_4 \neq \tau_1, \tau_4 \neq \tau_3, \tau_2 = \tau_1, w_{N,5,3} \neq 0, w_{1,3,1} = -w_{1,3,2}, w_{1,2,1} = 0, w_{1,3,2} \neq 0$

$$w_{N,5,1} = -w_{N,5,2} - w_{N,5,3}$$
$$\tau_3 = \tau_1$$
$$w_{N,5,4} = 0$$
$$n > 1,\ w_{2,3,1} = -w_{2,3,2}$$
$$n > 1,\ w_{2,2,1} = 0$$
$$n > 2,\ w_{3,3,1} = -\frac{w_{N,5,2}}{w_{N,5,3}} w_{3,2,1} - w_{3,3,2}$$

**0-1-1-0-1-0-1**

$\xi_{5,N} = 0, \tau_4 \neq \tau_1, \tau_4 \neq \tau_3, \tau_2 = \tau_1, w_{N,5,3} \neq 0, w_{1,3,1} = -w_{1,3,2}, w_{1,2,1} \neq 0$

$$w_{N,5,1} = -w_{N,5,3}$$
$$w_{N,5,2} = w_{N,5,4} = w_{1,3,2} = 0$$
$$\tau_3 = \tau_1$$
$$n > 1,\ w_{2,3,1} = w_{2,3,2} = 0$$
$$n > 2,\ w_{3,3,1} = -w_{3,3,2}$$

**0-1-1-0-1-1**

$\xi_{5,N} = 0, \tau_4 \neq \tau_1, \tau_4 \neq \tau_3, \tau_2 = \tau_1, w_{N,5,3} \neq 0, w_{1,3,1} \neq -w_{1,3,2}$

$$w_{N,5,1} = -w_{N,5,4}$$
$$w_{N,5,2} = -w_{N,5,3}$$
$$\tau_3 = \tau_1$$
$$w_{1,3,2} = w_{N,5,4} = 0$$
$$w_{1,3,1} = w_{1,2,1}$$
$$n > 1,\ w_{2,3,1} = w_{2,2,1}$$
$$n > 1,\ w_{2,3,2} = 0$$
$$n > 2,\ w_{3,3,1} = w_{3,2,1} - w_{3,3,2}$$

**0-1-1-1-0-0-0-0**

$\xi_{5,N} = 0, \tau_4 \neq \tau_1, \tau_4 \neq \tau_3, \tau_2 \neq \tau_1, \tau_4 = \tau_2, \tau_3 = \tau_1, w_{N,5,3} = 0, w_{N,5,4} = 0$

$$w_{N,5,1} = w_{N,5,2} = 0$$

**0-1-1-1-0-0-0-1-0-0**

$\xi_{5,N} = 0, \tau_4 \neq \tau_1, \tau_4 \neq \tau_3, \tau_2 \neq \tau_1, \tau_4 = \tau_2, \tau_3 = \tau_1, w_{N,5,3} = 0, w_{N,5,4} \neq 0, w_{1,4,3} = 0, w_{1,3,1} = -w_{1,3,2}$

$$\begin{aligned} w_{N,5,1} &= 0 \\ w_{N,5,2} &= -w_{N,5,4} \\ w_{1,2,1} &= w_{1,4,1} + w_{1,4,3} \\ w_{1,4,2} &= 0 \\ n > 1,\ w_{2,2,1} &= w_{2,4,1} + w_{2,4,3} \\ n > 1,\ w_{2,4,2} &= 0 \\ n > 2,\ w_{3,2,1} &= w_{3,4,1} + w_{3,4,2} + w_{3,4,3} \end{aligned}$$

**0-1-1-1-0-0-0-1-0-1**

$\xi_{5,N} = 0, \tau_4 \neq \tau_1, \tau_4 \neq \tau_3, \tau_2 \neq \tau_1, \tau_4 = \tau_2, \tau_3 = \tau_1, w_{N,5,3} = 0, w_{N,5,4} \neq 0, w_{1,4,3} = 0, w_{1,3,1} \neq -w_{1,3,2}$

$$\begin{aligned} w_{N,5,1} &= 0 \\ w_{N,5,2} &= -w_{N,5,4} \\ w_{1,2,1} &= w_{1,4,1} + w_{1,4,3} \\ w_{1,4,2} &= 0 \\ n > 1,\ w_{2,2,1} &= w_{2,4,1} \\ n > 1,\ w_{2,4,2} &= w_{2,4,3} = 0 \\ n > 2,\ w_{3,2,1} &= w_{3,4,1} + w_{3,4,2} + w_{3,4,3} \end{aligned}$$

**0-1-1-1-0-0-0-1-1**

$\xi_{5,N} = 0, \tau_4 \neq \tau_1, \tau_4 \neq \tau_3, \tau_2 \neq \tau_1, \tau_4 = \tau_2, \tau_3 = \tau_1, w_{N,5,3} = 0, w_{N,5,4} \neq 0, w_{1,4,3} \neq 0$

$$\begin{aligned} w_{N,5,1} &= 0 \\ w_{N,5,2} &= -w_{N,5,4} \\ w_{1,2,1} &= w_{1,4,1} + w_{1,4,3} \\ w_{1,4,2} &= w_{1,3,1} = w_{1,3,2} = 0 \\ n > 1,\ w_{2,2,1} &= w_{2,4,1} + w_{2,4,3} \\ n > 1,\ w_{2,4,2} &= 0 \\ n > 1,\ w_{2,3,1} &= -w_{2,3,2} n > 2,\ w_{3,2,1} = w_{3,4,1} + w_{3,4,2} + w_{3,4,3} \end{aligned}$$

**0-1-1-1-0-0-1-0**

$\xi_{5,N} = 0, \tau_4 \neq \tau_1, \tau_4 \neq \tau_3, \tau_2 \neq \tau_1, \tau_4 = \tau_2, \tau_3 = \tau_1, w_{N,5,3} \neq 0, w_{N,5,4} = 0$

$$\begin{aligned} w_{N,5,1} &= -w_{N,5,3} \\ w_{N,5,2} &= w_{1,3,1} = w_{1,3,2} = 0 \\ n > 1,\ w_{2,3,1} &= w_{2,3,2} = 0 \\ n > 2,\ w_{3,3,1} &= -w_{3,3,2} \end{aligned}$$

**0-1-1-1-0-0-1-1**

$\xi_{5,N} = 0, \tau_4 \neq \tau_1, \tau_4 \neq \tau_3, \tau_2 \neq \tau_1, \tau_4 = \tau_2, \tau_3 = \tau_1, w_{N,5,3} \neq 0, w_{N,5,4} \neq 0$

$$w_{N,5,1} = -w_{N,5,3}$$
$$w_{N,5,2} = -w_{N,5,4}$$
$$w_{1,3,1} = w_{1,3,2} = w_{1,4,2} = 0$$
$$w_{1,2,1} = w_{1,4,1} + w_{1,4,3}$$
$$n > 1,\ w_{2,3,1} = -w_{2,3,2}$$
$$n > 1,\ w_{2,4,2} = -\frac{w_{N,5,3}}{w_{N,5,4}} w_{2,3,2}$$
$$n > 1,\ w_{2,4,1} = w_{2,2,1} - w_{2,4,3} + \frac{w_{N,5,3}}{w_{N,5,4}} w_{2,3,2}$$
$$n > 2,\ w_{3,3,1} = -\frac{w_{N,5,4}}{w_{N,5,3}} \left( w_{3,4,1} + w_{3,4,2} + w_{3,4,3} - w_{3,2,1} \right) - w_{3,3,2}$$

**0-1-1-1-0-1-0**

$\xi_{5,N} = 0, \tau_4 \neq \tau_1, \tau_4 \neq \tau_3, \tau_2 \neq \tau_1, \tau_4 = \tau_2, \tau_3 \neq \tau_1, w_{N,5,4} = 0$

$$w_{N,5,1} = w_{N,5,2} = w_{N,5,3} = 0$$

**0-1-1-1-0-1-1-0**

$\xi_{5,N} = 0, \tau_4 \neq \tau_1, \tau_4 \neq \tau_3, \tau_2 \neq \tau_1, \tau_4 = \tau_2, \tau_3 \neq \tau_1, w_{N,5,4} \neq 0, w_{1,3,1} = -w_{1,3,2} + \frac{\tau_3 - \tau_1}{\tau_2 - \tau_1} w_{1,2,1}$

$$w_{N,5,1} = w_{N,5,3} = w_{1,4,2} = w_{1,4,3} = 0$$
$$w_{N,5,2} = -w_{N,5,4}$$
$$w_{1,4,1} = w_{1,2,1}$$
$$n > 1,\ w_{2,2,1} = w_{2,4,1} + w_{2,4,3} \frac{\tau_2 - \tau_3}{\tau_2 - \tau_1}$$
$$n > 1,\ w_{2,4,2} = -w_{2,4,3} \frac{\tau_3 - \tau_1}{\tau_2 - \tau_1}$$
$$n > 2,\ w_{3,2,1} = w_{3,4,1} + w_{3,4,2} + w_{3,4,3}$$

**0-1-1-1-0-1-1-1**

$\xi_{5,N} = 0, \tau_4 \neq \tau_1, \tau_4 \neq \tau_3, \tau_2 \neq \tau_1, \tau_4 = \tau_2, \tau_3 \neq \tau_1, w_{N,5,4} \neq 0, w_{1,3,1} \neq -w_{1,3,2} + \frac{\tau_3 - \tau_1}{\tau_2 - \tau_1} w_{1,2,1}$

$$w_{N,5,1} = w_{N,5,3} = w_{1,4,2} = w_{1,4,3} = 0$$
$$w_{N,5,2} = -w_{N,5,4}$$
$$w_{1,4,1} = w_{1,2,1}$$
$$n > 1,\ w_{2,2,1} = w_{2,4,1}$$
$$n > 1,\ w_{2,4,2} = w_{2,4,3} = 0$$
$$n > 2,\ w_{3,2,1} = w_{3,4,1} + w_{3,4,2} + w_{3,4,3}$$

**0-1-1-1-1-0-0**

$\xi_{5,N} = 0, \tau_4 \neq \tau_1, \tau_4 \neq \tau_3, \tau_2 \neq \tau_1, \tau_4 \neq \tau_2, \tau_3 = \tau_1, w_{N,5,3} = 0$

$$w_{N,5,1} = w_{N,5,2} = 0$$

**0-1-1-1-1-0-1**

$\xi_{5,N} = 0, \tau_4 \neq \tau_1, \tau_4 \neq \tau_3, \tau_2 \neq \tau_1, \tau_4 \neq \tau_2, \tau_3 = \tau_1, w_{N,5,3} \neq 0$

$$\begin{aligned} w_{N,5,1} &= -w_{N,5,3} \\ w_{N,5,2} &= w_{1,3,1} = w_{1,3,2} = w_{N,5,4} = 0 \\ n > 1,\ w_{2,3,1} &= w_{2,3,2} = 0 \\ n > 2,\ w_{3,3,1} &= -w_{3,3,2} \end{aligned}$$

**0-1-1-1-1-1-0**

$\xi_{5,N} = 0, \tau_4 \neq \tau_1, \tau_4 \neq \tau_3, \tau_2 \neq \tau_1, \tau_4 \neq \tau_2, \tau_3 \neq \tau_1, w_{N,5,3} = 0$

$$w_{N,5,1} = w_{N,5,2} = w_{N,5,4} = 0$$

**0-1-1-1-1-1-1**

$\xi_{5,N} = 0, \tau_4 \neq \tau_1, \tau_4 \neq \tau_3, \tau_2 \neq \tau_1, \tau_4 \neq \tau_2, \tau_3 \neq \tau_1, w_{N,5,3} \neq 0$

$$\begin{aligned} w_{N,5,1} &= w_{N,5,4} = w_{1,3,2} = 0 \\ w_{N,5,2} &= -w_{N,5,3} \\ w_{1,3,1} &= w_{1,2,1} \\ \tau_3 &= \tau_2 \\ n > 1,\ w_{2,3,1} &= w_{2,2,1} \\ n > 1,\ w_{2,3,2} &= 0 \\ n > 2,\ w_{3,2,1} &= w_{3,3,1} + w_{3,3,2} \end{aligned}$$

**1-0**

$\xi_{5,N} \neq 0, \tau_2 = \tau_1$

We now list the main operations performed to solve this system.

1) 4.5 to find $\tau_1$
2) 3.3 in 4.2
3) 4.3 to find $\tau_3$ and $\tau_4$
4) 4.4 to find contradiction

**1-1-0-0**

$\xi_{5,N} \neq 0, \tau_2 \neq \tau_1, \tau_3 = \tau_1, \tau_1 = 0$

We now list the main operations performed to solve this system.

1) multiply 4.10 $\tau_2 - \tau_1$
2) use 4.5 in 4.10
3) multiply 2.2 by $\tau_2 - \tau_1$
4) 2.1 in 2.2

5) 3.3 in 2.2
6) multiply 3.4 by $\tau_2 - \tau_1$
7) 3.3 in 3.4
8) 4.6 in 3.4
9) 3.1 in 4.1
10) 4.10 in 3.4
11) 4.2 in 3.4

Which yields

$$
\begin{aligned}
& w_{1,3,2} \neq 0 \\
& w_{N,5,1} = \xi_{5,N}^{N} \frac{N}{2(1+N)(2+N)} \left(1 + N - \frac{\xi_{5,N}}{w_{1,3,2}}\right) \\
& w_{N,5,2} = \xi_{5,N}^{N} \frac{N(3+N)^2}{2(1+N)^2(2+N)} \\
& w_{N,5,3} = \xi_{5,N}^{N+1} \frac{N}{2w_{1,3,2}(1+N)(2+N)} \\
& w_{N,5,4} = 2\frac{\xi_{5,N}^{N}}{(1+N)^2(2+N)} \\
& \tau_2 = \xi_{5,N} \frac{2}{3+N} \\
& \tau_4 = \xi_{5,N} \\
& w_{1,2,1} = \xi_{5,N} \frac{2}{3+N} \\
& w_{1,3,1} = -w_{1,3,2} \\
& w_{1,4,1} = \xi_{5,N} \frac{1+N}{4} \left(1 - \frac{\xi_{5,N}}{w_{1,3,2}}\right) - \xi_{5,N} N \\
& w_{1,4,2} = \xi_{5,N} 3 \frac{1+N}{4} \\
& w_{1,4,3} = \xi_{5,N}^{2} \frac{1+N}{4w_{1,3,2}} \\
n > 1,\ & w_{2,3,1} = \frac{w_{1,3,2}}{\xi_{5,N}} \left(\xi_{5,N}^{2} \frac{4}{3+N} - (3+N) w_{2,2,1}\right) - w_{2,3,2} \\
n > 1,\ & w_{2,4,1} = \xi_{5,N}^{2} \frac{(1+N)(3-N)}{2(3+N)} - \frac{N(3+N)}{2} w_{2,2,1} + \xi_{5,N} \frac{1+N}{4} N \frac{w_{2,3,2}}{w_{1,3,2}} - w_{2,4,3} \\
n > 1,\ & w_{2,4,2} = \xi_{5,N} \frac{1+N}{4} \left(2\xi_{5,N} - N \frac{w_{2,3,2}}{w_{1,3,2}}\right) \\
n > 2,\ & w_{3,4,1} = \frac{1+N}{2} \left(\xi_{5,N}^{3} \frac{6}{3+N} - \frac{N(3+N)^2}{2(1+N)} w_{3,2,1} - \xi_{5,N} \frac{N}{2w_{1,3,2}} (w_{3,3,1} + w_{3,3,2})\right) - w_{3,4,2} - w_{3,4,3}
\end{aligned}
$$

**1-1-0-1**

$\xi_{5,N} \neq 0, \tau_2 \neq \tau_1, \tau_3 = \tau_1, \tau_1 \neq 0$

We now list the main operations performed to solve this system.

1) multiply 4.10 $\tau_2 - \tau_1$
2) use 4.5 in 4.10
3) multiply 2.2 by $\tau_2 - \tau_1$
4) 2.1 in 2.2
5) 3.3 in 2.2
6) multiply 3.4 by $\tau_2 - \tau_1$
7) 3.3 in 3.4
8) 4.6 in 3.4
9) 3.1 in 4.1
10) 4.10 in 3.4
11) 4.2 in 3.4
12) expression of $\tau_1$ in 3.4, replace where fits
13) 4.2 for $\tau_2$
14) 4.1 for $\tau_4$
15) 3.1 for contradiction.

**1-1-1-0**

$\xi_{5,N} \neq 0, \tau_2 \neq \tau_1, \tau_3 \neq \tau_1, \tau_3 = \tau_2$

We now list the main operations performed to rewrite this system.

1) multiply 4.10 $\tau_2 - \tau_1$
2) use 4.5 in 4.10
3) multiply 2.2 by $\tau_2 - \tau_1$
4) 2.1 in 2.2
5) 3.3 in 2.2
6) multiply 3.4 by $\tau_2 - \tau_1$
7) 3.3 in 3.4
8) 4.6 in 3.4
9) $\tau_2$ in 4.6
10) $\tau_1$ in 4.2
11) 3.1 in 4.1 to find $\tau_4$

Which yields :

$$
\begin{aligned}
w_{1,3,2} &\neq 0 \\
w_{N,5,1} &= \xi_{5,N}^{N} \frac{N}{2(2+N)} \\
w_{N,5,2} &= \xi_{5,N}^{N} \frac{N(3+N)}{2(1+N)^2(2+N)} \left( 3 + N - \frac{\xi_{5,N}}{w_{1,3,2}} \right) \\
w_{N,5,3} &= \xi_{5,N}^{N+1} \frac{N(3+N)}{2w_{1,3,2}(1+N)^2(2+N)} \\
w_{N,5,4} &= \xi_{5,N}^{N} \frac{2}{(1+N)^2(2+N)} \\
\tau_1 &= 0 \\
\tau_2 &= \xi_{5,N} \frac{2}{3+N} \\
\tau_4 &= \xi_{5,N} \\
w_{1,2,1} &= \xi_{5,N} \frac{2}{3+N} \\
w_{1,3,1} &= \xi_{5,N} \frac{2}{3+N} - w_{1,3,2} \\
w_{1,4,1} &= \xi_{5,N} \frac{1-N}{4} \\
w_{1,4,2} &= \frac{\xi_{5,N}}{4} \left( 3 + N - \xi_{5,N} \frac{1+N}{w_{1,3,2}} \right) \\
w_{1,4,3} &= \xi_{5,N}^{2} \frac{1+N}{4w_{1,3,2}} \\
n > 1,\ w_{2,3,1} &= \frac{w_{1,3,2}}{\xi_{5,N}} \left( \xi_{5,N}^{2} \frac{4}{3+N} - \left( 3 + N - \frac{\xi_{5,N}}{w_{1,3,2}} \right) w_{2,2,1} \right) - w_{2,3,2} \\
n > 1,\ w_{2,4,1} &= \xi_{5,N}^{2} \frac{1-N}{2} + \xi_{5,N} \frac{N(3+N)}{4w_{1,3,2}} w_{2,3,2} \\
n > 1,\ w_{2,4,2} &= \xi_{5,N}^{2} \frac{1+N}{2} - \xi_{5,N} \frac{N(3+N)}{4w_{1,3,2}} w_{2,3,2} - w_{2,4,3} \\
n > 2,\ w_{3,4,1} &= \xi_{5,N}^{3} 3 \frac{1+N}{3+N} - \frac{N(3+N)}{4} \left( 3 + N - \frac{\xi_{5,N}}{w_{1,3,2}} \right) w_{3,2,1} \\
&\quad - \xi_{5,N} \frac{N(3+N)}{4w_{1,3,2}} \left( w_{3,3,1} + w_{3,3,2} \right) - w_{3,4,2} - w_{3,4,3}
\end{aligned}
$$

**1-1-1-1-0**

$\xi_{5,N} \neq 0, \tau_2 \neq \tau_1, \tau_3 \neq \tau_1, \tau_3 \neq \tau_2, \tau_4 = \tau_3$

We now list the main operations performed to solve this system.

1) multiply 4.10 $\tau_2 - \tau_1$
2) use 4.5 in 4.10
3) multiply 2.2 by $\tau_2 - \tau_1$
4) 2.1 in 2.2
5) 3.3 in 2.2
6) multiply 3.4 by $\tau_2 - \tau_1$
7) 3.3 in 3.4
8) 4.6 in 3.4
9) 4.7 for $\tau_3$
10) 4.3 for $\tau_1$
11) 4.1 for contradiction ($\tau_1 = \tau_2$)

**1-1-1-1-1-0**

$\xi_{5,N} \neq 0, \tau_2 \neq \tau_1, \tau_3 \neq \tau_1, \tau_3 \neq \tau_2, \tau_4 \neq \tau_3, \tau_4 = \tau_1, \tau_1 = 0$

We now list the main operations performed to solve this system.

1) multiply 4.10 $\tau_2 - \tau_1$
2) use 4.5 in 4.10
3) multiply 2.2 by $\tau_2 - \tau_1$
4) 2.1 in 2.2
5) multiply 2.2 by $\tau_3 - \tau_2$
6) use 3.1 in 2.2
7) 3.2 in 2.2
8) 4.4 in 2.2
9) 4.10 in 2.2
10) multiply 3.2 by $\tau_3 - \tau_1$
11) 4.4 in 3.2
12) multiply 3.3 by $\tau_3 - \tau_1$
13) 4.3 in 3.3
14) multiply 3.4 by $\tau_3 - \tau_1$
15) 4.7 in 3.4
16) 4.10 in 3.4
17) 3.3 in 3.4
18) if $\tau_1 = \xi_{5,N}$ it is not too difficult to use 2.2 and 4.1 to find $\tau_3 = \tau_1$ which is a contradiction, thus $\tau_1 \neq \xi_{5,N}$, which using 3.4 yields $\tau_1 = 0$, replace everywhere
19) 2.2 in 4.1
20) replace $w_{1,2,1}$ from 4.10 everywhere
21) 4.3 in 4.7 to find $w_{1,3,1} + w_{1,3,2} = \tau_3$, replace everywhere
22) 3.1 in 3.2
23) 4.2 in 4.6
24) 4.2 in 4.11
25) 4.3 in 4.8

26) 3.1 in 3.5
27) 3.5 in 4.13
28) 2.2 in 4.13
29) 4.8 in 4.12
30) 4.3 in 4.12 to find $w_{1,4,1} + w_{1,4,2} + w_{1,4,3} = 0$, replace everywhere
31) 2.2 in 4.9
32) multiply 2.1 by $\tau_3 - \tau_2$
33) 3.1 in 2.1
34) multiply 3.3 by $\tau_3 - \tau_2$
35) 3.1 in 3.3
36) multiply 4.3 by $\tau_3 - \tau_2$
37) 4.2 in 4.3
38) multiply 4.5 by $(\tau_3 - \tau_2)\tau_3$
39) 4.2 in 4.5
40) 3.3 in 4.5
41) 4.5 in 3.3
42) 2.2 in 3.3
43) we get $\tau_2 = \tau_1 = 0$ which is a contradiction.

**1-1-1-1-1-1-0**

$\xi_{5,N} \neq 0, \tau_2 \neq \tau_1, \tau_3 \neq \tau_1, \tau_3 \neq \tau_2, \tau_4 \neq \tau_3, \tau_4 \neq \tau_1, \tau_4 = \tau_2$

We now list the main operations performed to solve this system.

1) 4.4 in 4.9
2) 4.4 in 4.1
3) multiply 3.3 by $\tau_3 - \tau_2$
4) 4.3 in 3.3
5) multiply 3.4 by $\tau_3 - \tau_2$
6) 4.7 in 3.4
7) multiply 2.2 by $\tau_2 - \tau_1$
8) 2.1 in 2.2
9) multiply 2.2 by $\tau_3 - \tau_2$
10) 3.1 in 2.2
11) 3.2 in 2.2
12) multiply 4.10 by $\tau_2 - \tau_1$
13) 4.5 in 4.10
14) 4.10 in 3.4
15) 3.3 in 3.4
16) 4.4 in 4.1
17) if $\tau_1 = 0$ we can easily find $\tau_2 = \tau_1$ which is a contradiction, thus $\tau_1 \neq 0$, 4.1 to find $\tau_3 = 4\frac{\xi_{5,N}}{3+N} - \tau_2$
18) 4.1 in 4.4
19) we can use 4.4 to verify that $\xi_{5,N}\frac{2}{2+N} \neq \tau_2$, which can be used in 4.9 to find $w_{1,3,1} + w_{1,3,2} = \tau_3$, replace everywhere
20) 4.2 in 4.6
21) 4.6 in 4.2
22) 3.1 in 3.2

23) 3.2 in 3.1
24) 3.2 in 3.3
25) multiply 4.5 by $\tau_3 - \tau_2$
26) 4.6 in 4.5
27) 4.5 in 3.3
28) 4.1 in 3.3
29) 3.4 in 3.3
30) 4.1 in 4.2
31) 4.4 in 4.2
32) find contradiction

**1-1-1-1-1-1-1**

$\xi_{5,N} \neq 0, \tau_2 \neq \tau_1, \tau_3 \neq \tau_1, \tau_3 \neq \tau_2, \tau_4 \neq \tau_3, \tau_4 \neq \tau_1, \tau_4 \neq \tau_2$

We now list the main operations performed to solve this system.

1) multiply 2.2 by $\tau_2 - \tau_1$
2) 2.1 in 2.2
3) multiply 2.2 by $\tau_3 - \tau_2$
4) 3.1 in 2.2
5) 3.2 in 2.2
6) 4.4 in 2.2
7) multiply 2.1 by $\tau_3 - \tau_2$
8) 3.1 in 2.1
9) multiply 2.1 by $\tau_4 - \tau_2$
10) 4.1 in 2.1
11) multiply 3.1 by $\tau_4 - \tau_3$
12) 4.1 in 3.1
13) multiply 3.2 by $\tau_3 - \tau_4$
14) 4.4 in 3.2
15) multiply 3.3 by $\tau_3 - \tau_4$
16) 4.3 in 3.3
17) multiply 3.4 by $\tau_3 - \tau_4$
18) 4.7 in 3.4
19) multiply 4.10 by $\tau_2 - \tau_1$
20) 4.5 in 4.10
21) 4.10 in 3.4
22) 3.3 in 3.4
23) multiply 4.5 by $(\tau_3 - \tau_2)(\tau_3 - \tau_1)$
24) 4.2 in 4.3
25) multiply 3.3 by $(\tau_3 - \tau_2)(\tau_3 - \tau_1)$
26) 3.1 in 3.3
27) 3.3 + $\tau_4(3 + N)$ 4.5
28) 3.4 in 3.3
29) 4.5 in 3.3
30) 3.3 - 3.4
31) 3.3 yields $\tau_4 = \xi_{5,N}$, replace everywhere

32) multiply 3.2 by $\tau_3 - \tau_1$
33) 3.1 in 3.2
34) multiply 4.6 by $\tau_3 - \tau_1$
35) 4.2 in 4.6
36) 3.2 in 4.6
37) multiply 2.2 by $\xi_{5,N} - \tau_2$
38) 4.4 in 2.2
39) 4.1 in 2.2
40) 4.10 in 2.2
41) find $\tau_1 = 0$, replace everywhere
42) 4.3 in 4.7 to find $w_{1,3,1} + w_{1,3,2} = \tau_3$, replace everywhere
43) 4.2 in 4.11
44) 4.8 in 4.12
45) 4.1 in 4.4 to find $w_{1,4,1} + w_{1,4,2} + w_{1,4,3} = \xi_{5,N}$, replace everywhere
46) 4.3 in 4.8
47) 4.1 in 4.9
48) 4.1 in 4.13
49) multiply 3.5 by $\tau_3 - \xi_{5,N}$
50) 3.1 in 3.5
51) 4.1 in 3.5
52) multiply 4.3 by $\tau_3 - \tau_2$
53) 4.2 in 4.3
54) 4.14 in 4.16
55) 4.15 in 4.18
56) multiply 4.2 by $(\xi_{5,N} - \tau_3)(\xi_{5,N} - \tau_2)$
57) 4.1 in 4.2
58) multiply 4.3 by $\xi_{5,N} - \tau_2$
59) 4.1 in 4.3
60) 2.1 in 1.1
61) 3.1 in 1.1
62) 4.1 in 1.1
63) 4.1 in 4.19
64) 3.1 in 4.19
65) 2.1 in 4.19
66) isolate $w_{3,3,1}$ in 4.19
67) isolate $w_{N,5,3}(w_{2,3,1} + w_{2,3,2})$ in 3.6
68) 3.6 in 4.15
69) multiply 3.6 by $\tau_3 - \xi_{5,N}$
70) 4.15 in 3.6
71) multiply 4.17 by $w_{N,5,3}(\tau_3 - \xi_{5,N})$
72) 3.6 in 4.17
73) replace everything in 4.17 and 4.15
74) isolate $w_{2,4,1}$ in 4.15
75) multiply 3.6 by $\tau_2$
76) 4.14 in 3.6
77) replace everything in 3.6 and 4.14

Which yields :

$$\xi_{5,N}2 \neq \tau_2(3+N)$$

$$\xi_{5,N}\frac{2}{2+N}\left(\xi_{5,N}\frac{3}{3+N}-\tau_2\right) \neq \left(\xi_{5,N}\frac{2}{2+N}-\tau_2\right)\tau_3$$

$$\begin{aligned}
a &= \big(\xi_{5,N}2\big(\xi_{5,N}3-\tau_2(3+N)\big)-\big(\xi_{5,N}2-\tau_2(2+N)\big)\tau_3(3+N)\big) \neq 0\\
w_{N,5,1} &= \xi_{5,N}^{N}-\frac{\xi_{5,N}^{N}}{(1+N)(2+N)(3+N)}\left(\xi_{5,N}^{2}\frac{N}{\tau_3-\tau_2}\left(\frac{\tau_3(3+N)-\xi_{5,N}2}{\tau_2\big(\xi_{5,N}-\tau_2\big)}-\frac{\tau_2(3+N)-\xi_{5,N}2}{\tau_3\big(\xi_{5,N}-\tau_3\big)}\right)\right.\\
&\left.+\frac{a}{\big(\xi_{5,N}-\tau_3\big)\big(\xi_{5,N}-\tau_2\big)}\right)\\
w_{N,5,2} &= \xi_{5,N}^{2+N}N\frac{\tau_3(3+N)-\xi_{5,N}2}{(1+N)(2+N)(3+N)\tau_2(\tau_3-\tau_2)\big(\xi_{5,N}-\tau_2\big)}\\
w_{N,5,3} &= \xi_{5,N}^{2+N}N\frac{\xi_{5,N}2-\tau_2(3+N)}{(1+N)(2+N)(3+N)(\tau_3-\tau_2)\tau_3\big(\xi_{5,N}-\tau_3\big)}\\
w_{N,5,4} &= \xi_{5,N}^{N}\frac{a}{(1+N)(2+N)(3+N)\big(\xi_{5,N}-\tau_3\big)\big(\xi_{5,N}-\tau_2\big)}\\
\tau_1 &= 0\\
\tau_4 &= \xi_{5,N}\\
w_{1,2,1} &= \tau_2\\
w_{1,3,1} &= \tau_3-\frac{\xi_{5,N}(\tau_3-\tau_2)\tau_3}{\tau_2\big(\xi_{5,N}2-\tau_2(3+N)\big)}\\
w_{1,3,2} &= \frac{\xi_{5,N}(\tau_3-\tau_2)\tau_3}{\tau_2\big(\xi_{5,N}2-\tau_2(3+N)\big)}\\
w_{1,4,1} &= \xi_{5,N}-\xi_{5,N}^{2}\big(\xi_{5,N}-\tau_2\big)\frac{5\xi_{5,N}\tau_3+\xi_{5,N}\tau_2(3+N)-\tau_3(\tau_2+\tau_3)(3+N)-2\xi_{5,N}^{2}}{\tau_2\tau_3 a}\\
w_{1,4,2} &= \xi_{5,N}^{2}\big(\xi_{5,N}-\tau_2\big)\frac{\xi_{5,N}5\tau_3-(3+N)\tau_3^{2}-2\xi_{5,N}^{2}+\tau_2\xi_{5,N}N}{(\tau_3-\tau_2)\tau_2 a}\\
w_{1,4,3} &= \frac{\xi_{5,N}^{2}\big(\xi_{5,N}2-\tau_2(3+N)\big)\big(\xi_{5,N}-\tau_3\big)\big(\xi_{5,N}-\tau_2\big)}{(\tau_3-\tau_2)\tau_3 a}\\
n>1,\ w_{2,3,1} &= \frac{(\tau_3-\tau_2)\tau_3}{\xi_{5,N}^{2}N\big(\xi_{5,N}2-\tau_2(3+N)\big)\tau_2}\left[-\xi_{5,N}^{2}N\frac{\tau_3(3+N)-\xi_{5,N}2}{(\tau_3-\tau_2)}w_{2,2,1}\right.\\
&\left.+\frac{a}{\big(\xi_{5,N}-\tau_2\big)}\big(w_{2,4,2}\tau_2+w_{2,4,3}\tau_3\big)+\xi_{5,N}^{3}2\big(N\tau_2+\tau_3-\xi_{5,N}\big)\right]\\
n>1,\ w_{2,3,2} &= \frac{(\tau_3-\tau_2)\tau_3}{\xi_{5,N}^{2}N\big(\xi_{5,N}2-\tau_2(3+N)\big)\tau_2}\left(\xi_{5,N}^{3}2\big(\xi_{5,N}-\tau_3\big)-\frac{a}{\big(\xi_{5,N}-\tau_2\big)}\big(w_{2,4,2}\tau_2+w_{2,4,3}\tau_3\big)\right)
\end{aligned}$$

$$n>1,\ w_{2,4,1}=\frac{\xi_{5,N}^{2}}{a\tau_2}\big(2\tau_2\big(\xi_{5,N}-\tau_2\big)\big(\xi_{5,N}3-\tau_3(3+N)\big)+N\big(\tau_3(3+N)-\xi_{5,N}2\big)w_{2,2,1}\big)-w_{2,4,2}-w_{2,4,3}$$

$$n>2,\ w_{3,3,1}=\frac{(\tau_3-\tau_2)\tau_3\big(\xi_{5,N}-\tau_3\big)}{\xi_{5,N}^{2}N\big(\xi_{5,N}2-\tau_2(3+N)\big)}\Bigg(\xi_{5,N}^{3}6-\xi_{5,N}^{2}N\frac{\tau_3(3+N)-\xi_{5,N}2}{\tau_2(\tau_3-\tau_2)\big(\xi_{5,N}-\tau_2\big)}w_{3,2,1}$$

$$-\frac{a}{\big(\xi_{5,N}-\tau_3\big)\big(\xi_{5,N}-\tau_2\big)}\big(w_{3,4,1}+w_{3,4,2}+w_{3,4,3}\big)\Bigg)-w_{3,3,2}$$

### C.7 **Implicit method**, $s = 1$, $j = 1$, $v = 2$

The conditions are :

Order of consistency 1 :

$$w_{N,1,1} = \xi_{1,N}^{N}$$

Order of consistency 2 :

$$w_{N,1,1}\tau_1 = \xi_{1,N}^{N+1}\frac{1}{1+N}$$

$$w_{N,1,1}w_{1,1,1} = \xi_{1,N}^{N+1}\frac{1}{1+N}$$

We now list the main operations performed to solve this system.

Order of consistency 1 :

$$w_{N,1,1} = \xi_{1,N}^{N}$$

Order of consistency 2 :

$$\xi_{1,N}\left(\xi_{1,N}\frac{1}{1+N} - \tau_1\right) = 0$$

$$\xi_{1,N}\left(w_{1,1,1} - \tau_1\right) = 0$$

The solutions are :

**0**

$\xi_{1,N} = 0$

$$w_{N,1,1} = \xi_{1,N}^{N}$$

**1**

$\xi_{1,N} \neq 0$

$$w_{N,1,1} = \xi_{1,N}^{N}$$

$$\tau_1 = w_{1,1,1} = \xi_{1,N}\frac{1}{1+N}$$

### C.8 Implicit method, $s = 1$, $j = 2$, $v = 2$

This case is equivalent to the case $s = 1, j = 1, v = 2$ thanks to equivalence by permutation. The solutions hence are :

**0**

$\xi_{2,N} = 0$

$$w_{N,2,1} = 0$$

**1**

$\xi_{2,N} \neq 0$

$$\begin{aligned} w_{1,1,1} = \tau_1 = \xi_{2,N}\frac{1}{1+N} \\ w_{N,2,1} = \xi_{2,N}^{N} \end{aligned}$$

### C.9 Implicit method, $s = 2$, $j = 1$, $v = 2$

This case is equivalent to the case $s = 2, j = 2, v = 2$ thanks to equivalence by permutation. The solutions are :

**0-0**

$\xi_{1,N} = 0, \tau_1 = \tau_2$

$$w_{N,1,2} = -w_{N,1,1}$$

**0-1**

$\xi_{1,N} = 0, \tau_1 \neq \tau_2$

$$w_{N,1,2} = w_{N,1,1} = 0$$

**1-0**

$\xi_{1,N} \neq 0, \tau_1 = \tau_2$

$$w_{N,1,2} = \xi_{1,N}^N - w_{N,1,1}$$

$$\tau_1 = \frac{\xi_{1,N}}{1+N}$$

**1-1**

$\xi_{1,N} \neq 0, \tau_1 \neq \tau_2$

$$w_{N,1,1} = \frac{\xi_{1,N}^N}{\tau_1 - \tau_2}\left(\frac{\xi_{1,N}}{1+N} - \tau_2\right)$$

$$w_{N,1,2} = \frac{\xi_{1,N}^N}{\tau_1 - \tau_2}\left(\tau_1 - \frac{\xi_{1,N}}{1+N}\right)$$

### C.10 Implicit method, $s = 2, j = 2, v = 2$

The conditions are :

Order of consistency 1 :

$$w_{N,2,1} + w_{N,2,2} = \xi_{2,N}^{N}$$

Order of consistency 2 :

$$w_{N,2,1}\tau_1 + w_{N,2,2}\tau_2 = \xi_{2,N}^{N+1}\frac{1}{1+N}$$

$$w_{N,2,1}(w_{1,1,1} + w_{1,1,2}) + w_{N,2,2}(w_{1,2,1} + w_{1,2,2}) = \xi_{2,N}^{N+1}\frac{1}{1+N}$$

To simplify the system of equations, we assume that the method satisfies the conditions for order of consistency 1 at all stages and at rank 1.

Order of consistency 1 :

$$w_{N,2,2} = \xi_{2,N}^{N} - w_{N,2,1}$$

Order of consistency 2 :

$$w_{N,2,1}(\tau_1 - \tau_2) = \xi_{2,N}^{N}\left(\frac{\xi_{2,N}}{1+N} - \tau_2\right)$$

$$0 = 0$$

The solutions are :

**0-0**

$\xi_{2,N} = 0, \tau_1 = \tau_2$

$$w_{N,2,2} = -w_{N,2,1}$$

**0-1**

$\xi_{2,N} = 0, \tau_1 \neq \tau_2$

$$w_{N,2,2} = w_{N,2,1} = 0$$

**1-0**

$\xi_{2,N} \neq 0, \tau_1 = \tau_2$

$$w_{N,2,2} = \xi_{2,N}^{N} - w_{N,2,1}$$

$$\tau_2 = \frac{\xi_{2,N}}{1+N}$$

**1-1**

$\xi_{2,N} \neq 0, \tau_1 \neq \tau_2$

$$w_{N,2,2} = \frac{\xi_{2,N}^N}{\tau_1 - \tau_2}\left(\tau_1 - \frac{\xi_{2,N}}{1+N}\right)$$

$$w_{N,2,1} = \frac{\xi_{2,N}^N}{\tau_1 - \tau_2}\left(\frac{\xi_{2,N}}{1+N} - \tau_2\right)$$

### C.11 Implicit method, $s = 2, j = 3, v = 2$

The conditions are :

Order of consistency 1 :

$$w_{N,3,1} + w_{N,3,2} = \xi_{3,N}^{N}$$

Order of consistency 2 :

$$w_{N,3,1}\tau_1 + w_{N,3,2}\tau_2 = \xi_{3,N}^{N+1}\frac{1}{1+N}$$

$$w_{N,3,1}(w_{1,1,1} + w_{1,1,2}) + w_{N,3,2}(w_{1,2,1} + w_{1,2,2}) = \xi_{3,N}^{N+1}\frac{1}{1+N}$$

To simplify the system we assume that the method satisfy the conditions for order of consistency 1 at all stages and rank 1.

Order of consistency 1 :

$$w_{N,3,1} + w_{N,3,2} = \xi_{3,N}^{N}$$

Order of consistency 2 :

$$w_{N,3,1}\tau_1 + w_{N,3,2}\tau_2 = \xi_{3,N}^{N+1}\frac{1}{1+N}$$

$$0 = 0$$

The solutions are :

**0-0**

$\xi_{3,N} = 0, \tau_2 = \tau_1$

$$w_{N,3,2} = -w_{N,3,1}$$

**0-1**

$\xi_{3,N} = 0, \tau_2 \neq \tau_1$

$$w_{N,3,2} = w_{N,3,1} = 0$$

**1-0**

$\xi_{3,N} \neq 0, \tau_1 = \tau_2$

$$w_{N,3,2} = \xi_{3,N}^{N} - w_{N,3,1}$$

$$\tau_2 = \frac{\xi_{3,N}}{1+N}$$

**1-1**

$\xi_{3,N} \neq 0, \tau_1 \neq \tau_2$

$$w_{N,3,1} = \frac{\xi_{3,N}^N}{\tau_1 - \tau_2}\left(\frac{\xi_{3,N}}{1+N} - \tau_2\right)$$
$$w_{N,3,2} = \frac{\xi_{3,N}^N}{\tau_1 - \tau_2}\left(\tau_1 - \frac{\xi_{3,N}}{1+N}\right)$$

### C.12 Implicit method, $s = 2, j = 3, v = 3$

The conditions are :

Order of consistency 1 :

$$w_{N,3,1} + w_{N,3,2} = \xi_{3,N}^{N}$$

Order of consistency 2 :

$$w_{N,3,1}\tau_1 + w_{N,3,2}\tau_2 = \xi_{3,N}^{N+1}\frac{1}{1+N}$$

$$w_{N,3,1}(w_{1,1,1} + w_{1,1,2}) + w_{N,3,2}(w_{1,2,1} + w_{1,2,2}) = \xi_{3,N}^{N+1}\frac{1}{1+N}$$

Order of consistency 3 :

$$w_{N,3,1}\tau_1^2 + w_{N,3,2}\tau_2^2 = \xi_{3,N}^{N+2}2\frac{N!}{(2+N)!}$$

$$w_{N,3,1}\tau_1(w_{1,1,1} + w_{1,1,2}) + w_{N,3,2}\tau_2(w_{1,2,1} + w_{1,2,2}) = \xi_{3,N}^{N+2}2\frac{N!}{(2+N)!}$$

$$w_{N,3,1}(w_{1,1,1}\tau_1 + w_{1,1,2}\tau_2) + w_{N,3,2}(w_{1,2,1}\tau_1 + w_{1,2,2}\tau_2) = \xi_{3,N}^{N+2}\frac{N!}{(2+N)!}$$

$$w_{N,3,1}(w_{1,1,1}(w_{1,1,1} + w_{1,1,2}) + w_{1,1,2}(w_{1,2,1} + w_{1,2,2}))$$
$$+w_{N,3,2}(w_{1,2,1}(w_{1,1,1} + w_{1,1,2}) + w_{1,2,2}(w_{1,2,1} + w_{1,2,2})) = \xi_{3,N}^{N+2}\frac{N!}{(2+N)!}$$

$$w_{N,3,1}(w_{1,1,1} + w_{1,1,2})^2 + w_{N,3,2}(w_{1,2,1} + w_{1,2,2})^2 = \xi_{3,N}^{N+2}2\frac{N!}{(2+N)!}$$

$$\tilde{n} > 1 \; w_{N,3,1}(w_{2,1,1} + w_{2,1,2}) + w_{N,3,2}(w_{2,2,1} + w_{2,2,2}) = \xi_{3,N}^{N+2}2\frac{N!}{(2+N)!}$$

To simplify the system of equations we assume that the method satisfies the condition for order of consistency 1 at all stages at rank 1 and 2.

$$w_{1,1,1} = \tau_1 - w_{1,1,2}$$
$$w_{1,2,2} = \tau_2 - w_{1,2,1}$$

Order of consistency 1 :

$$w_{N,3,1} + w_{N,3,2} = \xi_{3,N}^{N}$$

Order of consistency 2 :

$$w_{N,3,1}\tau_1 + w_{N,3,2}\tau_2 = \xi_{3,N}^{N+1}\frac{1}{1+N}$$
$$0 = 0$$

Order of consistency 3 :

$$w_{N,3,1}\tau_1^2 + w_{N,3,2}\tau_2^2 = \xi_{3,N}^{N+2}2\frac{N!}{(2+N)!}$$

$$0 = 0$$
$$w_{N,3,1}\left(w_{1,1,1}\tau_1 + w_{1,1,2}\tau_2\right) + w_{N,3,2}\left(w_{1,2,1}\tau_1 + w_{1,2,2}\tau_2\right) = \xi_{3,N}^{N+2}\frac{N!}{(2+N)!}$$
$$0 = 0$$
$$0 = 0$$
$$0 = 0$$

We now list the main operations performed to solve this system.

1) each order of consistency 1 in 3.3
2) 3.1 in 3.3
3) isolate $w_{N,3,1}\tau_1$ in 2.1
4) 2.1 in 3.1
5) isolate $w_{N,3,1}$ in 1.1
6) 1.1 in 2.1, 3.3
7) isolate $w_{N,3,2}(\tau_2 - \tau_1)$ in 2.1
8) 2,1 in 3.1, 3.3

Which yields :

$$w_{1,1,1} = \tau_1 - w_{1,1,2}$$
$$w_{1,2,2} = \tau_2 - w_{1,2,1}$$

Order of consistency 1 :

$$w_{N,3,1} = \xi_{3,N}^{N} - w_{N,3,2}$$

Order of consistency 2 :

$$w_{N,3,2}(\tau_2 - \tau_1) = \xi_{3,N}^{N}\left(\frac{\xi_{3,N}}{1+N} - \tau_1\right)$$
$$0 = 0$$

Order of consistency 3 :

$$\xi_{3,N}\left[\left(\frac{\xi_{3,N}}{1+N} - \tau_1\right)\tau_2 - \frac{\xi_{3,N}}{1+N}\left(\xi_{3,N}\frac{2}{2+N} - \tau_1\right)\right] = 0$$
$$0 = 0$$
$$\xi_{3,N}\left[\left(\frac{\xi_{3,N}}{1+N} - \tau_2\right)w_{1,1,2} + \left(\frac{\xi_{3,N}}{1+N} - \tau_1\right)w_{1,2,1} - \xi_{3,N}^{2}\frac{N!}{(2+N)!}\right] = 0$$
$$0 = 0$$
$$0 = 0$$
$$0 = 0$$

The solutions ares :

**0-0**

$\xi_{3,N} = 0, \tau_2 = \tau_1$

$$\begin{aligned} w_{1,1,1} &= \tau_1 - w_{1,1,2} \\ w_{1,2,2} &= \tau_1 - w_{1,2,1} \\ w_{N,3,1} &= -w_{N,3,2} \end{aligned}$$

**0-1**

$\xi_{3,N} = 0, \tau_2 \neq \tau_1$

$$\begin{aligned} w_{1,1,1} &= \tau_1 - w_{1,1,2} \\ w_{1,2,2} &= \tau_2 - w_{1,2,1} \\ w_{N,3,1} &= w_{N,3,2} = 0 \end{aligned}$$

**1-0**

$\xi_{3,N} \neq 0, \tau_1 = \frac{\xi_{3,N}}{1+N}$

Impossible.

**1-1**

$\xi_{3,N} \neq 0, \tau_1 \neq \frac{\xi_{3,N}}{1+N}$

$$\begin{aligned} w_{N,3,2} &\neq 0 \\ \tau_2 &\neq \tau_1 \\ w_{1,1,1} &= \tau_1 - w_{1,1,2} \\ w_{1,2,2} &= \frac{1}{\xi_{3,N} - \tau_1(1+N)}\left[\xi_{3,N}\left(\frac{\xi_{3,N}}{2+N} - \tau_1\right) + \big(\xi_{3,N} - \tau_2(1+N)\big)w_{1,1,2}\right] \\ w_{N,3,1} &= \xi_{3,N}^N - \frac{\xi_{3,N}^N}{\tau_2 - \tau_1}\left(\frac{\xi_{3,N}}{1+N} - \tau_1\right) \\ w_{N,3,2} &= \frac{\xi_{3,N}^N}{\tau_2 - \tau_1}\left(\frac{\xi_{3,N}}{1+N} - \tau_1\right) \\ \tau_2 &= \frac{\xi_{3,N}}{\xi_{3,N} - \tau_1(1+N)}\left(\xi_{3,N}\frac{2}{2+N} - \tau_1\right) \\ w_{1,2,1} &= \frac{1}{\xi_{3,N} - \tau_1(1+N)}\left(\frac{\xi_{3,N}^2}{2+N} - \big(\xi_{3,N} - \tau_2(1+N)\big)w_{1,1,2}\right) \end{aligned}$$

### C.13 Implicit method, $s = 2$, $j = 3$, $v = 4$

The conditions are :

Order of consistency 1 :

$$w_{N,3,1} + w_{N,3,2} = \xi_{3,N}^{N}$$

Order of consistency 2 :

$$w_{N,3,1}\tau_1 + w_{N,3,2}\tau_2 = \xi_{3,N}^{N+1}\frac{1}{1+N}$$

$$w_{N,3,1}(w_{1,1,1} + w_{1,1,2}) + w_{N,3,2}(w_{1,2,1} + w_{1,2,2}) = \xi_{3,N}^{N+1}\frac{1}{1+N}$$

Order of consistency 3 :

$$w_{N,3,1}\tau_1^2 + w_{N,3,2}\tau_2^2 = \xi_{3,N}^{N+2}2\frac{N!}{(2+N)!}$$

$$w_{N,3,1}\tau_1(w_{1,1,1} + w_{1,1,2}) + w_{N,3,2}\tau_2(w_{1,2,1} + w_{1,2,2}) = \xi_{3,N}^{N+2}2\frac{N!}{(2+N)!}$$

$$w_{N,3,1}(w_{1,1,1}\tau_1 + w_{1,1,2}\tau_2) + w_{N,3,2}(w_{1,2,1}\tau_1 + w_{1,2,2}\tau_2) = \xi_{3,N}^{N+2}\frac{N!}{(2+N)!}$$

$$w_{N,3,1}(w_{1,1,1}(w_{1,1,1} + w_{1,1,2}) + w_{1,1,2}(w_{1,2,1} + w_{1,2,2}))$$

$$+w_{N,3,2}(w_{1,2,1}(w_{1,1,1} + w_{1,1,2}) + w_{1,2,2}(w_{1,2,1} + w_{1,2,2})) = \xi_{3,N}^{N+2}\frac{N!}{(2+N)!}$$

$$w_{N,3,1}(w_{1,1,1} + w_{1,1,2})^2 + w_{N,3,2}(w_{1,2,1} + w_{1,2,2})^2 = \xi_{3,N}^{N+2}2\frac{N!}{(2+N)!}$$

$$n > 1,\ w_{N,3,1}(w_{2,1,1} + w_{2,1,2}) + w_{N,3,2}(w_{2,2,1} + w_{2,2,2}) = \xi_{3,N}^{N+2}2\frac{N!}{(2+N)!}$$

Order of consistency 4 :

$$w_{N,3,1}\tau_1^3 + w_{N,3,2}\tau_2^3 = \xi_{3,N}^{N+3}6\frac{N!}{(3+N)!}$$

$$w_{N,3,1}(w_{1,1,1}\tau_1^2 + w_{1,1,2}\tau_2^2) + w_{N,3,2}(w_{1,2,1}\tau_1^2 + w_{1,2,2}\tau_2^2) = \xi_{3,N}^{N+3}2\frac{N!}{(3+N)!}$$

$$w_{N,3,1}\tau_1(w_{1,1,1}\tau_1 + w_{1,1,2}\tau_2) + w_{N,3,2}\tau_2(w_{1,2,1}\tau_1 + w_{1,2,2}\tau_2) = \xi_{3,N}^{N+3}3\frac{N!}{(3+N)!}$$

$$w_{N,3,1}\tau_1^2(w_{1,1,1} + w_{1,1,2}) + w_{N,3,2}\tau_2^2(w_{1,2,1} + w_{1,2,2}) = \xi_{3,N}^{N+3}6\frac{N!}{(3+N)!}$$

$$w_{N,3,1}(w_{1,1,1}(w_{1,1,1}\tau_1 + w_{1,1,2}\tau_2) + w_{1,1,2}(w_{1,2,1}\tau_1 + w_{1,2,2}\tau_2))$$

$$+w_{N,3,2}(w_{1,2,1}(w_{1,1,1}\tau_1 + w_{1,1,2}\tau_2) + w_{1,2,2}(w_{1,2,1}\tau_1 + w_{1,2,2}\tau_2)) = \xi_{3,N}^{N+3}\frac{N!}{(3+N)!}$$

$$w_{N,3,1}(w_{1,1,1}\tau_1(w_{1,1,1}+w_{1,1,2})+w_{1,1,2}\tau_2(w_{1,2,1}+w_{1,2,2})) \\ +w_{N,3,2}(w_{1,2,1}\tau_1(w_{1,1,1}+w_{1,1,2})+w_{1,2,2}\tau_2(w_{1,2,1}+w_{1,2,2}))=\xi_{3,N}^{N+3}2\frac{N!}{(3+N)!}$$

$$w_{N,3,1}\tau_1(w_{1,1,1}(w_{1,1,1}+w_{1,1,2})+w_{1,1,2}(w_{1,2,1}+w_{1,2,2})) \\ +w_{N,3,2}\tau_2(w_{1,2,1}(w_{1,1,1}+w_{1,1,2})+w_{1,2,2}(w_{1,2,1}+w_{1,2,2}))=\xi_{3,N}^{N+3}3\frac{N!}{(3+N)!}$$

$$w_{N,3,1}(w_{1,1,1}+w_{1,1,2})(w_{1,1,1}\tau_1+w_{1,1,2}\tau_2)+w_{N,3,2}(w_{1,2,1}+w_{1,2,2})(w_{1,2,1}\tau_1+w_{1,2,2}\tau_2)=\xi_{3,N}^{N+3}3\frac{N!}{(3+N)!}$$

$$w_{N,3,1}\tau_1(w_{1,1,1}+w_{1,1,2})^2+w_{N,3,2}\tau_2(w_{1,2,1}+w_{1,2,2})^2=\xi_{3,N}^{N+3}6\frac{N!}{(3+N)!}$$

$$w_{N,3,1}( \\ w_{1,1,1}(w_{1,1,1}(w_{1,1,1}+w_{1,1,2})+w_{1,1,2}(w_{1,2,1}+w_{1,2,2})) \\ +w_{1,1,2}(w_{1,2,1}(w_{1,1,1}+w_{1,1,2})+w_{1,2,2}(w_{1,2,1}+w_{1,2,2})) \\ ) \\ +w_{N,3,2}( \\ w_{1,2,1}(w_{1,1,1}(w_{1,1,1}+w_{1,1,2})+w_{1,1,2}(w_{1,2,1}+w_{1,2,2})) \\ +w_{1,2,2}(w_{1,2,1}(w_{1,1,1}+w_{1,1,2})+w_{1,2,2}(w_{1,2,1}+w_{1,2,2})) \\ )=\xi_{3,N}^{N+3}\frac{N!}{(3+N)!}$$

$$w_{N,3,1}\left(w_{1,1,1}(w_{1,1,1}+w_{1,1,2})^2+w_{1,1,2}(w_{1,2,1}+w_{1,2,2})^2\right) \\ +w_{N,3,2}\left(w_{1,2,1}(w_{1,1,1}+w_{1,1,2})^2+w_{1,2,2}(w_{1,2,1}+w_{1,2,2})^2\right)=\xi_{3,N}^{N+3}2\frac{N!}{(3+N)!}$$

$$w_{N,3,1}(w_{1,1,1}+w_{1,1,2})^3+w_{N,3,2}(w_{1,2,1}+w_{1,2,2})^3=\xi_{3,N}^{N+3}6\frac{N!}{(3+N)!}$$

$$w_{N,3,1}(w_{1,1,1}+w_{1,1,2})(w_{1,1,1}(w_{1,1,1}+w_{1,1,2})+w_{1,1,2}(w_{1,2,1}+w_{1,2,2})) \\ +w_{N,3,2}(w_{1,2,1}+w_{1,2,2})(w_{1,2,1}(w_{1,1,1}+w_{1,1,2})+w_{1,2,2}(w_{1,2,1}+w_{1,2,2}))=\xi_{3,N}^{N+3}3\frac{N!}{(3+N)!}$$

$$n>1,\ w_{N,3,1}(w_{2,1,1}\tau_1+w_{2,1,2}\tau_2)+w_{N,3,2}(w_{2,2,1}\tau_1+w_{2,2,2}\tau_2)=\xi_{3,N}^{N+3}2\frac{N!}{(3+N)!}$$

$$n>1,\ w_{N,3,1}(w_{1,1,1}+w_{1,1,2})(w_{2,1,1}+w_{2,1,2})+w_{N,3,2}(w_{1,2,1}+w_{1,2,2})(w_{2,2,1}+w_{2,2,2})=\xi_{3,N}^{N+3}6\frac{N!}{(3+N)!}$$

$$n>1,\ w_{N,3,1}(w_{2,1,1}(w_{1,1,1}+w_{1,1,2})+w_{2,1,2}(w_{1,2,1}+w_{1,2,2})) \\ +w_{N,3,2}(w_{2,2,1}(w_{1,1,1}+w_{1,1,2})+w_{2,2,2}(w_{1,2,1}+w_{1,2,2}))=\xi_{3,N}^{N+3}2\frac{N!}{(3+N)!}$$

$$n>1,\ w_{N,3,1}(w_{1,1,1}(w_{2,1,1}+w_{2,1,2})+w_{1,1,2}(w_{2,2,1}+w_{2,2,2})) \\ +w_{N,3,2}(w_{1,2,1}(w_{2,1,1}+w_{2,1,2})+w_{1,2,2}(w_{2,2,1}+w_{2,2,2}))=\xi_{3,N}^{N+3}2\frac{N!}{(3+N)!}$$

$$n > 1,\ w_{N,3,1}\tau_1\big(w_{2,1,1} + w_{2,1,2}\big) + w_{N,3,2}\tau_2\big(w_{2,2,1} + w_{2,2,2}\big) = \xi_{3,N}^{N+3} 6\frac{N!}{(3+N)!}$$

$$n > 2,\ w_{N,3,1}\big(w_{3,1,1} + w_{3,1,2}\big) + w_{N,3,2}\big(w_{3,2,1} + w_{3,2,2}\big) = \xi_{3,N}^{N+3} 6\frac{N!}{(3+N)!}$$

We define :

$$\tilde{\tau}_1 = w_{1,1,1} + w_{1,1,2}$$
$$\tilde{\tau}_2 = w_{1,2,1} + w_{1,2,2}$$

We rewrite the system as :

Order of consistency 1 :

$$w_{N,3,1} + w_{N,3,2} = \xi_{3,N}^{N}$$

Order of consistency 2 :

$$w_{N,3,1}\tau_1 + w_{N,3,2}\tau_2 = \xi_{3,N}^{N+1}\frac{1}{1+N}$$
$$w_{N,3,1}\tilde{\tau}_1 + w_{N,3,2}\tilde{\tau}_2 = \xi_{3,N}^{N+1}\frac{1}{1+N}$$

Order of consistency 3 :

$$w_{N,3,1}\tau_1^2 + w_{N,3,2}\tau_2^2 = \xi_{3,N}^{N+2} 2\frac{N!}{(2+N)!}$$
$$w_{N,3,1}\tau_1\tilde{\tau}_1 + w_{N,3,2}\tau_2\tilde{\tau}_2 = \xi_{3,N}^{N+2} 2\frac{N!}{(2+N)!}$$
$$w_{N,3,1}\big(w_{1,1,1}\tau_1 + w_{1,1,2}\tau_2\big) + w_{N,3,2}\big(w_{1,2,1}\tau_1 + w_{1,2,2}\tau_2\big) = \xi_{3,N}^{N+2}\frac{N!}{(2+N)!}$$
$$w_{N,3,1}\big(w_{1,1,1}\tilde{\tau}_1 + w_{1,1,2}\tilde{\tau}_2\big) + w_{N,3,2}\big(w_{1,2,1}\tilde{\tau}_1 + w_{1,2,2}\tilde{\tau}_2\big) = \xi_{3,N}^{N+2}\frac{N!}{(2+N)!}$$
$$w_{N,3,1}\tilde{\tau}_1^2 + w_{N,3,2}\tilde{\tau}_2^2 = \xi_{3,N}^{N+2} 2\frac{N!}{(2+N)!}$$
$$n > 1,\ w_{N,3,1}\big(w_{2,1,1} + w_{2,1,2}\big) + w_{N,3,2}\big(w_{2,2,1} + w_{2,2,2}\big) = \xi_{3,N}^{N+2} 2\frac{N!}{(2+N)!}$$

Order of consistency 4 :

$$w_{N,3,1}\tau_1^3 + w_{N,3,2}\tau_2^3 = \xi_{3,N}^{N+3} 6\frac{N!}{(3+N)!}$$
$$w_{N,3,1}\big(w_{1,1,1}\tau_1^2 + w_{1,1,2}\tau_2^2\big) + w_{N,3,2}\big(w_{1,2,1}\tau_1^2 + w_{1,2,2}\tau_2^2\big) = \xi_{3,N}^{N+3} 2\frac{N!}{(3+N)!}$$
$$w_{N,3,1}\tau_1\big(w_{1,1,1}\tau_1 + w_{1,1,2}\tau_2\big) + w_{N,3,2}\tau_2\big(w_{1,2,1}\tau_1 + w_{1,2,2}\tau_2\big) = \xi_{3,N}^{N+3} 3\frac{N!}{(3+N)!}$$
$$w_{N,3,1}\tau_1^2\tilde{\tau}_1 + w_{N,3,2}\tau_2^2\tilde{\tau}_2 = \xi_{3,N}^{N+3} 6\frac{N!}{(3+N)!}$$

$$w_{N,3,1}\big(w_{1,1,1}\big(w_{1,1,1}\tau_1+w_{1,1,2}\tau_2\big)+w_{1,1,2}\big(w_{1,2,1}\tau_1+w_{1,2,2}\tau_2\big)\big)$$

$$+w_{N,3,2}\big(w_{1,2,1}\big(w_{1,1,1}\tau_1+w_{1,1,2}\tau_2\big)+w_{1,2,2}\big(w_{1,2,1}\tau_1+w_{1,2,2}\tau_2\big)\big)=\xi_{3,N}^{N+3}\frac{N!}{(3+N)!}$$

$$w_{N,3,1}\big(w_{1,1,1}\tau_1\tilde{\tau}_1+w_{1,1,2}\tau_2\tilde{\tau}_2\big)+w_{N,3,2}\big(w_{1,2,1}\tau_1\tilde{\tau}_1+w_{1,2,2}\tau_2\tilde{\tau}_2\big)=\xi_{3,N}^{N+3}2\frac{N!}{(3+N)!}$$

$$w_{N,3,1}\tau_1\big(w_{1,1,1}\tilde{\tau}_1+w_{1,1,2}\tilde{\tau}_2\big)+w_{N,3,2}\tau_2\big(w_{1,2,1}\tilde{\tau}_1+w_{1,2,2}\tilde{\tau}_2\big)=\xi_{3,N}^{N+3}3\frac{N!}{(3+N)!}$$

$$w_{N,3,1}\tilde{\tau}_1\big(w_{1,1,1}\tau_1+w_{1,1,2}\tau_2\big)+w_{N,3,2}\tilde{\tau}_2\big(w_{1,2,1}\tau_1+w_{1,2,2}\tau_2\big)=\xi_{3,N}^{N+3}3\frac{N!}{(3+N)!}$$

$$w_{N,3,1}\tau_1\tilde{\tau}_1^2+w_{N,3,2}\tau_2\tilde{\tau}_2^2=\xi_{3,N}^{N+3}6\frac{N!}{(3+N)!}$$

$$w_{N,3,1}\big(w_{1,1,1}\big(w_{1,1,1}\tilde{\tau}_1+w_{1,1,2}\tilde{\tau}_2\big)+w_{1,1,2}\big(w_{1,2,1}\tilde{\tau}_1+w_{1,2,2}\tilde{\tau}_2\big)\big)$$

$$+w_{N,3,2}\big(w_{1,2,1}\big(w_{1,1,1}\tilde{\tau}_1+w_{1,1,2}\tilde{\tau}_2\big)+w_{1,2,2}\big(w_{1,2,1}\tilde{\tau}_1+w_{1,2,2}\tilde{\tau}_2\big)\big)=\xi_{3,N}^{N+3}\frac{N!}{(3+N)!}$$

$$w_{N,3,1}\big(w_{1,1,1}\tilde{\tau}_1^2+w_{1,1,2}\tilde{\tau}_2^2\big)+w_{N,3,2}\big(w_{1,2,1}\tilde{\tau}_1^2+w_{1,2,2}\tilde{\tau}_2^2\big)=\xi_{3,N}^{N+3}2\frac{N!}{(3+N)!}$$

$$w_{N,3,1}\tilde{\tau}_1^3+w_{N,3,2}\tilde{\tau}_2^3=\xi_{3,N}^{N+3}6\frac{N!}{(3+N)!}$$

$$w_{N,3,1}\tilde{\tau}_1\big(w_{1,1,1}\tilde{\tau}_1+w_{1,1,2}\tilde{\tau}_2\big)+w_{N,3,2}\tilde{\tau}_2\big(w_{1,2,1}\tilde{\tau}_1+w_{1,2,2}\tilde{\tau}_2\big)=\xi_{3,N}^{N+3}3\frac{N!}{(3+N)!}$$

$$n>1,\ w_{N,3,1}\big(w_{2,1,1}\tau_1+w_{2,1,2}\tau_2\big)+w_{N,3,2}\big(w_{2,2,1}\tau_1+w_{2,2,2}\tau_2\big)=\xi_{3,N}^{N+3}2\frac{N!}{(3+N)!}$$

$$n>1,\ w_{N,3,1}\tilde{\tau}_1\big(w_{2,1,1}+w_{2,1,2}\big)+w_{N,3,2}\tilde{\tau}_2\big(w_{2,2,1}+w_{2,2,2}\big)=\xi_{3,N}^{N+3}6\frac{N!}{(3+N)!}$$

$$n>1,\ w_{N,3,1}\big(w_{2,1,1}\tilde{\tau}_1+w_{2,1,2}\tilde{\tau}_2\big)$$

$$+w_{N,3,2}\big(w_{2,2,1}\tilde{\tau}_1+w_{2,2,2}\tilde{\tau}_2\big)=\xi_{3,N}^{N+3}2\frac{N!}{(3+N)!}$$

$$n>1,\ w_{N,3,1}\big(w_{1,1,1}\big(w_{2,1,1}+w_{2,1,2}\big)+w_{1,1,2}\big(w_{2,2,1}+w_{2,2,2}\big)\big)$$

$$+w_{N,3,2}\big(w_{1,2,1}\big(w_{2,1,1}+w_{2,1,2}\big)+w_{1,2,2}\big(w_{2,2,1}+w_{2,2,2}\big)\big)=\xi_{3,N}^{N+3}2\frac{N!}{(3+N)!}$$

$$n>1,\ w_{N,3,1}\tau_1\big(w_{2,1,1}+w_{2,1,2}\big)+w_{N,3,2}\tau_2\big(w_{2,2,1}+w_{2,2,2}\big)=\xi_{3,N}^{N+3}6\frac{N!}{(3+N)!}$$

$$n>2,\ w_{N,3,1}\big(w_{3,1,1}+w_{3,1,2}\big)+w_{N,3,2}\big(w_{3,2,1}+w_{3,2,2}\big)=\xi_{3,N}^{N+3}6\frac{N!}{(3+N)!}$$

We define the relations :

$$a_{1,1} : w_{1,1,1} = \tilde{\tau}_1 - w_{1,1,2}$$
$$a_{1,2} : w_{1,2,2} = \tilde{\tau}_2 - w_{1,2,1}$$
$$a_{2,1} : w_{2,1,1} = \tilde{\tau}_1' - w_{2,1,2}$$
$$a_{2,2} : w_{2,2,2} = \tilde{\tau}_2' - w_{2,2,1}$$

We now list the main operations performed to rewrite this system.

1) $a_{1,1}, a_{1,2}$ in 4.10 twice
2) 4.11 in 4.10
3) $a_{1,1}, a_{1,2}$ in 4.5 twice
4) 4.6 in 4.5,4.6
5) $a_{1,1}, a_{1,2}$ in 4.2
6) 4.4 in 4.2
7) $a_{1,1}, a_{1,2}$ in 4.3
8) 4.4 in 4.3
9) $a_{1,1}, a_{1,2}$ in 4.6
10) 4.9 in 4.6
11) $a_{1,1}, a_{1,2}$ in 4.7
12) 4.9 in 4.7
13) $a_{1,1}, a_{1,2}$ in 4.8
14) 4.9 in 4.8
15) $a_{1,1}, a_{1,2}$ in 4.11
16) 4.12 in 4.11
17) $a_{1,1}, a_{1,2}$ in 4.13
18) 4.12 in 4.13
19) $a_{1,1}, a_{1,2}$ in 3.3
20) 3.2 in 3.3
21) $a_{1,1}, a_{1,2}$ in 3.4
22) 3.5 in 3.4
23) 3.4 in 4.10, 4.13, 4.7, 4.11
24) 3.3 in 4.3, 4.5, 4.2, 4.8
25) Use identity $2(\tau_1\tilde{\tau}_1 - \tau_2\tilde{\tau}_2) = (\tau_1 - \tau_2)(\tilde{\tau}_1 + \tilde{\tau}_2) + (\tilde{\tau}_1 - \tilde{\tau}_2)(\tau_1 + \tau_2)$ in 4.6
26) 3.3 in 4.6
27) 3.4 in 4.6
28) 4.2 in 4.6
29) 4.11 in 4.6
30) 3.1 in 4.1
31) isolate in $w_{N,3,1}\tau_1^2$ 3.1
32) 3.1 in 4.1
33) isolate $w_{N,3,2}\tau_2\tilde{\tau}_2$ in 3.2
34) 3.2 in 4.4
35) isolate $w_{N,3,1}\tilde{\tau}_1^2$ in 3.5
36) 3.5 in 4.9,4.12
37) isolate $w_{N,3,2}\tau_2$ in 2.1
38) 2.1 in 3.1
39) isolate $w_{N,3,2}\tilde{\tau}_2$ in 2.2

40) 2.2 in 3.2
41) 3.5, 4.13 in 4.10
42) 4.8 in 4.13
43) 4.7 in 4.8
44) 4.7 in 4.5
45) 4.3 in 4.7
46) 4.8 in 4.10
47) 4.8 in 4.11
48) 4.2 in 4.11
49) 4.8 in 4.12
50) 4.7 in 4.13
51) 4.13 in 4.7
52) 4.9 in 4.12
53) 4.4 in 4.9
54) 4.1 in 4.4
55) 4.4 in 4.9
56) 4.9 in 4.4
57) 4.12 in 3.5
58) 4.8 in 3.5
59) 3.2 in 3.5
60) 3.1 in 4.1
61) 3.1 in 3.2
62) 3.2 in 4.9
63) isolate $w_{N,3,2}$ in 1.1
64) 1.1 in 2.1
65) 2.2, 4.8 in 2.2
66) 2.1 in 2.2
67) 2.2 in 4.12
68) 2.2 in 4.13
69) 2.1 in 3.1
70) 3.1 in 4.1
71) 4.2 in 4.1
72) 1.1 in 3.3, 3.4
73) 2.1 in 3.2
74) 4.11 in 3.2
75) 2.1 in 4.3
76) 2.2 in 3.4
77) 3.3 in 3.4
78) 2.1 in 3.3
79) $a_{1,2}$ in 4.5
80) 4.8 in 4.5
81) 4.5 in 3.3 twice
82) 3.3 in 4.3
83) 3.1 in 4.3
84) 4.8 in 3.4
85) 4.11 in 3.4

86) replace $\tilde{\tau}_1', \tilde{\tau}_2'$ in second rank equations
87) $a_{2,1}, a_{2,2}$ in 4.14, 4.18 in 4.14
88) $a_{2,1}, a_{2,2}$ in 4.16, 4.15 in 4.16
89) $a_{1,1}, a_{1,2}$ in 4.17, 4.15 in 4.17
90) isolate $w_{N,3,2}\tilde{\tau}_2'$ in 3.6
91) 3.6 in 4.15, 4.18
92) 2.2 in 4.15
93) 4.8 in 4.15
94) 4.18 in 4.15
95) 1.1 in 3.6
96) 1.1 in 4.14, 4.16,4.17, 4.19
97) 2.2 in 4.16
98) 4.8 in 4.16
99) 4.11 in 4.16
100) 4.14 in 4.16
101) 2.1 in 4.18
102) 3.6 in 4.17
103) 4.5 in 4.17
104) 2.1 in 4.14
105) 4.2 in 3.1

Which yields :

Order of consistency 1 :

$$w_{N,3,2} = \xi_{3,N}^N - w_{N,3,1}$$

Order of consistency 2 :

$$w_{N,3,1}(\tau_1 - \tau_2) = \xi_{3,N}^N \left( \frac{\xi_{3,N}}{1+N} - \tau_2 \right)$$

$$w_{N,3,1}[(\tilde{\tau}_1 - \tilde{\tau}_2) - (\tau_1 - \tau_2)] = 0$$

Order of consistency 3 :

$$0 = \xi_{3,N} \left( \tau_2 - 2\frac{\xi_{3,N}}{3+N} \left( 1 + \sqrt{\frac{1+N}{2(2+N)}} \right) \right) \left( \tau_2 - 2\frac{\xi_{3,N}}{3+N} \left( 1 - \sqrt{\frac{1+N}{2(2+N)}} \right) \right)$$

$$0 = 0$$

$$\xi_{3,N} \left[ w_{1,1,2}(\tau_1 - \tau_2) + \frac{\xi_{3,N}}{3+N} \left( \frac{\xi_{3,N}}{(2+N)} - 2\tau_1 \right) \right] = 0$$

$$0 = 0$$

$$0 = 0$$

$$n > 1,\ \ w_{N,3,1}(\tilde{\tau}_1' - \tilde{\tau}_2') = \xi_{3,N}^N \left( \xi_{3,N}^2 2\frac{N!}{(2+N)!} - \tilde{\tau}_2' \right)$$

Order of consistency 4 :

$$0 = 0$$

$$\xi_{3,N}\left(\tau_1+\tau_2-\xi_{3,N}\frac{4}{3+N}\right)=0$$
$$0=0$$
$$0=0$$
$$\xi_{3,N}\left(w_{1,1,2}+w_{1,2,1}-\xi_{3,N}\frac{2}{3+N}\right)=0$$
$$0=0$$
$$0=0$$
$$\xi_{3,N}(\tau_2-\tilde{\tau}_2)=0$$
$$0=0$$
$$0=0$$
$$\xi_{3,N}(\tilde{\tau}_1-\tau_1)=0$$
$$0=0$$
$$0=0$$
$$n>1,\ \xi_{3,N}\left(\frac{\xi_{3,N}}{1+N}-\tau_2\right)w_{2,1,2}+\xi_{3,N}w_{2,2,1}\left(\frac{\xi_{3,N}}{1+N}-\tau_1\right)=\xi_{3,N}^4 4\frac{N!}{(3+N)!}$$
$$0=0$$
$$0=0$$
$$n>1,\ \xi_{3,N}\left[w_{1,2,1}(\tilde{\tau}_1'-\tilde{\tau}_2')+\xi_{3,N}\frac{2}{3+N}\tilde{\tau}_2'\right]=0$$
$$n>1,\ \xi_{3,N}\left[\left(\frac{\xi_{3,N}}{1+N}-\tau_2\right)\tilde{\tau}_1'-\xi_{3,N}^2 2\frac{N!}{(2+N)!}\left(\xi_{3,N}\frac{3}{3+N}-\tau_2\right)\right]=0$$
$$n>2,\ w_{N,3,1}\left[w_{3,1,1}+w_{3,1,2}-\left(w_{3,2,1}+w_{3,2,2}\right)\right]=\xi_{3,N}^N\left(\xi_{3,N}^3 6\frac{N!}{(3+N)!}-\left(w_{3,2,1}+w_{3,2,2}\right)\right)$$

The solutions are :

**0-0**

$\xi_{3,N}=0, w_{N,3,1}=0$

$$w_{N,3,2}=0$$

**0-1**

$\xi_{3,N}=0, w_{N,3,1}\neq 0$

$$w_{N,3,2}=-w_{N,3,1}$$
$$\tau_1=\tau_2$$
$$\tilde{\tau}_1=\tilde{\tau}_2$$
$$n>1,\ \tilde{\tau}_1'=\tilde{\tau}_2'$$
$$n>2,\ w_{3,1,1}=w_{3,2,1}+w_{3,2,2}-w_{3,1,2}$$

**1-0**

$\xi_{3,N} \neq 0, \tau_2 = 2\frac{\xi_{3,N}}{3+N}\left(1+\sqrt{\frac{1+N}{2(2+N)}}\right)$

3.1 implies that $\tau_2 \neq \frac{\xi_{3,N}}{1+N}$, which implies thanks to 2.1 that $w_{N,3,1} \neq 0$ and $\tau_2 \neq \tau_1$.

$$\tau_1 = 2\frac{\xi_{3,N}}{3+N}\left(1-\sqrt{\frac{1+N}{2(2+N)}}\right)$$

$$\tau_2 = 2\frac{\xi_{3,N}}{3+N}\left(1+\sqrt{\frac{1+N}{2(2+N)}}\right)$$

$$w_{1,1,1} = \frac{\xi_{3,N}}{3+N}\left(1+\frac{1-N}{2(1+N)}\sqrt{\frac{1+N}{2(2+N)}}\right)$$

$$w_{1,1,2} = \frac{\xi_{3,N}}{3+N}\left(1-(5+3N)\frac{\sqrt{2(2+N)(1+N)}}{4(2+N)(1+N)}\right)$$

$$w_{1,2,1} = \frac{\xi_{3,N}}{3+N}\left(1+(5+3N)\frac{\sqrt{2(2+N)(1+N)}}{4(2+N)(1+N)}\right)$$

$$w_{1,2,2} = \frac{\xi_{3,N}}{3+N}\left(1+\frac{N-1}{2(1+N)}\sqrt{\frac{1+N}{2(2+N)}}\right)$$

$$w_{N,3,1} = \frac{\xi_{3,N}^N}{4}\left[2+(N-1)\sqrt{2\frac{2+N}{(1+N)^3}}\right]$$

$$w_{N,3,2} = \frac{\xi_{3,N}^N}{4}\left[2+(1-N)\sqrt{2\frac{2+N}{(1+N)^3}}\right]$$

$$n>1,\ w_{2,1,1} = \left(2\frac{\xi_{3,N}}{3+N}\left(1-\sqrt{\frac{1+N}{2(2+N)}}\right)\right)^2 - \frac{\xi_{3,N}^2 4\frac{\sqrt{2+N}}{2+N} - w_{2,2,1}\left[(1-N)\sqrt{2+N}+(1+N)\sqrt{2(1+N)}\right]}{(1-N)\sqrt{2+N}-(1+N)\sqrt{2(1+N)}}$$

$$n>1,\ w_{2,1,2} = \frac{\xi_{3,N}^2 4\frac{\sqrt{2+N}}{2+N} - w_{2,2,1}\left[(1-N)\sqrt{2+N}+(1+N)\sqrt{2(1+N)}\right]}{(1-N)\sqrt{2+N}-(1+N)\sqrt{2(1+N)}}$$

$$n>1,\ w_{2,2,2} = \left(2\frac{\xi_{3,N}}{3+N}\left(1+\sqrt{\frac{1+N}{2(2+N)}}\right)\right)^2 - w_{2,2,1}$$

$$n>2,\ w_{3,1,1} = w_{3,2,1}+w_{3,2,2}-w_{3,1,2}+4\frac{\sqrt{(1+N)^3}}{(N-1)\sqrt{2(2+N)}+2\sqrt{(1+N)^3}}\left(\xi_{3,N}^3 6\frac{N!}{(3+N)!}-\left(w_{3,2,1}+w_{3,2,2}\right)\right)$$

**1-1**

$\xi_{3,N} \neq 0, \tau_2 \neq 2\frac{\xi_{3,N}}{3+N}\left(1+\sqrt{\frac{1+N}{2(2+N)}}\right)$

3.1 implies that $\tau_2 \neq \frac{\xi_{3,N}}{1+N}$, which implies thanks to 2.1 that $w_{N,3,1} \neq 0$ and $\tau_2 \neq \tau_1$

$$\tau_1 = 2\frac{\xi_{3,N}}{3+N}\left(1+\sqrt{\frac{1+N}{2(2+N)}}\right)$$

$$\tau_2 = 2\frac{\xi_{3,N}}{3+N}\left(1-\sqrt{\frac{1+N}{2(2+N)}}\right)$$

$$w_{N,3,1} = \frac{\xi_{3,N}^N}{4}\left(2+(1-N)\sqrt{2\frac{2+N}{(1+N)^3}}\right)$$

$$w_{N,3,2} = \frac{\xi_{3,N}^N}{4}\left(2+(N-1)\sqrt{2\frac{2+N}{(1+N)^3}}\right)$$

$$w_{1,1,1} = \frac{\xi_{3,N}}{3+N}\left(1+\frac{N-1}{2(1+N)}\sqrt{\frac{1+N}{2(2+N)}}\right)$$

$$w_{1,1,2} = \frac{\xi_{3,N}}{3+N}\left(1+\frac{5+3N}{4(2+N)(1+N)}\sqrt{2(2+N)(1+N)}\right)$$

$$w_{1,2,1} = \frac{\xi_{3,N}}{3+N}\left(1-\frac{5+3N}{4(2+N)(1+N)}\sqrt{2(2+N)(1+N)}\right)$$

$$w_{1,2,2} = \frac{\xi_{3,N}}{3+N}\left(1+\frac{1-N}{2(1+N)}\sqrt{\frac{1+N}{2(2+N)}}\right)$$

$$n>1,\ w_{2,1,1} = \left(2\frac{\xi_{3,N}}{3+N}\left(1+\sqrt{\frac{1+N}{2(2+N)}}\right)\right)^2 - \frac{\xi_{3,N}^2 4\frac{\sqrt{2+N}}{(2+N)} - w_{2,2,1}\left(\sqrt{2+N}(1-N)-(1+N)\sqrt{2(1+N)}\right)}{\sqrt{2+N}(1-N)+(1+N)\sqrt{2(1+N)}}$$

$$n>1,\ w_{2,1,2} = \frac{\xi_{3,N}^2 4\frac{\sqrt{2+N}}{(2+N)} - w_{2,2,1}\left(\sqrt{2+N}(1-N)-(1+N)\sqrt{2(1+N)}\right)}{\sqrt{2+N}(1-N)+(1+N)\sqrt{2(1+N)}}$$

$$n>1,\ w_{2,2,2} = \left(2\frac{\xi_{3,N}}{3+N}\left(1-\sqrt{\frac{1+N}{2(2+N)}}\right)\right)^2 - w_{2,2,1}$$

$$n>2,\ w_{3,1,1} = w_{3,2,1}+w_{3,2,2}-w_{3,1,2}+4\frac{\sqrt{1+N}^3}{2\sqrt{1+N}^3+(1-N)\sqrt{2(2+N)}}\left(\xi_{3,N}^3 6\frac{N!}{(3+N)!}-(w_{3,2,1}+w_{3,2,2})\right)$$

# D Proof of generalized solved system conditions

We here prove the theorem of Section III.1.39.

**D.1 Theorem : Generalized solved system conditions**

Let $(n, s, v) \in \mathbb{N}^* \times \mathbb{N}^* \times \mathbb{N}^*$, $M = (\tau, w, \xi) \in \mathrm{MORK}_{n,s}$, $v' \in \mathbb{N}^*$, and $(j, N) \in [\![1, s+1]\!] \times [\![1, n]\!]$. If $M$ is at least of order of consistency $v'$ for length $n$ at stage $j$ and rank $N$, then :

$$\forall k \in [\![0, v'-1]\!],\ \forall \lambda \in \mathbb{N}^{[\![1,n]\!]}$$

$$|\lambda| \leq v' - 1 - k \Rightarrow \sum_{j'=1}^{s} w_{N,j,j'} \tau_{j'}^{k} \left( \prod_{N=1}^{n} \xi_{j',N}^{\lambda_N} \right) = \xi_{j,N}^{N+k+\,|\lambda|} \frac{N!(k+|\lambda|)!}{(k+|\lambda|+N)!}$$

The proof is more involved than the solved system conditions and require a somewhat complex formula [17], which uses multi-index notations.

**D.2 Theorem : Multivariate Faa di Bruno formula**

Let $(d, m, k) \in \mathbb{N}^* \times \mathbb{N}^* \times \mathbb{N}^*$. Let $X \in \mathcal{P}(\mathbb{R})$ a non empty set, and $x \in X$ such that there exists a neighborhood $U$ of $x$ contained in $X$. Let $Y \in \mathcal{P}(\mathbb{R}^m)$ a non empty set and $y \in Y$ such that there exists a neighborhood $V$ of $y$ contained in $Y$. Let $g \in \{X \to Y\}$ such that $g(x) = y$ and $g$ is of class $C^k$ on $U$. Let $f \in \{Y \to \mathbb{R}^d\}$ such that $f$ is of class $C^k$ on $V$.
$f \circ g$ is $k$ times differentiable at $x$ and :

$$\frac{\mathrm{d}^k f \circ g}{\mathrm{d}t^k}(x) = \sum_{\substack{\lambda \in \mathbb{N}^{[\![1,m]\!]} \\ 1 \leq\ |\lambda|\ \leq k}} \frac{\partial^{|\lambda|} f}{\prod_{k'=1}^{m} \partial g_{k'}^{\lambda_{k'}}}(y) \sum_{\substack{\alpha \in \mathbb{N}^{[\![1,k]\!] \times [\![1,m]\!]} \\ \sum_{k'=1}^{k} \alpha_{k'} = \lambda \\ \sum_{k'=1}^{k} k'\ |\alpha_{k'}|\ = k}} k! \prod_{k'=1}^{k} \frac{1}{(\alpha_{k'}!)(k'!)^{|\alpha_{k'}|}} \left( \frac{\mathrm{d}^{k'} g}{\mathrm{d}t^{k'}}(x) \right)^{\alpha_{k'}}$$

To prove this theorem we use the equation for order of consistency, which we can do because all multi-order Runge-Kutta methods satisfy, for all $(j', N') \in [\![1, s]\!] \times [\![1, n]\!]$, $\tilde{w}_{N',N',j'} = 1$ and $M$ is at least of order of consistency 0 at stage $j$ and rank $N$.

We use Faa di Bruno's formula to find the derivatives of $f \circ \overline{\mathcal{J}}^{n-1}\hat{\mathrm{y}}$, then we find the differentials of $f$. The next steps are simply a process of identification of the terms. In order, we identify the partial derivatives of $f$ in respect to $t$, the partial derivatives of $f$ in respect to $\hat{\mathrm{y}}$ and its derivatives, we will remove terms which involve partial derivatives of $f$, then we will identify by the powers of $h$.

In this proof we omit the arguments of the functions since they always are the initial values. Faa di Bruno's formula yields :

$$\frac{\mathrm{d}^{k''} f \circ \overline{\mathcal{J}}^{n-1}\hat{\mathrm{y}}}{\mathrm{d}t^{k''}} = \sum_{\substack{\lambda \in \mathbb{N}^{[\![0,nd]\!]} \\ 1 \leq |\lambda| \leq k''}} \frac{\partial^{|\lambda|} f}{\partial t^{\lambda_0} \prod_{N=1}^{n} (\partial y_N)^{\lambda_N}} \sum_{\substack{\alpha \in \mathbb{N}^{[\![1,k'']\!] \times [\![0,nd]\!]} \\ \sum \alpha_l = \lambda \\ \sum |\alpha_l|\ l = k''}} k''! \prod_{l=1}^{k''} \frac{\left( \frac{\mathrm{d}^l \overline{\mathcal{J}}^{n-1}\hat{\mathrm{y}}}{\mathrm{d}t^l} \right)^{\alpha_l}}{\alpha_l! l!^{|\alpha_l|}}$$

$$= \sum_{\substack{\lambda \in \mathbb{N}^{\llbracket 0, nd \rrbracket} \\ 1 \leq |\lambda| \leq k''}} \frac{\partial^{|\lambda|} f}{\partial t^{\lambda_0} \prod_{N=1}^{n} (\partial y_N)^{\lambda_N}} \sum_{\substack{\alpha' \in \mathbb{N}^{\llbracket 1, k'' \rrbracket} \\ \alpha \in \mathbb{N}^{\llbracket 1, k'' \rrbracket \times \llbracket 1, n \rrbracket \times \llbracket 1, d \rrbracket} \\ \sum [\alpha'_l, \alpha_{l,1}, .., \alpha_{l,n}] = [\lambda_0 \ \lambda_1 \ \dots \ \lambda_{nd}] \\ \sum (\alpha'_l + |\alpha_l|) l = k''}} \Bigg[$$

$$k''! \prod_{l=1}^{k''} \left( \frac{\mathrm{d}^l t}{\mathrm{d} t^l} \right)^{\alpha'_l} \frac{1}{\alpha'_l ! \alpha_l ! l!^{\alpha'_l + |\alpha_l|}} \prod_{N=1}^{n} \left( \frac{\mathrm{d}^{l+n-N} \hat{\mathrm{y}}}{\mathrm{d} t^{l+n-N}} \right)^{\alpha_{l,N}} \Bigg]$$

$$= \sum_{\substack{k \in \mathbb{N} \\ \lambda \in \mathbb{N}^{\llbracket 1, n \rrbracket \times \llbracket 1, d \rrbracket} \\ 1 \leq k + |\lambda| \leq k''}} \frac{\partial^{k+|\lambda|} f}{\partial t^{k} \prod_{N=1}^{n} (\partial y_N)^{\lambda_N}} \sum_{\substack{\alpha' \in \mathbb{N}^{\llbracket 1, k'' \rrbracket} \\ \alpha \in \mathbb{N}^{\llbracket 1, k'' \rrbracket \times \llbracket 1, n \rrbracket \times \llbracket 1, d \rrbracket} \\ \sum \alpha'_l = k \\ \sum \alpha_l = \lambda \\ \sum (\alpha'_l + |\alpha_l|) l = k''}} \Bigg[$$

$$\frac{k''!}{\alpha'!} \prod_{l=1}^{k''} \left( \frac{\mathrm{d}^l t}{\mathrm{d} t^l} \right)^{\alpha'_l} \frac{1}{\alpha_l ! l!^{\alpha'_l + |\alpha_l|}} \prod_{N=1}^{n} \left( \frac{\mathrm{d}^{l+n-N} \hat{\mathrm{y}}}{\mathrm{d} t^{l+n-N}} \right)^{\alpha_{l,N}} \Bigg]$$

$$= \sum_{\substack{k \in \mathbb{N} \\ \lambda \in \mathbb{N}^{\llbracket 1, n \rrbracket \times \llbracket 1, d \rrbracket} \\ 1 \leq k + |\lambda| \leq k''}} \frac{\partial^{k+|\lambda|} f}{\partial t^{k} \prod_{N=1}^{n} (\partial y_N)^{\lambda_N}} \sum_{\substack{\alpha' \in \mathbb{N} \\ \alpha \in \mathbb{N}^{\llbracket 1, k'' \rrbracket \times \llbracket 1, n \rrbracket \times \llbracket 1, d \rrbracket} \\ \alpha' = k \\ \sum \alpha_l = \lambda \\ \alpha' + \sum |\alpha_l| \ l = k''}} \frac{k''!}{\alpha'!} \prod_{l=1}^{k''} \frac{1}{\alpha_l ! l!^{|\alpha_l|}} \prod_{N=1}^{n} \left( \frac{\mathrm{d}^{l+n-N} \hat{\mathrm{y}}}{\mathrm{d} t^{l+n-N}} \right)^{\alpha_{l,N}}$$

$$= \sum_{k=0}^{k''} \sum_{\substack{\lambda \in \mathbb{N}^{\llbracket 1, n \rrbracket \times \llbracket 1, d \rrbracket} \\ 1 \leq k + |\lambda| \leq k''}} \frac{\partial^{k+|\lambda|} f}{\partial t^{k} \prod_{N=1}^{n} (\partial y_N)^{\lambda_N}} \sum_{\substack{\alpha \in \mathbb{N}^{\llbracket 1, k'' \rrbracket \times \llbracket 1, n \rrbracket \times \llbracket 1, d \rrbracket} \\ \sum \alpha_l = \lambda \\ k + \sum |\alpha_l| \ l = k''}} \frac{k''!}{k!} \prod_{l=1}^{k''} \frac{1}{\alpha_l ! l!^{|\alpha_l|}} \prod_{N=1}^{n} \left( \frac{\mathrm{d}^{l+n-N} \hat{\mathrm{y}}}{\mathrm{d} t^{l+n-N}} \right)^{\alpha_{l,N}}$$

The expressions of the differentials of $f$ are easier to find.

$$\mathrm{d}^{k'} f = \sum_{k=0}^{k'} \sum_{\substack{\lambda \in \mathbb{N}^{\llbracket 1, n \rrbracket \times \llbracket 1, d \rrbracket} \\ |\lambda| = k' - k}} \frac{k'!}{k! \lambda!} \frac{\partial^{k'} f}{\partial t^{k} \prod_{N=1}^{n} (\partial y_N)^{\lambda_N}} dt^k \prod_{N=1}^{n} (\mathrm{d} y_N)^{\lambda_N}$$

The equation for order of consistency is :

$$\sum_{j'=1}^{s} w_{N,j,j'} \sum_{k'=0}^{v-1} \frac{1}{k'!} \mathrm{d}^{k'} f - \sum_{k''=0}^{v-1} \xi_{j,N}^{N+k''} h^{k''} \frac{N!}{(k''+N)!} \frac{\mathrm{d}^{k''} f \circ \overline{\mathcal{J}}^{n-1} \hat{\mathrm{y}}}{\mathrm{d} t^{k''}} = \underset{h \to 0}{\mathcal{O}} (h^v)$$

Let $y_n$ the approximations of $M$. We replace these expressions in the equation for order of consistency.

$$
\begin{aligned}
&\sum_{j'=1}^{s} w_{N,j,j'} f - \sum_{k''=0} \xi_{j,N}^{N} f \\
&\quad + \sum_{j'=1}^{s} w_{N,j,j'} \sum_{k'=1}^{v-1} \frac{1}{k'!} \sum_{k=0}^{k'} \sum_{\substack{\lambda\in\mathbb{N}^{\llbracket 1,n\rrbracket\times\llbracket 1,d\rrbracket}\\ |\lambda| =k'-k}} \frac{k'!}{k!\lambda!} \frac{\partial^{k'} f}{\partial t^k \prod_{N=1}^{n} (\partial y_N)^{\lambda_N}} (\tau_{j'} h)^k \prod_{N=1}^{n} \left(y_{n,j',N} - y_{n,0,N}\right)^{\lambda_N} \\
&\quad - \sum_{k''=1}^{v-1} \xi_{j,N}^{N+k''} h^{k''} \frac{N!}{(k''+N)!} \sum_{k=0}^{k''} \sum_{\substack{\lambda\in\mathbb{N}^{\llbracket 1,n\rrbracket\times\llbracket 1,d\rrbracket}\\ 1\leq k+ |\lambda|\leq k''}} \frac{\partial^{k+ |\lambda|} f}{\partial t^k \prod_{N=1}^{n} (\partial y_N)^{\lambda_N}} \sum_{\substack{\alpha\in\mathbb{N}^{\llbracket 1,k''\rrbracket\times\llbracket 1,n\rrbracket\times\llbracket 1,d\rrbracket}\\ \sum\alpha_l=\lambda\\ k+\sum|\alpha_l|\, l=k''}} \Bigg[ \\
&\quad \frac{k''!}{k!} \prod_{l=1}^{k''} \frac{1}{\alpha_l!\, l!^{|\alpha_l|}} \prod_{N=1}^{n} \left(\frac{\mathrm{d}^{l+n-N}\hat{\mathrm{y}}}{\mathrm{d}t^{l+n-N}}\right)^{\alpha_{l,N}} \Bigg] = \underset{h\to 0}{\mathcal{O}}(h^v)
\end{aligned}
$$

The first terms are the only one with a derivative of order 0, they hence must be equal. The sum $\sum_{k''=1}^{v-1} \sum_{k=0}^{k''}$ is equivalent to $\sum_{k=0}^{v-1} \sum_{k''=\max(1,k)}^{v-1}$, therefore the system is equivalent to :

$$
\sum_{j'=1}^{s} w_{N,j,j'} = \xi_{j,N}^{N},
$$

$$
\begin{aligned}
&\sum_{k=0}^{v-1} \sum_{k'=\max(1,k)}^{v-1} \sum_{j'=1}^{s} w_{N,j,j'} \sum_{\substack{\lambda\in\mathbb{N}^{\llbracket 1,n\rrbracket\times\llbracket 1,d\rrbracket}\\ |\lambda| =k'-k}} \frac{1}{k!\lambda!} \frac{\partial^{k'} f}{\partial t^k \prod_{N=1}^{n} (\partial y_N)^{\lambda_N}} (\tau_{j'} h)^k \prod_{N=1}^{n} \left(y_{n,j',N} - y_{n,0,N}\right)^{\lambda_N} \\
&\quad -\sum_{k=0}^{v-1} \sum_{k''=\max(1,k)}^{v-1} \xi_{j,N}^{N+k''} h^{k''} \frac{N!}{(k''+N)!} \sum_{k'=1}^{k''} \sum_{\substack{\lambda\in\mathbb{N}^{\llbracket 1,n\rrbracket\times\llbracket 1,d\rrbracket}\\ k+ |\lambda| =k'}} \frac{\partial^{k+ |\lambda|} f}{\partial t^k \prod_{N=1}^{n} (\partial y_N)^{\lambda_N}} \sum_{\substack{\alpha\in\mathbb{N}^{\llbracket 1,k''\rrbracket\times\llbracket 1,n\rrbracket\times\llbracket 1,d\rrbracket}\\ \sum\alpha_l=\lambda\\ k+\sum|\alpha_l|\, l=k''}} \Bigg[ \\
&\quad \frac{k''!}{k!} \prod_{l=1}^{k''} \frac{1}{\alpha_l!\, l!^{|\alpha_l|}} \prod_{N=1}^{n} \left(\frac{\mathrm{d}^{l+n-N}\hat{\mathrm{y}}}{\mathrm{d}t^{l+n-N}}\right)^{\alpha_{l,N}} \Bigg] = \underset{h\to 0}{\mathcal{O}}(h^v)
\end{aligned}
$$

We identify by the partial derivatives of $f$ in respect to $t$.

$$
\begin{aligned}
&\sum_{j'=1}^{s} w_{N,j,j'} = \xi_{j,N}^{N},\ \forall k \in \llbracket 0, v-1\rrbracket, \quad \sum_{k'=\max(1,k)}^{v-1} \sum_{\substack{\lambda\in\mathbb{N}^{\llbracket 1,n\rrbracket\times\llbracket 1,d\rrbracket}\\ |\lambda| +k=k'}} \frac{1}{k!\lambda!} \frac{\partial^{|\lambda| +k} f}{\partial t^k \prod_{N=1}^{n} (\partial y_N)^{\lambda_N}} \sum_{j'=1}^{s} \Bigg[ \\
&\quad w_{N,j,j'} (\tau_{j'} h)^k \prod_{N=1}^{n} \left(y_{n,j',N} - y_{n,0,N}\right)^{\lambda_N} \Bigg] \\
&\quad - \sum_{k''=\max(1,k)}^{v-1} \sum_{k'=\max(1,k)}^{k''} \sum_{\substack{\lambda\in\mathbb{N}^{\llbracket 1,n\rrbracket\times\llbracket 1,d\rrbracket}\\ k+ |\lambda| =k'}} \xi_{j,N}^{N+k''} h^{k''} \frac{N!}{(k''+N)!} \frac{\partial^{k+ |\lambda|} f}{\partial t^k \prod_{N=1}^{n} (\partial y_N)^{\lambda_N}} \sum_{\substack{\alpha\in\mathbb{N}^{\llbracket 1,k''\rrbracket\times\llbracket 1,n\rrbracket\times\llbracket 1,d\rrbracket}\\ \sum\alpha_l=\lambda\\ k+\sum|\alpha_l|\, l=k''}} \Bigg[ \\
&\quad \frac{k''!}{k!} \prod_{l=1}^{k''} \frac{1}{\alpha_l!\, l!^{|\alpha_l|}} \prod_{N=1}^{n} \left(\frac{\mathrm{d}^{l+n-N}\hat{\mathrm{y}}}{\mathrm{d}t^{l+n-N}}\right)^{\alpha_{l,N}} \Bigg] = \underset{h\to 0}{\mathcal{O}}(h^v)
\end{aligned}
$$

Let $a \in \mathbb{R}^{\mathbb{N}^*\times\mathbb{N}^*}$. We have :

$$\sum_{k''=\max(1,k)}^{v-1} \sum_{k'=\max(1,k)}^{k''} a_{k'',k'} = \sum_{k''=0}^{v-1-\max(1,k)} \sum_{k'=0}^{k''} a_{k''+\max(1,k),k'+\max(1,k)}$$

$$= \sum_{k'=0}^{v-1-\max(1,k)} \sum_{k''=k'}^{v-1-\max(1,k)} a_{k''+\max(1,k),k'+\max(1,k)} = \sum_{k'=\max(1,k)}^{v-1} \sum_{k''=k'}^{v-1} a_{k'',k'}$$

We deduce that the system is equivalent to :

$$\sum_{j'=1}^{s} w_{N,j,j'} = \xi_{j,N}^{N},$$

$$\forall k \in [\![0, v-1]\!], \quad \sum_{k'=\max(1,k)}^{v-1} \sum_{\substack{\lambda \in \mathbb{N}^{[\![1,n]\!]\times[\![1,d]\!]} \\ |\lambda| + k = k'}} \frac{1}{\lambda!} \frac{\partial^{|\lambda| + k} f}{\partial t^k \prod_{N=1}^{n} \left(\partial y_N\right)^{\lambda_N}} \sum_{j'=1}^{s} \Bigg[ w_{N,j,j'} \left(\tau_{j'} h\right)^k \prod_{N=1}^{n} \left(y_{n,j',N} - y_{n,0,N}\right)^{\lambda_N} \Bigg]$$

$$- \sum_{k'=\max(1,k)}^{v-1} \sum_{\substack{\lambda \in \mathbb{N}^{[\![1,n]\!]\times[\![1,d]\!]} \\ k + |\lambda| = k'}} \sum_{k''=k'}^{v-1} \xi_{j,N}^{N+k''} h^{k''} \frac{N!k''!}{(k''+N)!} \frac{\partial^{k+|\lambda|} f}{\partial t^k \prod_{N=1}^{n} \left(\partial y_N\right)^{\lambda_N}} \sum_{\substack{\alpha \in \mathbb{N}^{[\![1,k'']\!]\times[\![1,n]\!]\times[\![1,d]\!]} \\ \sum \alpha_l = \lambda \\ k + \sum |\alpha_l| l = k''}} \Bigg[ \prod_{l=1}^{k''} \frac{1}{\alpha_l! l!^{|\alpha_l|}} \prod_{N=1}^{n} \left( \frac{\mathrm{d}^{l+n-N} \hat{\mathrm{y}}}{\mathrm{d}t^{l+n-N}} \right)^{\alpha_{l,N}} \Bigg] = \underset{h\to 0}{\mathcal{O}}\left(h^v\right)$$

We merge the first two sums of each side :

$$\sum_{j'=1}^{s} w_{N,j,j'} = \xi_{j,N}^{N},$$

$$\forall k \in [\![0, v-1]\!], \quad \sum_{\substack{\lambda \in \mathbb{N}^{[\![1,n]\!]\times[\![1,d]\!]} \\ \max(1,k) \leq |\lambda| + k \leq v-1}} \frac{1}{\lambda!} \frac{\partial^{|\lambda| + k} f}{\partial t^k \prod_{N=1}^{n} \left(\partial y_N\right)^{\lambda_N}} \sum_{j'=1}^{s} \Bigg[ w_{N,j,j'} \left(\tau_{j'} h\right)^k \prod_{N=1}^{n} \left(y_{n,j',N} - y_{n,0,N}\right)^{\lambda_N} \Bigg]$$

$$- \sum_{\substack{\lambda \in \mathbb{N}^{[\![1,n]\!]\times[\![1,d]\!]} \\ \max(1,k) \leq |\lambda| + k \leq v-1}} \frac{\partial^{k+|\lambda|} f}{\partial t^k \prod_{N=1}^{n} \left(\partial y_N\right)^{\lambda_N}} \sum_{k''=k+|\lambda|}^{v-1} \xi_{j,N}^{N+k''} h^{k''} \frac{N!k''!}{(k''+N)!} \sum_{\substack{\alpha \in \mathbb{N}^{[\![1,k'']\!]\times[\![1,n]\!]\times[\![1,d]\!]} \\ \sum \alpha_l = \lambda \\ k + \sum |\alpha_l| l = k''}} \Bigg[ \prod_{l=1}^{k''} \frac{1}{\alpha_l! l!^{|\alpha_l|}} \prod_{N=1}^{n} \left( \frac{\mathrm{d}^{l+n-N} \hat{\mathrm{y}}}{\mathrm{d}t^{l+n-N}} \right)^{\alpha_{l,N}} \Bigg] = \underset{h\to 0}{\mathcal{O}}\left(h^v\right)$$

We identify the partial derivatives of $f$ and replace $\left(y_{n,j',N} - y_{n,0,N}\right)^{\lambda_N}$ with its expression. To make the notations shorter we note $F_j$ the evaluation $f\left(t_0 + \tau_j h, y_{n,j}\right)$.

$$\sum_{j'=1}^{s} w_{N,j,j'} = \xi_{j,N}^{N},$$

$$\forall k \in \llbracket 0, v-1 \rrbracket,\ \ \forall \lambda \in \mathbb{N}^{\llbracket 1,n \rrbracket \times \llbracket 1,d \rrbracket}, \max(1,k) \leq |\lambda| + k \leq v-1$$

$$\frac{1}{\lambda!} \sum_{j'=1}^{s} w_{N,j,j'} \left(\tau_{j'} h\right)^{k} \prod_{N=1}^{n} \prod_{a=1}^{d} \left( \sum_{N'=1}^{N-1} \frac{\left(\xi_{j',N} h\right)^{N'}}{N'!} y_{n,0,N-N',a} + \frac{h^{N}}{N!} \sum_{j''=1}^{s} w_{N,j',j''} F_{j'',a} \right)^{\lambda_{N,a}}$$

$$- \sum_{k''=k+|\lambda|}^{v-1} \xi_{j,N}^{N+k''} h^{k''} \frac{N!k''!}{(k''+N)!} \sum_{\substack{\alpha \in \mathbb{N}^{\llbracket 1,k'' \rrbracket \times \llbracket 1,n \rrbracket \times \llbracket 1,d \rrbracket} \\ \sum \alpha_l = \lambda \\ k+\sum |\alpha_l|\ l = k''}} \prod_{l=1}^{k''} \frac{1}{\alpha_l! l!^{|\alpha_l|}} \prod_{N=1}^{n} \left( \frac{\mathrm{d}^{l+n-N} \hat{\mathrm{y}}}{\mathrm{d}t^{l+n-N}} \right)^{\alpha_{l,N}} = \underset{h \to 0}{\mathcal{O}}\left(h^{v}\right)$$

We develop the product of the first term.

$$\sum_{j'=1}^{s} w_{N,j,j'} = \xi_{j,N}^{N},$$

$$\forall k \in \llbracket 0, v-1 \rrbracket,\ \ \forall \lambda \in \mathbb{N}^{\llbracket 1,n \rrbracket \times \llbracket 1,d \rrbracket}, \max(1,k) \leq |\lambda| + k \leq v-1$$

$$\sum_{j'=1}^{s} \frac{w_{N,j,j'}}{\lambda!} \left(\tau_{j'} h\right)^{k} \prod_{N=1}^{n} \prod_{a=1}^{d} \sum_{\substack{\alpha \in \mathbb{N}^{\llbracket 1,N \rrbracket} \\ |\alpha| = \lambda_{N,a}}} \frac{\lambda_{N,a}!}{\alpha!} \left( \frac{h^{N}}{N!} \sum_{j''=1}^{s} w_{N,j',j''} F_{j'',a} \right)^{\alpha_N} \prod_{l=1}^{N-1} \left( \frac{\left(\xi_{j',N} h\right)^{l}}{l!} y_{n,0,N-l,a} \right)^{\alpha_l}$$

$$- \sum_{k''=k+|\lambda|}^{v-1} \xi_{j,N}^{N+k''} h^{k''} \frac{N!k''!}{(k''+N)!} \sum_{\substack{\alpha \in \mathbb{N}^{\llbracket 1,k'' \rrbracket \times \llbracket 1,n \rrbracket \times \llbracket 1,d \rrbracket} \\ \sum \alpha_l = \lambda \\ k+\sum |\alpha_l|\ l = k''}} \prod_{l=1}^{k''} \frac{1}{\alpha_l! l!^{|\alpha_l|}} \prod_{N=1}^{n} \left( \frac{\mathrm{d}^{l+n-N} \hat{\mathrm{y}}}{\mathrm{d}t^{l+n-N}} \right)^{\alpha_{l,N}} = \underset{h \to 0}{\mathcal{O}}\left(h^{v}\right)$$

We then permute the product and the sum.

$$\sum_{j'=1}^{s} w_{N,j,j'} = \xi_{j,N}^{N},$$

$$\forall k \in \llbracket 0, v-1 \rrbracket,\ \ \forall \lambda \in \mathbb{N}^{\llbracket 1,n \rrbracket \times \llbracket 1,d \rrbracket}, \max(1,k) \leq |\lambda| + k \leq v-1$$

$$\frac{1}{\lambda!} \sum_{j'=1}^{s} w_{N,j,j'} \left(\tau_{j'} h\right)^{k} \sum_{\substack{\alpha \in \prod_{N=1}^{n} \prod_{a=1}^{d} \mathbb{N}^{\llbracket 1,N \rrbracket} \\ \forall (N,a)\ |\alpha_{N,a}| = \lambda_{N,a}}} \Bigg[$$

$$\prod_{N=1}^{n} \prod_{a=1}^{d} \frac{\lambda_{N,a}!}{\alpha_{N,a}!} h^{\sum_{l=1}^{N} l \alpha_{N,a,l}} \left( \sum_{j''=1}^{s} \frac{w_{N,j',j''}}{N!} F_{j'',a} \right)^{\alpha_{N,a,N}} \prod_{l=1}^{N-1} \left( \frac{\xi_{j',N}^{l}}{l!} y_{n,0,N-l,a} \right)^{\alpha_{N,a,l}} \Bigg]$$

$$- \sum_{k''=k+|\lambda|}^{v-1} \xi_{j,N}^{N+k''} h^{k''} \frac{N!k''!}{(k''+N)!} \sum_{\substack{\alpha \in \mathbb{N}^{\llbracket 1,k'' \rrbracket \times \llbracket 1,n \rrbracket \times \llbracket 1,d \rrbracket} \\ \sum \alpha_l = \lambda \\ k+\sum |\alpha_l|\ l = k''}} \prod_{l=1}^{k''} \frac{1}{\alpha_l! l!^{|\alpha_l|}} \prod_{N=1}^{n} \left( \frac{\mathrm{d}^{l+n-N} \hat{\mathrm{y}}}{\mathrm{d}t^{l+n-N}} \right)^{\alpha_{l,N}} = \underset{h \to 0}{\mathcal{O}}\left(h^{v}\right)$$

Which is equivalent to :

$$\sum_{j'=1}^{s} w_{N,j,j'} = \xi_{j,N}^{N},$$

$$\forall k \in [\![0, v-1]\!],\ \forall \lambda \in \mathbb{N}^{[\![1,n]\!]\times[\![1,d]\!]}, \max(1,k) \leq |\lambda| + k \leq v-1$$

$$\frac{1}{\lambda!} \sum_{\substack{\alpha \in \prod_{N=1}^{n} \prod_{a=1}^{d} \mathbb{N}^{[\![1,N]\!]} \\ \forall (N,a)\ |\alpha_{N,a}| = \lambda_{N,a}}} \frac{\lambda!}{\alpha!} h^{k + \sum_{N=1}^{n} \sum_{a=1}^{d} \sum_{l=1}^{N} l \alpha_{N,a,l}} \sum_{j'=1}^{s} w_{N,j,j'} \tau_{j'}^{k} \Bigg[ \prod_{N=1}^{n} \prod_{a=1}^{d} \left( \sum_{j''=1}^{s} \frac{w_{N,j',j''}}{N!} F_{j'',a} \right)^{\alpha_{N,a,N}} \prod_{l=1}^{N-1} \left( \frac{\xi_{j',N}^{l}}{l!} y_{n,0,N-l,a} \right)^{\alpha_{N,a,l}} \Bigg]$$

$$- \sum_{k''=k+|\lambda|}^{v-1} \xi_{j,N}^{N+k''} h^{k''} \frac{N!k''!}{(k''+N)!} \sum_{\substack{\alpha \in \mathbb{N}^{[\![1,k'']\!]\times[\![1,n]\!]\times[\![1,d]\!]} \\ \sum \alpha_l = \lambda \\ k + \sum |\alpha_l|\ l = k''}} \prod_{l=1}^{k''} \frac{1}{\alpha_l! l!^{|\alpha_l|}} \prod_{N=1}^{n} \left( \frac{\mathrm{d}^{l+n-N} \hat{\mathrm{y}}}{\mathrm{d}t^{l+n-N}} \right)^{\alpha_{l,N}} \underset{h \to 0}{=} \mathcal{O}(h^v)$$

Since $k + \sum_{N=1}^{n} \sum_{a=1}^{d} \sum_{l=1}^{N} l\alpha_{N,a,l} \geq k + \sum_{N=1}^{n} \sum_{a=1}^{d} \sum_{l=1}^{N} \alpha_{N,a,l} \geq k + |\lambda|$ and every power of $h$ superior to $v-1$ disappear, we can clamp $k''$ to these values.

$$\sum_{j'=1}^{s} w_{N,j,j'} = \xi_{j,N}^{N},$$

$$\forall k \in [\![0, v-1]\!],\ \forall \lambda \in \mathbb{N}^{[\![1,n]\!]\times[\![1,d]\!]}, \max(1,k) \leq |\lambda| + k \leq v-1$$

$$\sum_{k''=k+|\lambda|}^{v-1} h^{k''} \sum_{\substack{\alpha \in \prod_{N=1}^{n} \prod_{a=1}^{d} \mathbb{N}^{[\![1,N]\!]} \\ \forall (N,a)\ |\alpha_{N,a}| = \lambda_{N,a} \\ k + \sum_{N=1}^{n} \sum_{a=1}^{d} \sum_{l=1}^{N} l\alpha_{N,a,l} = k''}} \frac{1}{\alpha!} \sum_{j'=1}^{s} \Bigg[ w_{N,j,j'} \tau_{j'}^{k} \prod_{N=1}^{n} \prod_{a=1}^{d} \left( \sum_{j''=1}^{s} \frac{w_{N,j',j''}}{N!} F_{j'',a} \right)^{\alpha_{N,a,N}} \prod_{l=1}^{N-1} \left( \frac{\xi_{j',N}^{l}}{l!} y_{n,0,N-l,a} \right)^{\alpha_{N,a,l}} \Bigg]$$

$$- \sum_{k''=k+|\lambda|}^{v-1} \xi_{j,N}^{N+k''} h^{k''} \frac{N!k''!}{(k''+N)!} \sum_{\substack{\alpha \in \mathbb{N}^{[\![1,k'']\!]\times[\![1,n]\!]\times[\![1,d]\!]} \\ \sum \alpha_l = \lambda \\ k + \sum |\alpha_l|\ l = k''}} \prod_{l=1}^{k''} \frac{1}{\alpha_l! l!^{|\alpha_l|}} \prod_{N=1}^{n} \left( \frac{\mathrm{d}^{l+n-N} \hat{\mathrm{y}}}{\mathrm{d}t^{l+n-N}} \right)^{\alpha_{l,N}} \underset{h \to 0}{=} \mathcal{O}(h^v)$$

We now identify the terms which do not involve the partial derivatives of $f$ and leave out the other terms. For the first part it implies replacing $\sum_{j''=1}^{s} w_{N,j',j''} F_{j'',a}$ with $\sum_{j''=1}^{s} w_{N,j',j''} f(t_0, y_{n,0})$. Let $y_{n,0,0} = f(t_0, y_{n,0})$. The first condition gives $\sum_{j''=1}^{s} w_{N,j',j''} = \xi_{j',N}^{N}$, we hence replace $\sum_{j''=1}^{s} w_{N,j',j''} F_{j'',a}$ with $\xi_{j',N}^{N} y_{n,0,0}$. For the second part it is necessary and sufficient that $l + n - N \geq n$ (thus $l > N$) implies $\alpha_{l,N} = 0_d$ :

$$\sum_{j'=1}^{s} w_{N,j,j'} = \xi_{j,N}^{N},$$

$$\forall k \in [\![0, v-1]\!],\ \forall \lambda \in \mathbb{N}^{[\![1,n]\!]\times[\![1,d]\!]}, \max(1,k) \leq |\lambda| + k \leq v-1$$

$$\sum_{k''=k+|\lambda|}^{v-1} h^{k''} \sum_{\substack{\alpha\in\prod_{N=1}^{n}\prod_{a=1}^{d}\mathbb{N}^{[\![1,N]\!]} \\ \forall(N,a)\ |\alpha_{N,a}|=\lambda_{N,a} \\ k+\sum_{N=1}^{n}\sum_{a=1}^{d}\sum_{l=1}^{N-1} l\alpha_{N,a,l}=k''}} \frac{1}{\alpha!}\sum_{j'=1}^{s}\Bigg[ w_{N,j,j'}\tau_{j'}^{k}\prod_{N=1}^{n}\prod_{a=1}^{d}\prod_{l=1}^{N}\left(\frac{\xi_{j',N}^{l}}{l!}\right)^{\alpha_{N,a,l}}\left(y_{n,0,N-l}^{[a]}\right)^{\alpha_{N,a,l}}\Bigg]$$

$$-\sum_{k''=k+|\lambda|}^{v-1}\xi_{j,N}^{N+k''}h^{k''}\frac{N!k''!}{(k''+N)!}\sum_{\substack{\alpha\in\mathbb{N}^{[\![1,k'']\!]\times[\![1,n]\!]\times[\![1,d]\!]} \\ \sum\alpha_l=\lambda \\ k+\sum|\alpha_l|\ l=k'' \\ l>N\Rightarrow\alpha_{l,N}=0_d}}\prod_{l=1}^{k''}\frac{1}{\alpha_l!l!^{|\alpha_l|}}\prod_{N=l}^{n}y_{n,0,N-l}^{\alpha_{l,N}} \underset{h\to 0}{=} \mathcal{O}(h^v)$$

We identify by the powers of $h$ and swap the coordinates of $\alpha$.

$$\sum_{j'=1}^{s} w_{N,j,j'} = \xi_{j,N}^{N},$$

$$\forall k\in[\![0,v-1]\!],\ \forall\lambda\in\mathbb{N}^{[\![1,n]\!]\times[\![1,d]\!]}, \max(1,k)\leq|\lambda|+k\leq v-1,\ \forall k''\in[\![k+|\lambda|,v-1]\!]$$

$$\sum_{\substack{\alpha\in\left(\prod_{N=1}^{n}\mathbb{N}^{[\![1,N]\!]}\right)^{[\![1,d]\!]} \\ \forall(N,a)\sum_{l=1}^{N-1}\alpha_{N,l,a}=\lambda_{N,a} \\ k+\sum_{N=1}^{n}\sum_{a=1}^{d}\sum_{l=1}^{N-1}l\alpha_{N,l,a}=k''}} \frac{1}{\alpha!}\sum_{j'=1}^{s}w_{N,j,j'}\tau_{j'}^{k}\prod_{N=1}^{n}\prod_{a=1}^{d}\prod_{l=1}^{N}\left(\frac{\xi_{j',N}^{l}}{l!}\right)^{\alpha_{N,l,a}}\left(y_{n,0,N-l}^{[a]}\right)^{\alpha_{N,l,a}}$$

$$=\xi_{j,N}^{N+k''}\frac{N!k''!}{(k''+N)!}\sum_{\substack{\alpha\in\mathbb{N}^{[\![1,k'']\!]\times[\![1,n]\!]\times[\![1,d]\!]} \\ \sum\alpha_l=\lambda \\ k+\sum|\alpha_l|\ l=k'' \\ l>N\Rightarrow\alpha_{l,N}=0_d}}\prod_{l=1}^{k''}\frac{1}{\alpha_l!l!^{|\alpha_l|}}\prod_{N=l}^{n}y_{n,0,N-l}^{\alpha_{l,N}}$$

Which is equivalent to :

$$\sum_{j'=1}^{s} w_{N,j,j'} = \xi_{j,N}^{N},$$

$$\forall k\in[\![0,v-1]\!],\ \forall\lambda\in\mathbb{N}^{[\![1,n]\!]\times[\![1,d]\!]}, \max(1,k)\leq|\lambda|+k\leq v-1,\ \forall k''\in[\![k+|\lambda|,v-1]\!]$$

$$\sum_{\substack{\alpha\in\left(\prod_{N=1}^{n}\mathbb{N}^{[\![1,N]\!]}\right)^{[\![1,d]\!]} \\ \forall(N,a)\sum_{l=1}^{N}\alpha_{N,l,a}=\lambda_{N,a} \\ k+\sum_{N=1}^{n}\sum_{l=1}^{N}l\ |\alpha_{N,l}|=k''}} \frac{1}{\alpha!}\sum_{j'=1}^{s}w_{N,j,j'}\tau_{j'}^{k}\left(\prod_{N=1}^{n}\prod_{a=1}^{d}\prod_{l=1}^{N}\xi_{j',N}^{l\alpha_{N,l,a}}\right)\left(\prod_{N=1}^{n}\prod_{a=1}^{d}\prod_{l=1}^{N}\frac{1}{l!^{\alpha_{N,l,a}}}\right)\prod_{N=1}^{n}\prod_{l=1}^{N}y_{n,0,N-l}^{\alpha_{N,l}}$$

$$=\xi_{j,N}^{N+k''}\frac{N!k''!}{(k''+N)!}\sum_{\substack{\alpha\in\mathbb{N}^{[\![1,k'']\!]\times[\![1,n]\!]\times[\![1,d]\!]} \\ \sum\alpha_l=\lambda \\ k+\sum|\alpha_l|\ l=k'' \\ l>N\Rightarrow\alpha_{l,N}=0_d}}\left(\prod_{l=1}^{k''}\frac{1}{\alpha_l!l!^{|\alpha_l|}}\right)\prod_{l=1}^{k''}\prod_{N=0}^{n-l}y_{n,0,N}^{\alpha_{l,N+l}}$$

Let $a\in\mathbb{R}^{\mathbb{N}^*\times\mathbb{N}^*}$. We have :

$$\prod_{l=1}^{k''}\prod_{N=0}^{n-l} a_{l,N} = \prod_{l=n-k''}^{n-1}\prod_{N=0}^{l} a_{n-l,N} = \prod_{N=0}^{n-1}\prod_{l=\max(N,n-k'')}^{n-1} a_{n-l,N} = \prod_{N=0}^{n-1}\prod_{l=1}^{\min(n-N,k'')} a_{l,N}$$

We deduce :

$$\sum_{j'=1}^{s} w_{N,j,j'} = \xi_{j,N}^{N},$$

$$\forall k \in [\![0, v-1]\!],\ \forall \lambda \in \mathbb{N}^{[\![1,n]\!]\times[\![1,d]\!]}, \max(1,k) \leq |\lambda| + k \leq v-1,\ \forall k'' \in [\![k+|\lambda|, v-1]\!]$$

$$\sum_{\substack{\alpha\in\left(\prod_{N=1}^{n}\mathbb{N}^{[\![1,N]\!]}\right)^{[\![1,d]\!]} \\ \forall N \sum_{l=1}^{N}\alpha_{N,l}=\lambda_N \\ k+\sum_{N=1}^{n}\sum_{l=1}^{N} l\ |\alpha_{N,l}|\ =k''}} \frac{1}{\alpha!}\sum_{j'=1}^{s} w_{N,j,j'}\tau_{j'}^{k}\left(\prod_{N=1}^{n}\xi_{j',N}^{\sum_{l=1}^{N} l\ |\alpha_{N,l}|}\right)\left(\prod_{N=1}^{n}\prod_{l=1}^{N}\frac{1}{l!^{|\alpha_{N,l}|}}\right)\prod_{l=1}^{n}\prod_{N=l}^{n} y_{n,0,N-l}^{\alpha_{N,l}}$$

$$= \xi_{j,N}^{N+k''}\frac{N!k''!}{(k''+N)!}\sum_{\substack{\alpha\in\mathbb{N}^{[\![1,k'']\!]\times[\![1,n]\!]\times[\![1,d]\!]} \\ \sum\alpha_l=\lambda \\ k+\sum|\alpha_l|\ l=k'' \\ l>N\Rightarrow\alpha_{l,N}=0_d}}\frac{1}{\alpha!}\left(\prod_{l=1}^{k''}\frac{1}{l!^{|\alpha_l|}}\right)\prod_{N=0}^{n-1}\prod_{l=1}^{\min(n-N,k'')} y_{n,0,N}^{\alpha_{l,N+l}}$$

We add zero coordinates to $\alpha$.

$$\sum_{j'=1}^{s} w_{N,j,j'} = \xi_{j,N}^{N},$$

$$\forall k \in [\![0, v-1]\!],\ \forall \lambda \in \mathbb{N}^{[\![1,n]\!]\times[\![1,d]\!]}, \max(1,k) \leq |\lambda| + k \leq v-1,\ \forall k'' \in [\![k+|\lambda|, v-1]\!]$$

$$\sum_{\substack{\alpha\in\mathbb{N}^{[\![1,n]\!]\times[\![1,n]\!]\times[\![1,d]\!]} \\ \forall N \sum_{l=1}^{N}\alpha_{N,l}=\lambda_N \\ k+\sum_{N=1}^{n}\sum_{l=1}^{N-1} l\ |\alpha_{N,l}|\ =k'' \\ l>N\Rightarrow\alpha_{N,l}=0_d}} \frac{1}{\alpha!}\sum_{j'=1}^{s} w_{N,j,j'}\tau_{j'}^{k}\left(\prod_{N=1}^{n}\xi_{j',N}^{\sum_{l=1}^{n} l\ |\alpha_{N,l}|}\right)\left(\prod_{N=1}^{n}\prod_{l=1}^{n}\frac{1}{l!^{|\alpha_{N,l}|}}\right)\prod_{l=1}^{n}\prod_{N=0}^{n-l} y_{n,0,N}^{\alpha_{N+l,l}}$$

$$= \xi_{j,N}^{N+k''}\frac{N!k''!}{(k''+N)!}\sum_{\substack{\alpha\in\mathbb{N}^{[\![1,k'']\!]\times[\![1,n]\!]\times[\![1,d]\!]} \\ \sum\alpha_l=\lambda \\ k+\sum|\alpha_l|\ l=k'' \\ l>N\Rightarrow\alpha_{l,N}=0_d}}\frac{1}{\alpha!}\left(\prod_{l=1}^{k''}\frac{1}{l!^{|\alpha_l|}}\right)\prod_{N=0}^{n-1} y_{n,0,N}^{\sum_{l=1}^{\min(n-N,k'')}\alpha_{l,N+l}}$$

Let $a \in \mathbb{R}^{\mathbb{N}^*\times\mathbb{N}^*}$. We have :

$$\prod_{l=1}^{n}\prod_{N=0}^{n-l} a_{l,N} = \prod_{l=0}^{n-1}\prod_{N=0}^{l} a_{n-l,N} = \prod_{N=0}^{n-1}\prod_{l=N}^{n-1} a_{n-l,N} = \prod_{N=0}^{n-1}\prod_{l=1}^{n-N} a_{l,N}$$

Therefore :

$$\sum_{j'=1}^{s} w_{N,j,j'} = \xi_{j,N}^{N},$$

$$\forall k \in [\![0, v-1]\!],\ \ \forall \lambda \in \mathbb{N}^{[\![1,n]\!]\times[\![1,d]\!]}, \max(1,k) \leq |\lambda| + k \leq v-1,\ \ \forall k'' \in [\![k+|\lambda|, v-1]\!]$$

$$\sum_{\substack{\alpha\in\mathbb{N}^{[\![1,n]\!]\times[\![1,n]\!]\times[\![1,d]\!]} \\ \forall N \sum_{l=1}^{n} \alpha_{N,l}=\lambda_N \\ k+\sum_{N=1}^{n}\sum_{l=1}^{n} l\ |\alpha_{N,l}|\ =k'' \\ l>N \Rightarrow \alpha_{N,l}=0_d}} \frac{1}{\alpha!} \sum_{j'=1}^{s} w_{N,j,j'} \tau_{j'}^{k} \left( \prod_{N=1}^{n} \xi_{j',N}^{\sum_{l=1}^{n} l\ |\alpha_{N,l}|} \right) \left( \prod_{l=1}^{n} \frac{1}{l!^{\sum_{N=1}^{n} |\alpha_{N,l}|}} \right) \prod_{N=0}^{n-1} y_{n,0,N}^{\sum_{l=1}^{n-N} \alpha_{N+l,l}}$$

$$= \xi_{j,N}^{N+k''} \frac{N!k''!}{(k''+N)!} \sum_{\substack{\alpha\in\mathbb{N}^{[\![1,k'']\!]\times[\![1,n]\!]\times[\![1,d]\!]} \\ \sum \alpha_l=\lambda \\ k+\sum|\alpha_l|\ l=k'' \\ l>N \Rightarrow \alpha_{l,N}=0_d}} \frac{1}{\alpha!} \left( \prod_{l=1}^{k''} \frac{1}{l!^{|\alpha_l|}} \right) \prod_{N=0}^{n-1} y_{n,0,N}^{\sum_{l=1}^{\min(n-N,k'')} \alpha_{l,N+l}}$$

We swap the first two coordinates of $\alpha$.

$$\sum_{j'=1}^{s} w_{N,j,j'} = \xi_{j,N}^{N},$$

$$\forall k \in [\![0, v-1]\!],\ \ \forall \lambda \in \mathbb{N}^{[\![1,n]\!]\times[\![1,d]\!]}, \max(1,k) \leq |\lambda| + k \leq v-1,\ \ \forall k'' \in [\![k+|\lambda|, v-1]\!]$$

$$\sum_{\substack{\alpha\in\mathbb{N}^{[\![1,n]\!]\times[\![1,n]\!]\times[\![1,d]\!]} \\ \sum_{l=1}^{n} \alpha_l=\lambda \\ k+\sum_{l=1}^{n} l\ |\alpha_l|\ =k'' \\ l>N \Rightarrow \alpha_{l,N}=0_d}} \frac{1}{\alpha!} \sum_{j'=1}^{s} w_{N,j,j'} \tau_{j'}^{k} \left( \prod_{N=1}^{n} \xi_{j',N}^{\sum_{l=1}^{n} l\ |\alpha_{l,N}|} \right) \left( \prod_{l=1}^{n} \frac{1}{l!^{|\alpha_l|}} \right) \prod_{N=0}^{n-1} y_{n,0,N}^{\sum_{l=1}^{n-N} \alpha_{l,l+N}}$$

$$= \sum_{\substack{\alpha\in\mathbb{N}^{[\![1,k'']\!]\times[\![1,n]\!]\times[\![1,d]\!]} \\ \sum \alpha_l=\lambda \\ k+\sum|\alpha_l|\ l=k'' \\ l\geq N \Rightarrow \alpha_{l,N}=0_d}} \frac{\xi_{j,N}^{N+k''}}{\alpha!} \frac{N!k''!}{(k''+N)!} \left( \prod_{l=1}^{k''} \frac{1}{l!^{|\alpha_l|}} \right) \prod_{N=0}^{n-1} y_{n,0,N}^{\sum_{l=1}^{\min(n-N,k'')} \alpha_{l,N+l}}$$

For both $\alpha$, if there exists $l > k'' - k$ such that $\alpha_{l,N,a} > 0$, then $k + \sum_{l=1}^{n} l\ |\alpha_l| > k''$, we thus clamp the coordinates of $\alpha$ to $k''$. For the second $\alpha$ we have $l > N$ implies $\alpha_{l,N} = 0_d$, so we clamp its coordinates to $n$.

$$\sum_{j'=1}^{s} w_{N,j,j'} = \xi_{j,N}^{N},$$

$$\forall k \in [\![0, v-1]\!],\ \ \forall \lambda \in \mathbb{N}^{[\![1,n]\!]\times[\![1,d]\!]}, \max(1,k) \leq |\lambda| + k \leq v-1,\ \ \forall k'' \in [\![k+|\lambda|, v-1]\!]$$

$$\sum_{\substack{\alpha\in\mathbb{N}^{\llbracket 1,\min(n,k''-k)\rrbracket\times\llbracket 1,n\rrbracket\times\llbracket 1,d\rrbracket}\\ \sum\alpha_l=\lambda\\ k+\sum l\ |\alpha_l|\ =k''\\ l>N\Rightarrow\alpha_{l,N}=0_d}} \frac{1}{\alpha!}\sum_{j'=1}^{s} w_{N,j,j'}\tau_{j'}^{k}\Bigg[$$

$$\left(\prod_{N=1}^{n}\xi_{j',N}^{\sum_{l=1}^{\min(n,k''-k)} l\ |\alpha_{l,N}|}\right)\left(\prod_{l=1}^{\min(n,k''-k)}\frac{1}{l!^{|\alpha_l|}}\right)\prod_{N=0}^{n-1} y_{n,0,N}^{\sum_{l=1}^{\min(n-N,k''-k)}\alpha_{l,l+N}}\Bigg]$$

$$=\sum_{\substack{\alpha\in\mathbb{N}^{\llbracket 1,\min(n,k''-k)\rrbracket\times\llbracket 1,n\rrbracket\times\llbracket 1,d\rrbracket}\\ \sum\alpha_l=\lambda\\ k+\sum|\alpha_l|\ l=k''\\ l>N\Rightarrow\alpha_{l,N}=0_d}} \frac{\xi_{j,N}^{N+k''}}{\alpha!}\frac{N!k''!}{(k''+N)!}\left(\prod_{l=1}^{\min(n,k''-k)}\frac{1}{l!^{|\alpha_l|}}\right)\prod_{N=0}^{n-1} y_{n,0,N}^{\sum_{l=1}^{\min(n-N,k''-k)}\alpha_{l,N+l}}$$

We introduce a new variable $\beta$ such that $|\beta| = \sum_{N=1}^{n}\sum_{l=1}^{\min(n-N,k'')}|\alpha_{l,l+N}| = |\alpha| = |\lambda|$.

$$\sum_{j'=1}^{s} w_{N,j,j'} = \xi_{j,N}^{N},$$

$$\forall k\in\llbracket 0,v-1\rrbracket,\ \forall\lambda\in\mathbb{N}^{\llbracket 1,n\rrbracket\times\llbracket 1,d\rrbracket},\max(1,k)\leq|\lambda|+k\leq v-1,\ \forall k''\in\llbracket k+|\lambda|,v-1\rrbracket$$

$$\sum_{\substack{\beta\in\mathbb{N}^{\llbracket 0,n-1\rrbracket\times\llbracket 1,d\rrbracket}\\ |\beta|\,=\,|\lambda|}}\prod_{N=0}^{n-1} y_{n,0,N}^{\beta_N}\sum_{\substack{\alpha\in\mathbb{N}^{\llbracket 1,\min(n,k''-k)\rrbracket\times\llbracket 1,n\rrbracket\times\llbracket 1,d\rrbracket}\\ \sum\alpha_l=\lambda\\ k+\sum l\ |\alpha_l|\ =k''\\ l>N\Rightarrow\alpha_{l,N}=0_d\\ \forall N,\sum_{l=1}^{\min(n-N,k''-k)}\alpha_{l,l+N}=\beta_N}} \frac{1}{\alpha!}\sum_{j'=1}^{s}\Bigg[$$

$$w_{N,j,j'}\tau_{j'}^{k}\left(\prod_{N=1}^{n}\xi_{j',N}^{\sum_{l=1}^{\min(n,k''-k)} l\ |\alpha_{l,N}|}\right)\left(\prod_{l=1}^{\min(n,k''-k)}\frac{1}{l!^{|\alpha_l|}}\right)\Bigg]$$

$$=\sum_{\substack{\beta\in\mathbb{N}^{\llbracket 0,n-1\rrbracket\times\llbracket 1,d\rrbracket}\\ |\beta|\,=\,|\lambda|}}\prod_{N=0}^{n-1} y_{n,0,N}^{\beta_N}\sum_{\substack{\alpha\in\mathbb{N}^{\llbracket 1,\min(n,k''-k)\rrbracket\times\llbracket 1,n\rrbracket\times\llbracket 1,d\rrbracket}\\ \sum\alpha_l=\lambda\\ k+\sum|\alpha_l|\ l=k''\\ l>N\Rightarrow\alpha_{l,N}=0_d\\ \forall N,\sum_{l=1}^{\min(n-N,k''-k)}\alpha_{l,l+N}=\beta_N}} \frac{\xi_{j,N}^{N+k''}}{\alpha!}\frac{N!k''!}{(k''+N)!}\left(\prod_{l=1}^{\min(n,k''-k)}\frac{1}{l!^{|\alpha_l|}}\right)$$

We identify by the initial values :

$$\sum_{j'=1}^{s} w_{N,j,j'} = \xi_{j,N}^{N},$$

$$\forall k\in\llbracket 0,v-1\rrbracket,\ \forall\lambda\in\mathbb{N}^{\llbracket 1,n\rrbracket\times\llbracket 1,d\rrbracket},\max(1,k)\leq|\lambda|+k\leq v-1,$$
$$\forall k''\in\llbracket k+|\lambda|,v-1\rrbracket,\forall\beta\in\mathbb{N}^{\llbracket 0,n-1\rrbracket\times\llbracket 1,d\rrbracket},|\beta|=|\lambda|$$

$$\sum_{\substack{\alpha\in\mathbb{N}^{[\![1,\,\min(n,k''-k)]\!]\times[\![1,n]\!]\times[\![1,d]\!]}\\ \sum\alpha_l=\lambda\\ k+\sum l\ |\alpha_l|\ =k''\\ l>N\Rightarrow\alpha_{l,N}=0_d\\ \forall N,\sum_{l=1}^{\min(n-N,k''-k)}\alpha_{l,l+N}=\beta_N}} \frac{1}{\alpha!}\sum_{j'=1}^{s} w_{N,j,j'}\tau_{j'}^{k}\left(\prod_{N=1}^{n}\xi_{j',N}^{\sum_{l=1}^{\min(n,k''-k)} l\ |\alpha_{l,N}|}\right)\left(\prod_{l=1}^{\min(n,k''-k)}\frac{1}{l!^{|\alpha_l|}}\right)$$

$$= \sum_{\substack{\alpha\in\mathbb{N}^{[\![1,\,\min(n,k''-k)]\!]\times[\![1,n]\!]\times[\![1,d]\!]}\\ \sum\alpha_l=\lambda\\ k+\sum|\alpha_l|\ l=k''\\ l>N\Rightarrow\alpha_{l,N}=0_d\\ \forall N,\sum_{l=1}^{\min(n-N,k''-k)}\alpha_{l,l+N}=\beta_N}} \frac{\xi_{j,N}^{N+k''}}{\alpha!}\frac{N!k''!}{(k''+N)!}\left(\prod_{l=1}^{\min(n,k''-k)}\frac{1}{l!^{|\alpha_l|}}\right)$$

We clear things up.

$$\sum_{j'=1}^{s} w_{N,j,j'}=\xi_{j,N}^{N},$$
$$\forall k\in[\![0,v-1]\!],\ \ \forall\lambda\in\mathbb{N}^{[\![1,n]\!]\times[\![1,d]\!]},\max(1,k)\leq|\lambda|+k\leq v-1,$$
$$\forall k''\in[\![k+|\lambda|,v-1]\!],\forall\beta\in\mathbb{N}^{[\![0,n-1]\!]\times[\![1,d]\!]},|\beta|=|\lambda|$$
$$\sum_{\substack{\alpha\in\mathbb{N}^{[\![1,\,\min(n,k''-k)]\!]\times[\![1,n]\!]\times[\![1,d]\!]}\\ \sum\alpha_l=\lambda\\ k+\sum l\ |\alpha_l|\ =k''\\ l>N\Rightarrow\alpha_{l,N}=0_d\\ \forall N,\sum_{l=1}^{\min(n-N,k''-k)}\alpha_{l,l+N}=\beta_N}} \left(\prod_{l=1}^{\min(n,k''-k)}\frac{1}{l!^{|\alpha_l|}}\right)\frac{1}{\alpha!}\Bigg[$$
$$\sum_{j'=1}^{s} w_{N,j,j'}\tau_{j'}^{k}\left(\prod_{N=1}^{n}\xi_{j',N}^{\sum_{l=1}^{\min(n,k''-k)} l\ |\alpha_{l,N}|}\right)-\xi_{j,N}^{N+k''}\frac{N!k''!}{(k''+N)!}\Bigg]=0$$

We now prove the equivalence of the generalized solved system conditions and this system.
In the previous equation, if we consider the case $k''=k+|\lambda|$, we have :

$$\sum_{j'=1}^{s} w_{N,j,j'}=\xi_{j,N}^{N},$$
$$\forall k\in[\![0,v-1]\!],\ \ \forall\lambda\in\mathbb{N}^{[\![1,n]\!]\times[\![1,d]\!]},\max(1,k)\leq|\lambda|+k\leq v-1,\ \ \forall\beta\in\mathbb{N}^{[\![0,n-1]\!]\times[\![1,d]\!]},|\beta|=|\lambda|$$
$$\sum_{\substack{\alpha\in\mathbb{N}^{[\![1,\,\min(n,|\lambda|)]\!]\times[\![1,n]\!]\times[\![1,d]\!]}\\ \sum\alpha_l=\lambda\\ \sum l\ |\alpha_l|\ =\ |\lambda|\\ l>N\Rightarrow\alpha_{l,N}=0_d\\ \forall N,\sum_{l=1}^{\min(n-N,|\lambda|)}\alpha_{l,l+N}=\beta_N}} \left(\prod_{l=1}^{\min(n,k+\ |\lambda|-k)}\frac{1}{l!^{|\alpha_l|}}\right)\frac{1}{\alpha!}\Bigg[$$
$$\sum_{j'=1}^{s} w_{N,j,j'}\tau_{j'}^{k}\left(\prod_{N=1}^{n}\xi_{j',N}^{\sum_{l=1}^{\min(n,|\lambda|)} l\ |\alpha_{l,N}|}\right)-\xi_{j,N}^{N+k+\ |\lambda|}\frac{N!(k+|\lambda|)!}{(k+|\lambda|+N)!}\Bigg]=0$$

Since $\sum\alpha_l=\lambda$ we have $\sum|\alpha_l|=|\lambda|$. If there exists $l>1$ such that $\alpha_l\neq 0_d$ then $\sum l\ |\alpha_l|>|\lambda|$. $\sum l\ |\alpha_l|=|\lambda|$, for all $l>1,\alpha_l=0_{n,d}$, and $\alpha_1=\lambda$, which implies that the sum has a single term, and :

$$\sum_{j'=1}^{s} w_{N,j,j'} = \xi_{j,N}^{N},$$

$$\forall k \in [\![0, v-1]\!],\ \ \forall \lambda \in \mathbb{N}^{[\![1,n]\!]\times[\![1,d]\!]}, \max(1,k) \leq |\lambda| + k \leq v-1$$

$$\sum_{j'=1}^{s} w_{N,j,j'} \tau_{j'}^{k} \left( \prod_{N=1}^{n} \xi_{j',N}^{|\lambda_N|} \right) = \xi_{j,N}^{N+k+\ |\lambda|} \frac{N!(k+|\lambda|)!}{(k+|\lambda|+N)!}$$

The first condition is equivalent to the second one in the case $k = |\lambda| = 0$, we can thus replace $\max(1,k)$ with $k$. We take $d = 1$.

$$\forall k \in [\![0, v-1]\!],\ \ \forall \lambda \in \mathbb{N}^{[\![1,n]\!]}, |\lambda| \leq v-1-k$$

$$\sum_{j'=1}^{s} w_{N,j,j'} \tau_{j'}^{k} \left( \prod_{N=1}^{n} \xi_{j',N}^{\lambda_N} \right) = \xi_{j,N}^{N+k+\ |\lambda|} \frac{N!(k+|\lambda|)!}{(k+|\lambda|+N)!}$$

The generalized solved system conditions are hence satisfied.

Assume the generalized solved system conditions are satisfied, if we use them in the system of equations we get :

$$\sum_{j'=1}^{s} w_{N,j,j'} = \xi_{j,N}^{N},$$

$$\forall k \in [\![0, v-1]\!],\ \ \forall \lambda \in \mathbb{N}^{[\![1,n]\!]\times[\![1,d]\!]}, \max(1,k) \leq |\lambda| + k \leq v-1$$

$$\forall k'' \in [\![k+|\lambda|, v-1]\!], \forall \beta \in \mathbb{N}^{[\![0,n-1]\!]\times[\![1,d]\!]}, |\beta| = |\lambda|$$

$$\sum_{\substack{\alpha \in \mathbb{N}^{[\![1,\, \min\left(n,k''-k\right)]\!]\times[\![1,n]\!]\times[\![1,d]\!]} \\ \sum \alpha_l = \lambda \\ k + \sum l\ |\alpha_l| = k'' \\ l > N \Rightarrow \alpha_{l,N} = 0_d \\ \forall N, \sum_{l=1}^{\min\left(n-N,k''-k\right)} \alpha_{l,l+N} = \beta_N}} \left( \prod_{l=1}^{\min(n,k''-k)} \frac{1}{l!^{|\alpha_l|}} \right) \frac{1}{\alpha!} \left[ \frac{N! \left( k + \sum_{N=1}^{n} \sum_{l=1}^{\min(n,k''-k)} l\ |\alpha_{l,N}| \right)!}{\left( k + \sum_{N=1}^{n} \sum_{l=1}^{\min(n,k''-k)} l\ |\alpha_{l,N}| + N \right)!} - \xi_{j,N}^{N+k''} \frac{N! k''!}{(k''+N)!} \right] = 0$$

Since $\sum_{N=1}^{n} \sum_{l=1}^{\min(n,k''-k)} l\ |\alpha_{l,N}| = k'' - k$ the terms simplify and we get $0 = 0$.